\documentclass[reqno]{memo-l}

\usepackage{amsfonts, amsbsy, amsmath, amsthm, amssymb, latexsym, verbatim, enumerate}
\usepackage{mathrsfs}
\usepackage[top=30mm,right=30mm,bottom=30mm,left=30mm]{geometry}
\usepackage{bm}
\usepackage[pdftex]{color,graphicx}

\makeatletter
\newcommand{\imod}[1]{\allowbreak\mkern4mu({\operator@font mod}\,\,#1)}

\def\@tocline#1#2#3#4#5#6#7{\relax
  \ifnum #1>\c@tocdepth 
  \else
    \par \addpenalty\@secpenalty\addvspace{#2}%
    \begingroup \hyphenpenalty\@M
    \@ifempty{#4}{%
      \@tempdima\csname r@tocindent\number#1\endcsname\relax
    }{%
      \@tempdima#4\relax
    }%
    \parindent\z@ \leftskip#3\relax \advance\leftskip\@tempdima\relax
    \rightskip\@pnumwidth plus4em \parfillskip-\@pnumwidth
    #5\leavevmode\hskip-\@tempdima
      \ifcase #1
       \or\or \hskip 1em \or \hskip 2em \else \hskip 3em \fi%
      #6\nobreak\relax
    \dotfill\hbox to\@pnumwidth{\@tocpagenum{#7}}\par
    \nobreak
    \endgroup
  \fi}

\makeatother

\def\hal{\unskip\nobreak\hfill\penalty50\hskip10pt\hbox{}\nobreak

\hfill\vrule height 5pt width 6pt depth 1pt\par\vskip 2mm}

\usepackage{mathdots}
\usepackage{amssymb,enumerate,graphicx,verbatim,calc}

\def\esixrt#1#2#3#4#5#6{\textstyle{
\hbox to4pt{$\hfil\scriptstyle{#1}\hfil$}
\hbox to4pt{$\hfil\scriptstyle{#3}\hfil$}
\hbox to4pt{$\hfil\scriptstyle{#4}\hfil$}
\hbox to4pt{$\hfil\scriptstyle{#5}\hfil$}
\hbox to4pt{$\hfil\scriptstyle{#6}\hfil$}
\atop
\raise2pt\hbox to4pt{$\hfil\scriptstyle{#2}\hfil$}}}

\def\esevenrt#1#2#3#4#5#6#7{\textstyle{
\hbox to4pt{$\hfil\scriptstyle{#1}\hfil$}
\hbox to4pt{$\hfil\scriptstyle{#3}\hfil$}
\hbox to4pt{$\hfil\scriptstyle{#4}\hfil$}
\hbox to4pt{$\hfil\scriptstyle{#5}\hfil$}
\hbox to4pt{$\hfil\scriptstyle{#6}\hfil$}
\hbox to4pt{$\hfil\scriptstyle{#7}\hfil$}
\atop
\raise2pt\hbox to4pt{$\hfil\scriptstyle{#2}\hfil$}
\raise2pt\hbox to4pt{$\hfil\scriptstyle{\phantom{0}}\hfil$}}}

\def\dfourrt#1#2#3#4{\textstyle{
\raise2pt\hbox to4pt{$\hfil\scriptstyle{#1}\hfil$}
\raise2pt\hbox to4pt{$\hfil\scriptstyle{#2}\hfil$}
{\hbox to4pt{$\hfil\scriptstyle{#3}\hfil$}
\atop
\raise2pt\hbox to4pt{$\hfil\scriptstyle{#4}\hfil$}}}}

\def\a{\alpha}
\def\b{\beta}
\def\e{\epsilon}
\def\d{\delta}

\def\g{\gamma}

\def\Z{\mathbb Z}
\def\C{\mathbb  C}

\def\do{\downarrow}
\def\ra{\rangle}
\def\la{\langle}

\def\Sym{{\rm Sym}}

\def\l{\lambda}

\def\om{\omega}
\def\o{\omega}

\def\pf{\noindent {\bf Proof $\;$ }}

\def\no{\noindent}

\def\hal{\unskip\nobreak\hfil\penalty50\hskip10pt\hbox{}\nobreak
  \hfill\vrule height 5pt width 6pt depth 1pt\par\vskip 2mm}

\newtheorem{theorem}{Theorem}
\newtheorem{thm}{Theorem}[section]
\newtheorem{theor}{Theorem}[chapter]
\newtheorem{prop}[thm]{Proposition}
\newtheorem{lem}[thm]{Lemma}
 \newtheorem{lemma}[theor]{Lemma}
 \newtheorem{propn}[theor]{Proposition}
 
\newtheorem{cor}[thm]{Corollary}
\newtheorem{coroll}[theorem]{Corollary}
 \newtheorem{coro}[theor]{Corollary}

\newtheorem{defn}[thm]{Definition}

\numberwithin{section}{chapter}
\numberwithin{equation}{chapter}

\usepackage{chngcntr}
\counterwithin{table}{chapter}

\begin{document}

\frontmatter

\parskip 1mm

\author{Martin W. Liebeck }
\address{M.W. Liebeck, Imperial College, London SW7 2AZ, UK} 
\email{m.liebeck@imperial.ac.uk}
 
\author{ Gary M. Seitz}
\address{G.M. Seitz, University of Oregon, Eugene, Oregon  97403, USA} 
\email{seitz@uoregon.edu}

\author{Donna M. Testerman}
\address{D.M. Testerman, EPFL, Lausanne,  CH-105  Switzerland}
\email{donna.testerman@epfl.ch}
  

\title{Multiplicity-free representations of algebraic groups}

 \maketitle

  \begin{abstract}
Let $K$ be an algebraically closed field of characteristic zero, and let $G$ be a connected reductive algebraic group over $K$. We address the problem of 
classifying triples $(G,H,V)$, where $H$ is a proper connected subgroup of $G$, and $V$ is a finite-dimensional irreducible $G$-module such that the restriction
of $V$ to $H$ is multiplicity-free -- that is, each of its composition factors appears with multiplicity 1. A great deal of classical work, going back to Dynkin, Howe, Kac, Stembridge, Weyl and others, 
and also more recent work of the authors, can be set in this context. In this paper we determine all such triples in the case where $H$ and $G$ are both simple algebraic groups of type $A$, and $H$ is embedded irreducibly in $G$. While there are a number of interesting familes of such triples $(G,H,V)$, the possibilities for the highest weights of the representations
defining the embeddings $H<G$ and $G<GL(V)$ are very restricted. For example, apart from two exceptional cases, both weights can only have support on at most two fundamental weights; and in many of the examples, 
one or other of the weights corresponds to the alternating or symmetric square of the natural module for either $G$ or $H$. 
 \end{abstract}

  \footnotetext{
Received by the editor: 



\vspace{2mm}
The authors would like to thank the referees for their many comments leading to improvements in the manuscript.

The authors acknowledge the support of the EPSRC Platform Grant  EP/I019111/1, and the  Swiss Science Foundation grants 200021-156583 and 200021-146223. 
The authors also acknowledge the support of the National Science Foundation under Grant No. DMS-1440140 while they were in residence at the Mathematical Sciences Research Institute in Berkeley, California, during the Spring 2018 semester. }

\tableofcontents

\mainmatter



 \chapter{Introduction}\label{intro}

Let $K$ be an algebraically closed field of characteristic zero. A finite-dimensional module for a connected reductive algebraic group over $K$ is said to be {\it multiplicity-free} if each of its composition factors appears with multiplicity 1. There is a great deal of work in the literature over many years that falls under the following rather general program: study triples $(G,H,V)$ satisfying the following properties:
\begin{itemize}
\item[(1)] $G$ is a connected reductive group over $K$, and $V$ is a finite-dimensional irreducible $G$-module such that 
the action of $G$ on $V$ does not contain a full classical group $SL(V)$, $Sp(V)$ or $SO(V)$;
\item[(2)] $H$ is a proper connected reductive subgroup of $G$;
\item[(3)] the restriction $V \downarrow H$ is multiplicity-free. 
\end{itemize}
Note that condition (3) is equivalent to saying that the endomorphism algebra ${\rm End}_H(V)$ is commutative. 

There are a number of instances of well-known work that fall into this program, which we now briefly discuss. 

First, there is an interesting collection of pairs $(G,H)$ with $H<G$ such that $V \downarrow H$ is multiplicity-free for {\it every} irreducible $KG$-module $V$. Under this condition, $H$ is called a {\it multiplicity-free} subgroup of $G$. 
As a consequence of well-known branching rules for the restriction of irreducible representations of $GL_n$ and $SO_n$ (see for example \cite[Chapter 8]{GW}), the pairs $(G,H) = (SL_n,GL_{n-1})$, $(SO_n,SO_{n-1})$ and $(SO_8,\, Spin_7)$ all satisfy this condition. In the case where $G$ is simple, multiplicity-free subgroups have been classified by Kr\"amer in \cite{kramer}, showing that the three pairs in the previous sentence form a complete list of examples. 

Next we mention the celebrated work of Dynkin \cite{dynk} on the maximal subgroups of classical algebraic groups over $\C$: much of this work amounts to classifying the triples $(G,H,V)$ having the property that $H<G$ and $V$ is an irreducible $G$-module such that $V\downarrow H$ is {\it irreducible} (which obviously implies that it is multiplicity-free). 

Our problem also has connections with a classical notion of multiplicity-freeness for an $H$-variety (see \cite[Sec. 12.2]{GW}). If $H$ is a reductive group over $K$, then an affine algebraic $H$-variety $X$ is said to be multiplicity-free if $K[X]$ is multiplicity-free as an $H$-module in the sense we have defined above. When $X$ is itself a finite-dimensional $H$-module, this amounts to saying that the symmetric algebra $S(X^*)$ of the dual $X^*$ is multiplicity-free as an $H$-module. The irreducible $H$-modules $X$ that are multiplicity-free in this sense were classified by Kac in \cite[Theorem 3]{Kac}; this was extended to arbitrary $H$-modules in \cite{BR, Br, Leahy}. One of Kac's examples is $H = GL_n(K)$ with $X = S^2(W)$, the symmetric square of the natural module $W = K^n$ (or its dual). This implies that for the pair $(G,H)$ where 
\[
H = GL_n = GL(W) < G = GL(X) = GL(S^2(W)),
\]
for any integer $k\ge 2$ the irreducible $G$-module $S^k(X)$ has multiplicity-free restriction to $H$. This gives an interesting and nontrivial collection of examples of triples satisfying the above properties (1)--(3). Kac's classification includes as well the family of examples $S^k(X)$ with $X = \wedge^2(W)$, the alternating square of $W$. The multiplicity-freeness of these symmetric algebras is closely related to famous work of Weyl \cite{weyl} on invariant theory; see for example \cite[Chapter 3]{Howe} and \cite[Sec. 5.7]{GW}, where connections with the First Fundamental Theorem are discussed. 

For the examples in the previous paragraph, the irreducible constituents of the restrictions $S^k(X)\downarrow H$ were classified in \cite{HU}, and a detailed discussion of these and other such examples is given in Chapters 3 and 4 of \cite{Howe}, which we shall use as a basic reference at various points in our work.

To avoid any possible confusion, we reiterate that there are two definitions for the multiplicity-freeness of an $H$-module $X$: ours, and the one described above from \cite[Sec. 12.2]{GW}. Our definition is much weaker -- it requires only that the single module $X$ is multiplicity-free, whereas the other definition requires this for a whole collection of modules, namely the symmetric powers of $X$.


Next we discuss some results of Stembridge \cite{stem} on tensor products that can also be interpreted in the framework (1)-(3). Stembridge addresses the following question: if $H$ is a simple algebraic group over $K$, for which pairs $V_1$, $V_2$ of irreducible $H$-modules is the tensor product $V_1\otimes V_2$ multiplicity-free as an $H$-module? (To set this in the above context, take $V = V_1\otimes V_2$ and $G = GL(V_1) \times GL(V_2)$ acting in the natural way on $V$.) Stembridge  solves this problem in \cite{stem}; we shall make great use of his results for the case where $H$ is of type $A_n$ (see Proposition \ref{stem1.1.A}).

Finally we mention our work on overgroups of distinguished unipotent elements \cite{LST}, which also fits into the above program. Let $G$ be a simple algebraic group over $K$, and recall that a unipotent element $u\in G$ is {\it distinguished} if $C_G(u)^0$ is a unipotent group. In $SL_n$, the distinguished unipotent elements are those with a single Jordan block; in $Sp_n$ (resp. $SO_n$) they are the elements with Jordan blocks of distinct even (resp. odd) sizes (see \cite[3.5]{book}). In \cite{LST} the following problem is addressed. Let $Cl(V)$ denote one of the classical groups $SL(V)$, $Sp(V)$, $SO(V)$. What are the connected simple irreducible subgroups $G$ of $Cl(V)$ that contain a distinguished unipotent element $u$ of $Cl(V)$? To set this in the above framework, observe that by the Jacobson-Morozov theorem, there is a subgroup $A \cong A_1$ of $G$ containing $u$. Since $u$ is distinguished in $Cl(V)$, it acts on $V$ with Jordan blocks of distinct sizes, and so the restriction $V \downarrow A$ is multiplicity-free. 

Hence to solve our problem we classified in \cite{LST} the triples $(A,G,V)$ where $A < G < Cl(V)$, $G$ is simple and irreducible on $V$, $A\cong A_1$ is $G$-irreducible (i.e. does not lie in a parabolic subgroup of $G$), and $V \downarrow A$ is multiplicity-free. 

This leads naturally to the problem of classifying triples $(G,H,V)$ where $H<G<Cl(V)$, $G$ is irreducible on $V$, $H$ is a $G$-irreducible subgroup, and $V \downarrow H$ is multiplicity-free. In this paper we solve this problem in the case where all of the groups $G$, $H$ and $Cl(V)$ are of type $A$. Perhaps surprisingly, while there are a number of very interesting families of examples of such triples (including of course those found by Dynkin, Weyl and Howe discussed above), the possibilities are rather restricted, and it turns out to be possible to completely classify them.
Our result will be a fundamental tool for future work on other cases of the general problem.

We now state our main result. This requires a little notation. In the statement, $X = A_{l+1}$ and $Y = A_n$ are simple algebraic groups, and $\d$ and $\l$ are dominant weights for $X$ and $Y$ respectively. We denote by $V_X(\d)$ the irreducible $X$-module of highest weight $\d$, and write $\d = \sum_1^{l+1} d_i\o_i$, where $\o_1,\ldots ,\o_{l+1}$ are fundamental dominant weights for $X$ and $d_i$ are non-negative integers. Similar notation applies for $Y$, writing $V_Y(\l)$ and $\l = \sum_1^n c_i\l_i$. Define $L(\d)$ to be the number of nonzero coefficients $d_i$, and define $L(\l)$  analogously. If $L(\d) = m$, we say that $\d$ has {\it support} on $m$ fundamental weights. If $W = V_X(\d)$ and $Y = SL(W)$, we consider $X$ as embedded in $Y$ via the representation of highest weight $\d$.

\begin{theorem}\label{MAINTHM}
Let $X = A_{l+1}$ with $l\ge 0$, let $W = V_X(\d)$ and $Y = SL(W) = A_n$. Suppose $V = V_Y(\l)$ is an irreducible $Y$-module such that $V \downarrow X$ is multiplicity-free, and assume $\l \ne \l_1,\l_n$ and $\d \ne \o_1,\o_{l+1}$.

Then $\l,\d$ are as in Tables $\ref{TAB1}-\ref{TAB4}$ below, listed up to duals. Conversely, for each possibility in the tables, $V \downarrow X$ is multiplicity-free.
\end{theorem}

Tables \ref{TAB1}--\ref{TAB4} are organized according to the values of $L(\d)$ and $L(\l)$. In the tables, $a$ and $c$ denote positive integers.


\begin{table}
\parbox{.45\linewidth}{
\centering
\begin{tabular}{|l|l|}
\hline
$\l$ & $\d$  \\
\hline
$2\l_1,\,\l_2$ & $\o_1+c\o_i$   \\
                   & $c\o_1+\o_i$  \\
                   & $\o_i+c\o_{i+1}$  \\
                   & $c\o_i+\o_{i+1}$   \\
\hline
$\l_2$ & $2\o_1+2\o_{l+1}$  \\
       & $2\o_1+2\o_2$   \\
       & $\o_2+\o_l$ \\
       & $\o_2 + \o_4$ \\
\hline
$\l_3$ & $\o_1+\o_{l+1}$   \\
\hline
$3\l_1$ & $\o_1+\o_2\,(l=1)$  \\
\hline
\end{tabular}
\vspace{2mm}
\caption{Examples with $L(\d)\ge 2$} \label{TAB1}
}
\hfill
\parbox{.45\linewidth}{
\centering
\begin{tabular}{|l|l|}
\hline
$\l$ & $\d$  \\
\hline
$\l_1+\l_n$ & $c\o_i$  \\
\hline
$\l_1+\l_i \,(2\le i \le 7)$ & $2\o_1,\,\o_2$    \\
$\l_1+\l_{n+2-i}\,(2\le i\le 7)$ &     \\
$\l_2+\l_3$ &  \\
$2\l_1+\l_n$ &   \\
$3\l_1+\l_n$ &   \\
$\l_2+\l_{n-1}$ & \\
$2\l_1+\l_2$ &   \\
$3\l_1+\l_2$ &  \\
\hline
$\l_1+\l_2$ & $3\o_1,\,\o_3$   \\
\hline
\end{tabular}
\vspace{2mm}
\caption{Examples with $L(\d)=1,\,L(\l)\ge 2$ for arbitrary ranks} \label{TAB2}
}
\vspace{1 cm}

\parbox{.45\linewidth}{
\centering
\begin{tabular}{|l|l|}
\hline
$\l$ & $\d$ \\
\hline
$2\l_1,\,\l_2$ & $c\o_i$  \\
\hline
$\l_3$ &  $c\o_1\,(c\le 6)$   \\
       &   $\o_i \,(i\le 6)$ \\
       & $2\o_2$  \\
\hline
$\l_4$ & $c\o_1\,(c\le 4)$  \\
       &   $\o_i \,(i\le 4)$  \\
\hline
$\l_5$ & $c\o_1\,(c\le 3)$  \\
       & $\o_2$   \\
\hline
$\l_i\,(i>5)$ & $2\o_1,\,\o_2$   \\
\hline
$3\l_1$  & $c\o_1\,(c\le 5)$  \\
 & $\o_i \,(i\le 5)$ \\
\hline
$4\l_1$  & $c\o_1\,(c\le 3)$   \\
 &  $\o_i\,(i\le 3)$  \\
\hline
$5\l_1$ & $2\o_1$   \\
         & $\o_i\,(i\le 3)$  \\
\hline
$c\l_1\,(c>5)$ & $2\o_1, \,\o_2$   \\
\hline
$2\l_2,\,3\l_2$ & $2\o_1,\,\o_2$ \\
\hline
\end{tabular}
\vspace{2mm}
\caption{Examples with $L(\d)=1$, $L(\l) = 1$ for arbitrary ranks}\label{TAB3}
}
\hfill
\parbox{.45\linewidth}{
\centering
\begin{tabular}{|l|l|l|}
\hline
$X$ & $\l$ & $\d$  \\
\hline
$A_1$ & $a\l_1+\l_2$ & $2\o_1$   \\
       & $5\l_1$ & $3\o_1$  \\
       & $\l_3$ & $7\o_1$  \\
\hline
$A_3$ & $a\l_i$ & $\o_2$  \\
       & $a\l_i+\l_j$ &  \\
        & $\l_1+\l_2+\l_3$ &  \\
       & $\l_1+\l_2+\l_5$ &  \\
\hline
$A_4$ & $a\l_1+\l_9$ & $\o_2$   \\
       & $a\l_1+\l_2$ &   \\
       & $4\l_2,\,5\l_2$ &   \\
       & $\l_1+2\l_2$ &   \\
       & $2\l_3$ & \\
       & $2\l_4$ & \\
\hline
$A_5$ & $\l_i$ & $\o_3$    \\
        & $\l_1+\l_{18}$ &    \\
\hline
$A_6$ & $\l_5$ & $\o_3$    \\
   & $\l_6$ & \\ 
\hline
$A_7$ & $\l_5$ & $\o_3$ \\
\hline
$A_{13}$ & $\l_3$ & $\o_7$ \\
\hline
\end{tabular}
\vspace{2mm}
\caption{Exceptional examples with $L(\d)=1$ for small ranks}\label{TAB4}
}
\end{table}


We next state some consequences picking out various striking aspects of Theorem \ref{MAINTHM}.

\begin{coroll}\label{cor1}
Assume the hypothesis of Theorem $\ref{MAINTHM}$.
\begin{itemize}
\item[{\rm (i)}] Then $\d$ has support on at most $2$ fundamental weights.
\item[{\rm (ii)}] Also $\l$ has support on at most $2$ fundamental weights, unless $X=A_3$, $\d = \o_2$, $Y=A_5$ and $\l = \l_1+\l_2+\l_3$, $\l_1+\l_2+\l_5$ or the dual of one of these.
\end{itemize}
\end{coroll}

It can be seen from the tables in Theorem \ref{MAINTHM} that many of the multiplicity-free examples have either $\l\in \{\l_2,\,2\l_1\}$ (i.e. $V_Y(\l) = \wedge^2W$ or $S^2W$) or $\d \in \{\o_2,\,2\o_1\}$. The next two corollaries pick out aspects of this phenomenon.

\begin{coroll}\label{cor1.5}
Assume the hypothesis of Theorem $\ref{MAINTHM}$.
\begin{itemize}
\item[{\rm (i)}]  Suppose $\l$ is not $\l_2$, $2\l_1$, $\l_1+\l_n$ or the dual of one of these. Then $\d$ is either $\o_1+\o_{l+1}$ or $c\o_i$ for some $c\le 7$ and some $i$.
\item[{\rm (ii)}] Suppose $\d$ is not $\o_2, 2\o_1$ or the dual of one of these. Then up to duals, $\l$ is $c\l_1\,(c\le 5)$, $\l_i\,(i\le 6)$, $\l_1+\l_n$, $\l_1+\l_2\,(\d=3\o_1 \hbox{ or }\o_3$), or $\l_1+\l_{18}\,(l=4,\d=\o_3)$.
\end{itemize}
\end{coroll}

Another notable fact visible from the tables is that provided $l\ge 4$, the lists of possible weights $\l$ for $\d=\o_2$ and $\d=2\o_1$ are identical. To state this 
formally, we introduce some special notation for dominant weights of $A_n$: 
let $\l_i\,(1\le i\le n)$ be fundamental dominant weights, and for $1\le i < \frac{n+1}{2}$ write $\l_i^*$ instead of $\l_{n+1-i}$ . 
With this notation, write any dominant weight in the form
\begin{equation}\label{nota}
\l = \sum a_i\l_i +  \sum b_i\l_i^*.
\end{equation}

\begin{coroll}\label{equi}
Let $X = A_{l+1}$ with $l\ge 4$, let $W_1=V_X(\o_2)$, $W_2=V_X(2\o_1)$, and embed $X$ in $Y_i = SL(W_i)$ for $i=1,2$. Let $\l$ be a weight written in the form $(\ref{nota})$, 
and consider $\l$ as a weight of both $Y_1$ and $Y_2$. Then $V_{Y_1}(\l)\downarrow X$ is multiplicity-free if and only if  $V_{Y_2}(\l)\downarrow X$ is multiplicity-free.
\end{coroll}

We now describe the layout of the paper, and outline some of the methods used to prove Theorem \ref{MAINTHM}. As in the hypothesis, let $X=A_{l+1}$ and suppose $X<Y = SL(W)$, where $W = V_X(\d)$. Let $V=V_Y(\l)$, where $\l$ is a dominant weight for $Y$. 

In Chapter \ref{notation} we define notation used throughout the paper, and Chapters \ref{levelset} and \ref{cflev} contains basic results needed for ``level analysis" of modules (discussed below), one of our main methods. 

Chapter \ref{MFfam} contains the proofs of the multiplicity-freeness of $V_Y(\l)\downarrow X$ for the weights $\d,\l$ listed in Tables \ref{TAB1}--\ref{TAB4}. Many of the proofs involve very detailed use of the Littlewood-Richardson Rule (Theorem \ref{LR}) for decomposing tensor products of representations, and also of the Carr\'e-Leclerc method of domino tilings (Theorem \ref{dominoes}) which gives the composition factors of the symmetric and alternating parts of the tensor square of a representation. The Littlewood-Richardson and Carr\'e-Leclerc rules are introduced in Chapter \ref{quote}. Special difficulties occur for some of the low rank examples in Table \ref{TAB4}; for example, the case where $X = A_3$ or $A_4$ and $\d=\o_2$ involves substantial effort in Section \ref{taby4}. 

The rest of the paper is devoted to showing that for pairs $\d,\l$ not in Tables \ref{TAB1}-\ref{TAB4}, the restriction $V_Y(\l)\downarrow X$ is not multiplicity-free. 
Our main approach for this is based on what we call ``level analysis". We can choose parabolic subgroups $P_X=Q_XL_X$ of $X$ (with unipotent radical $Q_X$ and Levi factor $L_X$) and $P_Y=Q_YL_Y$ of $Y$, such that $Q_X\le Q_Y$ and $L_X\le L_Y$. For $V = V_Y(\l)$, we define
\[
V^1(Q_Y) = \frac{V}{[V,Q_Y]}, \;\;V^2(Q_Y) = \frac{[V,Q_Y]}{[[V,Q_Y],Q_Y]},
\]
and so on, and call $V^{d+1}(Q_Y)$ the $dth$ level of $V$. We write $V^{d+1}$ to denote the restriction $V^{d+1}(Q_Y)\downarrow L_X'$. 

Now assume that $V\downarrow X$ is multiplicity-free, say $V\downarrow X = V_1+\cdots +V_s$, where each $V_i$ is an irreducible $X$-module. A preliminary result 
(Proposition \ref{induct}) shows that $V^1$ must be multiplicity-free, and $V^1 = \sum_{i\in S} V_i^1(Q_X)$, where $S$ is a certain subset of $\{1,\ldots ,s\}$ and $V^1_i(Q_X)
= V_i/[V_i,Q_X]$ as above. It also gives
\[
V^2 = \sum_{i\in S} V_i^2(Q_X) + M,
\]
where the summand $M$ must be multiplicity-free. The approach is then to attempt to obtain a contradiction by producing several $L_Y'$-composition factors of $V^2(Q_Y)$, 
 restricting them to $L_X'$, and hence finding an $L_X'$-composition factor that must appear more than once in $M$. A sketch of the procedure for finding such composition factors is given in Section \ref{levanal}. It works fairly widely, but not always, and in many cases we have to analyze deeper levels $V^3, V^4,\ldots$ before obtaining a contradiction.

This approach is carried out for the low rank case $X=A_2$ in Chapter \ref{casea2}, and for a series of specific weights $\d$ in Chapters \ref{rkge2}-\ref{o1ol1}. The general proof for arbitrary weights $\d$ is then carried out in Chapters \ref{pfstart}-\ref{pfend}. Along the way, a large amount of technical information about tensor products and symmetric and alternating powers of representations  is required, and this is collected in Chapter \ref{initiallemmas}. 

Throughout the paper, we make substantial use of the Lie-theoretic representation theory packages in Magma \cite{Magma} that enable one to decompose tensor products and  symmetric and alternating powers of representations of groups of type $A$ with given highest weights. Our code is very simple and makes use of just a few commands; so that the reader can reproduce our computations if they wish, we have included the code used at the first point that Magma is quoted in the paper, namely the proof of Proposition \ref{2om12om2}. The Magma code used in the rest of the paper is very similar to this, and we give some comments on this just after the code in Proposition \ref{2om12om2}.


\chapter{Notation}\label{notation}

In this chapter we introduce notation which will be used throughout the paper.
Let $Y = SL(W)$  where $W$ is a finite-dimensional vector space over
an algebraically closed field $K$ of characteristic $0.$ Let $V$ be a nontrivial irreducible $KY$-module of high
weight $\lambda$, so that $V = V_Y(\lambda)$. Denote by $\l^*$ the highest weight of the dual $V^*$. Let $X = A_{l+1}$ and assume that $X < Y = SL(W)$, where 
the restriction $W\downarrow X$ is the irreducible module $V_X(\d)$.  Our goal is to determine all possible $(X,Y,\l,\d)$ for which $V \downarrow X$ is multiplicity-free.  For $X$ of rank $1$ this was
 settled in  \cite{LST}.   For higher rank groups we will use inductive methods based on  parabolic embeddings.

Before describing our specific notation for $X$ and $Y$, we give a general definition of the $S$-{\it value} of a module for a semisimple algebraic group $G$. Let $\mu_1,\ldots ,\mu_n$ be fundamental dominant weights for $G$, and for a dominant weight $\mu = \sum_{i=1}^nc_i\mu_i$, define $S(\mu) = \sum c_i$. For a $G$-module $Z$, set $S(Z)$ to be the maximum of $S(\mu)$ over all irreducible summands $V_G(\mu)$ of $Z$. Also, we shall sometimes denote $\mu$ by the sequence of integers $c_1c_2\ldots$.

\vspace{2mm}
\no \textbf{Notation for $X$. }  We set up notation for $X$ as follows.  Let $P_X = Q_XL_X$ be a maximal parabolic subgroup of $X$ with unipotent radical $Q_X$ and Levi factor $L_X$ of type $A_l$. Write $L_X = L_X'T$, where $T$
is the $1$-dimensional central torus. We assume that $Q_X$ is a product of root subgroups for negative roots.
Let $\Sigma(X)$ denote the root system of $X$ with respect to a 
maximal torus $T_X < L_X$, and let $\Pi(X)$ be a fundamental system in $\Sigma(X)$. Write  
$\Pi(L_X') = \{\alpha_1, \ldots, \alpha_l\}$ and  $\Pi(X) = \Pi(L_X') \cup \{\alpha\}$, where  $\alpha= \alpha_{l+1}$. 
Let $\omega_1, \ldots, \omega_{l+1}$ be the corresponding fundamental dominant weights of $X$.
Write $T_X = S_XT$, where $S_X = T_X \cap L_X'$.  
 For a dominant weight $\mu = \sum_{i = 1}^l s_i\omega_i$ of $L_X'$, write $L(\mu)$ for the number of $i$ such that $s_i\ne 0$.

We have $W\downarrow X = V_X(\d)$, and we write $\d = \sum_{i = 1}^{l+1} d_i\omega_i$.  
With $V = V_Y(\l)$ as above, write $V \downarrow X = V_X(\theta_1) + \cdots + V_X(\theta _s)$  and assume that the restriction is multiplicity-free.  Write $V_i = V_X(\theta_i)$ for $1\le i\le s$.

We shall denote the Lie algebra of an algebraic group $G$ by $L(G)$.
For a positive root $\a$ in $\Sigma(X)$, let $e_\a$ and $f_\a = e_{-\a}$ be corresponding root elements in $L(X)$, and for $1\le i\le l+1$,  let $f_i = f_{\alpha_i + \cdots + \alpha_{l+1}}$. Then $\{f_1, \ldots, f_{l+1} \}$ is a basis of  commuting root elements  of $L(Q_X)$ for negative roots.

\vspace{2mm}
\no \textbf{Notation for $Y$. }  Let $T_Y$ be a fixed maximal torus of $Y,$ $\Sigma(Y)$ the corresponding root system,    $\Pi(Y) = \{\b_1, \ldots, \b_n\}$ a fundamental system of positive roots, and  $ \l_1, \ldots, \l_n $  the corresponding set of fundamental dominant weights for $Y$.  Let $P_Y = Q_YL_Y$ be a parabolic subgroup of $Y$ with  Levi factor $L_Y$ containing $T_Y$, and 
unipotent radical $Q_Y$  a product of root groups for negative roots. Write $L_Y' = C^0 \times \cdots \times C^k$ with each $C^i$ simple. For each $i$,  $\Pi(C^i) =  \la \b_1^i, \ldots, \b_{r_i}^i \ra$ is a string of  fundamental roots of 
$\Pi(Y)$ and we  order the factors so that the string for $\Pi(C^i)$ comes before the string for $\Pi(C^{i+1})$. 
For each $j$ with $1 \le j \le r_i$,  let 
$\l_j^i$ denote the fundamental dominant weight corresponding to $\b_j^i$.
 
 In situations to follow  there will be just
one fundamental root separating  $C^{i-1}$ and  $C^i$ which is designated $\gamma_i$ (see Corollary \ref{levelsnontrivial}). There are also possibly fundamental roots $\g_0$ before $C^0$, and $\g_{k+1}$ after $C^k$.

Finally, if $\l$ is a dominant weight for $Y$,  set $\mu^i = \l \downarrow (T_Y \cap C^i)$, 
so that $V_{C^0}(\mu^0) \otimes \cdots \otimes V_{C^k}(\mu^k)$ is a composition factor of $L_Y$ on $V_Y(\l)$.

\vspace{2mm}
\no \textbf{Notation for Levels. } 
The notion of levels was introduced in \cite{mem1}, and here and in the next section we shall make use of a few basic results from \cite{mem1}: specifically (1.2), (2.3) and (3.6). Although these results are proved there under the assumption of positive characteristic, their proofs make no use of this assumption, and thus the results also hold in our situation of characteristic zero.
Set $[V,Q_Y^0]=V$, define $[V, Q_Y^1]$ to be the commutator space $[V,Q_Y]$, and for  $d>1$  inductively define  $[V, Q_Y^d] = [[V,Q_Y^{d-1}],Q_Y]$.  Now set $V^{d+1}(Q_Y) = [V,Q_Y^d]/[V,Q_Y^{d+1}]$. This quotient is $L_Y$-invariant and will be called the $d$th level.  Thus $V^1(Q_Y) = V/[V, Q_Y]$ is level $0$; 
$V^2(Q_Y) = [V,Q_Y]/[V,Q_Y^2]$ is level $1$, and so on. Similarly, for fixed $i$ and $d \ge 0$ we write  $V_i^{d+1}(Q_X) = [V_i,Q_X^d]/[V_i,Q_X^{d+1}]$, where $V_i$ is as above.   By \cite[(1.2)]{mem1}, $V^1(Q_Y)$ is an irreducible $L_Y'$-module, and $V_i^1(Q_X) = V_i/[V_i,Q_X]$ is an irreducible
$L_X'$-module for each $i$.

  We have $V^1(Q_Y) = V/[V,Q_Y]$ and $V^2(Q_Y) = [V,Q_Y]/[V,Q_Y^2].$  It follows from \cite[(2.3)(ii)]{mem1} that  $V^2(Q_Y) = [V,Q_Y]/[V,Q_Y^2]$ can be regarded as the direct sum of weight spaces 
  of $V$ corresponding to weights of the form $\l -\psi -\gamma_j$, where $\psi$ is a sum of positive roots in $\Sigma (L_Y')$ and
 $\g_j$ is as defined above for $0 \le j \le k+1$.
  Therefore we can write $V^2(Q_Y) = \sum_{j=1}^k V_{\g_j}^2,$ where $V_{\gamma_j}^2 = V_{\gamma_j}^2(Q_Y)$ is the sum of such weight spaces for a fixed value of $j$. This is an $L_Y'$-module, and we set 
$S_j^2 = S(V_{\g_j}^2(Q_Y))$, the $S$-value of this module as defined above.

If  $\mu = \lambda - \sum c_i\b_i$ is a weight  of $V$, then the $Q_Y$-level of $\mu$
is defined to be $\sum c_j$, where the sum ranges over only those $j $ corresponding to fundamental roots in 
$\Pi(Y) \setminus \Pi(L_Y')$.

Next, we introduce an important piece of notation that will be used throughout. Suppose as above that $X = A_{l+1}$ and $X<Y=SL(W)$.
Assume that $P_X = Q_XL_X$ and $P_Y=Q_YL_Y$ are parabolic subgroups of $X$ and $Y$ as above, such that $Q_X \le Q_Y$ and $L_X\le L_Y$. 
With $V$ as above, for each level $V^i(Q_Y)$ we shall write $V^i$ to denote the restriction of the $L_Y$-module $V^i(Q_Y)$ to $L_X'$; that is, 
\begin{equation}\label{videf}
V^i = V^i(Q_Y)\downarrow L_X'.
\end{equation}

Finally, here are a few general pieces of terminology to be used throughout the paper. For a semisimple algebraic group $G$ and a dominant weight $\l$, we shall often denote the module $V_G(\l)$ just by the weight $\l$ -- and if there is any danger of confusion, we shall write $\l \oplus \mu$ for the sum $V_G(\l)+V_G(\mu)$ (rather than just $\l+\mu$, which could refer to the module $V_G(\l+\mu)$). 
If $V$ is an arbitrary $G$-module, we write $V \supseteq \l$ to mean that $V$ has $V_G(\l)$ as a composition factor. 
If we have modules $A = B\oplus C$ we shall write $C = A-B$. And lastly, throughout the rest of the paper we shall abbreviate the term ``multiplicity-free" by the initials MF.

\chapter{Level set-up}\label{levelset}

In this chapter we establish a number of basic  results which will be used throughout the paper.
Notation is as in  Chapter \ref{notation}. In particular, $X = A_{l+1} < Y = SL(W)$ where $W = V_X(\d)$. Also  $P_X = Q_XL_X$ is a maximal parabolic subgroup of $X$ with $L_X = L_X'T$, where $L_X' = A_l$ and $T$ a 1-dimensional 
central torus. Further, we assume that $V = V_Y(\l)$ is such that $V\downarrow X = V_1+\cdots +V_s$ is MF, and $V_i = V_X(\theta_i)$.

\begin{lemma}\label{parabembed}  There is a parabolic subgroup $P_Y = Q_YL_Y$  containing $P_X$ such that the following conditions hold:
\begin{itemize}
\item[{\rm (i)}] $Q_X < Q_Y$;
\item[ {\rm (ii)}] $L_X \le L_Y = C_Y(T)$ and $T_X \le T_Y$, a maximal torus of $L_Y$;
 \item[{\rm (iii)}] $\eta \downarrow T = \alpha \downarrow T$ for each fundamental root $\eta \in  \Pi(Y) \setminus \Pi(L_Y')$;
 \item[{\rm (iv)}] $[W, Q_X^d] = [W,Q_Y^d]$ for $d  \ge 0$.
\end{itemize}
\end{lemma}

\pf  Parts (i)-(iii) follow from (3.6) of \cite{mem1}, and (iv) follows from the proof of that result. \hal

Recall the notation $V^{d+1}(Q_Y) = [V,Q_Y^d]/[V,Q_Y^{d+1}]$, and similar notation for $X$.

\begin{lemma} \label{Twts} {\rm (i)} For $d \ge 0$, $T$ induces scalars on $V^{d+1}(Q_Y)$  via the weight $(\lambda - d\alpha) \downarrow T$.

{\rm (ii)}  For $1\le i \le s$ and $k \ge 0$, $T$ induces scalars on $V_i^{k+1}(Q_X)$  via the weight $(\theta _i - k\alpha) \downarrow T$.

\end{lemma} 

\pf  Here we refer to (2.3)((ii) of \cite{mem1}.  Applying that result to the action of $Y$ on $V$  shows that
$V^{d+1}(Q_Y)$ is isomorphic as a vector space to the direct sum of weight spaces of $V$ having
$Q_Y$-level $d$ (as defined in Chapter \ref{notation}).   As $T$ centralizes $L_Y$, part (i) follows from Lemma \ref{parabembed}(iii). Part (ii) is similar. \hal


\begin{lemma} \label{locateVi}  Let $1 \le i \le s.$
\begin{itemize}
\item[{\rm (i)}] There is a unique $n_i \ge 0$ such that $\theta_i \downarrow T =  (\lambda - n_i \alpha) \downarrow T.$
\item[{\rm (ii)}] $n_i$ is maximal subject to $V_i \le [V,Q_Y^{n_i}].$
\item[{\rm (iii)}] $\left(V_i + [V,Q_Y^{n_i + 1}]\right)/ [V,Q_Y^{n_i + 1}]$ is irreducible under the action of $L_X'$.
\end{itemize}
\end{lemma} 

\pf    Lemma \ref{Twts}(ii) shows that  the weights of $T$
on $V_X(\theta _i)$ have the form $\theta _i - k\alpha$ for $0\le k \le l_i$ where $l_i$ is
maximal among values $j$ with  $[V_i, Q_X^j] \ne 0.$  On the other hand, from Lemma \ref{Twts}(i) we see that all weights
of $T$ on $V$ have form $(\lambda - s\alpha)\downarrow T$ for some $s\ge 0$.  Therefore there is a unique value, say $n_i$,
such that $\theta_i \downarrow T =  (\lambda - n_i \alpha) \downarrow T.$  This gives (i) and (ii).

For (iii), first note that $(V_i + [V,Q_Y^{n_i + 1}])/ [V,Q_Y^{n_i + 1}] \cong V_i/(V_i \cap  [V,Q_Y^{n_i + 1}] )$
as $L_X'$-modules. Now  $Q_X \le Q_Y$ implies that $[V_i, Q_X] \le V_i \cap [V,Q_Y^{n_i + 1}] < V_i.$ The result follows since we know that $V_i/[V_i,Q_X]$ is an irreducible $L_X'$-module. \hal

\begin{lemma}\label{resttoLX'}  Assume $i,j \in \{1, \ldots, s\}$ satisfy $i\ne j$ and $n_i = n_j = m$. Then  
\begin{itemize}
\item[{\rm (i)}] $\theta_i  \downarrow S_X \ne \theta _j \downarrow S_X$;
\item[{\rm (ii)}] $(V_i + [V,Q_Y^{m + 1}])/ [V,Q_Y^{m + 1}] \not  \cong (V_j + [V,Q_Y^{m + 1}])/ [V,Q_Y^{m + 1}]$
as $L_X'$-modules.
\end{itemize}
 \end{lemma}
 
 \pf  Lemma \ref{locateVi} implies that both $(V_i + [V,Q_Y^{m + 1}])/ [V,Q_Y^{m + 1}]$ and  
$(V_j + [V,Q_Y^{m + 1}])/ [V,Q_Y^{m + 1}]$ are irreducible
$L_X'$-summands of $V^{m+1}(Q_Y)$. The highest weights of these summands 
are $\theta_i \downarrow S_X$ and $\theta_j \downarrow S_X$, respectively.  Now $T_X = S_XT$
and by hypothesis and Lemma \ref{locateVi}, $\theta_i \downarrow T = \theta_j \downarrow T$.  
As $\theta_i \ne \theta_j$ we conclude that they have different restrictions to $S_X$ and the
assertions follow. \hal  

We now combine some of the above results to obtain the following result which provides a basis for an inductive
approach to the main theorem, Theorem \ref{MAINTHM}. Recall from (\ref{videf}) our definition 
\[
V^i = V^i(Q_Y)\downarrow L_X'.
\]
For example, $V^1 = \left(V/[V,Q_Y]\right) \downarrow L_X'$.

\begin{propn}\label{induct}  The following assertions hold.
\begin{itemize}
\item[{\rm (i)}] $V^1 =  \sum _{i: n_i = 0} V_i/[V_i,Q_X]$ is MF.
\item[{\rm (ii)}] $V^2 = \sum _{i: n_i = 0} V_i^2(Q_X) + \sum _{i: n_i = 1} V_i/[V_i,Q_X]$.  Moreover, the second sum is MF.
\item[{\rm (iii)}] For any $d\ge 0$,
\[
V^{d+1}  = \sum _{i\,:\,0\le n_i\le d-1} V_i^{d+1-n_i}(Q_X) + \sum _{i\,:\, n_i = d} V_i/[V_i,Q_X],
\]
and the last summand is MF.
\end{itemize}
\end{propn}

\pf  Note that (i) and (ii) are the special cases $d=0,1$ of (iii), so it suffices to prove (iii). 
Fix $d \ge 0$. As mentioned in the proof of Lemma \ref{Twts}, $V^{d+1}(Q_Y) = [V,Q_Y^d]/[V,Q_Y^{d+1}]$ is isomorphic as a vector space to the direct sum of weight spaces of $V$ having
$Q_Y$-level $d.$ By Lemma \ref{Twts}(iii) these weights restrict to $T$ as does $\lambda - d\alpha.$  On the other hand,
$V \downarrow X = \sum V_i$ and for each $i$ the weights of $T$ on $V_i$ have the form $\theta _i - k\alpha$ for  
non-negative integers $k$.  Now we use a weight space comparison and Lemma \ref{locateVi} to complete the proof:  indeed, in order to get  weight $\lambda - d\alpha$ we can start with $\theta_j$ with $n_j \le d$
and take the $T$-weight space for weight $\theta_j - (d-n_j)\alpha$ of $V_j.$  Lemma \ref{Twts} shows that this weight space is $L_X$-isomorphic to $V_j^{d-n_j}(Q_X)$.  Furthermore, the sum of such weight spaces
yields  $V^{d+1}(Q_Y)$.  This together with Lemma \ref{resttoLX'} gives (iii).   \hal

\begin{coro}\label{conseq}  Let $\rho$ be a dominant weight for the maximal torus $S_X$  of $L_X'$.
\begin{itemize}
\item[{\rm (i)}]  If $V_{L_X'}(\rho)$ appears with
multiplicity $k$ in $\sum _{i: n_i = 0} V_i^2(Q_X),$  then $V_{L_X'}(\rho)$ appears with multiplicity
at most $k+1$ in $V^2$.  
\item[{\rm (ii)}] If $V_{L_X'}(\rho)$ appears with
multiplicity $k$ in $\sum _{i: n_i = 0} V_i^3(Q_X) +\sum _{j: n_j= 1} V_j^2(Q_X),$  then $V_{L_X'}(\rho)$ appears with multiplicity at most $k+1$ in $V^3$.
\end{itemize}
\end{coro}

\pf  This follows from Proposition \ref{induct}(ii),(iii). \hal

The next few results concern $S$-values of levels. We extend the definition of the $S$-function introduced in Chapter \ref{notation}, as follows. If  $\rho = \sum s_i\o_i$ is a weight for $X$, then $S(\rho) = \sum s_i$; also for an $X$-module $Z$, we define $S(Z)$ to be the maximal $S$-value over all weights of $Z$. Note that for a dominant weight $\rho$, the $S$-value of $V_X(\rho)$ is equal to $S(\rho)$, since $S(\rho-\sum a_i\a_i) \le S(\rho)$ for any $a_i\ge 0$.

 \begin{lemma}\label{weightsum}     
  {\rm (i)} Suppose that $V_i^j(Q_X)$ is a summand of $V^c$ as in Lemma $\ref{induct}(iii)$.  Then $V_i^{j+1}(Q_X)$ is a summand of 
$V^{c+1}$,  and 
\[
S(V_i^{j+1}(Q_X)) \le S(V_i^j(Q_X)) + 1.
\]
   {\rm (ii)} Suppose $x$ is the maximum $S$-value among irreducibles appearing in $V^c$. 
Then $x+1$ is an upper bound for the $S$-values of highest weights of irreducibles appearing
  within the image of  $[[V,Q_Y^{c-1}],Q_X] $ in $V^{c+1}(Q_Y)$.
    \end{lemma} 
  
\pf Part (ii) follows from (i), so we prove (i). We have $V_i^j(Q_X) =  [V_i, Q_X^{j-1}]/[V_i, Q_X^j]$, which by hypothesis  appears within $V^c(Q_Y)$. Taking commutators with $Q_X,$  $V_i^j(Q_X)$  gives rise to composition factors within $V_i^{j+1}(Q_X) = [V_i, Q_X^j]/[V_i, Q_X^{j+1}] = [V_i, Q_X^j/[[V_i, Q_X^j], Q_X]$, and these appear within 
$V^{c+1}(Q_Y)$. This gives the first assertion of (i).

Any weight of $V_i^{j+1}(Q_X)$ arises from $[[V_i, Q_X^j], Q_X]$, and hence is of the form $\g-\sum_{k=t}^{l+1}\a_k$ for some $t$, where $\g$ is a weight of $V_i^j(Q_X)$. We have $\la -\a_{l+1},\a_l\ra = 1$ and $\la -\a_{l+1},\a_k\ra = 0$ for $k<l$. And for fixed $k,m \le l,$ we have $\la -\alpha_k, \alpha_m\ra  = -2$ if   $m = k$, or 1 if $ m \in \{k+1,k-1\}$,  
and otherwise the inner product is $0$.    The $S$-value assertion in (i) follows.  \hal

The following consequence will be much used in our proofs.

\begin{propn}\label{sbd} Let $d\ge 1$, and suppose $V_{L_X'}(\nu)$ is a composition factor of multiplicity at least $2$ in $V^{d+1}(Q_Y)\downarrow L_X'$. Then 
\[
S(\nu) \le S(V^{d}) + 1.
\]
\end{propn}

\pf As  $V_X(\nu)$ has multiplicity at least $2$, it follows from 
Proposition \ref{induct} that it appears in $V'/[V_j,Q_X]$ for some composition factor $V'$ of $V^{d}$. The result follows, since $S(\nu) \le S(V')+1$ by Lemma \ref{weightsum}(i). \hal

The following general result  will be used to obtain estimates on $S$-values when passing from one level to another.  Note that this result has no multiplicity-free hypothesis.

\begin{lemma}\label{p.29}  Let $G,H$ be simple algebraic groups over $K$ with $G<H$, and let $\delta_1, \delta_2$ be dominant weights for $H$.  Then the following hold.
\begin{itemize}
\item[{\rm (i)}] There exists an irreducible summand $V_G(\nu)$ of $V_H(\delta_2) \downarrow G$ such that $S(V_H(\delta_1 + \delta_2) \downarrow G) \ge S(V_H(\delta_1) \downarrow G) + S(\nu)$.
\item[{\rm (ii)}] If $V_H(\delta_2) \downarrow G$ is irreducible, then 
\[
S(V_H(\delta_1 + \delta_2) \downarrow G) = S(V_H(\delta_1) \downarrow G) + S(V_H(\delta_2) \downarrow G).
\] 
\end{itemize}
\end{lemma}

\pf (i) Embed a Borel subgroup $B_G$ in a Borel subgroup $B_H$ of $H$ with a corresponding embedding of maximal tori $T_G < T_H.$ First note that $V_H(\delta_1 + \delta_2)$ is a direct summand of
  $V_H(\delta_1) \otimes V_H(\delta_2)$ and $A\otimes B$ is a maximal vector for $B_H$, where $A,B$ are maximal vectors for the tensor factors.  Let $V_G(\gamma)$  be an irreducible summand of $V_H(\delta_1) \downarrow G$ such that $\gamma$ has maximal $S$-value.  
There is a maximal vector in $V_H(\delta_1) \downarrow G$ affording $T_G$-weight $\g$, but this vector might not be $A$. However, if it is not then there is a nonzero vector affording $T_G$-weight $\g$ having the form $nA$, where $n$ is a product of terms $f_\b$ for $\b$ a positive root in $L(H)$, and the maximal vector is a linear combination of such
vectors.
Now $B$ is a maximal vector for $G$ so it does  generate an irreducible summand of $V_H(\delta_2) \downarrow G,$  say $V_G(\nu),$ although $\nu$  might not have maximal  $S$-value in this restriction.
 
If $A$ affords $T_G$-weight $\g$, then $A\otimes B$ affords $T_G$-weight $\g+\nu$. If instead we have $n$ as above, then   $n(A\otimes B ) = nA\otimes B + \cdots,$  where  terms following the first term have first coordinate $n_1A$ with $n_1$  a proper subproduct of the terms in $n$; hence these terms are independent of the first term.  The left side is a vector in $V_H(\delta_1+ \delta_2)$ and the first term on the right side is nonzero and affords $T_G$-weight $\gamma + \nu.$  In either case it follows that $\gamma + \nu$ is a dominant weight of $V_H(\delta_1 + \delta_2) \downarrow T_G$, and hence it is subdominant to the highest weight of an irreducible summand for $G.$  However, by considering the effect of subtracting roots we see that a subdominant weight of an irreducible summand has $S$-value at most that of the highest weight.  Therefore,  $S(V_H(\delta_1 + \delta_2) \downarrow G)  \ge  S(\gamma + \nu) = S(\gamma)  + S(\nu)$ and the result follows.  

(ii)  Assume $V_H(\delta_2) \downarrow G$ is irreducible.  Then from (i) we have
\[
\begin{array}{ll}
S(V_H(\delta_1) \downarrow G) + S(\nu) & \le S(V_H(\delta_1 + \delta_2) \downarrow G) \\
   & \le S(V_H(\delta_1 \otimes \delta_2) \downarrow G) \\
    &  = S(V_H(\delta_1) \downarrow G) + S(V_H(\delta_2) \downarrow G) \\
    & = S(V_H(\delta_1) \downarrow G) +S(\nu).
\end{array}
\]
  The assertion follows.  \hal


 
\chapter{Results from the Literature}\label{quote}

In this chapter we describe several results from the literature which will be important in our proofs.
We begin with the well-known result of Littlewood-Richardson on decompositions of tensor products of irreducible representions.  This is followed by a result of Carr\'{e}-Leclerc on the decomposition of the tensor square of a representation into its symmetric and skew symmetric parts.   In addition we quote results
of Stembridge  which provide information of multiplicity-free tensor products.  Finally we present a result
of Cavallin on dimensions of weight spaces.  

\section{Littlewood-Richardson theorem}\label{LRres}

Let $G = GL_n(K)$ and let $V$ be an irreducible  polynomial representation of $G$.  It is well-known that $V$ corresponds to a partition $\rho =(\rho_1, \ldots, \rho_n)$ with at most $n$ nonzero parts (where $\rho_i\ge \rho_{i+1}$ for all $i$). The character of the representation  is given by the Schur function $s_{\rho}$ which is a homogeneous polynomial consisting of the sum of monomials each corresponding to a semi-standard tableau of shape $\rho$; that is, a labelling of
the tableau  by integers $1, \ldots, n$ which is strictly increasing along columns and weakly increasing along rows. Each monomial appearing in $s_{\rho}$ is a product of $| \rho | = \rho_1 + \cdots + \rho_n$ terms.

 Restricting this representation to $X = SL_n(K)$ we obtain the irreducible representation of highest weight  $\l = \sum_{i=1}^{n-1} (\rho_i-\rho_{i+1})\om_i$, where $\o_1,\ldots,\o_{n-1}$ are fundamental dominant weights. 
Of course there are  many possibilities for $\rho$ which yield  $\l$. However, $\rho$ is determined if $\l$ and $| \rho | $ are known.
 
 Going in the other direction, if $\l = \sum_{i=1}^{n-1}c_i\om_i$ is a dominant weight, and we set $c(i) = c_i + \cdots + c_{n-1},$ then $\rho = (c(1), \ldots, c(n-1), 0)$ is one partition which yields $\l$.  In particular, in this way 
 $\om_i$ corresponds to the partition $(1, \ldots, 1, 0, \ldots, 0)$, where 1 appears  $i$ times.
 
 We will state the Littlewood-Richardson theorem in terms of Schur functions. A reference for this is \cite{HL}. We first require some  terminology and notation.  Given $\rho = (\rho_1, \ldots, \rho_n)$  as above, the {\it weight} of $\rho$ is the
 sequence $(1^{\rho_1}, 2^{\rho_2}, \ldots).$ 
 
 Now suppose $\e$ and $\nu$ are partitions.  If the tableau of shape
 $\e$ is an initial part, row by row,  of the tableau of shape $\nu$, then we can form $\nu/\e$.  Here
 we leave blank the cells corresponding to the $\e$ tableau and consider  labellings of the remaining
 cells. 
 
 A labelling of $\nu/\e$ is a {\it Littlewood-Richardson skew tableau} if the following conditions are satisfied:
 \begin{itemize}
\item[1.] Ignoring blank cells, within each row the labels are weakly increasing.
 \item[2.] Ignoring blank cells,  within each column the labels are strictly increasing. 
 \item[3.]  In addition the Y-condition (short for Yamanouchi) is  satisfied. 
\end{itemize}
   The Y-condition means  that the sequence obtained by adjoining the non-blank labels of the reversed rows (top to bottom) has the property that in every initial part of the sequence the number $i$ occurs at least as often as $i+1$.  We note that in particular, this forces all  non-blank cells in the first row of $\nu/\e$ to have label 1. Define the weight of $\nu/\e$ as above: it is $(1^{x_1},2^{x_2},\ldots)$, where $x_i$ is the number of $i$'s in the labelling.

\begin{thm}\label{LR}  Let $\d$ and $\e$ be partitions with at most $n$ nonzero parts. Then 
\[
s_{\d}s_{\e} =  \sum_{\nu} c_{\d,\e}^{\nu}s_{\nu},
\]
where $c_{\d,\e}^{\nu}$ is equal to the number of Littlewood-Richardson skew tableaux of shape $\nu/\e$ and weight equal to that of $\d$.
\end{thm}

Note that one requirement for $c_{\d,\e}^{\nu}\ne 0$ is that $|\nu| = |\e| + |\d|.$ 
\medskip 
 
Here is an  example.  Let $n = 3$ and consider the tensor product of the irreducible $SL_3$ modules of
highest weights $23$ and $22$.  These can arise from the partitions $\d = (5,3,0)$ and $\e = (4,2,0)$, respectively.  To decompose $23 \otimes 22$ we instead consider $s_{\d}s_{\e}.$  We claim
that the irreducible with highest weight $23$ appears with multiplicity 3.  By the above this could arise
from any of the partitions $(5,3,0), (6,4,1), (7,5,2), \ldots.$  However, $|\d| + |\e| = 14$, so the
only possibility is $\nu = (7,5,2).$  We list below the three possible labellings of $\nu/\e$ that satisfy
the conditions. 

\vspace{4mm}

\begin{picture}(300,30)

\put(0,0){\framebox(10,10){2}}
\put(10,0){\framebox(10,10){2}}
\put(0,10){\framebox(10,10){}}
\put(10,10){\framebox(10,10){}}
\put(20,10){\framebox(10,10){1}}
\put(30,10){\framebox(10,10){1}}
\put(40,10){\framebox(10,10){2}}
\put(0,20){\framebox(10,10){}}
\put(10,20){\framebox(10,10){}}
\put(20,20){\framebox(10,10){}}
\put(30,20){\framebox(10,10){}}
\put(40,20){\framebox(10,10){1}}
\put(50,20){\framebox(10,10){1}}
\put(60,20){\framebox(10,10){1}}


\put(100,0){\framebox(10,10){1}}
\put(110,0){\framebox(10,10){2}}
\put(100,10){\framebox(10,10){}}
\put(110,10){\framebox(10,10){}}
\put(120,10){\framebox(10,10){1}}
\put(130,10){\framebox(10,10){2}}
\put(140,10){\framebox(10,10){2}}
\put(100,20){\framebox(10,10){}}
\put(110,20){\framebox(10,10){}}
\put(120,20){\framebox(10,10){}}
\put(130,20){\framebox(10,10){}}
\put(140,20){\framebox(10,10){1}}
\put(150,20){\framebox(10,10){1}}
\put(160,20){\framebox(10,10){1}}


\put(200,0){\framebox(10,10){1}}
\put(210,0){\framebox(10,10){1}}
\put(200,10){\framebox(10,10){}}
\put(210,10){\framebox(10,10){}}
\put(220,10){\framebox(10,10){2}}
\put(230,10){\framebox(10,10){2}}
\put(240,10){\framebox(10,10){2}}
\put(200,20){\framebox(10,10){}}
\put(210,20){\framebox(10,10){}}
\put(220,20){\framebox(10,10){}}
\put(230,20){\framebox(10,10){}}
\put(240,20){\framebox(10,10){1}}
\put(250,20){\framebox(10,10){1}}
\put(260,20){\framebox(10,10){1}}

\end{picture}

\vspace{4mm}

The following corollaries of the theorem will be useful.

\begin{cor}\label{tensorwithoml}  Assume that $X = SL_n$  and  let $V_X(\delta)$ be irreducible with highest weight
$\delta = \sum_{i=1}^{n-1} c_i\omega_i$.  Then $V_X(\d) \otimes V_X(\om_{n-1})$ is MF, and the coefficients of the highest weights of the composition factors are as follows:
\begin{itemize}
\item[{\rm (i)}] $(c_1,\ldots, c_{i-2}, c_{i-1}+1, c_i-1, c_{i+1}, \ldots,  c_{n-1})$ for $1<i \le n-1$,
\item[{\rm (ii)}] $(c_1-1,c_2, \ldots,  c_{n-1})$,
\item[{\rm (iii)}] $(c_1, \ldots,  c_{n-2}, c_{n-1}+1)$,
 \end{itemize}
although  the term in $(i)$ (respectively  $(ii)$) does not occur  if $c_i = 0$ (respectively  $c_1 = 0)$. 
\end{cor}

\begin{cor}\label{LR_om2} Let $X = SL_n$, $n\geq 4$, and let $V_X(\delta)$ be irreducible with
 highest weight
$\delta=\sum_{i=1}^{n-1}c_i\omega_i$. Then $V_X(\delta)\otimes V_X(\omega_2)$ is MF and the 
highest weights of the summands are of the form $\delta+\omega_2-\beta$, for some $\beta\in\Z\Phi^+$. The highest weights and the corresponding $\beta$ are given below:
\begin{itemize}
\item[{\rm (i)}] $(c_1,c_2+1,c_3,\dots,c_{n-1})$; here $\beta=0$,
\item[{\rm (ii)}] $(c_1+1,c_2,\dots,c_{j-1},c_j-1,c_{j+1}+1,c_{j+2},\dots,c_{n-1})$, for 
$2\leq j\leq n-1$ and if  $c_j\ne 0$; here 
$\beta = \alpha_2+\cdots+\alpha_j$,
\item[{\rm (iii)}]  $(c_1-1,c_2,c_3+1,c_4,\dots,c_{n-1})$, if $c_1\ne 0$; here $\beta=\alpha_1+\alpha_2$,
\item[{\rm (iv)}]  $(c_1-1,c_2+1,c_3,\dots,c_{j-1},c_j-1,c_{j+1}+1,c_{j+2},\dots,c_{n-1})$, for $3\leq j\leq n-1$ and if $c_1c_j\ne 0$; here $\beta=\alpha_1+\alpha_2+\cdots+\alpha_j$,
\item[{\rm (v)}] $(c_1,c_2,\dots,c_{j-2},c_{j-1}-1,c_j,c_{j+1}+1,c_{j+2},\dots,c_{n-1})$ for $3\leq j\leq n-1$ and if $c_{j-1}\ne 0$; here $\beta=\alpha_1+2\alpha_2+\cdots+2\alpha_{j-1}+\alpha_j$,
\item[{\rm (vi)}] $(c_1,c_2,\dots,c_{i-2},c_{i-1}-1,c_i+1,c_{i+1},\dots,c_{j-1},c_j-1,c_{j+1}+1,c_{j+2},\dots c_{n-1})$ for $3\leq i<j\leq n-1$ and if $c_{i-1}c_j\ne 0$; here $\beta=\alpha_1+2\alpha_2+\cdots+2\alpha_{i-1}+\alpha_i+\cdots+\alpha_j$.
\end{itemize}
\end{cor}

\pf This follows by using Theorem \ref{LR}, noting that the partition for $\omega_2$ is $(1,1,0,\dots,0)$,  of weight $(1^1,2^1)$. One checks that the conditions on $c_i$ correspond precisely to the ways in which two cells labelled $1$ and $2$ can be added to the ``empty'' tableau corresponding to $\delta$. 

Another proof is to use \cite[15.25]{fulthar}. Here we note that the above weights
correspond to the choice of cardinality $2$ subsets of $\{1,2,\dots,n\}$:
 $\{1,2\}$, $\{1,j+1\}$, $\{2,3\}$, $\{2,j+1\}$, $\{j,j+1\}$, respectively $\{i,j+1\}$.
\hal

The following result  is called the Pieri formula. It is a consequence of Theorem \ref{LR} (see  the discussion in 
\cite[15.25(i)]{fulthar}). 

For convenience we write $(a_1,\ldots ,a_{n-1})$ for the irreducible $SL_n$-module of highest weight $\sum_1^{n-1}a_i\o_i$.

\begin{prop}\label{pieri} Let $X = SL_n.$ Then 
\[
(k,0,\ldots ,0) \otimes (a_1, \ldots, a_{n-1}) = \sum (b_1, \ldots, b_{n-1})
\]
is MF, where the sum is over all  weights $(b_1, \ldots, b_{n-1})$ for which  there exist
non-negative integers $c_1, \ldots, c_n$ such that $\sum c_i = k$, $c_{i+1} \le a_i$ 
and $b_i = a_i + c_i - c_{i+1}$ for $1\le i\le n-1$.
\end{prop}

\section{Decomposing the tensor square}\label{decom}

In this subsection we describe results of Carr\'{e}-Leclerc \cite{domi} on the decomposition of the tensor square of a representation into its symmetric and skew-symmetric parts.  In other words, suppose $V$ is an irreducible
module for  $X = SL_n.$  Then $V \otimes V = S^2(V)+ \wedge^2(V)$ and the goal is to decompose
these summands into  irreducible summands under the action of $X.$

Here too we will consider labellings of certain tableaux, but the situation is slightly different from the Littlewood-Richardson  considerations above.  Here there is just one module $V$, say of highest weight $\l$.  We are going to describe labellings of tableaux by $1\times 2$ and $2\times 1$ arrays, which we refer to as dominoes.

To be consistent with  \cite{domi} we will consider tableaux that are the opposite of those in the previous
subsection.  That is, we use tableaux where the longest rows are at the base.  

A {\it domino  tableau}  T is a tiling of a given tableau by dominoes, as above, such that the labels
are weakly increasing along rows and strictly increasing along columns (from bottom to top).  
(Here in each row or column we take the dominoes lying wholly or partly within it, 
and read the label of each such domino only once.)
We evaluate the weight of the
tiling by counting the number of 1's, 2's, etc.  As we are working with $SL_n$ we want the weight to correspond to  a partition of $n,$ so we also require that all labels are at most $n$.

The {\it column reading} of T is the sequence of numbers  $w_1 \ldots w_r$ obtained by reading the successive columns from top to bottom and left to right.  In this reading horizontal dominoes which belong to two columns are read only once in the leftmost column.  

Here again we have a sort of opposite version of the Y-condition in the previous section.  Namely, the word $w_1 \ldots w_r$ satisfies the Y-condition if for each $i$ and $j$  the sequence $w_i \ldots w_r$ has at least as many $j$'s as $j+1$'s.   In this case we say T is a {\it Yamanouchi domino tableau}.  If this condition holds, then the weight of T  corresponds to a partition.

We can now state the main result.  Suppose $V$ has highest weight $\l$ and let $\rho = (\rho_1, \ldots, \rho_{n-1},0)$ be a corresponding partition. Consider the partition $2\rho = (2\rho_1, 2\rho_1, 2\rho_2, 2\rho_2, \ldots, 2\rho_{n-1},2\rho_{n-1})$   and form the corresponding dual tableau.   

\begin{thm}\label{dominoes}  With notation as above, $S^2(s_{\rho}) = \sum_J a_Js_J$ and $\wedge^2(s_{\rho}) = \sum_J b_Js_J$, 
where $a_J$ (respectively $b_J$) is the number of Yamanouchi domino tableaux of shape $2\rho$ and weight
corresponding to J, whose number of horizontal dominoes is a multiple of $4$ (respectively not  a multiple of $4$).
\end{thm}

As an example let $X = SL_3$ and consider $31 \otimes 31$.  We will use the theorem to show
that the irreducible of highest weight 13 appears with multiplicity 2 in the tensor product with a summand in each of $S^2(31)$ and
$\wedge^2(31).$  Here we let $\rho = (4,1,0)$ so that   $2\rho = (8,8,2,2)$ and the corresponding
tableau has $20$ cells.  Therefore, we look for Yamanouchi domino tableaux with weight $(1^5,2^4,3^1).$
One checks that there are precisely two of these which are illustrated below.

\vspace{4mm}

\begin{picture}(200,40)

\put(0,0){\framebox(10,20){1}}
\put(0,20){\framebox(10,20){2}}
\put(10,0){\framebox(10,20){1}}
\put(10,20){\framebox(10,20){3}}
\put(20,0){\framebox(20,10){1}}
\put(20,10){\framebox(20,10){2}}
\put(40,0){\framebox(20,10){1}}
\put(40,10){\framebox(20,10){2}}
\put(60,0){\framebox(20,10){1}}
\put(60,10){\framebox(20,10){2}}


\put(120,0){\framebox(10,20){1}}
\put(130,0){\framebox(10,20){1}}
\put(120,20){\framebox(20,10){2}}
\put(120,30){\framebox(20,10){3}}
\put(140,0){\framebox(20,10){1}}
\put(140,10){\framebox(20,10){2}}
\put(160,0){\framebox(20,10){1}}
\put(160,10){\framebox(20,10){2}}
\put(180,0){\framebox(20,10){1}}
\put(180,10){\framebox(20,10){2}}

\end{picture}

\vspace{4mm}
\noindent The first of these corresponds to a copy of $13$ in $\wedge^2(31)$ while the second corresponds to
a copy of $13$ in $S^2(31).$

\section{Results of Stembridge and Cavallin}

We begin this subsection by stating  two results that follow from work of Stembridge \cite{stem}.  The first is a general result concerning 
conditions under which  the tensor product of two irreducible modules is MF.  The second provides more detailed information.


\begin{prop}\label{tensorprodMF}  Assume that $D$ is a simple algebraic group over $K$, and $\mu,\nu$ are dominant weights such that $V_D(\mu)\otimes V_D(\nu)$ is MF.  Then either $\mu$ or $\nu$ is a multiple of a fundamental dominant weight.
\end{prop} 

Fot the next result we need the following notation.  For dominant weights
$\mu$ and $\nu$ we write $\mu << \nu$ if $\nu - \mu$ is dominant.  Note that this is very different 
from  saying that $\mu$ is sub-dominant to $\nu$.


\begin{prop}\label{stem1.1.A} Assume $D = A_n$ with fundamental dominant weights $\{ \omega_1, \ldots, \omega_n\}$.  Then $V_D(\mu)\otimes V_D(\nu)$ is MF if and only if for some integers $m > 0$ and $1 \le i,j,k \le n$, one of the following holds (interchanging $\mu$ and $\nu$ if necessary):
\begin{itemize}
\item[{\rm (i)}] $\mu = 0$,  $\omega_i$, $m\omega_1$, or $m\omega_n$;
\item[{\rm (ii)}] $\mu = 2\omega_i$ and $\nu << m\omega_j + m\omega_k$;
\item[{\rm (iii)}] $\mu <<m\omega_2$ or $m\omega_{n-1}$, and $\nu << m\omega_j+m\omega_k$;
\item[{\rm (iv)}] $\mu << m\omega_i$ and $\nu <<\omega_j + m\omega_k$;
\item[{\rm (v)}] $\mu <<m\omega_i$, $\nu << m\omega_j + m\omega_k$ and $k \in \{1, j+1, n\}$.
\end{itemize}
\end{prop}

Next we state two useful results of Cavallin on  weight space dimensions
for representations of  Lie algebras or algebraic groups over the complex numbers.
Let $G$ be of type $A_n$ with fundamental system of roots $\{\a_1, \ldots, \a_n\}$ and let $\l = a_1\l_1 + \cdots + \a_n\l_n$ be  the highest weight of an irreducible representation of $G$. For a weight $\mu$, let $m_\l(\mu)$ be the multiplicity of $\mu$ as a weight of $V_G(\l)$.

The first result is \cite[Proposition 1]{cavallin}.

\begin{prop}\label{cav1}  Let $\mu = \l - \sum_1^nc_i\a_i$ where each $c_i$ is a non-negative integer.
Assume $J$ is a nonempty subset of $\{1, \ldots, n\}$ such that $c_j \le a_j$ for each $j \in J$.
Set $\l' = \l + \sum_{j\in J}(c_j-a_j)\l_j$ and $\mu' = \l' - (\l-\mu)$. Then 
\[
m_{\l}(\mu) = m_{\l'}(\mu').
\] 
\end{prop}

The following result is \cite[Proposition 3]{cavallin}.  
Let $I_{\l} = \{i : a_i\ne 0\}$, and write $I_\l = \{r_1,  r_2,  \ldots, r_{N_{\l}}\}$ with $r_1 < \cdots < r_{N_{\l}}$. 

\begin{prop}\label{cav2} With notation as above, let $\mu = \l - (\a_1 + \cdots +\a_n)$.  If $N_{\l} =1$, then $m_{\l}(\mu) = 1$.  Otherwise
\[
m_{\l}(\mu) = \prod_{i =2}^{N_{\l}}(r_i - r_{i-1}+1).
\]
\end{prop}


 \chapter {Composition Factors In Levels} \label{cflev}

In this chapter we establish a result which is closely related to the classical result on branching from
$GL_n$ to $GL_{n-1}$ mentioned in the Introduction.  Here $GL_{n-1}$ is  the stabilizer of vector in the natural module for $GL_n$.   Proofs that  branching here is always MF can be found in \cite[\S 6]{car}, \cite{NM} and \cite[5.4.1]{Howe}.  

In our situation let $X = A_{l+1}$ (an image of $SL_{l+1}$) and as before  let $P_X = Q_XL_X$ be a parabolic subgroup with $L_X = L_X’T$ where  $L_X'= A_l$, $T$ is a 1-dimensional torus and $Q_X$ is generated by negative root groups with respect to a maximal torus $T_X \le L_X$.  We establish a general result (Theorem \ref{LEVELS}) on the composition factors appearing in the various  $Q_X$-levels for a module $W = V_X(\delta)$, where $\delta$ is a dominant weight and \[ \d = d_1\omega_1 + \cdots +d_l\omega_l + d_{l+1}\omega_{l+1}.
\]
This result and the techniques used in the proof are tailored to our level setup and will play an important role in many arguments to follow. The result is clearly related to the classical results mentioned above but  our situation and the applications are  sufficiently different that we give a proof for completeness. 

In  Section \ref{leva2} we illustrate the main result for the case $X = A_2.$


\section{The main result on levels} \label{mainreslev}

Let $B_X$ be the Borel subgroup of $X$ generated by $T_X$ and positive root groups, 
and let $v$ a maximal vector for the action of $B_X$ on $W = V_X(\d)$.
In addition let $N$ denote the algebra of endomorphisms of  $W$ generated by  the actions of negative root elements $f_{\beta}$ for $\beta \in  \Sigma (L_X')^+$.  Then $W^1(Q_X) = Nv $ is the irreducible $L_X'$-module with highest weight $d_1\omega_1 + \cdots + d_l\omega_l$.

Recall from Chapter 2 that  for each $1\le i\le l+1$, we  let 
\[
f_i = f_{\alpha_i + \cdots + \alpha_{l+1}},
\]
and that $\{f_1, \ldots, f_{l+1} \}$ is a basis of  commuting elements of $L(Q_X)$.  
For $x > 0$, the $x$th level  $W^{x+1}(Q_X)$ is the sum of (images modulo $[W,Q_X^{x+1}]$ of) the 
 spaces $(f_1)^{a_1} \cdots (f_{l+1})^{a_{l+1}}Nv,$  where $a_1 + \cdots + a_{l+1} = x$.
Henceforth we shall refer to such spaces (and the vectors in them) as lying in the  level $W^{x+1}(Q_X)$, with the understanding that 
this really means their images modulo $[W,Q_X^{x+1}]$.

 \begin{thm}\label{LEVELS}  Let $X = A_{l+1}$, $L_X' = A_l$, and let $W=V_X(\delta)$ be the irreducible module for $X$ of highest weight  $\delta = d_1\omega_1 + \cdots +d_l\omega_l + d_{l+1}\omega_{l+1}$.  
\begin{itemize}
\item[{\rm (i)}] For each level $i$, the action of $L_X'$ on $W^{i+1}(Q_X)$ is multiplicity-free, and the irreducible modules that occur have highest weights 
\[
(d_1-a_1+a_2) \omega_1 + \cdots + (d_l-a_l+a_{l+1})\omega_l,
\]
 one for each sequence $(a_1, \ldots, a_{l+1})$ of  integers such that $a_1 + \cdots + a_{l+1} = i$ and $0 \le a_j \le d_j$ for all $j$.
\item[{\rm (ii)}] The highest weights in $(i)$ are afforded by the vectors $(f_1)^{a_1} \cdots (f_{l+1})^{a_{l+1}}v$.
\end{itemize} 
\end{thm}

A more detailed statement than that of part (ii) can be found in Lemma \ref{maxvectors}, giving precise maximal vectors for the composition factors of $W^{i+1}(Q_X)$.

Before beginning the proof we will establish a number of corollaries.

The following special case of Theorem \ref{LEVELS}  gives the composition factors of $W^2(Q_X)$
under the action of $L_X'$.

\begin{cor}\label{V^2(Q_X)}  Let $X = A_{l+1}$, $L_X'=A_l$  and $W=V_X(\delta)$ with 
$\delta = \sum_1^{l+1} d_i\omega_i$, as above. The irreducible $L_X'$-summands of $W^2(Q_X)$  occur with multiplicity $1$, and the coefficients of the corresponding highest weights are as follows:
\begin{itemize}
\item[{\rm (i)}] $(d_1,\ldots, d_{i-2}, d_{i-1}+1, d_i-1, d_{i+1}, \ldots,  d_l)$ for $1<i < l$,
\item[{\rm (ii)}] $(d_1-1,d_2, \ldots,  d_l)$,
\item[{\rm (iii)}] $(d_1, \ldots,  d_{l-1}, d_l+1)$,
 \end{itemize}
although  the terms in $(i),(ii),(iii)$ do not occur  if $d_i = 0$, $d_1 = 0$, $d_{l+1} = 0$, respectively.
\end{cor}

The next corollary shows that the levels described in Theorem \ref{LEVELS} are nontrivial
 modules for $L_X'$, with the possible exception of the first and last level.
 
 \begin{cor}\label{levelsnontrivial}  Let $X = A_{l+1}$, $L_X' = A_l$, and let $W=V_X(\delta)$ be the irreducible module for $X$ of highest weight  $\delta = d_1\omega_1 + \cdots +d_l\omega_l + d_{l+1}\omega_{l+1}$. If  $0 \le i \le d_1+ \cdots +d_{l+1}$, then the action of $L_X'$ on $W^{i+1}(Q_X)$ is nontrivial unless either $ \d = d_{l+1}\om_{l+1}$ with  $i = 0$   or $\d = d_1\om_1$ with $i = d_1$.
\end{cor} 

\pf  Fix a level $i$ and in accordance with Theorem \ref{LEVELS} consider a composition factor
of highest weight 
\[
(d_1-a_1+a_2) \omega_1 + \cdots + (d_l-a_l+a_{l+1})\omega_l,
\]
Suppose this weight is $0$.  Then $d_j - a_j + a_{j+1} = 0$ for $j = 1, \ldots, l$.
Since $d_j \ge a_j$ this forces $a_r = 0$ for $2 \le r \le l+1$ which then implies
$d_r = 0$ for $2 \le r \le l$.  At this point we have shown that $\d = d_1\om_1 + d_{l+1}\om_{l+1},$ and the highest weight of the composition factor is $(d_1 - a_1)\om_1$, 
and so $a_1 = d_1.$
If $d_1 = 0$, then $a_1 = 0 = i$ and we have the first possibility listed in the result.  
Now assume  $d_1 > 0$.  If also $d_{l+1} > 0$ then there is a second composition factor
at level $i$ with highest weight $\om_1 + \om_l$ arising from the  sequence $(a_1-1,0,\ldots , 0,1)$.  This composition factor is nontrivial so that $W^{i+1}(Q_X)$ is
nontrivial. Hence the only other possibility is when
$a_1 = d_1$, $\d = d_1\om_1$ and $i = a_1$, as asserted.  \hal

The following easy consequence of Theorem \ref{LEVELS}  provides information on $S$-values.  In the statement, for $\d = \sum_{i=1}^{l+1}d_i\om_i$, we write $\d' = \sum_{i=1}^{l}d_i\om_i$.

\begin{cor}\label{Svalues}  With notation as in Theorem $\ref{LEVELS} $, the $S$-value of the composition
factor at level $i$ corresponding to the sequence $(a_1, \ldots, a_{l+1})$ is $S(\delta') - a_1 + a_{l+1}$.
Moreover, $S(\delta') + i$ is an upper bound for the $S$-values of all $L_X'$-composition factors at level $i$.
\end{cor}

In the next four corollaries, we write $Y = SL(W)$, $V = V_Y(\l)$, and assume that $V\downarrow X = V_1+\cdots +V_s$ is MF, as in Chapter \ref{levelset}. Let $P_Y=Q_YL_Y$ be the parabolic subgroup of $Y$ given by Lemma \ref{parabembed}, and let $n_i$ be as in Lemma \ref{locateVi}.
Recall our notation 
\[
V^i = V^i(Q_Y) \downarrow L_X'.
\]  
The first corollary shows that part of $V^2$
is ``covered" by a specific tensor product involving $V^1$.  This will
be very useful in the sequel as we will also produce additional summands of  $V^2$.

\begin{cor}\label{cover} {\rm (i)} The $L_X'$-module $\sum_{i:n_i=0}V_i^2(Q_X)$ is isomorphic to a submodule of 
$V^1\otimes V_{L_X'}(\omega_l)$.

 {\rm (ii)} Any irreducible summand of $V^2$ that does not appear in $V^1\otimes V_{L_X'}(\omega_l)$ has multiplicity $1$.
\end{cor}

\pf  (i) Let $J$ be an irreducible $X$-module with highest weight 
$\mu = \sum_{i=1}^{l+1} a_i\omega_i$. Set $\mu' = \sum_{i=1}^l a_i\omega_i$.
Then Corollary~\ref{V^2(Q_X)} 
describes precisely the possible highest weights of the irreducible summands of $J^2(Q_X)$.
Comparing this with the Littlewood-Richardson rule (see Corollary \ref{tensorwithoml}) for decomposing 
$V_{L_X'}(\mu')\otimes V_{L_X'}(\omega_l)$, we see that all of these summands occur in 
$V_{L_X'}(\mu')\otimes V_{L_X'}(\omega_l)$. Hence 
$J^2(Q_X)$ is a submodule of $ V_{L_X'}(\mu')\otimes V_{L_X'}(\omega_l)$.

 Now apply this to each of the irreducible summands $V_i$ of $V\downarrow X$ with $n_i=0$. 
Recall that $V^1 = \sum_{i,n_i=0} V_i^1(Q_X)$ (see Proposition \ref{induct}).
Part (i)  now follows from the previous paragraph. 

Part (ii) is an immediate consequence of (i) together with Proposition \ref{induct}(ii). \hal

The following is a version of the previous result for higher levels.

\begin{cor}\label{icover} With the above notation, for $d \ge 1$,
 \[
V^{d+1}  \subseteq (V^d \otimes V_{L_X'}(\o_l)) + \sum _{i\,:\, n_i = d+1} V_i/[V_i,Q_X]. 
\]
\end{cor}

\pf Write $V \downarrow X = V_1 + \cdots + V_s$ as above.  In view of Proposition \ref{induct} we will work with the irreducible summands $V_j = V_X(\theta_j)$ individually.  Fix $j$ and first assume that $n_j = 0.$ Write $\theta_j = d_1\om_1 + \cdots + d_{l+1}\om_{l+1}$ and consider an irreducible summand $\xi$ of $V_j$ appearing at level $d$.  By Theorem \ref{LEVELS}  there is a sequence $(a_1, \ldots, a_{l+1})$ such that $a_1 + \cdots+ a_{l+1} = d$ and  $\xi$ has highest weight
$(d_1-a_1+a_2) \omega_1 + \cdots + (d_l-a_l+a_{l+1})\omega_l.$  Fix $k$ such that $0 < a_k.$  Then the sequence $(a_1, \ldots ,a_{k-1}, a_k -1, a_{k+1}, \ldots, a_{l+1})$ for $1< k < l+1,$ $(a_1-1, a_2, \ldots, a_{l+1})$ in case $k = 1$, or $(a_1, \ldots, a_l, a_{l+1}-1)$ for $k = l+1$ has terms summing to $d-1.$ The corresponding irreducible summand $\nu$ at level $d-1$ has highest weight which  differs from $\xi$ only for the coefficients of $\om_{k-1}$ and $\om_k$ (just at $\om_1$ if $k = 1$ or just at $\om_l$ if  $k = l+1).$ Indeed if $1 < k < l+1$ then $\nu$ has highest weight
  $$(\ldots, (d_{k-1}-a_{k-1}+(a_k-1))\om_{k-1}, (d_k-a_k+1 +d_{k+1})\om_k \ldots).$$
  And if $k = 1$ or  $k = l+1$ then $\nu$ has highest weight $((d_1-a_1+1)\om_1, \ldots)$ or  $(\ldots, (d_l-a_l+a_{l+1}-1)\om_l)$, respectively.  It now follows from Corollary \ref{tensorwithoml} that $\xi \subseteq \nu \otimes V_{L_X'}(\om_l)$.
  
  To complete the proof we must consider  summands $V_j$ for which $d > n_j > 0.$ Then
Proposition  \ref{induct}(iii) shows that $V_j^{d+1-n_j}$ is one of the summands of $V^{d+1}$.  We therefore use the same argument as above
  where we consider level $d+1-n_j$ of $V_j$ rather than level $d+1.$    The result follows by combining terms from all $V_j$ for which $n_j < d.$ \hal

The next two results compare the $S$-values of successive levels.

\begin{cor}\label{nextlevelSvalue} Let $X = A_{l+1}$, $P_X=Q_XL_X$ and $W = V_X(\d)$ be as above. Then for each integer $d \ge 0$, $S(W^{d+2}(Q_X)) \le S(W^{d+1}(Q_X)) + 1$.
\end{cor}

\pf  In view of Theorem \ref{LEVELS}  and Corollary \ref{Svalues}, for each irreducible summand, say $J,$ of $S(W^{d+2}(Q_X))$ there is a sequence $(a_1, \ldots, a_{l+1})$ such that $\sum a_i = d+1$, $a_j \le d_j$ for each $j$,  and the $S$-value of the irreducible summand is $S(\delta) - a_1 + a_{l+1}.$  By reducing precisely one of the $a_j$ by 1, we obtain a corresponding irreducible summand, say $K$, of  $S(W^{d+1}(Q_X))$.    Then $S(J) \le S(K) + 1$, with
equality holding only if $j = 1.$  \hal

\begin{cor}\label{Ylev}  Let $d\ge 1$. If $J$ is an irreducible $L_X'$-summand of $V^{d+2}$ which occurs with multiplicity at least $2$, then $S(J) \le S(V^{d+1}) + 1$.
\end{cor}

\pf Proposition \ref{induct}(iii) gives a decomposition  of $V^{d+2}$.  Since $J$
occurs with multiplicity at least $2$, it occurs within $V_j^a(Q_X)$ for some $j$ and some $a \ge 2.$  But
then the previous corollary shows that $S(J) \le S(V_j^{a-1}(Q_X)) + 1 \le S(V^{d+1}) + 1$.  \hal

\section{Proof of Theorem \ref{LEVELS} }\label{pflevel}

At this point we begin the proof of Theorem \ref{LEVELS}  with a series of lemmas. Let $X = A_{l+1}$, $P_X=Q_XL_X$ with $L_X'=A_l$, and $W = V_X(\d)$ with $\d = \sum_1^{l+1}d_i\o_i$ be as in the statement of the theorem.

Recall that for $1\le i\le l+1$ we define $f_i = f_{\a_i+\cdots +\a_{l+1}}$. 
We  introduce a total  (lexicographic) order of the monomials $(f_1)^{a_1} \cdots (f_{l+1})^{a_{l+1}}$ 
with all $a_i\ge 0$ and a given value of $\sum a_i$, as follows.  Assuming $a_1 + \cdots + a_{l+1} = b_1 + \cdots + b_{l+1} = i$, we say
 $(f_1)^{a_1} \cdots (f_{l+1})^{a_{l+1}} < (f_1)^{b_1} \cdots (f_{l+1})^{b_{l+1}}$ if $a_{l+1} < b_{l+1};$ or if $a_{l+1} = b_{l+1}$ and $a_l < b_l$; and so on. 

Recall also that $v$ is a maximal vector for the action of the Borel subgroup $B_X$ on $W = V_X(\d)$, and that 
$N$ is the algebra of endomorphisms of  $W$ generated by  the actions of negative root elements $f_{\beta}$ for $\beta \in  \Sigma (L_X')^+$,  so that $W^1(Q_X) = Nv$ is irreducible for $L_X'$ with highest weight $\sum_1^ld_i\o_i$.

 \begin{lem}\label{fivalue}  Let $\eta$ be a weight vector of $W=V_X(\d)$ whose $S_X$-weight is 
 $\sum_1^l c_j\omega_j$. Let $i\ge 1$ with $f_i\eta \ne 0$, and let the $S_X$-weight of $f_i\eta$ be 
$\sum_1^{l} r_j\omega_j$.
 \begin{itemize}
\item[{\rm (i)}] If $i = 1$, then $\sum_1^{l} r_j\omega_j =  \sum_1^{l}c_j\om_j - \omega_1$,  and 
$\sum_1^{l} r_j = \sum_1^{l} c_j - 1$.
 \item[{\rm (ii)}]  If $2 \le i \le l,$ then $\sum_1^{l} r_j\omega_j =  \sum_1^{l} c_j\omega_j +\omega_{i-1} - \omega_i$, and $\sum_1^{l} r_j = \sum _1^{l}c_j$.
 \item[{\rm (iii)}]  If $i = l+1$, then  $\sum_1^{l} r_j\omega_j  = \sum_1^{l} c_j\omega_j +\omega_l$, and    
$\sum_1^{l} r_j =  \sum_1^{l} c_j + 1$.
\end{itemize}
  \end{lem}
  
  \pf  This is immediate from the definition of the root vector $f_i$.  \hal

\begin{lem}\label{fiweight}   Let $a_1,\ldots,a_{l+1}$ be non-negative integers, and write 
$f = (f_1)^{a_1} \cdots (f_{l+1})^{a_{l+1}}$.
\begin{itemize}
\item[{\rm (i)}] Then  $fv$ affords the weight $\sum_{i=1}^l (d_i-a_i+a_{i+1})\o_i$
 of $S_X$.
\item[{\rm (ii)}]  If $fv$ is nonzero and affords a dominant weight of
  $S_X$, then $a_{l+1} \le d_{l+1}$ and $a_i \le d_i+a_{i+1}$ for $1\le i\le l$. 
  \item[{\rm (iii)}]  There are only finitely many monomials $f = (f_1)^{a_1} \cdots (f_{l+1})^{a_{l+1}}$ for which
  the vector $fv$  is nonzero and affords a dominant weight of $S_X$.
  \end{itemize}
  \end{lem}
  
  \pf The statement of (i) is immediate from the definition of the terms $f_i$.  
As $f_{l+1}^{a_{l+1}}v = f_{\a_{l+1}}^{a_{l+1}}v$,
  this vector is nonzero only if $a_{l+1} \le d_{l+1}.$  For the weight of $fv$ to restrict to a dominant weight of $S_X$ it follows from (i) that   $d_j-a_j+a_{j+1} \ge 0$ for $1 \le j \le l.$  This establishes (ii), and (iii) follows.  \hal
  
  \begin{lem}\label{W_i}
   {\rm (i)} Let  $a_1,\ldots,a_{l+1}$ be non-negative integers. Then
\[
(f_1)^{a_1} \cdots (f_{l+1})^{a_{l+1}}Nv \subseteq \sum N(f_1)^{b_1} \cdots (f_{l+1})^{b_{l+1}}v,
\]
 where the sum ranges over  terms $(f_1)^{b_1} \cdots (f_{l+1})^{b_{l+1}} \le (f_1)^{a_1} \cdots (f_{l+1})^{a_{l+1}}$ (in the above ordering) with  $\sum b_j = \sum a_j$.
 
 {\rm (ii)}  For $x\ge 1$, level $x$ can be expressed as 
\[
W^{x+1}(Q_X)  =  \sum N(f_1)^{b_1} \cdots (f_{l+1})^{b_{l+1}}v,
\]
 where the sum  is over all non-negative sequences with $b_1 + \cdots + b_{l+1} = x$.
 \end{lem}
  
  \pf The statement of (ii) follows from (i) since $W^{x+1}(Q_X)$ is the sum of terms of the form  $(f_1)^{a_1} \cdots (f_{l+1})^{a_{l+1}}Nv$ for which $a_1 + \cdots +a_{l+1} = x$. For (i), consider the expression $(f_1)^{a_1} \cdots (f_{l+1})^{a_{l+1}}Nv.$  Bring  terms   $f_{\beta}$ for $\beta \in  \Sigma (L_X')^+$ to the left past  terms $f_k$.  
For a given $\b$ and $k$, either $f_\b$ commutes with $f_k$, or there exists $j<k$ such that $\b = \a_j+\cdots +\a_{k-1}$, in which case $f_kf_{\beta} = f_{\beta}f_k \pm f_j$. As $f_j < f_k$ in the ordering, the assertion follows. \hal
  
 \begin{lem}\label{eifi^r}  Fix $j \le l$ and  $r > 0.$ Then $e_{\alpha_j}(f_j)^r = (f_j)^re_{\alpha_j} + \e r(f_j)^{r-1}f_{j+1},$  where $\e$ satisfies $[e_{\alpha_j},f_j] = \e f_{j+1}.$
\end{lem} 

\pf  The proof is by induction on $r$.  For $r = 1$ this is just the definition of $\e$.  Suppose true for $r$.
Then 
\[
\begin{array}{ll}
e_{\alpha_j}(f_j)^{r+1} & = (e_{\alpha_j}(f_j)^r)f_j \\
& = ((f_j)^re_{\alpha_j} + \e r(f_j)^{r-1}f_{j+1})f_j \\
& = (f_j)^re_{\alpha_j}f_j + \e r(f_j)^rf_{j+1} \\
& = (f_j)^{r+1}e_{\alpha_j} + \e (f_j)^rf_{j+1} + \e r(f_j)^rf_{j+1} \\
& = (f_j)^{r+1}e_{\alpha_j} + \e (r+1)(f_j)^rf_{j+1},
\end{array}
\]
completing the induction.  \hal

Notice that the second monomial on the right side of the equation above is greater than $(f_j)^r$ in the ordering, and the first term yields $0$ when applied to the maximal vector $v$.

Fix $x\ge 1$, and let $w_1 > w_2 > \cdots > w_z$ be the ordering of the finite number of monomials $(f_1)^{a_1} \cdots (f_{l+1})^{a_{l+1}}$ for which $\sum a_i = x$ and $(f_1)^{a_1} \cdots (f_{l+1})^{a_{l+1}}v$ is nonzero and affords a dominant weight of $S_X$. Set  $X_0 = [W,Q_X^{x+1}]$ and for $j\ge 1$,
\[
X_j = \left(\sum_{c \le j}L_X'w_cv\right) + [W,Q_X^{x+1}].
\]

\begin{lem}\label{maxvectors} The following hold.
\begin{itemize}
\item[{\rm (i)}] For each $j \ge 0$, either $w_{j+1}v\in X_j$, or $w_{j+1}v + X_j$ is a maximal vector in $W^{x+1}(Q_X) /X_j$.
\item[{\rm (ii)}]  $W^{x+1}(Q_X) = X_z$.
\end{itemize}
\end{lem}

\pf  (i) To simplify notation we will proceed as if $[W,Q_X^{x+1}] = 0$. We first claim that any $e_\a$ for $\a \in \Sigma (L_X')^+$ annihilates $w_1v.$   Suppose false. Then there is 
a fundamental root $\alpha_i$  for $i \le l$  such that $e_{\a_i}w_1v \ne 0.$ Write $w_1 = (f_1)^{a_1} \cdots (f_{l+1})^{a_{l+1}}.$   If $j \ne i$, then  $e_{\alpha_i}f_j = f_je_{\alpha_i}$ and otherwise 
$e_{\alpha_i}f_i = f_ie_{\alpha_i} \pm f_{i+1}$.  It follows from Lemma \ref{eifi^r} that  $e_{\alpha_i}(f_1)^{a_1} \cdots (f_{l+1})^{a_{l+1}}$ is the sum of   $(f_1)^{a_1} \cdots (f_{l+1})^{a_{l+1}}e_{\alpha_i}$ and a nonzero multiple of   $(f_1)^{b_1} \cdots (f_{l+1})^{b_{l+1}} $, where $b_i = a_i-1, b_{i+1} = a_i+1,$ and otherwise $b_j = a_j, $ 
Then $(f_1)^{b_1} \cdots (f_{l+1})^{b_{l+1}} > (f_1)^{a_1} \cdots (f_{l+1})^{a_{l+1}}$.  Now applying   the sum to the maximal vector $v$ we get  a nonzero multiple of  $(f_1)^{b_1} \cdots (f_{l+1})^{b_{l+1}}v \ne 0.$     Either this is a maximal vector or we can apply $e_{\a_j}$ for some $j \le l$ and get a nonzero vector.  In the first case
 this contradicts the choice of $w_1.$  In the second case apply $e_{\a_j}$ and get a nonzero multiple of $(f_1)^{c_1} \cdots (f_{l+1})^{c_{l+1}}v$, where $(f_1)^{c_1} \cdots (f_{l+1})^{c_{l+1}} > (f_1)^{b_1} \cdots (f_{l+1})^{b_{l+1}} > (f_1)^{a_1} \cdots (f_{l+1})^{a_{l+1}}.$
Continuing this process  we eventually get a maximal vector and contradict the choice of $w_1.$
This shows that $w_1v$ is a maximal vector.
 The same argument  shows that any $e_\a$ ($\a \in \Sigma (L_X')^+$) annihilates $w_2v + X_1$ in $W^{x+1}(Q_X)/X_1$, and continuing we obtain the result.  
 
 (ii) Suppose $W^{x+1}(Q_X)$ properly contains $X_z.$ Then Lemma \ref{W_i}(ii) implies that there exists a monomial  $(f_1)^{b_1} \cdots (f_{l+1})^{b_{l+1}}$ such that such that $b_1 + \cdots + b_{l+1} = x$ and  $(f_1)^{b_1} \cdots (f_{l+1})^{b_{l+1}}v$ is not contained in $X_z$.  But then this vector 
 affords a dominant weight, a contradiction. \hal


We now aim to determine precisely which $w_{j+1}v$ do not lie in $X_j$  in the above lemma.

\begin{lem}\label{a_jled_j} Let $w_{s+1} = (f_1)^{a_1} \cdots (f_{l+1})^{a_{l+1}}$ be one of the monomials in  Lemma $\ref{maxvectors}$ with $a_1 + \cdots + a_{l+1} = x.$  Then  $w_{s+1}v + X_s = 0$ unless
$a_k \le d_k$ for $1\le k \le  l+1.$
\end{lem}

\pf  Assume false, and let $X= A_{l+1}$, $W=V_X(\d)$ be a counterexample of minimal rank and with $S(\d)$ minimal for this rank. Choose $k$ maximal such that $a_k > d_k$ and $w_{s+1}v + X_s \ne 0$.

As in the last proof we simplify notation by working modulo $[W,Q_X^{x+1}]$.  That is, we proceed as
if $[W,Q_X^{x+1}] = 0.$  By definition, we have   
$X_s = \sum_{c \le s}L_X'w_cv.$  We claim that $X_s = \sum_{c \le s}Nw_cv.$   To see
this first note that Lemma \ref{maxvectors} shows that  $X_1 = L_X'w_1v$ is irreducible with maximal vector $w_1v$ and hence it can be
written as $Nw_1v.$  Similarly $X_2 /X_1 = (L_X'w_2v + X_1)/X_1 = (Nw_2v + X_1/X_1)$ and
so  $X_2 = Nw_1v + Nw_2v.$  Continuing in this way we obtain the claim.

Note that $k<l+1$, since otherwise $a_{l+1}>d_{l+1}$, hence 
 $(f_{l+1})^{a_{l+1}}v = (f_{\a_{l+1}})^{a_{l+1}}v = 0$, and so $w_{s+1}v + X_s = 0.$
Suppose $k > 1$.  Then $X$ has rank at least $3$ and we consider the $A_l$ Levi subgroup $S$ with fundamental system $\alpha_2, \ldots, \alpha_{l+1}$ and the irreducible module $Sv$ for this Levi subgroup, which has highest weight $\sum_2^{l+1}d_i\omega_i$.
 By minimality,  the conclusion of the lemma holds for the restriction to the Levi subgroup $L_S' = L_X' \cap S$ with base   $\alpha_2, \ldots, \alpha_l$  at  level $t = \sum_2^{l+1}a_i$.  Our supposition  $w_{s+1}v + X_s \ne 0$ implies that $(f_2)^{a_2} \cdots (f_{l+1})^{a_{l+1}}v \ne 0.$  Therefore $(f_2)^{a_2} \cdots (f_{l+1})^{a_{l+1}}$ is one of the  monomials in Lemma \ref{maxvectors} for level $t$, and again by minimality,
\[
(f_2)^{a_2} \cdots (f_{l+1})^{a_{l+1}}v = \sum e_{b_2, \ldots, b_{l+1}}(f_2)^{b_2}\cdots (f_{l+1})^{b_{l+1}}v, 
\]
 where each of the monomials $(f_2)^{b_2}\cdots (f_{l+1})^{b_{l+1}}$ appearing is larger than $(f_2)^{a_2} \cdots (f_{l+1})^{a_{l+1}}$ in the ordering, and the terms $e_{b_2, \ldots, b_{l+1}}$ are in the algebra generated by$f_\b$ for $\b \in \Sigma (L_S)^+$.  We then obtain a contradiction by multiplying both sides by $(f_1)^{a_1},$  noting that $f_1$ commutes with each of the terms   $e_{b_2, \ldots, b_{l+1}}$. Therefore we now assume $k = 1.$

 Lemma \ref{fiweight}(i) shows that the $S_X$-weight afforded by $(f_1)^{a_1} \cdots (f_{l+1})^{a_{l+1}}v$ is 
$\rho = \sum_1^l (d_i - a_i + a_{i+1})\omega_i$. By Lemma \ref{maxvectors}, in order
for the lemma to be false there must exist an irreducible submodule for $L_X'$ of this highest weight at level $x$.
Recall that the level is determined by the action of $T$. 

We shall use shorthand notation $(d_1, \ldots, d_{l+1})$ for $W = V_X(\d)$ with $\d = \sum_1^{l+1}d_i\o_i$, and similar notation for $L_X'$-modules. Observe that
\begin{equation}\label{tenso}
(d_1, \ldots, d_{l+1}) \subseteq (d_1, 0,\ldots, 0) \otimes (0, d_2, \ldots, d_{l+1}).
\end{equation}
  We will look for an irreducible $L_X'$-submodule of the tensor product, of the above highest weight $\rho$, 
and for which $T$ has the appropriate action. 
 We know that the action of $L_X'$ on the $X$-module  
 $(d_1, 0,\ldots, 0) = S^{d_1}(\om_1)$ has composition factors $(d_1, 0, \ldots, 0)$, $(d_1-1,0,\ldots, 0)$, $\ldots$, 
$(0, \ldots, 0)$, at levels $0,1,\ldots,d_1$ respectively. Fix a  level $a$ in this action, so that $L_X'$ acts as $(d_1-a,0,\ldots,0)$, and consider  level $b = x-a$ for the second tensor factor in (\ref{tenso}).
 By Lemma \ref{maxvectors} and the minimality of $S(\d)$,  the action of $L_X'$ at this level is given as a sum of terms corresponding to monomials $(f_2)^{j_2}\cdots (f_{l+1})^{j_{l+1}}$
 where $b = j_2 + \cdots + j_{l+1},$ which afford irreducibles for $L_X'$
  of highest weights $(j_2, d_2-j_2+j_3, \ldots, d_l-j_l+j_{l+1}).$
 
 Consider a term $(d_1-a,0, \ldots, 0) \otimes (j_2, d_2-j_2+j_3, \ldots, d_l-j_l+j_{l+1})$ and
 apply Pieri's formula \ref{pieri}.  In order for $\rho = (d_1 - a_1 + a_2, \ldots,d_l-a_l+a_{l+1})$ to appear as a summand,  there must exist non-negative integers $c_1, \ldots, c_{l+1}$ whose sum
 is $d_1-a$ and such that the following equalities hold:
 \[
\begin{array}{rl}
j_2 + c_1-c_2 & = d_1-a_1+a_2 \\
 d_2-j_2+j_3 + c_2-c_3 & = d_2-a_2+a_3 \\
 d_3 -j_3 + j_4 + c_3-c_4 & = d_3-a_3+a_4 \\
 & \vdots \\
 d_l-j_l+j_{l+1} + c_l-c_{l+1} & = d_l-a_l+a_{l+1}.
\end{array}
\]
  Call these equations $1,\ldots ,l$, and combine them as follows.  First cancel  the terms $d_i$ appearing on both sides of equations $2, \ldots, l.$
  Add the first two of the resulting equations, then the first three, etc.  Listing the first equation followed by the results of the additions we obtain:
 \[
\begin{array}{rl}
j_2 + c_1-c_2 &= d_1-a_1+a_2 \\
j_3 + c_1-c_3 & = d_1-a_1+a_3 \\
j_4 + c_1-c_4 & = d_1 -a_1+a_4 \\
  & \vdots \\
 j_{l+1}+c_1-c_{l+1} & = d_1-a_1 + a_{l+1}.
 \end{array}
\]
  Add these equations to get 
 \[
 (j_2 + \cdots + j_{l+1}) + lc_1 -(c_2 + \cdots +c_{l+1}) = l(d_1-a_1) + (a_2 + \cdots + a_{l+1}).
\]
 As $b =\sum_{i=2}^{l+1} j_i$, $\sum_1^{l+1}c_i = d_1-a$, and $\sum_1^{l+1}a_i = x$,  this reduces to
\[
b + lc_1 - (d_1-a -c_1) = l(d_1-a_1) + (x-a_1).
\]
And as $a+b = x$, this simplifies to $(l+1)c_1 = (l+1)(d_1-a_1),$  
 which  is a contradiction since $c_1 \ge 0$ while the right side is negative (by our assumption $a_1>d_1$). \hal

 Our goal at this point is to show that each of the monomials $w_{s+1} = (f_1)^{a_1} \cdots (f_{l+1})^{a_{l+1}}$ for which
$a_j \le d_j$ for all $j$, leads to a maximal vector.  Towards this end, consider such a term
$w_{s+1}$ and write $w_{s+1} = (f_1)^{a_1} \cdots (f_k)^{a_k},$  where $k$ is maximal such that $a_k \ne 0.$  Set $e_k = e_{\alpha_k + \cdots + \alpha_{l+1}}$, $h_k = h_{\alpha_k + \cdots + \alpha_{l+1}},$ and for $t\ge i$ let  $f_{i\cdots t} = f_{\alpha_i + \cdots + \alpha_t}.$  For $a_i \ne 0$, 
set 
\[
w_{s+1}/f_i = (f_1)^{a_1} \cdots (f_{i-1})^{a_{i-1}}(f_i)^{a_i-1}(f_{i+1})^{a_{i+1}} \cdots (f_k)^{a_k}.
\]

\begin{lem}\label{e_kw_{s+1}} With notation as above, the following hold.
\begin{itemize}
\item[{\rm (i)}] If $k > j$, then $e_k(f_j)^{a_j} = (f_j)^{a_j}e_k \pm a_j(f_{j \cdots  (k-1)})(f_j)^{a_j-1}.$
\item[{\rm (ii)}]  For $1 \le b \le a_k$ we have $e_k(f_k)^bv  = b\left(\d (h_k) - (b -1)\right)(f_k)^{b-1}v \ne 0.$ 
\item[{\rm (iii)}]We have
\[
e_kw_{s+1}v = c(f_1)^{a_1} \cdots (f_{k-1})^{a_{k-1}}(f_k)^{a_k -1}v + \sum_{i<k, \,a_i \ne 0}c_i
f_{i \cdots (k-1)}(w_{s+1}/f_i)v,
\]
where $c,\,c_i$ are integers and $c \ne 0$.
\end{itemize}
\end{lem}

\pf  The proof of (i) is just as in the proof of Lemma \ref{eifi^r}, 
except that $e_kf_j = f_je_k +  \e f_{\alpha_j + \cdots + \alpha_{k-1}}$  so that the term $f_{j+1}$ in \ref{eifi^r} is replaced by $f_{\alpha_j + \cdots + \alpha_{k-1}} = f_{j \cdots  (k-1)}$, noting that $f_j$ commutes with $f_{j \cdots  (k-1)}$.

For (ii), we argue by induction on $b.$  For $b = 1$,  $e_kf_kv = f_ke_kv +h_kv = \d(h_k)v.$  Note also, that 
$\d(h_k) = a_k+ \cdots +a_{l+1} \ge a_k > 0.$ So $\d(h_k)-(b-1) \ne 0$. 
Now consider the induction step.  Suppose the assertion holds for  $1 \le b < a_k.$
Then by induction we have
\[
\begin{array}{ll}
e_k(f_k)^{b+1}v & = (f_ke_k + h_k)(f_k)^bv \\
 & = f_k\left(b(\d (h_k) - (b -1))\right)(f_k)^{b-1}v + h_k(f_k)^bv \\
 & = b\left(\d (h_k) - (b -1)\right)(f_k)^bv + (\d (h_k) - 2b)(f_k)^bv \\
 &  = (b+1)( \d(h_k) -b)(f_k)^bv,
\end{array}
\]
 completing the induction. 
  
For (iii), consider $e_kw_{s+1}v = e_k(f_1)^{a_1} \cdots (f_k)^{a_k}v.$  We first compute $e_k(f_1)^{a_1} \cdots (f_{k-1})^{a_{k-1}}$ by  using  (i) repeatedly to bring $e_k$
 past the terms $(f_1)^{a_1}, \ldots, (f_{k-1})^{a_{k-1}}$ in turn.  The result has the form 
\[
 (f_1)^{a_1} \cdots (f_{k-1})^{a_{k-1}}e_k + \sum_{i<k}c_i
f_{i \cdots (k-1)}(f_1)^{a_1} \cdots (f_{i-1})^{a_{i-1}}(f_i)^{a_i-1}(f_{i+1})^{a_{i+1}} \cdots (f_{k-1})^{a_{k-1}}.
\]
 Now apply this to $(f_k)^{a_k}v. $ The conclusion follows from (ii).  \hal

\begin{lem}\label{nonzero}  Suppose $a_1 + \cdots +  a_{l+1} = x$ and $w_{s+1} = (f_1)^{a_1} \cdots (f_{l+1})^{a_{l+1}}$ with  each $a_j \le d_j$.  Then $w_{s+1}v \not \in X_s.$
\end{lem}

\pf  Suppose false, and among all such representations of $X=A_{l+1}$ where the assertion fails, choose $W = V_X(\d)$ such that 
the number of nonzero labels $d_h$  in $\d = \sum d_h\o_h$ is minimal.  As in other lemmas there is no harm in proceeding as if
$[W,Q_X^{x+1}] = 0.$

Assume first that there is just one nonzero $d_h.$ In this case at each level $x \le d_h$ there is just one associated nonzero monomial, namely  $f_h^x.$  So here it is only necessary to verify that $f_h^xv \ne 0.$  For this, first note that   $x \le d_h$ implies $(f_{\a_h})^xv \ne 0$ and therefore 
$f_h^xv = \pm ((f_{\a_h})^xv)^{s_{h+1} \cdots s_{l+1}} \ne 0$  as well (where $s_i$ is the reflection in the fundamental root $\a_i$). Therefore  from now on we assume there are at least two  nonzero labels. 

Let $x$  be minimal such that the assertion is false at level $x$ for $W$.  And at level $x$ choose the monomial $(f_1)^{a_1} \cdots (f_{l+1})^{a_{l+1}}$ to be as large as possible in the ordering such that the assertion fails.  That is, we take $s$ minimal. Set $k$  maximal with $a_k \ne 0$, and write 
\begin{equation}\label{disp}
w_{s+1}v = (f_1)^{a_1} \cdots (f_k)^{a_k}v = \sum n_{c_1 \cdots c_{l+1}}(f_1)^{c_1} \cdots (f_{l+1})^{c_{l+1}}v,
\end{equation}
 where each monomial  $(f_1)^{c_1} \cdots (f_{l+1})^{c_{l+1}} > 
(f_1)^{a_1} \cdots (f_k)^{a_k}$ in the ordering, and $n_{c_1 \cdots c_{l+1}}$ is in the algebra generated by $f_\b$ for $\b \in\Sigma (L_X')^+$. 
It follows that either $(f_1)^{c_1} \cdots (f_{l+1})^{c_{l+1}} = (f_1)^{c_1} \cdots (f_k)^{c_k}$ with $c_k \ge a_k$,
or $c_j > 0$ for some $j > k.$
All monomials appearing in (\ref{disp}) are at level $x.$  A weight space comparison shows that  terms on the right side which afford a weight  not equal to that afforded by $w_{s+1}v$ must sum to zero, so we delete all such terms.

{\bf Case (I)}
First assume that $c_j = 0$ for all $j > k$, for each of the nonzero terms on the right side of (\ref{disp}).  
This forces   $x > 1$, as
otherwise $w_{s+1}v = f_kv \ne 0$ but there is no possible nonzero term on the right side.   
Note that the assumption covers the possibility that $w_{s+1}v = 0$ and all terms on the right side are $0$. Therefore, all  nonzero terms  on the right side have the form  $n_{c_1 \cdots c_k}(f_1)^{c_1} \cdots (f_k)^{c_k} $ with $c_k \ge a_k$. Suppose some $c_k = a_k$ and choose $j$ maximal such that $c_j \ne a_j.$ In order for such a monomial to be larger than $(f_1)^{a_1} \cdots (f_k)^{a_k}$  we must have $c_j > a_j.$
Moreover, the equality $x = a_1 + \cdots + a_k = c_1 + \cdots + c_k$ implies that there exists $b< j < k$ such that $a_b > c_b.$  

Using  Lemma \ref{e_kw_{s+1}}, we apply $e_k$ to both sides of (\ref{disp}). Expanding $e_k(f_1)^{a_1} \cdots (f_k)^{a_k}v$ we get a linear combination of terms
such that the term with the smallest monomial is a   multiple of \linebreak 
 $(f_1)^{a_1} \cdots (f_{k-1})^{a_{k-1}}(f_k)^{a_k - 1}v$ by a nonzero integer (see  Lemma \ref{e_kw_{s+1}}(iii)) . Indeed  all other terms which appear  involve
monomials which end  with $(f_k)^{a_k}$ and  are larger than $(f_1)^{a_1} \cdots (f_{k-1})^{a_{k-1}}(f_k)^{a_k -1}$ in the ordering.
 
Similarly, apply  $e_k$ to terms on the right side of (\ref{disp}). All terms on the right side of (\ref{disp}) have $c_j = 0$ for $j > k.$ So the monomials have the form $(f_1)^{c_1} \cdots (f_k)^{c_k}$ with $c_k \ge a_k.$ We claim that  the term $n_{c_1 \cdots c_k}$ only involves negative roots in the span of the fundamental roots $\alpha_1, \ldots, \alpha_{k-1}$.
Indeed, if $k = l+1$, then clearly $n_{c_1\cdots c_k}$ only involves negative roots in the span of $\a_1, \ldots, \a_{k-1} = \a_1, \ldots , \a_l. $ And if $k \le l$  we obtain the claim  by comparing the $\om_t$ coefficient of  $w_{s+1}v$ with that on the right side of  (\ref{disp}) for $t \ge k,$ noting that $c_k \ge a_k.$
 Therefore the term $n_{c_1 \cdots c_k}$ commutes with $e_k$. Hence 
$e_kn_{c_1 \cdots c_k}(f_1)^{c_1} \cdots (f_k)^{c_k}v = $ 
$n_{c_1 \cdots c_k}e_k(f_1)^{c_1} \cdots (f_k)^{c_k}v$,
and we  expand as in the previous paragraph.  If $c_k > a_k,$ then all monomials appearing in the expansion are larger than $(f_1)^{a_1} \cdots (f_{k-1})^{a_{k-1}}(f_k)^{a_k -1}$ as they involve a larger power of $f_k$.  Whereas if  $c_k = a_k$, there are terms involving the monomial  $(f_1)^{c_1} \cdots (f_{k-1})^{c_{k-1}}(f_k)^{a_k-1}$, but the comments in the fourth paragraph imply  that these monomials are  larger than $(f_1)^{a_1} \cdots (f_{k-1})^{a_{k-1}}(f_k)^{a_k -1}.$ 
 
  The above procedure has now produced  an equation
 at level $x-1.$  Bring all but the term $(f_1)^{a_1} \cdots (f_{k-1})^{a_{k-1}}(f_k)^{a_k - 1}v$ to the right side of the equation to get a dependence relation where all monomials on the right are larger than $(f_1)^{a_1} \cdots (f_{k-1})^{a_{k-1}}(f_k)^{a_k - 1}$, which contradicts the minimality of $x$.

{\bf Case (II)} We may now assume that there exists a  nonzero term 
\[
 n_{c_1 \cdots c_{l+1}}(f_1)^{c_1} \cdots (f_{l+1})^{c_{l+1}}v
\]
 on the right side of (\ref{disp}), such that $c_j > 0$ for some $j > k.$  Choose $j$ maximal for this monomial.  
First note that $d_j > 0$.
Indeed if $j \le l$, then the $\om_j$ coefficient of the dominant weight afforded by
 $(f_1)^{c_1} \cdots (f_{l+1})^{c_{l+1}}v$ is $d_j - c_j$ so $d_j \ge c_j > 0.$  And if  $j = l+1$, then $(f_j)^{c_j}v = (f_{\a_{l+1}})^{c_j}v \ne 0$ so that again  $0 < c_j \le d_j.$
Consider the  tensor product 
\[
V_X(d_1\omega_1 + \cdots + d_k\omega_k) \otimes   V_X(d_{k+1}\omega_{k+1} +\cdots +d_{l+1}\omega_{l+1}).
\]
  Let $A,B$ be maximal vectors of the respective tensor factors.  Then under the action of $X$ and $L(X)$, $A\otimes B$ is a maximal vector and hence generates
an irreducible module of highest weight $(d_1,d_2, \ldots, d_{l+1}).$  So we can regard $A\otimes B = v$ and hence
\[
(f_1)^{a_1} \cdots (f_k)^{a_k}(A\otimes B) = \sum n_{c_1 \cdots c_{l+1}}(f_1)^{c_1} \cdots (f_{l+1})^{c_{l+1}}(A\otimes B).
\]
Now $V_X(d_1\omega_1 + \cdots + d_k\omega_k)$ has fewer nonzero labels than $W=V_X(\d)$, and hence the lemma holds for this module.  Expand both sides of the above equation and equate terms  of the form $dA\otimes B$. The term $((f_1)^{a_1} \cdots (f_k)^{a_k}A)\otimes B$ appears on the left side and it is nonzero since the lemma holds for $V_X(d_1\omega_1 + \cdots + d_k\omega_k)$.   
Consider a term $n_{c_1 \cdots c_{l+1}}(f_1)^{c_1} \cdots (f_{l+1})^{c_{l+1}}(A\otimes B),$ as above, for which $c_j >0$ for some $j > k.$   With $j$ maximal for this term we have 
\[
n_{c_1 \cdots c_{l+1}}(f_1)^{c_1} \cdots (f_{l+1})^{c_{l+1}}(A\otimes B) = (n_{c_1 \cdots c_{l+1}}(f_1)^{c_1} \cdots (f_{j-1})^{c_{j-1}})(A \otimes (f_j)^{c_j}B)
\]
 and there is no contribution 
of form $dA \otimes B$.  Therefore we obtain an equation
\[
((f_1)^{a_1} \cdots (f_k)^{a_k}A)\otimes B = (\sum n_{c_1 \cdots c_k}(f_1)^{c_1} \cdots (f_k)^{c_k}A)\otimes B
\]
and hence
\[
(f_1)^{a_1} \cdots (f_k)^{a_k}A = \sum n_{c_1 \cdots c_k}(f_1)^{c_1} \cdots (f_k)^{c_k}A.
\]
But this is a dependence relation in $V_X(d_1\omega_1 + \cdots + d_k\omega_k)$ at level $x$, where the monomials on the right side are larger
than the one on the left (see (\ref{disp})), contradicting the assumption that $W = V_X(\d)$ is a 
counterexample with a minimal  number of nonzero labels.  \hal

\noindent \textbf{Completion of proof of Theorem \ref{LEVELS} }
At this point we can complete the proof. Fix a level $i$. By Lemma \ref{W_i}(ii), 
\[
W^{i+1}(Q_X)  =  \sum N(f_1)^{a_1} \cdots (f_{l+1})^{a_{l+1}}v,
\]
 where the sum  is over all non-negative sequences with $a_1 + \cdots + a_{l+1} = i$. By Lemmas \ref{maxvectors}, 
\ref{a_jled_j} and \ref{nonzero}, all such vectors $(f_1)^{a_1} \cdots (f_{l+1})^{a_{l+1}}v$ contribute a  maximal vector in a composition factor of $W^{i+1}(Q_X)$, provided $a_j\le d_j$ for all $j$. Finally, by Lemma \ref{fiweight}, the weight afforded by such a maximal vector is $\sum_{j=1}^l (d_j-a_j+a_{j+1})\o_i$. This proves both parts of Theorem \ref{LEVELS} .  \hal

\section {Levels for $X = A_2$ }\label{leva2}

In this section we apply the results of the previous sections to study the embedding $L_X' < L_Y'$ for the case $X = A_2$, where the result can be displayed in a rectangle. 

In the next result $L_X' = A_1$, and we write just the positive integer $r$ for the irreducible module $V_{L_X'}(r\o_1)$.
 
 \begin{lem}\label{rectangle} Assume that $X = A_2$ and $W=V_X(\delta)$ with $\delta = r\omega_1 + s\omega _2$ and 
 $r \ge s \ge 0$. 
\begin{itemize}
\item[{\rm (i)}] There are $r+s + 1$ levels for the action of $T$ on $W$, namely $W_i = W^{i+1}(Q_X)$ for 
$0 \le i \le r+s$.
 \item[ {\rm (ii)}]  If $i \le s,$ then $W_i \downarrow L_X' = (r+i)\oplus (r+i-2) \oplus  \cdots \oplus  (r-i+2) \oplus  (r-i)$.
 \item[{\rm (iii)}]  If $s < i =s+j \le r,$ then $W_i \downarrow L_X' = (r+s-j)\oplus (r+s-j-2) \oplus  \cdots  \oplus  (r-s-j)$.
 \item[{\rm (iv)}] If $r< i = r+j$, then  $W_i \downarrow L_X' = (2s-j)+(2s-j-2) \oplus  \cdots  \oplus  j$.
 \end{itemize}
 \end{lem}
 
 \pf This is a direct application of Theorem \ref{LEVELS} .  Let  $\Pi(X) =  \{\alpha_1, \alpha_2 \}$ so that
  $\Pi(L_X') = \{\alpha_1\}.$ Write $f_1 = f_{\alpha_1 + \alpha_2}$ and $f_2 = f_{\alpha_2}.$
 The highest weights of the composition factors at level $i$ are afforded by the vectors 
$f_1^xf_2^yv$ for $x+y = i$ subject to the conditions $0 \le x \le r$ and $0 \le y \le s$. Part (i) follows.
 
If  $i \le s,$ then  the possibilities  are $(x,y) = (0,i), (1,i-1), \ldots, (i,0)$ and this leads to (ii).
Now suppose $s < i \le r$  and write $i = s + j.$  Then the possible choices are $(x,y) = (j,s), (j+1,s-1), \ldots, (i,0)$
and these yield the summands in (iii).  Finally assume $r < i = r+j $ and $j \le s$. In this case the possibilities are
$(x,y) = (i-s,s), (i-s+1,s-1), \ldots, (r,j)$, giving (iv).   \hal

 \vspace{4mm}
 The following is the explicit rectangle for the representation $\d = 63$.

	 \[
 \begin{array}{cccccccccc}
	&&&6&&&&&&\\
	&&7&&5&&&&&\\
	&8&&6&&4&&&&\\
	9&&7&&5&&3&&&\\
	&8&&6&&4&&2&&\\
	&&7&&5&&3&&1&\\
	&&&6&&4&&2&&0\\
	&&&&5&&3&&1&\\
	&&&&&4&&2&&\\
	&&&&&&3&&&\\
	
\end{array}
\]
Below is the general case for $\d = rs$.
{ \scriptsize
  \[
 \begin{array}{ccccccccccccccc}
	&&&&&r&&&&&&&&&\\
	&&&&r+1&&r-1&&&&&&&&\\
	&&&r+2&&r&&r-2&&&&&&&\\
	&&r+3&&r+1&&r-1&&r-3&&&&&&\\
	&\iddots &\vdots&&&&\vdots&&\vdots&\ddots&&&&&\\
	r+s&&\cdots&&&&&&\cdots&&r-s&&&&\\
	&r+s-1&&\cdots &&&\cdots&&&\cdots&&r-s-1&&\\
	&& \ddots&&&&&&&&&&\ddots&\\
	&&&2s&&&\cdots&&&\cdots&&&&&0\\
	&&&&2s-1&&&\cdots&&&\cdots&&&1&\\
	&&&&&\ddots&&&&&&&\iddots&&\\
	&&&&&&&&&&&&&&\\
	&&&&&&s+2&&s&&s-2&&&&\\
	&&&&&&&s+1&&s-1&&&&&\\
	&&&&&&&&s&&&&&&\\
\end{array}
\]	}

The above result provides the precise embedding of $L_X'$ in $L_Y'$,
which  is of considerable importance in view of Proposition \ref{induct}(i).
Indeed, if $Y = SL(W)$ then $L_Y'$	is a product of simple factors of type $A$ and $L_X'$ projects
to each factor as a sum of irreducibles of highest weights as presented. For example in
the $63$ example above, $L_Y' = A_6 \times A_{13} \times A_{20} \times \cdots \times A_{11} \times A_7 \times A_3$.

  Note that in the special case $s = 0$ the above array reduces to just the top right-slant row which is
 $r, \ldots, 0$.  Here $L_Y' = A_r \times A_{r-1} \times \cdots \times A_1$.

 \section{$Y$-levels}

We conclude this section with the following useful result on levels.  We use notation as in Chapter \ref{notation}.  Let $Y = SL(W)$ and let $P_Y = Q_YL_Y$ be a parabolic subgroup such that $L_Y' = C^0 \times \cdots \times C^k$  and $Q_Y$ is a product of root groups for negative roots.  By Corollary \ref{levelsnontrivial}, there is a single fundamental root $\gamma_1$ between $C^0$ and $C^1.$ Let $\l$ be a dominant weight for $Y$, let $V = V_Y(\l)$ and let $\mu^i$ denote the restriction of $\l$ to to $T_Y\cap C^i$.

\begin{prop}\label{v2gamma1} With the above notation,  assume that $\la \l, \gamma_1 \ra > 0$. Then 
\[
V^2_{\gamma_1}(Q_Y) \supseteq V_{C^0}(\mu^0) \otimes V_{C^0}(\l_{r_0}^0) \otimes V_{C^1}(\l_1^1 + \mu^1) \otimes V_{C^2}(\mu^2) \otimes \cdots \otimes  V_{C^k}(\mu^k).
\]
\end{prop} 

\pf  Let $J$ be the  parabolic subgroup of $Y$ with fundamental system $ \la \b_1^0, \ldots, \b_{r_0}^0, \gamma_1 \ra$. There is a maximal parabolic subgroup $P_J = Q_JL_J$ of $J$ such that $L_J'=C^0$.  Let $\l_J$ denote the restriction of
$\l$ to  $T_Y \cap J$.  Theorem \ref{LEVELS}  gives  the composition factors of $C^0$ on $V^2(Q_J)$.
Moreover, using the hypothesis  $\la \l, \gamma_1 \ra > 0$, the proof of Corollary \ref{cover} shows that 
$V^2(Q_J) = V_{C^0}(\mu^0) \otimes V_{C^0}(\l_{r_0}^0).$

The next step is to lift this to $Y$. 
Write $\l_J = \sum_1^{r_0+1} a_i\l_i$. List the subscripts $j$ for which $a_j\ne 0$ in decreasing order $i_1>\cdots >i_r$ (so that $i_1 = r_0+1$). 
Then $f_{i_1} > \cdots > f_{i_r}$ in the ordering of monomials,
where $f_{i_j} = f_{ \b_{i_j}^0 + \cdots +\b_{r_0}^0 +\gamma_1}$.

Lemmas \ref{maxvectors}, \ref{a_jled_j}, and \ref{nonzero} show that $V^2(Q_J) = \sum_j R_{i_j}$ where 
each $R_{i_j}$ is an irreducible module
for $C^0$ having a maximal vector $v_{i_j}$ of form $(f_{i_j} + \sum_{s < j}n_sf_{i_s})v$ where the terms $n_s$ are in the algebra generated by $f_\b$ for $\b \in \Sigma(C^0)^+$.  Then for each $j,$ $v_{i_j}$ is  a maximal vector for $C^1 \times  \cdots \times C^k$ affording weights $(\l_1^1+\mu^1), \mu^2, \ldots, \mu^k$, for the respective factors. 
Therefore 
\[
L_Y'v_{i_j} = R_{i_j} \otimes V_{C^1}(\l_1^1 + \mu^1) \otimes V_{C^2}(\mu^2) \otimes \cdots \otimes  V_{C^k}(\mu^k).
\] 
Adding these terms it follows that 
\[
V^2_{\gamma_1}(Q_Y) \supseteq  \sum_j L_Y'( v_{i_j}) = V_{C^0}(\mu^0) \otimes V_{C^0}(\l_{r_0}^0) \otimes V_{C^1}(\l_1^1 + \mu^1) \otimes V_{C^2}(\mu^2) \otimes \cdots \otimes  V_{C^k}(\mu^k),
\]
as required. \hal

\section{Method of Proof  - Level Analysis }\label{levanal}

Let $X = A_{l+1}$ be embedded in $Y = A_n$ via an irreducible $X$-module $W$  of highest weight $\d$.
In this subsection we describe a method for proving that most representations of $Y$ are not MF upon restriction to $X$.
  This method and more elaborate variations of it will be used a great many times in this paper.  

Choose parabolic subgroups $P_X = Q_XL_X$ and $P_Y = Q_YL_Y$, such that Lemma \ref{parabembed} holds.  Assume that $V = V_Y(\l)$ is such that $V\downarrow X$ is MF. Then Proposition \ref{induct}(ii) shows that 
\begin{equation}\label{sumf}
V^2 = V^2(Q_Y) \downarrow L_X' =  \sum _{i: n_i = 0} V_i^2(Q_X) + \sum _{i: n_i = 1} V_i/[V_i,Q_X].
\end{equation}
Here the first sum consists of terms that arise from $V^1$, and the second sum is MF.

 The key idea is to show  that the second sum in  (\ref{sumf}) is actually not MF, 
which is a contradiction. This is accomplished  by producing a number of composition factors in 
$V^2(Q_Y)$ for $L_Y'$, restricting them to  $L_X'$  and finding  composition factors that appear with multiplicity at least two greater than what can possibly arise from $V^1$.

We have $L_X' < L_Y' = C^0 \times \cdots \times C^k$.  From  Chapter \ref{notation} we have $V^2(Q_Y) = \sum_{j=1}^k V_{\g_i}^2(Q_Y)$, where $\g_i$ is the node between $C^{i-1}$ and $C^i$ (see Corollary \ref{levelsnontrivial}).
Fix $i\ge 1$, let $\g = \g_i$, and let  $v^+$ be a maximal vector for $Y$ in $V$.  If $\la \l, \g \ra > 0, $ then a consideration of positive roots in the Lie algebra of $C^0 \times \cdots \times C^k$ shows that  $f_{\g}v^+$ is a maximal vector for $L_Y'$ and affords an irreducible module, say $F_{\g},$ in $V_{\g}^2(Q_Y).$ 
The highest weight of $F_{\g}$  is $\nu_{\g} = (\l - \g) \downarrow (C^0 \times \cdots \times C^k)$, and $\nu_{\g}$  restricts to $\mu^i+ \l_1^i$ for $C^i,$ to $\mu^{i-1} + \l_{r_{i-1}}^{i-1}$ for $C^{i-1},$ and  to $\mu^s$ for $s \ne i-1,i. $  If $\la \l, \g \ra = 0,$ then set $F_{\g} = 0.$  
 
We can obtain a further $L_Y'$-composition factor in $V^2_\g(Q_Y)$ if $\mu^i \ne 0$. To see this,  let $t$ be minimal such that $\mu^i$ has a nonzero label for fundamental root $\b_t^i$. Consider the weight $\nu_i = \l - \g - \b_1^i - \cdots - \b_t^i.$ If $\la \l,\g\ra >0$,   the $\nu_i$-weight space has dimension $t+1$ in $V_{\g}^2(Q_Y)$  and dimension $t$ in $F_{\g}$; 
and if $\la \l,\g\ra = 0$ this weight space has dimension 1. Therefore there is a vector $w_i$ in the weight space but not in $F_{\g}$, and one checks that $w_i$ is a maximal vector modulo $F_{\g}$. Consequently there is a composition factor, $F_i,$ of $V_{\g}^2(Q_Y)$ of highest weight $\nu_i$, and we can easily determine its restriction to $C^{i-1}$ and $C^i$.
 
We can also get a new composition factor  in $V^2_\g(Q_Y)$ if $\mu^{i-1}\ne 0$. In this case, let $j$ be maximal such $\mu^{i-1}$ has a nonzero label for fundamental root 
$\b_j^{i-1}$.   Then as above there is a composition factor, $F_{i-1}$, with highest weight $\nu_{i-1} = \l  - \b_j^{i-1} - \cdots - \b_{r_{i-1}}^{i-1} - \g.$
 
 This idea can be pushed still further if both $\mu^{i-1}$ and $\mu^i$ are nonzero.  Indeed, here there is a new composition factor $F_{i-1,i}$  of highest weight $\nu_{i-1,i} = \l - \b_j^{i-1} - \cdots -  \b_{r_{i-1}}^{i-1} - \g - \b_1^i - \cdots - \b_t^i.$  The justification varies  according to whether or not $\la \l,\g\ra$ is 0.  To simplify notation we set $r = r_{i-1}.$
 
 First suppose $\la \l, \g \ra= 0$.  Then Proposition \ref{cav2} shows that  the dimension of the $\nu_{i-1,i}$-weight space  in $V^2_\g(Q_Y)$ is $r-j+2+t$, whereas this weight space has dimension $r-j+1$ in $F_i$ and dimension $t$ in $F_{i-1}$.  So the usual argument produces a new composition factor in $V^2_\g(Q_Y)$ of  highest weight $\nu_{i-1,i}$.
 
 The argument in the case $\la \l, \g \ra\ne 0$ is similar but the count is a bit more complicated.  Here Proposition  \ref{cav2} shows that the  $\nu_{i-1,i}$-weight space dimensions 
 in $V_{\g}^2(Q_Y)$, $F_{\g}$, $F_i$ and $F_{i-1}$   are 
   $(r-j+2)(t+1)$,$ (r-j+1)t$, $t$ and  $r-j+1$, respectively.  It follows that the weight space in $V^2_\g(Q_Y)$ has dimension 1 more than the
  sum of the other three dimensions and hence there is a composition factor of highest weight $\nu_{i-1,i}$.
  
   To conclude, we have the following $L_Y'$-composition factors in $V^2_\g(Q_Y)$:
   \begin{itemize}
   \item[(1)]  $F_{\g}$,  nontrivial only if  $\la \l, \g \ra > 0$,
   \item[(2)] $F_i$,  nontrivial only if $\mu^i\ne 0$,
   \item[(3)] $F_{i-1}$,  nontrivial only if  $\mu^{i-1}\ne 0$,
   \item[(4)] $F_{i-1,i}$,  nontrivial only if  $\mu^{i-1}\ne 0$ and $\mu^i\ne 0$.
\end{itemize}
As indicated above, at this point we consider the restrictions of these composition factors to $L_X'$ and try to show that the second sum in the expression (\ref{sumf}) for $V^2$ is not MF. 

In many cases the above approach gives a contradiction;  
when it does not (for example, when the second sum in (\ref{sumf}) turns out to be MF), we need to analyze higher levels $V^{d+1}$ for 
$d\ge 2$, and aim to contradict the fact that the second summand in the expression for $V^{d+1}$ in Proposition \ref{induct}(iii) is MF. Here the analysis is often much more complicated than the above, as composition factors in $V^{d+1}(Q_Y)$ can be much harder to find than those in $V^2(Q_Y)$.


\chapter{Multiplicity-free families}\label{MFfam}

In this chapter we show that the  configurations in Tables \ref{TAB1}--\ref{TAB4} of Theorem \ref{MAINTHM} are indeed MF.  
We concentrate on cases where $X$ has rank at least 2, since \cite{LST} settles the case where $X =A_1.$
Specifically, we prove the following result.

\vspace{4mm}

\begin{theor}\label{yesMF}
Let $X = A_{l+1}$ with $l\ge 1$, let $W = V_X(\d)$ and $Y = SL(W) = A_n$. Let $V = V_Y(\l)$, and suppose $\l,\d$ are as in Tables $\ref{TAB1}-\ref{TAB4}$ of Theorem $\ref{MAINTHM}$. Then $V\downarrow X$ is multiplicity-free.
\end{theor}

In Section \ref{cnbn} we study 
wedge and symmetric powers of  $W$ and show that their restrictions to orthogonal groups
are MF.  In the following sections we work through Tables \ref{TAB1}-\ref{TAB4}, showing that in each
case the configuration is indeed MF.

\section{Restrictions of $SL_n$ representations to $SO_n$}\label{cnbn}
  
In this subsection we record some results on how the wedge powers and symmetric
powers of the natural module for $SL(W)$ restrict to $SO(W)$.  
Let $n= \dim W$ and write $n=2m$ or $2m+1$.  
It is well known (see \cite[19.2, 19.14]{fulthar}) that $\wedge^k(W) \downarrow SO(W)$ is irreducible for all 
$k < m$.  The following result settles the case of symmetric powers. Write $\o_1$ for the highest weight of the natural module $W$. Recall that $SO(4) \cong A_1A_1$ and $SO(3)\cong A_1$, and we denote the highest weight of an $A_1$-module simply by a positive integer $k$. 

\begin{thm}\label{SOpowers}  Write $D = SO(W)$. Then for any positive integer $c$,
\[
S^c(W) \downarrow D = \left\{\begin{array}{l}
\sum_{i=0}^{[c/2]} V_D((c-2i)\o_1), \hbox{ if }\dim W \ge 5 \\
\sum_{i=0}^{[c/2]} V_D(c-2i, c-2i), \hbox{ if }\dim W =4 \\
\sum_{i=0}^{[c/2]} V_D(2(c-2i)), \hbox{ if }\dim W =3.
\end{array}
\right.
\]




\end{thm}

\pf  We will use the result in \cite[p.427]{fulthar} which gives the following formula for restricting representations from $SL(W)$ to $SO(W)$: 
 for a partition $\g = (\g_1\ge \g_2\ge \cdots \ge \g_n\ge 0)$ with $\g_i=0$ for $i>\frac{n}{2}$,   
\begin{equation}\label{fhform}
V_{GL(W)}(\g)\downarrow SO(W) = \sum_{\xi} N_{\g,\xi}V_{SO(W)}(\xi),
\end{equation}
where the sum is over all partitions $\xi = (\xi_1 \ge \xi_2 \ge \cdots \ge \xi_m \ge 0)$, and
\[
N_{\g,\xi}  = \sum_{\e} c_{\e,\xi}^{\g},
\]
the sum  over all partitions $\e = \e_1 \ge \e_2 \ge \cdots \ge 0$ with even parts and the
terms $c_{\e,\xi}^{\g}$ in the sum are Littlewood-Richardson coefficients. 

We sketch the argument, noting that in  Subsection \ref{a3a5} we will give further details on using this formula.
Here we have  $V = S^c(W)$ so we take the partition $\gamma = (c,0,\ldots,0).$
The strategy is as follows.  Fix $\g$.  
Given  $\xi$, we show that there is at most one even  partition $\e$  
such that $c_{\e,\xi}^{\g} \ne 0$  and that this coefficient is  $1$.  This is particularly simple in this case.
 Indeed, the relevant tableaux will have just one row of length $c$.
So we can only use partitions $\e = (\e_1 \ge \e_2 \ge \cdots )$ with $\e_1 = 2d$ for some $d \ge 0$.  Then the row
has $2d$ blank entries followed by $c - 2d$ 1's.  Therefore, $\xi = (c-2d,0, \ldots ,0)$ and the corresponding
representation has highest weight $(c-2d)\o_1$, where $\o_1$ is the highest weight of the natural module $W$. The result follows.  \hal

 \section{Table \ref{TAB1} configurations}\label{keys}

 In this section will show that the configurations in Table \ref{TAB1} of Theorem \ref{MAINTHM} are indeed MF.   Most of
 the proofs are based on the domino method described in Section \ref{decom}.
 
We begin with two lemmas which will be used, often implicitly, in several of the proofs to follow.  Suppose $\d = \sum d_i \om_i$ is  the highest weight of an irreducible  representation of $A_n$.  Let $(c_1, \ldots , c_{n},0)$ be a corresponding partition.
As explained in Section \ref{decom}, we can use Yamanouchi domino tilings of shape $T = (2c_1,2c_1,\ldots ,2c_n,2c_n)$ to decompose the tensor product $\d \otimes \d$.The next lemma shows that distinct domino tilings correspond to distinct composition factors.
 
 \begin{lem}\label{unique}  If  $(1^{a_1}, \ldots, (n+1)^{a_{n+1}}) $ and $(1^{b_1}, \ldots, (n+1)^{b_{n+1}})$ are the weights of two domino tilings
of shape $T$  which correspond to the same irreducible representation for $A_n,$ then $a_i = b_i$ for all $i.$
  \end{lem}
  
  \pf  The hypothesis forces $a_i - a_{i+1} = b_i - b_{i+1}$ for $1\le i\le n$.  We may take it
  that  $a_1 \ge b_1.$  Then  from  $a_1 - a_2 = b_1 - b_2$ we conclude that $a_2 \ge b_2.$  Continuing
  we have $a_i \ge b_i$ for all $i$.
  
   On the other hand the total number of dominoes used in each tiling is the same. Therefore $\sum a_i = \sum b_i$ and it follows
 that $a_i = b_i$ for all $i.$ \hal
 
 \begin{lem} \label{chain}  Fix a Yamanouchi domino tableau, and let its weight be $(1^{a_1}, \ldots, k^{a_k})$ with  $a_k >0$. Then the following hold.
 \begin{itemize}
\item[{\rm (i)}]  $(a_1, \ldots, a_k)$ is a partition.  
 \item[{\rm (ii)}]  Let $w_1\ldots w_r$ be the column reading of the tableau and let $i\in \{1,\ldots ,r\}$. If 
$j = w_i$, then $j, j-1, \ldots, 1$ all appear in  $w_i\ldots w_r$.
\end{itemize} 
 \end{lem}
 
 \pf    (i) The Y-condition applied to the column reading $w_1\ldots w_r$ of the domino tableau implies that $a_1\ge a_2\ge \cdots$. Hence $(a_1,\ldots ,a_k)$ is a partition.

(ii)  This follows by applying the Y-condition to  $w_i\ldots w_r$.   \hal
  
 \subsection{Weights $c\o_i + \o_{i+1}$ and  $\o_i + c\o_{i+1}$}

Here we establish the following result.

\begin{prop}\label{comi+omi+1}  Let $W$ be the irreducible module for $A_n$ with highest weight either $c\o_i + \o_{i+1}$ or 
$\o_i + c\o_{i+1}$ for some $i$, where $c \ge 1$. Then  $\wedge^2(W)$ and $S^2(W)$ are both MF.
\end{prop}

\pf  Replacing $W$ by its dual we can assume that the highest weight is $c\o_i + \o_{i+1}$ for fixed $i$. Using the domino approach we consider the sequence $(2c+2, \ldots, 2c+2, 2,2)$ where there are $2i$ terms $2c+2$.  

We will show that there are at most two domino tilings affording a given weight $(1^{k_1},2^{k_2}, \ldots, )$
 and if there are two, then one corresponds to a symmetric composition factor of $W \otimes W$ and the other
corresponds to an alternating composition factor.
Towards this end assume we have a domino tiling with the fixed weight above.
 
 The tableau has the shape of a $ 2i \times (2c+2)$ matrix on top of which is a $2 \times 2$ square. 
We will first analyze possible tilings of the large rectangle, leaving the small square until later.
The ${\it base}$ of the array refers to the bottom two rows.

The bottom two rows of a tiling consist of $t$ vertical 1-dominoes
followed by  a sequence of $s$ horizontal 2-dominoes on top of horizontal 1-dominoes.  Thus $2s+t = 2c+2.$   Then $t$ is even and 1 appears $s+t$ times  in the tiling and this  determines both $s$ and $t$.  Therefore, since $c$ is fixed, $k_1$ determines both $s$ and $t$.

Suppose $s > 0$.  Then at the right end of the base there is a horizontal $2$-domino lying above a horizontal $1$-domino.  At this point we state a lemma which will be used at several points in this section in slightly
different settings.  The arguments in those other  applications are only slight variations of the argument below.

\begin{lem}\label{horizontal}  The remaining dominoes in the last two columns are all horizontal with labels $3, \ldots, 2i.$  
\end{lem}

\pf  Suppose false, and write $x=2s+t=2c+2$. Then in columns $x$ and $x-1$ there is a sequence 
 of horizontal dominoes with labels $1,2,\ldots, d$ and above this there is a vertical domino in column $x$, say with label $r.$

For the purposes of column reading (see Section \ref{decom}), the labels on the sequence of horizontal dominoes contribute to column $x-1$ in the column reading.  Therefore $r$ is the smallest label in the column reading in column $x$ and the Y-condition forces $r = 1.$  But this violates the decreasing condition on columns and establishes the lemma. \hal
 
Continuing with the proof of Proposition \ref{comi+omi+1}, 
suppose $s > 0$ and apply Lemma \ref{horizontal}. Set $b_0 = 2i.$  As horizontal dominoes have the fastest vertical numerical growth,   $b_0$ is the largest label in the large rectangle. Note that $b_0 = 2i \le n+1.$   If $s > 1$, then the same argument shows that the next two columns also consist of horizontal dominoes.  This continues and determines  columns $t+1$ through $t + 2s.$ They each consist of a sequence of horizontal dominoes with labels $1, \ldots, b_0$.   These $2s$ columns contribute $(1^{s}, \ldots , b_0^{s})$ to the weight of the array.  As mentioned above $b_0$ is the largest possible
label of the large rectangle,  it appears with multiplicity $s$, and the pair $(b_0,s)$ determines the labelling and tiling of the 
$2s$ columns on the right side of the large rectangle.

The above analysis assumed $s > 0.$ On the other hand if $s = 0$, then there are no horizontal dominoes on the right side of the base.   If $t = 0$ then the labelling  and tiling of the large rectangle is determined by the above.  

So assume $t > 0$ and consider column $t$, where there is a vertical 1-domino covering the bottom two rows.   Using the Y-condition we see that column $t$ 
begins with a sequence of vertical dominoes  having labels $1, \ldots, a_1$, allowing the possibility
 $a_1 = 1.$  Vertical dominoes have the slowest vertical numerical growth, so the non-decreasing condition
on rows implies that each of the columns $1, \ldots, t$ begins with a sequence of vertical dominoes
with  labels $1, \ldots, a_1$.

 If $a_1 < i,$ then above the vertical dominoes in column $t$ with labels $1, \ldots, a_1$ we argue as in
 Lemma \ref{horizontal}, which is based on the Y-condition, to see that there is a sequence of horizontal dominoes  with labels $a_1+1, \ldots, b_1$ going to the top of the large rectangle. Note that $b_1 < b_0.$ This same pattern occurs for some even number, say $2s_1,$ of columns ending at  column $t$.  So the tiling of  these $2s_1$ columns is determined and it contributes  $(1^{2s_1}, \ldots, a_1^{2s_1},(a_1+1)^{s_1},\ldots, (b_1)^{s_1})$  to the weight of the array.  On the other hand, if $a_1 = i$, then the vertical dominoes
 go to the top of the array, each of the columns $1, \ldots, t$ consist of vertical dominoes with labels
 $1, \ldots, a_1,$ and these columns in the large rectangle  contribute $(1^t, \ldots, i^t)$ to the weight. 
 
At this point we continue filling out more and more of the large array.  For example if $2s_1 < t$, then  column $t-2s_1$ begins with vertical dominoes having labels $1, \ldots, a_2$, where $a_2 > a_1$ and then continues
with horizontal dominoes to the top of the large rectangle. There will be an even number, say  $2s_2,$ of such columns and these columns contribute $(1^{2s_2}, \ldots, a_2^{2s_2}, (a_2+1)^{s_2},\ldots, (b_2)^{s_2}).$ This process continues for say $k$ terms until all columns of the large rectangle have been tiled.   In the process we
have determined triples $(s_1,a_1,b_1), (s_2,a_2,b_2),\ldots, (s_k,a_k,b_k).$  Note that $a_i$ determines $b_i$
in each case.  Moreover, $a_k > a_{k-1} > \cdots > a_2 >a_1.$

Counting, we find that the contribution  of the large rectangle to the weight of the array is as follows
\begin{equation}\label{list2}
\begin{array}{l}
(1^{s+t},\ldots,a_1^{s+t},(a_1+1)^{s+t-s_1},\ldots, a_2^{s+t-s_1},(a_2+1)^{s+t-s_1-s_2}, \ldots, a_k^{s+t-s_1-\cdots -s_{k-1}},\ldots, \\
(a_k+1)^{s+s_1+\cdots +s_k}, \ldots , b_k^{s+s_1+\cdots s_k}, (b_k+1)^{s+s_1+\cdots s_{k-1}}, \ldots , b_{k-1}^{s+s_1\cdots s_{k-1}}, \ldots, (b_1+1)^s, \ldots, b_0^s).
\end{array}
\end{equation}
  Next consider the possible tilings of the $2 \times 2$ square. The square is tiled with either two vertical dominoes with labels $x \le y$ or  a horizontal $y$-domino lying over a horizontal $x$-domino in which case $x < y.$  In either case we have $x,y > a_k$ by the column decreasing condition so that $x,y > a_i$  for all $i.$     Therefore from the form of (\ref{list2}) we see that the weight of the array determines
  each of $s, t, a_1,s_1,a_2, s_2, \ldots, s_{k-1}, a_k$ (recall that $k_1$ determines $s$ and $t$).  
Moreover, $t = \sum_{i=1}^k 2s_i$, so
  $s_k$ is also determined.  It follows that the weight of the array completely determines the tiling of the large
  rectangle. 
  
 The weight of the array also determines the labels $x \le y.$  If $x = y$  the only possibility is
 where both tiles are vertical with label $x.$  So in this case the tiling is uniquely determined from
 the weight.  If $x < y$, then there are potentially two tilings satisfying the various conditions.
 If indeed this occurs  then one will correspond to an alternating summand of $W \otimes W$ and the other is symmetric.  \hal

   \subsection {Weights $c\o_1 + \o_i$ }

In this subsection we prove

\begin{prop}\label{clambda1+lambdai}   Let $W$ be the irreducible module for $A_n$ with highest weight $c\o_1 + \o_i$ for some $i>1$, where $c \ge 1$. Then  $\wedge^2(W)$ and $S^2(W)$ are both MF.
\end{prop}

\pf We use the domino approach as in the previous result.
First note that  $\d = c\o_1 + \o_i$ corresponds to the partition $(c+1, 1, \ldots, 1),$ where 1 appears $i-1$ times.  Now to convert the partition to dominoes
we double each entry and  study domino tableaux satisfying the Y-condition of shape $ (2c+2,2c+2,2,\ldots, 2),$    where the number of 2's is $2(i-1).$  

As in the previous result we will show that there are at most two domino tilings of any given weight, and if there are two then one corresponds to summand of $S^2(W)$ and the other to a summand of $\wedge^2(W)$.  Therefore we now consider a domino tableau with a fixed weight.

The base consists of a $2\times (2c+2)$ rectangle.  In view of the Y-condition the base must be labelled with  $t$ vertical 1-dominoes followed by $s$ horizontal 2-dominoes on top of horizontal 1-dominoes.  Thus $2s+t = 2c+2.$   Then $t$ is even and 1 appears $s+t$ times  in the labelling and this determines both $s$ and $t$.  Note also that the number of  $1 \times 1$ squares in each of columns 1 and 2
is $2i$.  Above the base, all remaining layers have width $2$.  

\vspace{2mm}
\noindent{\bf A.} Assume first that the base
consists of only vertical 1-dominoes so that the labelling has $1$ occurring with multiplicity $2c+2.$  The labellings that occur here have the largest possible number of 1's and
will not occur in other situations to follow.  It will follow from our considerations that for (A), multiplicities are not
possible.

Assume that above the base there are  two vertical dominoes.
The strictly decreasing condition on columns implies  that these dominoes
have labels greater than $1.$  And the Y-condition implies that the second vertical domino has  label $2$.
As the rows are weakly increasing we conclude that both vertical dominoes in
the second layer are labelled by $2$.  If the next level has two vertical dominoes, then
the same considerations show that they are both labelled by $3.$  Now continue in this way. 

\vspace{2mm}
{\bf A(i)} Assume that all remaining  dominoes are vertical.  Then the weight of the tiling is 
 \begin{equation} 
 (1^{2c+2},2^2, \ldots i^2).
 \end{equation}
 Also, Theorem \ref{dominoes} shows that this factor is in the symmetric square.  Of course, this is obvious, anyway, since
 the corresponding highest weight is $2\d$.

\vspace{2mm}
{\bf A(ii)} Assume that, including the base there is a sequence  of pairs of vertical dominoes.  The conditions imply that at each level the vertical dominoes have the same label, so say the pairs have labels $1, 2,3,\ldots, a.$  Following this sequence   
there is a horizontal domino which must be labelled $a+1$ by the Y-condition. Labels above this are strictly larger than $a+1$ as columns decrease. This begins a sequence of horizontal dominoes, which would have to be labelled   $a+1, \ldots, a+b$ for some $b \ge 1.$

A variation of the argument of Lemma \ref{horizontal} shows that there cannot be a pair of vertical dominoes, say with 
labels $j$ in column 1 and $k$ in column 2 above the horizontal domino with label $a+b$.  To review that argument note that for purposes of column reading  the labels on the horizontal dominoes contribute only  to column 1.  Then the Y-condition forces $k = a+1$ and this contradicts the fact that
labels in columns are strictly decreasing. Therefore the horizontal domino with label $a+b$ is in the top row of the array.

If $a = 1$ the weight is 
\begin{equation} 
 (1^{2c+2},2^1, \ldots (b+1)^1),
 \end{equation}
 where $b = 2i-2$.  And if $a \ge 2$, then the weight is 
\begin{equation}\label{eight}
(1^{2c+2},2^2, \ldots a^2,(a+1)^1,\ldots ,(a+b)^1).
\end{equation}
In either case there is a unique tiling for a given weight. Moreover,  the corresponding composition factor is symmetric or alternating according as $i-a$ is even or odd.
\vspace{4mm}

Note  in both {\bf A(i)} and {\bf A(ii)}  the weight of the tiling determines $t = 2c+2$ (i.e. the multiplicity of 1).  The largest label in {\bf A(i)} is strictly greater than that in {\bf A(ii)}, so the weights
differ in these two cases.  Moreover we see from (\ref{eight}) that the weight in  {\bf A(ii)} determines $a$, so that the tiling is uniquely determined by the weight.

\vspace{4mm}
From now on we assume that the base does have horizontal dominoes.  There are several configurations to consider.  As mentioned above none of the weights to follow can be of type {\bf (A)} as they have smaller multiplicity of  $1$.
\vspace{4mm}

\noindent{\bf B.}  Assume next that the base has only horizontal dominoes.  That is there are two
layers of such dominoes, each with $c+1$ entries.  The bottom row is labelled with 1's and the top row is labelled with 2's.  So the weight begins with $1^{c +1}2^{c+1}.$ Note that the weights that occur here have the fewest possible number of 1's and will therefore not occur in types to follow.  And as in (A)  we will see that there is a unique  tiling.

  The remaining layers in the first two columns can be described as in (A).  There are vertical dominoes for
a while followed by horizontal ones, but it is possible  that there are either no vertical dominoes or
no horizontal dominoes.  Say the vertical dominoes have labels $3, \ldots, a$ and the horizontal ones
$a+1, \ldots, a+b.$  So in the general case the  weight  is 
\begin{equation}\label{list1}
(1^{c+1},2^{c+1},3^2, \ldots ,a^2,(a+1), \ldots ,(a+b))
\end{equation}
and there is a unique tiling for a given weight as (\ref{list1}) determines $a$ and $b$.
We allow the special cases where the terms $3^2, \ldots ,a^2$ or $(a+1), \ldots ,(a+b)$ do not appear.

\vspace{4mm}
\noindent{\bf C.}  Finally, suppose that $t,s > 0$  and recall that $t$ is even. The  fixed weight determines both $t$ and $s.$  We shall see that multiplicities  do occur here.

Above the base in columns 1 and 2 there can be either  two vertical dominoes or a horizontal domino.  In the horizontal case the Y-condition implies that this domino has label either $2$ or  $3.$  

\vspace{4mm}
{\bf C(i)} Assume that the layer above the base is a horizontal 3-domino. The variation of Lemma \ref{horizontal} given in {\bf (A)(ii)} above shows  that all remaining dominoes must be horizontal with labels $4, \ldots, 2i.$  Indeed if some level has a vertical domino, then going to the first such level we contradict the Y-condition working in the second column.
So in this special
case the weight is 
\begin{equation}\label{ten}
(1^{t+s},2^s,3^1, \ldots ,(2i)^1).
 \end{equation}
Morover, given the weight as in (\ref{ten}) we must be in case {\bf C(i)} as otherwise the multiplicity of 2 would be
greater than $s$.
\vspace{4mm}

{\bf C(ii)}  Assume that there are two vertical dominoes above the base.  The Y-condition allows for the labelling of the  vertical dominoes to be $(2,2), (2,3),$ or $(3,3),$  although the  third case requires $s \ge 2$.  We then have a string of vertical dominoes  above the base  with labels   $(2,2), \ldots,(a,a); $  $(2,3), \ldots, (a,a+1);$ or $(3,3), \ldots,(a,a), $  respectively. 
Above this there can be a sequence of horizontal dominoes with labels $(a+1, \ldots , a+b)$, 
$(a+2, \ldots, a+b),$ or $(a+1, \ldots, a+b),$ respectively, extending to the top of the array.  In the second case we will consider separately the situation where $a =2$, in which case the first horizontal domino has label 4. 

In the general case where there are horizontal dominoes present and excluding the special case mentioned at the end of the previous paragraph the weights are
as follows
 \begin{equation}\label{11}
(1^{t+s},2^{s+2},3^2, \ldots, a^2,(a+1)^1,\ldots, (a+b)^1)
 \end{equation}
\begin{equation}\label{12}
( 1^{t+s},2^{s+1},3^2,\ldots, a^2,(a+1)^1,\ldots, (a+b)^1)\;(a>2)
 \end{equation}
 \begin{equation}\label{13}
( 1^{t+s},2^s,3^2, \ldots ,a^2,(a+1)^1,\ldots, (a+b)^1).
  \end{equation}
 Note that if there are only two vertical dominoes above the base, then we delete the terms $3^2\ldots a^2$ in (\ref{11}) and (\ref{12}) and the terms $4^2\ldots a^2$  in (\ref{13}).  And if there are no horizontal
 dominoes we delete the terms $(a+1)^1\ldots (a+b)^1$.  
 In the special case of (\ref{12})  mentioned above with $a = 2$  the weight is
 \begin{equation}\label{14}
(1^{t+s},2^{s+1},3^1,4^1,\ldots, (2i-1)^1).
  \end{equation} 
 Consideration of equations  (\ref{11})--(\ref{14}) shows that the weight determines the tiling of the array. Indeed the multiplicity of $2$ indicates which of the three cases $(2,2), (2,3),$ or $(3,3)$ we are in and the form of the weight determines $a$ and $b$.

  \vspace{4mm} 
 {\bf C(iii)} Assume that the level above the base has a  horizontal 2-domino.
 If there are additional horizontal dominoes then, as in Lemma \ref{horizontal}, the Y-condition implies that there are no
 vertical dominoes and that above the base there are only horizontal dominoes
 with labels  $2, \ldots, a$ where $a = 2i-1.$  In this case the weight is 
  \begin{equation}\label{15}
(1^{t+s},2^{s+1},3^1, \ldots, (2i-1)^1).
 \end{equation}
Otherwise, it is possible to have vertical dominoes  labelled as $(3,3), \ldots, (a,a)$ followed possibly by horizontal
 dominoes with labels $a+1, \ldots, a+b.$  Here the weight is
 \begin{equation}\label{16}
(1^{t+s},2^{s+1},3^2, \ldots, a^2,(a+1)^1, \ldots, (a+b)^1).
 \end{equation}
In this case we again see that the weight determines $a$ and $b$.
 
 \vspace{4mm}
  We now complete the proof of Proposition \ref{clambda1+lambdai}.  In view of Lemma \ref{unique} we need only  look for coincidences among the above labellings.   By earlier remarks
 we need only consider coincidences among  types {\bf (C)}.
 
 From the multiplicity of $2$ we see that the only possible coincidences
 occur among (\ref{12}), (\ref{14}), (\ref{15}) and (\ref{16}).  In cases (\ref{12}) and (\ref{16}) we have $a >2$ so neither of these could be of type (\ref{14}) or (\ref{15}).   And comparing (\ref{12}) and (\ref{16}) we see that the latter has two more horizontal tiles than the former.  So one of these corresponds to an alternating composition factor and the other a symmetric composition factor.  The same holds for   (\ref{14}) and (\ref{15}), completing the proof.  \hal
 
    \subsection {Weights $\o_1 + c\o_i$}

In this subsection we prove
  
\begin{prop}\label{lambda1+clambdai}    Let $W$ be the irreducible module for $A_n$ with highest weight $\o_1 + c\o_i$ for some $i>1$, where $c \ge 1$. Then  $\wedge^2(W)$ and $S^2(W)$ are both MF.
\end{prop}

\pf As in the previous cases we will work with domino tilings, but there are more  configurations  possible in this case.  
 Taking duals we will work with $\d = c\o_j + \o_n,$ with $1 \le j < n.$    The highest weight $\d$  corresponds to the partition 
$(c+1, \ldots, c+1, 1, \ldots, 1)$  where $c+1$ appears $j$ times and $1$ appears $n-j$ times.  We therefore consider a domino tiling of shape $(2c+2,\ldots, 2c+2,2, \ldots, 2)$,  where the term $2c+2$ occurs $2j$ times and $2$ appears $2(n-j)$ times.  This determines the shape
of the array.  It has two sections.  The lower part is a $2j \times (2c+2) $ rectangle.  On top of
this there is a  $2(n-j) \times 2$ rectangle which lies above the first two columns in the lower part.

We will show that for a fixed weight there are at most two tilings of the array affording the weight and if there are two, then one must correspond to an alternating composition factor and the other a symmetric composition factor of $\d \otimes \d.$

The base consists of $t$ vertical 1-dominoes followed by  $s$ horizontal 1-dominoes on top of which are $s$ horizontal 2-dominoes.  Then $t + 2s = 2c+2$. 
For a given tiling, the number of 1's is $s+t$, so as $c$ is fixed this determines both $s$ and $t$.  Therefore, two different tilings yielding the same weight must have the same base.

The argument here is similar to, but more complicated than, Proposition \ref{comi+omi+1}.  We begin with the tiling of the lower part.  If $s > 0$ then  using Lemma \ref{horizontal} we find that columns $t+1, \ldots, 2c+2$ must be tiled with $s$ sequences of horizontal dominoes with labels $1, \ldots, 2j.$  Since all labels are at most $n+1$ this forces $j \le \frac{n+1}{2}$ in the case where $s > 0.$

Now consider the first $t$ columns of the lower part of the array, where the base consists of vertical $1$-dominoes.  Here we can apply the analysis of Proposition \ref{comi+omi+1} to find sequences $a_k > \cdots > a_1$ and $s_k > \cdots > s_1$ such that the contribution to the weight from the lower part is
\begin{equation}\label{17}
\begin{array}{l}
(1^{s+t},\ldots,a_1^{s+t},(a_1+1)^{s+t-s_1},\ldots, a_2^{s+t-s_1},(a_2+1)^{s+t-s_1-s_2}, \ldots, a_k^{s+t-s_1-\cdots -s_{k-1}},\ldots, \\
(a_k+1)^{s+s_1+\cdots +s_k}, \ldots , b_k^{s+s_1+\cdots s_k}, (b_k+1)^{s+s_1+\cdots s_{k-1}}, \ldots , b_{k-1}^{s+s_1\cdots s_{k-1}}, \ldots, (b_1+1)^s, \ldots, b_0^s).
\end{array}
\end{equation}
Here $a_k$ is the label of the largest vertical domino which appears in  columns 1 and 2 of the lower part
of the array.  The column decreasing condition implies that all dominoes in columns 1 and 2 in the upper part of the array are greater than $a_k.$
Therefore it follows from the above expression for the weight that the sequences $a_1, \ldots, a_k$ and $s_1, \ldots, s_{k-1}$ are uniquely determined by the weight.  Also $s_k$ is determined since $t = 2s_1 + \cdots + 2s_k.$  

We have therefore shown that if two tilings yield the same weight than the tilings must agree on  the lower part of the array.    As the tilings of the lower part agree the  contributions to the weight from the upper part of the array must also agree.

 The  first two columns each have height $2n$ and the largest possible label for any tile is $n+1.$  Suppose that for each possible label $1\le i \le n+1$ the multiplicity of tiles in columns 1 and 2 with label $i$ is $c_i$.  Then  $ c_i \le 2$ and we have an equation $2n = \sum_{i = 1}^{n+1} c_i.$  It follows that there are at most two labels appearing with multiplicity 1.  Moreover, in view of the column decreasing condition, any horizontal tile with label $i$ must satisfy $c_i = 1.$
  
First suppose that $c_i =0$ for some $i.$  Then $c_k = 2$ for all $k \ne i$ and the conditions imply that the only possible tiling has all vertical tiles with labels
$(1,1), \ldots, (i-1,i-1), (i+1,i+1), \ldots $ where we allow the case $i = n+1$ with the sequence ending at $(n,n)$.  So  in this case there is a unique tiling.
 
 Now assume $c_i \ne 0$ for all $i$. This implies that  there are precisely two labels appearing with multiplicity 1 and all others appear with multiplicity 2.  Moreover, since there are an even number of rows in both the top and bottom parts, the exceptional tiles are either both in the top section or the bottom section. 

 Suppose the latter occurs. Then the analysis in Proposition \ref{comi+omi+1} implies that columns 1 and 2 of the lower part  have pairs of vertical tiles with labels  $(1, 1) \ldots, (j-1,j-1)$ followed by a horizontal $(j+1)$-tile lying over a horizontal $j$-tile.  This forces the upper part of the array to be
 tiled with pairs of vertical dominoes with labels $(j+2,j+2), \ldots, (n+1,n+1).$  Here again there is a unique labelling.
 
 Otherwise, columns 1 and 2 of the lower part of the array are tiled with pairs of vertical dominoes 
 with equal labels $(1,1), \ldots, (j,j)$ and there exist integers $ j < x < y \le n+1$ such that $x,y$ are the only 
 labels appearing with multiplicity 1 in the top part of the array.  The tiling can be done in two ways.
 
 One possibility is that the upper part is tiled
 $$(j+1,j+1), \ldots,(x-1,x-1), (x,x+1),(x+1,x+2), \ldots, (y-1,y), (y+1,y+1), \ldots, (n+1,n+1),$$
 where we allow the possibilities that  $x = j+1,$ omitting  the first sequence of tiles,
 and $y =n+1,$ omitting the last sequence of tiles.
 
 Alternatively, the upper part begins with pairs of  vertical tiles with labels $(j+1,j+1) \ldots, (x-1,x-1)$, then
 a horizontal $x$-tile, then pairs of  vertical tiles with labels $(x+1,x+1) \ldots, (y-1,y-1)$, then a horizontal
 $y$-tile, then pairs of equal vertical tiles going to the top of the array.  Here we again allow the case $x = j+1$, omitting the first sequence of tiles.
 
  It follows from the above  that when there are two tilings giving the  same weight then one tiling has precisely two more horizontal tiles than the other.  Therefore one tiling  corresponds to an alternating composition factor and the other corresponds to a symmetric composition factor of $\d \otimes \d$.  \hal

\section{Remaining Table \ref{TAB1} configurations}\label{remain}

 \begin{prop}\label{2om12om2}  Let $W$ be the irreducible module for $A_n$ ($n\ge 2$) with highest weight $\d = 2\o_1 + 2\o_2$. 
 Then  $\wedge^2(W)$ is MF.
\end{prop}

\pf  If $n \le 5$ we  check the assertion using Magma.  Since this is the first time (of many) that we are using Magma, and so that the reader can reproduce the computations if they wish, we give the Magma code for this computation in the case of $A_4$:

\begin{tabular}{l}
A4:=  RootDatum(``A\_4" : isogeny:= ``SC"); \\
W :=   LieRepresentationDecomposition(A4, [2,2,0,0]); \\
X:= AlternatingPower(W,2); \\
Weights(X); 
\end{tabular}

\noindent The output of this is a list of weights which are the highest weights of the composition factors of $\wedge^2(W)$, together with a list of numbers giving the multiplicities of these composition factors. Since this list consists of a string of 1's, it follows that $\wedge^2(W)$ is MF, as claimed.

Note that the Magma code used in the rest of the paper is very similar to the example above; in particular the commands used are the ones given above, together with the ``SymmetricPower" and  ``TensorProduct" commands.

Now assume $n \ge 5$.   As before we proceed by the method of domino tilings.  
The highest weight $\d$  corresponds to the partition $(4,2,0, \ldots , 0)$ so the tableau to be tiled is $(8,8,4,4).$  Consider a fixed tiling. All $1 \times 1$ squares on the bottom row are labelled by 1's.  A column has at most 4 rows  and the labels must be strictly decreasing.  

It follows that the highest label of a tile is 4 and hence the weight of the tiling
has the form $(1^a,2^b,3^c,4^d)$ and the highest weight of the corresponding irreducible module is $(x,y,z,w,0\ldots ,0).$  
Thus the possible tilings are identical for $n =5$ and for $n > 5$.  The former
is MF by a Magma computation, and hence so is the latter. \hal

\begin{prop}  Let $W$ be the irreducible module for $A_n$ ($n\ge 2$) with highest weight $\d = 2\o_1 + 2\o_n$. 
 Then  $\wedge^2(W)$ is MF.
 \end{prop}

\pf  Small cases can be settled with a Magma computation.  So  assume $n \ge 6.$  The shape of the tableau to be tiled by dominoes   is $(8,8,4,\ldots, 4)$ where $4$ occurs $2n-2$ times.  Say the base has $t$ vertical 1-dominoes followed by $s$ horizontal 1-dominoes on top of which are
horizontal 2-dominoes.  Therefore $t + 2s = 8$ and $t+s$ is the multiplicity of 1 in a given weight. Note that $t$ is even. As  the group in question is $A_n$, the largest possible labelling is $n+1.$  

First assume that $s \ge 3$.  Then the base in columns 3 and 4 consists of a horizontal 2-domino lying over a horizontal 1-domino.  In view of the various conditions,  
 above the base in columns 3 and 4  there must exist pairs of vertical
dominoes with labels $3, \ldots, d$ followed by horizontal dominoes with labels $d+1, \ldots d+e$,
where we allow the possibility that there are no vertical dominoes or no horizontal dominoes.  Then we have equations $2 + 2(d-2) + e = 2n$
and $d+e \le n+1.$  The only possibility is that $d = n+1$ and above the base in columns 3 and 4  the dominoes are all vertical.

If $t = 0$, then columns 1 and 2 must be the same as columns 3 and 4.  So the tiling is determined, there are 8 horizontal dominoes, and the corresponding composition factor is symmetric.

Next assume $t = 2$  and consider columns $1$ and $2$.  Above the base there could be two vertical dominoes with labels $(2,2), (2,3),(3,3)$.  Or there could be a horizontal 3-domino on top of a horizontal 2-domino.  In fact, we argue that the latter cannot occur.  For otherwise, the weakly increasing and Y conditions rule out any possibility for tiles in the next two rows.  

 If there are two vertical $3$-dominoes then above the base these columns have pairs of vertical dominoes with labels $3, \ldots, n+1$ which  corresponds to an alternating factor since there are 6 horizontal dominoes.  Similarly if there is a vertical $2$-domino
 followed by a vertical $3$-domino, then above the base there are only vertical dominoes with labels $(2,3), \ldots, (n,n+1).$  Again the corresponding factor is alternating.
 
 Finally assume that above the base there is 
 a pair of vertical $2$-dominoes.   Then columns $1$ and $2$  begin with pairs of vertical dominoes labelled by
$1,2, \ldots, d$ followed by horizontal dominoes labelled $d+1, \ldots , d+e$ (allowing $e = 0$).
We then have $2d + e = 2n$ and $d+e \le n+1.$  Therefore, $d \ge n-1$ so that either $(d,e) = (n, 0)$
or $(n-1,2).$  In the former case we again get an alternating summand, whereas in the latter case the summand
is symmetric.

At this point we can assume $s \le 2,$ so that $t = 4$, 6 or 8.  In particular the base of columns 3 and 4 consists 
of vertical 1-dominoes.

Consider  columns $3$ and $4$. The possibilities (consistent with the various conditions) for what can occur   above the base are as follows: 

(1) two vertical $2$-dominoes. 

(2) a vertical $2$-domino followed by a vertical $3$-domino   ($s \ge 1$)

(3) two vertical $3$-dominoes  ($s \ge 2$)

(4) a horizontal $2$-domino followed by two vertical $3$-dominoes ($s \ge 1$)

(5) column $4$ is like one of the columns of (3) and columns $2$ and $3$ are as in (4)  ($s \ge 2$).

\noindent If (1) holds then columns $3$ and $4$  begin with pairs of vertical dominoes labelled by
$1,2, \ldots, d$ followed by horizontal dominoes labelled $d+1, \ldots , d+e$ (allow $e = 0$).
We then have $2d + e = 2n$ and $d+e \le n+1.$  Therefore, $d \ge n-1$ so that either $(d,e) = (n, 0)$
or $(n-1,2).$

If (2) holds then there is a sequence of pairs of vertical dominoes labelled $(1,1), (2,3), \ldots, $ \linebreak $(d,d+1)$
followed by horizontal dominoes with labels $d+2, \ldots , d+e$ (allowing $e = 0$).  We  find that
$d = n$ is the only possibility.

If (3) holds, the only possibility is  a sequence of pairs of vertical dominoes with labels $1,3,\ldots, n+1.$

If (4) holds, then the only possibility is that columns $3$ and $4$ begin with a pair of vertical $1$-dominoes, followed by a horizontal $2$-domino, followed by pairs of vertical dominoes with labels $3, \ldots, n$,
followed by a horizontal $n+1$-domino.

And if (5) holds, column $4$ has all vertical dominoes with labels $1,3, \ldots, n+1$ and columns $2$ and $3$
are as described in the last paragraph.

Suppose $t = 8$. Then (1) must occur.  The weakly increasing condition implies that each of the columns begins with a sequence of pairs of vertical dominoes with labels $1, \ldots, d$.  If $d = n$, then
the tiling is determined, there are no horizontal dominoes, and the weight corresponds to a  symmetric factor.   Suppose
$d = n-1$.  Then  the sequence of vertical dominoes in columns  1 and 2 ends at either $n$ or $n-1$. Only
the former will  yield an alternating summand.

Now consider $t = 4$ or $6.$   Suppose columns $3$ and $4$ have a pair of vertical dominoes just above the base as described in (1), (2) or (3). If $t = 4$ these can have labels  $(2,2), (2,3),$ or $(3,3). $  If $t = 6$ only the first two are possible. 

First suppose that (1) holds.  Then the weakly increasing condition implies that the level above the base
consists of only vertical $2$-dominoes and the analysis of (1) applies equally to the pairs of columns $(1,2)$  and $(3,4).$   Therefore the possible configurations are clear from the
analysis of (1) above.  The tiling is entirely determined by $t$ and the multiplicity of $n+1$ and we
see that there is precisely one alternating possibility.

Now suppose that (2) holds. Again the weakly increasing condition implies that there are vertical $2$-dominoes above the base in columns $1$ and $2$.  The analysis of (2) shows that above the base columns $3$ and $4$ consist of pairs of vertical dominoes with labels $(1,1), (2,3), \ldots, (n,n+1).$  And applying the analysis of (1)
to  columns 1 and 2 we see that these columns consist of vertical dominoes with labels $1, \ldots, n$. So there is an alternating possibility only for $t = 6$.

Suppose (3) holds.  Here the Y-condition forces columns 1 and 2 to have vertical 2-dominoes  above the base,
so the analysis of these columns is as for (1).   And the analyis of (3) shows that columns $3$ and $4$ consist of pairs of vertical dominoes with labels $(1,1), (3,3) \ldots, (n+1,n+1).$  Once again there
is a unique alternating possibility.

If (4) holds, then  columns $3$ and $4$ are uniquely determined as described above.  Columns $1$ and $2$
could either be of type (1) or of type (4). In the latter case $s \ge 2$ by the Y-condition, so $s = 2$ and $t = 4.$
In this case the tiling corresponds to a symmetric factor.  Otherwise,  the remaining pairs of columns are of type (1)  and the weight uniquely determines whether the factor is symmetric or alternating.

Finally suppose that (5) holds. The Y-condition implies that $s \ge 2$ in this situation and hence $s = 2$ and $t = 4.$  Then the tiling in columns $2,3$ and $4$ is  determined. Column $1$ must be labelled with vertical dominoes $1,2, \ldots, n$ and the resulting composition factor is alternating.

At this point it remains to rule out multiplicities of configurations which correspond to alternating summands.   For two tilings with the same weight the values of $s$ and $t$ are determined.  If $t=8$ we have seen that there is only one possible configuration. If $t=6$ the above analysis shows that there are four possibilities, with the following weights:
\[
\begin{array}{l}
1^7,2^5,3^4,\ldots, (n-1)^4,n^2,(n+1)^2 \\
1^7,2^5,3^4,\ldots, n^4 \\
1^7,2^4,\ldots, n^4,n+1 \\
1^7,2^4,\ldots, (n-1)^4,n^3,(n+1)^2.
\end{array}
\]
And if $t=4$ there are also four possible weights:
\[
\begin{array}{l}
1^6,2^6,3^4,\ldots,(n-1)^4, n^3,n+1 \\
1^6,2^4,\ldots,(n-1)^4, n^3,(n+1)^3 \\
1^6,2^5,3^4,\ldots, n^4,n+1 \\
1^6,2^4,3^4,\ldots, n^4,(n+1)^2.
\end{array}
\]
From this it is apparent that each possible weight determines a unique summand, 
and this establishes the result.  \hal

\begin{prop}  Let $W$ be the irreducible module for $A_n$ ($n\ge 4$) with highest weight $\d = \o_2 + \o_{n-1}$. 
 Then  $\wedge^2(W)$ is MF.
\end{prop}

\pf This is easily verified by a Magma computation for small values of $n$.  We will assume $n \ge 7$  and proceed via a domino argument.  
The highest weight $\d$ corresponds to the partition $(2,2,1,\ldots,1,)$.  Doubling and repeating entries
we obtain the shape of the tableau $(4,4,4,4,2,\ldots, 2)$, where 2 occurs  $2(n-3)$ times.  Therefore the array will have $2n-2$ rows and must  be tiled with $2n+2$ dominoes.  The largest possible label for a domino is $n+1$ and in order to have an alternating representation the number of horizontal dominoes must be twice an odd number.

The bottom four rows of the array each have four $1 \times 1$ squares.  Above these rows there are $2n-6$ rows each with two $1 \times 1$ squares.  We will work through the possible domino tilings.  The possibilities
for the bottom four rows are  straightforward.  The remaining part of the array is a $2 \times (2n-6)$ matrix.
The labelling of this part is constrained by the fact  that the largest possible is $n+1$ and the fact that the irreducible summands must be alternating. Typically we will see that the base of the $2 \times (2n-6)$ array may begin with one or two horizontal tiles and may end with a few horizontal tiles, but otherwise it consists of vertical tiles.  This will follow from the argument  of Lemma \ref{horizontal}, which is sometimes applied implicitly.  Details will be provided in the various cases.

 For future reference we first list  the highest weights of the composition factors followed by the weight of a corresponding tiling:

1.  $(020\ldots 0101) -  (1^4,2^4,3^2,\ldots, (n-2)^2,(n-1)^1,n^1)$

2.  $(1010\ldots 0 20) - (1^4,2^3,3^3,4^2,\ldots ,(n-2)^2,(n-1)^2)$

3.  $(110\ldots 0100) - (1^4,2^3,3^2,\ldots, (n-2)^2,(n-1)^1,n^1(n+1)^1)$

4.  $(0010\ldots 0100) - (1^3,2^3,3^3,4^2,\ldots, (n-2)^2,(n-1)^1,n^1,(n+1)^1)$

5.  $(10 \ldots 01) - (1^3,2^2, \ldots, n^2,(n+1)^1)$

6.  $(1010 \ldots 01000) - (1^4,2^3,3^3,4^2,\ldots, (n-3)^2,(n-2)^1,(n-1)^1,n^1,(n+1)^1)$

7.  $(110 \ldots 011) - (1^4,2^3,3^2,\ldots, (n-1)^2,n^1)$

8.  $(0010\ldots 011) - (1^3,2^3,3^3,4^2,\ldots, (n-1)^2,n^1)$

9.  $ (00010 \ldots 0101) -  (1^3,2^3,3^3,4^3,5^2,\ldots, (n-2)^2,(n-1)^1,n^1)$

10.  $ (010 \ldots 010) - (1^3,2^3,3^2,\ldots, (n-1)^2,n^1,(n+1)^1)$

11.  $ (20 \ldots 010) - (1^4,2^2,\ldots, (n-1)^2,n^1,(n+1)^1)$

12.  $ ( 010 \ldots 02) - (1^3,2^3,3^2,\ldots, n^2)$.

At this point we must indicate why these occur and show that there are no additional configurations.
The possibilities for the base (bottom two rows) are as follows: 

(a) 4 vertical 1-dominoes; 

(b) two vertical 1-dominoes followed by a horizontal 2-domino lying over a horizontal 1-domino; 

(c) two pairs of horizontal 2-dominoes lying over horizontal 1-dominoes.

\noindent Suppose that (a) occurs.  Then there are three possibilities for the next two rows:  4 vertical 2-dominoes; 2 vertical 2-dominoes followed by a horizontal 3-domino above a horizontal 2-domino; two pairs of a horizontal 3-domino lying over a horizontal 2-domino.

Continuing with (a), assume that the next layer  has four vertical 2-dominoes.  If there is a horizontal domino at the base of the $2 \times (2n-6)$ array, then Lemma \ref{horizontal} implies that all the dominoes of this array are horizontal, but
this contradicts the fact that the labelling of dominoes is bounded by $n+1.$   Indeed the same argument shows that above the bottom four rows of the array the remaining labelling consists of pairs of vertical dominoes with labels $(3,3), \ldots, (n-2,n-2)$ and above this there is a horizontal $(n-1)$-domino followed by a horizontal $n$ domino.  This gives example 1. 

  Next assume that
the layer above the base contains two vertical 2-dominoes followed by a horizontal pair of a 3-domino above
a horizontal 2-domino.  Above this there could be a horizontal 3-domino or a pair of vertical dominoes.  In the first case the 3-domino must be followed by a pair of vertical dominoes with labels $(4,4), \ldots, (n-2,n-2),$ 
and this is followed by horizontal dominoes with labels $n-1, n, n+1.$ We note that Y-condition in this case is satisfied due to the existence of   the 3-domino above the base. This gives example 3.  Otherwise there is a sequence of pairs of vertical dominoes  with labels $(3,3), \ldots, (n-1,n-1),$ or $(3,3), \ldots, (n-3,n-3),$  or $(3,4), \ldots,  (n-1,n).$ The vertical dominoes continue
to the top of the array  in the first case and the last case.  In the other case the  array  ends with 4 horizontal dominoes having labels $n-2, n-1, n, n+1.$  This yields the examples 2, 6 and 7. 

 If the two rows above the base have two pairs of horizontal dominoes then the argument
of Lemma \ref{horizontal} shows that what follows are pairs of vertical dominoes with labels $(4,4), \ldots (n-1,n-1)$ followed by horizontal dominoes with labels $n$ and $n+1.$  This is example 11.

Now suppose that  (b) holds. Then the Y-condition implies that above the base in columns 3 and 4 there must be a horizontal 4-domino above a horizontal 3-domino.  In the first two columns there could be two vertical 2-dominoes or possibly a horizontal 3-domino above a horizontal 2-domino.  Suppose the latter case holds.  Then above this level there must be a sequence of pairs of vertical dominoes with labels $(4,5), \ldots, (n,n+1).$  This yields example 5.
Now assume that the former case holds.

First assume that above the vertical 2-dominoes there is a pair of vertical dominoes.  Then there must be a string of vertical dominoes with  labels $(3,3), \ldots , (n-1,n-1)$, or $(3,3), \ldots , (n-2,n-2)$    or $(3,3), (4,5), \ldots (n-2,n-1),$  and these are  followed either by a single horizontal $n$-domino, or by two horizontal dominoes with labels $(n-1,n)$ or $(n,n+1)$, respectively.  These give examples 8, 9 and 4, respectively.  

Now assume that above the 2-dominoes there is a horizontal domino. If there is a single horizontal 3-domino, then above this there must be a string of vertical dominoes with labels $(4,5), \ldots, (n-1,n)$  followed by a horizontal $n+1$-domino.  This gives example 10. Otherwise there must be a pair of horizontal dominoes with labels 3 and 4.  Following this there are  pairs of vertical dominoes with labels $(5,5), \ldots, (n,n).$ and this gives example 12.  
 
 Finally we claim (c) cannot occur. Arguing with the Y-condition we see that above the base there must be two pairs of a horizontal 4-domino lying over a horizontal 3-domino.  If there is a horizontal domino above this level in columns 1 and 2, then the argument of Lemma \ref{horizontal} shows that  all the remaining dominoes would be  horizontal, which is not possible since the largest possible label of a domino is $n+1.$  The only other possibility consistent with the conditions is a sequence of pairs of vertical dominoes with labels $(5,5), \ldots, (n+1,n+1).$
But this yields the trivial representation which is symmetric.  This establishes the claim and completes the proof.
 \hal

\begin{prop}  Let $W$ be the irreducible module for $A_n$ ($n\ge 5$) with highest weight $\d = \o_2 + \o_4$. 
 Then  $\wedge^2(W)$ is MF.
\end{prop}

\pf  This is a routine Magma computation for $n \le 7$, so assume $n \ge 8.$ We start with the partition $(2,2,1,1,0,\ldots,0)$ which corresponds to $\d$ and then double and repeat
exponents to get the tableau $(4,4,4,4,2,2,2,2)$ to be tiled with dominoes.  
With such an easy array it is  a simple matter to list all the possible tilings.
The list below gives all possible composition factors and corresponding weights.
 
1.  $(001110\ldots 0) -  (1^3,2^3,3^3,4^2,5^1)$

2.  $(0101010\ldots 0) - (1^3,2^3,3^2,4^2,5^1,6^1)$

3.  $(00001010\ldots 0) - (1^2,2^2,3^2,4^2,5^2,6^1,7^1)$

4.  $(0020010\ldots 0) - (1^3,2^3,3^3,4^1,5^1,6^1)$

5.  $(10120 \ldots 0) - (1^4,2^3,3^3,4^2)$

6.  $(12000010 \ldots 0) - (1^4,2^3,3^1,4^1,5^1,6^1,7^1)$

7.  $(010020 \ldots 0) - (1^3,2^3,3^2,4^2,5^2)$

8.  $ (1110010 \ldots 0) - (1^4,2^3,3^2,4^1,5^1,6^1)$

9.  $(11011\ldots 0) - (1^4,2^3,3^2,4^2,5^1)$

10.  $ (2001010 \ldots 0) - (1^4,2^2,3^2,4^2,5^1,6^1)$

11.  $ (10010010 \ldots 0) - (1^3,2^2,3^2,4^2,5^1,6^1,7^1)$

12.  $ (1000110 \ldots 0) - (1^3,2^2,3^2,4^2,5^2,6^1)$

13.  $ (101000010 \ldots 0) - (1^3,2^2,3^2,4^1,5^1,6^1,7^1,8^1)$

14.  $ (021010 \ldots 0) - (1^4,2^4,3^2,4^1,5^1)$

15.  $ (01100010 \ldots 0) - (1^3,2^3,3^2,4^1,5^16^1,7^1)$.  \hal

\begin{prop}\label{wedgeadj}  Let $W$ be the irreducible module for $A_n$ ($n\ge 2$) with highest weight $\d = \o_1 + \o_{n}$. 
 Then  $\wedge^3(W)$ is MF.
\end{prop}

\pf This is easily checked by a Magma computation for $n \le 6.$
So assume $n \ge 7.$ The highest weight of any composition factor of  $\wedge^3(W) \downarrow X$ is subdominant to $\l = 3\o_1+3\o_n=(30\ldots 03).$  However the dominant
weights $\l, \l-\a_1$, and $\l - \a_n$ clearly cannot occur in in $\wedge^3(W).$
The remaining dominant weights that are subdominant to $\l$ are listed below (where as usual $\a_1,\ldots ,\a_n$ form a system of fundamental roots for $X = A_n$):

 \vspace{2mm}
1. $(110\ldots 011) = \l - \a_1 - \a_n$

2. $(30\ldots 0100) = \l - \a_{n-1}-2\a_n$

3. $(0010\ldots 03) = \l - 2\a_1-\a_2$

4. $(110\ldots 0100) =  \l - \a_1-\a_{n-1}-2\a_n$

5.  $(0010\ldots 011) = \l - 2\a_1-\a_2 -\a_n$

6. $(0010\ldots 0100) = \l - 2\a_1-\a_2 -\a_{n-1} - 2\a_n$

7. $(20\ldots  02) = \l - (\a_1 + \cdots +\a_n)$

8.  $(10\ldots  01) = \l - 2(\a_1 + \cdots +\a_n)$

9.  $(010\ldots  02) = \l - (2\a_1 + \a_2 + \cdots  +\a_n)$

10. $(20\ldots  010) = \l - (\a_1 + \a_2 + \cdots +\a_{n-1} + 2\a_n)$

11. $(010\ldots  010) = \l - (2\a_1 + \a_2 + \cdots + \a_{n-1} + 2\a_n)$

12. $(0 \ldots  0) = \l - 3(\a_1 + \cdots +\a_n)$.

\vspace{2mm}

Weights 1,2,3 each occur with multiplicity one in $\wedge^3(W)$ and are not subdominant
to any other dominant weight.  Therefore there is a  composition factor of $\wedge^3(W)$
for each of these highest weights appearing with multiplicity 1.  Next consider weight 4.
This occurs as the wedge of weight vectors  with the following weights: $(\d, \d - \a_1 -\a_n, \d - \a_{n-1}-\a_n)$ and
$(\d-\a_n,  \d - \a_1, \d - \a_{n-1}-\a_n),$ so this weight  space has dimension 2  in $\wedge^3(W).$  But it appears in the irreducible of highest weight 1 with multiplicity 2.  So there is no composition factor with highest weight 4.  Similarly for 5.

Now consider weight 6. This occurs as the wedge of weight vectors as follows:
\[
\begin{array}{l}
(\d, \d-\a_1-\a_n, \d-\a_1-\a_2-\a_{n-1}-\a_n),\; (\d, \d-\a_1-\a_2 -\a_n, \d-\a_1-\a_{n-1}-\a_n), \\
(\d-\a_1, \d-\a_{n-1} -\a_n, \d-\a_1-\a_2-\a_n),\; (\d-\a_n, \d-\a_1 -\a_2, \d-\a_1-\a_{n-1}-\a_n), \\
(\d-\a_1-\a_n, \d-\a_1-\a_2, \d-\a_{n-1}-\a_n),\; (\d-\a_1, \d-\a_n, \d - \a_1-\a_2 - \a_{n-1} - \a_n).
\end{array}
\]
It follows that this weight space has multiplicity 6.  It appears in the composition factor with highest weight 1 with multiplicity 4 and in the composition factors of highest weights 2 and 3, each with multiplicity 1.  Therefore there is no composition factor with highest weight 6 in $\wedge^3(W)$.

At this point we aim to show that for the remaining weights, 7-12, there is a composition factor in $\wedge^3(W)$ appearing with multiplicity 1. For 9 and 10 this follows from Lemma \ref{wedge3adj} in the next chapter which is proved independently of results in this chapter.  And a simple counting argument   shows that 7 also appears with multiplicity 1.  Therefore we must consider 8, 11, 12. Suppose the multiplicities of these in $\wedge^3(W)$ are $a,b,c$, respectively.  The multiplicity of 12 in the tensor cube of $\d$ is 2, hence $c\le 2$.

To settle these cases we first decompose $(110\ldots  0) \otimes (0 \ldots  011).$  Of course the multiplicity of the highest weight
is 1.  Easy applications of Theorem \ref{LR} show that the  irreducibles with highest weights  7,9,10,11, and 12 also occur with multiplicity 1, while the irreducible with highest weight   8 appears with multiplicity 2.   Moreover,  there are no composition factors of highest weights 4, 5 or 6.  This covers all subdominant weights so that, denoting irreducible modules by their highest weights as listed in 1--12, 
\[
\wedge^3(\d) - \left((110\ldots  0) \otimes (0 \ldots  011)\right) = 2 + 3 + 8^{a-2} + 11^{b-1} + 12^{c-1}
\]
(here we are using the notation $A-B$ for modules $A\subseteq B$, to denote a module with the same composition factors as the quotient $A/B$).  

We  now use dimension arguments.  We obtain the dimension of $(110\ldots  0)$ by viewing it as $((10\ldots  0) \otimes 
(010\ldots  0)) - (0010\ldots  0).$  Squaring we have the dimension of $(110\ldots  0) \otimes (0 \ldots  011).$  Similarly we can obtain the dimension of 2 and 3 by viewing them as alternating tensor products.
For example 2 can be viewed as 
\[
\left((30\ldots  0) \otimes (0\ldots  0100)\right) - \left((20\ldots  0) \otimes (0\ldots  010)\right) + \left((10\ldots  0) 
\otimes (0\ldots  01)\right) - 0.
\]
At this point we can  compute the dimension  
of $\wedge^3(\d)$ and subtract the dimension of $(110\ldots  0) \otimes (0 \ldots  011)$ as well as the dimensions of 2 and 3. The result is $-n^2-2n$.  Therefore we have the equation
\[
-n^2-2n = (a-2)(n^2+2n) + (b-1)\dim (010\ldots  010) + (c-1).
\] 
Now $\dim (010\ldots  010)$ is ${{n+1}\choose 2}^2- (n+1)^2$ and this forces $b = 1$. In turn this forces
$c = 1$ and hence $a = 1$ as well, completing the proof.  \hal

\begin{prop}  If $X = A_2$ and $W = V_X(\om_1 + \om_2),$ then $S^3(W)$ is MF.
\end{prop}
\pf This is immediate using a Magma calculation. \hal

\section{Table \ref{TAB2} configurations}

We next consider entries in Table \ref{TAB2} of Theorem \ref{MAINTHM}. We use the following notation:
\[
\begin{array}{l}
X = A_{l+1}, \\
W = V_X(\d), \\
Y = SL(W) = A_n.
\end{array}
\]   
As usual write $\o_i,\a_i$ ($1\le i \le l+1$) and $\l_i,\b_i$ ($1\le i\le n$) for the fundamental dominant weights and roots of $X$ and $Y$ respectively.

\begin{prop}  Let $\d = c\om_i$. Then $V_Y(\l_1+\l_n) \downarrow X$ is MF.
\end{prop}

\pf  Observe that $V_Y(\l_1+\l_n)$ is a direct summand of $W \otimes W^*.$  
Since Proposition \ref{stem1.1.A} shows that  $(W \otimes W^*) \downarrow X$ is MF, the result follows.  \hal

\begin{prop}\label{l2l3,etc} {\rm (i)} If $\d = 2\om_1$ or $\om_2$,  then $V_Y(\l) \downarrow X$
is MF for the following highest weights $\l$:
\[
\begin{array}{l}
\l = \l_1 + \l_i \, (2 \le i \le 7), \\
\l = \l_2 + \l_3, \\
\l = 2\l_1+\l_2, \\
\l = 3\l_1 + \l_2.
\end{array}
\]
{\rm (ii)}  If $\d = 3\om_1$ or $\om_3,$  then $V_Y(\l_1+\l_2) \downarrow X$ is MF.
\end{prop}

\pf  We first note that
\[
V_Y(\l_1 )\otimes V_Y(\l_i) \cong V_Y(\l_1 + \l_i) + V_Y(\l_{i+1}).
\]
We shall need a number of isomorphisms like this, and abbreviate notation by writing just $\l$ instead of the module $V_Y(\l)$, and also using subtraction notation $A-B$, for modules $B\subseteq A$, to denote a module with the same composition factors as $A/B$. With this notation the above isomorphism can be written $\l_1+\l_i= (\l_1\otimes \l_i )-\l_{i+1}$. We list some similar such isomorphisms:
\begin{equation}\label{18}
\l_1 + \l_i = (\l_1\otimes \l_i) - \l_{i+1}
\end{equation}
\begin{equation}\label{19}
\l_2 + \l_3 = (\l_2 \otimes \l_3) - (\l_1 \otimes \l_4)
\end{equation}
\begin{equation}\label{20}
2\l_1 + \l_2 = (2\l_1 \otimes \l_2) - (\l_1 \otimes \l_3) + \l_4
\end{equation}
\begin{equation}\label{21}
3\l_1 + \l_2 = (3\l_1 \otimes \l_2) - (2\l_1 \otimes \l_3) + (\l_1 \otimes \l_4)- \l_5
\end{equation}
These can each be checked using Theorem \ref{LR}.
Using the method described in Section 4.4 of \cite{Howe} we determine explicitly the decompositions of $\wedge^j(\om_2)$ and 
$\wedge^j(2\om_1)$ for $j \le 7$.
These are all quite simple. For each $j \le 7$ there is a small number of irreducible summands each with support  on the first $j+1$ roots and of form  $j\d-\sum_{i\le j}c_i\a_i$.  The entries are unchanged for $n > j$. We illustrate with the most complicated case $j =7$:
$$\wedge^7(\om_2) = (6\om_1 +\om_8) \oplus (4\om_1 +\om_3 +\om_7) \oplus (3\om_1 
+\om_2 +\om_4 +\om_6 ) \oplus (2\om_2 +2\om_5)\oplus (\om_1 +2\om_4 + \om_5),$$
$$ \wedge^7(2\om_1) = (7\om_1 +\om_7)\oplus (4\om_1 +2\om_2 +\om_6)\oplus (2\om_1 
+2\om_2 + \om_3 + \om_5)\oplus (\om_1 +3\om_3+\om_4) \oplus (3\om_2 +2\om_4).$$
   At this point  a Magma computation can be easily used to compute the tensor product with the  other wedge powers appearing in (6.15) and (6.16), and thereby application of Theorem \ref{LR} settles the  cases $\l = \l_1+\l_i$  and $\l = \l_2+\l_3$  in (i).

   For the  cases $\l = 2\l_1+\l_2$ and $\l =3\l_1+\l_2$ we use (6.17) and (6.18) where we must decompose  the restrictions of symmetric powers $2\l_1$ and $3\l_1$.  These  are summands of $\bigotimes^2\d$ and $\bigotimes^3\d$, respectively,  and these tensor products can be decomposed using Theorem \ref{LR}. As above we find each irreducible summand has a simple expression with
   support on the first $6$ roots and is of the form $a\delta-\sum_{i\leq 5}c_i\alpha_i$. Therefore  the same holds for the
   irreducible summands of $S^2(2\om_1),$ $S^3(2\om_1)$, $S^2(\om_2),$ $S^3(\l_1)$.
   At this point we are in a position to use Magma
   to tensor the terms with the restriction to $X$ of  $\l_2, \l_3$ and proceed to settle the cases $2\l_1+\l_2$ and $3\l_1+\l_2$ as above. Therefore (i) holds.

Finally  consider (ii).  From (6.15) we have $\l_1 + \l_2 = (\l_1 \otimes \l_2)-\l_3$.
Using Theorem \ref{dominoes} we compute the irreducible summands of $\wedge^2(\d)$ and we find that each summand has a simple expression with support on $\a_1, \ldots , \a_6$. The possibilities are independent of $n\ge 6$. Given that the support of the weight $\delta$ is on $\alpha_i$ for $i\le 3$, we can now establish (ii) using Magma checks for $n \le 9$, as in the previous cases. \hal


  

\begin{prop}\label{l1l(n+2-i)} Let $\d = 2\om_1$ or $\om_2.$  Then $V_Y( \l_1 + \l_{n+2-i}) \downarrow X$
is MF for $2 \le i \le 7.$
\end{prop}

\pf  Taking duals it will suffice to show that $V_Y(\l_{i-1} + \l_n) \downarrow X$ is MF.  To ease notation we
set $j = i-1$, so that $j \le 6.$ Then
\begin{equation}\label{23}
\l_j + \l_n = (\l_j \otimes \l_n) - \l_{j-1}.
\end{equation}
As in the proof of Proposition \ref{l2l3,etc}, we see that $\wedge^j(\om_2)$ and $\wedge^j(2\om_1)$ have support on the first $7$ roots, and the irreducible summands are of a very simple nature and they are unchanged for  $n >j$.
Now consider $(V_Y(\l_j) \otimes V_Y(\l_n)) \downarrow X$.  To decompose this tensor product we use Theorem \ref{LR} in the usual way.  Let $\e = \l_n \downarrow X$ with corresponding partition $(1,\ldots ,1,0)$ or $(2,\ldots ,2)$, depending on 
whether $\d = \om_2$ or $\d = 2\om_1.$  Now consider possible Littlewood-Richardson skew tableaux of shape $\nu/\e$ corresponding to composition factors of $V_Y(\l_j \otimes \l_n) \downarrow X.$  In such a tableau the crossed out squares correspond to $1^1\ldots  l^1$ or $1^2 \ldots  (l+1)^2,$  respectively.
In either case they form a simple rectangle of width 1 or 2 and length $l$ or $l+1$, respectively.  
Consider a labelling of such a tableau corresponding to a composition factor.  In view of the previous paragraph
all of the labels appear either in the first 7 rows or the last row. It therefore follows that all possible configurations are already apparent when $l+1 = 8.$  Consequently, we can  use Magma to decompose the tensor product and then (\ref{23}) yields the result.  \hal 

\begin{prop}\label{l2l(n-1_}  If $\d = 2\om_1$ or $\d =\om_2$, then $V_Y(\l_2+\l_{n-1}) \downarrow X$ is MF.
\end{prop}

\pf  For these cases we decompose the restriction using Theorem \ref{LR}.  We begin with the observation that
\begin{equation}\label{24}
\l_2 + \l_{n-1} = (\l_2 \otimes \l_{n-1}) - (\l_1 \otimes \l_n)
\end{equation}
We will decompose the restrictions to $X$ of the tensor products in (\ref{24}) for the two cases $\d = \om_2$ and $\d = 2\om_1.$  

First assume $\d = \om_2$.  Then $\l_2 \downarrow X = 1010\ldots  0$ and $\l_{n-1} \downarrow X = 0 \ldots  0101.$  In terms of weights of partitions these become  $\g =(1^2,2^1,3^1)$  and $\e = (1^2, \ldots , (n-2)^2,(n-1)^1,n^1).$  To decompose the restriction of the first tensor product we consider all possible Littlewood-Richardson skew tableaux of shape $\nu/\e$ and weight $\g$.  This is relatively easy due to the simple nature of $\g.$  The possible weights corresponding to  $\nu$
are as follows: 

 \vspace{2mm}
$(1^4,2^3,3^3,4^2,\ldots ,(n-2)^2,(n-1)^1,n^1),$ \ \ \  $ (1^4,2^3,3^2,\ldots ,(n-1)^2,n^1),$ 

$(1^4,2^2, \ldots , (n-1)^2,n^1,(n+1)^1),$  \ \ \  $(1^4,2^3,3^2, \ldots , (n-2)^2,(n-1)^1,n^1,(n+1)^1),$

$(1^4,2^2,3^2,\ldots , n^2),$  \ \ \  $(1^3,2^3,3^2,\ldots ,(n-1)^2,n^1,(n+1)^1),$

$(1^3,2^3,3^3,4^2,\ldots ,(n-1)^2,n^1),$ \ \ \  $(1^3,2^3,3^3,4^2,\ldots ,(n-2)^2,(n-1)^1,n^1,(n+1)^1),$

$(1^3,2^3,3^2,\ldots , n^2),$ \ \ \ $(1^3,2^2, \ldots , n^2,(n+1)^1),$ \ \ \   $(1^2,\ldots , (n+1)^2).$
 
\noindent  In addition we check that the labelling for each $\nu$ is unique with the exception of \linebreak 
$(1^3,2^3,3^2,\ldots ,(n-1)^2,n^1,(n+1)^1)$ and $(1^3,2^2, \ldots , n^2,(n+1)^1).$  Each of these have
 two possible labellings.
 Therefore the first tensor factor in (\ref{24}) is
 \[
\begin{array}{ll}
(1010\ldots  0) \otimes (0 \ldots  0101) = &  (1010\ldots  0101) + (110 \ldots  011) + (20 \ldots  010) + \\
                                                                    &  (110 \ldots  0100) + (20 \ldots  02)  + (010 \ldots  010)^2 + (0010\ldots  011)+ \\
                                                                   &  (0010 \ldots  0100) + (010\ldots  02) + (10 \ldots  01)^2 + (0 \ldots  0).
\end{array}
\]
 The restriction of the second tensor product in (\ref{24}) is just $(010\ldots  0) \otimes (0 \ldots  010)$ and
 applying Corollary \ref{LR_om2} again shows that this is  $ (010 \ldots  010) + (10 \ldots  01) + (0\ldots  0).$
 Therefore $V_Y(\l_2+\l_{n-1}) \downarrow X$ is MF when $\d = \om_2.$
 
 Now assume $\d = 2\om_1$  so that $\l_2 \downarrow X = 210\ldots  0$ and $\l_{n-1} \downarrow X = (0 \ldots  012).$ In terms of weights of partitions we have $\g =(1^3,2^1)$  and $\e = (1^3, \ldots , (n-1)^3,n^2).$ Here again we consider
 Littlewood-Richardson skew tableaux of shape $\nu/\e$ and weight $\g.$  It is relatively easy to list the possible choices of $\nu$ due to the simple nature of  $\g.$ The possible weights are as follows:
 
  \vspace{2mm}
  $(1^6,2^4,3^3, \ldots , (n-1)^3,n^2),$ \ \ \  $(1^6,2^3, \ldots , n^3),$ \ \ \  $(1^6,2^3, \ldots , (n-1)^3, n^2,(n+1)^1),$
  
  $(1^5,2^3, \ldots , n^3,(n+1)^1),$ \ \ \  $(1^5,2^4,3^3, \ldots , n^3),$ \ \ \ $(1^5,2^3, \ldots , (n-1)^3,n^2,(n+1)^2),$
  
  $(1^5,2^4,3^3,\ldots , (n-1)^3,n^2,(n+1)^1),$  \ \ \ $(1^4,2^4,3^3, \ldots , (n-1)^3,n^2,(n+1)^2),$
  
 $(1^4,2^4,3^3, \ldots , n^3,(n+1)^1),$  \ \ \ $(1^4,2^3, \ldots , n^3,(n+1)^2),$ \ \ \ $(1^3, \ldots , (n+1)^3).$ 
 
  \vspace{2mm}
  
\noindent We also check that the labelling for each $\nu$ is unique with the exception of $(1^5,2^3, \ldots ,$ \linebreak
$ n^3,(n+1)^1)$
and $(1^4,2^3, \ldots , n^3,(n+1)^2),$  where in each case there are two possible labellings.
Therefore in this case the first tensor factor in (\ref{24}) is 
\[
\begin{array}{ll}
 \vspace{2mm}
(210 \ldots  0) \otimes (0 \ldots  012) = & (210\ldots  012) + (30 \ldots  03) + (30 \ldots  011) + \\
                                                                & (20 \ldots  02)^2 + (110 \ldots  03) + (20 \ldots  010) + (110 \ldots  11) + \\
                                                                & (010 \ldots  010) + (010 \ldots  02) + (10 \ldots  01)^2 + (0 \ldots  0).
\end{array}
\]
The restriction of the second tensor product in (\ref{24}) is $(20\ldots  0) \otimes (0 \ldots  02)$ and this
 is easily seen to be $(20 \ldots  02) + (10 \ldots  01) + (0 \ldots  0).$  We therefore
 conclude that $V_Y(\l_2+\l_{n-1}) \downarrow X$ is also MF when $\d = 2\om_1.$  \hal

\begin{prop}\label{2l1ln}   If $\d = 2\om_1$ or $\om_2$, then $V_Y(2\l_1+\l_n) \downarrow X$ is MF.
\end{prop}

\pf  We proceed along the lines of the last result. We first note that
\begin{equation}\label{25}
2\l_1 + \l_n = (2\l_1 \otimes \l_n) - \l_1,
\end{equation}
so we must work out the restriction of the tensor product to $X$.  We can assume $n \ge 5$ as otherwise
we obtain the result from a Magma computation.

First assume $\d = \om_2$ so that $V_Y(2\l_1) \downarrow X = (020\ldots  0) + (00010\ldots  0).$  Now
we use Theorem \ref{LR} to decompose each of the terms on the right side tensored with $\l_n \downarrow T_X = (0\ldots  010).$ Using Corollary \ref{LR_om2}, we then obtain
\[
\begin{array}{c}
 (020\ldots  0) \otimes (0\ldots  010) =  (020\ldots  010) + (010\ldots  0) +(110\ldots  01) \\
(00010\ldots  0) \otimes (0\ldots  010) =  (00010\ldots  010) + (010\ldots  0) + (0010\ldots  01).
\end{array}
\]
Therefore, $V_Y(2\l_1 \otimes \l_n) \downarrow X$ only fails to be MF because of a repeated summand $\d =(010\ldots  0).$
Using  (\ref{25}) we obtain the result in this case.

Now assume $\d = 2\om_1$ so that $V_Y(2\l_1) \downarrow X = (40\ldots  0) + (020\ldots  0).$  This
time $\l_n \downarrow T_X = \d^* = (0\ldots  02)$.  Using Proposition \ref{pieri} we find that 
\[
\begin{array}{c}
 (40\ldots  0) \otimes \d^* = (40\ldots  02) + (30\ldots 01) +(20\ldots  0) \\
(020\ldots 0) \otimes \d^* = (110\ldots  01) + (20\ldots  0) + (020\ldots  02).
\end{array}
\]
As above, the only repeated composition factor in $V_Y(2\l_1 \otimes \l_n) \downarrow X$ is $\d = (20\ldots  0)$,
so that $V_Y(2\l_1 + \l_n) \downarrow X$ is MF.  \hal

\begin{prop}\label{3l1ln}   If $\d = 2\om_1$ or $\om_2$, then $V_Y(3\l_1+\l_n) \downarrow X$ is MF.
\end{prop}

\pf  The argument here is very similar to the proof in Proposition \ref{2l1ln}. As above we may assume that $n\ge 5$. 
Here we have
\begin{equation}\label{30}
3\l_1 + \l_n = (3\l_1 \otimes \l_n) - 2\l_1,
\end{equation}
First assume $\d = \om_2$ so that 
$$V_Y(3\l_1) \downarrow X = (030\ldots  0) + (0000010\ldots  0)+ (01010\ldots 0).$$
Now we use Corollary \ref{LR_om2} to decompose each of the summands on the right side above tensored
with $\l_n \downarrow X = (0\ldots  010)$. The results are as follows
\[
\begin{array}{ll}
 (030\ldots  0) \otimes \d^* = & (030\ldots 010) + (120\ldots 01) +(020\ldots  0) \\
 (000001\ldots  0) \otimes \d^* = & (000001\ldots 010) + (00001\ldots 01) +(0001\ldots  0) \\
 (01010\ldots  0) \otimes \d^* = & (01010\ldots 010) + (10010\ldots  01) +(1010\ldots  0) +\\
             &  (020\ldots 0) + (00010\ldots  0) + (0110\ldots  01).
\end{array}
\]
  It follows that $V_Y(3\l_1 \otimes \l_n) \downarrow X$ has precisely two irreducible summands appearing
 with multiplicity 2, namely $(020\ldots  0)$ and $(0001\ldots  0).$  But these sum to $V_Y(2\l_1) \downarrow X$,
 so that (\ref{30}) implies that $V_Y(3\l_1+\l_n) \downarrow X$ is MF.
 
 Now assume that $\d = 2\om_1$ so that 
$$V_Y(3\l_1) \downarrow X = (60\ldots  0) + (2200\ldots  0)+ (0020\ldots 0).$$
Once again we use Proposition \ref{pieri} to decompose each of the summands on the right side above tensored
with $\l_n \downarrow X = \d^* = (0\ldots  02)$. The results are as follows
\[
\begin{array}{c}
 (60\ldots  0) \otimes \d^* = (60\ldots 02) + (50\ldots 01) +(40\ldots  0) \\
 (220\ldots  0) \otimes \d^* = (220\ldots 02) + (310\ldots 01) +(40\ldots  0) + (210\ldots 0)+ (020\ldots  0) + (120\ldots  01)\\
 (0020\ldots  0) \otimes \d^* = (0020\ldots 02) + (0110\ldots 01) +(020\ldots  0).
\end{array}
\]
Therefore  $V_Y(3\l_1 \otimes \l_n) \downarrow X$ has  two irreducible summands appearing
 with multiplicity 2, namely $(40\ldots  0)$ and $(020\ldots  0).$  But these sum to $V_Y(2\l_1) \downarrow X$
 and again we conclude that $V_Y(3\l_1+\l_n) \downarrow X$ is MF.  \hal
 
 \section{Table \ref{TAB3} configurations}

Continue with the notation of the previous section -- that is, $X = A_{l+1}$, $W = V_X(\d)$ and $Y = SL(W) = A_n$.
The following two results can be found in  Theorems 3.8.1, 4.7.1 and Section 3.1 of \cite{Howe}.
As mentioned in the Introduction, Theorem \ref{symc} also follows from \cite[Theorem 3]{Kac}. 
A shorter proof of Theorem \ref{wedgec} was given by Stembridge in \cite{stem2}.

\begin{thm}\label{wedgec}   If $\d = \o_2$ or $2\o_1$, then $\wedge^c(W)\downarrow X$  is MF for all $c \ge 1.$
\end{thm}

\begin{thm}\label{symc}   If $\d = \o_2$ or $2\o_1$, then $S^c(W)\downarrow X$  is MF for all $c \ge 1.$
\end{thm}

We also record the following, which is immediate from Proposition \ref{stem1.1.A}.

\begin{prop} If $\d=c\o_i$, then $S^2(W)\downarrow  X$ and $\wedge^2(W)\downarrow X$ are both MF.
\end{prop}

For later use, in the next result we record a few of the composition factors of $S^c(W)$ in the case where $\d = \o_2$. This follows from \cite[3.8.1]{Howe}.

\begin{lem}\label{cfssc}
If $\d=\o_2$, then $S^c(W)$ contains composition factors $c\o_2$, $(c-2)\o_2+\o_4$ (for $c\ge 2, l\ge 3$), $(c-4)\o_2+2\o_4$ (for $c\ge 4,l\ge 3$) and $(c-3)\o_2+\o_6$ (for $c\ge 3,l\ge 5$).
\end{lem}

The next five results concern certain entries in Table \ref{TAB3} of Theorem \ref{MAINTHM}.  In each case the proof
is achieved by reducing consideration to a configuration where $X$ has bounded rank and then using
Magma. To do this we use Theorem \ref{LR}.  To illustrate the idea, suppose we are considering the tensor product of irreducible $X$-modules of highest weights $\mu$ and $\om$.  Suppose  $\mu$ has support on $\a_1, \ldots , \a_k$ and $\om = ab0 \ldots 0,$ so that a partition corresponding to $\om $ is $(a+b,b,0,\ldots ,0)$.  Then it follows
from Theorem \ref{LR}  and the shape of the partition for $\om$ that any irreducible constitutent of $\mu \otimes \om$ must have support on $\a_1, \ldots , \a_{k+2}.$ In the special case  $b = 0$,  the support is on  $\a_1, \ldots , \a_{k+1}.$

\begin{prop}\label{lambda3(020..)}  Suppose $l\ge 2$ and $\d = 2\o_2 = 020\dots 0$.  Then $V_{Y}(\l_3)\downarrow X$ is MF.
\end{prop}

\pf  We have $V_Y(\l_1 \otimes \l_2) = V_Y(\l_1 + \l_2) + V_Y(\l_3)$  so we we can view $V_{Y}(\l_3)\downarrow X$ as a submodule
of $\l_1 \otimes \l_2  \downarrow X = (020\ldots  0) \otimes \wedge^2(020\ldots 0).$  Using Theorem \ref{dominoes}
one checks that $\wedge^2(020\ldots 0) = (1210\ldots  0) + (10110\ldots  0).$  Now use Theorem \ref{LR} to
study $(020\ldots  0) \otimes  (1210\ldots  0)$ and   $(020\ldots  0) \otimes (10110\ldots  0).$
In each case we find that the highest weights of irreducible summands have support on $\a_1, \a_2, \ldots , \a_6.$
Consequently the same holds for the highest weights of all irreducible summands of $V_{Y}(\l_3)\downarrow X$
and hence it suffices to establish the result for $n \le 7.$   
For $n \le 7$ the result follows from  a Magma computation.  \hal

\begin{prop}\label{lamdai(com1)}  Let $\l = \l_i$ with $i=3, 4$ or $5$, and let $\d = (c0 \ldots  0)$ with $c\le 6,4$ or $3$, respectively. Then  $V_Y(\l) \downarrow X$ is MF.
\end{prop}

\pf  For $k \ge 1$ we have $V_Y(\lambda_1)\otimes V_Y(\lambda_k) = V_Y(\l_1 + \l_k) + V_Y(\l_{k+1})$.  Restricting to $X$ the left hand side is
$(c0 \ldots  0) \otimes \wedge^k(c0 \ldots  0).$  Suppose $\mu$ is an irreducible summand of
$\wedge^k(c0 \ldots  0)$.  If $\mu$ has support on $\a_1, \ldots , \a_j$, then it follows from Proposition \ref{pieri} that $(c0 \ldots  0) \otimes \mu$ has support on $\a_1, \ldots , \a_j, \a_{j+1}.$ Applying this repeatedly to all the highest weights of irreducible summands of $(c0 \ldots  0) \otimes \wedge^k(c0\ldots  0)$ for $k = 2,3,4$ we find that the irreducible summands of $\wedge^{k+1}(c0\ldots 0)$ have highest weights with support on $\a_1, \ldots , \a_{k+1}.$  

The weights of $ \l_1 = (10\ldots  0)$ are $\l_1, \l_1 -\a_1, \l_1-\a_1-\a_2, \ldots  $ which correspond to the basis vectors $v_1, v_2, v_3, \ldots  $ of $V_Y(\l_1)$. It follows that each weight of $\d$ has the form $(c0 \ldots  0) - (r_1\a_1 + \cdots  + r_{l+1}\a_{l+1})$ such that $r_1 \le c$ and $r_i \ge r_{i+1}$ for $i \ge 1$.  Then  weights of
$\wedge^k(c0 \ldots  0)$ have the form $((kc)0 \ldots  0) - (s_1\a_1 + \cdots + s_{l+1}\a_{l+1})$ with $s_1\le kc$ and $s_i \ge s_{i+1}$ for $i \ge 1$.

 Now consider the highest weight of an irreducible summand $\mu$ of $\wedge^k(c0 \ldots  0)$ which by the above has support on $\a_1, \ldots , \a_k$. Then  $s_j = (l+2-j)s_{l+1}$ for all $j \ge k+2$.  Applying this to $j = k+2$ we find that $((l+2)-(k+2))s_{l+1} \le kc$. If $s_{l+1} = 0$ then $\mu =((kc)0 \ldots  0) - (s_1\a_1 + \cdots + s_{k+1}\a_{k+1})$.  And if $s_{l+1} \ge 1$ then we find that $l \le 21$ for all values of $c$ and $k$.  Therefore we can apply a Magma computation to see that $\wedge^k(c0 \ldots  0)$ is MF in each case.  \hal

\begin{prop}\label{lambda3(omi)}  Let $\l = \l_j$ with $i=3$ or $4$, and let $\d = \o_i$ with $i\le 6$ or $4$, respectively. 
 Then  $V_Y(\l) \downarrow X$ is MF.
\end{prop}

\pf  We will show that these cases reduce to Magma checks as in the two preceding lemmas.
First note that $\l_1 \otimes \l_1 = V_Y(2\l_1) + V_Y(\l_2)$  and $\l_1 \otimes \l_2 = V_Y(\l_1 + \l_2) + V_Y(\l_3)$.  
Using the first tensor product and Theorem \ref{LR} we see that all irreducible summands of $V_Y(\l_2) \downarrow X$ have highest weights with support on $\a_1, \ldots , \a_{2i}.$  Then using the
second tensor product and Theorem \ref{LR} we find that the highest weights of irreducible summands of $V_Y(\l_3) \downarrow X$ have support on $\a_1, \ldots , \a_{3i}.$  Therefore, we can verify the result using a Magma
check.  Similarly reasoning applies for  $V_Y(\l_4) \downarrow X$ with $i \le 4.$ \hal

\begin{prop}\label{3l1,4l1} Let $V = V_Y(k\l_1)$ where $k = 3,4,5$.

{\rm (i)}  If $\d = c\om_1$, then $V \downarrow X$ is MF, provided $c \le 5,3,2$, respectively.

 {\rm (ii)} If $\d = \om_i$, then $V \downarrow X$ is MF, provided $i \le 5,3,3$, respectively.
\end{prop}

\pf  An argument with Theorem \ref{LR} shows that the highest weights of irreducible summands of $V_Y(k\l_1) \downarrow X$ have support on $\a_1, \ldots , \a_k$ or $\a_1, \ldots , \a_{ki}$, according to whether $\d = c\om_1$ or $\om_i$.
Consequently the assertion reduces to a Magma check as for the preceding results.  \hal

\begin{prop}\label{2l2,3l2}  If $\d = 2\o_1$ or $\o_2$, then both $V_Y(2\l_2) \downarrow X$ and 
$V_Y(3\l_2) \downarrow X$ are MF. Moreover, if $\d=\o_2$, then $V_Y(2\l_2) \downarrow X$ contains $(\o_1+\o_2+\o_5)\oplus (2\o_2+\o_4)$, and  $V_Y(3\l_2) \downarrow X$ contains
$(2\om_2 + \om_3 + \om_5) \oplus (\om_1 + 2\om_2 + \om_3 + \om_4)$.
\end{prop}

\pf  We first note that $V_Y(\l_2) \downarrow X = V_X(2\om_1 + \om_2)$ or $V_X(\om_1 + \om_3)$ according
to whether $\d = 2\o_1$ or $\o_2$.  Moreover, $V_Y(2\l_2)$ is contained in $V_Y(\l_2) \otimes V_Y(\l_2)$
and $V_Y(3\l_2)$ is contained in $V_Y(\l_2) \otimes V_Y(\l_2) \otimes V_Y(\l_2).$  Restricting to $X$ and applying 
Theorem \ref{LR}, we
see that irreducible summands of $V_Y(2\l_2) \downarrow X$ have support on $\a_1,\ldots ,\a_4$ or
$\a_1,\ldots ,\a_6$, according to whether $\d = 2\om_1$ or $\om_2.$  Similarly, irreducible summands of $V_Y(3\l_2) \downarrow X$ have support on $\a_1,\ldots , \a_6$ or $\a_1,\ldots , \a_9.$  Consequently
we can verify the result using Magma calculations assuming that $X$ has bounded rank as indicated.

Now $V_Y(2\l_2) \downarrow X$ is contained in $S^2(\wedge^2(\d)) = S^2(210\ldots  0)$ or $S^2(1010\ldots 0)$
according to whether $\d = 2\om_1$ or $\om_2.$  In each case we use a Magma computation to show the full symmetric square is MF.

Next consider $V_Y(3\l_2)$.  We have $S^3V_Y(\l_2) = V_Y(3\l_2) + V_Y(\l_2 + \l_4) + V_Y(\l_6)$ and
$V_Y(\l_2) \otimes V_Y(\l_4) =  V_Y(\l_2 + \l_4) + V_Y(\l_1 + \l_5) +  V_Y(\l_6).$  We conclude that
$V_Y(3\l_2) = S^3V_Y(\l_2) - (V_Y(\l_2) \otimes V_Y(\l_4)) + (V_Y(\l_1) \otimes V_Y(\l_5)) - V_Y(\l_6).$
So at this point we can restrict the above terms to $X$ and use Magma to verify that $V_Y(3\l_2) \downarrow X$ is MF.

Finally, the composition factors claimed in the statement follow along the way from the Magma computations. \hal

\section{ Table \ref{TAB4} configurations}\label{taby4}

In this subsection we will consider the configurations of Table \ref{TAB4} of Theorem \ref{MAINTHM}.  We begin with the special case where $X = A_3 < A_5.$
The statements of many of the lemmas to follow involve parameters $a$ and $d$, which are always taken to be positive integers. 

 \subsection{Embedding $X = A_3$, $\d = \om_2 $}\label{a3a5}

Here we consider the case $X = A_3  < Y = A_5,$ where the embedding is given by $\d = \om_2.$  Therefore we can regard $X$ as  $D_3.$ Note that the graph automorphism acts on the orthogonal module so, adjusting by a scalar, we see that there is an element of  $Y$ that induces a graph automorphism on $X$. Therefore when restricting representations from $Y$ to $X$ the restrictions are self dual.  

We temporarily change our usual notation to coincide with the orthogonal group notation.
That is we write  $\om_1, \om_2, \om_3$ for the fundamental dominant weights
of $D_3$, regarding $\om_1$ as the orthogonal representation and $\om_2, \om_3$ as
spin representations.  

With this in mind we will use some information from \cite{Howe2}  (see also \cite{fulthar}). Subject to the ordering
above this gives the fundamental dominant weights as $\om_1 = L_1, \om_2 = (L_1 + L_2 + L_3)/2,
\om_3 = (L_1 + L_2 - L_3)/2,$  where  $\pm L_i$ are the weights of the standard representation.
Suppose $a_1\om_1 + a_2\om_2 + a_3\om_3$ is a dominant weight, with $a_2 \ge a_3.$ We will only be considering representations of the orthogonal group $X$, which forces $a_2+a_3$ to be even. Then writing this in
terms of the $L_i$ we get $(a_1 + \frac{a_2+a_3}{2})L_1 + \frac{a_2+a_3}{2}L_2 +  \frac{a_2-a_3}{2}L_3$.  Thus we have a partition
$\e = a_1 + \frac{a_2+a_3}{2} \ge \frac{a_2+a_3}{2} \ge \frac{a_2-a_3}{2} \ge 0.$  
For future reference we note that if we write this partition as
$a+b+c \ge a+b \ge a$, then the corresponding dominant weight
is $ (c, 2a+b,b)$ for $D_3 $ or $(2a+b,c,b)$ for $A_3.$
On the other hand if $a_3 > a_2$
we get the partition  where the last term is replaced by $\frac{a_3-a_2}{2}.$  Dual pairs of irreducible representations correspond to the two options and correspond to the same partition.

Results of Littlewood \cite{Littlewood} (see Theorem 1.1 of \cite{Howe2} and  Equation (25.37) on p. 427 of \cite{fulthar}) provide a formula for restricting certain representations of $GL_6$ to $Y$.  Fix a partition
$\g = (\g_1 \ge \g_2 \ge
 \g_3 \ge 0).$ 
This result is as follows:
\begin{equation}\label{fhform}
V_{GL_6}(\g)\downarrow SO_6 = \sum_{\bar \xi} N_{\g,\bar \xi}V_{SO_6}(\bar \xi),
\end{equation}
where the sum is  over all partitions $\bar \xi = (\xi_1 \ge \xi_2 \ge \xi_3 \ge 0)$ and
\begin{equation}\label{moreform}
N_{\g,\bar \xi}  = \sum_{\e} c_{\e,\bar \xi}^{\g},
\end{equation}
the sum  over partitions $\e = \e_1 \ge \e_2 \ge \e_3 \ge 0$ with even parts, and the
terms $c_{\e,\bar \xi}^{\g}$ in the sum are Littlewood-Richardson coefficients.  We can rewrite this as
\begin{equation}\label{genform}
V_{GL_6}(\g)\downarrow SO_6 = \sum_{\e}(\sum_{\bar \xi} N_{\g,\bar \xi}V_{SO_6}(\bar \xi)).
\end{equation}

The strategy is as follows.  Fix $\g$.  Given an even partition $\e$ we determine those 
partitions $\bar \xi$ such that $c_{\e,\bar \xi}^{\g} \ne 0$  and show that this coefficient is  $1$.  There are usually very few such partitions.    Conversely, we argue that a given $\bar \xi$ can only arise from a single even partition $\e$.

We will have occasion to go back and forth between the $D_3$ and $A_3$ notation.  When confusion is possible we will identify how $X$ is being viewed.
At various points in the subsection we will use the  abbreviation $(abc)^+ = (abc) + (cba)$ to write the sum of a pair
of dual representations of $X = A_3$.

\begin{lem}\label{00d00}  The restrictions of the $A_5$-modules $d0000, 0d000,$ and  $00d00$ to  $X$ are MF.
\end{lem}

\pf  For $d0000$ this follows from Theorem \ref{SOpowers}. 
We will provide details for the case of $00d00$ and just indicate the changes required for the other case $0d000$, which is  easier.  The partition corresponding to $00d00$ is  $\g = (d,d,d,0,0,0)$ with weight $1^d2^d3^d$.  Fix a partition $\e = \e_1 \ge \e_2 \ge \e_3 \ge 0$ with even terms and  consider a skew tableaux of shape $\g/\e$.  This tableaux has  three rows of length $d$.  The first row  has $\e_1$ blank entries and the remainder the entries must be 1's by the Y-condition.  Suppose there are $a$ 1's so that $\e_1 + a = d$.  

The second row begins with  $\e_2$ blank entries.  The $Y$-condition implies that only 1's and 2's can appear and there cannot exist more than $a$ 2's.  But as column entries increase there must
exist at least $a$ 2's.  Therefore there must exist $b$ 1's, where $\e_2 + b + a = d.$

Finally we consider the third row which starts with $\e_3$ blank entries.  As column entries increase there must exist at least $a$ 3's and the $Y$ condition implies that there must exist exactly $a$ 3's. 
Another application of column increasing implies that there must exist $b$ 2's below the 1's in the
second row.  And the $Y$ condition implies that there cannot exist additional 2's in this row.  So the remaining entries are $c$ 1's where $\e_3 + a + b + c = d.$
We note that $a \equiv d \hbox{ mod }2$ and $b,c \equiv 0  \hbox{ mod }2$. 

What we have shown is that the labelling of the skew tableaux is $1^{a+b+c}2^{a+b}3^a.$  Therefore,
$\e$ determines a unique partition $\bar \xi = (a+b+c) \ge a+b \ge a$.  The corresponding highest weight in the  $D_3$
ordering is  $(c(2a+b)b)).$  Correspondingly the weight is  $((2a+b)cb)$ in the $A_3$ ordering.

Conversely, the above argument shows that given a partition $\bar \xi$ with weight $1^r2^s3^t$ the only possible even partition is $\e = (d-t, d-s,d-r).$  We have shown that all the coefficients in (\ref{genform}) are at most 1 and  the result follows in this case.

We illustrate the above with the case $d = 4$.  Here the possible even partitions $\e$ are 
$$(4,4,4), (4,4,2), (4,4,0),(4,2,2),(4,2,0),(4,0,0),(2,2,2),(2,2,0),(2,0,0),(0,0,0)$$
and these yield the respective $A_3$ (rather than $D_3$)  summands
$$(000),(020),(040),(202),(222),(404),(400)^+,(420)^+,(602)^+,(800)^+.$$

Now consider $0d000$, where the relevant partition has weight  $1^d2^d$.  The corresponding even partition $\e$ has the form $\e = e_1 \ge \e_2$.
Using the above techniques we see  the labelling of the skew tableau is $1^{a+b}2^a$.
Therefore, the labelling is determined by the partition and conversely.  The assertion follows. \hal

\begin{lem}\label{01d00}  The restriction $01d00 \downarrow X$ is MF.
\end{lem}

\pf   The proof here is similar to the proof of Lemma \ref{00d00}.  The partition corresponding to $01d00$ is  $\g = (d+1,d+1,d,0,0,0)$ with weight $1^{d+1}2^{d+1}3^d$.  Fix a partition $\e = e_1 \ge e_2 \ge e_3 \ge 0$ with even terms and  consider a skew tableaux of shape $\g/\e$.  The tableaux has  two rows of length $d+1$ and one row of length $d$.  

The analysis of the first two rows is just as in Lemma \ref{00d00}. The first row  has $\e_1$ blank entries followed by $a$ 1's by the Y-condition and the second row has $\e_2$ blank entries followed by $b$ 1's and then $a$ 2's.
Therefore $d+1 = \e_1 + a = \e_2 + b + a.$ 

First assume $a > 0$  and consider the third row which has length $d$ and which begins with $\e_3$ blank entries.  As columns are strictly increasing there must exist $a-1$ 3's at the end of the row under the $a-1$ 2's in columns
$\e_1 +1, \ldots , \e_1 + a-1$ of row 1. By the Y-condition there are at most $a$ 3's in the third row. Therefore  under the $b$ 1's in the second row there are either $b$ 2's or $b-1$ 2's followed by a single $3$.  Finally the row contains $c$ 1's so that $d = \e_3 +c + b + (a-1)$ or  $d = \e_3 +c + (b-1) + a$  corresponding to partitions with evaluations $1^{a+b+c}2^{a+b}3^{a-1}$ or $1^{a+b+c}2^{a+b-1}3^a,$ respectively.  Note that the multiplicity of $3$ has the same parity as $d$ or $d+1$ respectively.

So there are two partitions associated with  $\e$.  In one case the labelling contains $3^a$  and in the other case $3^{a-1}$.  Now let's reverse things and start with a partition $\bar \xi$ with evaluation $1^x2^y3^z.$  Then the parity of
$d$ and $z$ determine which of the configurations at the end of the last paragraph occur, and this in turn determines $\e.$  The result follows in this case.

Finally assume $a = 0$.  Here the first row is blank and the second row has
$\e_2$ blank entries followed by $b$ 1's.  Now consider the third row which begins with $\e_3$ blank entries.  Under the $b$ 1's in the second row there are  $b-1$ 2's.  There remain $c =  d- \e_3 - (b-1)$ entries and the possible labellings   are either $1^c$ or $1^{c-1}2^1.$  The  corresponding weight of $\bar \xi$ is   $(1^{b+c},2^{b-1})$  or $(1^{b+c-1},2^b)$, respectively.  Again the result follows.  \hal

\noindent{\bf Example.}
As an example of the last result we consider the representation $01300.$  Restricting to $O_6$ we have even partitions
$$442, 440, 422, 420, 400, 222, 220, 200, 000$$
 which yield the  summands  of $A_3$ as follows
$$010, 030, 111, 212, 313, 131,(301)^+, (321)^+,  (503)^+, (511)^+ , (701)^+, (410)^+.$$

\begin{lem}\label{0d100}  The restriction $0d100 \downarrow X$ is MF.
\end{lem}

\pf  This case is similar to, but easier than Lemma \ref{01d00}. The partition corresponding to $0d100$ is  $\g = (d+1,d+1,1,0,0,0)$ with weight $1^{d+1}2^{d+1}3^1$.  Therefore we only  consider  partitions $\e = \e_1 \ge \e_2  \ge 0$ with even terms.  Fix such a partition  and  consider a skew tableaux of shape $\g/\e$.  The tableaux has  two rows of length $d+1$ and one row of length $1$.  

The first row has $\e_1$ blank entries followed by $a$ 1's and the second row has $\e_2$ blank entries
followed by $b$ 1's and then $a$ 2's.  The third row has a single $1,2,$ or $3$ subject to the following
condition.  If $\e_2 = 0$, then the entry of the third row must be $2$ or $3$ as columns are strictly increasing, while if $\e_2 > 0$, the entry might be $1,2$ or $3$. So the labelling of the skew tableaux is $1^{a+b}2^a3^1,$ 
$1^{a+b}2^{a+1}$ (if $b > 0$), or $1^{a+b+1}2^a.$

Now reverse the situation and start with a partition with evaluation $1^x2^y3^z.$  If $z = 1$, then $\e_1$ satisfies
$\e_1 + y = d+1$, so $\e_1$ is determined uniquely.  Similarly $\e_2$ satisfies $\e_2 + x = d+1.$  Now suppose
$z = 0$, so that the partition has evaluation $1^x2^y$.  This time $\e$ is determined by the parity of $y$.
If $y$ and $d+1$ have the same parity, then $\e_1 + y = d+1$ and $\e_2 + x-1 = d+1.$  Otherwise,
$\e_1 + y-1 = d+1$ and $\e_2 + x = d+1.$  In any case $\e$ is determined by the partition and we conclude
that the restriction is MF.  \hal

\begin{lem}\label{1d000}  The restriction $1d000 \downarrow X$ is MF.
\end{lem}

\pf  The partition corresponding to $1d000$ is $\g = (d+1,d,0^4)$ with weight $1^{d+1}2^d.$ So let $\e = \e_1 \ge \e_2 \ge 0$ be an even partition and consider $\g/\e$.  Then $d+1 = \e_1 +a$ and $d = \e_2 + b + (a-1)$ or $\e_2 + (b-1) + a$   corresponding to labellings  $1^{a+b}2^{a-1}$ or $1^{a+b-1}2^a,$ the two cases differing by
the parity of the multiplicity of 2.  Reversing the situation we chose a partition with evaluation $1^x2^y$.  If $y$ and $d+1$ have the same parity, then  $\e_1 = d+1 - y,$  $\e_2 = d - x,$ and the evaluation determines the
even partition.  While if
$y$ and $d+1$ have opposite parity, then  $\e_1 = d+1 - (y+1)$
and $\e_2 = d - (x-1)$ and the partition is again determined.   The result follows.  \hal

\begin{lem}\label{d0100}  The restriction of  $d0100$  to  $X$ is MF. 

\end{lem}

\pf  Here the partition is $\g = (d+1,1,1,0^3)$ with weight $1^{d+1}2^13^1.$  So here even partitions have
the form $\e = \e_1 \ge 0 \ge 0.$  If $\e_1 = 0$, then the only possible labelling is $1^{d+1}2^13^1$  and for $\e_1 > 0$ there are two possible labellings, namely $1^a2^13^1$ and $1^{a+1}2^1$, where
in each case $\e_1 + a = d+1.$  
Conversely given a partition $\bar \xi$ with weight $1^x2^y3^z$  we see that there is at most 1 even partition $\e$ corresponding to it. Necessarily $y = 1$.  If $z = 0$, then the partition is  $\e_1 = d+2 - x.$  And if $z = 1$, 
then $\e_1 = d+1 -x.$  The result follows.  \hal

\begin{lem}\label{10d00}  The restriction of  $10d00$  to  $X$ is MF. 
\end{lem}

\pf  This time the partition is $\g = (d+1,d,d,0^3)$ with weight $1^{d+1}2^d3^d.$  Let $\e = \e_1 \ge \e_2 \ge \e_3$ be an even partition and consider $\g/\e$.  In row 1 there are $\e_1$ blank cells followed by $a$ 1's.
Row two has length $d$.  There are $a-1$ 2's below the first $a-1$ 1's in row $1$. These must be  preceded by either   $b$ 1's or $b-1$ 1's and a 2.
Then $d+1 = \e_1 +a$ and $d = \e_2 + b + (a-1)$ or $\e_2 + (b-1) + a$   so that the labellings  of just these two rows is given by$1^{a+b}2^{a-1}$ or $1^{a+b-1}2^a.$   

Now consider row 3.  In the first case  the Y-condition implies that this row has $\e_3$ blank cells followed by $c$ 1's, then $b$ 2's, then $a-1$ 3's, so that  $d = \e_3 + c + b + (a-1).$  This yields
the labelling $1^{a+b+c}2^{a+b-1}3^{a-1}.$  In the second case  the third row has $\e_3$ blank cells followed by $c$ 1's, then $b-1$ 2's, then $a$ 3's, so that  $d = \e_3 + c + (b-1) + a.$  This yields
the labelling $1^{a+b+c-1}2^{a+b-1}3^a.$   The parity of the multiplicity of $3$ distinguishes between the two
possible partitions. 

Conversely, given a partition $\bar \xi = 1^x2^y3^z$, suppose it arises from an even partition $\e$. Then it corresponds to the first or second case above depending on whether or not  $z \equiv d  \hbox{ mod }2$.  And we see that
$\bar \xi$ determines $\e$.   The result follows.  \hal

\begin{lem}\label{d1100}  The restrictions of  $d1100$ and $1d100$ to  $X$ are not MF for $d \ge 2.$
The restriction of $11100$ to $X$ is MF.  More precisely the restriction is $210 +012+ 311+113+101+202+020+121.$
\end{lem}

\pf  First consider $d1100$, where the corresponding partition  $\g = (d+2,2,1,0^3)$ with weight $1^{d+2}2^23^1.$  Let $\e = \e_1 \ge \e_2 \ge 0$ be an even partition.  Note that $\e_2 = 0$ or $2$.  First assume $d$ is even and take $ \e_1 = d+2$  and $\e_2 = 0.$  Then we get the labelling $1^22.$   But now take $\e_1 = d$ and $\e_2 = 2$.  Here too we get
a labelling $1^22$.  Therefore the labelling $1^22$ arises from two distinct even partitions  and we conclude that  the restriction is not MF.  If $d\ge 2$ is odd we have a similar argument by considering the partitions $\e_1 = d+1 > \e_2 = 0$ and  $\e_1 = d-1 > \e_2 = 2.$  In both cases there is a labelling $1^32.$  The result follows for this case.  We leave it to the reader to work out $11100 \downarrow X.$

Now consider the case $1d100$ where the corresponding partition is $\g = (d+2, d+1,1,0^3)$ with weight 
$1^{d+2}2^{d+1}3^1.$  First take $d$ even.  Then the partitions $d+2 \ge 0$ and $d \ge 2 \ge 0$ can yield
a labelling $1^{d+1}2^1.$  And if $d$ is odd, the partitions $d+1 \ge 0$ and $d-1 \ge 2 \ge 0$ both provide a
labelling $1^{d+1}2^2.$  Again we obtain the assertion.  \hal

The following lemma covers certain additional cases when $X = A_3$ and $\d = \om_2$, including
some cases where the above restriction techniques do not apply.

\begin{lem}\label{a3om2}  The restrictions of the $A_5$-modules $a0001$, $a0010$, $a0100$ and $a1000$ to $X$ are all MF.  More precisely
\begin{itemize}
\item[{\rm (i)}] $(a0001) \downarrow X = (0(a+1)0) + (1(a-1)1) + (0(a-1)0) +(1(a-3)1) +(0(a-3)0) +   \cdots.$
\item[{\rm (ii)}] $(a0010)\downarrow X = ((1a1) + (1(a-2)1) + \cdots) + ((2(a-1)0)^+ + (2(a-3)0)^+ \cdots)$.
\item[{\rm (iii)}] $(a0100)\downarrow X = ((1(a-1)1) + (1(a-3)1) + \cdots) + (2a0)^+ + (2(a-2)0)^+ + \cdots).$
\item[{\rm (iv)}] $(a1000)\downarrow X= ((1a1) + (1(a-2)1) + \cdots) + ((0a0)+(0(a-2)0) + \cdots)$, 
except that $(000)$ does not occur in the latter sum if $a$ is even.
\end{itemize}
\end{lem}

\pf  We begin with some results on tensor products for $A_3$ which follow from Littlewood-Richardson arguments (see Theorem \ref{LR}):
\begin{equation}\label{arr1}
\begin{array}{l}
(0x0) \otimes (010) = (0(x+1)0) + (1(x-1)1) + (0(x-1)0) \\
(0y0) \otimes(101) = (1y1) +  (2(y-1)0)^+  + (0y0) +(1(y-2)1) \\
(0z0) \otimes ((200) + (002)) = (2z0)^+  + (1(z-1)1)^2+(2(z-2)0)^+,
\end{array}
\end{equation}
although certain terms do not occur for small values of $x,y,z.$  

(i) First note that $(a0001) = (a\l_1 \otimes \l_5) - (a-1)\l_1.$  Viewing
$\om_2$ as $(010)$,  Theorem \ref{SOpowers} shows that   $S^a(010) = (0a0) + (0(a-2)0) + \cdots.$
Therefore, applying the first of the above tensor products we see that $S^a(010) \otimes (010) = 
(0(a+1)0) + (1(a-1)1) + (0(a-1)0) + (0(a-1)0) + (1(a-3)1) +(0(a-3)0) +   \cdots,$   so that
$(0(a-1)0), (0(a-3)0), \ldots $ each appear with multiplicity 2.  But subtracting
$S^{a-1}(010)$ we obtain 
$(a0001) \downarrow X = (0(a+1)0) + (1(a-1)1) + (0(a-1)0) +  (1(a-3)1) +(0(a-3)0) +   \cdots.$ 

(ii) We first note that $(a0010) =  (a\l_1 \otimes \l_4)) - ((a-1)\l_1 + \l_5))$. 
Restricting to $X$ we have $ a\l_1 \downarrow X = S^a(010) = (0a0) + (0(a-2)0) + \cdots$ and $\l_4 \downarrow X = \wedge^4(010) = (101).$
Consequently we use the second tensor product formula in (\ref{arr1}) 
to see that $(a\l_1 \otimes \l_4)) \downarrow X $ is the sum of the following terms:
$$(0a0) \otimes (101) = (1a1) +  (2(a-1)0)^+  + (0a0) +(1(a-2)1),$$
$$(0(a-2)0) \otimes (101) = (1(a-2)1) +  (2(a-3)0)^+  + (0(a-2)0) +(1(a-4)1),$$
$$(0(a-4)0) \otimes (101) = (1(a-4)1) +  (2(a-5)0)^+  + (0(a-4)0) +(1(a-6)1),$$
$$\vdots$$
Now by (i), $((a-1)0001) \downarrow X = (0a0) + (1(a-2)1) + (0(a-2)0) +(1(a-4)1) +(0(a-4)0) +   \cdots.$
The assertion follows.

(iii)  We have $(a0100) =  (a\l_1 \otimes \l_3)) - ((a-1)\l_1 + \l_4))$.  Restricting the first tensor product to $X$
we obtain
$$((0a0) + (0(a-2)0) + \cdots) \otimes ((200) + (002))$$
 To evaluate this we use the third tensor product formula in (\ref{arr1}) to obtain 
 $$(0a0) \otimes  ((200) + (002)) = (1(a-1)1)^2 + (2a0)^+ + (2(a-2)0)^+$$ 
$$(0(a-2)0) \otimes  ((200) + (002)) = (1(a-3)1)^2 + (2(a-2)0)^+ + (2(a-4)0)^+$$
$$\vdots$$
By (ii), we have
$$((a-1)0010)\downarrow X = ((1(a-1)1) + (1(a-3)1) + \cdots) + ((2(a-2)0)^+ + (2(a-4)0)^+ \cdots).$$
Subtracting this from the above we obtain the result.

(iv)  We have $(a1000) =  (a\l_1 \otimes \l_2)) - ((a-1)\l_1+ \l_3))$. 
Restricting the first tensor product to $X$ gives
 $$((0a0) + (0(a-2)0) + \cdots) \otimes (101).$$
 Using the second tensor product formula in (\ref{arr1}) we get the sum of the terms
 $$(1a1) + (1(a-2)1) + (0a0) + (2(a-1)0)^+$$
 $$(1(a-2)1) + (1(a-4)1) + (0(a-2)0) + (2(a-3)0)^+$$
$$\vdots$$
By (iii), we have
$$((a-1)0100)\downarrow X = ((1(a-2)1) + (1(a-4)1) + \cdots) + (2(a-1)0)^+ + (2(a-3)0)^+ + \cdots).$$
The result follows.  \hal

\begin{lem}\label{11001} The restriction $11001\downarrow X$ is MF. 
\end{lem}

\pf  This is a straightforward Magma computation.  \hal

\begin{lem}\label{0a001}  For all $a$, the restriction $0a001\downarrow X$ is MF. 
\end{lem}

\pf  We first observe that $0a000 \otimes 00001 = 0a001 + 1(a-1)000$  so we can use the analysis
in the proofs of Lemmas \ref{00d00} and \ref{1d000} to assist with the proof here.
We begin with some  results on tensor products which follow from Littlewood-Richardson arguments (see Theorem \ref{LR}):
\[
\begin{array}{ll}
(0x0) \otimes (010) = & (0(x+1)0) + (1(x-1)1) + (0(x-1)0) \\
(yxy) \otimes (010) = & (y(x+1)y) + ((y+1)(x-1)(y+1)) + (y(x-1)y) + \\
                      & ((y-1)(x+1)(y-1)) + ((y+1)x(y-1)) + ((y-1)x(y+1)),
\end{array}
\]
where certain terms are deleted if $x=1$ or $y=1.$ The weight of the partition associated with $0a000$ is $1^a2^a$.  We will first assume $a$ is even and later indicate the changes required for the odd case. 

We will list  possible
even partitions $\e = \e_1 \ge \e_2 \ge 0$ by simply writing $(\e_1,\e_2)$.  As $a$ is even, these pairs
are as follows: $(a,a),(a,a-2), \ldots , (a,0),(a-2,a-2), (a-2,a-4)\ldots , (a-2,0), \ldots , (0,0)$
and each gives rise to a unique labelling; they correspond to $A_3$ representations 
$$(000), (020), \ldots , (0a0),$$
$$ (202), (222), \ldots , (2(a-2)2),$$
$$(404), (424), \ldots , (4(a-4)4)$$
$$\vdots$$
$$(a0a).$$

Next we tensor each of the above irreducibles with $(010)$ using the above equations. We obtain a
multiplicity in  precisely two ways.  Namely, the tensor product of $(tst)$ and $(t(s+2)t)$ with $(010)$ each
contain a copy of $(t(s+1)t),$ while the tensor product of $(tst)$ and $((t+2)(s-2)(t+2))$ with $(010)$ each
contain a copy of $((t+1)(s-1)(t+1))$. Consequently we get  multiplicity 2 for the following terms:
$(010), (030), \ldots , (111), (131), \ldots , (212), (232), \ldots .$ 

Now consider the restriction of $1(a-1)000$ to $X$, still assuming that $a$ is even. The weight of the partition associated with $1(a-1)000$ is $1^a2^{a-1}$.  Using the notation for even partitions as above we have
partitions $(a,a-2), (a,a-4) \ldots , (a,0),(a-2,a-2), (a-2,a-4)\ldots , (a-2,0), \ldots , (0,0).$
Fix an even partition $\e = \e_1 \ge \e_2 \ge 0$  and write $a = \e_1 + t+1.$ Then we can write
$a-1 = \e_2 + v + t$ or $a-1 = \e_2 +(v-1) + (t+1)$ corresponding to the two possible labellings
$1^{v+t+1}2^t$ and $1^{t+v}2^{t+1}$ both of which are consistent with the Y-condition, etc.  These
yield representations $(t,v+1,t)$ and $((t+1)(v-1)(t+1)),$ respectively.  It follows that $0a001 \downarrow X$
is indeed MF.

Now suppose $a$ is odd.  The situation is very similar to the above. The possible partitions (in the above notation) are $(a-1,a-1), (a-1,a-3) \ldots , (a-1,0),(a-3,a-3), (a-3,a-5)\ldots , (a-3,0), \ldots , (0,0)$ and these
correspond to the representations
$$(101), (121), \ldots , (1(a-1)1),$$
$$ (303), (322), \ldots , (3(a-3)3),$$
$$(505), (525),\ldots , (5(a-5)5)$$
$$\vdots$$
$$(a0a).$$
We get multiplicities just as indicated in the third paragraph and we argue as above that each irreducible
appearing with multiplicity 2 appears within the restriction of $1(a-1)000$ to $X$.  \hal

\begin{lem}\label{abc00,c>1} If $c > 1$ and $ab \ne 0$, then  $abc00\downarrow X$ is not MF.
\end{lem}
\pf  Assume $c > 1.$  We will work through the various possibilities for $a,b,c$.  For each case other than $11200$
we indicate  an even partition and  two different labellings with the same weight.   This implies that the
case is not MF.  For the exceptional case we find that there are two distinct even partitions giving the same
labelling. To indicate  that a row begins with $x$ blank cells we write $(-)^x$.   
 
\vspace{2mm}

$ a=b=1, c$ odd. $\e = (c+1, c-1, c-3):$

\vspace{2mm}

$(-)^{c+1}1$; $(-)^{c-1}11$; $(-)^{c-3}122$

 $(-)^{c+1}1$; $(-)^{c-1}12$; $(-)^{c-3}112$
 
 \vspace{2mm}
 
 $a=b=1,c > 2$ even.  $\e = (c, c-2, c-4):$
 
  \vspace{2mm}

 $(-)^c11$; $(-)^{c-2}112$; $(-)^{c-4}1223$

 $(-)^c11$; $(-)^{c-2}122$; $(-)^{c-4}1123$
 
\vspace{2mm}
$a=b=1,c = 2$

Below we see that each of the  partitions  $(4,0,0)$ and $(2,2,0)$  yields a labelling with weight $1^32^2.$

\vspace{2mm}

$(-)^4; 111; 22$

$(-)^211; (-)^22; 12.$

\vspace{2mm}

From now on we assume that either $b > 1$ or $a > 1$.  In listing the third row of the partitions below
we will write $c^*rs$.  Here $c^*rs = (-)^{c-2}rs$ or  $(-)^{c-3}1rs$ according to whether $c$ is even or $c$ is odd.
Also in describing the third row of a partition we  write $c-y$ to mean $c-2$ or $c-3$ according to whether $c$ is even or odd.
\vspace{2mm}

$b \ge 2$, $a+b+c$ odd, $b+c$ even.  $\e = (a+b+c-1, b+c-2, c-y):$

\vspace{2mm}

$(-)^{a+b+c-1}1; (-)^{b+c-2}11; c^*12$

$(-)^{a+b+c-1}1; (-)^{b+c-2}12; c^*11$

\vspace{2mm}

$b \ge 2$, $a+b+c$ odd, $b+c$ odd.  $\e = (a+b+c-1, b+c-1, c-y):$

\vspace{2mm}

$(-)^{a+b+c-1}1; (-)^{b+c-1}1; c^*12$

$(-)^{a+b+c-1}1; (-)^{b+c-1}2; c^*11$

\vspace{2mm}
$b \ge 2$, $a+b+c$ even, $b+c$ even $a > 1.$  $\e = (a+b+c-2, b+c-2, c-y):$

$(-)^{a+b+c-2}11; (-)^{b+c-2}11; c^*12$

$(-)^{a+b+c-2}11; (-)^{b+c-2}12; c^*11.$

\vspace{2mm}
$b \ge 2$, $a+b+c$ even, $b+c$ odd, $a > 1.$ $\e = (a+b+c-2, b+c-1, c-y):$

$(-)^{a+b+c-2}11; (-)^{b+c-1}1; c^*12$

$(-)^{a+b+c-2}11; (-)^{b+c-1}2; c^*11.$

\vspace{2mm}
$b \ge 2$,  $b+c$ odd, $a = 1.$  $\e = (a+b+c-2, b+c-3, c-y):$

$(-)^{a+b+c-2}11; (-)^{b+c-3}122; c^*12$

$(-)^{a+b+c-2}11; (-)^{b+c-3}112; c^*22.$

\vspace{2mm}
$b \ge 2$,  $b+c$ even, $a = 1.$  $\e = (a+b+c-1, b+c-2, c-y):$

$(-)^{a+b+c-1}1; (-)^{b+c-2}12; c^*11$

$(-)^{a+b+c-1}1; (-)^{b+c-2}11;c^*12.$

\vspace{2mm}

It remains to consider those cases where $b = 1$ and $a > 1$.

\vspace{2mm}
$b = 1, a>1, a+b+c$ even, $b+c$ odd. $\e = (a+b+c-2, b+c-1, c-y):$

$(-)^{a+b+c-2}11; (-)^{b+c-1}1; c^*12$

$(-)^{a+b+c-2}11; (-)^{b+c-1}2; c^*11$

\vspace{2mm}

$b = 1, a>1, a+b+c$ even, $b+c$ even.  $\e = (a+b+c-2, b+c-2, c-y):$

$(-)^{a+b+c-2}11; (-)^{b+c-2}12; c^*12$

$(-)^{a+b+c-2}11; (-)^{b+c-2}11; c^*22.$

\vspace{2mm}

$b = 1, a>1, a+b+c$ odd, $b+c$ odd.  $\e = (a+b+c-1, b+c-1, c-y):$

$(-)^{a+b+c-1}1; (-)^{b+c-1}1; c^*12$

$(-)^{a+b+c-2}11; (-)^{b+c-1}2; c^*11.$

\vspace{2mm}

$b = 1, a>1, a+b+c$ odd, $b+c$ even.  $\e = (a+b+c-1, b+c-2, c-y):$

$(-)^{a+b+c-1}1; (-)^{b+c-2}11; c^*12$

$(-)^{a+b+c-2}11; (-)^{b+c-2}12; (-)^{c^*}11.$ \hal

\vspace{2mm}

The final case of this subsection is by far the most complicated.

\begin{lem}\label{0a010}  Let $X= A_3$ embedded in $A_5$ via $\d = 010.$ Then $0a010\downarrow X$ is MF. 
\end{lem}

\pf The proof is more  complicated than the previous results.  
To illustrate some of the 
ideas we will provide explicit  details for the special case $a = 5$  as we proceed 
through the proof, first for the case $a$ odd.

We use the decomposition $(0(a+1)000)\otimes (01000) = (0(a+2)000) + (1a1000) + (0a010)$.

As in the proof of Lemma \ref{0a001}, we have the following $A_3$-irreducible sumands of $(0(a+1)000)$:
 $$(000), (020), \ldots , (0(a+1)0),$$
$$ (202), (222), \ldots , (2(a-1)2),$$
$$(404), (424), \ldots , (4(a-3)4)$$
$$\vdots$$
$$((a+1)0(a+1)).$$
Note that all the entries in the above highest weights are even.  Also, it follows from the
 construction that 
if $(yxy)$ occurs, then $x+y \le a+1.$ We have a similar decomposition of $(0(a+2)000)\downarrow X$:
 $$(101), (121), \ldots , (1(a+1)1),$$
$$ (303), (323), \ldots , (3(a-1)3),$$
$$(505), (525), \ldots , (5(a-3)5)$$
$$\vdots$$
$$((a+2)0(a+2)).$$

In order to restrict the tensor product $(0(a+1)000)\otimes (01000)$ to $X$ we begin by recording the 
result of certain tensor products where
one of the factors is $(101)$.  The results are all checked using Theorem \ref{LR}.

If $x >0$ we have
$$(0x0) \otimes (101) = (1x1) + (1(x-2)1) +(0x0) + (2(x-1)0)^+.$$

If $y > 0$ we have
\[
\begin{array}{ll}
(y0y) \otimes (101) = & ((y+1)0(y+1)) + (y0y)^2 + ((y-1)0(y-1)) + \\
                      &  ((y-1)2(y-1)) + ((y+1)1(y-1)^+  + (y1(y-2))^+. 
\end{array}
\]

If  $x,y>0$ we obtain 
\[
\begin{array}{ll}
(yxy) \otimes (101) = &((y+1)x(y+1)) + ((y+1)(x-2)(y+1))+ \\
                      & (yxy)^3 + ((y-1)x(y-1)) + ((y-1)(x+2)(y-1)) + \\
                      & ((y+2)(x-1)y)^+  +  ((y+1)(x-1)(y-1))^+ + \\
                      & (y(x+1)(y-2))^+ + ((y+1)(x+1)(y-1))^+,
\end{array}
\]
although some terms are missing if $x$ or $y$ equals 2. 

\vspace{2mm}

In the following we illustrate the result of tensoring $(0(a+1)000)\downarrow X$ with $(101)$
for the special case $a = 5$.  The patterns are already clear for this case.  We list the 
irreducible followed by its 
tensor product with $(101)$.  

\vspace{2mm}

$(000): (101)$

$(020): (121), (101), (020), (210)^+$

$(040): (141), (121), (040), (230)^+$

$(060): (161), (141), (060), (250)^+$

$(202): (303), (202)^2, (121), (101), (311)^+, (210)^+$

$(222): (323), (222)^3, (141), (303), (121), (412)^+, (331)^+, (311)^+, (230)^+$

$(242): (343), (242)^3, (161), (323), (141), (432)^+, (351)^+, (331)^+, (250)^+$

$(404): (505), (404)^2, (323), (303), (513)^+, (412)^+$

$ (424): (525), (424)^3, (343), (323), (505), (614)^+, (533)^+, (513)^+, (432)^+$

$ (606):  (707), (606)^2, (525), (505), (715)^+, (614)^+.$

\vspace{2mm}

We study the repeated summands, beginning with the
irreducible summands that are not self-dual.  We claim that a fixed irreducible summand $(rst)$ with $r>t$ 
occurs with multiplicity at most $2$ in $(0(a+1)000)\downarrow X \otimes (101)$.
In the following we indicate all coincidences of summands $(rst)$ with $r > t$ that appear when tensoring 
irreducibles of the form $(yxy)$ with $(101).$  Here $x$ and $y$ are even.

If $x,y >0, $ then $(yxy) \otimes (101)$ and $(y(x+ 2)y) \otimes (101)$ (respectively 
$(y(x- 2)y) \otimes (101)$) both contain $((y+1)(x+1)(y-1))$ (respectively $(y+1)(x-1)(y-1))$) 
(note that the former requires  $x+y+2 \le a$). And $(yxy) \otimes (101)$ and 
$((y+2)(x-2)(y+2)) \otimes (101)$ both contain $((y+2)(x-1)y)$.  
There are additional cases when $xy = 0$.  Indeed  $(0x0) \otimes (101)$ and
$(2(x-2)2) \otimes (101)$ both contain $(2(x-1)0).$  Also $(y0y) \otimes (101)$ contains
both $((y+1)1(y-1))$ and $(y1(y-2))$.  These also appear in $(y2y) \otimes (101)$ 
and $((y-2)2(y-2)) \otimes (101)$, respectively.  It follows that the multiplicity of  
$(rst)$ for $r>t$ in $(0(a+1)000)\downarrow X \otimes (01000)$ is at most $2$.
We need to show that each such summand occurs in the module $(1a100)\downarrow X$ (as none occurs
in $(0(a+2)000)\downarrow X$).

We separate into cases as follows:\begin{enumerate}[1.]
\item $(2(x-1)0)$, where $2\leq x\leq a+1$ is even;
\item $((y+1)1(y-1))$, where $2\leq y\leq a+1$ is even;
\item $((y+1)(x-1)(y-1))$, where $2\leq y\leq a+1$ is even, $4\leq x\leq a+1$ is even and $x+y\leq a+1$;
\item $((y+2)(x-1)y)$ where $2\leq y\leq a+1$ is even, $2\leq x\leq a+1$ is even and $x+y\leq a+1$;
\item $((y+1)(x+1)(y-1))$, same conditions as in previous case;
\item $(y(x+1)(y-2))$, same conditions as in previous case.
\end{enumerate}

We now give the partitions in each case which show that these occur as summands of $(1a100)\downarrow X$. 
Note that the weight of the associated partition is $1^{a+2}2^{a+1}3^1$.
\begin{enumerate}[1.]
\item Here we take the even partition $a+1\geq a-x+1$, with labelled skew tableau \linebreak
$(-)^{a+1}1;(-)^{a-x-1}1^{x+1}2;3$, giving the weight $1^x2^13^1$;

\item  Here we take the even partition $a+1-y\geq a+1-y$, with labelled skew tableau \linebreak
$(-)^{a+1-y}1^{y+1};(-)^{a+1-y}2^y;3$, giving the weight $1^{y+1}2^y3^1$;

\item  Here we take the even partition $a+2-(y+1)\geq a+1-(x+y-2)$,  with labelled skew tableau \linebreak
$(-)^{a+2-(y+1)}1^{y+1};(-)^{a+1-(x+y-2)}1^{x-2}2^y;3$, giving the weight $1^{x+y-1}2^y3^1$;

\item  Here we take the even partition $a+2-(y+1)\geq a+1-(x+y)$, with labelled skew tableau \linebreak
$(-)^{a+2-(y+1)}1^{y+1};(-)^{a+1-(x+y)}1^{x-1}2^{y+1};3$, giving the weight $1^{x+y}2^{y+1}3^1$;

\item  Here we take the even partition $a+2-(y+1)\geq a+1-(x+y)$, with labelled skew tableau \linebreak
$(-)^{a+2-(y+1)}1^{y+1};(-)^{a+1-(x+y)}1^x2^y;3$, giving the weight $1^{x+y+1}2^y3^1$;

\item  Here we take the even partition $a+2-(y-1)\geq a+1-(x+y)$, with labelled skew tableau \linebreak
$(-)^{a+2-(y-1)}1^{y-1};(-)^{a+1-(x+y)}1^{x+1}2^{y-1};3$, giving the weight $1^{x+y}2^{y-1}3^1$.\end{enumerate}

Now look at the repeated summands of the form $(rsr)$ where $r$ is even. Note that none of these occur in $((0(a+2)000)\downarrow X$.

These are of the form $(y0y)$ for $2\leq y\leq a+1$ even,  occurring with multiplicity $2$, and
$(yxy)$ for $2\leq y\leq a+1$, even, and $2\leq x\leq a+1$ even with $x+y\leq a+1$, occurring with 
multiplicity $3$. So we need to produce one summand $(y0y)$ and two summands $(yxy)$ of 
$(1a100)\downarrow X$ for all relevant values of $x$ and $y$.

For $(y0y)$, take the partitions $\epsilon_1\geq\epsilon_2$ for 
$(\epsilon_1,\epsilon_2) = (a+1-2j,a+1-2j-2)$, for $0\leq j\leq \frac{a-1}{2}$. 
For $(yxy)^2$, we will need to take two different partitions, each of which gives rise to a 
labelled skew tableau whose weight is $1^{x+y}2^y$, namely the partitions $a+2-(y-1)\geq a+1-x-y$ 
and $a+2-(y+1)\geq a+1-(y+x-2)$.
For the first partition, we  label the tableaux with $(-)^{a+3-y}1^{y-1};(-)^{a+1-x-y}1^{x+1}2^{y-1};2$. 
For the second 
partition, we label via $(-)^{a+1-y}1^{y+1};(-)^{a+3-x-y}1^{x-2}2^y;1$.

Now look at the self-dual summands of the form $(rsr)$ where $r$ is odd. These occur as 
follows:
\begin{enumerate}[1.]
\item $(1s1)$ for $0\leq s\leq a+1$ even, with multiplicity $3$.
\item $(r0r)$, with $r\geq 1$ odd, with multiplicity $3$
\item $(rsr)$, with $r\geq 3$ odd and $s\geq 2$ even, $r+s\leq a+2$, with multiplicity $4$ 
(except a limit case which will be discussed below).
\end{enumerate}

All of these summands occur in $(0(a+2)000)\downarrow X$, so we need to show that $(1s1)$, 
$(r0r)$, and $(rsr)^2$ occur in $(1a100)\downarrow X$.
For $(1s1)$ take the partitions $(\epsilon_1,\epsilon_2)$ of the form $(a+1,2j)$ for 
$0\leq j\leq \frac{a+1}{2}$.
For $(r0r)$, take the partitions $(a+1-2j,a+1-2j)$ for $0\leq j\leq \frac{a-1}{2}$.
Finally for $(rsr)$, we take the partition $(a+2-r, a+1-(r+s-1))$, and indicate two 
different
labelled skew tableaux, each of which will have weight $1^{r+s}2^r$. But first we 
point out 
that if $r+s=a+2$ this summand occurs only twice in the tensor product 
$(0(a+1)000)\downarrow X\otimes (101)$ and so we do not need to consider this case.
The two labellings are: $(-)^{a+2-r}1^r;(-)^{a+2-r-s}1^s2^{r-1};2$ and 
$(-)^{a+2-r}1^r;(-)^{a+2-r-s}1^{s-1}2^{r};1$.

We now indicate how the analysis of the case $a$ even goes through. (There are no particular difficulties.) As above, we study the repeated summands in the tensor product 
$(0(a+1)000)\downarrow X\otimes (101)$, beginning with the
irreducible summands that are not self-dual. As before, these can occur with multiplicity at most $2$ and do not occur in the summand $(0(a+2)000)\downarrow X$. 
We need to show that each such summand occurs in the module $(1a100)\downarrow X$.

We separate into cases as follows:\begin{enumerate}[1.]
\item $((y+1)1(y-1))$, where $1\leq y\leq a+1$ is odd;
\item $(y1(y-2))$  where $1\leq y\leq a+1$ is odd;
\item $((y+2)(x-1)y)$ where $1\leq y\leq a+1$ is odd, $4\leq x\leq a$ is even and $x+y\leq a+1$;
\item $((y+1)(x-1)(y-1))$, same conditions as in previous case.
\item $((y+1)(x+1)(y-1))$, where $1\leq y\leq a+1$ is odd, $2\leq x\leq a$ is even and $x+y\leq a+1$.
\end{enumerate}

We now give the partitions in each case which show that these occur as summands of 
$(1a100)\downarrow X$. 
\begin{enumerate}[1.]
\item Here we take the even partition $a+2-(y+1)\geq a-y+1$, with labelled skew tableau \linebreak
$(-)^{a+2-(y+1)}1^{y+1};(-)^{a+1-y}2^y;3$, giving the weight $1^{y+1}2^y3^1$;

\item  Here we take the even partition $a+2-(y-1)\geq a+1-y$, with labelled skew tableau \linebreak
$(-)^{a+2-(y-1)}1^{y+1};(-)^{a+1-y}12^{y-1};3$, giving the weight $1^{y}2^{y-1}3^1$;

\item  Here we take the even partition $a+2-(y+1)\geq a+1-(x+y)$, with labelled skew tableau 
$(-)^{a+2-(y+1)}1^{y+1};(-)^{a+1-(x+y)}1^{x-1}2^{y+1};3$, giving the weight $1^{x+y}2^{y+1}3^1$;

\item Here we take the even partition $a+2-(y+1)\geq a+1-(x+y-2)$, with labelled skew tableau 
$(-)^{a+2-(y+1)}1^{y+1};(-)^{a+1-(x+y-2)}1^{x-2}2^{y};3$, giving the weight $1^{x+y-1}2^{y}3^1$;

\item  Here we take the even partition $a+2-(y+1)\geq a+1-(x+y)$, with labelled skew tableau \linebreak
$(-)^{a+2-(y+1)}1^{y+1};(-)^{a+1-(x+y)}1^{x}2^{y};3$, giving the weight $1^{x+y+1}2^{y}3^1$;
\end{enumerate}

Now look at the repeated summands of the form $(rsr)$ where $r$ is odd. Note that none of these occur in $((0(a+2)000)\downarrow X$.

These are of the form $(y0y)$ for $1\leq y\leq a+1$ odd,  occurring with multiplicity $2$, and
$(yxy)$ for $1\leq y\leq a+1$, odd, and $2\leq x\leq a+1$ even with $x+y\leq a+1$, occurring with 
multiplicity $3$. So we need to produce one summand $(y0y)$ and two summands $(yxy)$ of 
$(1a100)\downarrow X$ for all relevant values of $x$ and $y$.

For $(y0y)$, take the partitions $\epsilon_1\geq\epsilon_2$ for 
$(\epsilon_1,\epsilon_2) = (a+2-2j,a-2j)$, for $0\leq j\leq \frac{a}{2}$. 
For $(yxy)^2$, we will need to take two different partitions, each of which gives rise to a 
labelled skew tableau whose weight is $1^{x+y}2^y$, namely the partitions $a+2-(y+1)\geq a+1-(x+y-2)$ 
and $a+2-(y-1)\geq a+1-(y+x)$.
For the first partition, we  label the tableaux with $(-)^{a+2-(y+1)}1^{y+1};(-)^{a+1-(x+y-2)}1^{x-2}2^{y};1$. 
For the second 
partition, we label via $(-)^{a+2-(y-1)}1^{y-1};(-)^{a+1-x-y}1^{x+1}2^{y-1};2$.

Now turn to repeated summands of the form $(rsr)$, $r$ and $s$ both even. Each of these will occur 
in $(0(a+2)000)\downarrow X$, and as in the case of $a$ odd, those of the form $(0s0)$, for $s\geq 2$ 
occur with multiplicity $2$, those of the form $(r0r)$ for $2\leq r\leq a$ with multiplicity $3$, and 
those of the form $(rsr)$ with $2\leq r$, $s\geq 2$ even and $r+s\leq a+2$, with multiplicity $4$ in the 
tensor product $(0(a+1)000)\downarrow X\otimes (101)$ (or multiplicity 2 if $r+s=a+2$).
So we need to show that 
$(r0r)$ (for $2\leq r\leq a$) and $(rsr)^2$ (for $r+s\leq a$) occur in $(1a100)\downarrow X$.
For $(r0r)$, where $2\leq r\leq a$ is even, take the partitions $(a-2j,a-2j)$ for $0\leq j\leq \frac{a-2}{2}$.
Finally for $(rsr)$, we take the partition $(a+2-r, a+1-(r+s-1))$, and indicate two 
different
labelled skew tableaux, each of which will have weight $1^{r+s}2^r$.
The two labellings are: $(-)^{a+2-r}1^r;(-)^{a+2-r-s}1^s2^{r-1};2$ and 
$(-)^{a+2-r}1^r;(-)^{a+2-r-s}1^{s-1}2^{r};1$. \hal


 \subsection{$X = A_4, \d = \om_2$}\label{A4}

Let $X = A_4 < Y = SL(W) = A_9$, where $W = V_X(\o_2)$. 
In this subsection we consider two infinite families in Table \ref{TAB4}, namely the restrictions to $X$ of the 
$Y$-modules with highest weights $a\l_1+\l_9$ or $a\l_1+\l_2$.

 \begin{lem}\label{a4(a0...01)}  For any $a\ge 1$, the restriction $V_Y(a\l_1 + \l_9))\downarrow X$ is MF.
\end{lem}

\pf We first observe that  $a\l_1 + \l_9 =  (S^a(\l_1) \otimes \l_9) - S^{a-1}(\l_1).$    Restricting
to $X$ this becomes $(S^a(0100) \otimes (0010)) - S^{a-1}(0100).$

Applying \cite[3.8.1]{Howe}, we find that 
\begin{equation}\label{35}
S^a(\l_1) \downarrow X = (0a00) + (0(a-2)01) + (0(a-4)02) + \cdots
\end{equation}
Indeed,  it follows from  \cite[3.8.1]{Howe}  that composition factors have the form $0x0y$ and correspond to partitions  of the form $(1^{x+y},2^{x+y},3^y,4^y)$, subject to the condition
$2x+4y=2a$. This condition reflects the fact that the center of $GL_5$  acts with weight $2a$ on  $S^a(0100)$. All such composition factors occur and have
the form $0a00 -( t\a_1 +2t\a_2 +t\a_3)$ as indicated. Similarly for
$S^{a-1}(\l_1) \downarrow X.$

Next we apply Corollary \ref{LR_om2} to  verify that 
\begin{equation}\label{36}
\begin{array}{ll}
(0x0y) \otimes (0010) = & (0(x-1)0y) + (0(x+1)0(y-1)) + (0x1y) + \\
   &   (1(x-1)1(y-1)) +(1(x-1)0(y+1)),
\end{array}
\end{equation}
noting that there are some missing terms if $xy = 0.$
Applying (\ref{36}) to the terms in (\ref{35}) we see that we only obtain a multiplicity  from  consecutive
terms in (\ref{35}) and these occur from the first two summands in (\ref{36}).  For example,
both $(0(a-2)01)$ and $(0(a-4)02)$ yield a term $(0(a-3)01)$.  Therefore we get
a series of terms appearing with multiplicity 2.  Indeed, these are $(0(a-1)00), (0(a-3)01),
\ldots .$  But each such term
appears in $S^{a-1}(0100).$  Indeed the repeated terms exhaust $S^{a-1}(0100)$ except when $a = 2k+1$
in which case $ (000k)$ appears in $S^{a-1}(0100)$ but it only appears with multiplicity 1 in $(S^a(0100) \otimes (0010))$.  This establishes the lemma.  \hal

\begin{lem}\label{a4(a10...}  For any $a \ge 1$, the restriction $V_Y(a\l_1+\l_2) \downarrow X$ is  MF.
 \end{lem}
 
\pf  This proof is very similar to the previous one. Note that 
$((a+1)0\ldots  0)\otimes (10\ldots 0) = ((a+2)0\ldots  0)+(a10\ldots 0)$.
Restricting to $X$ we then have 
$V_Y(a\lambda_1+\lambda_2)\downarrow X = S^{a+1}(0100))\otimes(0100)-S^{a+2}(0100)$.
Now we apply Corollary \ref{LR_om2} to verify that 
\[
\begin{array}{ll}
(0x0y)\otimes(0100) = & (0(x+1)0y)+(1(x-1)1y)+ (0(x-1)0(y+1))+(1x0(y-1))+ \\
  & (0(x-1)1(y-1)).
\end{array}
\]
 Again, using (\ref{35}), 
 we see that there can only 
be repetitions between two successive terms of $S^{a+1}(0100)$ tensored with $(0100)$ and then only one repeated summand 
in the sum of two terms. For example, the tensor products $(0(a+1)00)\otimes(0100)$ and 
$(0(a-1)01)\otimes(0100)$ each have a summand $(0a01)$, and the tensor products  
$(0(a-1)01)\otimes(0100)$ and $(0(a-3)02)\otimes(0100)$ each have a summand $(0(a-2)02)$.
Each of these repeated summands occurs in $S^{a+2}(0100)$, establishing the result. \hal

  \subsection{Remaining Table \ref{TAB4} configurations}

The remaining configurations in Table 4 occur when either $(X,\d)=(A_4,\o_2)$ or  $(A_c,\om_3)$ with $c = 5,6,7$, and $\l$ is one of a few possible weights.   The result is as follows

\begin{lem}\label{a5,a6,a7}  Assume $(X,\d)=(A_4,\o_2)$, $(A_{13},\o_7)$ or  $(A_c,\om_3)$ with $c = 5,6,7$. Then $V_Y(\l)\downarrow X$ is MF for each of the following weights $\l$:
\begin{itemize}
\item[{\rm(i)}] $X=A_4$, $\l = 4\l_2$ $5\l_2$, $2\l_3$, $2\l_4$ or $\l_1+2\l_2$;
\item[{\rm (ii)}] $X = A_5$, $\l = \l_i$ or $\l_1 + \l_{18}$;
\item[{\rm (iii)}] $X = A_6$, $\l = \l_5$ or $\l_6$;
\item[{\rm (iv)}] $X = A_7$, $\l = \l_5$;
\item[{\rm (v)}] $X = A_{13}$, $\l = \l_3$.
\end{itemize}
\end{lem}

\pf   These are all verified using Magma. The computations are straightforward in all cases except (i) with $\l = a\l_2$, so we give some details for this case. 

We have $X=A_4$, and $X < SL(W) = A_9$, where $W = V_X(\o_2)$. Using Magma for representations of $A_9$, we compute that $V_{A_9}(a\l_2)  = V^+-V^-$, where $V^+,V^-$ are as in Table \ref{vplusminus}.

\begin{table}[h]
\caption{}\label{vplusminus}
\[
\begin{array}{|l|l|l|}
\hline
\l & V^+ & V^- \\
\hline


4\l_2 & S^4\l_2 +(\l_4\otimes \l_4)+(\l_1\otimes \l_2\otimes \l_5) & (S^2\l_2\otimes \l_4)+S^2\l_4+ 
(S^2\l_1\otimes \l_6) +  \\
  &  & (\l_3\otimes \l_5) \\

5\l_2 & S^5\l_2 + (S^2\l_2 \otimes \l_1\otimes \l_5) + & (S^3\l_2\otimes \l_4) + (S^2\l_4\otimes \l_2)+\\
        &  (S^2\l_1\otimes \l_1\otimes \l_7) + (\l_2\otimes \l_2\otimes \l_6) + &  
(S^2\l_2\otimes \l_6) +(S^2\l_1\otimes \l_2\otimes \l_6) +  \\
        & (\l_2\otimes \l_4\otimes \l_4) + (\l_3\otimes \l_7)^2 + & S^2\l_5 + (\l_2\otimes \l_3\otimes \l_5) + \\
        & (\l_1\otimes \l_9)^2      & (\l_1\otimes \l_2\otimes \l_7) + (\l_1\otimes \l_1\otimes \l_8) + 0^2 \\
\hline
\end{array}
\]
\end{table}

We now use Magma for representations of $A_4$ to restrict each of the modules $V^+$ and $V^-$  to $X$, using the fact that $\l_i \downarrow X = \wedge^iW$ for all $i$. For later use we give the complete decompositions $V_Y(a\l_2)\downarrow X$ for all $2\le a\le 5$ in Table \ref{comple}; in particular, these are all MF. \hal

\begin{table}[h]
\caption{}\label{comple}
\[
\begin{array}{|l|l|l|}
\hline
\l & V_Y(\l) \downarrow X \\
\hline
2\l_2 & 0 0 0 2+  0 2 0 1+  2 0 2 0 + 1 1 0 0 \\

3\l_2  & 0 1 0 0+ 1 2 1 1+ 2 1 1 0+ 1 1 0 1+ 0 2 1 0+ 3 0 3 0+ 1 0 1 2 \\

4\l_2 & 3 1 2 0+  2 0 0 1+ 2 2 2 1+ 2 2 0 0+ 1 2 2 0+ 2111+ \\
 &   0 0 0 4+ 0 2 0 3+ 2 0 2 2+ 0 2 1 1+ 1 1 0 2+\\
   &   4 0 4 0+ 1 1 1 0 + 0 0 2 0 + 1 3 0 1+ 0 4 0 2 \\

5\l_2 &  0411+ 1103+ 5050+ 2311+ 1014+ 2120+  \\
    &    1030+  1200+ 2112+ 0102+ 0212+  3032+ \\
    &     0301+ 4130+ 2230+ 1310+ 1213+ 1412+ \\
     &     3011+ 1111+ 1221+ 2201+ 3121+ 1302+ \\
     &     2010+ 3210+ 3231+ 0000 \\
\hline
\end{array}
\]
\end{table}

\chapter{Initial Lemmas}\label{initiallemmas}

In this chapter we establish a range of preliminary lemmas that will be required in later chapters.
The chapter is divided into three subsections.  In the first subsection we prove  lemmas  which determine particular summands in tensor products or in alternating or symmetric powers of modules.   These results are established using either Littlewood-Richardson or domino techniques, as described in Chapter \ref{quote}.   The second subsection establishes a number of  results showing that various modules are not multiplicity-free.  The final subsection gives lower bounds for the  $L$-values of various modules.

\section{Summands of Tensor Products}\label{summands}  In the next several lemmas, let $G = A_n$  with fundamental dominant weights $\omega_1, \ldots , \omega_n$ and fundamental roots $\a_1,\ldots,\a_n$. For a dominant weight $\psi$ we shall often just write $\psi$ to denote the irreducible $G$-module $V_G(\psi)$; similarly $\psi\otimes \xi$ denotes $V_G(\psi) \otimes V_G(\xi)$.

 \begin{lem}\label{abtimesce}  Let $G = A_n$, $\psi = \sum a_j\omega_j $ and $\xi = \sum b_k\omega_k.$  Assume that $i,j$ are such that $i<j$ and  $a_i,b_j \ne 0$, and let $\mu= (\psi + \xi) - (\alpha_i + \dots + \alpha_j).$ Then the following hold.
 \begin{itemize}
\item[{\rm (i)}]   $\psi \otimes \xi$ has a composition factor of highest weight $\mu.$
 \item[{\rm (ii)}]  If $b_{i} = b_{i+1} = \cdots = b_{j-1} = 0$, then $\mu$ appears with multiplicity $1$ in $\psi \otimes \xi$.
 \item[{\rm (iii)}] Assume $n \ge 3$, $\psi = a_1\om_1 + a_2\om_2,$  $\xi = b_{n-1}\omega_{n-1} + b_n\om_n,$  and $a_1a_2b_{n-1}b_n \ne 0.$  Then  $(\psi + \xi) - (\alpha_1 + \dots + \alpha_n)$ appears with multiplicity $1$ in $\psi \otimes \xi$.
\end{itemize}
\end{lem}
 
\pf  (i) Working in a proper Levi factor, if necessary, we can assume that $i = 1$ and $j = n.$ Then $\mu = (a_1+b_1-1, a_2+b_2, \ldots  ,a_{n-1}+b_{n-1}, a_n+b_n-1).$ 
 
We will use Littlewood-Richardson skew tableaux as in Theorem \ref {LR}.  In order to establish (i)  it will suffice to show that there exists at least one labelling of certain tableaux for which the corresponding partition gives rise to  the dominant weight  $\mu $.   For (ii) and (iii) we must show that the labelling is unique.

To simplify notation we use the following notation.  For $k =1, \ldots , n$ let $a(k) = a_k + \cdots + a_n$, $b(k) = b_k + \cdots + b_n,$ and
$ab(k) = a(k) + b(k).$
With this notation the   partitions   $\d_{\psi} = (a(1), a(2), \ldots , a(n),0)$
and $\d_{\xi} =  (b(1), b(2), \ldots , b(n),0)$ correspond to $\psi$ and $\xi$, respectively.

We  need a partition $\nu$ that
 corresponds to the dominant weight $\mu$.  One partition corresponding to $\mu$ is $(ab(1)-2, ab(2)-1, \ldots , ab(n)-1,0).$  But we need a partition $\nu$ corresponding to $\mu$   such that $|\nu| = |\d_{\psi}| + |\d_{\xi}|.$
 Therefore we take  $\nu = (ab(1)-1, ab(2), \ldots ,  ab(n),1).$
 
 We now construct  a Littlewood-Richardson tableau of shape $\nu/\d_{\psi}$ for which the labelling  equals the weight of $\d_{\xi}$.  Let $r_1, \ldots , r_{n+1}$ denote the rows of the tableaux.  In the following  we indicate the entries of these rows, letting ``$x$" be a placeholder for a blank cell:
  \vspace{2mm}
  
  $r_1:  (x^{a(1)}, 1^{b(1)-1})$
  
  $r_2:  (x^{a(2)},1^1, 2^{b(2)-1})$
  
  $r_3:  (x^{a(3)},2^1, 3^{b(3)-1})$
  
  \indent \indent $ \vdots$
  
  $r_{n-1}:  (x^{a(n-1)},(n-2)^1, (n-1)^{b(n-1)-1})$
  
  $r_n: (x^{a(n)}, (n-1)^1, n^{b(n)-1})$
  
  $r_{n+1} : (n^1)$
  
   \vspace{2mm}
\no This tableau satisfies the various conditions and it has labelling equal to the weight of $\d_{\xi}$.  So this establishes (i).

(ii) Again we can assume $i = 1$ and $j = n$, so by hypothesis  $\xi =(0\ldots  0b_n).$  We must show that the above labelling is the only possible labelling of the tableaux satisfying the conditions and with the weight  that of  $\d_{\xi}$. Note that $b(j) = b_n$ for $j = 1, \ldots , n.$  The Y-condition will give the assertion.  In the first row the non-blank cells are necessarily labelled by 1's, so $r_1$ is as above.  Now consider $r_2$ where there are $b(2) = b(1)$ non-blank cells.  The  Y-condition implies that these cells are labelled by 1's followed by 2's.  But there is a single $1$  available and the Y-condition implies that there are at most $(b(1)-1)$ 2's.  Therefore $r_2$ must be as above.
 Continuing in this way we see that the above labelling is the only one possible.

(iii)  Assume  $n \ge 3$, $\psi = a_1\om_1 + a_2\om_2,$  $\xi = b_{n-1}\omega_{n-1} + b_n\om_n,$  and $a_1a_2b_{n-1}b_n \ne 0.$ Again we aim to show that the above labelling is the only one possible.
 We have $\d_{\psi} = (a(1), a(2), 0 \ldots , 0),$  $b(1) = \cdots = b(n-1) = b_{n-1}+b_n$, and $b(n) = b_n.$  The argument of (ii) shows that rows $r_1, \ldots , r_{n-1}$ must be labelled as above.  So consider  $r_n$. As $n \ge 3$, $a(n) = 0$ so there are no  blank entries in the row. The remaining entries  in the tableau are $(n-1)^1$ and $n^{b(n)}$.  The column decreasing condition implies that  $r_n$ cannot be labelled as $n^{b(n)}$, so the only possibility is as above.  \hal

 \begin{lem}\label{aibjnotzero}  Assume $1 < i \le j < n$.  Let $\l = \sum a_k\om_k$ and $\xi =\sum b_k\om_k.$ Assume that $a_i \ne 0 \ne b_j.$   Set $\mu = (\l + \xi) - (\a_{i-1} + 2\a_i + \cdots + 2\a_j + \a_{j+1}).$
\begin{itemize}
\item[{\rm (i)}]  Then $\l \otimes \xi \supseteq \mu$.  
\item[{\rm (ii)}]  If, in addition, $a_{i-1} \ne 0 \ne b_{j+1}$, then $\l \otimes \xi \supseteq (\mu)^2$. 
\end{itemize}
\end{lem}

\pf  Working in a proper Levi factor, if necessary, we can assume that $i = 2$ and $j = n-1.$  We will use Littlewood-Richardson skew tableaux as in Theorem \ref {LR}.  In order to establish (i) (respectively (ii)) it will suffice to show that there exists at least 1 (respectively 2) labellings of certain tableaux for which the corresponding partition gives rise to  the dominant weight $\mu $. 

Using the notations of the previous lemma we have
partitions   $\d_{\l} = (a(1), a(2), \ldots , a(n))$
and $\d_{\xi} =  (b(1), b(2), \ldots , b(n))$ corresponding to $\l$ and $\xi$, respectively.
We next need a partition $\nu$ that
 corresponds to the dominant weight $\mu$ and such that $|\nu| = |\d_{\l}| + |\d_{\xi}|.$ We take $\nu = (ab(1)-1, ab(2)-1,ab(3), \ldots , ab(n-1), ab(n)+1,1).$
 
 We now construct  a Littlewood-Richardson tableau of shape $\nu/\d_{\xi}$ for which the labelling  equals the weight of $\d_{\l}$.
  We start with first row of the tableau for $\nu.$  After ignoring
  the $b(1) = b_1 +\cdots +b_n$ cells corresponding to $\d_{\xi}$ there remain $a(1)-1$ cells which we label with 1's.  Similarly for row $r_2$ where there are $a(2)-1$ cells available which we label with 2's.  In rows $r_k$ ($3\le k\le n-1$), there are are an additional $a(k)$ cells available per row and we label each of these cells with $k$.  We note that to satisfy the Y-condition for $3$ we require our assumption $a_2 \ne 0.$
  
  Now consider row $r_n$ where there are $ab(n)+1 = a_n+b_n+1$ cells.  The first $b_n$ are ignored as these are correspond to $\d_{\xi}$.  We can label the remaining ones as $1, n, \ldots  ,n$, but in order to satisfy the
  strictly increasing condition on columns we need our assumption that $b_{n-1} \ne 0$ as otherwise, the $1$ in $r_n$ would lie below an $n-1$ in the tableau.
  Row $r_{n+1}$ has size $1$ and we label the blank cell with a $2$.  This yields a LR tableaux and establishes (i).
  In the following we  illustrate the above tableaux by listing for each row $r_1, \ldots , r_{n+1}$ the labellings of the cells from left to right.  We use ``$x$" as a place holder for an empty cell.
  
  \vspace{2mm}
  
  $r_1:  (x^{b(1)}, 1^{a(1)-1})$
  
  $r_2:  (x^{b(2)}, 2^{a(2)-1})$
  
  $r_3:  (x^{b(3)}, 3^{a(3)})$
  
  \indent \indent $ \vdots$
  
  $r_{n-1}:  (x^{b(n-1)}, (n-1)^{a(n-1)})$
  
  $r_n: (x^{b(n)}, 1^1, n^{a(n)})$
  
  $r_{n+1} : (2^1)$
  
   \vspace{2mm}
   
  To establish (ii) we make only a minor adjustment to the above. Label rows $r_1, \ldots , r_{n-1}$ as above. If we also label rows $r_n$ and $r_{n+1}$ as above then we get  the first possible labelling.  For the second possibility   label the blank cells in row $r_n$ as  $2, n, \ldots  ,n$ and  the blank cell in row $r_{n+1}$  with a $1$.  The resulting array must be column strictly increasing and this forces $b_n \ne 0.$  In addition applying the Y-condition  to $2$ we must have  $a_1 \ne0.$  This completes the proof.  \hal

 \begin{lem}\label{aibj12}  Assume $1 \le i \le j \le n$.  Let $\l = \sum a_k\om_k$ and $\xi =\sum b_k\om_k.$  If $a_i, b_j \ge 2$, then $\l \otimes \xi \supseteq (\l + \xi) - (2\a_i + \cdots + 2\a_j).$ 
\end{lem}

\pf  The proof here is similar to that  of Lemma \ref{aibjnotzero}. 
 Working within a Levi factor, if necessary, we reduce to the case $i = 1$ and $j = n.$ We  use  notation  as before where for $k =1, \ldots , n$ let $a(k) = a_k + \cdots + a_n$, $b(k) = b_k + \cdots + b_n,$ and $ab(k) = a(k) + b(k).$ 

The   partitions   $\d_{\l} = (a(1), a(2), \ldots , a(n))$
and $\d_{\xi} =  (b(1), b(2), \ldots , b(n))$  correspond to $\l$ and $\xi$, respectively.
 We next need a partition $\nu$ that corresponds to the dominant weight $(\l + \xi) - (2\a_1 + \cdots + 2\a_n)$ and such that $|\nu| = |\d_{\l}| + |\d_{\xi}|.$ We take $\nu = (ab(1)-2, ab(2), ab(3), \ldots , ab(n),2).$
 
  We now describe  a Littlewood-Richardson tableau of shape $\nu/\d_{\l}$ for which the labelling  equals the weight of $\d_{\l}$.
  We start with first row of the tableau for $\nu.$  This row begins with    $a(1) = a_1 +\cdots +a_n$ blank cells. There remain $b(1)-2$ cells which we label with 1's.  Continuing, in rows $r_k$ for $2 \le k \le n$  there are there are $a(k)$  blank cells followed by $(k-1)^2$ and then $k^{b(k)-2}.$   
 Finally, row $r_{n+1}$  is labelled $n^2.$  This yields a LR tableau and establishes the result. \hal

 \begin{lem}\label{july4/17} Let $G = A_n$ and assume $\l$ and $\mu$ are dominant weights with $L(\l), L(\mu) \ge 2.$ Then  $\l \otimes \mu \supseteq \nu^2$, where $\nu$ is a dominant weight such that $S(\nu) \ge  ) + S(\mu) -2.$
  \end{lem}
  
\noindent{\bf Remark.} The weight $\nu$ will be given explicitly in the course of the proof.
  
  \vspace{2mm}
 
 \pf
We will make use of the following fact.  Let $L$ be a standard Levi subgroup of
$G$ of type $A$ and let $\l_L$ and $\mu_L$ be the respective restrictions of $\l$ and $\mu$ to $L$. If $\l_L \otimes \mu_L \supseteq ((\l_L \otimes \mu_L)- \sum c_i\a_i)^2$ where
the $\a_i$ are fundamental roots in  $\Pi(L)$, then $\l \otimes \mu \supseteq ((\l \otimes \mu)- \sum c_i\a_i)^2.$

 Write  $\l = c_{i_1}\l_{i_1} + c_{i_2}\l_{i_2} + \cdots ,$  where $i_1 < i_2 < \cdots$ and each coefficient is nonzero. 
Similarly  write $\mu =  d_{j_1}\l_{j_1} + d_{j_2}\l_{j_2} + \cdots $ with $j_1 < j_2 < \cdots.$  We will use the above to reduce to a Levi subgroup of $G$.  There are three cases.

 (i) Up to interchanging the roles of $\l$ and $\mu,$ there exist $i_k$ and $j_l$ such that $i_k \le j_l <j_{l+1} \le i_{k+1}.$ Applying the first paragraph and changing notation we  reduce to  $\l = a\l_1 + b\l_n$ and $\mu = c_1\l_{i_1} + \cdots + c_k\l_{i_k},$ where $abc_1 \ldots  c_k \ne 0$ and  $1 \le i_1 < \cdots < i_k \le n.$  Set $\nu = \l + \mu - (\a_1 + \cdots + \a_n).$

(ii)  Here (i) does not hold and  there exist $i_k$ and $j_l$ such that $i_k < j_l < i_{k+1} < j_{l+1}.$ 
This time we can reduce to the case  $\l =   a\l_1 + b\l_j$ and $\mu = c\l_k + d\l_n$, where $1 < k < j < n$ 
and $abcd \ne 0.$ We again set  $\nu = \l + \mu - (\a_1 + \cdots + \a_n).$

  (iii) Neither (i) nor (ii) hold and there exist $i_k$ and $j_l$ such that $i_k < i_{k+1} \le j_l < j_{l+1}$ and $j_l$ is minimal for this. In this case  a change of notation reduces us  to $\l = a\l_1 + b\l_i$ and  $\mu = c\l_j + d\l_n$, where $i\le j$.  Here we set $\nu = \l + \mu - (\a_1 + \cdots \a_{i-1} +2\a_i + \cdots + 2\a_j + \a_{j+1} \cdots + \a_n).$
 
 \vspace{2mm}
 
 We will work through the above cases  using  Theorem \ref{LR}, starting with (ii).
 \vspace{2mm}
 
{\bf Case (ii) } In this case we take  partitions for $\l$ and $\mu$  with weights $1^{a+b}2^b \ldots  j^b$ and $1^{c+d}2^{c+d} \ldots  k^{c+d}(k+1)^d \ldots  n^d$ respectively.  
Here  $\nu =  (a-1)\l_1 +c\l_k +b\l_j + (d-1)\l_n$  and we take the  partition with weight:
 $$1^{a+b+c+d-1}2^{b+c+d} \ldots  k^{b+c+d}(k+1)^{b+d} \ldots  j^{b+d}(j+1)^d\ldots  n^d(n+1)^1.$$
 
 We construct two Littlewood-Richardson tableaux of shape $\nu/\l$.  These have rows as in the following, where in the $r_i$ entry, we list the $ith$ row of the first tableau, then the $ith$ row of the second:
\[ 
\begin{array}{ccc} 
r_1:  & (x^{a+b},1^{c+d-1}) & (x^{a+b},1^{c+d-1}) \\
r_2: & (x^b,1^1,2^{c+d-1}) & (x^b,1^1, 2^{c+d-1}) \\
\vdots & \vdots & \vdots \\
r_k: & (x^b,(k-1)^1,k^{c+d-1}) &  (x^b,(k-1)^1, k^{c+d-1}) \\
r_{k+1}: & (x^b,k^1,(k+1)^{d-1}) &  (x^b, (k+1)^d) \\
\vdots & \vdots & \vdots \\
r_j :&  (x^b,(j-1)^1,j^{d-1}) & (x^b, j^d) \\
r_{j+1}: & (j^1, (j+1)^{d-1}) & (k^1, (j+1)^{d-1}) \\
r_{j+2}: & ((j+1)^1,(j+2)^{d-1}) & ((j+1)^1,(j+2)^{d-1}) \\
\vdots & \vdots & \vdots \\
r_n :&  ((n-1)^1, n^{d-1}) & ((n-1)^1, n^{d-1}) \\
r_{n+1}: &  (n^1) &  (n^1) 
\end{array}
\]
where the $r_{j+2}$ row does not occur if $j+1=n$.

  \vspace{2mm}
   
{\bf Case (i) } Here we take the partition for $\l$ with weight $1^{a+b}2^b \ldots  n^b$.  Recalling the notation
$c(t) = c_t+ \cdots +c_k$ for $1 \le t \le k$ we take the partition for $\mu$ with weight $1^{c(1)} \ldots  i_1^{c(1)}(i_1+1)^{c(2)} \ldots  i_2^{c(2)} \ldots  i_k^{c(k)}.$ 
Then $\nu = (a-1)\l_1 + c_1\l_{i_1} + \cdots +c_k\l_{i_k} + (b-1)\l_n$   and corresponding to $\nu$ we use the partition with weight
 $$1^{a+b +c(1)-1}2^{b+c(1)} \ldots  i_1^{b+c(1)}(i_1+1)^{b+c(2)} \ldots  i_2^{b+c(2)} \ldots  i_k^{b+c(k)}(i_k+1)^b \ldots  n^b(n+1)^1,$$ 
where the terms $(i_k+1)^b,\ldots , n^b$ do not occur if $i_k = n$, and the terms $2^{b+c(1)},\ldots ,i_1^{b+c(1)}$ do not occur if $i_1=1$.
 
  We construct two Littlewood-Richardson tableaux of shape $\nu/\l$.  The first has rows as follows:
  
  \vspace{2mm}
   
   $r_1: (x^{a+b},1^{c(1)-1})$
   
   $r_2: (x^b,1^1,2^{c(1)-1})$
 
  \indent \indent $ \vdots$

$r_{i_1}: (x^b,(i_1-1)^1,i_1^{c(1)-1})$

$r_{i_1+1}: (x^b, i_1^1,(i_1+1)^{c(2)-1})$

  \indent \indent $ \vdots$
  
 $r_{i_2}: (x^b, (i_2-1)^1, i_2^{c(2)-1})$
 
 $r_{i_2+1}: (x^b, i_2^1,(i_2+1)^{c(3)-1})$
 
 \indent \indent $ \vdots$

$r_{i_k} : (x^b,(i_k-1)^1,i_k^{c(k)-1})$

 $r_{i_k+1} : (x^b)$
 
 \indent \indent $ \vdots$
 
 $r_n : (x^b)$
 
 $r_{n+1} : (i_k^1)$
 
 \vspace{2mm}
 
\noindent  Note that some collapsing takes place if either $i_1 = 1$ or $i_k = n$.  
The second tableau is obtained by interchanging the terms $i_k^1$ and $(i_k-1)^1$ which appear in rows $r_{n+1}$ and $r_{i_k}$, respectively.

  \vspace{2mm} 
{\bf Case (iii) } We take partitions for $\l$ and $\mu$  with weights $1^{a+b}2^b \ldots  i^b$ and $1^{c+d}2^{c+d} \ldots  j^{c+d}(j+1)^d \ldots  n^d$ respectively.   We will construct two Littlewood-Richardson tableaux of shape $\nu/\mu$.
 There are special cases where $i = 2$, $i = j$, or $j = n-1.$   For the general case  $2 < i < j < n-1, $ 
  $$\nu = (a-1)\l_1 + \l_{i-1}+ (b-1)\l_i + (c-1)\l_j + \l_{j+1} + (d-1)\l_n$$
  and we use the partition with weight
 $$1^{a+b+c+d-1}2^{b+c+d} \ldots  (i-1)^{b+c+d}i^{b+c+d-1}(i+1)^{c+d} \ldots  j^{c+d}(j+1)^{d+1}(j+2)^d\ldots  n^d(n+1)^1.$$
 
  We now construct the  tableaux. The first tableau has
   rows as follows:
   
    \vspace{2mm}
   
   $r_1: (x^{c+d},1^{a+b-1})$
   
   $r_2: (x^{c+d},2^b)$
 
  \indent \indent $ \vdots$

$r_{i-1}: (x^{c+d},(i-1)^b)$

$r_i: (x^{c+d},i^{b-1})$

$r_{i+1}: (x^{c+d})$

  \indent \indent $ \vdots$
  
  $r_j : (x^{c+d})$
  
  $r_{j+1}: (x^d,i^1)$
  
  $r_{j+2} :(x^d)$
  
    \indent \indent $ \vdots$
    
   $r_n : (x^d)$
   
   $r_{n+1}: (1^1)$
   
    \vspace{2mm}
  
\noindent  For the second tableau we simply interchange the terms $i^1$ and $1^1$ which appear in rows $r_{j+1}$ and $r_{n+1},$ respectively.  The tableaux for the special cases are entirely similar and we leave these to the reader. \hal

The following is an immediate corollary of three cases in the proof of Lemma \ref{july4/17}.

\begin{lem}\label{7.5cor} Let $G = A_n$ and assume $\l$ and $\mu$ are dominant weights with $L(\l), L(\mu) \ge 2$ satisfying one of the conditions (i), (ii), (iii) below. Then  $\l \otimes \mu \supseteq \nu^2$, where $\nu$ is as indicated and  $S(\nu) \ge S(\l) + S(\mu) -2.$ 
\begin{itemize} 
\item[{\rm (i)}] $\l = a\l_1 + b\l_n$, $\mu = \sum c_i\l_i$ with $\nu = \l + \mu - (\a_1 + \cdots + \a_n).$  \item[{\rm (ii)}] $\l = a\l_1 + b\l_j$, $\mu = c\l_k +d\l_n$ with $1 \le k < j <n$  and $\nu = \l + \mu - (\a_1 + \cdots + \a_n).$ 
 \item[{\rm (iii)}] $\l = a\l_1 + b\l_i$, $\mu = c\l_j +d\l_n$ with $1 <i \le j <n$  and $\nu = \l + \mu -(\a_1 + \cdots \a_{i-1} +2\a_i + \cdots +2\a_j + \a_{j+1} + \cdots +\a_n).$ 
\end{itemize} 
\end{lem}

\begin{lem}\label{abtimescd}  Let $G = A_n$ with $n \ge 2$ and let $\nu_1 = (a,b,0 \dots 0)$,
$\nu_2 = (0,\dots, 0, c,d)$  and $\lambda = \nu_1 + \nu_2$ be dominant weights for $G$.  
\begin{itemize}
\item[{\rm (i)}] Suppose $a,b,c,d \ne 0$. Then the $G$-composition factor with highest weight $\lambda - (\alpha_1 + 2\alpha_2 + \dots + 2\alpha_{n-1} + \alpha_n)  = (a,b-1,0, \dots ,0,c-1,d)$ appears in $\nu_1 \otimes \nu_2$ with multiplicity  $2$ and has $S$-value $S(\lambda)-2.$ 
\item[{\rm (ii)}] Suppose $b=c=0$, $a\ge 3$ and $d\ge 2$. Then the $G$-composition factor with highest weight $(2a-2,0,\ldots ,0,2d-2)$ appears in $\wedge^2(\nu_1) \otimes \wedge^2(\nu_2)$ with multiplicity  $2$.
\item[{\rm (iii)}]  Suppose that $abcd \ne 0,$ $G$ is a proper Levi subgroup of $Y = A_{n+t}$, $G$ has base $\{\alpha_1, \dots, \alpha_n\}$, and $\bar \nu_1$, $\bar \nu_2$ are dominant weights for $Y$ which restrict to $G$ as $\nu_1$, $\nu_2$, respectively. Set $ \bar \lambda =  \bar \nu_1 + \bar \nu_2$. Then $V_Y(\bar \lambda_1) \otimes V_Y(\bar \lambda_2) $ has a composition factor of highest weight $\bar \lambda - (\alpha_1 + 2\alpha_2 + \dots + 2\alpha_{n-1} + \alpha_n)$ which appears with  multiplicity at least $2$ and has $S$-value at least $S(\lambda)-1$.
\end{itemize}
\end{lem}

 \pf   (i) For $n = 2$, the weight in question is  $\lambda - (\alpha_1 + \alpha_2)$.  Then the tensor product  has summands of highest weights $\lambda, \lambda -\alpha_1,$ and $ \lambda - \alpha_2$.  The multiplicity of weight $\lambda - (\alpha_1 + \alpha_2)$ in  these summands is $2, 1, 1,$ respectively.  On the other hand, in the tensor product this weight  has multiplicity $6$ and the assertion follows.
 
 For $n \ge 3$,  Lemma \ref{7.5cor} shows that an irreducible with the indicated weight does occur with multiplicity at least 2 in the tensor product.  Moreover in the special situation of this lemma the proof of that lemma shows that the two composition factors  produced are the  only ones possible.  So this gives the assertion about the multiplicity  being precisely 2.  The assertion about $S$-values is obvious.

(ii)  The $G$-module $\wedge^2(d\omega_n)$ has a summand of highest weight $\omega_{n-1}+(2d-2)\omega_n$, and $\wedge^2(a\omega_1)$ has a summand of highest weight $(2a-2)\omega_1+\omega_2$. Then Lemma~\ref{aibjnotzero} shows that the tensor product of these modules contains that $V_G((2a-2)\omega_1+(2d-2)\omega_n)$ with multiplicity 2.

 (iii) This follows from (i); note that the increased $S$-value holds since $G$ is a proper Levi subgroup of $Y$.  \hal
 
 \begin{lem}\label{compfactor} Let $G = A_n$ and $\lambda = (r0 \ldots  0s)$ with $rs \ne 0.$ Let $\mu = (c0\ldots  0d0 \ldots  0e)$.
\begin{itemize}
\item[{\rm (i)}]  If $cde \ne 0,$ then $ \lambda \otimes \mu \supseteq ((\lambda + \mu) - (\alpha_1 + \cdots+  \alpha_n))^3.$ \item[{\rm (ii)}]  If $cd \ne 0,$ $ de \ne 0$, or $ce\ne 0,$ then $\lambda \otimes \mu \supseteq ((\lambda + \mu) - (\alpha_1 + \cdots+  \alpha_n))^2.$ \end{itemize}
   \end{lem}
 
 \pf  This follows from the proof of Lemma \ref{july4/17}. The form of $\l$ and $\mu$ imply that we are in  case (i) of the proof of Lemma \ref{july4/17}. Therefore part (ii) is immediate from that proof.  For part (i) suppose label $d$ occurs at node $j$.  Then we have $i_1 = 1,$ $i_2 = j,$ and $i_3 = n.$ The first two tableaux are the ones given in the proof of Lemma \ref{july4/17}.  To get the  third tableau replace the entries $1^1, \ldots , (n-1)^1$ appearing in rows $r_2, \ldots , r_n$ by $2^1, \ldots , n^1$  and then  replace  the entry in $r_{n+1}$ by $1^1.$  \hal

\begin{lem}\label{lemma1.5(i)(ii)}  Let $G = A_n$ and let $\lambda = (a0\ldots  0b0 \ldots  0c)$
 with $a,c \ne 0.$ 
\begin{itemize} 
\item[{\rm (i)}]  If $b \ne 0$, then  both $\wedge^2(\lambda)$ and $S^2(\lambda)$ contain $2\lambda - (\alpha_1 + \cdots + \alpha_n)$ with multiplicity at least $2.$
 \item[{\rm (ii)}]  If $b = 0$, then both $\wedge^2(\lambda)$ and $S^2(\lambda)$ contain $2\lambda - (\alpha_1 + \cdots + \alpha_n).$
\end{itemize}
   \end{lem}
   
  \pf   (i) We proceed using the domino technique to study the tensor square of $\l$. 
Assuming that $b$ occurs at node $k$ we use the partition
 $$(1^{a+b+c},2^{b+c}, \ldots , k^{b+c},(k+1)^c, \ldots , n^c).$$
 Doubling and repeating each exponent we obtain the sequence
 $$(2(a+b+c), 2(a+b+c), 2(b+c), \ldots , 2(b+c), 2c, \ldots , 2c),$$
  where the terms $2(b+c)$ and $2c$ occur $2(k-1)$ times and $2(n-k)$ times, respectively. 
 Forming the half sum of the terms of the sequence we obtain $2a+2kb+2nc$, so this is the
 number of dominoes in the array.
 
 Now let $\d = 2\l - (\a_1+ \cdots \a_n) = ((2a-1)0 \ldots  0\,2b\,0\ldots  0(2c-1)).$  We need a partition
 corresponding to $\d$ for which the sum of the exponents is the number of dominoes.  Consequently we
 use the  partition with weight
 $$(1^{2a+2b+2c-1},2^{2b+2c},\ldots ,k^{2b+2c},(k+1)^{2c}, \ldots , n^{2c},(n+1)^1).$$
We now describe tilings of the array.  Assume the base consists of $s$ pairs of a horizontal $2$-domino lying above a horizontal $1$-domino and $t$ vertical dominoes.  Then we have equations $2s+t = 2(a+b+c)$ and
$s+t = 2a+2b+2c-1.$   Therefore $s = 1$ and $t = 2a+2b+2c-2.$  The base of the diagram is now determined.

The level above  the base has $2b+2c-1$ vertical $2$-dominoes
followed by $1$ vertical $3$-domino.  At this point all the required 1's and 2's are accounted for.  At the
next level there are $2b+2c-1$ vertical $3$-dominoes
followed by $1$ vertical $4$-domino.   We continue in this way until we get a level with  $2b+2c-1$ vertical $k$-dominoes followed by $1$ vertical $(k+1)$-domino. 

The next level has $2c-1$ vertical $(k+1)$-dominoes
followed by $1$ vertical $(k+2)$-domino.  Continue  until at the top we have $2c-1$ vertical $n$-dominoes
followed by $1$ vertical $(n+1)$-domino.  The tiling satisfies all the conditions and corresponds to $\d$.
There are precisely $2$ horizontal dominoes so this corresponds to an alternating summand of the tensor
square of $\l$.

Slight changes yield the additional summands required.  First note that if we change the top level to $2c-2$ vertical $n$-dominoes followed by a horizontal $(n+1)$-domino lying above a horizontal $n$ domino, then we also have
a labelling corresponding to $\d$ and this time it has $4$ horizontal dominoes.  Hence this yields a symmetric
summand.

To obtain additional summands we return to the level where there were $2b+2c-1$ vertical $k$-dominoes
followed by $1$ vertical $(k+1)$-domino. In each of the two summands described above replace that level by $2b+2c-2$ vertical $k$-dominoes  followed by a horizontal $(k+1)$-domino lying above a horizontal $k$ domino. This yields additional summands.  The alternating (respectively symmetric) summand above is converted to a symmetric (respectively alternating) summand.   So this yields the result.

(ii)  Here we assume $b = 0.$ Following through the above argument we see that the rows above the base all have the same size and we produce two composition factors of highest weight $2\lambda - (\alpha_1 + \cdots + \alpha_n)$.  One is alternating and the other is symmetric.   \hal
 
   \begin{lem}\label{wedgesquareab}  Assume that $G = A_n$ for $n \ge 2$ and that $\lambda = (a0 \ldots  0b).$
\begin{itemize}
\item[{\rm (i)}]  If $a\ge3$, $b \ge 2,$ then  $\wedge^2(\lambda)$ is not MF.  Indeed   $\wedge^2(\lambda)\supseteq ((2a-4)10\ldots  0(2b-2))^2$ or $((2a-4)(2b-1))^2$ according as $n \ge 3$ or $n =2.$ 
\item[{\rm (ii)}] If $a \ge b \ge 2,$ then  $S^2(\lambda)$ is not MF.  Indeed $S^2(\lambda) \ge ((2a-2)0\ldots  0(2b-2))^2.$
\end{itemize}
\end{lem}

\pf  We prove both parts using the domino technique.  The proofs are quite similar.

 (i) We will establish the assertion under the assumption  $n \ge 3.$ The $n = 2$ case is almost identical. Let $\mu =   ((2a-4)10\ldots  0(2b-2)). $ We begin
with the partition having weight $(1^{a+b},2^b, \ldots ,n^b)$ which corresponds to $\l$.  Doubling and repeating each exponent we have the sequence $((2a+2b),(2a+2b),2b, \ldots , 2b)$ where the term $2b$ appears $2(n-1)$ times.
Adding the terms of this sequence we find that number of $1 \times 1$ tiles in the array is $4a+4nb,$ so that there will be a total of $2a+2nb$ dominoes.  
 
 Therefore we require a partition corresponding to $\mu$ such  that the sum of the exponents in its weight is $2a+2nb$.  
One partition corresponding to $\mu$ is $(1^{2a+2b-5},2^{2b-1},3^{2b-2}, \ldots , n^{2b-2}).$  The sum of the exponents here is $2a+2nb-2-2n$, so we replace this partition by  $$(1^{2a+2b-3},2^{2b+1},3^{2b}, \ldots , n^{2b},(n+1)^2).$$
  We then look for the corresponding tilings of the above array.
 The number of dominoes in the bottom two rows is $2a+2b$.  So if there are $t$ vertical 1-tiles and  $s$ pairs of a horizontal 2-tile above a horizontal 1-tile, then we must have $2s+t = 2a+2b$ and  since all 1-dominoes must be in the bottom two rows we also have $s+t = 2a+2b-3$.  Therefore $s = 3,$ $t = 2a+2b-6,$ and the base is determined.
 
We now describe two tilings.  For the first tiling just above the base we set  $2b-2$ vertical 2-dominoes followed by 2 vertical 3-dominoes. Similarly, above this there are   $2b-2$ vertical 3-dominoes followed by 2 vertical 4-dominoes. Continue in this way, until finally there are $2b-2$ vertical $n$-dominoes followed by 2 vertical $(n+1)$-dominoes.  

In the second tiling the first $2b-4$ columns are exactly as in the first labelling.  Consider columns $2b-3, 2b-2, 2b-1, 2b.$ 
At the base these have 4 vertical $1$-dominoes.  Above these we place
2 horizontal $2$-dominoes.  Then we begin a series of  4 vertical dominoes with labels $3, \ldots , n$ followed by 2 horizontal 
$(n+1)$-dominoes. The tilings have $6$ and $10$ horizontal dominoes, respectively.  
Therefore both are alternating, completing the proof of (i).

(ii) For this proof $\l$, its corresponding partition, and the array are all as above.  Set $\nu =   ((2a-2)0\ldots  (2b-2)). $  The weight of one partition corresponding to $\nu$ is $(1^{2a+2b-4},2^{2b-2},3^{2b-2}, \ldots , n^{2b-2})$ and here the sum of the exponents is again equal to $2a + 2nb-2n-2.$ Therefore we replace this partition by
$$(1^{2a+2b-2},2^{2b},3^{2b}, \ldots , n^{2b}, (n+1)^2).$$

We  first tile the base.  If there are $t$ vertical dominoes and $s$ pairs of a horizontal 2-tile above a horizontal 1-tile, then we get the equations $2s+t = 2a+2b$ and $s + t = 2a+2b-2.$  Therefore $s = 2$ and $t = 2a+2b-4.$  This determines the tiling of the base and we move up the diagram as before.

The tilings here are very similar to those in (i).  For the first tiling just above the base we use $2b-2$ vertical $2$-dominoes followed by $2$ vertical $3$-dominoes.  Above this we have $2b-2$ vertical $3$-dominoes followed by $2$ vertical $4$-dominoes and so on. We continue  until just below the top of the array we have $2b-2$ vertical $n$-dominoes followed by $2$ vertical $(n+1)$ dominoes.  In the second tiling the first $2b-4$ columns are exactly as in the first tiling.  And we label
columns $2b-3, 2b-2,2b-1,2b$ precisely as in the second tiling in (i).     Once again all the conditions are satisfied but this time there are either $4$ or $8$ horizontal dominoes, respectively.  So both tilings yield symmetric composition factors, completing the proof.  \hal


\begin{lem}\label{nonMF_tensors} 
Let $X = A_{l+1}$ ($l\ge 1$), and let $x,y,z\geq 2$. Then the following hold:
\begin{itemize}
\item[{\rm (i)}] $V_X((2z+x-3)\omega_1+\omega_2+(y-1)\omega_{l+1})$ occurs with multiplicity two in $V_X(x\omega_1)\otimes V_X(y\omega_{l+1})\otimes \wedge^2V_X(z\omega_1)$.
\item[{\rm (ii)}] $V_X((2z+x+y-4)\omega_1+2\omega_2)$ occurs with multiplicity two in $V_X(x\omega_1)\otimes V_X(y\omega_1)\otimes \wedge^2V_X(z\omega_1)$.
\item[{\rm (iii)}]  $V_X((x+y-1)\omega_1+(2z-1)\omega_{l+1})$ occurs with multiplicity two in $V_X(x\omega_1)\otimes V_X(y\omega_1)\otimes \wedge^2V_X(z\omega_{l+1})$.
\item[{\rm (iv)}] $V_X(\omega_1+\omega_{l}+2\omega_{l+1})$ occurs with multiplicity two in 
$(V_X(\omega_1+2\omega_{l+1})\oplus V_X(2\omega_1+3\omega_{l+1}))\otimes V_X(\omega_l)$.
\item[{\rm (v)}] Suppose $2\le i\le 4$ and $l\geq i$. Then $V_X(\omega_1+\omega_{l+1-i}+2\omega_{l+1})$ occurs with multiplicity two in 
$(V_X(\omega_1+2\omega_{l+1})\oplus V_X(2\omega_1+3\omega_{l+1}))\otimes V_X(\omega_{l+1-i})$.
\end{itemize}
\end{lem}

\pf This is a straightforward application of Lemmas~\ref{abtimesce} and \ref{compfactor}(ii), as well as 
some easy calculations with exterior squares.\hal

\begin{lem}\label{wedge3adj}  Let $G = A_n$ for $n \ge 2$ and let $\lambda = 10\ldots 01$.  Then $\wedge^3(\lambda) $ contains an irreducible summand of highest weight $3\lambda -(2\alpha_1 + \alpha_2 + \cdots +\alpha_n)$ and multiplicity $1$. 
\end{lem}

\pf  This is checked using Magma for $n \le 4$, so assume $n \ge 5.$  The dominant weights in $\wedge^3(\lambda) $ strictly above $\sigma = 3\lambda -(2\alpha_1 + \alpha_2 + \cdots +\alpha_n)$ are $3\lambda - \alpha_1-\alpha_n, 3\lambda - 2\alpha_1-\alpha_2$, and $3\lambda - (\alpha_1+ \cdots +\alpha_n)$ and the multiplicities of $\sigma$ in the  irreducibles with these highest weights
are $2(n-2), n-2, 1$, respectively.  Easy counting arguments show that each of these irreducibles appears
with multplicity $1$ in $\wedge^3(\lambda) $, so this gives a total of $3n-5$ appearances of $\sigma.$ 
On the other hand $\sigma$ can occur in $\wedge^3(\lambda) $ using the wedge of the following triples

$(\lambda, \lambda - \alpha_1, \lambda - (\alpha_1+ \cdots +\alpha_n)),$

$(\lambda,  \lambda - \alpha_1-\alpha_n, \lambda - (\alpha_1+ \cdots +\alpha_{n-1})),$

$(\lambda, \lambda - \alpha_1-(\alpha_{n-1} + \alpha_n), \lambda -(\alpha_1+ \cdots +\alpha_{n-2})) ,$

{$\vdots$} 

$(\lambda, \lambda - \alpha_1-(\alpha_3 + \cdots +\alpha_n), \lambda -(\alpha_1+ \alpha_2)), $

$(\lambda - \alpha_1, \lambda - \alpha_n, \lambda - (\alpha_1 + \cdots + \alpha_{n-1}))$

$(\lambda - \alpha_1, \lambda - (\alpha_{n-1} + \alpha_n), \lambda - (\alpha_1 + \cdots + \alpha_{n-2}))$

$\vdots$

$ (\lambda - \alpha_1, \lambda - (\alpha_3 + \cdots + \alpha_n), \lambda - (\alpha_1 +  \alpha_2)).$

The first item on the list appears with multiplicity $n$.  All the rest occur with multiplicity $1$ and they occur in two
batches of size $n-2.$  This gives a total multiplicity of $3n-4$   which then implies that $\wedge^3(\lambda) $
contains a irreducible summand of highest weight $\sigma$ and multiplicity $1$.   \hal

\begin{lem}\label{wedge2adj} The following hold for $G = A_n$ with $n \ge 3$:
\begin{itemize}
\item[{\rm (i)}] $\wedge^2(\om_1+\om_n) = (\om_1+\om_n) \oplus  (\om_2+2\om_n)\oplus  (2\om_1+\om_{n-1})$;
\item[{\rm (ii)}] $S^2(\om_1+\om_n) = (2\om_1+2\om_n)\oplus  (\om_1+\om_n)\oplus  (\om_2+\om_{n-1}) \oplus  0.$
\end{itemize}
\end{lem}
\pf  This is easily established using the domino technique as in preceding proofs.  We leave the details
to the reader.  \hal

\section{Some non-MF representations}\label{techlemmas}

In this section we establish a number of results showing that certain modules 
are not MF.  The results come in five basic flavours:
\begin{itemize}
\item[(1)] restrictions $V_Y(\l)\downarrow X$, where $X = A_{l+1} < Y = SL(W)$ and $W = V_X(\d)$ for $\d = \o_2$,
\item[(2)] as for (1), but with $\d = 2\o_1$,
\item[(3)] modules of the form $S^2(V_X(\d))$, $\wedge^2(V_X(\d))$ for $X = A_{l+1}$ and various highest weights $\d$,
\item[(4)] modules for $A_{l+1}$ with $l$ small,
\item[(5)] tensor products of various modules for $X = A_{l+1}$.
\end{itemize}

We divide the section into five subsections accordingly.


 \subsection{Non-MF modules for $\d = \o_2$}\label{nonmfo2}

In this subsection, we adopt the following notation:
\begin{equation}\label{o2not}
\begin{array}{l}
X=A_{l+1} \hbox{ with } l\ge 4,\\
W = V_X(\o_2), \\
X < Y = SL(W) = A_n.
\end{array}
\end{equation}
We shall need some of the notation of  Chapter \ref{levelset}. By Theorem \ref{LEVELS}, 
there are two levels $W^1(Q_X) \cong V_{L_X'}(\o_2)$ and $W^2(Q_X) \cong V_{L_X'}(\o_1)$,
so $L_X' < L_Y' = C^0 \times C^1,$  where $C^0 = A_{r_0}$ with $r_0 = \frac{(l+1)l}{2}-1$, and $C^1 = A_l.$ 
For a dominant weight $\l$ of $T_Y$, let $\mu^i$ be the restriction of $\l$ to $T_Y\cap C^i$, so that 
$V^1(Q_Y) = V_{C^0}(\mu^0) \otimes V_{C^1}(\mu^1)$. By Proposition \ref{induct}(i), if $V\downarrow X$ is MF, then also 
$V^1$ is MF, where we recall the notation from Chapter \ref{notation}:
\[
V^i = V^i(Q_Y) \downarrow L_X'.
\]

\begin{lem}\label{om2a1000}  Adopt the notation in $(\ref{o2not})$, and let $V = V_Y(a\l_1+\l_2)$ with $a \ge 4. $ Then $V\downarrow X$ is not MF.
\end{lem}

\pf By the preceding remarks, it suffices to show that $V^1$ is not MF, and inductively it then suffices to consider the case $l=4$.
So assume that $l=4$. Now $Y = A_{14}$ and we have the decomposition of $Y$-modules
$$((a+1)0\ldots 0)\otimes (10\ldots 0) = ((a+2)0\ldots  0)+(a10\ldots  0).$$
Since $S^b((01000))$ is MF by \cite[3.8.1]{Howe}, it suffices to show that 
for $b\geq 5$, $S^b((01000))\otimes (01000)$ has a multiplicity three summand. 
Now we will use \cite[3.8.1]{Howe}, as in the proof of Lemma~\ref{a4(a0...01)}, to obtain three 
particular summands of $S^b((01000))$. Here \cite[3.8.1]{Howe} implies that
 the composition factors have the form $(0x0y0)$, corresponding to partitions with weight of the form $(1^{x+y+z}, 2^{x+y+z}, 3^{y+z}, 4^{y+z}, 5^z,6^z)$ satisfying $2x+4y+6z = 2b$. Taking $z=1$, $y=0$ and $x=b-3$ gives $(0 (b-3) 0 0 0)$, further $z=1=y$ and $x=b-5$ gives $(0 (b-5) 0 1 0)$, and $z=0$, $y=2$ and $x=b-4$ gives $(0(b-4)0 2 0)$.

Now tensoring each of these with $(0 1 0 0 0)$, using Corollary \ref{LR_om2}, we obtain three summands $(0(b-4) 0 1 0)$, giving the desired conclusion. \hal


\begin{lem}\label{nolab1} Adopt the notation in $(\ref{o2not})$, and let $\l = \l_1 + \l_i + \l_n$ with $2 \le i \le 7$. Then $V_Y(\l)  \downarrow X$ is not MF.  
\end{lem}

\pf  For $l=4,5,6$ the assertion can be checked using Magma, so assume $l\ge 7$. 
It is convenient to replace $\l$ by the dual $\l_1+\l_{i'}+\l_n$, where $i' = n-i+1 > r_0+1$. 
By way of contradiction assume that  $V \downarrow X$ is  MF, where $V = V_Y(\l)$.  

Now $V^1 = \om_2 \otimes (\om_j + \om_l)$, where $j = i' - (r_0+1)$,  and  $V^2(Q_Y)$ contains the following $L_Y'$-summands:

\vspace{2mm}

(1)  $(\l_1^0 + \l_{r_0}^0) \otimes (\l_{j+1}^1 + \l_l^1)$

(2)  $(\l_1^0 + \l_{r_0}^0) \otimes  \l_j^1$

(3) $0 \otimes (\l_{j+1}^1 + \l_l^1)$

(4) $0 \otimes \l_j^1$.

\vspace{2mm}
 \no The sum of (1) and (3) is  $(\l_1^0 \otimes \l_{r_0}^0) \otimes (\l_{j+1}^1 + \l_l^1)$ and 
 the sum of (2) and (4) is $(\l_1^0 \otimes \l_{r_0}^0) \otimes \l_j^0.$  Therefore the sum
 of all four terms is $(\l_1^0 \otimes \l_{r_0}^0) \otimes (\l_{j+1} \otimes \l_l^1)$,  and
 restricting to $L_X'$ this becomes $\om_2 \otimes \om_{l-1} \otimes \om_{j+1} \otimes \om_l.$
 The tensor product of the middle two terms contains $(\om_{j+1} + \om_{l-1}) \oplus (\om_j + \om_l)$
 so the full tensor product contains
\begin{equation}\label{2nds}
(\om_2 \otimes (\om_{j+1} +\om_{l-1}) \otimes \om_l) \oplus (\om_2 \otimes (\om_j + \om_l) \otimes \om_l).
\end{equation}
 By Corollary \ref{cover} the second summand in (\ref{2nds}) contains $\sum_{i,n_i=0}V_i^2(Q_X)$.  By Proposition \ref{stem1.1.A} the first summand is  not MF unless $l-1 = j+1$.  We claim that this summand is not MF in the latter case as well. Indeed, here the first summand becomes $\om_2 \otimes 2\om_{l-1} \otimes \om_l.$   Writing highest weights as sequences, this contains $((010\ldots  20) + (10\ldots  011))\otimes(0\ldots  01)$ which contains $(10\ldots  20)^2$. Thus we have a contradiction to Proposition \ref{induct} in all cases, completing the proof. \hal

\begin{lem}\label{nolab2}  With notation as in $(\ref{o2not})$, let  $\l = a\l_1 + b\l_n$ with $a\ge b \ge 2$, and $V = V_Y(\l)$. Then $V \downarrow X$ is not MF.
\end{lem}

\pf We will show that $V \downarrow X \supseteq \nu^3$, where $\nu = (a-1)\om_2 + (b-1)\om_l$.  

We begin by claiming that  $S^a(010\ldots 0) \otimes S^b(0\ldots  010) \supseteq ((a-1)\om_2 + (b-1)\om_l)^4.$
First note that $S^a(010\ldots 0) \supseteq (0a0\ldots  0) + (0(a-2)010\ldots 0)$ and similarly
$S^b(0\ldots 010) \supseteq (0\ldots 0b0) + (0\ldots 010(b-2)0).$  

Set $\gamma_1 = (0a0\ldots  0),$  $\gamma_2 = (0\ldots 0b0)$,$\gamma_3 = (0\ldots 010(b-2)0)$, and
$\gamma_4 = (0(a-2)01,0\ldots 0).$ Then it follows from Lemma \ref{aibjnotzero} that $\g_1 \otimes \g_2, \g_1 \otimes \g_3, \g_2 \otimes \g_4$, and $\g_3 \otimes \g_4$ each contain $\nu.$  So this establishes the claim.

We can now complete the proof.  Using  Proposition \ref{pieri} we have $a\l_1 \otimes b\l_n = (a\l_1 + b\l_n) + ((a-1)\l_1 + (b-1)\l_n) + \cdots + (a-b)\l_1$.  Since $\l_1 \downarrow T_X = \om_2$ and  $ \l_n\downarrow T_X =  \om_l$ it follows that $\nu$ can only appear in the restriction to $X$ of the first two summands and   it appears with multiplicity 1 in $((a-1)\l_1 + (b-1)\l_n)$.  Therefore $V \downarrow X \supseteq ((a-1)\om_2 + (b-1)\om_l)^3$ as asserted.
This completes the proof.  \hal



In the proof of the next result, we use a result of Howe  (see 4.4.4 of \cite{Howe}) which shows
how to produce maximal vectors in the module  $\wedge^i(\om_2)$ for $A_{l+1}$.  
Let $v_1, v_2, \ldots , v_{l+2}$ be a basis for the natural module of $A_{l+1}$, chosen such that the vectors afford weights $\om_1, \om_1-\a_1, \om_1-\a_1-\a_2, \ldots , \om_1-\a_1-\a_2- \cdots  -\a_{l+1}.$  For $i<j$ set $e_{ij} = v_i \wedge v_j,$  so that these elements form a basis for $V_{A_{l+1}}(\om_2)$.  We now consider an array as follows:

\vspace{2mm}

$e_{12} \ e_{13} \ e_{14} \ e_{15} \ \ldots  \ e_{1,l+2}$

\ \ \ \ \ $e_{23}  \ e_{24} \ e_{25} \ \ldots  \ e_{2,l+2}$

\ \ \ \ \ \ \ \ \ \ $e_{34} \ e_{35} \ \ldots  \ e_{3,l+2}$

\ \ \ \ \ \ \ \ \ \ $\vdots$  \ \ \ \ \ \ \ \ \ \ \ $\vdots$

\noindent Next we list the weights of the above:

$\om_2, \  \om_2 -\a_2, \  \om_2 -\a_2-\a_3, \ldots $

\ \ \ \ \ \  $\om_2-\a_1-\a_2, \  \om_2 -\a_1-\a_2-\a_3, \ldots $

\ \ \ \ \ \ \ \ \ \ \ \ \ \ \ \ \ \ \ \ \ \ \ \ \  $\om_2-\a_1-2\a_2-\a_3, \ \om_2-\a_1-2\a_2-\a_3-\a_4, \ldots  $

\ \ \ \ \ \ \ \ \ \ $\vdots$  \ \ \ \ \ \ \ \ \ \ \ $\vdots$

We say a set $S\subseteq \{e_{rs}| 1\leq r\leq s\leq l+2\}$ is increasing if,
 whenever $e_{ij}$ lies in $S$, 
all $e_{rs}$ with $r\leq i$ and $s\leq j$ also lie in $S$ -- that is, $S$ is increasing if all $e_{rs}$ above
or to the left of any element in $S$ also lie in $S$.
In \cite[4.4.4]{Howe}, it is shown that the wedge
product of an increasing set of size $j$ yields a maximal vector of $\wedge^j(\o_2)$, and that this process yields
all irreducible summands of the exterior algebra of $V_{A_{l+1}}(\o_2)$.

\begin{lem}\label{(s+2)cge2} With notation as in $(\ref{o2not})$, let $\l = \l_l+c\l_n$ with $c\ge 1$, and $V = V_Y(\l)$. If $c=1$, suppose further that $l\ge 7$. Then $V \downarrow X$ is not MF.
\end{lem}

\pf  First consider the case where $c\ge 2$. Assume that $l \ge 5.$  
We shall use level analysis for the parabolic of $X$ with Levi factor $L_X' = A_l$, with notation as in Chapter \ref{levelset}.
The result will follow provided we can show that the $A_l$-module $V^1 = \wedge^l(\om_2)  \otimes c\o_l$ is not MF.

We now use the theory described in the preamble to the lemma to find some composition factors of $\wedge^l(\om_2)$.
Consider the two following increasing subsets of length $l$:
$e_{12}, e_{13}, \ldots , e_{1l}, e_{23}$ and $e_{12}, e_{13}, \ldots , e_{1(l-1)}, e_{23}, e_{24}.$
The wedges of these elements have weights $\xi_1 = (l-3)\om_1+\om_3+\om_l$ and $\xi_2 = (l-5)\om_1+\om_2 + \om_4 + \om_{l-1}$, respectively.

Now consider $\wedge^l(\om_2) \otimes c\om_l.$ In the first tensor factor there are irrreducible summands
with the highest weights $\xi_1$ and $\xi_2$.  Using Proposition \ref{pieri}, we see that each of these tensors 
with $c\om_l$ and yields an irreducible of highest weight $\nu = ((l-4)\om_1 + \om_3 + \om_{l-1} + (c-2)\om_l.$   
So this establishes the result for $l \ge 5.$

Now assume $l = 4$, still with $c\ge 2$. Replace $V$ with the dual, which has highest weight $c\l_1 + \l_{r_0+2}.$  Then 
$V^1 = S^c(0100) \otimes (1000).$ We have $V^2 \supseteq (V_{C^0}(c\l_1^0 + \l_{r_0}^0) \otimes V_{C^1}(\l_2^1)) + (V_{C^0}((c-1)\l_1^0) \otimes V_{C^1}(\l_2^1)) = V_{C^0}(c\l_1^0 ) \otimes V_{C^0}(\l_{r_0}^0)  \otimes V_{C^1}(\l_2^1).$  Restricting to $L_X'$ this becomes $S^c(0100) \otimes (0010) \otimes (0100).$ Now $(0100) \otimes (0010) = (0110) + (1001) + (0000) = (0110) + ( (1000) \otimes (0001)).$  Therefore
$$V^2  \supseteq (S^c(0100) \otimes (0110)) +
(S^c(0100) \otimes (1000) \otimes (0001)).$$
Corollary \ref{cover} shows that the second summand contains $\sum_{i,n_i=0}V_i^2(Q_X).$  If
$c \ge 3$  then $S^c(0100) \supseteq (0(c-2)01)$ so it follows from Proposition \ref{tensorprodMF} that the
first summand is not MF.  Therefore the result holds in this case.  Finally assume $c = 2.$
Then the first term becomes $S^2(0100) \otimes (0110) \supseteq (1010)^2.$  So the result also holds
here.

Now assume that $c=1$, so $\l=\l_l+\l_n$ and $l\ge 7$ by hypothesis.
We will show that  $V\downarrow X = (\wedge^l(\om_2) \otimes \om_l) - \wedge^{l-1}(\om_2)$ is not MF.  We know that
the second term is MF, so it will suffice to show that the tensor product contains a composition factor
of multiplicity at least 3.

Towards this end we take three increasing sequences, $(l-3,3),(l-3,2,1),(l-2,2)$, of $l$ terms from the above $e_{ij}$ array. The notation indicates that for the first sequences we
take the first $l-3$ terms in the first row of the array and the first 3 terms of the second row.  Similarly for
the other sequences.  Taking wedge products of the terms of each sequence we get a maximal vector
in $\wedge^l(\om_2)$ and one checks that these maximal vectors have weights  $\g_1 = ((l-7)20010\ldots 01000), \g_2 = ((l-6)0020\ldots 01000), \g_3 = ((l-5)1010\ldots 0100),$ respectively.

At this point we use Corollary \ref{LR_om2} to see that  $ \g_i \otimes \om_l \supseteq ((l-6)1010 \ldots  01000)$ for $i = 1,2,3$.
This completes the proof. \hal


\begin{lem}\label{cl1+ls} With notation as in $(\ref{o2not})$, let $\l = \l_{l+2} + c\l_n$ with $c\ge 1$, and $V = V_Y(\l)$.  
\begin{itemize}
\item[{\rm (i)}] If $l\ge 5$ and $1\le c\le 3$, then $V \downarrow X$ is not MF.
\item[{\rm (ii)}] If $l=4$ and $c\ge 2$, then $V \downarrow X$ is not MF.
\end{itemize}
\end{lem}

\pf  (i) Assume $l\ge 5$ and $1\le c\le 3$. 
Here $V^1 = \wedge^{l+2}(\om_2) \otimes c\om_l.$  

 We begin with an analysis of $\wedge^{l+2}(\om_2)$  and
as in the last result we will consider three increasing subsets of length $l+2$ in the  array given in the proof of 
Lemma \ref{(s+2)cge2}. These are $(l,2), (l-1,3),(l-1,2,1),$ where as before the  notation means that for $(l,2)$ we take the 
first $l$ entries in row 1 of the array and the first 2 entries in the  row 2 of the array, etc.  We then wedge the entries
in each case, obtaining maximal vectors in $\wedge^{l+2}(\om_2)$ affording 
 the following weights of $X$:  $\g_1 = ((l-3)1010 \ldots  01)$, $\g_2 = ((l-5)20010 \ldots  010)$,
$\g_3 = ((l-4)0020\ldots 010),$ respectively.  Let $\nu_1 = ((l-3)1010 \ldots  0)$, $\nu_2 =  ((l-5)20010 \ldots  01)$,
 $\nu_3 = ((l-4)0020 \ldots  01)$ be the corresponding weights of $L_X’$.

First assume $c \ge 2.$  Here we claim that $V^1 = \wedge^{l+2}(\om_2) \otimes
 c\om_l$ is not MF.  Indeed an application of Proposition \ref{pieri} shows that 
$\nu_1 \otimes c\om_l$ and $\nu_2 \otimes c\om_l$ both contain $(l-4)1010\ldots  0(c-1)$. This establishes the 
lemma for $c\ge 2$  and $l \ge 5.$

Now assume that $c=1$. Observe that $V = (\wedge^{l+2}(\om_2) \otimes \om_l) - \wedge^{l+1}(\om_2).$  
A Magma check shows that $V \downarrow X$ is not MF for $l = 5,6,$ so  assume $l \ge 7.$
Now Corollary \ref{LR_om2} shows that $\g_i\otimes \o_l$ contains $((l-4)1010 \ldots  010)$ for $i=1,2,3$.
The latter therefore  occurs with multiplicity 3 in $\wedge^{l+2}(\om_2) \otimes \om_l$.
This establishes the result because we know that $\wedge^{l+1}(\om_2)$ is MF.

(ii) Here we consider the case $l = 4$, with $c\ge 2$. Replace $V$ by $V^*$.
Then  $V^2(Q_Y) $ contains
$V_{C^0}((c-1)\l_1^0 + \l_9^0) \otimes V_{C^1}(\l_1^1)$ and $ V_{C^0}(c\l_1^0 + \l_8^0) \otimes V_{C^1}(\l_1^1)$.  These sum to $V_{C^0}(c\l_1^0 ) \otimes V_{C^0}(\l_8^0 ) \otimes V_{C^1}(\l_1^1)$ and restricting to
$L_X'$ this becomes $S^c(0100) \otimes (0101) \otimes (1000) = S^c(0100) \otimes ((1101) +(0011) + (0100))
= (S^c(0100) \otimes (1101)) + (S^c(0100) \otimes (0010) \otimes (0001)).$ 

By Corollary \ref{cover}, the first summand must be MF.  If $c \ge 3$, then  $S^c(0100)$ contains $(0(c-2)01)$ so  that the first summand is not MF.  
And if $c = 2$ a Magma calculation shows that $S^2(0100) \otimes (1101)$ is not MF.  \hal

\begin{lem} \label{newestlem}
 With notation as in $(\ref{o2not})$, let $\l = \l_{l+3}+\l_{n-1}$, and $V = V_Y(\l)$. Then $V \downarrow X$ is not MF.
\end{lem}

\pf By way of contradiction assume $V \downarrow X$ is MF.  We consider  $V^*$ where the label of $C^0 $ is $ (010 \ldots  010)$ and the label of $C^1$ is $(0\ldots 0).$ Then $V^2(Q_Y)$ contains  $V_{C^0}(\l_1^0 +\l_{r_0-1}^0) \otimes V_{C^1}(\l_1^1)$ and
$V_{C^0}(\l_2^0 +\l_{r_0-2}^0) \otimes V_{C^1}(\l_1^1).$ These sum to $(V_{C^0}(\l_2^0) \otimes V_{C^0}(\l_{r_0-2}^0) - V_{C^0}(\l_{r_0}^0)) \otimes V_{C^1}(\l_1^1).$ Restricting to $L_X'$ we find that 
\[
V^2 \supseteq \left(\wedge^2(010\ldots 0) \otimes \wedge^3(0\ldots 010) - (0\ldots  010)\right) \otimes (10\ldots 0).
\]
 Now $\wedge^2(010\ldots 0) = (1010\ldots  0)$ and $\wedge^3(0\ldots 010) = (0\ldots  0200)+ (0 \ldots  0102).$ First assume $l \ge 5.$ Using Theorem \ref{LR} we find that $\wedge^2(010\ldots 0) \otimes \wedge^3(0\ldots 010) \supseteq (010\ldots  012)^1 + (10\ldots  011)^3 + (010\ldots  0101)^3 +(010\ldots 020)^1$.  Tensoring with $(10\ldots 0)$ we see that  $V^2 \supseteq (010\ldots 011)^8$, which is a 
contradiction, since viewing $V^1 \subseteq \wedge^2(\om_2) \otimes  \wedge^2(\om_{l-1})$ we see that only five such summands can arise from $V^1$.
For $l = 4$ the situation is slightly different.  Here $(\wedge^2(0100) \otimes \wedge^3(0010) \supseteq (0112)^1 + (1011)^3 + (0201)^2 +(0120)^1$ and $V^2 \supseteq (0111)^7$, which is again a contradiction.  \hal

\begin{lem}\label{nolab3} With notation as in $(\ref{o2not})$, let $\l = \l_1+\l_i$ with $r_0+1 \le i \le n-6$, and $V = V_Y(\l)$. Then $V \downarrow X$ is not MF.
\end{lem}

\pf Replace $\l$ by the dual $\l^* = \l_j+\l_n$, so that $7\le j\le l+1$.
Observe that $\l_j+\l_n = (\l_j\otimes \l_n)-\l_{j-1}$, so $V\downarrow X = (\wedge^j(\o_2)\otimes \o_l)-\wedge^{j-1}(\o_2)$. We know that $\wedge^{j-1}(\o_2)$ is MF for $X$ by Theorem \ref{wedgec}, so it will suffice to show that 
$(\wedge^j(\o_2)\otimes \o_l)$ has a summand of multiplicity at least 3.

To do this we use the $e_{ij}$ array as in the previous proofs. First we check the assertion for $l=6,7,8$ using Magma, so assume $l\ge 9$. 

Suppose first that $7<j<l+1$. Here we take three increasing sequences $(j-2,2)$, $(j-3,3)$ and $(j-3,2,1)$ of $j$ terms from the $e_{ij}$ array, where as before the notation indicates that for the sequence $(j-2,2)$ we take first $j-2$ terms from the first row and the first 2 terms from the second, and similarly for the other sequences. The three maximal vectors given by taking wedge products of the terms in each of these sequences have the following highest weights:
\begin{itemize}
\item[{\rm (1)}] $(j-5)\o_1+\o_2+\o_4+\o_{j-1}$,
\item[{\rm (2)}] $(j-7)\o_1+2\o_2+\o_5+\o_{j-2}$,
\item[{\rm (3)}] $(j-6)\o_1+2\o_4+\o_{j-2}$.
\end{itemize}
Hence each of these summands appears in $\wedge^j(\o_2)$. Using Corollary \ref{LR_om2}, 
we see that the tensor product of each summand with $\o_l$ has a summand of highest weight 
$(j-6)\o_1+\o_2+\o_4+\o_{j-2}$. Hence this appears with multiplicity 3 in $(\wedge^j(\o_2)\otimes \o_l)$, as required.

Now assume $j=l+1$. In this case we take the three increasing sequences $(l-2,3)$, $(l-3,4)$ and $(l-3,3,1)$, and check that the corresponding maximal vectors have highest weights $(l-6)\o_1+2\o_2+\o_5+\o_{l-1}$, 
$(l-8)\o_1+3\o_2+\o_6+\o_{l-2}$ and $(l-7)\o_1+\o_2+\o_4+\o_5+\o_{l-2}$. The tensor product of each of these with $\o_l$ has a summand $(l-7)\o_1+2\o_2+\o_5+\o_{l-2}$.

Finally, assume $j=7$. Here the increasing sequences $(5,2)$, $(4,3)$ and $(4,2,1)$ yield summands of $\wedge^7(\o_2)$  
of highest weights $2\o_1+\o_2+\o_4+\o_6$, $2\o_2+2\o_5$ and $\o_1+2\o_4+\o_5$, and the tensor product of each of these with $\o_l$ has a summand of highest weight $\o_1+\o_2+\o_4+\o_5$. This completes the proof. \hal 

\begin{lem}\label{nolab4} With notation as in $(\ref{o2not})$, let $\l = \l_1+\l_{r_0+2-i}$, where $2\le i\le 7$. 
Then $V_Y(\l) \downarrow X$ is not MF if either $l \ge 5$ or $l = 4$ and $i = 3.$ \end{lem}

\pf  Replacing $V$ by $V^*$ it will suffice to show that $(\wedge^{l+j}(010\ldots  0) \otimes (0\ldots  010)) - \wedge^{l+j-1}(010\ldots  0)$ fails to be MF, where $2 \le j \le 7$. The subtracted term is MF, so it will suffice to show that the tensor product contains an irreducible summand appearing with multiplicity at least 3.  We use Magma to show this for 
$l = 4,5$. So from now on asssume $l \ge 6$.

For each value of $j$ we produce three sequences which yield composition factors of $\wedge^{l+j}(010\ldots  0)$ with the property that they have a common composition factor upon tensoring with $(0 \ldots  010)$.  In Table \ref{sequ} we list the various values of $j$, followed by the sequences, the corresponding composition factors, and the repeated composition factor.  We leave it to the reader to verify the details. \hal 

\begin{table}[h]
\caption{}\label{sequ}
\[
\begin{array}{|l|l|l|l|}
\hline
j & \hbox{sequences} & \hbox{comp. factors of} & \hbox{ common comp. factor} \\
  &                                &  \wedge^{l+j}(010\ldots  0) &          \\
\hline
2 &  (l,2), (l-1,3),( l-1,2,1) & ((l-3)1010 \ldots  01),  & ((l-4)1010 \ldots  10) \\
   &          & ((l-5)20010 \ldots  010), & \\
   &                                      & ((l-4)0020 \ldots  010)  &  \\
3 &  (l+1,2), (l,3), (l,2,1) & ((l-2)1010 \ldots  0),  & ((l-3)1010 \ldots  01) \\
  &  & ((l-4)20010 \ldots  01), & \\
   &   & ((l-3)0020 \ldots  01) & \\
4 &  (l+1,3), (l,4), (l,3,1) & ((l-3)20010 \ldots  0),  & ((l-4)20010 \ldots  01) \\
  & & ((l-5)300010 \ldots  01), & \\
  &  & ((l-4)10110 \ldots  01) &  \\
5 &  (l+1,3,1), (l,4,1), (l,3,2) & ((l-3)10110 \ldots  0),  & ((l-4)10110 \ldots  01) \\
 & & ((l-5)201010 \ldots  01), & \\
  &  & ((l-4)01020 \ldots  01) & \\
6 &  (l,4,2), (l+1,3,2), (l,3,2,1) & ((l-5)110110 \ldots  01), & ((l-4)01020 \ldots  01) \\
 & & ((l-3)01020 \ldots  0),  & \\
  &  & ((l-4)00030 \ldots  01) & \\
7 &  (l+1,4,2), (l,4,3), (l,4,2,1) & ((l-4)110110 \ldots  0), & ((l-5)110110 \ldots  01) \\
 & &  ((l-5)020020 \ldots  01), & \\
  &  & ((l-5)100210 \ldots  01) & \\
\hline
\end{array}
\]
\end{table}

 \subsection{Non-MF modules for $\d = 2\o_1$}\label{nonmf2o1}

In this subsection, we adopt the following notation:
\begin{equation}\label{2o1not}
\begin{array}{l}
X=A_{l+1} \hbox{ with } l\ge 1,\\
W = V_X(2\o_1), \\
X < Y = SL(W) = A_n.
\end{array}
\end{equation}
Write $\d = 2\o_1$. Again we shall need some of the notation of  Chapter \ref{levelset}.
In this case there are two levels $W^1(Q_X)$ and $W^2(Q_X)$ on which $L_X'$ acts irreducibly with highest weights $2\o_1$ and $\o_1$, respectively, so $L_Y' = C^0\times C^1 \cong A_{r_0}\times A_l$, where $r_0 = \frac{(l+1)(l+2)}{2}-1$. 

We shall need some information about the restrictions of certain modules $V_Y(\nu)$ for the smallest rank case where $X = A_2 < Y = A_5$, recorded in the next result.

\begin{lem}\label{A2info}  Let $X = A_2$ be embedded in $Y = A_5$ via $(20)$, and let $\nu$ be a dominant weight for $Y$ as in Table $\ref{resty}$.  Then the restrictions $V_Y(\nu)\downarrow X$ are as in the table.
\end{lem}

\begin{table}[h]
\caption{}\label{resty}
\[
\begin{array}{|l|l|}
\hline
\nu & V_{A_5}(\nu)\downarrow A_2 \\
\hline
\lambda_1&(20)\\
\lambda_2&(21)\\
\lambda_3&(30)+ (03)\\
\lambda_1+\lambda_2&(41)+ (22)+ (11)\\
\lambda_1+\lambda_3&(50)+ (12)+ (23)+ (31)+ (01)\\
\lambda_1+\lambda_4&(32)+ (13)+ (21)+ (10)\\
\lambda_1+\lambda_5&(22)+ (11)\\
2\lambda_2&(42)+ (31)+ (04)+ (20)\\
3\lambda_2&(63)+ (52)+ (25)+ (60)+ (33)+ (41)+ (22)+ (30)+ (03)\\
\lambda_1+2\lambda_5&(24)+ (40)+ (13)+ (21)+ (02)\\
\lambda_2+\lambda_4&(33)+ (41)+ (14)+ (22)+ (30)+ (03)+ (11)\\
\lambda_2+\lambda_3&(51)+ (24)+ (32)+ (40)+ (13)+ (21)+ (02)\\
\lambda_1+3\lambda_5&(26)+ (15)+ (42)+ (23)+ (04)+ (31)+ (12)+ (20)\\
2\lambda_1+\lambda_2&(20)+ (12)+ (42)+ (23)+ (31)+ (61)\\
3\lambda_1+\lambda_2&(43)+ (32)+ (21)+ (13)+ (02)+ (51)+ (62)+ (40)+ (24)+ (81)\\
a\l_1+\l_2&  (2a+2\,1) + \cdots \\
\hline
\end{array}
\] 
\end{table}

\pf  With the exception of the last row,  these can all be checked using Magma by regarding the
given representation as an alternating sum of tensor products.  For the last row, 
identifying an $A_5$ irreducible representation with its highest weight, we have  
$a\l_1+\l_2 =(a\l_1 \otimes \l_2) - ((a-1)\l_1 \otimes \l_3) + ((a-2)\l_1 \otimes \l_4)
- ((a-3) \l_1 \otimes \l_5) + (a-4)\l_1,$  noting that some of the terms will not appear for small values of $a$.  
On restriction to $A_2$, the irreducible summand whose highest weight has   largest $S$-value
occurs in the first summand $a\l_1\otimes \l_2$ and has highest weight $(2a+2\,1)$. \hal


\begin{lem}\label{nonMF} With notation as in $(\ref{2o1not})$, the following $X$-modules are not MF:
\begin{enumerate}[]
\item{\rm (1)}  $V_Y(2\lambda_2)\downarrow X\otimes (a\omega_{l+1})$, for 
$a=2,3$;
\item{\rm (2)} $V_Y(3\lambda_2)\downarrow X\otimes (a\omega_{l+1})$, for 
$a\geq 1$;
\item{\rm (3)}  $V_Y(\lambda_1+\lambda_j)\downarrow X\otimes (a\omega_{l+1})$, 
for $2\leq j\leq l+2$ and for $a\geq 1$;
\item{\rm (4)} $V_Y(\lambda_1+\lambda_j)\downarrow X\otimes (\omega_k)$, for 
$2\leq j\leq l+2$, $1\leq k\leq l+1$;
\item{\rm (5)} $V_Y(a\lambda_2)\downarrow X\otimes (\omega_k)$, for 
$a=2,3$ and $1\leq k\leq l+1$;
\item{\rm (6)} $V_Y(2\lambda_1+\lambda_{n-1})\downarrow X$;
\item{\rm (7)} $V_Y(\lambda_2+\lambda_3)\downarrow X\otimes (a\omega_{l+1})$, for $a\geq 1$;
\item{\rm (8)} $V_Y(\lambda_1+\lambda_j)\downarrow X\otimes (\omega_{l+1})$,  for $j\geq {\rm max}\{2,n-5\}$.  
\end{enumerate}
\end{lem}


\pf 
For (1), first note that $S^2(V_Y(\lambda_2)) = V_Y(2\lambda_2)\oplus V_Y(\lambda_4)$, and so
$$V_Y(2\lambda_2)\downarrow X\oplus \wedge^4(2\omega_1) = S^2(2\omega_1+\omega_2).$$ From this we 
deduce that 
$V_Y(2\lambda_2)\downarrow X$ has summands $4\omega_1+2\omega_2$ and $4\omega_2$ (applying
\cite[4.4.2]{Howe} to see that these do not occur in $\wedge^4(2\omega_1)$). Now 
tensoring with 
$a\omega_{l+1}$ and applying Lemma \ref{aibj12} gives a repeated summand 
$(2\omega_1+2\omega_2+(a-2)\omega_{l+1})$.

The case (2) is straightforward for $l=1$, using the explicit decomposition of 
$V_Y(3\lambda_2)\downarrow X$ given in Lemmas \ref{A2info} and \ref{abtimesce}; so we now assume that 
$l\geq 2$. Note that $S^3(V_Y(\lambda_2)) = V_Y(3\lambda_2)\oplus V_Y(\l_2+\l_4)\oplus V_Y(\l_6)$, 
and  any $X$-summand of $V_Y(\l_2+\l_4)\oplus V_Y(\l_6)$ occurs as an irreducible summand of 
$\wedge^2(2\omega_1)\otimes\wedge^4(2\omega_1)$.
Recall that $\wedge^2(2\omega_1) = (2\omega_1+\omega_2)$. Hence
 $V_Y(3\lambda_2)\downarrow X$ has a summand $(6\omega_1+3\omega_2)$. In addition, there is a summand 
$(5\omega_1+2\omega_2+\omega_3)=(6\delta-4\alpha_1-\alpha_2)$, since this occurs in 
$S^3(2\omega_1+\omega_2)$,
but does not  occur as a summand of $\wedge^2(2\omega_1)\otimes\wedge^4(2\omega_1)$ 
(since any weight occurring here must have the form $6\delta-\sum_i m_i\alpha_i$ with $\sum_im_i\geq 6$). 

We then deduce using again Lemma ~\ref{abtimesce} that
$(5\omega_1+3\omega_2+(a-1)\omega_{l+1})$ occurs with multiplicity at least two in 
$V_Y(3\lambda_2)\downarrow X\otimes V_X(a\omega_{l+1})$.

The case (3) is a straightfoward check for $l=1,2$, so assume that 
$l\geq 3$. 
We use the fact that $V_Y(\lambda_1)\otimes V_Y(\lambda_j) = 
V_Y(\lambda_1+\lambda_j)\oplus V_Y(\lambda_{j+1})$.
Restricting the tensor product to $X$ gives $\delta\otimes\wedge^j(\delta)$ (recall $\d=2\o_1$). In 
$\wedge^j(\delta)$ there is a summand
$V_X(j\delta-(j-1)\alpha_1-(j-2)\alpha_2-\cdots-\alpha_{j-1}) = V_X(j\omega_1+\omega_j)$,
where $\omega_{l+2}$ should be interpreted as the zero weight, 
 so that $\delta\otimes\wedge^j(\delta)$ has summands $(j+2)\omega_1+\omega_j$ and 
$j\omega_1+\omega_2+\omega_j$. 
Now one checks that 
$(j+2)\omega_1+\omega_j = (j+1)\delta-(j-1)\alpha_1-(j-2)\alpha_2-\cdots-\alpha_{j-1}$ and
$j\omega_1+\omega_2+\omega_j = (j+1)\delta-j\alpha_1-(j-2)\alpha_2-\cdots-\alpha_{j-1}$ do 
not occur as 
summands in $\wedge^{j+1}(\delta)$, using \cite[4.4.2]{Howe}, and hence these summands occur in 
$V_Y(\lambda_1+\lambda_j)\downarrow X$.
Finally, tensoring with $a\omega_{l+1}$ and using Lemma \ref{abtimesce} gives two summands 
$V_X((j+1)\omega_1+\omega_j+(a-1)\omega_{l+1})$,
establishing the result.

For (4), a Magma check handles the cases $l=1,2$, so assume $l\geq 3$.
As in case (3), $V_Y(\lambda_1+\lambda_j)\downarrow X$ has 
summands $((j+2)\omega_1+\omega_j)$ and 
$(j\omega_1+\omega_2+\omega_j)$. 
Now tensoring these with $\omega_k$, and using the Littlewood-Richardson rules Theorem \ref{LR}, we see that 
there are two summands $((j+1)\omega_1+\omega_{k+1}+\omega_j)$.

For (5) in case $l=1$, this is a Magma check, so we assume $l\geq 2$. For the case $a=2$, 
note that 
$S^2(V_Y(\lambda_2))=V_Y(2\lambda_2)\oplus \wedge^4(\lambda_1)$.
Using this we can show that
 $V_Y(2\lambda_2)\downarrow X$ has irreducible summands of highest weights $4\omega_1+2\omega_2$ and 
$3\omega_1+\omega_2+\omega_3$. Now using Lemma \ref{abtimesce}, we see that 
$((4\omega_1+2\omega_2)\oplus (3\omega_1+\omega_2+\omega_3))\otimes \omega_k$ is not MF. For the case 
$a=3$, we use the summands $(6\omega_1+3\omega_2)$ and $(5\omega_1+2\omega_2+\omega_3)$ 
of $V_Y(3\lambda_2)\downarrow X$, as discussed in the consideration of (2).
 Then tensoring with $\omega_k$ is not MF.

Next, (6) is straightfoward: $V_Y(2\lambda_1+\lambda_{n-1})\oplus 
V_Y(\lambda_1+\lambda_n) = 
V_Y(2\lambda_1)\otimes V_Y(\lambda_{n-1})$. One checks that $S^2(2\omega_1)\otimes 
\wedge^2(2\omega_{l+1})$ has a summand $(2\omega_1+\omega_{l})$ occurring with 
multiplicity two, while $V_Y(\lambda_1+\lambda_n)\downarrow X = (2\omega_1+2\omega_{l+1})
\oplus (\omega_1+\omega_{l+1})$.

Next we prove (7). This is straightforward for $l=1$, using the explicit decomposition of
 $V_Y(\lambda_2+\lambda_3)\downarrow X$ given in Lemma \ref{A2info}, together with Theorem \ref{LR}. So now 
assume $l\geq 2$. Note that $V_Y(\lambda_2)\otimes V_Y(\lambda_3) = 
V_Y(\lambda_2+\lambda_3)\oplus V_Y(\lambda_1+\lambda_4)\oplus V_Y(\lambda_5)$. Hence, the 
$X$-modules 
$$V_Y(\lambda_2+\lambda_3)\downarrow X \oplus ((2\omega_1)\otimes \wedge^4(2\omega_1))$$ and
$$\wedge^2(2\omega_1)\otimes\wedge^3(2\omega_1)$$ have the same set of irreducible summands.

Recall that $\wedge^2(2\omega_1) = (2\omega_1+\omega_2)$, and $2\omega_1+\omega_2 = 
2\delta-\alpha_1$.
 Then we find a summand $V_X(5\delta-4\alpha_1)$, in the tensor product 
$\wedge^2(2\omega_1)\otimes\wedge^3(2\omega_1)$ but 
not in the summand
$(2\omega_1)\otimes\wedge^4(2\omega_1)$. Hence $V_Y(\lambda_2+\lambda_3)\downarrow X$  has a summand 
$(2\omega_1+4\omega_2)$.

We now argue that $V_X(5\delta-4\alpha_1-\alpha_2) = V_X(3\omega_1+2\omega_2+\omega_3)$ 
occurs with 
multiplicity two in $\wedge^2(2\omega_1)\otimes\wedge^3(2\omega_1)$. Indeed, 
$\wedge^3(2\omega_1)$ has summands 
$(3\omega_1+\omega_3)$ and $(3\omega_2)$, and tensoring each of these with 
$(2\omega_1+\omega_2)$ 
gives rise to a summand $(3\omega_1+2\omega_2+\omega_3)$. 
We now claim that $V_X(5\delta-4\alpha_1-\alpha_2) = V_X(3\omega_1+2\omega_2+\omega_3)$ occurs
 as a summand of $(2\omega_1)\otimes\wedge^4(2\omega_1)$ with multiplicity exactly one, and 
so  $V_X(3\omega_1+2\omega_2+\omega_3)$ is a summand of 
$V_Y(\lambda_2+\lambda_3)\downarrow X$. To see this, note that any weight 
$\mu=4\delta-\sum_{i=1}^{l+1}a_i\alpha_i$ occurring in
 $\wedge^4(\delta)$ satisfies: \begin{enumerate}
\item $a_1\geq 3$,
\item if $a_3=0$, then $a_j=0$ for all $j\geq 3$, and
\item $\sum_{i=1}^{l+1} a_i \geq 5$.
\end{enumerate}
So we see that the weight $5\delta-4\alpha_1-\alpha_2$ has multiplicity one in the tensor 
product 
$\delta\otimes\wedge^4(\delta)$, and as it is not subdominant to any other weight occurring 
in the tensor product, 
the claim follows.

Now tensoring $(2\omega_1+4\omega_2)$ with $(a\omega_{l+1})$ gives rise to a summand 
$(3\omega_1+3\omega_2+(a-1)\omega_{l+1})$,
and the same is true for the tensor product 
$(3\omega_1+2\omega_2+\omega_3)\otimes (a\omega_{l+1})$.
This completes the proof of (7). 

Finally, we prove (8). 
This is a Magma check for $l=1,2$, so we assume $l\geq 3$ and so $j\geq n-5$. For 
$j\in\{n,n-1\}$, it is easy
to see that $V_X(\omega_1+\delta_{j,n-1}\omega_l+2\omega_{l+1})$ occurs with 
multiplicity at least two
 in the tensor product. 

Consider now the case $j=n-2$, where $V_Y(\lambda_1+\lambda_j)\downarrow X\oplus 
(\omega_l+2\omega_{l+1})$ has the same irreducible summands as 
$\delta\otimes\wedge^3(\delta^*)$. Using that 
$\wedge^3(\delta^*)$ has summands
$(\omega_{l-1}+3\omega_{l+1})$ and $(3\omega_l)$, we deduce then that 
$V_Y(\lambda_1+\lambda_j)\downarrow X$ has irreducible summands $(\omega_1+
\omega_{l-1}+2\omega_{l+1})$ and $(\omega_1+2\omega_l+\omega_{l+1})$.
Now using Theorem \ref{LR}, tensoring each of these with $(\omega_{l+1})$ produces a 
summand $(\omega_1+\omega_{l-1}+\omega_l+\omega_{l+1})$.

We turn now to the case $j=n-3$, and assume $l\geq 4$, the case $l=3$ 
being a straightforward check.
Here the $X$-modules $V_Y(\lambda_1+\lambda_j)\downarrow X\oplus \wedge^3(\delta^*)$ and 
$\delta\otimes \wedge^4(\delta^*)$ have the same
set of irreducible summands. Now $\wedge^4(\delta^*)$ has a summand 
$(\omega_{l-2}+4\omega_{l+1})$.
Tensoring this with $\delta$ we obtain summands $(2\omega_1+\omega_{l-2}+
4\omega_{l+1})$, $(\omega_1+\omega_{l-1}+4\omega_{l+1})$ and 
$(\omega_1+\omega_{l-2}+3\omega_{l+1})$. None 
of these summands occurs in $\wedge^3(\delta^*)$, and hence they occur in 
$V_Y(\lambda_1+\lambda_j)\downarrow X$. 
Now tensoring with $\omega_{l+1}$ produces more than one summand 
$\omega_1+\omega_{l-2}+4\omega_{l+1}$.

For the case $j=n-4$, we proceed as above, handling the cases $l=3,4$ with a Magma computation.
We note that the $X$-module $\delta\otimes \wedge^5(\delta^*)$ has summands 
$(\omega_1+\omega_{l-2}+5\omega_{l+1})$ and $(\omega_1+\omega_{l-3}+4\omega_{l+1})$.
Neither of these summands occurs in $\wedge^4(\delta^*)$.
Tensoring each of these with $\omega_{l+1}$ produces a summand 
$\omega_1+\omega_{l-3}+5\omega_{l+1}$.

Finally, we turn to the case $j=n-5$, handling the cases $l=3,4,5$ with Magma. 
This is entirely similar to the previous case; here we use the summands 
$(\omega_1+\omega_{l-3}+6\omega_{l+1})$ and 
$(\omega_1+\omega_{l-4}+5\omega_{l+1})$ of $\delta\otimes \wedge^6(\delta^*)$.\hal

In the proof of the next lemma, it will be useful
to adopt a notation used in \cite{Howe} for certain calculations within $\wedge^j(2\omega_1)$.
(This is quite similar to the notation used for analysing $\wedge^j(\o_2)$ described in the preamble to Lemma \ref{(s+2)cge2}.)
Let $X=A_{l+1}$ and let  $e_1,\dots,e_{l+2}$ be a basis of the natural $X$-module such that the weight 
of $e_1$ is $\omega_1$ and the weight of $e_i$ is $\omega_1-\alpha_1-\cdots-\alpha_{i-1}$, for $i\geq 2$. 
Then a basis of $V_X(2\omega_1)$ is given by
the symmetric tensors $e_i\otimes e_i$, $1\leq i\leq l+2$ and $e_i\otimes e_j-e_j\otimes e_i$, for $1\leq i<j\leq l+2$. We will note these by $e_{ij}$.

Then we arrange them in the following tableau:
$$e_{11}\quad e_{12}\quad e_{13}\quad e_{14}\quad\dots$$
$$\phantom{e_{11}}\quad e_{22}\quad e_{23}\quad e_{24}\quad\dots$$
$$\phantom{e_{11}}\quad \phantom{e_{12}}\quad e_{33}\quad e_{34}\quad\dots$$
The corresponding weights are as follows:
\[
\begin{array}{llll}
2\o_1 & 2\o_1-\a_1 & 2\o_1-\a_1-\a_2 & \cdots \\
   & 2\o_1-2\a_1 & 2\o_1-2\a_1-\a_2 & \cdots \\
   & & 2\o_1-2\a_1-2\a_2 & \cdots \\
   &&& \cdots 
\end{array}
\]
As before, we say a set $S\subseteq \{e_{rs}| 1\leq r\leq s\leq l+2\}$ is increasing if
all $e_{rs}$ above
or to the left of any element in $S$ also lie in $S$.
Then it is clear that for any increasing set $\{v_1,\dots,v_t\}\subseteq \{e_{rs}| 1\leq r\leq s\leq l+2\}$, the vector 
$v_1\wedge v_2\wedge\cdots\wedge v_t$ is a maximal vector.
 Hence these vectors provide irreducible $X$-summands in the tensor algebra 
$\wedge(V_X(2\omega_1))$. Moreover, \cite[Thm.4.4.2]{Howe} shows that all irreducible summands
occur in this fashion, and with multiplicity 1. We will use this often in what follows.

\begin{lem}\label{nonMF12} Assume $l\geq 5$ and let $7\leq j\leq l+8$. Then $V_Y(\lambda_j+\lambda_n) \downarrow X$ is not MF.
\end{lem}

\pf Note that $V_Y(\lambda_j)\otimes V_Y(\lambda_n)=V_Y(\lambda_j+\lambda_n)\oplus V_Y(\lambda_{j-1})$ 
and since the restriction of the second summand to $X$ is MF by Theorem \ref{wedgec}, 
it suffices to show that the restriction of the tensor product has a multiplicity three $X$-summand. This can be checked using Magma for $l=5,6$, so assume that $l\ge 7$. 

As in the proof of Lemma \ref{nolab4}, for each value of $j$ we produce three sequences which yield summands of $\wedge^j(2\o_1)$ which have have common composition factor when tensored with $2\o_{l+1}$. The details are recorded in Table \ref{howetab}. This completes the proof. \hal
\begin{table}[h!]
\caption{}\label{howetab}
\[
\begin{array}{|l|l|l|l|}
\hline
j & \hbox{sequences} & \hbox{weight} &  \hbox{common comp. factor} \\
\hline
\le l+3 & (j-2,2) & (j-5)\omega_1+2\omega_2+\omega_3+\omega_{j-2} & (j-6)\omega_1+2\omega_2+\omega_3+\omega_{j-3} \\
           & (j-3,2,1)  & (j-6)\omega_1+3\omega_3+\omega_{j-3} & \\
           &  (j-3,3) & (j-7)\omega_1+3\omega_2+\omega_4+\omega_{j-3} & \\
\hline
l+4 & (l+2,2) & (l-1)\omega_1+2\omega_2+\omega_3 &  (l-2)\omega_1+2\omega_2+\omega_3+\omega_{l+1} \\
          & (l+1,3) & (l-3)\omega_1+3\omega_2+\omega_4+\omega_{l+1} & \\
          & (l+1,2,1) & (l-2)\omega_1+3\omega_3+\omega_{l+1} & \\
\hline
l+2+m, & (l+2,m) & (l-m+1)\omega_1+m\omega_2+\omega_{m+1} &  (l-m)\omega_1+m\omega_2+\omega_{m+1}+\omega_{l+1} \\
3\le m\le 6  & (l+1,m+1) & (l-m-1)\omega_1+(m+1)\omega_2+ & \\
&& \omega_{m+2}+\omega_{l+1} & \\
        & (l+1,m,1) & (l-m)\omega_1+(m-2)\omega_2+2\omega_3+ & \\
&& \omega_{m+1}+\omega_{l+1} & \\
\hline
\end{array}
\]
\end{table}

\begin{lem}\label{a1000}  With notation as in $(\ref{2o1not})$, let 
  $\l = a\l_1+\l_2$ with  $a \ge 4$, and $V = V_Y(\l)$. Then $V\downarrow X$ is not MF.
   \end{lem}

\pf Suppose the assertion is false. Using induction it will suffice to obtain a contradiction when   $X = A_2$.   Then
$Y = A_5$ and $L_X'$ is embedded in $L_Y' = C^0 \times C^1 = A_2 \times A_1.$   Write $\pi(Y) = \{\b_1, \ldots , \b_5\}$ with  $\pi(C^0) = \{\b_1, \b_2\}$  and  $\pi(C^1) = \{ \b_4 \}.$  Take the base of the module with highest weight $\d = (20)$ to be $\{\d, \d-\a_1, \d-2\a_1, \d-\a_1-\a_2, \d-2\a_1-\a_2, \d-2\a_1-2\a_2 \}$, corresponding to the weights 
$\l_1,\l_1-\b_1,\ldots, \l_1-\b_1-\cdots -\b_5$. 
 We then find that $\b_1,  \b_2$ and $ \b_4$ all restrict to $T_X$ as $\a_1$, whereas $\b_3 \downarrow T_X = \a_2-\a_1$ and $\b_5 \downarrow T_X = \a_2.$   Therefore, $\l_1 \downarrow T_X = (2,0)$ and $\l_2 \downarrow T_X = (2,1)$, so that  $\l \downarrow T_X = (2a+2,1).$

We will use level analysis to obtain the result. 
The contradiction will come from the multiplicity of the irreducible $(2a-6)$ in $V^3.$

We begin with the top level  $V^1(Q_Y) = (a,1)\otimes (0).$  Now $(a,1) = ((a,0) \otimes (0,1)) - (a-1,0)$  and restricting to $L_X'$ 
using Theorem \ref{SOpowers}, we obtain
$(2a+ (2a-4) \oplus   \cdots ) \otimes 2) - ((2a-2) \oplus  (2a-6) \cdots).$   Expanding we have 
$$V^1 = (2a+2) \oplus  (2a) \oplus  (2a-2) \oplus  (2a-4) \oplus  (2a-6) \oplus  \cdots$$
Each of the summands in the above expression corresponds to level 1 of a certain irreducible 
$X$-module appearing as a composition factor of $V \downarrow X$.  
Since $\l \not \in [V,Q_Y]$, and $\l \downarrow T_X = (2a+2, 1)$, one such composition factor is $\xi_1 = (2a+2,1).$  
Write the others as  $\xi_2 = (2a,x_2), \xi_3 = (2a-2,x_3), \ldots .$   Let $T(c) = h_1(c)h_2(c^2),$ so that
the elements $T(c)$ generate the center of $L_Y$.  Then $T(c)$ induces scalars on $V^1(Q_Y)$ and it induces
$c^{(2a+2)+2} = c^{2a+4}$ on the irreducible with highest weight $(2a+2, 1)$.  Using this we find $x_i = i$
for $i = 2,3,\ldots .$  For example, for $x_3$ we must have $2a+4 = 2a-2 + 2x_3$, so that $x_3 = 3.$

Consider which of the $\xi_i$ can contribute a term $(2a-6)$ in $V^3.$  First
note that these will occur in $V^3(\xi_i)$ and  all irreducibles at the level, in particular $(2a-6)$, have
multiplicity 1. Let $v_{\xi_i}$ denote a maximal vector of $\xi_i$ The irreducibles $\xi_4  = (2a-4, 4),$ 
$\xi_5  = (2a-6, 5),$  and $\xi_6  = (2a-8, 6)$  each contribute  a factor $(2a-6).$  
Indeed, modulo other composition factors at the level, these are afforded by $f_{12}^2v_{\xi_4},$ $f_{12}f_2v_{\xi_5},$ $f_2^2v_{\xi_6}$, 
respectively.  No other $\xi_i$ can contribute a term $(2a-6)$. Therefore there are at most 3 such terms arising from $V^1(Q_Y).$

 Of course $(2a-6)$ may also arise in $V^3$ from irreducible summands $(2a-5,x)$ and $(2a-7,y)$ of $V \downarrow X$ which do not contribute a term at level 1 but do contribute a term $(2a-5)$ or $(2a-7)$, respectively, in $V^2$.  As we are assuming that $V \downarrow X$ is MF, 
Proposition \ref{induct} implies that there can exist at most 1 of each.  Therefore, in view of the MF assumption, 
the total multiplicity of $(2a-6)$ in $V^3$ is at most $3 + 2 + 1 = 6.$

Now consider $V^3(Q_Y).$ Listed below are weights that are highest weights of $L_Y'$-composition factors. 

\vspace{2mm}
(1) $\l - \b_2-\b_3-\b_4-\b_5 = (a+1,0) \otimes (0)$

(2) $\l - \b_1 - \b_2-\b_3 -\b_4-\b_5 = (a-1,1) \otimes (0)$

(3) $\l -2\b_1- 2\b_2- 2\b_3 = (a-2,1) \otimes (2)$

(4) $\l -\b_1- 2\b_2- 2\b_3 = (a,0) \otimes (2)$

\vspace{2mm}

\noindent Now  (1) and (2) sum to $((a,0) \otimes (1,0)) \otimes (0)$, and (3) and (4) sum to $((a-1,0)\otimes (1,0))\otimes (2)$.
Using this information  
we list the retrictions to $L_X'$ of the above composition factors:

\vspace{2mm}
(1) + (2):  $(2a \oplus  (2a-4) \oplus  (2a-8) \cdots) \otimes (2) \supseteq (2a-6)^2$

(3) + (4) :  $((2a-4) \oplus  (2a-6))\otimes 2\otimes 2 \supseteq (2a-6)^5$. 

\vspace{2mm}

\noindent We conclude that $(2a-6)$ appears with multiplicity at least 7 in $V^3$,
which contradicts the above.  \hal

In the final lemma of this subsection we change notation to $X = A_{l+1} < Y = SL(W)=A_n$, with $W = V_X(r\o_1)$, $r\ge 3$.

\begin{lem}\label{nonMF2} Let $l\geq 1$ and $r\geq 3$. Then $V_Y(a\lambda_1+\lambda_{n-1})\downarrow X$ ($a=1,2$), 
and $V_Y(\lambda_2+\lambda_{n-1})\downarrow X$, are not MF.
\end{lem}

\pf We first note that 
\[
\begin{array}{l}
V_Y(a\lambda_1)\otimes V_Y(\lambda_{n-1}) = 
V_Y(a\lambda_1+\lambda_{n-1})\oplus V_Y((a-1)\lambda_1+\lambda_n), \\ 
V_Y(\lambda_2)\otimes V_Y(\lambda_{n-1}) = 
V_Y(\lambda_2+\lambda_{n-1})\oplus( V_Y(\lambda_1)\otimes V_Y(\lambda_n)).
\end{array}
\]
Hence we must show that the following $X$-modules are not MF:
\[
\begin{array}{l}
((S^a(r\omega_1)\otimes \wedge^2(r\omega_{l+1}))\oplus \delta_{a,2}V_X(0))\;-\;
(S^{a-1}(r\omega_1)\otimes (r\omega_{l+1}))\;(a=1,2), \hbox{and } \\ 
(\wedge^2(r\omega_1)\otimes \wedge^2(r\omega_{l+1}))\;-\;(r\omega_1\otimes r\omega_{l+1}).
\end{array}
\] 
Note that $S^2(r\omega_1)$ has summands $2r\omega_1$ and $((2r-4)\omega_1+2\omega_2)$, 
 and $\wedge^2(r\omega_1)$ has summands $((2r-2)\omega_1+\omega_2)$ and 
$((2r-6)\omega_1+3\omega_2)$.
Then using Lemma \ref{aibj12}, we check that $((ar-2)\omega_1+\omega_l+(2r-4)\omega_{l+1})$ occurs 
with multiplicity at least 2 in $V_Y(a\lambda_1+\lambda_{n-1})\downarrow X$ for $a=1,2$, and 
$((2r-4)\omega_1+\omega_2+\omega_l+(2r-4)\omega_{l+1})$ with multiplicity at least 2 in 
$V_Y(\lambda_2+\lambda_{n-1})\downarrow X$. This completes the proof of the lemma.\hal

 \subsection{Non-MF symmetric and wedge squares}\label{nonmfs2}
This subsection consists of lemmas showing that $\wedge^2W$  and $S^2W$ are not MF for various $X$-modules $W$, where $X = A_{l+1}$. 

\begin{lem}\label{donna1}
Let $X = A_{l+1}$ with $l\ge 2$, and let $c,d>0$ and $e\ge 0$.

{\rm (i)} The modules $\wedge^2(V_X(\d))$ and $S^2(V_X(\d))$ are not MF for the following weights: 
\[
\begin{array}{ll}
\d = & \om_1+c\om_i+d\om_{l+1},\,c\om_1+\om_i+d\om_{l+1},\\
        & \om_i+c\om_{i+1}+d\om_{l+1}, \,c\om_i+\om_{i+1}+d\om_{l+1}.
\end{array}
\]

{\rm (ii)} The module $\wedge^2(V_X(\d))$ is not MF for the following weights:
\[
\begin{array}{ll}
\d = & 2\om_1+2\om_2+e\om_{l+1}\, (e \ge 1), \\  & 2\om_1+2\om_l+e\om_{l+1},\, ((l,e) \ne (2,0)). \\  \end{array} \]


\end{lem}

\pf This follows from   Lemma \ref{lemma1.5(i)(ii)} with the exception of the case  $\delta = 2\om_1 + 2\om_l$ ($e = 0$ in (ii)).  For this case we will use dominoes to show that $\mu = 2\om_1 + \om_{l-1} + \om_l + \om_{l+1}$ appears with multiplicity $2$ in $\wedge^2 ( \delta).$  Note that $\mu = 2\delta - (2\a_1 + \cdots + 2\a_{l-1} + 3\a_l + \a_{l+1}).$   The weight of a partition corresponding to $\delta$ is $1^42^2 \ldots  l^2.$  

According to Theorem \ref{dominoes} we double and repeat all exponents, getting $8844\ldots  44$, where there are $2(l-1)$ 4's. Then the array will have a total of $4l+4$ tiles.  Consequently we consider the partition of weight $1^62^4 \ldots  (l-1)^4l^3(l+1)^2(l+2)^1$ which corresponds to $\mu$ and exponents summing to $4l+4.$

We will indicate two labellings giving this weight.  In each case the base has 4 vertical tiles followed by a pair of horizontal 2-tiles lying over horizontal 1-tiles.  We indicate two options for the remaining part of the tiling.  

Just above the base in  one tiling  there is a sequence of  4 vertical dominoes with labels $(2,2,3,3), \ldots$, 
$(l-1,l-1,l,l)$.  Above this there is a vertical $l$-domino, followed by a vertical $(l+1)$-domino, followed by a horizontal $(l+2)$-domino above a horizontal $(l+1)$-domino.  
For the other tiling just above the base there are 2 horizontal 2-dominoes.
This is followed by a sequence of 4 vertical dominoes with labels $(3,3,3,3), \ldots , (l-1,l-1,l-1,l-1).$ At the next level there are 2 vertical $l$-dominoes followed by a horizontal $(l+1)$-domino above a horizontal $l$-domino.  And above this there is  a horizontal $(l+1)$-domino followed by a horizontal $(l+2)$-domino.  

For these tilings there are $6$ (respectively  $10$) horizontal dominoes,  so both summands are alternating. \hal

\begin{lem}\label{0c0...010} Let  $X = A_{l+1}$  with $l\ge 4$ and let $\d = c\om_2 + \om_l.$ 

{\rm (i)}  If $c  > 1$, then  $\wedge^2(\d)$ is not MF.

{\rm (ii)}  $S^2(\d)$ is not MF.

\end{lem}

\pf  We verify this using the domino technique.  First assume $c > 1.$  For $\d$ we use the partition with weight  $1^{c+1}2^{c+1}3^1 \ldots  l^1$.  Doubling and repeating exponents we have $2c+2, 2c+2, 2c+2, 2c+2, 2, \ldots , 2$ where  2 occurs $2(l-2)$ times.  So the bottom part of the array is a $4 \times (2c+2)$ matrix on top of which is a $(2l-4) \times 2$ matrix. 

 We claim that the composition factor of highest weight $\mu = (1(2c-3)10 \ldots  010)$ appears with multiplicity at least 2 in each of $\wedge^2(\d)$ and $S^2(\d).$  We will use the partition with weight $1^{2c+1}2^{2c}3^34^2 \ldots  l^2(l+1)^1(l+2)^1$ to correspond to $\mu.$

We first describe two tilings for  $\wedge^2(\d)$    yielding the above partition.  In both cases the base has $2c$ vertical dominoes followed by a horizontal 2-domino above a horizontal 1-domino and the next layer begins with $2c-2$ vertical 2-dominoes, then a horizontal 3-domino above a horizontal 2-domino,  and this is followed by a horizontal 4-domino lying over a horizontal 3-domino.  

Consider the $(2l-4) \times 2$ matrix. In one case  the tiling is a sequence of vertical dominoes with labels  $(3,5), (4,6), \ldots , (l,l+2).$ In the other case the tiling begins with a horizontal 4-domino above a horizontal 3-domino and  is then followed by pairs of vertical dominoes with labels $(5,5), \ldots , (l,l)$, which in turn is followed by horizontal dominoes with labels $l+1$ and $l+2$.  These tilings have $6$ and $10$ horizontal dominoes, respectively, so they both correspond to alternating summands.

Now consider $S^2(\d),$ where we will again produce two tilings.  In both cases the bottom four rows are as above.  Consider the tiling of the $(2l-4) \times 2$ matrix.
One tiling begins with a single horizontal 3-domino which is followed by a sequence of  pairs of vertical dominoes with labels $(4,5), (5,6), \ldots , (l,l+1)$ and then a horizontal $(l+2)$-domino.
The second labelling has a sequence of pairs of vertical dominoes with labels $(3,4), (5,5), \ldots , (l,l)$
 and this is followed by horizontal dominoes with labels $l+1$ and $l+2.$   Both labellings have 8 horizontal dominoes so correspond to symmetric summands.
 
 Finally we must show that $S^2(\d)$ is not MF when $c = 1.$ We will show that  $S^2(\d) \supseteq (\d)^2$   by describing two symmetric tilings with weight $1^32^33^2 \ldots  l^2(l+1)^1(l+2).$  The resulting partition corresponds to $\d$ so this will give the
 assertion.
 
 In each case the base is labelled with two vertical 1-dominoes followed by a horizontal 2-domino above a 
horizontal 1-domino.  The next two rows are tiled with two vertical 2-dominoes  followed by a horizontal 4-domino above a 
horizontal 3-domino.  Now consider the labelling of the $(2l-4) \times 2$ matrix.

 In the first tiling this part of the labelling begins with a horizontal 4-domino above a horizontal 3-domino. This is followed by a sequence of pairs of vertical dominoes with labels $(5,5), \ldots , (l,l)$, which in turn is followed by horizontal dominoes with labels $l+1$ and $l+2.$  The second tiling has  a sequence of pairs of vertical dominoes with labels $(3,5), (4,6), \ldots , (l+1,l+2).$ 
 Both labellings correspond to symmetric summands, so this completes the proof.   \hal

 \begin{lem}\label{010...0100} Let  $X = A_{l+1}$  with $l \ge 6$ and let $\d = \om_2 + \om_{l-1}.$
  Then $\wedge^2(V_X(\d))$ is not MF.
  \end{lem}

  \pf  
We  claim that $\mu = \om_2 + \om_{l-2}$ appears in $\wedge^2(\d)$ with multiplicity 2.    We again apply the domino technique.
  Starting with the partition  of weight $1^22^23^1\ldots  (l-1)^1$ we double and repeat exponents to  obtain
  the sequence $4,4,4,4,2, \ldots , 2$ where 2 appears $2(l-3)$ times.  We will use the partition with weight $1^32^33^2\ldots (l-2)^2(l-1)^1l^1(l+1)^1(l+2)^1$ to correspond to $\mu.$
  
  We will indicate two tilings of the array which correspond to the above weight.  In each
  case the base consists of two vertical 1-dominoes followed by a horizontal 2-domino above a horizontal 1-domino.  And the layer above the base has  two vertical 2-dominoes followed by a horizontal 4-domino
  above a horizontal 3-domino.  We now describe the tilings of the $(2l-6)\times 2$ array.

 The first tiling of the $(2l-6)\times 2$ array has a sequence of vertical dominoes with labels $(3,5), \ldots , (l-2,l)$ followed by a horizontal $l+2$-domino lying above a horizontal $l+1$-domino.  This tiling has 6 horizontal dominoes and hence is alternating.

For the second possibility the tiling begins with a horizontal 4-domino above a horizontal 3-domino
  and is followed by a sequence of pairs of vertical dominoes with labels $(5,5), \ldots , (l-2,l-2)$  (omit
  this sequence if $l = 6$) and this followed by horizontal dominoes with labels $l-1, l, l+1, l+2.$  The second
  tiling has 10 horizontal dominoes and so it also yields an alternating composition factor.  \hal 
 
 \begin{lem}\label{S2(01010...0)}  Let  $X = A_{l+1}$  with $l \ge 4$ and let $\d = \om_2 + \om_4.$
  Then $S^2(V_X(\d))$ is  not MF.
    \end{lem}
  
  \pf  This follows from a Magma computation for $l = 4$ and then an induction argument using a parabolic subgroup gives  
the result for  $l > 4.$   \hal

\begin{lem}\label{om3oml}  Assume $X = A_{l+1}$ and let $\d = \om_3 + \om_l$.
\begin{itemize}
\item[{\rm (i)}] $\wedge^2(\d)$ is not MF  if $l \ge 6.$
\item[{\rm (ii)}]  $S^2(\d)$ is not MF if $l \ge 5.$
\end{itemize}
 \end{lem}

\pf  We will apply the domino technique. Starting with the weight  $1^22^23^24^1 \ldots  l^1$ we double and repeat exponents getting $4,4,4,4,4,4,2,2 \ldots , 2,2$ where 2 appears $2(l-3)$ times.  The tiling will therefore involve $2l+6$ tiles.

(i) First consider $\wedge^2(\d)$ where we will show that
$\om_4 + \om_l$ occurs with multiplicity 2.  In order to get the proper number of tiles we use a partition with weight $1^32^33^34^35^2 \ldots  l^2(l+1)^1(l+2)^1.$ The base of the tiling has two vertical 1-dominoes followed by a horizontal 2-domino above a horizontal 1-domino.  Columns 3 and 4 consist entirely of horizontal dominoes with labels $1,\ldots , 6.$  One tiling has columns 1 and 2 beginning with pairs of vertical dominoes with labels $(1,1), (2,2), (3,3), (4,4)$ followed by pairs of vertical dominoes with labels 
 $(5,7), \ldots , (l,l+2)$. The first two columns of the second tiling begin with pairs of vertical dominoes with labels $(1,1), (2,2), (3,3), (4,4)$ as before and  these are followed by a horizontal 6-domino lying above a horizontal 5-domino.  If $l \ge 7,$ this is followed by a sequence of vertical dominoes with labels $(7,7), \ldots , (l,l)$ and this sequence (if it occurs) is followed by a horizontal $(l+2)$-domino lying above an $(l+1)$-domino.  These tilings have 6 and 10 horizontal dominoes, respectively, so both correspond to alternating summands.

(ii) For $S^2(\d)$ we claim that $\om_1 + \om_3 + \om_l$ occurs with multiplicity 2.
This time we use a partition with weight $1^42^33^34^2 \ldots  l^2(l+1)^1(l+2)^1.$ Both tilings have base consisting of 4 vertical 1-dominoes.  Above the vertical 1-dominoes in columns
3 and 4 there are horizontal dominoes with labels 2,3,4,5.  In one tiling columns 1 and 2 begin with pairs of vertical dominoes with labels $(1,1), (2,2), (3,3)$ followed by vertical dominoes with labels 
$(4,6), (5,7), \ldots , (l,l+2)$. The first two columns of the second tiling also begin with vertical dominoes with labels $(1,1),(2,2), (3,3)$ and these are followed by a horizontal $5$-domino lying above a horizontal $4$-domino.
Then comes a sequence of pairs of vertical dominoes with labels $(6,6), \ldots , (l,l)$ and the sequence is followed by a horizontal $(l+2)$-domino lying over a horizontal $(l+1)$-domino.  So this time the tilings have 4 or 8 horizontal dominoes and hence correspond to symmetric summands. \hal

 \subsection{Low rank cases}\label{nonmflowrk}

This subsection consists of some non-MF results for groups of low rank.

\begin{lem}\label{newprop1}
{\rm (i)} Let $X=A_4$. Then for any $c\ge 2$, the following $X$-modules are not MF:
\[
\begin{array}{l}
S^c(\om_2)\otimes \om_2, \\
S^c(\om_2)\otimes \om_3, \\
S^c(\om_2)\otimes  S^2(\om_2),  \\
S^c(\om_2)\otimes S^2(\om_3).
\end{array}
\]

{\rm (ii)} Let $X = A_5$. For $c\ge 6$, the $X$-module $S^c(\om_3)$ is not MF.  
\end{lem}

 \pf  The last two cases of (i) are  easy and we settle these first. First note that $ S^c(\om_2) \supseteq 0c00 + 0(c-2)01.$  
Similarly $S^2(\om_2) = 0200 + 0001$ and  $S^2(\om_3) = 0020 + 1000.$  It is then clear  that  
$S^c(\om_2)\otimes  S^2(\om_2) \supseteq (0c01)^2.$ Also an application of Lemma \ref{abtimesce} shows  that  
$S^c(\om_2)\otimes S^2(\om_3) \supseteq (0(c-1)10)^2$.
  
 This leaves the first two cases of (i) to deal with, and also (ii).  
These are all   handled by the same argument.  We claim that $S^c(\om_2)\otimes \om_2 \supseteq (0(c-1)01)^2$, 
$S^c(\om_2)\otimes \om_3 \supseteq (0(c-1)00)^2$ (first two cases of (i), resp.) and
 $S^c(\om_3) \supseteq (00(c-4)00)^2$ (case (ii)).  Note that if $\g$ denotes the highest weight of the module on the left hand side of each of 
these expressions, the highest weights of the irreducible modules indicated on the right hand sides are 
$\g -1210$, $\g -1221$ and $\g -24642,$ respectively.
 
We illustrate the argument with (ii). Using Magma we verify that  $S^6(\om_3) \supseteq (00200)^2,$ 
so that the claim holds for $c = 6.$  Now consider $S^c(\om_3)$ for $c > 6.$  This has a basis of symmetric products
of $c$ terms of the form $\om_3 - \a,$ where 
either $\a = 0$ or $\a$ is a positive root with nonzero coefficient of $\a_3.$  Therefore there can be at most $6$ of the latter.  
It follows that the number of ways of achieving 
$\g -24642$ is precisely the same as for the $c = 6$ case. Similarly for all dominant weights above $\g -24642$.  
The claim follows for this case and a similar argument works for the remaining cases of (i), where it suffices to consider the case $c=2$. \hal

\begin{lem}\label{donnathreeseven}  Assume $X = A_{l+1}$, $W = V_X(\d)$ and $Y = SL(W) = A_n$, and let $a \ge 2.$
\begin{itemize}
\item[{\rm (i)}] If $\delta = a\omega_1+2\omega_{l+1}$, then $V_Y(\lambda_1+\lambda_{n-1})\downarrow X$ is not MF.
\item[{\rm (ii)}] If $l=2$ and $\delta = a\omega_1+2\omega_3$, then $V_Y(\lambda_1+\lambda_{n-3})\downarrow X$ is not MF.
\item[{\rm (iii)}] If $l=2$ and $\delta = a\omega_1+2\omega_3$, then $V_Y(\lambda_1+\lambda_{n-4})\downarrow X$ is not MF.
\end{itemize}
\end{lem}


\pf (i) Let $V= V_Y(\lambda_1+\lambda_{n-1})$. 
Then $V$ is isomorphic to the quotient 
$(V_Y(\lambda_1)\otimes V_Y(\lambda_{n-1}))/V_Y(\lambda_n)$. So we exhibit a repeated summand of $V_X(\delta)\otimes 
\wedge^2(V_X(\delta^*)$ which is not the irreducible $V_X(\delta^*)$.

Note that $\wedge^2(\delta^*)$ has summands 
$(4\omega_1+\omega_l+(2a-2)\omega_{l+1})$ and 
$(2\omega_1+\omega_2+2a\omega_{l+1})$.
Tensoring with $\delta$, we obtain a repeated summand 
$((a+2)\omega_1+\omega_2+\omega_l+2a\omega_{l+1})$, as desired.

(ii) We proceed as in (i). Here $V=V_Y(\lambda_1+\lambda_{n-3})$ is isomorphic to the 
quotient $(V_Y(\lambda_1)\otimes V_Y(\lambda_{n-3}))/V_Y(\lambda_{n-2})$. 
Consider the tensor product $\delta\otimes \wedge^4\delta^*$; the second tensor factor has summands 
$(5\omega_1+\omega_2+(4a-1)\omega_3)$ and $(7\omega_1+\omega_2+(4a-3)\omega_3)$. Then tensoring each of these with 
$\delta$, we obtain a repeated summand $((a+5)\omega_1+2\omega_2+(4a-1)\omega_3)$. Now
an $S$-value comparison shows that the weight $(a+5)\omega_1+2\omega_2+(4a-1)\omega_3$ is not subdominant to $3\delta^*$ 
and so the repeated summand cannot lie in $\wedge^3(\delta^*)$.

(iii) We proceed as in the previous proof. Here $V=V_Y(\lambda_1+\lambda_{n-4})$ is isomorphic to the 
quotient $(V_Y(\lambda_1)\otimes V_Y(\lambda_{n-4}))/V_Y(\lambda_{n-3})$. 
We find summands $(5\omega_1+3\omega_2+(5a-3)\omega_3)$ and 
$(7\omega_1+3\omega_2+(5a-5)\omega_3)$ in  $\wedge^5(\delta^*)$.
 Then tensoring each of these with $\delta$, we obtain
a repeated summand $((a+5)\omega_1+4\omega_2+(5a-3)\omega_3)$. An $S$-value comparison shows that the weight 
$(a+5)\omega_1+4\omega_2+(5a-3)\omega_3$ is not subdominant to $4\delta^*$ and so the repeated summand cannot
lie in $\wedge^4(\delta^*)$.\hal

 \begin{lem}\label{lemma1.5(iii)(iv)} Let $X = A_{l+1}$ and $c\ge 0$.
 \begin{itemize}
\item[{\rm (i)}]  If $l =1$ and $c > 1$, then the $X$-module  $\wedge^3(c\omega_1 + \omega_2)$ is not MF.
\item[{\rm (ii)}] For $l = 2$ the module $S^3(\omega_1+ \omega_2 + c\omega_3)$ is not MF.
\item[{\rm (iii)}]  For $l = 2$ and $c > 0$, the module $\wedge^3(c\omega_1 + \omega_2)$ is not MF.
\end{itemize}
   \end{lem}

   
   \pf (i) Assume $\d = c\om_1 + \om_2$ and consider $V = \wedge^3(\d).$ The weight $3\d-2\a_1-\a_2$ appears with multiplicity 3 in $\wedge^3(\d)$. However, the only dominant weight $\nu$ with $3\d-2\a_1-\a_2<\nu \le 3\d$ giving rise to a summand of $\wedge^3(\d)$ is $3\d-\a_1-\a_2$, and the weight $3\d-2\a_1-\a_2$ occurs with multiplicity 1 in this summand. Part (i) follows.

 (ii) Let $\d = \omega_1 + \omega_2 + c\omega_3$ and $\a_0 = \alpha_1+\alpha_2 +\alpha_3.$  
If $c = 0$ we use Magma to obtain the result, so assume $c \ne 0.$
 The weight  $3\d-\a_0$ can be obtained in $S^3(\d)$ as the sum of the following $3$-tuples: 
$(\d,\d-\alpha_2,\d-\alpha_1-\alpha_3),   (\d-\alpha_1,\d-\alpha_2,\d-\alpha_3), 
(\d,\d-\alpha_1,\d-\alpha_2-\alpha_3),  (\d,\d-\alpha_1-\alpha_2,\d-\alpha_3),$ and 
$(\d,\d,\d-\alpha_1-\alpha_2-\alpha_3).$ These occur with multiplicities $1,1,2, 2,4$ for a total of $10.$  
 
 The irreducible summands
 of $S^3(\d)$ with highest weights strictly above $3\d-\a_0$ have highest weights
 $3\d, 3\d-\alpha_1-\alpha_2, 3\d-\alpha_2-\alpha_3,$ and $3\d-\alpha_1-\alpha_3.$  
These irreducibles each occur with multiplicity $1$ and dimensions of the weight spaces for weight 
$3\d-\a_0$ are $4,1,1,1,$ respectively.  Therefore, $S^3(\d) \supseteq (3\d-\a_0)^2.$
 
     (iii)  Let $\d = c\omega_1 + \omega_2.$ If $c =1$ the result follows from (ii).  
And  for $c > 1$ the result follows from (i) since $V^1$ is not MF.   \hal

\begin{lem}\label{alambda2} Let $X = A_4$, $W = V_X(\o_2)$ and $X < Y = SL(W) = A_{9}$.
If $\l = a\l_2$ with $a \ge 6,$ then $V_Y(\l)\downarrow X$ is not MF.
\end{lem}

\pf We will show that  $V_Y(\l)\downarrow X \supseteq ((a-5)3(a-5)1)^2.$ We have $C^0 = A_5,$ $C^0 = A_3$, and $V^1 = (0a000)\downarrow L_X'$.  

First consider $V^2(Q_Y)$.  It follows from Theorem \ref{LEVELS} that  $V^2(Q_Y)$ is irreducible and
$V^2 = ((1(a-1)000)\downarrow L_X') \otimes (100).$ We will consider the multiplicity of $ ((a-5)3(a-5))$ in $V^3$.  The composition factors in $V^2$ that can give rise to  $ ((a-5)3(a-5))$ are as follows:

\vspace{2mm}

 $a.$  $ ((a-5)3(a-6))$
 
 $b.$  $ ((a-5)2(a-4))$
 
  $c.$ $ ((a-6)4(a-5))$ 
  
   $d.$  $ ((a-4)3(a-5))$
   
   \vspace{2mm}

It follows from the discussion  in Subsection \ref{a3a5}  that  $(1(a-1)000)\downarrow L_X'$ has all of its composition factors being self-dual. Indeed, the partition corresponding to  $(1(a-1)000)$ is $1^a2^{a-1}$ so the array will have just 2 rows and the Y-condition implies  that  for each composition factor in the restriction, the partition has the form $1^{b+c}2^b$ and the composition factor is $(b,c,b).$

The self-dual irreducibles  of $L_X'$ for which tensoring with $(100)$ produces one of the irreducibles above are  listed below along with the irreducibles that they produce.

\vspace{2mm}

$1. \  ((a-4)3(a-4)): ((a-4)3(a-5)).$ 

$2. \  ((a-4)1(a-4)): ((a-5)2(a-4)).$

$3. \  ((a-5)3(a-5)):  ((a-5)3(a-6)), ((a-4)3(a-5)), ((a-5)2(a-4)), ((a-6)4(a-5))$. 

 $4. \  ((a-6)3(a-6): ((a-5)3(a-6)).$

$5. \  ((a-6)5(a-6)):  ((a-6)4(a-5)).$ 

 \vspace{2mm}

We next show how these 5 composition factors  occur in $(1(a-1)000)\downarrow L_X'$ using the method described in Subsection \ref{a3a5}.  The partition corresponding to $(1(a-1)000)$ is $1^a2^{a-1}$ so the arrays will have two rows with lengths $a$ and $a-1.$ All 2's must occur in the second row and each 1 appearing in the first row must have a 2 below it in the second row. Using these facts one shows  that in each case there is only one even partition giving rise to the composition factor.
For each case we indicate below the composition factor,  a corresponding partition,   the only possible even partition yielding the partition and the rows of the resulting array.  

\vspace{2mm}

$1.$ $((a-4)3(a-4))$, $1^{a-1}2^{a-4}$:

\ \ \ \  $(4,0)$:  $(xxxx1^{a-4}),$  $(1^32^{a-4})$ 

$2.$  $((a-4)1(a-4)),$  $1^{a-3}2^{a-4}:$

\ \ \ \ $(4,2):$  $(xxxx1^{a-4}),$  $(xx1^12^{a-4})$

$3.$  $((a-5)3(a-5)),$ $1^{a-2}2^{a-5}:$ 

\ \ \ \ $(4,2):$ $(xxxx1^{a-4})$, $(xx1^22^{a-5})$

$4.$ $((a-6)3(a-6)),$ $1^{a-3}2^{a-6}:$ 

\ \ \ \ $(6,2):$  $(xxxxxx1^{a-6}), $  $(xx1^32^{a-6})$

$5.$  $((a-6)5(a-6)),$   $1^{a-1}2^{a-6}:$ 

\ \ \ \ $(6,0):$ $(xxxxxx1^{a-6}),$ $(1^52^{a-6})$

\vspace{2mm}
It follows that each of the 5 summands appears with multiplicity 1 and so $((a-5)3(a-5))^8$ can potentially arise from  $V^2.$  We claim that in fact only $((a-5)3(a-5))^6$ can possibly arise in $V^3$.  To see this consider $V^1 = (0a000) \downarrow L_X'$.  Arguing as in Lemma \ref{0a001} we see that $V^1$ contains  $(a0a),$  $((a-4)2(a-4))$, and $(((a-6)4(a-6)). $  Then  $V_Y(a\l_2)\downarrow X$  contains irreducibles $ V_1,$ $V_2$, and and $V_3$ such that  $V_1 = (a0a0)$, $ V_2  = ((a-4)2(a-4)x)$, and $V_3 = ((a-6)4(a-6)y).$ Using the method of Lemma \ref{a4(a10...} find that  $x=2$ and $y=4.$ Note that $V_2^2 \supseteq ((a-5)2(a-4)) +((a-4)3(a-5))$ and these summands are (b) and (d) in the list above. Similarly $V_3^2 \supseteq ((a-6)4(a-5)) +((a-5)3(a-6))$ and these are (c) and (a). All four of these  give a summand  $((a-5)3(a-5))$ in $V^3.$  But since $V_2^3$  and $V_3^3$ are MF the claim follows.
 
 We now turn to $V^3$ which contains the summands
 $(2(a-2)000) \otimes (200)$ and $(0(a-1)000) \otimes (010).$  We will show that $((2(a-2)000) \otimes (200))\downarrow L_X' \supseteq ((a-5)3(a-5))^6$ and $(0(a-1)000) \otimes (010) \downarrow L_X'  \supseteq ((a-5)3(a-5))^2.$ 
 
 Consider $(2(a-2)000) \downarrow L_X'.$  We claim that this  contains $((a-6)4(a-6))^2,$ $((a-5)2(a-5))^1$, $((a-5)4(a-5))^1$,  and $((a-4)2(a-4))^2.$ Each of these tensors with $(200)$  to yield copies of $((a-5)3(a-5))$, so the claim will show that
  $((a-5)3(a-5))^6$ appears in the restriction of $(2(a-2)000) \otimes (200)$.
 
 In the following we list the partitions giving rise to these irreducibles, the relevant even partitions, and the labelling of the rows.  In each case there are two rows of lengths $a$ and $a-2$, respectively.
 
 \vspace{2mm}
 
 $((a-6)4(a-6)),   1^{a-2}2^{a-6}:$
 
 \ \ \ \   $(4,2):$ $(xxxx1^{a-4})$, $(xx1^22^{a-6})$ 
 
  \ \ \ \   $(6,0):$  $(xxxxxx1^{a-6})$, $(1^42^{a-6})$
  
 $((a-5)2(a-5)),   1^{a-3}2^{a-5}:$ 
 
 \ \ \ \  $(4,2):$ $(xxxx1^{a-4})$, $(xx1^12^{a-5})$ 
 
 $((a-5)4(a-5)),   1^{a-1}2^{a-5}:$ 
 
 \ \ \ \  $(4,0):$ $(xxxx1^{a-4})$, $(1^32^{a-5})$
 
 $((a-4)2(a-4)),   1^{a-2}2^{a-4}:$  

 \ \ \ \  $(4,0):$ $(xxxx1^{a-4})$, $(1^22^{a-4})$
 
  \ \ \ \  $(2,2):$ $(xx1^{a-2})$, $(xx2^{a-4})$
  
  \vspace{2mm}
  
  Similarly we claim that  $(0(a-1)000) \downarrow L_X'$ contains
  $((a-5)2(a-5))^1$  and $((a-5)4(a-5))^1.$ Each of these tensors with $(010)$ to yield $((a-5)3(a-5)).$ So together with the above
  this will produce $((a-5)3(a-5))^8.$ 
    
  As above we illustrate how the composition factors arise, noting
  that this time the arrays will have two equal rows of lengths $a-1$.
  
  \vspace{2mm}
  
  $((a-5)2(a-5)),   1^{a-3}2^{a-5}:$ 
 
 \ \ \ \  $(4,2):$ $(xxxx1^{a-5})$, $(xx1^22^{a-5})$ 

$((a-5)4(a-5)),   1^{a-1}2^{a-5}:$

 \ \ \ \  $(4,0):$ $(xxxx1^{a-5})$, $(1^42^{a-5})$
 
 \vspace{2mm}
 
 We have now shown that $((a-5)3(a-5))^8$ appears in $V^3$  and only $((a-5)3(a-5))^6$ can arise from $V^2$. It follows that  $V_Y(a\l_2)\downarrow X \supseteq ((a-5)3(a-5)x)^2$ and using the method of Lemma \ref{a4(a10...} we see that $x = 1.$ This completes the proof of the lemma.  \hal

\begin{lem}\label{00100}  Let $X = A_5$, $W = V_X(\o_3)$ and $X < Y = SL(W) = A_{19}$. Let $a \ge 2,\,c\ge 0$ and let $V = V_Y(\l)$, 
where $\l$ is one of the following weights:
\begin{itemize}
\item[{\rm (i)}]  $\l = a\l_1 + \l_2 + c\l_{10}$  or  $\l_8 + a\l_9+ c\l_{10}$
\item[{\rm (ii)}]  $\l = a\l_2  + c\l_{10}$ or $a\l_8+c\l_{10}$
 \item[{\rm (iii)}]  $\l = a\l_1 + \l_9 + c\l_{10}$ or $\l_1+a\l_9+c\l_{10}$.
\end{itemize} 
 Then $V \downarrow X$ is not MF.
\end{lem}

\pf  (i) By way of contradiction assume $V \downarrow X$ is MF.  We will consider levels in the usual way.  Here $C^0$ and $C^1$ are both of type $A_9$ and the embeddings
for $L_X'$ are via the representations $(0010)$ and $(0100),$ respectively.  Now $V^1 = 
V_{C^0}(a\l_1 + \l_2) \downarrow L_X'$ or its dual. 

First assume $V^1 = 
V_{C^0}(a\l_1^0 + \l_2^0) \downarrow L_X'.$  We will consider the multiplicity of the module $(10(a-2)0)$ in $V^2.$  
The only possible modules in $V^1$ that can give rise to this module are $(20(a-2)0), (01(a-2)0),$ and $(10(a-3)1).$ 
In view of our supposition that $V \downarrow X$ is MF, so is $V^1$, and so there are at most 3 of these. 
Therefore $(10(a-2)0)$ has multiplicity at most 4 in  $V^2.$

Now $V^2(Q_Y) $ contains the modules $V_{C^0}((a+1)\l_1^0) \otimes V_{C^1}(\l_1^1)$ and  $V_{C^0}((a-1)\l_1^0 + \l_2^0) \otimes  V_{C^1}(\l_1^1)$.  These add to $V_{C^0}(a\l_1^0) \otimes V_{C^0}(\l_1^0) \otimes 
V_{C^1}(\l_1^1).$  Restricting to $L_X'$ this becomes $S^a(0010) \otimes (0010) \otimes (0100)$.
Now $S^a(0010)$ contains $(10(a-2)0)$, by Lemma \ref{cfssc}.  So by Corollary \ref{LR_om2}, $S^a(0010) \otimes 0010 \otimes 0010$ contains $(10(a-2)0)^5$, a contradiction.

Now assume $V^1 = V_{C^0}(\l_8^0 + a\l_9^0) \downarrow L_X'.$   Then $V^2(Q_Y)$ contains the summands 
$V_{C^0}(2\l_8^0 + (a-1)\l_9^0) \otimes V_{C^1}(\l_1^1)$ and  $V_{C^0}(\l_7^0 + a\l_9^0) \otimes V_{C^1}(\l_1^1).$
Let $Z = (2\l_8^0 + (a-1)\l_9^0) \oplus  (\l_7^0 + a\l_9^0)$, so $V^2$ contains $(Z\downarrow L_X') \otimes \o_2$.

Using Theorem \ref{LR} we see that $Z = ((\l_8^0 + a\l_9^0) \otimes \l_9^0 -  (\l_8^0 + (a+1)\l_9^0)$ and 
$\l_8^0+b\l_9^0 = ((b+1)\l_9^0\otimes \l_9^0)-(b+2)\l_9^0$. Hence
\[
Z = ((a+1)\l_9^0\otimes \l_9^0\otimes \l_9^0) + (a+3)\l_9^0 - ((a+2)\l_9^0\otimes \l_9^0)^2.
\]
Now consider $Z\downarrow L_X'$. Using Corollary \ref{LR_om2}, we see that $S^{a+1}(\o_2)\otimes \o_2\otimes \o_2$ contains $(2,a,0,1) \oplus (2,a-1,2,0)$. Neither of these summands lies in $S^{a+2}(\o_2)\otimes \o_2$, since all composition factors of $S^{a+2}(\o_2)$ have $\o_1$-coefficient 0 (see \cite[3.8.1]{Howe}). Hence $Z\downarrow L_X'$ contains 
$(2,a,0,1) \oplus (2,a-1,2,0)$, and it follows that $V^2$ contains $((2,a,0,1) \oplus (2,a-1,2,0)) \otimes (0100)$. By Corollary \ref{LR_om2}, this contains $(3,a-1,1,1)^2$, with $S$-value $a+4$. However $V^1 = (\l_8^0+a\l_9^0)\downarrow L_X'
\subseteq \wedge^2(\o_2)\otimes S^a(\o_2)$ has $S$-value at most $a+2$. This is a contradiction, proving that $V\downarrow X$ is not MF.

(ii) Suppose $V \downarrow X$ is MF.  Then Lemma \ref{alambda2} shows that $a \le 5.$
We begin by showing that $c = 0$.  Suppose otherwise.   If $\l = a\l_2 + c\l_{10},$ then $V^2(Q_Y)$ contains $V_{C^0}(a\l_2^0)\otimes V_{C^0}(\l_9) \otimes V_{C^1}(\l_1^1)$ by Proposition \ref{v2gamma1}, and restricting to $L_X'$ we see that $V^2$ contains
$V^1\otimes (0100) \otimes (0100)$ which equals $V^1 \otimes ((0200)+(1010)+(0001)).$  Now we apply Lemmas \ref{twolabs} and \ref{cover} to obtain a contradiction.  A very similar argument works for $\l = a\l_8 + c\l_{10}.$ So from now on assume $c = 0.$

First consider the case $\l = a\l_8.$ Then $V^2(Q_Y) \supseteq V_{C^0}(\l_7^0 + (a-1)\l_8^0) \otimes V_{C^1}(\l_1^1).$ 
We claim that the restriction of the first tensor factor to $L_X'$ contains $((a+1)0(a-1)1)$ and $((a-1)1(a-1)1)$. This restriction is the same as that of the $A_9$-module $W = (0(a-1)10\ldots  0))$ to $D = A_4$ where the embedding is $(0100).$  Consider $W^1$
where we use the method of Subsection \ref{a3a5} to see that there are $L_D'$-composition factors of highest weights $((a+1)0(a-1))$ and $((a-1)1(a-1))$.  Indeed these arise from the even partitions $0 \ge 0 \ge 0$ and $2 \ge 0 \ge 0$, respectively.  Therefore  there exist $D$-composition factors of highest weights $((a+1)0(a-1)x),$  $((a-1)1(a-1)y),$ and $((a-1)0(a+1)z),$ where the existence of the last one follows from the fact that $W^1$ is self-dual.  Letting $T_D$ be the 1-dimensional torus centralizing $L_D'$ we see  $T_D$ acts on these composition factors via the weights $(a+1) +3(a-1) + 4x$,
$(a-1) + 2 + 3(a-1) + 4y$, and $(a-1) + 3(a+1) + 4z,$ respectively. As $T_D$ has the same action on all three factors, this implies that $x = y = z+1$. 

Now $W$ is contained
in $(0(a-1)0 \ldots  0) \otimes (0010\ldots  0)$ and 
$D$-composition factors of the tensor factors have $S$-value at most $2(a-1)$ and $3$, respectively. The $S$-value of the third factor above is $2a+z$, so that $z \le 1$. If $z = 1$, then $x = 2$ and the $S$-value of the first factor is $2a+2$, a contradiction.  Therefore  $z = 0,$ $x =y = 1$ and this establishes the claim.

Using Corollary \ref{LR_om2} one checks that each of
$((a+1)0(a-1)1) \otimes (0100)$ and $((a-1)0(a-1)1) \otimes (0100)$ contains $(a1(a-2)2)$, so that $V^2 \supseteq ((a1(a-2)2)^2.$ Finally, using the decomposition of $V^1$ given in the proof of Lemma \ref{a5,a6,a7} we see that this factor does not arise from a summand of $V^1$  and hence this is a contradiction.

The cases $\l = a\l_2$ is settled using Magma computations (recall that $a\le 5$), along the lines described in the proof of Lemma \ref{a5,a6,a7}:  we find that $V_{A_{19}}(a\l_2) = V^+-V^-$ with $V^+,V^-$ as in Table \ref{vplusminus}, with the exception that for $a=5$, $V^-$ has an extra term $\l_{10}^2$. Using this we restrict 
$V_{A_{19}}(a\l_2)$ to $X$ and find that this restriction fails to be MF.

(iii) First assume $\l = a\l_1 + \l_9 + c\l_{10}$ and assume, by way of contradiction, that $V \downarrow X$ is MF. We have $V^1 = V_{C^0}(a\l_1^0 + \l_9^0)\downarrow L_X'$ while $V^2$ contains $V_{C^0}((a-1)\l_1^0 + \l_9^0) \otimes V_{C^1}(\l_1^1)$ and
$V_{C^0}(a\l_1^0 + \l_8^0) \otimes V_{C^1}(\l_1^1).$ These sum to $V_{C^0}(a\l_1^0) \otimes V_{C^0}(\l_8^0) \otimes V_{C^1}(\l_1^1).$ Restricting to $L_X'$ we have $S^a(0010) \otimes (1010) \otimes (0100).$ The first tensor factor contains $(00a0) + (10(a-2)0)$ and the product of the other two factors is $(2001) + (0020) + (0101) + (1000) + (1110).$ Using Theorem \ref{LR} we find that full tensor product  contains $(10a0)^6.$  On the other hand $(10a0)$ can only arise from terms $(20a0), (01a0)$, and $(10(a-1)1)$ in $V^1$, so this contradicts our assumption that $V\downarrow X$ is MF.

	Now assume  $\l = \l_1 + a\l_9 + c\l_{10}$ and again assume $V \downarrow X$ is MF.  We claim that
$(1(a-1)10)^6$ appears in $V^2$ which will again be a contradiction.   This time $V^2$ contains $V_{C^0}(a\l_9^0) \otimes 
V_{C^1}(\l_1^1)$ and $V_{C^0}(\l_1^0 + \l_8^0 +(a-1)\l_9^0) \otimes V_{C^1}(\l_1^1).$ The situation is a little more complicated here as these terms sum to $(V_{C^0}(\l_1^0) \otimes V_{C^0}( \l_8^0 +(a-1)\l_9^0) \otimes V_{C^1}(\l_1^1))  -(V_{C^0}( \l_8^0 +(a-2)\l_9^0) \otimes V_{C^1}(\l_1^1)).$  Now $V_{C^0}( \l_8^0 +(a-1)\l_9^0)) \downarrow L_X'$ contains $(1(a-1)10)$, so the restriction of the first term contains $(1(a-1)10) \otimes (0010) \otimes (0100)$. Using Corollary \ref{LR_om2},  we find that this tensor product contains $(1(a-1)10)^7.$  To complete the argument we must consider the subtracted term $V_{C^0}( \l_8^0 +(a-2)\l_9^0) \otimes V_{C^1}(\l_1^1).$ The restriction of this term  is contained $S^{a-2}(0100) \otimes (1010) \otimes (0100).$ The first tensor factor is $(0(a-2)00) + (0(a-4)01) + (0(a-6)02) + \cdots $ with corresponding $S$-values $a-2, a-3, a-4, \ldots $. The term $(1(a-1)10)$ has $S$-value $a+1$ so this summand can only arise from $(0(a-2)00) \otimes (1010) \otimes (0100)$ and here it only occurs with multiplicity 1.  This establishes the claim and we have a contradiction. \hal



\begin{lem} \label{donnanewa5lem}
Let $X = A_5$, $W = V_X(\o_2)$ and $X < Y = SL(W) = A_{14}$.
If $\l = c\l_2$ with $c\ge 4$, then $V_Y(\l)\downarrow X$ is not MF.
\end{lem}


\pf We pass to $M=V^*$ and look for a contradiction in $M^4$.  Note that $M^4(Q_Y)\supseteq (00\dots 03)\otimes (00(a-3)3))\oplus (0\dots011)\otimes (00(a-2)1)$. Restricting these summands to $L_X'$, we see that each summand contributes $(01(a-4)3)^2$, for a total of four such summands. Now $(01(a-4)3)$ can only arise from summands in $M^3$ of the form $(01(a-4)2)$, $(11(a-4)3)$, $(01(a-5)4)$ and $(00(a-3)3)$.

All $L_Y'$-summands of $M^3(Q_Y)$ have highest weight of the form 
$\lambda^*-2\beta_{10}-x\beta_{11}-y\beta_{12}-z\beta_{13}-w\beta_{14}$ or 
$\lambda^*-\beta_9-2\beta_{10}-x\beta_{11}-y\beta_{12}-z\beta_{13}-w\beta_{14}$, 
for nonnegative integers $x,y,z,w$ satisfying:
\begin{itemize}
\item $x,y,z\geq 2$,
\item $y+2\geq 2x$,
\item $x+z\geq 2y$,
\item $y+w+a\geq 2z$,
\item $z\geq 2w$.
\end{itemize}

Henceforth let us write $\lambda^*-(abcdef)$ for the weight 
$\lambda^*-a\beta_9-b\beta_{10}-c\beta_{11}-d\beta_{12}-e\beta_{13}-f\beta_{14}$. It is easy to see that 
$\mu_1=\lambda^*-(022220)$  and $\mu_2=\lambda^*-(122221)$ afford the highest weights of $L_Y'$-summands of $M^3(Q_Y)$.
Note that these weights have respective restriction to $L_Y'$ as follows :
 $$(2\lambda_9+(a-2)\lambda_{13}+2\lambda_{14})|_{L_Y'}, \ \ (\lambda_8+(a-1)\lambda_{13})|_{L_Y'}.$$ 

We claim that these are the only weights of $L_Y'$-summands whose restriction gives rise to $L_X'$-summands of $S$-value
 at least $a-1$. To see this, we first note  the following two additional conditions, the first because of the lower bound on the 
$S$-value and the second because if $w=0$, then the corresponding weight has multiplicity exactly 1 in $M$ and occurs 
already in the $L_Y'$-summand afforded by $\mu_1$.
\begin{itemize}
\item $x+w\leq 5$,
\item $w\ne 0$.
\end{itemize}

We now work our way through the various pairs $(x,w)$ with $x+w\leq 5$, $x\geq 2$, $w\geq 1$, giving all of the
 details in the first few cases and leaving the details to the reader for the remaining cases. \medbreak

\emph{Case I. $w=1$, $x=2$.} We must determine if the weight $\lambda^*-(022yz1)$ or the weight $\lambda^*-(122yz1)$
affords the highest weight of a $L_Y'$-summand of $M^3(Q_Y)$. In the second case, we may assume $(y,z)\ne (2,2)$,
as this pair gives rise to the weight $\mu_2$. 

For the first case, we note that $z\leq a$ (this follows from the list of inequalities and $a\geq 4$). 
The weight $\lambda^*-(022yz1)$ is conjugate to the weight $\lambda^*-(002yz1)$, and the multiplicity of the 
latter in $M$ is the same as the multiplicity of the weight $\nu_1-(002yz1)$ in the $Y$-module with highest weight 
$\nu_1=z\lambda_{13}$ (Cavallin). If $z=2=y$, this  multiplicity is 1 and since $\lambda^*-(022yz1)$ already appears in the 
summand afforded by $\mu_1$, we assume henceforth that $(y,z)\ne (2,2)$. Now $\nu_1-(002yz1)$ is 
conjugate to $\nu_1-(0013(y+1)1)$, and for all allowed values of $y,z$ with $(y,z)\ne (2,2)$, we have 
$y+1\leq z$. Hence we again apply (Cavallin) and calculate the multiplicity by
considering the multiplicity of the weight $\nu_2 -(0013(y+1)1)$ in the module with highest weight
 $\nu_2=(y+1)\lambda_{13}$. The weight $\nu_2 -(0013(y+1)1)$ is conjugate to $\nu_2-(001241)$. If $y=2$, this 
weight has multiplicity 2,
while if $y\geq 3$, again using Proposition \ref{cav1} we see that this weight has multiplicity 3. 

Now we must consider the 
multiplicity of the weight $\lambda^*-(022yz1)$ in the irreducible afforded by $\mu_1$. Here we have 
$\lambda^*-(022yz1) = \mu_1-(000(y-2)(z-2)1)$. Since  $z-2\leq a-2$, we may 
calculate the multiplicity of the weight $\nu_1-(000(y-2)(z-2)1)$ in the module with highest weight 
$\nu_1 = (z-2)\lambda_{13}+2\lambda_{14}$. Here $\nu_1-(000(y-2)(z-2)1)$ is conjugate to 
$\nu_1-(000(y-2)(y-1)1)$. As
remarked above, $y+1\leq z$, so $y-1\leq z-2$, and we calculate the multiplicity of the weight 
$\nu_2-(000(y-2)(y-1)1)$ in the module with highest weight $\nu_2=(y-1)\lambda_{13}+2\lambda_{14}$. It is now 
straightforward (using Cavallin and treating the cases $y=2$ and $y\geq 3$ separately) to see that the 
multiplicity of this weight is $2$, respectively $3$, according to whether $y=2$ or $y\geq 3$. Hence, the 
weight does not afford a summand of $M^3(Q_Y)$.

Now consider the  weight  $\lambda^*-(122yz1)$. As above we have $z\leq a$ and
 so the multiplicity of this weight in $M$ is the same as the multiplicity of the weight $\sigma_1-(122yz1)$ in 
the module with highest weight $\sigma_1=z\lambda_{13}$. Here the weight $\sigma_1-(122yz1)$
is conjugate to $\sigma_1-(1223(y+1)1)$.  As $(y,z)\ne (2,2)$, we have $y+1\leq z$ and so the multiplicity of this 
weight is the same as the multiplicity of the weight $\sigma_2-(1223(y+1)1)$ in the module with highest weight 
$\sigma_2=(y+1)\lambda_{13}$. The weight is conjugate to $\sigma_2-(012341)$ and if $y=2$ its multiplicity is 3 and if $y\geq 3$ 
its multiplicity is 4.
Now we determine the multiplicity of the weight $\lambda^*-(122yz1)$ in the summands afforded by $\mu_1$ and $\mu_2$. 
In the first summand, this weight is of the form $\mu_1-(100(y-2)(z-2)1)$, where the multiplicity is 2, respectively 3, 
depending on whether $y=2$ or $y\geq 3$. In the summand afforded by $\mu_2$, this weight has the form 
$\mu_2-(000(y-2)(z-2)0)$, which has multiplicity 1 in both cases. So we see that $\lambda^*-(122yz1)$, for $(y,z)\ne (2,2)$ 
does not afford a $L_Y'$-summand of $M^3(Q_Y)$.\medbreak

\emph{Case II. $w=1$, $x=3$.} We must determine if the weight $\lambda^*-(023yz1)$ or the weight $\lambda^*-(123yz1)$
affords the highest weight of a $L_Y'$-summand of $M^3(Q_Y)$.

In both cases, we have $y\geq 4$, $z\geq 5$, and $z\leq a$. 
Now considering first the weight $\lambda^*-(023yz1)$, we see that its multiplicity in $M$ is the same as the 
multiplicity of the weight $\nu_1-(023yz1)$ in the module with highest weight $\nu_1=z\lambda_{13}$. 
This weight is then conjugate to the weight $\nu_1-(013y(y+1)1)$ and one checks that $y+1\leq z$ and so we replace 
the pair $(\nu_1,\nu_1-(013y(y+1)1)$ by $(\nu_2=(y+1)\lambda_{13},\nu_2-(013y(y+1)1))$, and subsequently 
by $(\nu_3=5\lambda_{13},\nu_3-(012341))$ and find that the weight has multiplicity 4 in $M$.
Now we must determine the multiplicity of this weight in the summand afforded by $\mu_1$ (it does not occur in the 
summand afforded by $\mu_2$). After a first application of (Cavallin), we replace the pair 
$(\mu_1,\mu_1-(001(y-2)(z-2)1)$
by the pair $(\sigma_1=(z-2)\lambda_{13}+2\lambda_{14}, \sigma_1-(001(y-2)(z-2)1)$.  Conjugating and using that 
$y+1\leq z$ as above, we reduce to the pair $(\sigma_2=(y-1)\lambda_{13}+2\lambda_{14},\sigma_2-(001(y-2)(y-1)1)$, 
and then finally by $(\sigma_3=3\lambda_{13}+2\lambda_{14},\sigma_3-(001231)$, where it is easy to see that the
 multiplicity is 4. So the weight $\lambda^*-(023yz1)$ does not afford the highest weight of a 
$L_Y'$-summand of $M^3(Q_Y)$.

Turn now to the second weight $\lambda^*-(123yz1)$. Arguing as above we see that the multiplicity of this weight in 
$M$ is 5, while the multiplicity in the summand afforded by $\mu_1$ is 4 and in the summand afforded by $\mu_2$
the weight occurs with multiplicity 1. So this weight as well does not afford a summand.\medbreak

\emph{Case III. $w=1$, $x=4$.} We must determine if the weight $\lambda^*-(024yz1)$ or the weight $\lambda^*-(124yz1)$
affords the highest weight of a $L_Y'$-summand of $M^3(Q_Y)$.

In both cases, we have $y\geq 6$ and $z\geq 8$. Then one checks that the inequalities imply that $z\leq a$. So then we 
proceed as in 
previous cases to determine the multiplicity of the weight in $M$, first replacing the pair $(\lambda^*,\lambda^*-(024yz1))$ by the pair 
$(\nu_1=z\lambda_{13},\nu_1-(024yz1))$. Conjugating and noting that $y+1\leq z$, we further replace by the pair 
$(\nu_2 = (y+1)\lambda_{13},\nu_2-(024y(y+1)1))$.
Then again using that $y\geq 6$, we are able to finally replace by the pair $(\nu_3=6\lambda_{13},\nu_3-(023561))$, and 
determine the 
multiplicity to be 4.  We must now find the multiplicity of the weight  $\lambda^*-(024yz1)$ in the summand afforded by $\mu_1$.
The weight here is of the form $\mu_1-(002(y-2)(z-2)1)$, and since $z-2\leq a-2$ and $y-1\leq z-2$, we apply Proposition \ref{cav1} several 
times to reduce this to determining the multiplicity of the weight $\sigma_1-(001241)$ in the $L_Y'$-module of highest weight 
$\sigma_1=4\lambda_{13}+2\lambda_{14}$. This is easily seen to be 4. Hence $\lambda^*-(024yz1)$ does not give rise to a summand
of $M^3(Q_Y)$.

Turn now to the second weight  $\lambda^*-(124yz1)$. Proceeding as above, we reduce to determining the multiplicity of the weight 
$\nu_1-(123451) $ in the irreducible with highest weight $\nu_1=5\lambda_{13}$. This multiplicity is 5. Now as just above, the multiplicity of this weight in the summand afforded by $\mu_1$ is 4, while the multiplicity in the summand afforded by $\mu_2$ is easily seen to be 
1. So
$\lambda^*-(124yz1)$ does not afford a summand.\medbreak

\emph{Case IV. $w=2$, $x=2$.} We must determine if the weight $\lambda^*-(022yz2)$ or the weight $\lambda^*-(122yz2)$
affords the highest weight of a $L_Y'$-summand of $M^3(Q_Y)$.

In both cases, we have $z\geq 4$. We start by determining the muultiplicity of the weight $\lambda^*-(022yz2)$ in $M$.
One first checks that $z\leq a$ and so the multiplicity is the same as the multiplicity of the weight $\nu_1-(022yz2)$ in the 
irreducible with highest weight $\nu_1=z\lambda_{13}$. The weight  $\nu_1-(022yz2)$ is conjugate to  $\nu_1-(022y(y+2)2)$.
Now if $(y,z)=(3,4)$, then it is straightfoward to determine the multiplicity, which is 4. In all other cases, $y+2\leq z$ and 
we may replace the pair $(\nu_1,\nu_1-(022y(y+2)2))$ by $(\nu_2=(y+2)\lambda_{13},\nu_2-(022y(y+2)2))$. Now conjugating, we are 
reduced to determining the multiplicity of the weight $\nu_2-(002462)$. 
There 
are two special cases, namely when $y=2$ and when $y=3$, where one checks directly that the multiplicity is 3, respectively 5. When $y\geq 4$, the
 multplicity is the same as the multiplicity of $\nu_3-(002462)$ in the irreducible with highest weight $\nu_3=6\lambda_{13}$, 
which is 6.

Now we must determine the multiplicity of the weight $\lambda^*-(022yz2)$ in the summand afforded by $\mu_1$. 
Here the weight is of 
the form $\mu_1-(000(y-2)(z-2)2)$ and the multiplicity is the same as the multiplicity of the weight 
$\sigma_1-(000(y-2)(z-2)2)$ in 
the $L_Y'$-module with highest weight $\sigma_1=(z-2)\lambda_{13}+2\lambda_{14}$. The special case $(y,z) = (3,4)$ 
gives rise to a 
multiplicity of 4, while otherwise, we may replace by the pair $(\sigma_2=y\lambda_{13}+2\lambda_{14},\sigma_2-(000(y-2)y2))$. 
Conjugating we reduce to the pair $(\sigma_3=y\lambda_{13}+2\lambda_{14},\sigma_3-(000242))$. We must again treat separately the cases 
$y=2$, respectively $ 3$, where we find that the multiplicity of the weight is 3, respectively 5. In the other cases, $y\geq 4$ and the multiplicity 
is the same as the multiplicity of the weight $\sigma_4-(000242)$ in the module with highest weight 
$\sigma_4=4\lambda_{13}+2\lambda_{14}$. This last multiplicity is 6. So in all cases, we see that there is no $L_Y'$-summand 
of $M^3(Q_Y)$ afforded by $\lambda^*-(022yz2)$.

Now turn to the second weight $\lambda^*-(122yz2)$.  As above we have $z\leq a$ and so the multiplicity of this weight in $M$ 
is the same as the multiplicity of the weight $\nu_1-(122yz2)$ in the irreducible with highest weight $\nu_1=z\lambda_{13}$. 
As above we must consider the special case where $(y,z)=(3,4)$ and find that the multiplicity is 6. 
The weight  $\nu_1-(122yz2)$ is conjugate to $\nu_1-(122y(y+2)2)$ and having treated the case $(y,z)=(3,4)$, we now may assume 
$y+2\leq z$ and calculate the multiplicity of $\nu_2-(122y(y+2)2)$ in the irreducible with highest weight $\nu_2=(y+2)\lambda_{13}$. 
The weight $\nu_2-(122y(y+2)2)$ is conjugate to $\nu_2-(122462)$. If $y=2$, this weight has multiplicity $4$, and multiplicity 
7 if $y=3$. When $y\geq 4$, the multiplicity is the same as the multiplicity of the weight $\nu_3-(122462)$ in the ireducible with
 highest weight $6\lambda_{13}$, which is 8.

Now we must determine the multiplicity of the weight $\lambda^*-(122yz2)$ in the summands afforded by $\mu_1$ and $\mu_2$. 
For $\mu_1$,
 we may argue precisely as in the consideration of the weight $\lambda^*-(022yz2)$ and find that the multiplicity is 4, 
respectively $3, 5$,
in the special cases $(y,z)=(3,4)$, respectively, $y=2$, $y=3$. And when $y\geq 4$, we find multiplicity 6. Now we must 
determine the multiplicity in the second summand. Here the weight has the form $\mu_2-(000(y-2)(z-2)1)$. Since $z-2\leq a-1$, 
this multiplicity is the same as the multiplicity of the weight $\sigma_1-(000(y-2)(z-2)1)$ in the $L_Y'$-module of highest 
weight $(z-2)\lambda_{13}$. Now conjugating and using that $y-1\leq z+2$, we are led to determine the multiplicity of the 
weight $\sigma_2-(0001(y-1)1)$ in the module with highest weight $(y-1)\lambda_{13}$. This  multiplicity is $2$ if $y\geq 3$ 
and 1 if $y=2$. Now combining this with the multiplicity in the summand afforded by $\mu_1$, we see that the weight 
$\lambda^*-(122yz2)$ affords no $L_Y'$-summand of $M^3(Q_Y)$.

The remaining cases, where $(x,w) = (3,2)$ and $(x,w) = (2,3)$ are entirely similar and we omit the details. \hal


\begin{lem}\label{a2one} 
For $X = A_2$, the following $X$-modules are not MF for $c>0$ and $i=1,2$:
\[
\wedge^2(c\om_1+\om_2)\otimes \om_i,\; S^2(c\om_1+\om_2)\otimes \om_i,\; S^3(c\om_1+\om_2)\otimes \om_i.
\]
\end{lem}

\pf Write $\d = c\om_1+\om_2$. Now $\wedge^2(\d)$ has summands of highest weights $2\d-\a_1$, $2\d-\a_2$ and $2\d-\a_1-\a_2$. 
The tensor products $(2\d-\a_2) \otimes \om_1$ and $(2\d-\a_1-\a_2) \otimes \om_1$ both have summands 
$(2c\om_1+\om_2)$, while $(2\d-\a_2) \otimes \om_2$ and $(2\d-\a_1-\a_2) \otimes \om_2$ both have $2c\om_1$. Hence 
$\wedge^2(\d) \otimes \om_i$ is not MF for $i=1,2$.

Similarly, $S^2(\d) \otimes \om_i$ has a multiplicity 2 summand $(2c\om_1+\om_2)$ (if $i=1$) or $((2c-1)\om_1+2\om_2)$ (if $i=2$); 
and $S^3(\d) \otimes \om_i$ has a multiplicity 2 summand $(3c\om_1+2\om_2)$ (if $i=1$) or $((3c-1)\om_1+3\om_2)$ (if $i=2$). \hal

\begin{lem}\label{mags}
{\rm (i)} For $X=A_2$ and $c\le 3$, the module $S^3(\om_1+\om_2) \otimes S^c(\om_2)$ is not MF.

{\rm (ii)} For $X = A_3$, $\d' = \om_1+\om_3$ and $\d'' = \om_2$, each of the following $L$-modules is not MF:
\[
\wedge^3(\d') \otimes \d'',\; \wedge^3(\d') \otimes S^2(\d''),\; \wedge^3(\d') \otimes \wedge^2(\d''),\; 
\wedge^3(\d') \otimes \wedge^3(\d'').
\]
\end{lem}

\pf This was checked using Magma. \hal

\begin{lem}\label{nonMF11} Let $X=A_3$, $W = V_X(2\o_1)$ and $X<Y = SL(W) = A_9$. Then $V_Y(c\lambda_5)\downarrow X$ is non-MF for $c\geq 2$.
\end{lem}

\pf  For $c=2$ this is a Magma check, so assume $c\geq 3$. Note that $V^1 = S^c(2\omega_2)$, while 
$V^2 = V_{C^0}(\lambda^0_4+(c-1)\lambda^0_5)\downarrow {L_X'}\otimes \omega_1$.
We first claim that $V^2 \supseteq (2\omega_1+(2c-3)\omega_2)^3$.  
In order to do so, we show that $V_{C^0}(\lambda^0_4+(c-1)\lambda^0_5)\downarrow {L_X'}$ has summands 
$(\omega_1+(2c-3)\omega_2)$, $(2\omega_1+(2c-2)\omega_2)$, and $(3\omega_1+(2c-4)\omega_2)$.
Upon tensoring these with $\omega_1$ we obtain the three summands $(2\omega_1+(2c-3)\omega_2)$.

To prove the claim, observe first that $c\l_5^0\otimes \l_5^0 = (\l_4^0+(c-1)\l_5^0) \oplus (c+1)\l_5^0$. 
Hence $(\l_4^0+(c-1)\l_5^0)\downarrow L_X' = S^c(02)\otimes (02)-S^{c+1}(02)$. The first summand contains 
$((0,2a) + (2,2c-4))\otimes (02)$, which gives $(1,2c-3)+(2,2c-2)^2+(3,2c-4)$. Of these only $(2,2c-2)$ occurs in 
$S^{c+1}(02)$, with multiplicity 1. This proves the claim. 

Thus $V^2 \supseteq (2\omega_1+(2c-3)\omega_2)^3$.  
We next show that only one of these summands can arise from an irreducible summand of $V^1$, which 
will give the result. Using the standard monomial theory, we see that an irreducible summand of 
$V^1$ which gives rise to a summand $(2\omega_1+(2c-3)\omega_2)$ 
has highest weight in $\{2\omega_1+(2c-4)\omega_2, \omega_1+(2c-2)\omega_2, 3\omega_1+(2c-3)\omega_2\}$.
The first of these is indeed the highest weight of an irreducible summand of $S^c(2\omega_2)$.
It is straightforward to see that the second one is not, while the third
is not even subdominant to the weight $2c\omega_2$. Hence we have the result.\hal

 \subsection{Tensor products, symmetric and exterior powers}\label{nonmftensor}

The section contains a selection of results showing that various symmetric powers, exterior powers and tensor products are non-MF.


\begin{lem}\label{wed3ex}
For $X = A_{l+1}\,(l\ge 2)$, the module $\wedge^3(\omega_1 + \omega_l + c\omega_{l+1})$ is not MF for $c \ge 0.$
\end{lem}

\pf      Let $\d = \omega_1 + \omega_l + c\omega_{l+1}$ and assume first that $c>0$.
Suppose by way of contradiction that $\wedge^3(\d) \downarrow X$ is MF.  First assume $l = 2.$ 
If $c = 1$ we can use Magma to obtain a contradiction.  So assume $c > 1.$  Then  replacing $\delta$ by $\delta^*$ and
 using Lemma \ref{lemma1.5(iii)(iv)}(i) we obtain a contradiction.  Therefore we may now assume $l \ge 3.$
   We will show that $V^2$ contains $(20\ldots 011)^5$ and  only $3$ copies of this irreducible arise from 
$V^1.$
   
  The MF supposition implies that $V^1$ is MF.  Using
Corollary \ref{V^2(Q_X)} one sees that the 
only possible irreducibles that can yield $(20\ldots  011)$ are $(20\ldots  02)$, $(20\ldots  010)$, $(30\ldots  011)$, and $(110\ldots  011)$ (or $(121)$ in the last case if $l = 3$). 
However the third weight here is $(30\ldots  03) - \a_l$, so this cannot occur in $V^1=\wedge^3(\d).$ So this verifies the last statement
 of the above paragraph.
 
 Now $V^2 \supseteq \wedge^2(10\ldots 01) \otimes ((10\ldots  010) + (0\ldots  01) +(10\ldots  02)).$ 
One checks that $\wedge^2(10\ldots 01) \supseteq (20\ldots 010) + (10\ldots 01) + (010\ldots 02)$ (this is actually an equality).
 By Lemma \ref{compfactor}(ii), $(20\ldots 010) \otimes (10\ldots 010)  \supseteq (20\ldots 011)^2$
 and clearly $(10\ldots 01) \otimes (10\ldots 010) \supseteq (20\ldots 011).$
 Also $(20\ldots 10) \otimes (0\ldots 01)$  $\supseteq (20\ldots 011).$ Finally, an easy application of Theorem \ref{LR} shows that
 $(010\ldots 02) \otimes (10\ldots 010) \supseteq (20\ldots 011).$ 
 Therefore $V^2 \supseteq (20\ldots 011)^5$, giving a contradiction.
 
Now suppose $c=0$, so $\d= \om_1 + \om_l.$ We can assume $l \ge 4$ since a Magma computation gives the result for smaller values of $l.$  
We will work with the dual module $M =\wedge^3(V_X(\om_2 + \om_{l+1}).$ As usual define $M^i = M^i(Q_Y)\downarrow L_X'$.  
Note that the embeddings of $L_X'$ in $C^0$ and $C^1$ are given by the representations $\om_2$ and $(\om_1 \oplus  (\om_2 + \om_l))$, 
respectively.  Then $M^1 = \wedge^3(\om_2) = 2\om_3 \oplus  (2\om_1 + \om_4).$   And $M^2 = \wedge^2(\om_2) \otimes (\om_1 \oplus  (\om_2 + \om_l)).$  
Now $ \wedge^2(\om_2) = \om_1 + \om_3$ and Lemma \ref{aibjnotzero} shows that $ (\om_1+ \om_3) \otimes (\om_2 + \om_l) \supset
(\om_1 + \om_4)^2. $  Also $ (\om_1+ \om_3) \otimes (\om_1) \supseteq  \om_1 + \om_4$, so  that $M^2 \supseteq (\om_1 + \om_4)^3. $  
At most one such summand can arise from $M^1$, so $M \downarrow X$ is not MF, as required. \hal


\begin{lem}\label{MFtensors} Let $X=A_{l+1}$ ($l\ge 1$),  and let $W = V_X(r\o_1)$ with $r\ge 3$. Suppose $Y = SL(W)=A_n$ and $\l$ is as in Tables $\ref{TAB2}-\ref{TAB4}$ of Theorem $\ref{MAINTHM}$ for $\d = r\o_1$. Suppose also that $a,b\ge 1$.
\begin{itemize}
\item[{\rm (i)}]  If $V_Y(\l)\downarrow X\otimes V_X(a\omega_1)$ is MF then $(\l,a)\in\{(2\lambda_1,1),(\lambda_2,1)\}$.
\item[{\rm (ii)}] If $V_Y(\l)\downarrow X\otimes V_X(b\omega_{l+1})$ is MF then $(\l,b)\in\{(2\lambda_1,1),(\lambda_2,1)\}$.
\end{itemize}
\end{lem}

\pf Write $T_a = V_Y(\l)\downarrow X\otimes V_X(a\omega_1)$, $T_b=V_Y(\l)\downarrow X\otimes V_X(b\omega_{l+1})$.
In Table \ref{tensy} below we give, for each possible $\l$, some summands of the restriction $V_Y(\l)\downarrow X$, and repeated summands of $T_a$ and $T_b$. In some columns of the table we delete terms $\om_3$ and $\om_4$ when $l = 2, 3$, respectively.
The multiplicities follow by application of Proposition \ref{pieri}.

\begin{table}[h]
\caption{} \label{tensy}
{\small 
\[
\begin{array}{|l|l|l|l|}
\hline
\l &  V_Y(\l)\downarrow X \supseteq & \hbox{Repeated summand} & \hbox{Repeated summand} \\
  & & \hbox{of }T_a & \hbox{of }T_b \\
\hline
\l_1+\l_n & \sum_{i=0}^{r-1} ((r-i)\o_1+(r-i)\o_{l+1}) & (a+r-1)\o_1+(r-1)\o_{l+1} &  \\
2\l_1\,(a,b\ge 2) & 2r\o_1 \oplus ((2r-4)\o_1+2\o_2) & (2r+a-4)\o_1+2\o_2 & (2r-2)\o_1+\\
&&& (b-2)\o_{l+1}  \\
\l_2\,(a,b\ge 2) & ((2r-2)\o_1+\o_2)\oplus  & (2r+a-6)\o_1+3\o_2 & (2r-4)\o_1+\o_2+\\
    & ((2r-6)\o_1+3\o_2) & & (b-2)\o_{l+1}  \\
\l_3 & ((3r-6)\o_1+3\o_2)\oplus & (3r+a-7)\o_1+2\o_2+\o_3 & (3r-7)\o_1+3\o_2+ \\
        & ((3r-7)\o_1+2\o_2+\o_3)  & & (b-1)\o_{l+1} \\
3\l_1 & ((3r-4)\o_1+2\o_2) \oplus & (a+3r-6)\o_1+3\o_2 & (3r-5)\o_1+2\o_2+\\
  & ((3r-6)\o_1+3\o_2) & & (b-1)\o_{l+1}  \\
\l_4 & ((4r-7)\o_1+2\o_2+\o_3) \oplus & (a+4r-9)\o_1+3\o_2+\o_3 & (4r-8)\o_1+2\o_2+\o_3 +\\
& ((4r-9)\o_1+3\o_2+\o_3) && (b-1)\o_{l+1},\,l>1  \\
       &  & & (4r-8)\o_1+ \\
&&& (b+1)\o_2,\,l=1 \\
4\l_1\,(r=3) & (8\o_1+2\o_2)\oplus (6\o_1+3\o_2) & (a+6)\o_1+3\o_2 & 7\o_1+2\o_2+(b-1)\o_{l+1} \\
\l_5\,(r=3,l\ge 2) & (7\o_1+2\o_2+\o_4)\oplus & (a+5)\o_1+3\o_2+\o_4 & 6\o_1+2\o_2+\o_4+ \\
  & (5\o_1+3\o_2+\o_4)  & & (b-1)\o_{l+1}\\
\l_1+\l_2\,(r=3) & (5\o_1+2\o_2)\oplus (3\o_1+3\o_2) & (a+3)\o_1+3\o_2 & 4\o_1+2\o_2+(b-1)\o_{l+1} \\
\hline
\end{array}
\]}
\end{table}

One case is omitted from the table: $\l=\l_5$ with $r=3,\,l=1$. This is handled easily using Theorem \ref{LR}. \hal

\begin{lem}\label{tensor_1} Let $X=A_{l+1}$ ($l\ge 1$),  let $W = V_X(r\o_1)$ with $r\ge 3$, and let 
$Y = SL(W)=A_n$. For $b\ge 1$ the restriction  
$V_Y(\lambda_{n-1}+b\lambda_n)\downarrow X$ has summands $(\omega_l+((b+2)r-2)\omega_{l+1})$ and  
$(2\omega_l+((b+2)r-4)\omega_{l+1})$.
\end{lem}

\pf Note that $V_Y((b+1)\lambda_n)\otimes V_Y(\lambda_n) = V_Y((b+2)\lambda_n)\oplus V_Y(\lambda_{n-1}+b\lambda_n)$. We show that $(\omega_l+((b+2)r-2)\omega_{l+1})$  occurs
as a summand of  $S^{b+1}(r\omega_{l+1})\otimes (r\omega_{l+1})$ and does not occur in 
$S^{b+2}(r\omega_{l+1})$ and that $(2\omega_l+((b+2)r-4)\omega_{l+1})$ occurs with 
multiplicity two in  $S^{b+1}(r\omega_{l+1})\otimes (r\omega_{l+1})$ and with multiplicity one in $S^{b+2}(r\omega_{l+1})$. This will establish the result.

Since $r\geq 3$ and $b\geq 1$, $S^{b+1}(r\omega_{l+1})$ has summands $((b+1)r\omega_{l+1})$ and $(2\omega_l+((b+1)r-4)\omega_{l+1})$. Tensoring the first of these with 
$(r\omega_{l+1})$ yields a summand $(\omega_l+((b+2)r-2)\omega_{l+1})$. It is 
straightforward to see that this does not occur as a summand of $S^{b+2}(r\omega_{l+1})$.
The same tensor product also has a summand $(2\omega_l+((b+2)r-4)\omega_{l+1})$. The latter
also occurs as a summand in
$(2\omega_l+((b+1)r-4)\omega_{l+1})\otimes (r\omega_{l+1})$. So it remains to see that 
$(2\omega_l+((b+2)r-4)\omega_{l+1})$ occurs with multiplicity one in $S^{b+2}(r\omega_{l+1})$.
Again this is a straightforward check.\hal

\begin{lem}\label{newprop}
Let $X=A_{l+1}\,(l\ge 3)$, let $2\le c \le 5$ and let $3\le i \le \frac{l+3}{2}$. 
If $c=3$, assume $i\le 6$; if $c=4$ assume $i\le 4$; and if $c=5$ assume $i=3$.
Then none of the following modules are  MF:
\begin{itemize}
\item[{\rm (i)}] $S^c(\om_i) \otimes (\om_{i-1})^*$
\item[{\rm (ii)}] $\wedge^c(\om_i) \otimes (\om_{i-1})^*$, except for $(c,i,l) = (2,3,3)$
\item[{\rm (iii)}] $S^c(\om_i) \otimes \om_{i-1}$ for $i=\frac{l+3}{2}$
\item[{\rm (iv)}] $\wedge^c(\om_i) \otimes \om_{i-1}$ for $i=\frac{l+3}{2}$, except for $(c,i,l) = (2,3,3)$
\item[{\rm (v)}] $S^c(\om_3) \otimes S^2(\om_2^*)$ for $l \ge 4$
\item[{\rm (vi)}] $\wedge^c(\om_3) \otimes S^2(\om_2^*)$ for $l\ge 4$
\item[{\rm (vii)}] $((\om_3\otimes \wedge^2(\om_3))/\wedge^3(\om_3))\otimes U$, where $U$ is  $\om_2^*$  or $S^2(\om_2)^*$, for $l\ge 4$
\item[{\rm (viii)}] $\wedge^6(\o_4)\otimes \o_j$ for $l=5$, $j=3,4$
\item[{\rm (ix)}] $((\om_4\otimes \wedge^2(\om_4))/\wedge^3(\om_4))\otimes \o_j$, for  $l=5$, $j=3,4$.
\end{itemize}
\end{lem}

\pf All cases where $l \le 9$ can easily be settled using Magma. In particular, the exceptional cases in (ii) and (iv) appear here, as well as cases (viii), (ix). So from now on assume $l \ge 10$.

Next suppose that $i = \frac{l}{2}$, $\frac{l+1}{2}$,  $\frac{l+2}{2}$ or $\frac{l+3}{2}$. Then  $\om_i, \om_{i-1}$ and $(\om_{i-1})^*$ are fundamental weights associated with fundamental roots at or near the center of the Dynkin diagram.  The two nodes have distance at most 4 from each other.  Now consider the Levi factor of $X$ with base $\a_{i-2},\a_{i-1},\a_i,\ldots , a_{i+5},$ of rank  $8.$ The assertions hold for this group  by the above.  It follows that each tensor product appearing in (i)-(iv) has a summand of multiplicity at least two where the high
weight is a rational combination of roots in the root system of the Levi subgroup.  Consequently the result
holds of $X$ as well.  So we now asssume $i < \frac{l}{2}.$  In particular (iii) and (iv) have been settled.

We now consider parts (i) and (ii). Let $\om_k = \om_{i-1}^*$ so that $k = l-i+3.$ In each case we will produce  
summands of highest weights $\nu_1, \nu_2$ such that $\nu_1 \otimes \om_k$ and $\nu_2\otimes \om_k$ both contain a common summand of highest weight $\nu.$  The weights are as follows, where we omit certain terms such as $\om_0$ or $\om_{k+2}$ when $k+2 > l+1$:
\[
\begin{array}{llll}
\hline
\hbox{module} & \nu_1 & \nu_2 & \nu \\
\hline
\wedge^2(\om_i) &  \om_{i-1}+\om_{i+1} & \om_{i-3}+\om_{i+3} & \om_{i-3}+\om_{i+1}+\om_{k+2} \\
S^2(\om_i) &   2\om_i & \om_{i-2}+\om_{i+2} & \om_{i-2}+\om_i+\om_{k+2} \\
\wedge^3(\om_i) &  2\om_{i-1} + \om_{i+2} & \om_{i-2} + 2\om_{i+1} &  \om_{i-2}+\om_{i-1}+\om_{i+1}+\om_{k+2}\\
 S^3(\om_i) &  3\om_i & \om_{i-2}+\om_i + \om_{i+2} & \om_{i-2} + 2\om_i + \om_{k+2}\\
\wedge^4(\om_i) &  \om_{i-2} + \om_{i-1} + \om_{i+1} + \om_{i+2} & 3\om_{i-1} + \om_{i+3} &  \om_{i-2} + 2\om_{i-1} + \om_{i+2} + \om_{k+2}\\
S^4(\om_i) &  4\om_i &  \om_{i-2} + 2\om_i + \om_{i+2} &  \om_{i-2} + 3\om_i + \om_{k+2}\\
\wedge^5(\om_i) (i=3) & 4\om_2 + \om_7 &  \om_1 + 2\om_2 + \om_4 + \om_6 &  \om_1 +3\om_2 + \om_6+
\o_{k+2}\\
S^5(\om_i) (i = 3) & \om_1 + 3\om_3 + \om_5 &  2\om_1 + \om_3 + 2\om_5 &  2\om_1 + 2\om_3 + \om_5 + \o_{k+2} \\
\hline
\end{array}
\]
All the assertions in (i) and (ii) follow from Lemma \ref{aibjnotzero} or Theorem \ref{LR}.  

Now consider parts (v) and (vi). If $c \ge 2$, then $S^c(\om_3) $ contains $\nu_1 = c\om_3$ and $\nu_2 = \om_1 + (c-2)\om_3 + \om_5$.  As $S^2(\om_2)^* $ contains $2\om_{l}$,  it follows from Lemma \ref{aibjnotzero} that the tensor products of $\nu_1$ and $\nu_2$ with $2\om_{l}$ each contain $\nu = \om_1 + (c-1)\om_3 + \om_{l}.$ Therefore $S^c(\om_3) \otimes 2\om_l$ is not MF, giving the conclusion for these cases.

 Similarly $\wedge^4(\om_3)$  contains $\nu_1 =  \om_2 + \om_3 + \om_7$ and $\nu_2 = \om_1 +\om_4+\om_7.$  Lemma \ref{aibjnotzero} shows that tensoring each of these with $2\om_{l}$ produces a summand
 $\om_1 +\om_2 + \om_7 + \om_{l}.$  Next, $\wedge^5(\om_3)$ 
contains $\nu_1 = 2\om_2 + \om_3 + \om_8$ and $\nu_2 = \om_2 + \om_5 + \om_8$, and tensoring 
 each with $2\om_{l}$ produces a summand $\nu = \om_2 + \om_3 +\o_8+ \om_{l}.$ And  $\wedge^3(\om_3)$ contains $\om_1 + 2\om_4$ and $2\om_2 + \om_5$ and we argue with Theorem \ref{LR} that each of these tensored with $2\om_{l}$ contains $\om_2 + \om_4 + \om_{l+1}$. 
 Finally $\wedge^2(\om_3) \otimes S^2(\om_2)^* \supseteq 
((\om_2+\om_4)\oplus \om_6) \otimes (2\om_{l})$ and  Lemma \ref{aibjnotzero} shows that this contains $\om_4+\om_{l}$ with multiplicity 2. 

Finally consider part (vii). The first tensor factor contains both $\nu_1 = \om_2 + \om_3 + \om_4$ and $\nu_2 = 
\om_1 + \om_3 + \om_5$.  If $U = S^2(\om_2)^*$, then the result follows from Proposition \ref{stem1.1.A}  since $\nu_1 \otimes 2\om_{l}$ is not MF.  For $U = (\om_2)^* = \om_{l}$,  set 
$\nu =  \om_1 + \om_3 + \om_4+ \om_{l+1}$.  Then using Lemma \ref{abtimesce} we find
that tensoring $\nu_1$ and $\nu_2$ with $U$ each yield a summand $\nu.$  \hal

\begin{lem}\label{morenonmf}
Let $X=A_{l+1}$ with  $l\ge 1$. For $c>0$, the following $X$-modules are not MF:
\begin{itemize}
\item[{\rm (i)}] $\wedge^2(c\om_1+\om_2)\otimes \om_{l+1},\, \wedge^2(c\om_1+\om_2)\otimes \wedge^2(\om_{l+1})$
\item[{\rm (ii)}]  $S^2(c\om_1+\om_i)\otimes \om_{l+3-i},\; S^2(c\om_1+\om_i)\otimes S^2(\om_{l+3-i})$ for $i=2,3,l+1$
\item[{\rm (iii)}] $S^2(c\om_1+2\om_2)\otimes 2\om_{l+1},\; S^2(c\om_1+2\om_2)\otimes S^2(2\om_{l+1})$
\item[{\rm (iv)}]  $S^2(c\om_1+\om_2)\otimes \wedge^2(\om_{l+1})$.
\end{itemize}
\end{lem}

\pf  For (i), we can assume by Lemma \ref{a2one} that $l\ge 2$.  Let $\d = c\om_1+\om_2$. Observe that $\wedge^2V(\d)$ has summands $V(2\d-\a_1)$ and $V(2\d-\a_1-\a_2)$. By Lemma \ref{abtimesce}, the tensor product of each of these with $V(\om_{l+1})$ (resp. $V(\om_{l})$) has a summand $V((2c-1)\om_1+2\om_2)$ (resp. $V((2c-1)\om_1+2\om_2+\om_{l+1})$). Part (i) follows.

For the remaining parts, we similarly compute (using Lemma \ref{abtimesce} and  results in Section \ref{LRres}) a multiplicity 2 summand $V(\mu)$ in the relevant tensor product, where $\mu$ is as in Table \ref{relev}. \hal

\begin{table}[h]
\caption{}\label{relev}
\[
\begin{array}{|ll|}
\hline
\hbox{tensor product} & \mu \\
\hline
S^2(c\om_1+\om_i)\otimes \om_{l+3-i} & (2c-1)\om_1+2\om_2\,(i=2, l\ge 2) \\
                                                                 & (2c-1)\om_1+2\om_3+\om_{l+1}\,(i=3, l\ge 3) \\
                                                                  & (2c-1)\om_1+\om_2+\om_{l+1}\,(i=l+1 \ge 3) \\
S^2(c\om_1+\om_i)\otimes S^2(\om_{l+3-i}) & (2c-1)\om_1+2\om_2+\om_{l+1}\,(i=2) \\
                                                                 & (2c-1)\om_1+2\om_3+\om_{l}+\om_{l+1}\,(i=3, l\ge 3) \\
                                                                  & (2c-1)\om_1+2\om_2+\om_{l+1}\,(i=l+1\ge 3) \\
S^2(c\om_1+2\om_2)\otimes 2\om_{l+1} & (2c-1)\om_1+4\om_2+\om_{l+1} \\
S^2(c\om_1+2\om_2)\otimes S^2(2\om_{l+1}) & (2c-1)\om_1+4\om_2+3\om_{l+1} \\
S^2(c\om_1+\om_2)\otimes \wedge^2(\om_{l+1}) & (2c-1)\om_1+2\om_2+\om_{l+1}\,(l\ge 2) \\
\hline
\end{array}
\]
\end{table}

\begin{lem}\label{nonemf}
Let $X=A_{l+1}$, $l\ge 1$, and let $a\ge b\ge 2$, $(a,b)\ne (2,2)$. 
Let $Z$ denote either of the $L$-modules $S^2(2\om_1)$ or $\wedge^3(2\om_1)\,(l=1)$.
Then the following $X$-modules are not MF:
\[
\begin{array}{l}
(a\om_1) \otimes S^c(b\om_1)\;\; (2 \le c \le 4) \\
(a\om_1) \otimes \wedge^c(b\om_1) \;\;(2\le c\le 5,\,(b,c,l) \ne (2,2,l), (2,4,1) \hbox{ or }(2,5,1)) \\
\wedge^c(a\o_1) \otimes \o_i\;(a\ge 3,\,3\le c\le 4,\,2\le i\le 5)  \\
\wedge^5(3\o_1) \otimes \o_i\;(2\le i\le 5)  \\
S^c(a\om_1) \otimes Z\;\; (2\le c \le 4) \\
\wedge^c(a\om_1) \otimes Z \;\;(2\le c\le 5) \\
(3\om_1) \otimes ((b\om_1)\otimes \wedge^2(b\om_1)/\wedge^3(b\om_1))\;\;(b=2 \hbox{ or } 3) \\
((3\om_1)\otimes \wedge^2(3\om_1)/\wedge^3(3\om_1)) \otimes Z \\
S^b(2\omega_1)\otimes (a\omega_{l+1}) \,(a=2,3,\,b\geq 2) \\
S^b(2\omega_1)\otimes (2\omega_1)\, (b\geq 2).
\end{array}
\]
\end{lem}

\pf We compute a summand $V(\mu)$ of multiplicity at least 2 in each tensor product, as in Table \ref{relev1}.
(Note that some terms in the table may not be present for small ranks.) \hal

\begin{table}[h]
\caption{}\label{relev1}
\[
\begin{array}{|ll|}
\hline
\hbox{tensor product} & \mu \\
\hline
(a\om_1) \otimes S^c(b\om_1) & (a+bc)\om_1-2\a_1 \\
(a\om_1) \otimes \wedge^2(b\om_1) & (a+2b)\om_1-3\a_1 \\
(a\om_1) \otimes \wedge^3(b\om_1) & (a+3b)\om_1-4\a_1-\a_2 \\
(a\om_1) \otimes \wedge^4(b\om_1) & (a+4b)\om_1-5\a_1-\a_2,\;b\ge 3 \\
                                                             & (a+8)\om_1-5\a_1-2\a_2-\a_3,\;b=2 \\
(a\om_1) \otimes \wedge^5(b\om_1) &(a+5b)\om_1 -8\a_1-2\a_2,\; l = 1\\
&(a+5b)\om_1 -7\a_1-3\a_2-\a_3,\; l \ge 2\\

\wedge^3(a\o_1)\otimes \o_i \,(2\le i\le 5) & (3a-7)\o_1+3\o_2+\o_{i+1} \\
\wedge^4(a\o_1)\otimes \o_i \,(2\le i\le 5) & (4a-9)\o_1+2\o_2+\o_3+\o_{i+2},\;l\ge 2 \\
           & (4a-9)\o_1+\o_2,\;l=1 \\
\wedge^5(3\o_1)\otimes \o_i \,(2\le i\le 5) & 4\o_1+\o_2+\o_3+\o_4+\o_{i+2},\;l\ge 2 \\
                  & 4\o_1+2\o_2,\,l=1 \\
S^c(a\om_1) \otimes Z & ca\om_1 +2\om_2, Z = S^2(2\om_1) \\
&(ca-1)\om_1 +2\om_2, Z = \wedge^3(2\om_1) \\
\wedge^c(a\om_1) \otimes Z &(2a-2)\om_1 + 3\om_2, c=2, Z = S^2(2\omega_1) \\
& (2a-3)\om_1 + 3\om_2, c=2, Z = \wedge^3(2\omega_1) \\
&(3a-2)\om_1 + 3\om_2 -\a_1-\a_2, c=3, Z = S^2(2\omega_1) \\
&(3a-3)\om_1 + 3\om_2+\o_3, c=3, Z = \wedge^3(2\omega_1) \\
&(4a-6)\om_1 + 5\om_2 -\a_1-\a_2, c=4, Z = S^2(2\omega_1) \\
&(4a-6)\om_1 + 3\om_2+2\o_3, c=4, Z = \wedge^3(2\omega_1) \\
&7\omega_1 + 6\omega_2 - 2\a_1-2\a_2, c = 5 (a = 3), Z = S^2(2\omega_1)\\
&  5\om_1 + 5\om_2+2\o_3,  c=5 (a = 3), Z = \wedge^3(2\omega_1) \\

(3\om_1) \otimes ((b\om_1)\otimes \wedge^2(b\om_1)/\wedge^3(b\om_1)) & 8\om_1 +2\om_2, b = 3  \\
&    5\om_1 +2\om_2, b = 2    \\
((3\om_1)\otimes \wedge^2(3\om_1)/\wedge^3(3\om_1)) \otimes 
Z & 5\om_1 + 4\om_2, Z=S^2(2\omega_1) \\
&  4\om_1 + 4\om_2+\o_3,  Z = \wedge^3(2\omega_1)\\
S^b(2\omega_1)\otimes (a\omega_{l+1})  & (2b-2)\o_1+(a-2)\o_{l+1} \\
S^b(2\omega_1)\otimes (2\omega_1)\, (b\geq 2) & (2b-2)\o_1+2\o_2 \\
\hline 
\end{array}
\]
\end{table}

\begin{lem}\label{manymany}
Let $X = A_{l+1}$ with $l\ge 2$,  let $\d = a\om_1+b\om_{l+1}$ with $a\ge b\ge 1$ and  let $W = V_X(\d)$.
\begin{itemize}
\item[{\rm (i)}] Suppose $a\ge 2$. Then $S^cW\; (c=3,4)$ and $\wedge^cW\; (c=3,4,5)$ are not MF.
\item[{\rm (ii)}] If $a=3$, then $(W\otimes \wedge^2W)/\wedge^3W$ is not MF.  
\end{itemize}
\end{lem}

\pf 
We first prove part (ii). Note that $\wedge^2(W) \supseteq ((2a-2)10 \ldots  (2b)) + ((2a)0\ldots  01(2b-2)) +((2a-1)0\ldots  0(2b-1))$, where we use Lemma \ref{lemma1.5(i)(ii)}(ii) for the last summand. The tensor product of $W$ with each of these summands yields a summand $((3a-1)0\ldots  0(3b-1)).$    Hence  $W \otimes \wedge^2(W) \supseteq ((3a-1)0\ldots  0(3b-1))^3.$  The corresponding highest weight can be written as $3\d - (\a_1 + \cdots + \a_n)$  and an easy count shows that there is only  one such summand in $\wedge^3(W).$ 

Now we prove (i). We proceed using level analysis where $X < Y = SL(W)$.  First consider $\wedge^c(W)$.  Then using the usual parabolic subgroup we have $ V^1(Q_Y) = V_{C^0}(\l_c^0)$ and $V^2(Q_Y) = V_{C^0}(\l_{c-1}^0) \otimes V_{C^1}(\l_1^1).$   The embedding of $L_X'$ in $C^1$ is given by the representation $(a0\ldots  01) + ((a-1)00 \ldots  0) = (a0\ldots  0) \otimes (0 \ldots  01).$  Therefore $V^2  = V_{C^0}(\l_{c-1}^0) \otimes (a0\ldots  0) \otimes (0 \ldots  01).$  Now $a0\ldots  0 = V_{C^0}(\l_1^0) \downarrow L_X'$ so we can write the restriction
as $((V_{C^0}(\l_{c-1}^0) \otimes  V_{C^0}(\l_1^0)) \downarrow L_X') \otimes (0\ldots  01).$  Moreover, $V_{C^0}(\l_{c-1}^0) \otimes  V_{C^0}(\l_1^0) = V_{C^0}(\l_c^0) + V_{C_0}(\l_1^0 + \l_{c-1}^0).$   

At this point we tensor each of the above  two summands with $(0 \ldots  01).$  By Corollary \ref{cover}(ii) it suffices to show that the second tensor product fails to be MF.
 We have $V_{C_0}(\l_1^0 + \l_{c-1}^0) \downarrow L_X' = (a0\ldots  0) \otimes \wedge^{c-1}(a0\ldots  0) / \wedge^c(a0\ldots  0).$  Now $ \wedge^{c-1}(a0\ldots  0)$ contains an irreducible summand of highest weight
 $(2a-2)10 \ldots  0), ((3a-3)010 \ldots 0), ((4a-7)210 \ldots  0)$, according as $c = 3,4,5,$ respectively.
 Tensoring with $(a0\ldots  0)$ we find that the tensor product contains $((3a-2)10 \ldots  0) +(3a-4)20 \ldots  0)), ( ((4a-3)010 \ldots 0) +  ((4a-5)110 \ldots 0)), (((5a-7)210 \ldots  0) + ((5a-9)310 \ldots  0)$, respectively.
It is easy to see these summands do not lie in $\wedge^c(a0 \ldots  0).$ For example the weights $ ((5a-7)210 \ldots  0)$ and $((5a-9)310 \ldots  0)$ have the form $(5a0 \ldots  0)- 4\a_1-\a_2 $ and
$(5a0 \ldots  0)- 5\a_1-\a_2 $, respectively and these weights do not occur in $\wedge^5(a0\ldots 0).$
Now tensoring with $(0\ldots  01)$, we obtain repeated summands $((3a-3)10\ldots  0)$, $((4a-4)01\ldots  0)$, $((5a-8)21\ldots  0)$, 
respectively.   This completes the proof for $\wedge^c(W).$
 
 Now consider $S^c(W)$ for $c = 3,4.$ Then $ V^1(Q_Y) = V_{C^0}(c\l_1^0)$ and $V^2(Q_Y) = V_{C^0}((c-1)\l_1^0) \otimes V_{C^1}(\l_1^1).$   Arguing as above we see that $V^2 = 
 S^c(a0\ldots  0) \otimes (0\ldots  01) + (V_{C^0}((c-2)\l_1^0 + \l_2) \downarrow L_X') \otimes (0\ldots  01).$
Applying Corollary \ref{cover}(ii)  we see that it will suffice to show that the second summand is not MF.
For $c = 3$ $(V_{C^0}((c-2)\l_1^0 + \l_2) \downarrow L_X') = ((a0\ldots  0) \otimes \wedge^2( a0\ldots  0)) -\wedge^3(a0\ldots 0)$ and arguing as above we see that this contains $((3a-2)10\ldots  0)$ and
$((3a-4)20\ldots  0)$.  Then tensoring with $(0 \ldots  01)$ we have a repeated composition factor $((3a-3)10 \ldots  0). $ And for $c = 4$, $V_{C^0}((c-2)\l_1^0 + \l_2^0) \downarrow L_X' = \wedge^2(S^2(a0\ldots  0))$.
This contains $((4a-2)10 \ldots  0)$ and $((4a-4)20 \ldots  0)$, and
tensoring with $(0\ldots  01)$ we have a repeated composition factor $((4a-3)10\ldots 0).$  \hal

 \section{$L(\nu )\ge 2$ results}\label{lge2}
 
 In this subsection we establish a result showing that almost all of the $Y$-modules in Tables \ref{TAB1}-\ref{TAB4} of Theorem \ref{MAINTHM} have restriction to $X$ containing  a summand of 
highest weight $\nu$ satisfiying $L(\nu) \ge 2.$  
 
\begin{lem}\label{twolabs} Let $X = A_{l+1}\,(l\ge 1)$, $W = V_X(\d)$, $X<Y = SL(W)$ and $V = V_Y(\l)$. 
Suppose $\l,\d$ are in Tables $\ref{TAB1}-\ref{TAB4}$ of Theorem $\ref{MAINTHM}$ (up to duals), and are not in the following list:
\[
\begin{array}{|l|l|l|}
\hline
\l & \d & l \\
\hline
2\l_1 & 2\omega_1,\,\omega_2 & \hbox{any} \\
\l_3 & 2\omega_1 & 1 \\
\l_3,\,c\l_1 & \omega_2 & 2 \\
\hline
\end{array}
\]
 Then $V\downarrow X$ has a composition factor $V_X(\nu)$ such that $L(\nu) \ge 2$.
\end{lem}

\pf First consider 
the examples in Table \ref{TAB1}. These satisfy $L(\d)\ge 2$. If $\l = 2\l_1$ then $V = S^2W$ and we can take $\nu = 2\d$. 
If $\l = \l_2$ then $V = \wedge^2W$ and we take $\nu$ as follows:
\[
\begin{array}{|ll|}
\hline
\d & \nu \\
\hline
\om_1+c\om_i & 2\d-\a_1 = \om_2+2c\om_i\;(i>2) \\
                        & 2\d-\a_2 = 3\om_1+(2c-2)\om_2+\om_3\;(i=2,\,l\ge2) \\
                        & 2\d-\a_2 = 3\om_1+(2c-2)\om_2\;(i=2,\,l=1,\,c\ge 2) \\
                        & 2\d-\a_1-\a_2 = \om_1+\om_2\;(i=2,\,l=1,\,c=1) \\
\hline
c\om_1+\om_i & 2\d-\a_1 = (2c-2)\om_1+\om_2+2\om_i \\
(c>1) & \\
\hline
c\om_i+d\om_{i+1} & 2\d-\a_i = \om_{i-1}+(2c-2)\om_i+(2d+1)\om_{i+1} \\
(1<i<l) & \\
\hline
2\om_1+2\om_{l+1} & 2\d-\a_1 = 2\om_1+\om_2+4\om_{l+1} \\
\hline
2\om_1+2\om_2 & 2\d-\a_1 = 2\om_1+5\om_2 \\
\hline
\om_2 + \om_l & 2\d-\a_2 = \om_1 + \om_3 + 2\om_l\\
\hline
\om_2 + \om_4 & 2\d-\a_2 = \om_1 + \om_3 + 2\om_4\\
\hline
\end{array}
\]
For example, when $\d$ has a nonzero coefficient of $\om_1$, then $W=V_X(\d)$ contains vectors of
 weights $\d$ and 
$\d-\a_1$, and the wedge product of these vectors affords the weight $2\d-\a_1$, which is not subdominant to 
any other
 weight of $V = \wedge^2W$; therefore $V_X(2\d-\a_1)$ occurs as a summand of $V$.

For the case $\l=\l_3$ we have  $\d = \om_1+\om_{l+1}$ and $V = \wedge^3W$, and we take $\nu = 3\d-\a_1-\a_{l+1} = 
\om_1+\om_2+\om_l+\om_{l+1}$ (afforded by a wedge of weight vectors $\d \wedge (\d-\a_1) \wedge (\d-\a_{l+1})$); 
and 
for the last case $\l = 3\l_1$, 
we take $\nu = 3\d$. 

Before turning to Table \ref{TAB2}, consider the cases $\delta=2\omega_1$ or $\omega_2$ and 
$\lambda=\lambda_i$, for $1<i<n$, so $V = \wedge^i(W)$. The cases $\delta = 2\omega_1$ and $l = 1,2$ and 
$\delta=\omega_2$ with 
$l = 2,3$ can  be checked directly, yielding the exceptions in the table as well as a base case for induction. So 
now assume $l>2$ if $\delta=2\omega_1$ and $l>3$ if $\delta=\omega_2$. Let $Q_XL_X$ be the maximal parabolic 
subgroup defined in Chapter \ref{notation}.  Then $V^1 = \wedge^i(\delta)$, 
provided $i \le  \frac{(l+1)(l+2)}{2}-1$, respectively $\frac{l(l+1)}{2} -1$, for $\delta=2\omega_1$, resp. 
$\omega_2$. So in this case induction shows that some irreducible summand has highest weight 
with at least 2 nonzero labels, if 
$i < \frac{(l+1)(l+2)}{2} -1$, respectively $\frac{l(l+1)}{2} -1 .$ Then the corresponding irreducible summand of 
$\wedge^i(W) $ also has at least 2 nonzero 
labels.  So now assume $i \geq \frac{(l+1)(l+2)}{2} -1$, respectively $\frac{l(l+1)}{2} -1 .$ Then 
$n= \frac{(l+2)(l+3)}{2} -1$, respectively $\frac{(l+2)(l+1)}{2} -1 $, and $V^*$ has highest weight $\lambda_{n-i+1}$. 
But here $n-i+1 < \frac{(l+1)(l+2}{2} -1$, respectively $\frac{l(l+1)}{2} -1$, and hence $V^*\downarrow X$ has an irreducible 
summand with at least two nonzero labels,  
and so the same must hold for $V\downarrow X$. 

Now turn to the configurations of Table ~\ref{TAB2}, and note that the following cases are straightfoward verifications:

\[
\begin{array}{|lll|}
\hline
\d & \l&\nu \\
\hline
c\om_i & \l_1+\l_n&c\om_i+c\om_{l+2-i}\;(i\ne\frac{l+2}{2}) \\
                       & & 2\d-\a_i = \om_{i-1}+2c\om_i+\om_{i+1}\;(i=\frac{l+2}{2}) \\
                       
\hline
3\om_1 & \l_1+\l_2&3\d-\a_1 = 7\om_1+\om_2 \\ 
\hline
\om_3& \l_1+\l_2&3\d-\a_3 = \om_2+\om_3+\om_4 \\ 
\hline
\end{array}
\]

For the remaining cases in Table~\ref{TAB2}, we have $\delta = 2\omega_1$ or $\omega_2$. For the weights  
$a\lambda_1+\l_2$, $a\leq3$, $\l_2+\l_3$, $\l_2+\l_{n-1}$, and 
$\l_1+\l_{n-1} = (\lambda_2+\lambda_n)^*$, we first note that $\wedge^2(\delta) = 2\omega_1+\omega_2$, 
respectively $\omega_1+\omega_3$. 
Now we argue inductively, looking at the action of $L_X$ on $V^1$, as in the above
considerations for $\l = \l_i$. For the weights $\l_2+\l_{n-1}$ and $\l_2+\l_n$, we will use the result 
established for $\l=\l_2$. For the other weights, the base cases for $l\leq3$ are easily handled using Magma.
For $\l = a\l_1+\l_n$, we use that 
$V\downarrow X = (S^a(\d)\otimes \d^*)/S^{a-1}(\d)$ to see that $\nu = a\d+\d^*$ occurs as the highest weight of an 
irreducible summand. Then $\nu$ has the required property unless $\d = \om_2$ and $l=2$, where one can directly check
the existence of an appropriate summand.

For Table~\ref{TAB2}, it remains to consider the weights of the form $\l_1+\l_i$, for $3\leq i\leq n-2$. 
We first note that for $r_0+2\leq i\leq n-2$, $V^1 = 
V_{L_X'}(\d)\otimes V_{L_X'}(\omega_i)$ for some $1\leq i\leq l$. Since there exists an
irreducible summand of this module having two nonzero labels, this corresponds to a summand of 
$V\downarrow X$ having two nonzero labels as well. Now for $3\leq i\leq r_0+1$, we use 
induction on $l$.
It is straightforward to check that the result holds for $\d = 2\om_1$, $l=1,2$, and for $\d = \om_2$, $l=2,3$. 
 So assume $l\geq 3$, when $\d = 2\omega_1$ and $l\geq 4$ when $\d=\om_2$. If $i\leq r_0$, the induction hypothesis 
together with the previously considered cases show that $V^1$ has an irreducible summand with two nonzero labels, 
which again implies the result.
If $i=r_0+1$, we replace $V$ by $M=V^*$; then  $M^1 = \wedge^{l+2}(\delta)$ or $\wedge^{l+1}(\delta)\otimes \omega_2^*$, respectively. Then using the result for $\lambda_i$ and recalling that $l\geq 3$, respectively $l\geq 4$, we see that $M^1$ has a summand with two nonzero labels and hence we get the result for $V$ as well.

We now turn to Table~\ref{TAB3}, and recall that we have already handled the cases $\d = 2\om_1$ or $\om_2$ 
when $\lambda=\l_i$.   The following cases can be easily handled, sometimes working in a group of a fixed rank in order 
to 
establish the existence of the weight $\nu$:
\[
\begin{array}{|lll|}
\hline
\d & \l&\nu \\
\hline
c\om_1, c=2 & a\l_2,a=2,3&2a\om_1+a\om_2 \\
\ \ \ \ \ \         c\geq 2               &a\l_1,ac>4&a\d-2\alpha_1 = (ac-4)\om_1+2\om_2\\
\ \ \ \ \ \         c>2                    &\lambda_2&(2c-2)\omega_1+\omega_2\\
\hline
\om_2,l\geq 3&a\l_1,a\geq3 & a\d-\alpha_1-2\alpha_2-\alpha_3 = (a-2)\om_2+\om_4 \\
\hline
\om_i, 3\leq i\leq l-1&a\l_1, a\geq 2&\omega_{i-2}+(a-2)\omega_i+\omega_{i+2}\\
\hline
\om_2&a\l_2,a=2,3 & a\om_1+a\om_3 \\
\hline
\om_i,i=3,4&\l_4&4\d-\alpha_{i-1}-3\alpha_i-\alpha_{i+1} = \\
i\leq\frac{l+2}{2}&&\ \ \ \om_{i-2}+\om_{i-1}+\om_{i+1}+\om_{i+2}\\
\hline
\om_i, i=3,4,5,6&\l_3&3\d-\alpha_{i-1}-2\alpha_i = \om_{i-2}+2\om_{i+1}\\
i\leq\frac{l+2}{2}&&\\
\hline
c\om_i, i>1, c>1&2\l_1&2\d-2\alpha_i = 2\om_{i-1}+(2c-4)\om_i+2\om_{i+1}\\
\hline
c\omega_i, i>1, c\geq 1&\l_2&2\d-\alpha_i = \om_{i-1}+(2c-2)\om_i+\om_{i+1}\\
\hline
\end{array}
\]

To complete the consideration of Table~\ref{TAB3}, for the  cases where  
$\d=\om_i$, $i=2,3,4,5$ and $\l=3\l_1$, we work in a group of large enough fixed rank (here when $\delta=\omega_2$, we have $l\geq3$ and in any case $l\geq 2i-2$)  and produce the desired
summand. Finally, we must consider the weights $\delta=c\omega_1$, $3\leq c\leq 6$, respectively, 
$c=3,4$, $c=3$, and $\lambda=\lambda_3$, respectively $\lambda_4$, $\lambda_5$. Then we use the weights $\nu = 3\delta-3\alpha_1$, respectively $\nu = 4\delta-4\alpha_1-\alpha_2$, $\nu = 5\delta-7\alpha_1-\alpha_2$. 

Finally, we turn to Table~\ref{TAB4}. If $l\geq 4$, i.e., if ${\rm rank}(X)\geq 5$, this is a straightforward Magma check.
For $X=A_3$, we refer to the proof of Lemma~\ref{lnu2}  for the restrictions of the $\l_i$ to $T_X$, $1\leq i\leq 5$. 
Taking $\nu=\l\downarrow T_X$, we reduce to the weights (up to duals) $\{a\lambda_3,a\lambda_1+\lambda_5\}$.
For the first family, we may assume $a\geq2$, in which case, again using the restrictions of the $\lambda_i$ and the 
$\beta_i$ to $T_X$, we deduce that $\l-\beta_3$ affords the highest weight of an irreducible $X$-summand, affording $\nu=2\om_1+(2a-2)\om_3$. For the final family, $\l=a\lambda_1+\lambda_5$,  $V = (S^a(\d)\otimes \d^*)/S^{a-1}(\d)$
and it is easy to see then that $\l\downarrow T_X$ affords $\nu = \om_1+a\om_2+\om_3$, the highest weight of an 
irreducible summand.

So finally we turn to the cases where $X=A_4$, and $\d=\omega_2$. Here we may choose an embedding of $X$ in $Y$ such that 
$\lambda\downarrow T_X$ affords the highest weight of an irreducible summand of $V\downarrow X$ and $\beta_1$ restricts to 
$\alpha_2$ and $\beta_2$ restricts to $\alpha_1$. Indeed, let $P$ be the parabolic subgroup of $X$,
containing $T_X$ and the Borel subgroup defined by $\Pi(X)$ and whose Levi factor has simple roots 
$\{\alpha_1,\alpha_2\}$, Then, we may assume $P$ lies in a parabolic subgroup 
$R$ of $Y$, containing $T_Y$, and corresponding to the base $\Pi(Y)$. Then considering the action of both groups on 
$L_X(\delta)$  shows that a Levi factor of $R$ has simple roots
$\{\beta_i\ i=1,2,4,5,7,8\}$. Then again considering the action on the natural module for $Y$, we see that 
we have the restrictions $\beta_j\downarrow T_X = \alpha_1$ for $j=2,4,7$ and $\beta_j\downarrow T_X = 
\alpha_2$ for $j=1,5,8$. 
Then $\lambda_2\downarrow T_X = \omega_1+\omega_3$, which settles the cases in rows 2-4 of the $A_4$ section
 of Table~\ref{TAB4}. For $\l=a\l_1+\l_9$, as in the previous paragraph, 
we use that $V\downarrow X = (S^a(\delta)\otimes \delta^*)/S^{a-1}(\delta)$ and take $\nu=a\om_2+\om_8$. For the last two rows of 
the $A_4$ part of Table~\ref{TAB4}, we can use Magma.\hal 


We note the following corollary of the proof.

\begin{cor}\label{cor_giantlemma} Let $\d = 2\om_1$ or $\om_2$ and $\l=\l_1+\l_i$, for $i\geq 2$. Then $V\downarrow X$ has an irreducible summand with highest weight $\nu$ such that $L(\nu)\geq 2$.

\end{cor}

We finish the section with a lemma on dimensions that will be needed later.

\begin{lem}\label{diminequal}
Let $X = A_{l}$ with $l \ge 2$, and let $2\le j\le l$ and $a,b \ge 1$. Then 
$\dim V_X(a\om_1+b\om_j) > \dim V_X(b\om_{j-1})+2$ for the following values of $b,j$:
\begin{itemize}
\item[{\rm (i)}] $b=1$, any $j$
\item[{\rm (ii)}] $j=2$, any $b$
\item[{\rm (iii)}] $b=2, j=3$.
\end{itemize}
\end{lem}

\pf (i) The dimension of $ V_X(a\om_1+\om_j)$ is at least the number of conjugates of  $a\om_1+\om_j$ under the action of the Weyl group $S_{l+1}$, which is $\frac{(l+1)!}{(j-1)!(l-j+1)!}$. One checks that this is greater than ${l+1 \choose j-1}+2$, which is equal to $\dim V_X(\om_{j-1})+2$, as required. 

(ii) Let $j=2$. The Weyl dimension formula shows that $\dim V_X(a\om_1+b\om_2)$ is divisible by
\[
d(a,b):=  \frac{b+1}{1} \frac{b+2}{2}\cdots   \frac{b+l-1}{l-1}\cdot  \frac{a+b+l}{l}.
\]
Hence 
\[
\dim V_X(a\om_1+b\om_2) - \dim V_X(b\om_1) \ge d(a,b) - {b+l\choose l}.
\]
The right hand side is greater than 2 unless either $b=1$ or $l \le 3$. In the latter cases the conclusion is easily checked.

(iii) Since the dimension of  $\dim V_X(a\om_1+2\om_3)$ grows with $a$, we can assume that $a=1$. So we need to show that $\dim V_X(\om_1+2\om_3) > \dim V_X(2\om_2)+2$. Using the Weyl dimension formula we have 
\[
\dim  V_X(2\om_2) = \frac{(l+2)(l+1)^2l}{12}.
\]
On the other hand, for $l\ge 4$, $V_X(\om_1+2\om_3)$ has a subdominant weight $\om_1+\om_2+\om_4$, and the number of conjugates of this under the Weyl group is $\frac{1}{2}\frac{(l+1)!}{(l-3)!}$, which is greater than $\dim  V_X(2\om_2)$. Finally, for $l=3$ the result holds by the dimension formula.
 \hal

\chapter{The case $X = A_2$}\label{casea2}

In this chapter we prove Theorem \ref{MAINTHM} in the case where $X = A_2$. Notation will be  as in  Chapter \ref{notation}. In particular $\Pi(X) = \{\a_1,\a_2\}$ and $\Pi(Y) = \{\b_1, \ldots , \b_n\}$  are fundamental systems of positive roots for $X$ and $Y$, respectively,  with corresponding fundamental weights $\{\om_1,\om_2\}$  and $\{\l_1, \ldots , \l_n\}.$ 

\vspace{4mm}
\section{Case $\delta = rs$ with $r,s>0$}

Let $X = A_2$, let $\d$ be the dominant weight $rs: = r\om_1+s\om_2$ and let $W = V_X(\d)$. Assume $r \ge s \ge 1$.
Then $X$ embeds in $SL(W)$. In this section we determine all irreducible $Y$-modules $V_Y(\l)$ such that $V_Y(\l)\downarrow X$ is multiplicity-free.

\begin{thm}\label{rgreater}
Let $X = A_2$, $\d = rs$ and assume $r\ge s\ge 1$. Let $W = V_X(\d)$ and take $X < Y = SL(W)$ as above. Suppose $\l$ is a dominant weight for $Y$ such that $\l$ is not $\l_1$ or its dual. Then $V_Y(\l)\downarrow X$ is multiplicity-free if and only if one of the following holds, where $\l$ is given up to duals:

{\rm (i)} $s=1$ and  $\l = \l_2$ or $2\l_1$

{\rm (ii)} $r=s=1$ and $\l = 3\l_1$ or $\l_3$ 

{\rm (iii)} $r=s=2$ and $\l=\l_2$.
\end{thm}

 \subsection{Preliminaries}\label{prel}

First note that all the examples listed in (i), (ii) and (iii) of Theorem \ref{rgreater} are indeed multiplicity-free, as was shown in Chapter \ref{MFfam}. 

We shall need the following in the proof.

\begin{lem}\label{A2rssquares}  Assume $X = A_2$ and assume
$r \ge s \ge 2.$  Then $\Sym^2(rs)$ is not MF and
$\wedge^2(rs)$ is MF only if $r = s = 2.$

\end{lem}

\pf  Lemma \ref{wedgesquareab}(ii) shows that $\Sym^2(rs)$ is not MF and 
Lemma \ref{wedgesquareab}(i) shows that $\wedge^2(rs)$ is not MF provided
$r \ge 3.$  Finally a Magma computation shows that  $\wedge^2(22)$ is MF.  \hal

\vspace{4mm}
Let $X = A_2$, $\om = rs$ and assume $r\ge s\ge 1$. Let $W = V_X(\om)$ and take $X < Y = SL(W).$ As in the statement of Theorem \ref{rgreater}, suppose $V_Y(\l)\downarrow X$ is multiplicity-free.
Write $V:= V_Y(\l)$, and let $V\downarrow X = V_1+\cdots +V_t$ where the $V_i$ are distinct irreducible $X$-modules.

As in Chapter \ref{levelset}, let $P_X = Q_XL_X = Q_XL_X'T$ be a maximal parabolic in $X$ with $L_X'$ corresponding to the root $\a_1$, and embed it in a parabolic $P_Y=Q_YL_Y$ of $Y$ satisfying the conditions of Lemma \ref{parabembed}. By Lemma \ref{rectangle},
\begin{equation}\label{factors}
L_Y' = A_r+A_{2r+1}+ A_{3r+2} +\cdots + A_{2s+1}+A_s.
\end{equation}
 In the notation of Chapter \ref{notation} we have $A_r = C^0, A_{2r+1} = C^1, \ldots , A_s = C^k$, with $k=r+s+1$.
Moreover, $L_X'$ projects to each factor as a sum of irreducibles of highest weights as in the array presented after 
Lemma \ref{rectangle}. By \cite[3.18]{book}, the 
$S_X$-labelling of the factors of $L_Y'$ is $22\ldots  2$ (for the $A_r$ factor), $2020\ldots  2$ (for the $A_{2r+1}$ factor), 
$20200200\ldots  202$ (for the $A_{3r+2}$ factor) and so on.
As in Chapter \ref{levelset}, define the $d^{th}$ level in $V$ to be 
\[
V^{d+1}(Q_Y) : = [V,Q_Y^d]/[V,Q_Y^{d+1}],
\]
and likewise define $V_i^{d+1}(Q_X) : = [V_i,Q_X^d]/[V_i,Q_X^{d+1}]$. Recall also $V^i = V^i(Q_Y)\downarrow L_X'$.
By Proposition \ref{induct}, 
\begin{equation}\label{levone}
V^1 = \sum_{i,n_i=0} V_i^1(Q_X)
\end{equation}
 is multiplicity-free. 

Recall that $L_X' = A_1$, and let $m$ denote the highest weight appearing in $V^1$.

\begin{lem}\label{bound}
Suppose level $1$ of $V$, namely $V^2$, has a composition factor of highest weight $t$.
\begin{itemize}
\item[{\rm (i)}] The composition factor $t$ has multiplicity at most $3$ in $V^2$.
\item[{\rm (ii)}] If $t> m+1$ then composition factor $t$  has multiplicity at most $1$ in $V^2$.
\item[{\rm (iii)}] If $t= m+1$ then composition factor $t$  has multiplicity at most $2$ in $V^2$.
\item[{\rm (iv)}] Suppose the second highest weight of a composition factor of $V^1$ is at most $m-4$. Then the composition factors $m-1$ and $m-3$ appear with multiplicity at most $2$ in $V^2$.
\end{itemize}
\end{lem}

\pf By Proposition \ref{induct}(ii) we have 
\begin{equation}\label{levonemore}
V^2 = \sum_{i,n_i=0} V_i^2(Q_X) + \sum_{j,n_j=1} V_j^1(Q_X),
\end{equation}
where the second sum on the right hand side is multiplicity-free. From Theorem \ref{LEVELS}, it follows that each composition factor $x$ in the multiplicity-free sum in (\ref{levone}) corresponds to at most two composition factors in the first sum in (\ref{levonemore}), with highest weights among $x\pm 1$. Therefore any given composition factor $t$ can appear no more than twice in the first sum, hence no more than three times in $V^2$, which proves (i). Moreover $m+1$ can appear at most once in the first sum in (\ref{levonemore}), and no weight greater than $m+1$ can appear at all, which gives (ii) and (iii). Finally, under the assumption in (iv), $m-1$ and $m-3$ can appear only once in the first sum, and part (iv) follows. \hal

\vspace{4mm}
We conclude this section by dealing with the case where the highest weight $\l$ has no nonzero labels on any factor of $L_Y'$, which is to say that $V^1(Q_Y)$ is trivial.

\begin{lem}\label{triv} Assume $(r,s) \ne (1,1)$. 
If $V^1(Q_Y)$ is trivial, then $s=1$, $Y=A_n$ and $\l = \l_{n-1}\, (= \l_2^*)$, as in Table  $\ref{TAB1}$ of Theorem $\ref{MAINTHM}$ .
\end{lem}

\pf In this case $m=0$. Pick $\g \in \Pi(Y)\setminus \Pi(L_Y')$ such that $\la \l, \g \ra \ne 0$. Suppose that $\g \ne \b_{n-1}$. Then $\g$ is adjacent to two factors $A_c,A_d$ of $L_Y'$ with $c,d>1$.
Then the natural modules for the two factors adjacent to $\g$ have $L_X'$-summands $x$ and $(x+1)\oplus (x-1)$, for some $x\ge 2$. Hence $\l-\g$ affords a summand $x \otimes ((x+1)\oplus (x-1))$ in $V^2$. This contains $(2x-1)^2$, contradicting Lemma \ref{bound}(ii) (since $m=0$). 

Hence $\g = \b_{n-1}$, which implies that $Y=SL(W)= A_n$, $s=1$ (so $r>1$) and $\l = c\l_{n-1}$ for some $c \ge 1$. Suppose $c\ge 2$. Replace $V$ by the dual $V^* = V_Y(c\l_2)$. The fact that $V^1$ is multiplicity-free implies that $r=2$ by \cite{LST}. Hence $\d = 21$ and now the highest weight of $L_X'$ on $V^1(Q_Y)$ is $m= 2c$. Then $\l-\a_3-\a_2$ affords the level 1 summand 
$(1,c-1) \otimes 10000$ for $A_2A_5 \le L_Y'$, for which the restriction to $L_X'$ contains $(2c\oplus (2c-2))\otimes (3\oplus 1)$, which contains $(2c+1)^3$. This contradicts Lemma \ref{bound}(iii). Hence $c=1$, completing the proof. \hal 

 \subsection{Proof of Theorem \ref{rgreater}}

We now embark on the proof of Theorem \ref{rgreater}. Adopt the assumptions of the previous section, assuming that $r\ge s\ge 1$,  that $Y = SL(W) = A_n$, so that the factors of $L_Y'$ are given by (\ref{factors}), and that $V_Y(\l)\downarrow X$ is MF. We also assume until Lemma \ref{oneone} that
\begin{equation}\label{not11}
(r,s) \ne (1,1).
\end{equation}

In view of Lemma \ref{triv}, we assume that  $V^1(Q_Y)$ is nontrivial, so that $\l$ has a nonzero label on at least one of the factors of $L_Y'$ in (\ref{factors}).

\begin{lem}\label{lab}
{\rm (a)} Suppose $\l$ has a nonzero label on a factor $A_d$ of $L_Y'$, where $d\ne r,s$. Then one of the following holds:

{\rm (i)} the $\l$-labelling of $A_d$ is $10\ldots  0$ or $00\ldots  1$

{\rm (ii)} $d=3, s=1$ and the $\l$-labelling is $200$ or $002$.

\no Moreover, there is no other factor of $L_Y'$ with a nonzero $\l$-label.

\vspace{2mm}
 {\rm (b)} Suppose $\l$ has nonzero labels on both of the factors $A_r$, $A_s$ of $L_Y'$. Then these labels are
$10\ldots  0$ or $00\ldots  1$.

 \vspace{2mm}
 {\rm (c)} If $\l$ has a nonzero label on just one of the factors $A_r$, $A_s$, this label is given by \cite[Thm.1]{LST}.
\end{lem}

\pf Since $V^1$ is multiplicity-free, the weight $m$ appears exactly once, and for each $c\ge 1$ the weight $m-2c$ appears with multiplicity at most $c+1$. 

(a) Assume as in (a) that $\l$ has a nonzero label on a factor $A_d$ of $L_Y'$, where $d\ne r,s$. The $S_X$-labelling of $A_d$ is $2020\ldots $. If $\l$ has a nonzero label over a root $\b_i$ with $S_X$-label 0, then the weight $m$ is afforded by both $\l$ and 
$\l-i$ (where $\l-i$ denotes $\l-\b_i$), so $m^2$ appears in $V^1$, a contradiction. If $\l$ has a nonzero label over a root $\b_i$ with $S_X$-label 2 which is not an end node of $A_d$, and $\b_i$ adjoins $\b_j,\b_k$, then $(m-2)^3$ appears in $V^1$ as it is afforded by $\l-i, \l-ij, \l-ik$, again a contradiction. 

Hence $\l$ can have nonzero labels only on the two end nodes of $A_d$; moreover, only one such label is possible, since otherwise $m-2$ appears in $V^1$ with multiplicity greater than 2.
Now assume $\l$ has a label $c>1$ over an end node $\b_i$ of $A_d$. 
If $d\ge 4$ let $ijkl$ be adjoining nodes of $A_d$ (with $S_X$-labelling $2020 \ldots  $); then the weight $m-4$ appears with multiplicity 4 in $V^1$ as it is afforded by $\l-i^2, \l-i^2j,\l-ijk, \l-ijkl$. Hence $d=3$, which forces $s=1$ (see the array after Lemma \ref{rectangle}). Let $ijk$ be the nodes of $A_3$. If $c>2$ then $m-6$ is afforded by $\l-i^3,\l-i^3j,\l-i^3j^2,\l-i^3j^3,\l-i^2jk$. Hence $c=2$.

We have now shown that (i) or (ii) of part (a) holds. For the final assertion of (a), observe that if $\l$ has a nonzero label over another factor of $L_Y'$, then $m-2$ has multiplicity 3 in $V^1$, so this cannot occur.

(b)  The weights $m-2$ and $m-4$ appear in $V^1$ with multiplicities at most $2$ and $3$, respectively. Also each of the fundamental roots $\b_i^0$ of $A_r$, respectively $A_s$, restricts to $S_X$ affording weight 2.  Let $x, y$ be the highest weights of the largest compositions factors of $\mu^0 \downarrow L_X'$ and $\mu^k\downarrow L_X’,$ respectively.

We claim that $\mu^0$ is the natural module or its dual. Suppose otherwise and recall that $r > 1.$  Then either there is an end node with label strictly greater than 1, a node with nonzero label which is not an end node, or two nodes with nonzero labels.  In the first two cases weights $x, x-2,(x-4)^2 $ occur in the $L_X'$ summand of $V^1$ afforded by $\mu^0$, and in the third case $x, (x-2)^2, (x-4)^2$ occur. Applying a similar analysis to the $A_s $ factor we have a contradiction in the third case and in the other cases we must have $s = 1$ and $\mu^k = \l_1^k.$

We now argue using the methods described in Section \ref{levanal}. 
If $\mu^0 = (a0\ldots  0)$ with $a > 1$, then  $V^1 = S^a(r) \otimes 1 = ((ar)\oplus (ar-4) \oplus  \cdots) \otimes (1)$ and $V^2 \supseteq S^{a-1}(r) \otimes ((r+1) \oplus  (r-1)) \otimes (1)$; this contains $(ar-2)^4$.  At most two such factors can arise from $V^1$, so this contradicts Corollary \ref{cover}. 

Next assume $\mu^0 = a\l_i^0$ for $i > 1$ and consider the dual module $V^*$.  By part (a),
$(V^*)^1 = r$ and there is a node $\g \in \Pi(Y)\setminus \Pi(L_Y')$  which is not adjacent to $A_r$ and such that $\la \l^*, \g \ra \ne 0$.   The natural modules for the two factors adjacent to $\g$ have $L_X'$-summands $x$ and $(x+1)\oplus (x-1)$, for some $x\ge 2$. Hence $\l-\g$ affords a summand $r \otimes x \otimes ((x+1)\oplus (x-1))$ in $V^2$. This contains $(r+ 2x-1)^2$, contradicting Lemma \ref{bound}(ii). 

We have now proved that  $\mu^0$ is the natural or dual module.  Applying this to $V^*$ we obtain the assertion.

(c) This follows from the fact that $V^1$ is multiplicity-free. \hal

\vspace{4mm}
We now handle separately cases (a), (b) and (c) of Lemma \ref{lab}, starting with (a).

\begin{lem}\label{casealab}
Assume that case {\rm (a)(i)} of Lemma $\ref{lab}$ holds. Then $\l = \l_k$ for some $k$.
\end{lem}

\pf In this case there is a factor $A_d$ of $L_Y'$ with $d\ne r,s$ such that the $\l$-labelling of $A_d$ is $10\ldots  0$ or $00\ldots  1$, and $\l$ has zero labels on all other factors of $L_Y'$. 

The conclusion of the lemma will follow if we show that $\la \l,\g\ra = 0$ for every root $\g \in \Pi(Y) \setminus \Pi(L_Y')$. So suppose $\la \l,\g\ra \ne 0$ for some $\g \in \Pi(Y) \setminus \Pi(L_Y')$. Then $\l-\g$ affords a summand $S$ of the first level $V^2(Q_Y)\downarrow L_Y'$. We investigate the possibilities for $S$.

Assume first that $\g$ is not adjacent to the factor $A_d$ of $L_Y'$. Let $A_u,A_v$ be the factors adjacent to $\g$. Then $S$ is isomorphic to the tensor product of natural or dual modules for $A_d,A_u,A_v$. Each of these natural modules restricts to $L_X'$ as a row of the $rs$-array exhibited after Lemma \ref{rectangle}. Hence for some $x \ge 1$, $S\downarrow L_X'$ contains
\[
x \otimes ((x+1)\oplus (x-1)) \otimes (m\oplus (m-2)).
\]
This contains the summand $2x+m-1$ with multiplicity at least 4, contradicting Lemma \ref{bound}(i).

Now assume $\g$ is adjacent to $A_d$, and let $A_e$ be the other factor of $L_Y'$ adjacent to $\g$. Then as a module for $A_dA_e$, $S$ is either $20\ldots  0 \otimes M,$  $0\ldots  02 \otimes M$, or $10\ldots  01 \otimes M$, where $M$ is the natural or dual module for $A_e$. Let $x$ be the highest weight of $M\downarrow L_X'$. 

Now we reason as in Section \ref{levanal}. 
If $S = 20\ldots  0 \otimes M$ or $0\ldots  02 \otimes M$, then $S\downarrow L_X'$ contains $S^2(m\oplus (m-2)) \otimes x$, which contains $(2m+x-4)^4$ unless $m=2$, $x=1$. Hence the latter holds by Lemma \ref{bound}. Then $d=3$ and $e=s=1$. Now $\g = \b_{n-1}$ and 
there is an additional irreducible in $V^2(Q_Y)$ for which the highest weight is  $\lambda - \gamma - \b_{n-2}$ restricted to the maximal torus of $L_Y'$, and hence the irreducible $010 \otimes 1$ appears as a further level 1 summand for $L_Y'$. Hence 
$V^2$ contains 
\[
(S^2(2\oplus 0) \otimes 1) + (\wedge^2(2\oplus 0)\otimes 1)
\]
which contains $1^5$, contradicting Lemma \ref{bound}.

If $S = 10\ldots  01 \otimes M,$ then $S\downarrow L_X'$ contains $((m\oplus (m-2)) \otimes (m\oplus (m-2)))^\dagger \otimes x$ (where $\dagger$ indicates that one trivial composition factor should be omitted). This contains $(2m+x-2)^4$, contrary to Lemma \ref{bound}. This completes the proof of the lemma. \hal

\begin{lem}\label{caseano}
Case {\rm (a)(i)} of Lemma $\ref{lab}$ does not occur.
\end{lem}

\pf Suppose false. By Lemma \ref{casealab} we have $\l = \l_k$, where $\b_k$ is an end node of the $A_d$ factor of $L_Y'$. Let $\g$ be the node in $\Pi(Y)\setminus \Pi(L_Y')$ adjoining $\b_k$, and let $A_e$ denote the factor of $L_Y'$ adjoining $\g$ (with $e\ne d$). Then $\l-\b_k-\g$ affords a summand $S:=00\ldots  10 \otimes 10\ldots  0$ (or the dual)  of level 1 for $L_Y'$. Let $x$ be the highest weight of the natural module for $A_e$ restricted to $L_X'$, so that $x = m\pm 1$. Then $S\downarrow L_X'$ contains $\wedge^2(m\oplus (m-2)) \otimes x$. If $m\ge 3$ this contains $(2m+x-6)^4$, contradicting Lemma \ref{bound}(i). So suppose $m=2$. If $x=3$ then  $\wedge^2(m\oplus (m-2)) \otimes x$ contains $5^2$, contrary to Lemma \ref{bound}(ii). 

This leaves the case where $m=2$, $x=1$. Then also $s=1$, $d=3$ and $k=n-2$. Let $\eta = \b_{n-5}$ be the other node (apart from $\g$) adjoining the $A_d$ factor of $L_Y'$, and let $A_f$ be the factor of $L_Y'$ adjoining $\eta$ (with $f=5$). Then $\l-\eta-\b_{n-4}-\b_{n-3}-\b_{n-2}$ affords the level 1 summand $00001 \otimes 000$ for $A_5A_3$. Now we see that level 1 for $L_X'$ contains
\[
(\wedge^2(2\oplus  0) \otimes 1) + (3 \oplus  1)
\]
which contains $3^3$, contradicting Lemma \ref{bound}(iii). \hal

\begin{lem}\label{caseaiino}
Case {\rm (a)(ii)} of Lemma $\ref{lab}$ does not occur.
\end{lem}

\pf Assume false, so that $d=3, s=1$ and the $\l$-labelling of the $A_3$ factor of $L_Y'$ is $200$ or $002$.
Let $\g = \b_{n-1}$ and $\eta = \b_{n-5}$ be the nodes adjoining the $A_3$ factor. 

Observe that $m=4$. If the $\l$-labelling of $A_3$ is $200$, then $\l-\eta-\b_{n-4}$ affords a level 1 summand $00001\otimes 110$ for $A_5A_3$; and if the $\l$-labelling is $002$ then $\l-\g-\b_{n-2}$ affords level 1 summand $011\otimes 1$ for $A_3A_1$. In either case the restriction of this summand to $L_X'$ contains $3^4$, contradicting Lemma \ref{bound}(i). \hal

\vspace{4mm}
We now move on to cases (b) and (c) of Lemma \ref{lab}.

\begin{lem}\label{casesbc}
In cases {\rm (b)} and {\rm (c)} of Lemma $\ref{lab}$, all $\l$-labels on $\Pi(Y)\setminus \Pi(L_Y')$ are $0$.
\end{lem}

\pf 
In these cases Lemmas \ref{caseano} and \ref{caseaiino} imply that the only nonzero $\l$-labels on $L_Y'$ are on the $A_r$ and $A_s$ factors in (\ref{factors}). 
Suppose that there is a nonzero $\l$-label on some $\g \in \Pi(Y)\setminus \Pi(L_Y')$. 

If $\g$ is not adjacent to either factor $A_r$ or $A_s$, let $A_d,A_e$ be the factors of $L_Y'$ adjacent to $\g$. Then $L_X'$ acts on the natural modules for these factors as $x\oplus (x-2)\oplus \cdots$, respectively $(x-1)\oplus (x-3)\oplus \cdots$ for some $x \ge 3$. Now $\l-\g$ affords a level 1 summand restricting to $L_X'$ as
\[
m \otimes (x\oplus (x-2)) \otimes ((x-1)\oplus (x-3)).
\]
This contains $(m+2x-3)^4$, contradicting Lemma \ref{bound}(i).

Now suppose $\g$ is adjacent to $A_x$ with $x \in \{r,s\}$. Using Lemma \ref{p.29} we see that  $\l-\g$ affords a level 1 summand restricting to $L_X'$ as $(m+x)\otimes ((m+1)\oplus (m-1))$. This contains $(m+2x-1)^2$, so by Lemma \ref{bound} we have $m+2x-1 \le m+1$, and so $x=1$. Therefore $s=1$ and $\g$ is adjacent to $A_s$ (so $\g = \b_{n-1}$).   In case (b) of Lemma \ref{lab}, $m = r+1$. If the $\lambda$-label on 
$A_s$ is nonzero, then  $\l-\g$ affords a level 1 $L_X'$-summand $(m+1)\otimes (2\oplus 0)$, 
while $\l-\g-\b_n$ affords $(m-1)\otimes (2\oplus 0)$. Hence level 1 for $L_X'$ contains 
$(m+1)^3$, contradicting Lemma \ref{bound}(iii).  So the $\lambda$-label on $A_s$ is zero and we are in case (c). Now 
$\l-\g$ contains  $m\otimes (2+0)\otimes 1$ which also contains $(m+1)^3$, producing the 
same contradiction.
\hal

\begin{lem}\label{caseblab}
Assume case {\rm (b)} of Lemma $\ref{lab}$ holds. Then $\l$ is $\l_1+\l_n$, $\l_r+\l_n$, $\l_1+\l_{n-s+1}$ or
$\l_r+\l_{n-s+1}$.
\end{lem}

\pf In this case the $\l$-labellings of the $A_r$ and $A_s$ factors in (\ref{factors}) are $10\ldots  0$ or $00\ldots  1$. Hence by Lemma \ref{casesbc}, $\l$ is one of the four possibilities in the statement of the lemma. \hal

\begin{lem}\label{caseblabno}
Case {\rm (b)} of Lemma $\ref{lab}$ does not occur.
\end{lem}

\pf Suppose false. Then $\l$ is one of the four possibilities in Lemma \ref{caseblab}.

Assume first that $\l = \l_1+\l_n$. Then $V_Y(\l)\downarrow X = (rs \otimes sr)^\dagger$. However for $r>s\ge 1$ this is not multiplicity-free by  Lemma \ref{compfactor}(ii).  

Now assume $\l = \l_r+\l_n$. We have $m=r+s$. Let $\g = \b_{r+1}$ be the root adjacent to the factor $A_r$. Then $\l-\b_r-\g$ affords the level 1 $L_X'$-summand 
\[
\wedge^2r \otimes ((r+1)\oplus (r-1)) \otimes s.
\]
This contains $(3r+s-3)^2$, so by Lemma \ref{bound}(ii) we have $3r+s-3 \le r+s+1$. Hence $r\le 2$ and so $(r,s) = (2,1)$ or $(2,2)$  (recall our assumption (\ref{not11}) that $(r,s)\ne (1,1)$). In this case the above level 1 summand is $\wedge^22 \otimes (3\oplus 1) \otimes s$, which contains $2^4$ or $3^4$ according to whether $s = 1$ or $2$. This contradicts Lemma \ref{bound}(i).

Next consider $\l = \l_1+\l_{n-s+1}$. If $r=s$ this is just the dual of the previous case, so assume $r>s$.
The dual $\l^* = \l_s+\l_n$ also has the property that $V_Y(\l^*) \downarrow X$ is multiplicity-free. Since $r>s$, by Lemma \ref{caseblab} this implies that $s=1$, so that $\l = \l_1+\l_n$. This case was dealt with above.

Finally, suppose $\l = \l_r+\l_{n-s+1}$. The argument for the $ \l_r+\l_n$ case forces $r\le 2$, hence also  $s \le 2$. The argument for the $\l = \l_r+\l_n$ above gives a contradiction here as well. \hal

\vspace{4mm}
At this point we are in case (c) of Lemma \ref{lab}, in which $\l$ has a nonzero labelling on just one factor $A_r$ or $A_s$ of $L_Y'$ in (\ref{factors}). If all nonzero labels are on the $A_r$ factor, replace $\l$ by its dual $\l^*$; then either all labels are on the $A_s$ factor of $L_Y'$ or we are in a previous case which has already been dealt with.
Hence we may assume that the nonzero $\l$-labelling is on the $A_s$ factor. We use the notation of Chapter \ref{notation}, so we write $\l_{A_s} = \sum_{i=1}^s c_i\l_i^k$ for the $\l$-labelling of $A_s$, etc.
 By Lemma \ref{casesbc}, the only nonzero $\l$-labels on $\Pi(Y)$ are those on the $A_s$ factor of $L_Y'$.

\begin{lem}\label{labas}
The possibilities for $\l_{A_s}$ are as follows, up to duals:
\begin{itemize}
\item[{\rm (i)}] $c\l^k_1$ (where $s\le 2$ if $c> 5$; $s\le 3$ if $c=4,5$; $s\le 5$ if $c=3$);
\item[{\rm (ii)}] $\l^k_2$ (where $s\ge 3$);
\item[{\rm (iii)}] $\l^k_1+\l^k_s$ (where $s\ge 2$);
\item[{\rm (iv)}] $\l^k_3$ (where $5\le s\le 7$);
\item[{\rm (v)}] $c\l^k_1+\l^k_2$ (where $s=2$, $c\ge 2$);
\item[{\rm (vi)}] $\l^k_1+\l^k_2$ (where $s=3$).
\end{itemize}
\end{lem}

\pf Since $V^1(Q_Y) = V_{A_s}(\l_{A_s})$ is multiplicity-free on restriction to $L_X'$, this follows immediately from the $A_1$ result in \cite{LST}. \hal

\vspace{4mm} 
We handle these cases one by one.

\begin{lem}\label{co1}
Suppose $\l_{A_s} = c\l_1^k$ or $c\l_s^k$. Then one of the following holds:

{\rm (i)} $s=1$, $c=2$ and $\l = 2\l_n = 2\l_1^*$

{\rm (ii)} $r=s=2$ and $\l = \l_{n-1} = \l_2^*$.

 In both cases $V_Y(\l)\downarrow X$ is multiplicity-free. 
\end{lem}

\pf First assume that $s\ge 2$ and $\l_{A_s} = c\l^k_1$. If $c=1$ and $s=2$ then $\l = \l_2^*$, and 
Lemma \ref{A2rssquares} shows that $V_Y(\l)\downarrow X$ is multiplicity-free if and only if $r=2$, as in conclusion (ii). So assume that $s\ge 3$ if $c=1$.

Observe that $m=cs$. The nonzero $\l$-label is over the root $\b_1^k$ of $A_s$. Let $\g_k$ be the adjacent root to $A_s$. Then $\l-\g_k-\b_1^k$ affords the level 1 $L_X'$-summand
\[
((s+1)\oplus (s-1)) \otimes V_{A_s}(c-1,1,\ldots ,0)\downarrow L_X'.
\]
This contains $((s+1)\oplus (s-1)) \otimes ((c+1)s-2)$, which contains the composition factor $(c+2)s-3$ with multiplicity 2. Hence by Lemma \ref{bound}(ii) we have $(c+2)s-3 \le m+1 = cs+1$, forcing $s\le 2$. Therefore $s=2$ (as we are assuming $s\ge 2$ at the moment), and so $c\ge 2$ by the observation in the first paragraph. Then 
$ V_{A_2}(c-1,1)\downarrow L_X'$ contains $2c \oplus (2c-2)$, so level 1 for $L_X'$ contains 
\[
(3\oplus 1) \otimes (2c \oplus (2c-2)).
\]
This contains $(2c+1)^3$, contradicting Lemma \ref{bound}(iii). 

Now assume that $\l_{A_s} = c\l^k_s$ and $s\ge 1$ (allowing $s=1$ here). Then $\l = c\l_n = c\l_1^*$, so replacing 
$V$ by the dual we may take $\l = c\l_1$. We have $c > 1$ by assumption (in the statement of Theorem \ref{rgreater}); and if $c=2$ then Lemma \ref{A2rssquares} shows that $s=1$, as in the statement of the lemma. Hence we may assume that $c\ge 3$. 

We have $m=cr$. The nonzero $\l$-label is over the root $\b_1^0$ in the $A_r$ factor of $L_Y'$. Let $\g_1$ be the root adjacent to $A_r$. Then $\l-\b_r^0-\cdots -\b_1^0-\g_1$ affords the level 1 $L_X'$-summand
\[
((r+1)\oplus (r-1)) \otimes V_{A_r}(c-1,0,\ldots ,0)\downarrow L_X'.
\]
Now $ V_{A_r}(c-1,0,\ldots ,0)\downarrow L_X'$ contains $r(c-1)\oplus (r(c-1)-4)$. If $r\ge 3$ it follows that the above tensor product contains $(cr-5)^4$, contrary to Lemma \ref{bound}(i). So finally assume that $r=2$.  Then level 0 for $L_X'$ is $S^c(2) = 2c\oplus (2c-4)\oplus \cdots$, while level 1 contains $(3\oplus 1)\otimes 
((2c-2)\oplus (2c-6))$, which contains $(2c-3)^3$. This contradicts Lemma \ref{bound}(iv). \hal

\begin{lem}\label{o2}
$\l_{A_s}$ is not $\l^k_2$ or $(\l^k_2)^*$ with $s\ge 3$.
\end{lem}

\pf If $\l_{A_s} = (\l^k_2)^* = \l_{s-1}^k$ then $\l=\l_2^*$ and $V\downarrow X$ is not multiplicity-free by Lemma \ref{A2rssquares}. 

Assume now that $\l_{A_s} = \l^k_2$. We have $s\ge 4$ by the previous paragraph. The top weight is $m = 2s-2$. The nonzero label is over the root $\b_{n-s+2}^k$. Let $\g_k$ be the root  adjoining the $A_s$ factor. Then $\l-\g_k-\b_{n-s+1}^k-\b_{n-s+2}^k$ affords the level 1 $L_X'$-summand
\[
((s+1)\oplus (s-1)) \otimes V_{A_s}(\l^k_3)\downarrow L_X'.
\]
This contains $((s+1)\oplus (s-1)) \otimes (3s-6)$, which contains $(4s-7)^2$. Hence Lemma \ref{bound}(ii) gives 
$4s-7 \le m+1 =2s-1$, forcing $s\le 3$, a contradiction. \hal

\begin{lem}\label{o1s}
$\l_{A_s}$ is not $\l^k_1+\l^k_s$.
\end{lem}

\pf Suppose $\l_{A_s} = \l^k_1+\l^k_s$. Then $m=2s$. If $s\ge 3$ then the weight $\l-\g_k-\b_1^k$ affords the level 1 $L_X'$-summand $((s+1)\oplus (s-1)) \otimes V_{A_s}(\l^k_2+\l^k_s)\downarrow L_X'$, and the latter factor contains $3s-2$. Hence this tensor product contains $(4s-3)^2$, forcing $4s-3 \le m+1 = 2s+1$ by Lemma \ref{bound}(ii), contradicting $s\ge 3$. Therefore $s=2$. Now  $\l-\g_k-\b_1^k$ affords $(3\oplus 1)\otimes (4\oplus 0)$ for $L_X'$, while 
 $\l-\g_k-\b_1^k-\b_2^k$ affords $(3\oplus 1)\otimes 2$. Hence level 1 for $L_x'$ contains $3^5$, contrary to Lemma \ref{bound}. \hal

\begin{lem}\label{o3}
$\l_{A_s}$ is not $\l^k_3$ or $(\l^k_3)^*$.
\end{lem}

\pf  Suppose $\l_{A_s} = \l^k_3$. Here $s=5,6$ or 7 (see Lemma \ref{labas}(iv)), and $m=3s-6$. The weight 
$\l-\g_k -\b_1^k -\b_2^k -\b_3^k$ affords the level 1 $L_X'$-summand $((s+1)\oplus (s-1)) \otimes 
V_{A_s}(\l^k_4)\downarrow L_X'$, which contains $((s+1)\oplus (s-1)) \otimes (4s-12)$, hence has $(5s-13)^2$. Therefore $5s-13\le m+1 = 3s-5$, forcing $s\le 4$, a contradiction.

If $\l_{A_s} = (\l^k_3)^*$ then again $m=3s-6$ and $\l-\g_k- \b_1^k -\cdots -\b_{s-2}^k$ affords level 1 $L_X'$-summand 
$((s+1)\oplus (s-1)) \otimes V_{A_s}(\l^k_{s-1})\downarrow L_X'$. This contains $(3s-3)^2$, contradicting Lemma \ref{bound}(ii). \hal

\begin{lem}\label{lc1}
$\l_{A_s}$ is not $c1$ or $1c$ with $s=2, c\ge 2$.
\end{lem}

\pf Suppose false. We have $m=2c+2$. If $\l_{A_s} = c1$, the weight 
$\l-\g_k-\b_1^k$ affords the level 1 $L_X'$-summand $(3\oplus 1) \otimes V_{A_2}(c-1,2)\downarrow L_X'$, which contains $(2c+1)^4$, contradicting Lemma \ref{bound}. And if $\l_{A_s}=1c$ then 
$\l-\g_k-\b_1^k$ affords the level 1 $L_X'$-summand $(3\oplus 1) \otimes V_{A_2}(0,c+1)\downarrow L_X'$, while 
$\l-\g_k-\b_1^k-\b_2^k$ affords $(3\oplus 1) \otimes V_{A_2}(1,c-1)\downarrow L_X'$; together, these contain $(2c+1)^4$, again a contradiction. \hal

\begin{lem}\label{l110}
$\l_{A_s}$ is not $110$ or $011$ with $s=3$.
\end{lem}

\pf Suppose false. We have $m=7$. If $\l_{A_s}=110$ then $\l-\g_k-\b_1^k$ affords the level 1 $L_X'$-summand $(4\oplus 2) \otimes V_{A_3}(020)\downarrow L_X'$, which contains $10^2$, contradicting Lemma \ref{bound}(ii). And if  
$\l_{A_s}=011$ then $\l-\g_k-\b_1^k-\b_2^k$ affords $(4\oplus 2) \otimes V_{A_3}(002)\downarrow L_X'$, which contains $4^4$, contrary to Lemma \ref{bound}(i). \hal

\vspace{4mm}
In view of Lemma \ref{labas}, we have now dealt with case (c) of Lemma \ref{lab}. 

It remains to deal with the case where $(r,s)=(1,1)$, excluded by assumption (\ref{not11}) until now.

\begin{lem}\label{oneone}
Suppose $(r,s)=(1,1)$ and $Y = SL(W) = SL_8$. Then up to duals, $\l$ is  $c\l_1\,(c\le 3)$, $\l_2$, or $\l_3$, as in Table $\ref{TAB1}$ of Theorem $\ref{MAINTHM}$.
\end{lem}

\pf  
Under the hypotheses of the lemma,  $\Pi(Y)\setminus \Pi(L_Y') = \{\g_1,\g_2\}$. Suppose first that $V^1(Q_Y)$ is trivial. 
If $\la \l,\g_1 \ra$ and $\la \l,\g_2 \ra$ are both nonzero, then $\l-\g_1$ and $\l-\g_2$ both afford level 1 $L_X$-summands 
$(2\oplus 0)\otimes 1 = 3\oplus 1^2$, so level 1 contains $1^4$, contrary to Lemma \ref{bound}. Hence, replacing $\l$ by the dual if necessary, we may take $\l = c\l_2$. If $c\ge 2$ then level 2 for $L_X'$ contains summands afforded by $\l-2\g_1$ and $\l-\b_1^0 -2\g_1-\b_1^1$. These summands are $2 \otimes S^2(2\oplus 0)$ and $ \wedge^2(2\oplus 0)$, which between them contain $2^6$. Since level 1 is only $3\oplus 1^2$, this leads to a contradiction using Proposition \ref{induct}. Hence $c=1$ and $\l=\l_2$, as in the conclusion of the lemma. 

Suppose now that $V^1(Q_Y)$ is nontrivial, so that $\l$ has a nonzero label on a factor of $L_Y' = A_1+A_3+A_1$. 
The proof of Lemma \ref{lab} implies that one of the following holds, replacing $\l$ by its dual if necessary:

(a) the $\l$-labelling of the $A_3$ factor of $L_Y'$ is $100$ or $200$, and the $A_1$ factors have label 0.

(b) the $\l$-labelling of the $A_3$ factor is $000$ and the $A_1$ factors have labels $c \ge d \ge 1.$

(c) the $A_3$ factor is labelled $000$, $\b_1^k$ is labelled 0 and $\b_1^0$ is labelled $m$ for some $m$.

Consider case (a). The proof of Lemma \ref{casealab} shows that the $\l$-labels of $\Pi(Y)\setminus \Pi(L_Y') $ are 0, so that 
$\l = c\l_3$ with $c=1$ or 2. If $c=1$ then $\l$ is as in the conclusion, so suppose $c=2$. Then $\l-\a_2-\a_3$ affords the level 1 $L_X'$-summand $1 \otimes V_{A_3}(110)\downarrow L_X'$. This contains $3^5$, contradicting Lemma \ref{bound}.

Now consider (b). We apply the method of Section \ref{levanal}. Here $\l - \b_1^0-\g_1$ and $\l - \g_2 - \b_1^k$ afford
$(c-1) \otimes (2\oplus 0) \otimes d$ and $c \otimes (2\oplus 0) \otimes (d-1),$ 
respectively. The first of these equals
$c-1\otimes 1\otimes1\otimes d = (c\oplus c -2)\otimes 1\otimes d = 
(V^1\otimes 1)+(c-2)\otimes 1\otimes d$,
where the final summand does not occur if $c=1$. So by 
Corollary~\ref{cover}, it suffices to show that
$c\otimes (2\oplus 0)\otimes d-1$ is not MF. This is a straightfoward check.

Finally, consider case (c). The proof of Lemma \ref{casesbc} shows that the $\l$-labels of $\Pi(Y)\setminus \Pi(L_Y') $ are 0, so that 
$\l = m\l_1$ and $V_Y(\l) = S^m(W)$. Suppose $m \ge 4.$  At level 1 the restriction to $L_X'$ is $(m-1) \otimes (2\oplus 0)$.  And at level 2, $\l-2\b_1^0-2\g_1 $ and $\l-\b_1^0-\g_1-(\b_1^1-\b_2^1-\b_3^1)-\g_2 $ afford
$(m-2) \otimes S^2(2\oplus 0)$ and $(m-1) \otimes 1$, respectively.  As $m-2 \ge 2$, together these contain $(m-2)^5$ whereas only three summands of $(m-2)$ can arise from level 1.  So this contradicts Corollary \ref{conseq}.  Therefore
$m \le 3$ as in Theorem \ref{rgreater}.
\hal

\vspace{2mm}
This completes the proof of Theorem \ref{rgreater}.

\section{Case $\delta =r0$}

We now prove 

\begin{thm}\label{a2r0}
Let $X = A_2$ and $\d = r\o_1$ with $r\ge 2$. Let $W = V_X(\d)$ and take $X < Y = SL(W)= A_n$. Suppose $\l$ is a dominant weight for $Y$ such that $\l$ is not $\l_1$ or its dual. Then $V_Y(\l)\downarrow X$ is multiplicity-free if and only if $r,\l$ are as in Table $\ref{tab:A2An}$, where $\l$ is given up to duals.
\end{thm}

\begin{table}[h]
\caption{MF pairs, $A_2\leq A_n$, $\d = r\omega_1$}\label{tab:A2An}
\[
\begin{array}{|l|l|}
\hline
r & \lambda\\
\hline
\hline
\hbox{all } r& \lambda_2, \,2\lambda_1,\, \lambda_1+\lambda_n \\
\hline
r\leq 6 & \lambda_3\\
\hline
r\leq 5 & 3\lambda_1 \\
\hline
r\leq 4  & \lambda_4 \\
\hline
r\leq 3 & \lambda_i\, (\hbox{all } i), \,4\lambda_1,\, \lambda_1+\lambda_2 \\
\hline
r=2 & a\lambda_1\,(a\geq 1), \\
           &  \lambda_i+\lambda_j \,(\hbox{all } i,j), \\ 
            & 2\lambda_2,\, 3\lambda_2,\\
             & 2\lambda_1+\lambda_5, \,3\lambda_1+\lambda_5,\\
            &2\lambda_1+\lambda_2,\,3\lambda_1+\lambda_2\\
\hline
\end{array}
\]
\end{table}

Note that the multiplicity-freeness of all the examples in Table \ref{tab:A2An} is given by Theorem \ref{yesMF}.

 \subsection{Case $r=2$}

Here we prove Theorem \ref{a2r0} in the case where $r=2$: 

\begin{prop}\label{r=2} Let $r=2$, let $V=V_Y(\l)$, and suppose that $V\downarrow X$ is multiplicity-free. Then $\lambda$
or its dual is as in Table $\ref{tab:A2An}$.
\end{prop}

\pf Note that $\dim V_X(2\o_1) = 6$, so $Y  =A_5$. 
Notation will be as in Chapter \ref{notation} although matters
simplify since $Y$ has small rank. In particular we note that $\g_1 = \b_3$ and $\g_2 = \b_5$.

Let $P_X=Q_XL_X$ be the parabolic subgroup of $X$ with Levi factor $L_X = T\langle U_{\pm\a_1}\rangle$
 and 
unipotent radical $Q_X =\langle U_{-\a_2},U_{-\a_1-\a_2}\rangle$.  Let $P_Y = Q_YL_Y$ be 
the parabolic subgroup of $Y$ constructed according to the $Q_X$-levels of $W$; so without loss of 
generality,
we have $L_Y' = \langle U_{\pm\b_i}\ |\ i=1,2,4\rangle$ and $P_Y$ contains the opposite Borel 
subgroup $B^-$.

Recall that we often simplify notation by writing $(a\l_1^0 + b\l_2^0) \otimes (d\l_1^1)$ instead of $V_{C^0}(a\l_1^0 + b\l_2^0) \otimes V_{C^1}(d\l_1^1).$  So this is the $L_Y'$ irreducible module with highest weight 
$(a\l_1^0 + b\l_2^0) + (d\l_1^1).$


We assume neither $\lambda$ nor its dual is as in 
Table \ref{tab:A2An}, and aim for a contradiction.
We assume also that $\l$ and its dual are not equal to any of the following:
\begin{equation}\label{listb}
\l_1+2\l_2,\,\l_1+2\l_4,\,2\l_1+\l_4,
\end{equation}
since for these weights it can be checked using Magma that $V_Y(\l)\downarrow X$ is not MF.

We will treat a series of cases. But first we make some general remarks. 
Say $ \lambda = a\lambda_1 + b\lambda_2 + c\lambda_3 + d\lambda _4 + e\lambda_5$ and 
$V\downarrow X$  is multiplicity-free.
Then $V^1 = ((ab)\downarrow L_X')\otimes d$ is multiplicity-free. By \cite[Theorem 1]{LST}, the fact that 
$(ab)\downarrow L_X'$ is multiplicity-free implies that $a \le 1$ or $b\le 1$. Moreover $(a1)\downarrow L_X'$ contains 
$(2a+2)\oplus (2a)$, and $(a0)\downarrow L_X'$ contains $(2a)\oplus (2a-4)$ (if $a\ge 2$); hence the following hold:
\begin{enumerate}[]
\item{ i)} $abd = 0$;
\item{ ii)}  if $a+b>1$, then $d = 0$ or $1$; 
\item{ iii)} if $ab\ne 0$, then $\{a,b\} = \{x,1\}$, for some $x\ne 0$.
\end{enumerate}
Similarly, applying this reasoning to $(V^*)^1(Q_Y)$, we see that
\begin{enumerate}[]
\item{iv)} $bde = 0$;
\item{v)} if $d +e>1$, then $b = 0$ or $1$;
\item{vi)} if $de\ne 0$, then $\{d,e\} = \{x, 1\}$, for some $x\ne 0$.
\end{enumerate}

\noindent{\bf{Case 1}}: Assume $V^1(Q_Y)$ is the trivial $L_Y'$-module.

In this case, $\lambda=x\lambda_3+y\lambda_5$. Then as $\lambda$ and its dual are not as in Table \ref{tab:A2An} 
or (\ref{listb}), we have $x\ne 0$ and one of $x$, $y$ is greater than $1$; hence
\begin{equation}\label{v2eq}
V^2=\left\{\begin{array}{l} 1\oplus1\oplus 3, \hbox{ if } y\ne 0,\\
 1\oplus 3, \hbox{ if } y=0.\end{array}\right.
\end{equation}
Now we consider $V^3$. The weight $\lambda-\g_1-\b_1^1-\g_2$ affords a summand $2$ for $L_X'$, 
occurring with multiplicity 3 if $y\ne0$. In addition, $\lambda-\b_2^0-2\g_1-\b_1^1$ 
affords 
a summand $2$. If $x>1$, $\lambda-2\g_1$ affords a summand 
$(2\l_2^0 \otimes 2\l_1^1)\downarrow L_X' = 6\oplus4\oplus2\oplus2$; and if $y>1$, $\lambda-2\g_2$ 
affords a summand $2$ as well. (Note that we have used the fact that when $x>1$, the 
multiplicity 
of the weight $\lambda-\b_2^0-2\g_1-\b_1^1$ is 2.)

Hence the summand 2 appears in $V^3$ with multiplicity at least 4 (at least 5 if $y\ne 0$). However, in view of (\ref{v2eq}), this contradicts Corollary~\ref{conseq}.


\noindent{\bf{Case 2}}: $\lambda=c\lambda_1+x\lambda_3+y\lambda_5$.

By duality, the considerations of Case 1, and the fact that $\lambda$ is not as in Table \ref{tab:A2An} or (\ref{listb}),
 we see that $cy\ne 0$. We apply the method of Section \ref{levanal}. Now $V^1$ consists of precisely
the summands $2c-4i$,
for $0\leq i\leq \lfloor\frac{c}{2}\rfloor$. We now consider $V^2(Q_Y)$. 
Since $y\ne 0$, $\lambda-\gamma_2$ affords a summand $V^1\otimes 1$ of $V^2$. Then Corollary~\ref{cover} implies that any other summands must be MF. So in particular $x=0$, else $\lambda-\gamma_1$ affords the summand $(c\lambda_1^0+\lambda_2^0)\otimes \lambda_1^1$, which is non-MF on restriction to $L_X'$.

Assume for the moment 
that $c\ge 2$ and $y\geq 2$. Then $V^2$ has summands $(2c-1)\oplus (2c-3)$,
afforded by $\lambda-\b_1^0-\b_2^0-\g_1$ and $2c+1\oplus 2c-1\oplus 2c-3\oplus 2c-5$,
afforded by $\lambda-\g_2$. All
other irreducible summands have lower weights. Now consider $V^3(Q_Y)$, where we have 
 \begin{enumerate}[(1)]
\item a summand $S^c(2)\otimes 2$, afforded by $\lambda-2\g_2$,
\item a summand $S^{c-2}(2)\otimes 2$, afforded by the weight 
$\lambda-2\b_1^0-2\b_2^0-2\g_1$, 
\item a summand $S^{c-1}(2)\otimes 2$, afforded by 
$\lambda-\b_1^0-\b_2^0-\g_1-\g_2$, 
\item a summand $V_{L_1}(c\lambda_1^0+\lambda_2^0)\downarrow L_X'$ afforded by 
$\lambda-\g_1-\b_1^1-\g_2$, and 
\item two summands
$S^{c-1}(2)$, afforded by two remaining basis vectors in the weight space for 
$\lambda-\b_1^0-\b_2^0-\g_1-\b_1^1-\g_2$.
\end{enumerate}
Note that the sum of the summand in (4) and one of the summands in (5) is isomorphic to $S^c(2)\otimes 2$. 
Now let us count the occurrences of the weight $2c-2$ in $V^3$:
the summand from (1), respectively (2), (3), affords 2, resp. 1, 1, and the summands in (4) and (5) 
provide another 3. So in all we have 7 such occurrences while only 5 are allowed by Corollary \ref{conseq}.
 Hence we deduce that one
of $c$ and $y$ must equal 1.

By duality, 
we may assume $c>1$, and so $y=1$. By our initial assumptions on $\l$, 
we have $c\geq 4$. In this case, $V^2 = (S^c(2)\oplus S^{c-1}(2))\otimes 1$,
while $V^3(Q_Y)$ consists of all of the above summands except the first one listed. Here we count 
the occurrences of the weight $2c-6$. At most five such occurrences are allowed in $V^3$ by 
Corollary \ref{conseq}, while 
the summands (2) and (3) afford 2, respectively 1, and the summands (4) and (5) together afford 3 more summands,
 which gives the final contradiction.
 This completes the consideration of Case 2.

\noindent{\bf{Case 3}}: $\lambda=c\lambda_2+x\lambda_3+y\lambda_5$.

By Case 1, $c\ge1$. Here $V^1(Q_Y)\downarrow L_Y' = c\lambda_2$ affords $L_X'$-summands $2c-4i$, for 
$0\leq i\leq\lfloor \frac{c}{2}\rfloor$. 
As in Case 2, we deduce that if $y\ne 0$, then $x=0$.

Suppose $y=0$ and $x\ne 0$. Then $\lambda-\gamma_1$ and $\lambda-\beta_2^0-\gamma_1$ afford summands of $V^1(Q_Y)$ whose sum affords the $L_X'$ module $2\otimes 1\otimes V^1 = (3\oplus 1)\otimes V^1$. So by Corollary~\ref{cover}, $s\otimes V^1 = 3\otimes S^c(2)$ must be MF. This shows that $c=1$. 
Moreover, by our initial assumptions on $\lambda$, we have $x>1$. We now consider $V^3(Q_Y)$. Now $V^1=2$ and
$V^2 = 5\oplus 3\oplus 1\oplus 3\oplus 1$. In particular, there exist at most four $L_X'$-summands of weight $4$ in $V^3$. But $\lambda-\gamma_1-\beta_1^1-\gamma_2$
affords $S^2(2)$, $\lambda-2\gamma_1$ affords $S^3(2)\otimes 2$ and $\lambda-\beta_2^0-2\gamma_1$ affords $(\lambda_1^0+\lambda_2^0)\downarrow L_X'\otimes 2$, in which we find a total of $5$ summands of weight $4$, giving the final contradiction.

We now suppose $x=0$ and $y\ne 0$. We claim that $c=1$. As before, $\lambda-\gamma_2$ affords $V^1\otimes 1$ and so any remaining $L_X'$ summands of $V^2$ must be MF. 
Now $\lambda-\beta_2^0-\gamma_1$
affords $(\lambda_1^0+(c-1)\lambda_2^0)\downarrow L_X'\otimes 1$, which is easily seen to be MF only if $c=1$.
 By our initial assumptions on $\lambda$, we 
have $y\geq 3$. Now consider $M: = V^\ast$, the irreducible $Y'$-module with highest 
weight $\mu= y\lambda_1+\lambda_4$.
The weights $\mu-\beta_1^1-\gamma_2$, $\mu-\gamma_1-\beta_1^1$, $\mu-\beta_1^0-\beta_2^0-\gamma_1$ and $\mu-\beta_1^0-\beta_2^0-\gamma_1-\beta_1^1$, afford summands of $V^1(Q_Y)$ the sum of which has restriction to $L_X'$ giving $(V^1\otimes 1)\oplus (S^{y-1}(2)\otimes 2)$. But since $y\geq 3$, the second summand is not MF, contradicting Corollary~\ref{cover}.

So finally in Case 3, we have reduced to $\lambda = c\lambda_2$.
By our initial assumptions on $\lambda$, we have $c\geq 4$. Let $M = V^*$, the $Y$-module
with highest weight $\mu =c\lambda_4$. Here $M^1 = c$, and 
$M^2 = c+1\oplus c-1\oplus c-3\oplus c-1$. Let us now consider $M^3(Q_Y)$. 
The weight
$(\mu-2\g_1-2\b_1^1)\downarrow L_Y'$
affords $L_X'$-summands  $c+2\oplus c\oplus c-2\oplus c-2$ of $M^3(Q_Y)$.
The weight
$(\mu-2\b_1^1-2\g_2)\downarrow L_Y'$
affords an $L_X'$-summand  $c-2$ of $M^3(Q_Y)$. In addition, the weight
$(\mu-\g_1-\b_1^1-\g_2)\downarrow L_Y'$
affords $L_X'$-summands  $c+2\oplus c\oplus c-2$ of $M^3(Q_Y)$. Finally, the weight 
$\mu-\g_1-2\b_1^1-\g_2$ occurs with multiplicity 2 in $M$, but
 with multiplicity 1 in 
the previous $L_X'$-summand, so we have a further summand 
which affords $L_X'$-summands $c\oplus c-2\oplus c-4$.
In particular, we have the weight $c-2$ occurring 5 times in $M^3(Q_Y)$, contradicting Corollary \ref{conseq}.

\noindent{\bf{Case 4}}: $\lambda = c\lambda_1+x\lambda_3+\lambda_4+y\lambda_5$, $c\geq1$. 

Here $V^1 = S^c(2)\otimes 1$. The summands of $V^1(Q_Y)$ afforded by $\lambda-\gamma_1-\beta_1^1$ and $\lambda-\beta_1^0-\beta_2^0-\gamma_1-\beta_1^1$ have sum affording the $L_X'$ module $S^c(2)\otimes 2$.
In addition, $\lambda-\beta_1^1-\gamma_2$ affords $S^c(2)$. The sum of these then gives $V^1\otimes 1$ and so by Corollary~\ref{cover}, the sum of all remaining summands must be MF. Now $\lambda-\beta_1^0-\beta_2^0-\gamma_1$ affords $S^{c-1}(2)\otimes 2$, which is MF only if $c\leq 2$.
If $y\ne 0$, $\lambda-\gamma_2$ affords $S^c(2)\otimes 2$ and the sum of these two is not MF.
So $y=0$ and as $\lambda$ is not as in Table~\ref{tab:A2An} nor in the list (\ref{listb}), we have $x\ne 0$. But then finally, $\lambda-\gamma_1$ affords a non-MF summand of $V^2$, giving the final contradiction.

\noindent{\bf{Case 5}}: $\lambda = c\lambda_2+x\lambda_3+\lambda_4+y\lambda_5$, $c\geq1$.

The argument here is similar to the previous case. The $L_Y'$-summands afforded by $\lambda-\beta_2^0-\gamma_1-\beta_1^1$ and $\lambda-\gamma_1-\beta_1^1$ sum to give the $L_X'$-module $S^c(2)\otimes 2$. Summing this with the summand afforded by $\lambda-\beta_1^1-\gamma_2$ gives $V^1\otimes 1$. Again by Corollary~\ref{cover}, all remaining summands should sum to an MF $L_X'$-module. The summand afforded by $\lambda-\beta_2^0-\gamma_1$ is MF only if $c=1$, and then the restriction is $2\otimes 2$. And as $\lambda$ is not in Table~\ref{tab:A2An}, $x+y>0$. But then we have either a summand $2\otimes 2$ or a summand $S^2(2)\otimes 2$, in addition to the summand $2\otimes 2$, and the sum of these is not MF.

\noindent{\bf{Case 6}}: $\lambda = \lambda_1+x\lambda_3+d\lambda_4+y\lambda_5$
 or $\lambda = \lambda_2+x\lambda_3+d\lambda_4+y\lambda_5$, with $d\geq 2$. 

Here $V^1=d+2\oplus d\oplus d-2$. The weights $\lambda-\g_1-\b_1^1$ and 
$\lambda-\b_1^1-\g_2$ each afford an $L_Y'$-summand of $V^2(Q_Y)$, and each contributes a 
$L_X'$-summand of weight $d+1$. In addition, the weight 
$\lambda-\b_1^0-\b_2^0-\g_1$ (or $\lambda-\b_2^0-\g_1$ in the 
second case) contributes another $L_X'$-summand $d+1$. This then implies that $y=0$, else 
$\lambda-\g_2$ affords an $L_Y'$-summand which gives rise to a fourth $d+1$ 
$L_X'$-summand, contradicting Corollary \ref{conseq}.
Now consider the module $M = V^*$ with highest weight $\lambda = 
d\lambda_2+x\lambda_3+\lambda_j$, where $j = 5$ or $j=4$, according to the
 choice of $\lambda$. These weights have been handled in Cases 3 and 5.

\noindent{\bf{Case 7}}: $\lambda = x\lambda_3+d\lambda_4+y\lambda_5$, $d\geq 1$. 

Here we note that for $M = V^*$, we have previously treated this configuration 
unless $\{y,d\} = \{1,c\}$ for some $c\geq 1$; so assume this is the case.
 Let $\mu$ be the highest weight of $M$, so 
$\mu\downarrow L_X'$ affords summands $2c+2\oplus 2c$ and all other $L_X'$-summands have lower highest weights. Moreover, it is straightforward to see that $x=0$,
 else there are too many $L_X'$-summands of $M^2(Q_Y)$ of  weight $2c+3$ (afforded by $\mu-\g_1$ and $\mu-\b_2^0-\g_1)$. 

Since we are assuming $\lambda$ and $-w_0\lambda$ are not as in Table \ref{tab:A2An}, the case $\lambda = \lambda_4+c\lambda_5$  is covered by Lemma \ref{a1000}.

Finally, consider the module $M =V^*$ with highest weight $\mu = \lambda_1+c\lambda_2$. By our assumptions on $\l$, we have  $c\geq 3$.  If $c=3$, a Magma computation shows that $V\downarrow X$ is not MF, so assume $c\ge 4$. Here $M^1$ has summands 
$2c+2\oplus 2c\oplus 2c-2\oplus 2c-4$ and all other
summands have lower highest weights. Now consider $M^2(Q_Y)$. The weight 
$(\mu-\b_2^0-\g_1)\downarrow L_Y'$
affords three $L_X'$-summands of weight $2c-3$ and the weight $\mu-\b_1^0-\b_2^0-\g_1$ affords a fourth such summand,
contradicting Corollary \ref{conseq}.

\noindent{\bf{Case 8}}: $V^1(Q_Y)\downarrow L_Y' = c\lambda_1+\lambda_2$ or 
$\lambda_1+c\lambda_2$, for $c\geq 1$.

Here we simply note that the module $M = V^*$ has been treated in one of the above cases.

The  proof of Proposition~\ref{r=2} is completed by applying the above cases to $V$ and $V^*$.\hal

 \subsection{General case $\d = r\omega_1$, $r\geq 3$}

Here we prove Theorem \ref{a2r0} under the assumption that $r\ge 3$. Suppose $V_Y(\l)\downarrow X$ is MF
and write $V = V_Y(\l)$.

Notation will be as in Chapter \ref{notation}. Let $Y = A_n$, so that $n+1 = \frac{(r+1)(r+2)}{2}$. 
We have $\Pi(Y) = \{\b_1, \ldots , \b_n\}$ with $ \{\l_1, \ldots , \l_n\}$ the corresponding fundamental dominant weights. Let 
$P_X=Q_XL_X$ and $P_Y = Q_YL_Y$ be as in Chapter \ref{notation}.
 Then $L_Y' = C^0 \times  \cdots \times C^{r-1} $ ($C^r = 0$), where $C^i$ is of 
type 
$A_{r-i}$. For $0\leq i\leq r-1$, let $\{\l_1^i, \ldots , \l_{r-i}^i\}$ be the fundamental 
dominant weights for $C^i$, and for convenience, we will write $(\l_j^i)^*$ for the weight 
$\l_{r-i-j+1}^i$,
that is, the highest weight of the dual to the $C^i$ module with highest weight $\l_j^i$. 
Recall that $\l_j^i$ corresponds to the simple root 
$\b_j^i\in \Pi(C^i).$ 
Now since the embedding of $L_X'$ in $C^i$ is via the irreducible representation with highest weight 
$r-i$, we have that each simple root $\b\in\Pi(L_Y)$ has restriction to $T\cap L_X'$
being $2$. Given that $V^1$ is multiplicity-free and using \cite{LST}, it is straightforward to deduce the
following lemma.

\begin{lem}\label{cases} One of the following holds.
\begin{itemize}
\item[{\rm (1)}] $V^1(Q_Y)$ is the trivial $L_Y'$-module.
\item[{\rm (2)}] There exists a unique $i$ such that $\mu^i\ne 0$. Moreover,
the pair $(C^i,\mu^i)$ (or the pair corresponding to the dual module)
 appears in Table~{\rm\ref{dist}} below.
\item[{\rm (3)}] There exists $0\leq i<j\leq r-1$ such that 
$\mu^i\ne 0\ne \mu^j$, 
and
$\mu^k=0$ for all $k\in\{0,\dots,r-1\}\setminus\{i,j\}$. Moreover precisely one of 
the 
following holds:
\begin{itemize}
\item[{\rm (a)}] $\mu^i = \l_1^i$ or $(\l_1^i)^*$ and 
$\mu^j = \l_1^j$ or $(\l_1^j)^*$;
\item[{\rm (b)}] $j=r-1$, $\mu^j = b\l_1^{r-1}$, $\mu^i = c\l_1^i$ 
or $c(\l_1^i)^*$, for some $cb\ne 0$, and exactly one of $c,b$ is greater than 1.
\item[{\rm (c)}] $j=r-1$, $\mu^j = \l_1^{r-1}$, $\mu^i = \l_k^i$ or 
$(\l_k^i)^*$ for some $1< k < r-i$.
\end{itemize}
\end{itemize}
\end{lem}

 \begin{table}[h]
\caption{}\label{dist}
\[
\begin{array}{|l|l|}
\hline
C^i & \mu^i  \\
\hline
A_m &  \l_1, \,\l_2,\, 2\l_1, \,\l_1+\l_m,  \\
& \l_3 \,(5 \le m \le 7), \\
&  3\l_1\,( m \le 5),\, 4\l_1 (m\le 3),\,5\l_1 (m\le 3)  \\
A_3 &  110  \\
A_2   &  c1, c0   \\
\hline
\end{array}
\]
\end{table}

We shall treat each of the cases (1), (2), (3) of Lemma \ref{cases} in turn. First we record the following, which is just a Magma check.

\begin{lem}\label{a2-examples-2} 
If $r,\l$ are as in the table below, then $V_Y(\l) \downarrow X$ is not MF.
\[
\begin{array}{|l|l|}
\hline
r & \l \\
\hline
3 & 5\l_1,\,2\l_1+\l_9,\,2\l_1+2\l_9,\\
   & \l_1+\l_3,\,\l_2+\l_3,\,\l_2+\l_9 \\
4 & \l_5,\,4\l_1 \\
5 & \l_4,\,4\l_1 \\
7 & \l_3 \\
\hline
\end{array}
\]
\end{lem}

\begin{lem}\label{v1-triv} Assume that  $V^1(Q_Y)$ is the trivial 
$L_Y'$-module. Then  
one of the following holds:
\begin{itemize}
\item[{\rm (a)}] $\lambda=a\lambda_n$ for $a\leq 2$;
\item[{\rm (b)}] $r=3$, $\lambda\in\{\lambda_4,\lambda_7,a\lambda_n\, (a\leq 4)\}$;
\item[{\rm (c)}] $r=4,5$, $\lambda=3\lambda_n$;
\item[{\rm (d)}] $r\leq 7$, $\lambda=\lambda_{n-2}$.
\end{itemize}
Hence $\l$ is as in Tables $\ref{TAB1}-\ref{TAB4}$ of Theorem $\ref{MAINTHM}$.
\end{lem}

\pf By Corollary \ref{conseq}, each irreducible $L_X'$-summand of $V^2(Q_Y)$ with highest weight different from 
1 occurs with multiplicity at most 1, and 1 can occur at most twice as a highest weight. 
This then implies that one of the following holds.
\begin{enumerate}[(i)]
\item There exists a unique $\gamma\in\Pi(Y)$ with $\langle\lambda,\gamma\rangle\ne 0$, or
\item $\langle\lambda,\b_n\rangle\ne 0$ and $\langle\lambda,\gamma\rangle\ne 0$ for a 
unique $\gamma\in\Pi(Y)\setminus\{\b_n\}$.
\end{enumerate}

Indeed, if there exist distinct $\gamma,\delta\in\Pi(Y)\setminus\Pi(L_Y)$, with 
$\gamma\ne\b_n\ne\delta$ and $\langle\lambda,\gamma\rangle\ne 0\ne\langle\lambda,\delta\rangle$,
 then there exists $s\ne t$, $s,t\ge2$ such that $\lambda-\gamma$ affords an $L_X'$-summand 
$s\otimes (s-1)$ and $\lambda-\delta$ affords an $L_X'$-summand $t\otimes(t-1)$ of $V^2(Q_Y)$, contradicting the above remarks.

For case (i), we first suppose $\gamma\not\in\{\b_{n-2}, \b_n\}$. Then there 
exists $s>2$ such that $V^2 = 2s-1\oplus 2s-3\oplus\cdots\oplus 1$.
Now consider $V^3(Q_Y)$; if $\gamma = \b_t$, then the weight 
$\lambda-\b_{t-1}-2\b_t-\b_{t+1}$ affords an $L_Y'$ summand of $V^3(Q_Y)$, and this in 
turn restricts to give an $L_X'$-summand 
$\wedge^2(s)\otimes \wedge^2(s-1)$. In addition, there exists an $L_Y'$-summand of $V^2(Q_Y)$ 
whose restriction to $L_X'$ affords a summand $s\otimes (s-2)$, and if $s<r$ an additional $L_X'$-summand
$(s+1)\otimes (s-1)$. Counting the occurrences of the $L_X'$-summands 2, and applying Corollary \ref{conseq},
we deduce that $s=r=3$. Moreover $\langle\lambda,\gamma\rangle=1$ else $\lambda-2\gamma$
 affords a summand $S^2(3)\otimes S^2(2)$ of $V^3$, yielding too many summands $2$. But now $V = \wedge^4(W)$
and $r,\l$ are as in (b).

If $\gamma = \b_{n-2}$, we claim that $\langle\lambda,\gamma\rangle =1$. Suppose otherwise.  We will obtain a contradiction from consideration of $V^3.$
We have $V^2 = 3\oplus 1$. The weight 
$(\lambda-2\b_{n-2})\downarrow L_Y'$ affords an $L_X'$-summand  containing
$2\oplus 2$. In addition, the weight $\lambda-\b_{n-3}-2\b_{n-2}-\b_{n-1}$ affords an $L_Y'$-summand 
which produces a third
$L_X'$-summand 2. Finally, the weight $(\lambda-\b_{n-2}-\b_{n-1}-\b_n)\downarrow L_Y'$
 affords a fourth $L_X'$-summand with highest weight 2, contradicting Corollary \ref{conseq}. Hence 
$\lambda = \lambda_{n-2}$. Now consider the module $M = V^*$ and deduce from Table~\ref{dist} that 
$r\leq 7$, as in (d).

To complete case (i), we must consider the case where $\gamma = \b_n$ and so 
$\lambda = a\l_n$ for some $a\geq 1$. Assume $\l$ is not as in the conclusion. Combined with Lemma \ref{a2-examples-2}, this implies that $a\ge 3$, and also that $a\ge 6$ if $r=3$, and $a\ge 5$ if $r=4,5$. 
Now the result of Table \ref{dist} applied to the dual $V^*$ gives a contradiction.

We now turn to case (ii). Here there exists $s>1$ such that $V^2(Q_Y)$ has an $L_X'$-summand 
$2s-1\oplus 2s-3\oplus\cdots\oplus 1$, in addition to the summand of highest weight 1 
afforded by $\lambda-\b_n$; moreover, there are no further summands of $V^2$.
Now consider $V^3(Q_Y)$. If $s>2$, then as in case (a), we have an $L_X'$-summand of the 
form $\wedge^2(s)\otimes\wedge^2(s-1)$, which affords at least two $L_X'$-summands of
highest weight 2. In addition,
we have a summand of the form $s\otimes (s-2)$, which affords another irreducible of highest weight 2. Finally, we have
a summand of the form $s\otimes (s-1)\otimes 1$ which has two $L_X'$-irreducible summands of 
highest weight 2. This contradicts Corollary \ref{conseq}. So finally, assume $s=2$, so that 
$\gamma = \b_{n-2}$. Considering the dual module $M = V^*$, Table~\ref{dist} implies that $r=3$, and $M$ has 
highest weight $\mu=\lambda_1+\lambda_3$. This contradicts  Lemma~\ref{a2-examples-2}.

This completes the proof of the Lemma.\hal

\begin{lem}\label{2fctrs-1wedge} $\lambda$ is not as in Lemma~$\ref{cases}{\rm (3c)}$. 
\end{lem}

\pf Assume $\lambda$ is as in Lemma~\ref{cases}(3c).  
Set $r-i = d$, so $C^i$ is of type $A_d$ and by assumption $d>2$. Moreover, $C^j$ is of type 
$A_1$.
Since $V^1$ is multiplicity-free, Table~\ref{dist} implies that either 
$\lambda\downarrow C^i = \l_2^i$ or $(\l_2^i)^*$ or $5\leq d\leq 7$ and 
$\lambda\downarrow C^i =  \l_3^i$ or $(\l_3^i)^*$. We first show that $d>3$; 
in particular $r>3$. Indeed, otherwise,
$V^1 = 5\oplus 3\oplus 1$. The weight $\lambda-\b_{n-1}-\b_n$ affords 
$L_X'$-summands 
$4\oplus 0$ of $V^2(Q_Y)$; $\lambda-\b_{n-2}-\b_{n-1}$ affords 
$6\oplus 4\oplus 2\oplus 2$; and
finally, $\lambda-\b_{n-7}-\b_{n-6}-\b_{n-5}$ affords an additional $L_X'$-summand $3\otimes2\otimes1$. 
Counting the highest 
weight 4 summands 
leads to a contradiction. Hence $d>3$ as claimed.

Note that $\lambda-\beta_{n-1}-\beta_n$ and $\lambda-\beta_{n-2}-\beta_{n-1}$ afford summands of $V^1(Q_Y)$ whose restriction to $L_X'$ has sum equal to $V^1\otimes 1$. So by Corollary~\ref{cover}, the quotient $V^2/V^1\otimes 1$ must be MF.
Now there exists a summand of $V^2(Q_Y)$ afforded by $\lambda-\beta-\gamma_{r-i+1}$ for some $\beta$ in the root system of $C^{r-i}$. The restriction to $L_X'$ is then one of $$d\otimes (d-1)\otimes 1,\ \wedge^2(d)\otimes (d-1)\otimes 1,\ \wedge^3(d)\otimes (d-1)\otimes 1,\ \wedge^4(d)\otimes (d-1)\otimes 1,$$ where $5\leq d\leq 7$ in the last case. None of these summands is MF, giving the desired contradiction.\hal

\begin{lem}\label{2fctrs-1sym} $\lambda$ is not as in Lemma~$\ref{cases}{\rm (3b)}$. 
\end{lem}

\pf Assume $\lambda$ is as in Lemma~\ref{cases}(3b). Here we have $j=r-1$ and $C^j = A_1$. 
Set $d = r-i$, so $C^i$ is of type $A_d$ and $d\geq 2$. 
First note that $\lambda\downarrow C^j = \l_1^{r-1}$, else for the module $M = V^*$,
$M^1(Q_Y)$ is not multiplicity-free, contradicting Corollary \ref{conseq}.
Hence, we have
$\lambda\downarrow C^j = \l_1^{r-1}$ and $\lambda\downarrow C^i = c\l_1^i$ or $(c\l_1^i)^*$ 
for some $c>1$.
So $V^1 = S^c(d)\otimes 1$; in particular, $V^1(Q_Y)$ has $L_X'$-summands 
$cd+1\oplus cd-1\oplus cd-3\oplus cd-5$, and all other summands have smaller highest weights. 

For the moment, assume that $d>2$. Then the weights $\lambda-\b_{n-1}-\b_n$ and 
$\lambda-\b_{n-2}-\b_{n-1}$ afford 
$L_Y'$-summands of $V^2(Q_Y)$, each of which contributes an $L_X'$-summand of highest weight $cd-4$. 
In addition, there exists a 
weight $\mu$ affording an $L_Y'$-summand of $V^2(Q_Y)$
such that $\mu\downarrow L_Y' = (c-1)\l_1^i + +\l_1^{i+1} + \l_1^{r-1}$, if 
$\lambda\downarrow C^i = c\l_1^i$, and a weight 
$\nu$
affording an $L_Y'$-summand of $V^2(Q_Y)$ such that $\nu\downarrow L_Y'= 
(\l_2^i+(c-1)\l_1^i)^* + \l_1^{i+1} + \l_1^{r-1}$,
if $\lambda\downarrow C^i = (c\l_1^i)^*$.
The summand afforded by $\mu$ gives two further $L_X'$-summands of highest weight $cd-4$, 
contradicting Corollary \ref{conseq}, while 
the summand afforded by $\nu$ affords two $L_X'$-summands of highest weight $cd+2d-4$, and 
hence $cd+2d-4\leq cd+2$, and so $d=3$. 
But in this last case, the $L_Y'$-summand afforded by $\nu$ gives three $L_X'$-summands of 
weight $3c+2=cd+2$, contradicting Corollary \ref{conseq}.

Hence we have now reduced to the case where $d=2$, so $V^1$ has summands
$2c+1\oplus 2c-1\oplus2c-3\oplus 2c-5$, and all other summands have lower highest weights.
 Recall that $r>2$, so $i\ne 0$. If $\lambda\downarrow C^i = c\l_1^i$, the weight 
$\mu = \lambda-\b_{n-5}-\b_{n-4}\in V^2(Q_Y)$
affords an $L_Y'$-summand
such that $\mu\downarrow L_Y'= (\l_1^{i-1})^* + ((c-1)\l_1^i  + \l_2^i) + \l_1^{r-1}$. The action 
of 
$L_X'$ on this summand affords three $L_X'$-summands of weight $2c+2$, contradicting Corollary \ref{conseq}.
So finally, we have $\lambda\downarrow L_Y' = (c\l_1^i)^* + \l_1^{r-1}$. Here we have 
$L_Y'$-summands of $V^2(Q_Y)$ afforded by weights $\lambda-\b_{n-5}-\b_{n-4}-\b_{n-3}$, 
$\lambda-\b_{n-3}-\b_{n-2}$ and $\lambda-\b_{n-2}-\b_{n-1}$. Restricting 
these to $L_X'$ produces
three summands of highest weight $2c+2$, again contradicting Corollary \ref{conseq}.\hal

\begin{lem}\label{2ntrl} $\lambda$ is not as in Lemma$~\ref{cases}\rm{(3a)}$.
\end{lem}

\pf Assume $\lambda$ is as in Lemma  \ref{cases}(3a). Here 
$\lambda\downarrow C^m = \l_1^m$ or $(\l_1^m)^*$, for $m\in\{i,j\}$. We first 
claim that 
$j= r-1$; that is, $C^j$ is of type $A_1$. Suppose otherwise. Set $d = r-i$ and $f = r-j$, so 
$C^i$ is of type $A_d$
and $C^j$ is of type $A_f$, where $d>f\geq 2$.  
Then $V^1 =  d\otimes f$. Now considering all possible cases for the various 
dispositions
of the factors $C^i$ and $C^j$ and the restrictions of the weight $\lambda$ to these factors, 
we see that
there are at least two $L_X'$-summands of $V^2(Q_Y)$ of highest weight $d+3f-5$, so $d+3f-5\leq d+f+1$, 
which implies that $f\leq 3$. If $f=3$, then the  usual
arguments show that there are at least four $L_X'$-summands of $V^2(Q_Y)$ with highest weight
 $d+2$, unless $d = f+1 = 4$, in which case one can obtain the precise decomposition of $V^2(Q_Y)$ and 
find three summands with highest weight $d+4$, contradicting Corollary \ref{conseq}. If $f=2$, where 
$V^1 = d+2\oplus d\oplus d-2$, 
in every configuration, we find that $V^2$ has at least 4 summands of 
highest weight $d+1$,
 contradicting Corollary \ref{conseq}. Hence $j=r-1$ and $C^j$ is of type $A_1$, as claimed.

Keeping $d$ as above, we have $V^1 = d+1\oplus d-1$. Note that $d>2$, otherwise
the weight $4=d+2$ occurs as the highest weight of three $L_X'$-summands of $V^2(Q_Y)$, 
contradicting Corollary \ref{conseq}. 

The weights $\lambda-\beta_{n-1}-\b_n$ and $\lambda-\beta_{n-2}-\beta_{n-1}$ afford summands of $V^2(Q_Y)$ whose sum restricts to $L_X'$ as $V^1\otimes 1$ and so as usual Corollary~\ref{cover} applies. 
If $\mu^i = (\lambda_1^i)^*$, we have
a summand of $V^2$ of the form $\wedge^2(d)\otimes (d-1)\otimes 1$, and if $\mu^i=\l_1^i$ with $d<r$ we have a summand $\wedge^2(d)\otimes (d+1)\otimes 1$.
Neither of these is MF, so we deduce that $\lambda\downarrow L_Y' = \lambda_1^0+\lambda_1^{r-1}$. Next we claim that $\langle\lambda,\gamma\rangle=0$ for all $\gamma\in\Pi(Y)\notin\Pi(L_Y)$, $\gamma\ne\beta_{n}$. For otherwise, we have a summand of $V^2$ of the form $S\in\{s\otimes(s-1)\otimes d\otimes 1, 1<s<r, d\otimes d\otimes (d-1)\otimes 1, d\otimes 2\otimes 2\}$, none of which is MF. So finally we have reduced to $\lambda=\lambda_1+\lambda_{n-1}+x\lambda_n$. If $x\ne 0$, comparing the dual module $V^*$ with Table~\ref{dist}, shows that $r=3=d$ and $x=1$. Now $V^1 = 3\otimes 1$ and $V^2(Q_Y)$ has summands afforded by $\lambda-\beta_1-\beta_2-\beta_3-\beta_4$, $\lambda-\beta_n$, $\lambda-\beta_{n-1}-\beta_n$ and $\lambda-\beta_{n-2}-\beta_{n-1}$,
whose restrictions to $L_X'$ are $2\otimes 1$, $3\otimes 2$, $3$, respectively $3\otimes 2$, which gives rise to four summands $3$, contradicting Corollary~\ref{conseq}.
So finally we have reduced to $\lambda = \lambda_1+\lambda_{n-1}$.

By Lemma~\ref{a2-examples-2}, we may assume $d=r\ge 4$.
Set $M = V^*$ of highest $\mu=\lambda_2+\lambda_{n}$, so $M^1 = \wedge^2(d)$.
 Now $M^2(Q_Y)$ has $L_X'$-summands 
$(2d-1)\oplus (2d-1)\oplus (2d-3)\oplus (2d-3)\oplus (2d-5)\oplus (2d-5)$, and all
other summands have lower highest weight. So we may now consider $M^3(Q_Y)$. Here we have 
$L_Y'$-summands afforded
by the weights $\nu_1 = \mu-\b_{n-2}-\b_{n-1}-\b_n$, 
$\nu_2 = \mu-\b_2-\cdots-\b_{r+1}-\b_n$, 
$\nu_3 = \mu-\b_2-\cdots-\b_{2r+1}$, and 
$\nu_4 = \mu-\b_1-2(\b_2+\cdots+\b_{r+1})-\b_{r+2}$.
These summands afford six $L_X'$-summands of highest weight $2d-4$, contradicting Corollary \ref{conseq}.

This completes the proof of the lemma.\hal

It remains to consider the configuration of Lemma~\ref{cases}(2), that is, when there is a unique $i$ with
$V^1(Q_Y)\downarrow C^i$ nontrivial. Set $d = r-i$, so that $C^i$ is of type $A_d$. First we 
collect some information in Table~\ref{tab:A1MF}. We will require the 
precise decompositions of the multiplicity free actions established in \cite{LST}. Most of these follow from Magma calculations. Items 2., 3. and 4. follow
from a direct weight count and items 14. and 15. are established in the proof 
of \cite[3.2]{LST}.
\begin{table}[h]
\caption{MF pairs for $A_1$ in $A_d$}\label{tab:A1MF}
\begin{tabular}{|l|l|l|}
\hline
reference&$(C^i,\lambda\downarrow C^i)$&$V^1$\\
number&&\\
\hline
\hline
1.&$(A_d,\l_1^i)$ & $d$\\
\hline
2.&$(A_d,\l_2^i)$ & $\oplus_{i=0}^{\lfloor (d-1)/2\rfloor} (2d-2-4i)$\\
\hline
3.&$(A_d,2\l_1^i)$ & $\oplus_{i=0}^{\lfloor d/2\rfloor} (2d-4i)$\\
\hline
4.&$(A_d,\l_1^i+\l_d^i)$&$\oplus_{i=0}^{d-1}(2d-2i)$\\
\hline
5.&$(A_5,\l_3^i)$&$9\oplus 5\oplus 3$\\
\hline
6.&$(A_6,\l_3^i)$ & $12\oplus 8\oplus 6\oplus 4\oplus 0$\\
\hline
7.&$(A_7,\l_3^i)$ &$15\oplus 11\oplus 9\oplus 7\oplus 5\oplus 3$\\
\hline
8.&$(A_3,3\l_1^i)$ & $9\oplus 5\oplus 3$\\
\hline
9.&$(A_4,3\l_1^i)$ &$12\oplus 8\oplus 6\oplus 4\oplus 0$\\
\hline
10.&$(A_5,3\l_1^i)$ &$15\oplus 11\oplus 9\oplus 7\oplus 5\oplus 3$\\
\hline
11.&$(A_3,4\l_1^i)$ &$12\oplus 8\oplus 6\oplus 4\oplus 0$\\
\hline
12.&$(A_3, 5\l_1^i)$ &$15\oplus 11\oplus 9\oplus 7\oplus 5\oplus 3$\\
\hline
13.&$(A_3, \l_1^i+\l_2^i)$ &$7\oplus 5\oplus 3\oplus 1$\\
\hline
14.&$(A_2, c\l_1^i)$ &$2c\oplus 2c-4\oplus\cdots\oplus 2c-4\lfloor\frac{c}{2}\rfloor$\\
\hline
15.&$(A_2, c\omega_1^i+\omega_2^i)$,  $c\geq 1$&$2c+2\oplus 2c\oplus\cdots\oplus 2$\\
\hline
\end{tabular}

\end{table}




\begin{lem}\label{1fctr-1} Let $\lambda$ be as in Lemma~$\ref{cases}(2)$. Then one of the following holds:\begin{itemize}
\item[{\rm (i)}] $i=0$, 
\item[{\rm (ii)}] $i=r-1$, or 
\item[{\rm (iii)}] $d=2$.
\end{itemize}
\end{lem}

\pf Suppose false, so that  $0<i<r-1$ and $2<d<r$. We must 
consider the various possibilities of Table~\ref{tab:A1MF}. 

Suppose the pair $(A_d,\lambda\downarrow C^i)$ is as in item 13  in Table~\ref{tab:A1MF}, or its dual.   
Then $d=3$ and one of the following holds:

\begin{enumerate}[a)]
\item $V^2(Q_Y)\downarrow L_Y'$ has summands with highest weights 
$(\l_1^{i-1})^*+ 2\l_2^i$ and  $ 2\l_1^i + \l_1^{i+1}$;
\item $V^2(Q_Y)\downarrow L_Y'$ has summands with highest weights 
$(\l_1^{i-1})^*+ (2\l_1^i)^*$ and $ 2\l_2^i + \l_1^{i+1}$.
\end{enumerate} 
In the first case, the restriction to $L_X'$ gives four summands of highest weight 8, and in the 
second
case five summands of highest weight 6. This contradicts Corollary \ref{conseq}.

This leaves us with $\lambda\downarrow C^i\in\{a\l_1^i,a\l_d^i,\l_2^i,\l_{d-1}^i,
\l_1^i+\l_d^i,\l_3^i,\l_{d-2}^i\}$. 

Suppose $\lambda\downarrow L_i = \l_1^i+\l_d^i$. Then $V^2(Q_Y)\downarrow L_Y'$ has a summand
with highest weight $(\l_1^{i-1})^*+ (\l_2^i + \l_d^i)$. Upon restriction to $L_X'$ we obtain 
two summands of highest weight $4d-3$; hence $4d-3\leq 2d+1$, contradicting $d>2$.

Consider now the case where $\lambda\downarrow C^i\in\{\l_2^i,\l_{d-1}^i\}$. 
Considering the two
cases separately, we see that if $d>3$, there are two $L_X'$-summands of $V^2(Q_Y)$ of highest 
weight $4d-9$ or $4d-11$. Hence $4d-11\leq 2d-1$ and so $d=4$ or 5. But if $d=4$, 
$V^2$ has summands $(5\otimes(6\oplus 2))\oplus (4\otimes 3)$ or 
$(5\otimes 4)\oplus ((6\oplus 2)\otimes3)$.
In each case, there are three $L_X'$-summands of highest weight 5, contradicting Corollary \ref{conseq}. If $d=5$,
 $V^2$ has summands $(6\otimes(9\oplus 5\oplus 3))\oplus (5\otimes 4)$ or 
$(6\otimes 5)\oplus ((9\oplus 5\oplus 3)\otimes4)$. In each case there are three $L_X'$-summands
of highest weight 9, again contradicting Corollary \ref{conseq}.

So finally, in case $\lambda\downarrow C^i\in\{\l_2^i,\l_{d-1}^i\}$, we have $d=3$. 
Here we have $V^1 = 4\oplus 0$. In addition, $\lambda-\b_{n-7}-
\b_{n-6}-\b_{n-5}$ and $\lambda -\b_{n-9}-\b_{n-8}-\b_{n-7}$ afford highest 
weights of $L_Y'$-summands of $V^2(Q_Y)$, which upon restriction to $L_X'$ give summands 
$(4\otimes 3)+(3\otimes 3)$. In particular, we deduce that 
$\langle\lambda,\gamma\rangle =0$ for all $\gamma\in\Pi(Y)\setminus \Pi(L_Y)$ (else there is a 
third $L_X'$-summand of $V^2(Q_Y)$ of highest weight 1). This then means that 
$V^2 = (4\otimes 3)+(3\otimes 3)$. We now consider 
$V^3(Q_Y)$. The weight $(\lambda-\b_{n-9}-\b_{n-8}-\b_{n-7}-\b_{n-6}-
\b_{n-5})\downarrow L_Y' = (\l_1^{i-1})^* + \l_2^i + \l_1^{i+1}$ affords an 
$L_Y'$-summand, which upon restriction to $L_X'$ gives a summand $4\otimes(4\oplus 0)\otimes 2$ 
which affords four summands of highest weight 6.
In addition, we have the weight $(\lambda-\b_{n-10}-2\b_{n-9}-2\b_{n-8}-
2\b_{n-7}-\b_{n-6})\downarrow L_Y' = (\l_2^{i-1})^*$ which affords a fifth 
$L_X'$-summand of highest weight 6, contradicting Corollary \ref{conseq}. This completes the consideration of the 
case $\lambda\downarrow C^i\in\{\l_2^i,\l_{d-1}^i\}$.

We now turn to the case where $\lambda\downarrow L_i = 2\l_1^i$ or $(2\l_1^i)^*$. 
Here there exists a weight $\mu\in V^2(Q_Y)$ such that $\mu\downarrow L_Y' =  (\l_1^{i-1})^*+
(\l_1^i + \l_2^i)$,
respectively $(\l_1^{i-1})^*+ (\l_1^i)^*$. As well we have $\nu\in V^2(Q_Y)$ such that 
$\nu\downarrow L_Y' = \l_1^i+ \l_1^{i+1}$,
respectively $(\l_1^i+\l_2^i)^* + \l_1^{i+1}$. Now in the first case, restriction of these two 
summands to
$L_X'$ affords two $L_X'$-summands of highest weight $4d-3$ and so $4d-3\leq 2d+1$, contradicting $d>2$. 
In the second
case, restriction of $\mu$ and $\nu$ to $L_X'$ affords two $L_X'$-summands of highest weight $4d-5$ and so 
$4d-5\leq 2d+1$; hence $d=3$. It is now a direct check to see that the weights $\nu$ and $\mu$ 
afford three $L_X'$-summands of $V^2(Q_Y)$ of highest weight $4d-5$, contradicting Corollary \ref{conseq}.

Consider the case where $5\leq d\leq 7$ and $\lambda\downarrow C^i = \l_3$, 
or $(\l_3^i)^*$. In the first case, we have an $L_Y'$-summand of highest weight 
$(\l_1^{i-1})^*+ \l_4^i$. This affords two $L_X'$-summands of highest weight $5d-15$, so 
$5d-15\leq 3d-5$ and hence $d=5$. But then $5d-15 = 10$ and we have an additional $L_X'$-summand 
of highest weight 10 afforded by 
an $L_Y'$-summand of highest weight $\l_2^i+\l_1^{i+1}$. This contradicts Corollary \ref{conseq}. In the 
second case, there is an $L_Y'$-summand of $V^2(Q_Y)$ with highest weight $(\l_1^{i-1})^*+(\l_2^i)^*$,
which upon restriction to $L_X'$ gives two summands with highest weight $3d-5$, 
and a further $L_Y'$-summand with highest weight $(\l_4^i)^*+\l_1^{i+1}$ affording a
third summand $3d-5$, again contradicting Corollary \ref{conseq}.

Now turn to the case where $3\leq d\leq 5$ and $\lambda\downarrow L_Y' = 3\l_1^i$ or
$(3\l_1^i)^*$. In the first case, there is an $L_Y'$-summand of $V^2(Q_Y)$ with highest weight 
$(\l_1^{i-1})^*+(2\l_1^i+\l_2^i)$. Upon restriction to $L_X'$ we obtain two summands with 
highest weight $5d-3$ which contradicts Corollary \ref{conseq} since $d>2$ and so $5d-3>3d+1$. In the second case,
 there is
an $L_Y'$-summand with highest weight $(\l_2^i+2\l_1^i)^*+\l_1^{i+1}$.
Upon restriction to $L_X'$ this affords two summands of highest weight $5d-5$, and hence $d=3$. But then 
$5d-5=3d+1$ and we have a third summand of highest weight $3d+1$ afforded by an $L_Y'$-summand with 
highest weight $(\l_1^{i-1})^*+(2\l_1^i)^*$, yielding the usual contradiction. 

In case $d=3$ and $\lambda\downarrow L_Y' = a\l_1^i$ or $(a\l_1^i)^*$, for $a\in\{4,5\}$, 
it is straightforward to produce the usual contradiction. We omit the details. Hence we are left with  
the case $\lambda\downarrow L_Y' = \l_1^i$ or $(\l_1^i)^*$. Here 
$V^1 = d$. 

Consider first the case where $\lambda\downarrow L_Y' = \l_1^i$; then there
exists an $L_Y'$-summand of $V^2(Q_Y)$ with highest weight $(\l_1^{i-1})^*+\l_2^i$. Upon restriction
to $L_X'$, this gives two summands of highest weight $3d-5$ and so $3d-5\leq d+1$ and $d=3$. In addition
 to the given $L_Y'$-summand, there is a second summand of highest weight $\l_1^{i+1}$, which affords
an $L_X'$-summand of highest weight 2. We now have $V^1 = 3$ and $V^2$
containing summands $8\oplus 6\oplus 4\oplus4\oplus 2\oplus 2\oplus 0$. It is now easy to see
that $\langle\lambda,\gamma\rangle = 0$ for all $\gamma\in\Pi(Y)\setminus\Pi(L_Y)$, else there is a third $L_X'$-summand of $V^2(Q_Y)$ of highest weight 2. So we now have $d=3$, $\lambda = \lambda_{n-8}$,
$V^2 =  8\oplus 6\oplus 4\oplus4\oplus 2\oplus 2\oplus 0$. We turn to $V^3(Q_Y)$. 
The weight $(\lambda-\b_{n-9}-\cdots -\b_{n-5})\downarrow 
L_Y' = ( \l_1^{i-1})^*+\l_1^i+\l_1^{i+1}$ affords an $L_X'$-summand $4\otimes 3\otimes 2$ and 
the weight $\lambda-\b_{n-10}-2\b_{n-9}-2\b_{n-8}-\b_{n-7})\downarrow L_Y' = 
(\l_2^{i-1})^*+(\l_1^i)^*$ affords an $L_X'$-summand $(6\oplus 2)\otimes 3$. Now one counts the number 
of
$L_X'$-summands of highest weight 5 and obtains the usual contradiction. 

Finally, consider the case where $\lambda\downarrow L_Y' = (\l_1^i)^*$. Here we have $L_Y'$-summands
 of $V^2(Q_Y)$ with highest weights $(\l_1^{i-1})^*$ and $(\l_2^i)^*+\l_1^{i+1}$. The second
summand upon restriction to $L_X'$ affords two summands of highest weight $3d-7$, so $3d-7\leq d+1$ and 
$d=3$ or 4. But if $d=4$, then $3d-7 = d+1$ and the first $L_Y'$-summand affords a third summand of 
highest  weight $d+1$, contradicting Corollary \ref{conseq}. Hence $d=3$. We now have $V^1 = 3$ and 
$V^2$
containing summands $6\oplus 4\oplus 4\oplus2\oplus 2$. As above it is easy to see
that $\langle\lambda,\gamma\rangle = 0$ for all $\gamma\in\Pi(Y)\setminus\Pi(L_Y)$, else there is a third
 $L_X'$-summand of $V^2(Q_Y)$ of highest weight 2. So we now have $d=3$, $\lambda = \lambda_{n-6}$,
$V^2 = 6\oplus 4\oplus 4\oplus2\oplus 2$.  We turn to $V^3(Q_Y)$. 
The weight $(\lambda-\b_{n-9}-\cdots -\b_{n-5})\downarrow 
L_Y' = (\l_1^{i-1})^*+(\l_1^i)^*+\l_1^{i+1}$ affords an $L_X'$-summand $4\otimes 3\otimes 2$; 
the weight $\lambda-\b_{n-7}-2\b_{n-6}-2\b_{n-5}-\b_{n-4})\downarrow L_Y' = 
\l_1^i+(\l_1^{i+1})^*$ affords an $L_X'$-summand $3\otimes 2$; the weight 
$(\lambda-\b_{n-6}-\cdots-\b_{n-2})\downarrow L_Y' = \l_2^i+\l_1^{r-1}$ affords an 
$L_X'$-summand $(4\oplus 0)\otimes 1$. Now one counts the number of
$L_X'$-summands of highest weight 5 and obtains the usual contradiction. This completes the proof of the lemma .\hal

\begin{lem}\label{1fctr-2} Let $\lambda$ be as in Lemma $\ref{cases}(2)$ with $i\not\in\{0,r-1\}$. Then either 
 $r=3$ and $\lambda=\lambda_5$ or $\lambda_6$, or $r=4$ and $\lambda=\lambda_{11}$. Hence $(\l,\d)$ are as in Table $\ref{TAB3}$ of Theorem $\ref{MAINTHM}$.
\end{lem}

\pf Assume false.
 Then the previous lemma implies that $d=2$, so $i=r-2$, and 
$\lambda\downarrow C^i \in\{c\l_1^i, c\l_2^i, c\l_1^i+\l_2^i, 
\l_1^i+c\l_2^i, c\geq 1\}$. 

Consider
first the case where $r=3$, so $i=1$.  If $\lambda\downarrow C^1 = 
c\l_1^1+\l_2^1$ or $\l_1^1+c\l_2^1$, for $c\geq 1$, then $V^1$ has 
summands $2c+2\oplus 2c$
and all other summands have lower highest weights. In the first case, $V^2(Q_Y)$ has
$L_Y'$-summands with highest weights $(\l_1^0)^*+((c-1)\l_1^1+2\l_2^1)$ and
$(c+1)\l_1^1+\l_1^2$. If $c>1$, the restriction of these summands to $L_X'$ affords 
three $L_X'$-summands of highest weight $2c+3$, contradicting Corollary \ref{conseq}. If 
$\lambda\downarrow C^1 = \l_1^1+c\l_2^1$, with $c>1$, then 
$V^2(Q_Y)$ has
$L_Y'$-summands with highest weights $(\l_1^0)^*+(c+1)\l_2^1$ and
$(2\l_1^1+(c-1)\l_2^1)+\l_1^2$. But then there are 4 $L_X'$-summands with 
highest weight $2c+1$, contradicting Corollary \ref{conseq}. So we have reduced to $\lambda = 
x\lambda_4+\lambda_5+\lambda_6+y\lambda_7+z\lambda_9$, and so $V^1 = 4\oplus 2$. 
If $x+y+z\ne 0$, we easily produce four $L_X'$-summands of $V^2(Q_Y)$ of highest weight 3, 
which gives the 
usual contradiction. Hence $x+y+z=0$. Now set $M= V^*$, of highest weight
$\mu=\lambda_4+\lambda_5$. Now $M^1 =2$, while $\mu-\b_4$ and 
$\mu-\b_4-\b_5$ 
produce 3 $L_X'$-summands of $M^2(Q_Y)$ of highest weight 3, contradicting Corollary \ref{conseq}.

Continuing with the case where $r=3$, we are left with $\lambda\downarrow C^1 = 
c\l_1^1$ or $c\l_2^1$. If $c>1$, we argue as above to reduce to the case $c=1$.
(There are at least three $L_X'$-summands of $V^2(Q_Y)$ of highest weight $2c-1$, whereas 
$V^1 = \oplus_{i=0}^{\lfloor c/2\rfloor} (2c-4i)$.)
So now $\lambda=x\lambda_4+\lambda_j+y\lambda_7+z\lambda_9$, for $j=5$ or 6, for some $x,y,z\geq 0$ and 
$V^1 = 2$  . 
It is completely straightforward to show that $x+y+z=0$, using the standard arguments, and
so the result holds.

For the remainder of the proof, we assume $d=2$ and $r>3$. A good part of the above analysis
goes through and 
we reduce to one of the following.
\begin{enumerate}[i.]
\item $\lambda\downarrow L_Y' = \l_1^i+\l_2^i$ and 
$\langle\lambda,\b_i\rangle=0$
for $n-8\leq i\leq n-5$ and for $i\geq n-2$;
\item $\lambda\downarrow L_Y' = \l_1^i$ and $\langle\lambda,\b_i\rangle=0$
 for  $n-8\leq i\leq n-5$ and for $i\geq n-3$; 
\item $\lambda\downarrow L_Y' = \l_2^i$ and $\langle\lambda,\b_i\rangle=0$ 
for $n-8\leq i\leq n-4$ and for $i\geq n-2$.
\end{enumerate}

Set $M = V^*$, with highest weight $\mu$. Now $M$ must also be as in Lemma~\ref{cases}. 
 Consulting the list of pairs in Table~\ref{tab:A1MF} for $(C^i,\mu\downarrow C^i)$,
 and applying Lemmas~\ref{v1-triv}, \ref{2fctrs-1wedge}, \ref{2fctrs-1sym}, 
\ref{2ntrl}, 
\ref{1fctr-1}, as well as the above arguments, we
see that $M$ multiplicity-free implies that one of the following holds.
\begin{enumerate}[i.]
\item $r=4$;
\item $r=5$ and $\mu\downarrow C^0 = \l_4^0$ or $\l_5^0$;
\item $r=6$ and $\mu\downarrow C^0 = \l_4^0$ or $\l_5^0$;
\item $r=7$ and $\mu\downarrow C^0 = \l_5^0$.
\end{enumerate}

If $5\leq r\leq 7$, Lemmas~\ref{2fctrs-1wedge},~\ref{2fctrs-1sym} and ~\ref{2ntrl}  show that 
$C^0$ is the unique factor of $L_Y'$ acting
nontrivially on $M^1(Q_Y)$.

In case $r=7$, we have $M^1 = 15\oplus 11\oplus 9\oplus 7\oplus 5\oplus 3$, 
while $M^2(Q_Y)\downarrow L_Y'$ has a summand with highest weight $\l_4^0+\l_1^1$, which affords
two $L_X'$-summands of highest weight 18, contradicting Corollary \ref{conseq}.

In case $r=6$,  we have $M^1= 10\oplus 6\oplus 2$, if $\mu\downarrow C^0 = \l_{5}^0$,
and $M^1 = 12\oplus 8\oplus 6\oplus 4\oplus 0$ if $\mu\downarrow C^0=\l_4^0$.
In the first case we find two $L_X'$-summands of $M^2(Q_Y)$ of highest weight 13,
and in the second case three summands of highest weight 11, giving the usual contradiction.

In case $r=5$, we have $M^1= 8\oplus 4\oplus 0$ if $\mu\downarrow C^0 = \l_4^0$ and
$M^1 = 5$ if $\mu\downarrow C^0=\l_5^0$. In the first case we find three 
$L_X'$-summands of $M^2(Q_Y)$ of highest weight 7,
and in the second case two summands of highest weight 8, giving the usual contradiction.

So finally, we have reduced to the case $r=4$. Here, by applying all previously established results
(including the earlier results in this proof) to both $V$ and $M$,
we deduce that either $\lambda=\lambda_{11}$ or one of the following holds:
\begin{enumerate}[a)]
\item $\mu = \lambda_4+\lambda_5$;
\item $\mu = \lambda_5+\lambda_{10}$;
\end{enumerate}
In the first case, $M^1=4$, while $\mu-\b_5$ and $\mu-\b_4-\b_5$ afford
three $L_X'$-summands of $M^2(Q_Y)$ of highest weight 5, contradicting Corollary \ref{conseq}. For $\mu$ as in b), 
$M^1 = 2$. Now $\mu-\b_5$ affords an $L_X'$-summand $4\otimes3\otimes2$ of 
$M^2(Q_Y)$, which produces two $L_X'$-summands
of highest weight 7, contradicting Corollary \ref{conseq}. \hal


\begin{lem}\label{1fctr-3} Let $\lambda$ be as in Lemma $\ref{cases}(2)$, with $i=0$ or $i=r-1$.
Then the pair $(r,\lambda)$ 
or $(r,\lambda^*)$ is as in Table $\ref{tab:A2An}$.
\end{lem}

\pf
\noindent{\bf{Case I.}} Assume $i=r-1$ and so $C^i = A_1$. 

In particular,
$V^1 = c$, for some $c\geq 1$. Recall that $r>2$. Consider $M = V^*$, the irreducible
 $Y$-module of highest weight $\mu$. Since $M^1$ is multiplicity-free, one of the following holds.
\begin{itemize}
\item[i.] $c=1$ and $\langle\lambda,\b_n\rangle=0=\langle\lambda,\b_{n-2}\rangle$;
\item[ii.] $r=3$, $c=1=\langle\lambda,\b_n\rangle$ and $\langle\lambda,\b_{n-2}\rangle=0$;
\item[iii.] $r=3$, $c=1=\langle\lambda,\b_{n-2}\rangle$ and $\langle\lambda,\b_{n}\rangle=0$.
\end{itemize} 

So in all cases, $V^1=1$.
In the first case, we claim that $\lambda = \lambda_{n-1}$. 
Otherwise, there exists 
$\gamma\in\Pi(Y)\setminus \Pi(L_Y)$ with $\langle\lambda,\gamma\rangle\ne 0$, $\gamma\not\in\{\b_{n-2},\b_n\}$, by assumption. But then 
$(\lambda-\gamma)\downarrow L_X'=s\otimes (s-1)\otimes 1$, for some $s\geq 3$,
 and $\lambda-\b_{n-2}-\b_{n-1}$ produce three 
$L_X'$-summands
of $V^2(Q_Y)$ of highest weight 2, contradicting Corollary \ref{conseq}. Hence $\lambda = \lambda_{n-1}$ and the result holds.

In the second case, the usual argument gives that $\lambda = \lambda_8+\lambda_9$
and the result holds. Finally, in the third case, the usual considerations show that $\lambda = \lambda_7+\lambda_8$, and the result follows from Lemma~\ref{a2-examples-2}. This completes Case I.

\noindent{\bf {Case II.}} Assume $i=0$ and $r=3$.  

If $\lambda = \sum_{j=1}^9 a_j\lambda_j$,
by considering the dual module $M = V^*$, and the previously treated cases, we may assume $a_j=0$ for $j=2,8$ and
$4\leq j\leq 6$. There are various possibilities for $\lambda$ arising from the list of pairs 
$(C^0,\lambda\downarrow C^0)$. 

If $\lambda\downarrow C^0=(a\l_1^0)\downarrow C^0$, for $2\leq a\leq 5$, then  $V^1(Q_Y)$ has $L_X'$-summands 
$3a\oplus 3a-4$, one summand of highest weight $3a-6$ if $a\geq 3$, and all other summands have lower highest weights. Now if $a_7\ne 0$, then 
$(\lambda-\b_7)\downarrow L_Y' = a\l_1^0+\l_2^1+\l_1^2$ affords three $L_X'$-summands of 
highest weight $3a-3$, contradicting Corollary \ref{conseq}.
Hence $\lambda = a\lambda_1+x\lambda_9$. By Lemma~\ref{a2-examples-2}, we may assume 
$a\geq 3$ or $x\geq 3$; without loss of generality, we assume $a\geq 3$.
 But now $V^2(Q_Y)$ has $L_Y'$-summands
with highest weights $(a-1)\l_1^0+\l_1^1$ and $a\l_1^0+\l_1^2$. These afford four 
$L_X'$-summands of highest weight $3a-5$, contradicting Corollary \ref{conseq}.

If $\lambda\downarrow C^0=(a\lambda_3^0)\downarrow C^0$, then the usual arguments show that 
$\lambda = a\lambda_3$. But now 
Lemma~\ref{v1-triv} applied to $M = V^*$ gives the result.

Consider now the case where $\lambda\downarrow C^0 = \lambda_1^0+\lambda_3^0$, so 
$V^1 = 6\oplus 4\oplus 2$. Then $(\lambda-\b_3-\b_4)\downarrow 
L_Y' = (\l_1^0+\l_2^0)+\l_1^1$, which affords
three $L_X'$-summands of $V^2(Q_Y)$ of highest weight 5. By Lemma~\ref{a2-examples-2},
we may assume $\langle\lambda,\b_7\rangle+\langle\lambda,\b_9\rangle\ne 0$, which gives
a fourth summand of weight 5, giving the usual contradiction.

So for Case II, we have reduced to $\lambda\downarrow C^0 = \lambda_1^0$, so 
$\lambda = \lambda_1+a_7\lambda_7+a_9\lambda_9$. If $a_7\ne 0$, the dual module, $M= V^*$,
has been handled above. So we may assume $a_7=0$. By Lemma~\ref{a2-examples-2}, we may assume $a_9\geq3$, and then the module $M = V^*$ is a configuration which
has been previously considered. This completes Case II.

\noindent{\bf{Case III.}} Assume $i=0$ and $r\geq 4$.

Consider first the case where $\lambda\downarrow L_Y' = 3\omega_1^0$, so by Table~\ref{tab:A1MF}, $r=4$ or 5. In each case,
it is straightforward to compare the $L_X'$-summands in $V^1(Q_Y)$ and those in $V^2(Q_Y)$ and applying Corollary \ref{conseq},
we see that $\lambda = 3\lambda_1$, in which case $r,\l$ are as in Table \ref{tab:A2An}.

The case $\lambda\downarrow C^0 = 3\l_r^0$ for $r=4,5$
 is ruled out by applying Lemma ~\ref{1fctr-2} to $V^*$. 
Now we turn to the configurations where $\lambda\downarrow C^0 = \l_3^0$ or $(\l_3^0)^*$; 
in particular,
by Table~\ref{tab:A1MF}, we have $5\leq r\leq 7$. Using the standard techniques, we show that 
if $\lambda= \l_3^0$, then $\lambda = \lambda_3$ and the result follows from Lemma~\ref{a2-examples-2}. 
If $\lambda\downarrow C^0 = (\l_3^0)^*$, so we may assume $r=6$ or 7, then $V^1(Q_Y)$ has 
$L_X'$-summands $3r-6\oplus 3r-10\oplus 3r-12 
\oplus 3r-14$ and all other summands have lower highest weights. There exists a weight 
$\mu\in V^2(Q_Y)$ affording an $L_Y'$-summand with highest weight 
$(\l_4^0)^*+\l_1^1$, which
in turn restricts to $L_X'$ to give two $L_X'$-summands of highest weight $5r-17$. 
Hence $5r-17\leq 3r-5$ and so $r=6$. But now we can be more precise: 
$V^1 = 12\oplus 8\oplus 6\oplus 4\oplus 0$ and the weight 
$\mu$ affords an $L_X'$-summand $(12\oplus 8\oplus 6\oplus 4\oplus 0)\otimes 5$,
which produces four summands of highest weight 9, contradicting Corollary \ref{conseq}. 

Consider now the case where $\lambda\downarrow C^0 = \l_1^0+(\l_1^0)^*$. 
Here $V^1 = 2r\oplus 2r-2\oplus\cdots\oplus 2$. 
The weight $\lambda-\b_r-\b_{r+1}$ affords two $L_X'$-summands of $V^2(Q_Y)$ of 
highest weight $4r-5$ and so $4r-5\leq 2r+1$ contradicting $r>3$.

Now consider the case where $\lambda\downarrow C^0 = 2\l_1^0$; in particular 
$V^1 = 2r\oplus 2r-4\oplus\cdots\oplus 2r-4\lfloor(r/2)\rfloor$.
Now using the usual aruments comparing the $L_X'$-summands of $V^1(Q_Y)$ and those of $V^2(Q_Y)$,
we deduce that $\lambda=2\lambda_1+x\lambda_n$, for some $x\geq 0$. Now setting $M = V^*$,
we deduce that $x\leq 2$, since we have already treated the case where $r\geq 4$ and $x\geq 3$
above. This then leaves us with $\lambda = 2\lambda_1$, which is in Table \ref{tab:A2An},  or $\lambda = 2\lambda_1+\lambda_n$, or  $\lambda = 2\lambda_1+2\lambda_n$. We must treat the latter two cases.
We can now determine precisely $V^2(Q_Y)$; $V^2(Q_Y)\downarrow L_Y'$ has exactly
two irreducible summands with highest weights $\l_1^0+\l_1^1$ and 
$2\l_1^0+\l_1^{r-1}$. The restriction of these summands to $L_X'$ affords 
$L_X'$-summands $r\otimes (r-1)$ and $2r\otimes 1$, $(2r-4)\otimes 1,\dots,
(2r-4\lfloor(r/2)\rfloor)\otimes 1$. Let us now consider $V^3(Q_Y)$. There is an 
$L_Y'$-summand with highest weight $2\l_1^1$, another with highest weight 
$2\l_1^0+\l_1^{r-2}$,
another with highest weight $\l_1^0+\l_1^2$, and another with highest weight
 $\l_1^0+\l_1^1+\l_1^{r-1}$. These summands produce six $L_X'$-summands of 
highest weight $2r-2$, which
exceeds the number allowed by Corollary \ref{conseq}. This completes the consideration of the case  
$\lambda\downarrow C^0 = 2\l_1^0$.

If $\lambda\downarrow C^0  = 2\l_r^0$, $V^1(Q_Y)$ is as in the previous case. However, 
$(\lambda-\b_r-\b_{r+1})\downarrow L_Y' = (\l_1^0+\l_2^0)^*+\l_1^1$ affords an 
$L_Y'$-summand of $V^2(Q_Y)$ whose restriction to $L_X'$ gives two $L_X'$-summands of highest weight $4r-5$.
Hence $4r-5\leq 2r+1$, contradicting the assumption that $r>3$.

Suppose now that $\lambda\downarrow C^0 = \l_2^0$, so 
$V^1 = 2r-2\oplus 2r-6\oplus\cdots\oplus 2r-4\lfloor(r-1)/2\rfloor$.
 Now comparing the $L_X'$-summands of $V^1(Q_Y)$ and those of $V^2(Q_Y)$, 
we deduce that $\lambda = \lambda_2+x\lambda_n$, for some $x\geq 0$. If $x=0$, we have 
one of the examples in Table \ref{tab:A2An}. If $x\geq 1$, then we consider 
the dual module $M=V^*$ and see that
we are in a case treated in Lemmas~\ref{2fctrs-1sym} or \ref{2ntrl}. Now if 
$\lambda\downarrow L_Y' = (\l_2^0)^*$, then $V^1$ is as in the
previous case. Here the weight $(\lambda-\b_{r-1}-\b_r-\b_{r+1})\downarrow L_Y' = 
(\l_3^0)^*+\l_1^1$, affords two $L_X'$-summands of highest weight $4r-11$. 
Hence $4r-11\leq 2r-1$, so $r=4$ or 5. If $r=5$,
then $M = V^*$ has been treated in Lemma \ref{1fctr-3}. If $r=4$, then $V^1 = 6\oplus 2$ and our usual considerations show that $\lambda = \lambda_3$. This is
one of the examples in Table \ref{tab:A2An}.

So finally it remains to consider the case where $\lambda\downarrow C^0 = 
\l_1^0$ or $(\l_1^0)^*$. In each case $V^1 = r$. In the second case, 
the weight $(\lambda-\b_r-\b_{r+1})\downarrow L_Y' = (\l_2^0)^*+\l_1^1$ affords
two $L_X'$-summands of $V^2(Q_Y)$ of highest weight $3r-7$. Hence $3r-7\leq r+1$ and so $r=4$.
Here the module $M = V^*$ has  been treated in an earlier case. Thus, 
we are left with $\lambda\downarrow L_Y' = \l_1^0$. Comparing the $L_X'$-summands
in $V^2(Q_Y)$ with $V^1(Q_Y)$, we deduce in the usual way that $\lambda = \lambda_1+x\lambda_n$,
for some $x\geq 0$. Note that $x \ne 0$ as $\l\ne \l_1$ by hypothesis, and if $x=1$, then we have a configuration of 
Table \ref{tab:A2An}; otherwise, we note that the module $M=V^*$ 
has been treated earlier.

This completes the proof of the lemma. \hal

\vspace{4mm}
This completes the consideration of cases 1, 2 and 3 of Lemma~\ref{cases}, and therefore we have established 
Theorem \ref{a2r0} in the case where $r \ge 3$.

\vspace{4mm}
We have now finished the proof of Theorems \ref{rgreater} and \ref{a2r0}, which 
constitute Theorem \ref{MAINTHM} in the case where $X = A_2$.

\chapter{The case $\d = r\omega_k$ with $r,k\ge 2$}\label{rkge2}

In this chapter  we consider the case where $X = A_{l+1}$ for $l \ge 2$ and $X$ is embedded in $Y = SL(W)=A_n$ via the
representation of highest weight $\d = r\omega_k$ with $r>1$ and $l+1 > k > 1.$  
Using the notation established in Chapter \ref{levelset}, 
we let $\Pi(L_X') = \{\alpha_1, \dots, \alpha_l\}$ and $\alpha = \alpha_{l+1}.$  Replacing  the embedding with  its dual, if necessary, we can assume $k \le \lbrack \frac{l+2}{2}\rbrack.$ Then $L_X'  < L_Y' = C^0 \times C^1 \times \dots \times C^r$  is a product of subgroups of type $A$, where for each $i$ the projection, $\pi_i(L_X'),$ of $L_X'$ to $C^i$ corresponds to the
action of $L_X'$ on the $i$th level of $W$. It follows from Theorem \ref{LEVELS}
that the action of $L_X'$ on the $i$th level is irreducible with highest weight $\delta^i
= i\omega_{k-1} + (r-i)\omega_k$.  Let $\lambda_1, \lambda_2, \dots, \lambda_n$ be the fundamental dominant weights of $Y.$

 Let $V = V_Y(\lambda)$ be a nontrivial irreducible $Y$-module.  We have  
$V^1(Q_Y) \downarrow L_Y' = V_{C^0}(\mu^0) \otimes \dots \otimes V_{C^r}(\mu^r) $, where $\mu ^i$ is the restriction of $\lambda$ to $T_Y \cap C^i.$   For each $i$  let the fundamental
system of $C^i$ be $\Pi(C^i) = \{\beta_1^i, \dots, \beta_{r_i}^i\}$ with corresponding fundamental dominant weights $\lambda_j^i.$  Let $\lambda^*$ and  $(\mu^i)^*$ denote  the dual of $V$ and
 $\mu^i,$ respectively. 

We recall the following notation from Chapter \ref{levelset}: for $\eta = \sum x_i \omega_i$ a dominant weight of a simple group algebraic group $G$, let $S(\eta) = \sum x_i.$

In this section we prove Theorem \ref{MAINTHM} for this case, under the inductive assumption that Theorem \ref{MAINTHM} holds in general for rank smaller than $l+1$.

\begin{theor}\label{main} Let $X = A_{l+1}$, $W = V_X(r\o_k)$, $Y = SL(W)$ and assume $V_Y(\l)\downarrow X$ is multiplicity-free, with $\l \ne 0,\l_1,\l_n$. Assume also that the conclusion of Theorem $\ref{MAINTHM}$ holds for groups $A_m$ of rank $m<l+1$. 
Then replacing $V$ by $V^*$ if necessary, either 
\begin{itemize}
\item[{\rm (i)}] $\lambda \in \{ \lambda_2, 2\lambda_1, \lambda_1+\lambda_n\}$, or 
\item[{\rm (ii)}] $k=r= 2$ and $\lambda = \lambda_3$.
\end{itemize}
\end{theor}

Note that in both cases (i) and (ii) of the theorem, $V \downarrow X$ is MF by Theorem \ref{yesMF}.

To establish the theorem we will separate the cases $l > 2$ and $l = 2$, beginning  
with $l > 2.$ Then we will briefly discuss the changes required to settle the case $l = 2$ with $r > 2.$ The
final case is $l = r = 2$ and this will require some extra effort.

The first lemma establishes a comparison between the ranks of $C^0$ and $C^r.$

\begin{lemma}\label{dimineq}	 Assume the hypotheses of Theorem $\ref{main}$.
\begin{itemize}
\item[{\rm (i)}] Suppose $l$ is even and $k = \frac{l+2}{2}$. Then 
$\dim V_{L_X'}(r\omega_k) = \dim V_{L_X'}(r\omega_{k-1}).$
\item[{\rm (ii)}]  If the assumption in {\rm (i)} does not hold, then 
$\dim V_{L_X'}(r\omega_k) > \dim V_{L_X'}(r\omega_{k-1}) + 4.$
\end{itemize}
 \end{lemma}
 
 \pf  (i)  If $l$ is even with $k =  \frac{l+2}{2},$ then $k + (k-1) = l+1$ so that  $V_{L_X'}(r\omega_k)$ and $ V_{L_X'}(r\omega_{k-1})$ are dual representations and hence have the same dimension.  

(ii) Suppose the assumption of (i) does  not hold.  Then $k \le \frac{l+1}{2}.$  Here we use the Weyl degree formula to compare dimensions.  Setting
 $x = \dim V_{L_X'}(r\omega_{k-1})$ we find that $\dim V_{L_X'}(r\omega_k) \ge \frac{r+k}{k}x.$ Therefore, (ii) holds provided $x > \frac{4k}{r}.$  But $x \ge \frac{l(l+1)}{2}$ so this
 holds. \hal

Assume the hypotheses of Theorem \ref{main}. Then Lemma \ref{induct} implies that for each $i$, $V_{C^i}(\mu^i)\downarrow L_X'$ is MF, and $L_X'$ embeds into $C^i$ via the highest weight $\d^i$ given above. Hence the inductive assumption in the statement of Theorem \ref{main} implies the following.

\begin{lemma}\label{induc} Assume $l>2$.  
\begin{itemize}
\item[{\rm (i)}] If $1 < i < r-1,$  then $\mu^i $ or $(\mu^i)^* \in \{ 0,\,\lambda_1^i, \,   \lambda_2^i \,(i = k=2, r=4)\}.$
\item[{\rm (ii)}] If $i = 1, r-1,$ then  $\mu^i$ or $(\mu^i)^* \in \{0, \,\lambda_1^i, \, \lambda_2^i, \,2 \lambda_1^i \}$.
\item[{\rm (iii)}] If $i = 0, r$, then $\mu^i$ or $(\mu^i)^* \in \{0,\, \lambda_1^i, \, \lambda_2^i, \,2 \lambda_1^i,  \, \lambda_1^i+ \lambda_{r_i}^i, \, \lambda_3^i \,(k=r=2)\}.$
\end{itemize}
\end{lemma}

 \pf  This follows from the induction hypothesis except for (iii) when $i = r$ and $k = 2$, where there are a number of additional possibilities  for $\mu^r$ or its dual.  These include
$c\l_1^r,$ $\l_i^r$,  $\l_1^r+\l_2^r$ ($r = 3$), and all the Table 2 possibilities when $r = 2.$ The extra cases are ruled out using  the induction hypothesis applied to $(\mu^*)^0 \downarrow L_X'$  along with an application of Lemma \ref{dimineq}. \hal

   The next lemma provides some composition factors of certain wedge and symmetric powers of $r\omega_k.$

\begin{lemma}\label{wedgefactors}  Let $X =   A_{l+1}$ and $k \le \frac{1}{2}(l+1).$ 
\begin{itemize}
\item[{\rm (i)}] $\wedge^2(V_X(r\omega_k))$ has a summand $V_X(\omega_{k-1} +(2r-2)\omega_k + \omega_{k+1})$.
\item[{\rm (ii)}] $\wedge^3(V_X(2\omega_k))$ contains $V_X(\omega_{k-2} + \omega_{k-1} +2\omega_k + \omega_{k+1}+\omega_{k+2})$ as a summand (noting that the first term does not appear if $k = 2$ and the last term does not occur if $l = 2$).
\item[{\rm (iii)}] $S^2(V_X(r\omega_k))$ has a summand $V_X(2r\omega_k) + V_X(2\omega_{k-1} + (2r-4)\omega_k + 2\omega_{k+1})$. 
\end{itemize}
\end{lemma}

\pf Parts (i) and (iii) follow by noting that $\wedge^2(V_X(r\omega_k)) \supseteq 2r\omega_k-\a_k$ and $S^2(V_X(r\omega_k))\supseteq 2r\omega_k-2\a_k.$ For (ii), first use Magma to see that the assertion
holds  for the case $k = 2$ and $l = 3,$  where we find that there is a composition factor of highest weight
$6\omega_k -(\a_{k-1}+3\a_k+\a_{k+1}).$  A consideration of Levi factors implies that the same holds in larger rank configurations.  This yields (ii). \hal

Now we handle a useful special case with $l=2$.

  \begin{lemma}\label{020}  Let $X = A_3$, $W = V_X(020)$ and $X < Y = SL(W) = A_{19}$. Let $a \ge 1$ and let $V = V_Y(\l)$ with $\l = a\l_1 + \l_2$ or  $\l_4 + a\l_5.$ Then $V \downarrow X$ is not MF.
\end{lemma}

\pf  If $a\ge 4$, then $V^1$ is not MF for $L_X' = A_2$, by Theorem \ref{a2r0}. The case where $\l = a\l_1+\l_2$ with  $a=1,2,3$ is handled using Magma.  Finally, assume $\l =\l_4 + a\l_5$ with $a \le 3.$  The decomposition of $V^1$ is given in Lemma \ref{A2info}. We can take  the embedding of $L_X'$ in $C_0$  so that $\b_1^0, \b_2^0, \b_3^0, \b_4^0,  \b_5^0$ restrict to $S_X$ as $\a_2, \a_2, \a_1-\a_2, \a_2, \a_1,$  respectively. Now $V^2$ contains 
$((0002(a-1)) \downarrow L_X') \otimes 11$ and $((0010a) \downarrow L_X') \otimes 11$, and it follows from the above that $\l_5^0, \l_4^0, \l_3^0$ restrict to $L_X'$ as $(20), (21), (30)$, respectively.  Then one checks that 
$V^2 \supseteq ((2a+2,2) + (2a+3,0)) \otimes (11)$ and this contains $(2a+2,2)^3$.
But at most one such factor can arise from $V^1.$
\hal

\section{Case $l>2$ } We now begin the proof of Theorem \ref{main}. In this subsection we prove the theorem under the assumption that $l>2$. Assume then that $l>2$.

In what follows, we shall refer to the inductive assumption in the statement of Theorem \ref{main} as the Inductive Hypothesis.
   
 \begin{lem}\label{atmost1}  There exists at most one  value of $i$ such that  $ 0 < i < r$ and $\mu^i \ne 0.$
\end{lem} 

\pf  Suppose both $i \ne j$ satisfy the conditions $ 0 < i,j < r$  and $\mu^i \ne 0 \ne \mu^j.$  Then $\delta^i = i\omega_{k-1} + (r-i)\omega_k$ and  $\delta^j = j\omega_{k-1} + (r-j)\omega_k. $   It follows from Lemma \ref{induc}, Lemma \ref{lemma1.5(i)(ii)}, and Proposition \ref{tensorprodMF}  that  $(V_{C^i}(\mu^i) \otimes V_{C^j}(\mu^j)) \downarrow L_X'$ is not MF and hence neither is $V^1,$ a contradiction.  \hal
 
 \begin{lem}\label{no1<i<r-1}  
 \begin{itemize}
\item[{\rm (i)}]  $\mu^i = 0$ if $1 < i < r-1.$
 \item[{\rm (ii)}]  $\la \lambda, \gamma_i \ra = 0$  if $1 < i \le r-1.$
\end{itemize}
 \end{lem}
 
 \pf  (i)  By way of contradiction assume $1 < i < r-1$ but  $\mu^i \ne 0$.  Note that this forces $r \ge 4.$ The possibilities for $\mu^i$ are given by Lemma \ref{induc}.

First assume $\mu^i = \lambda_1^i$.  Then Lemma \ref{atmost1} implies that  
 $$V^1(Q_Y) \downarrow L_Y'  = V_{C^0}(\mu^0) \otimes V_{C^i}(\lambda_1^i)\otimes V_{C^r}(\mu^r).$$
   Also
  $\lambda - \gamma_i -\beta_1^i$ affords the $L_Y'$-module 
  $$V_{C^0}(\mu^0) \otimes V_{C^{i-1}}(\lambda_{r_{i-1}}^{i-1})\otimes V_{C^i}(\lambda_2^i) \otimes V_{C^r}(\mu^r)$$
  in $V^2(Q_Y)$.  
 
We consider the restriction to $L_X'$ of the tensor product of the middle two terms.  We have $(V_{C^{i-1}}(\lambda_{r_{i-1}}^{i-1})\otimes V_{C^i}(\lambda_2^i)) \downarrow L_X' = (V_{L_X'}(\delta^{i-1}))^* \otimes
\wedge^2(V_{L_X'}(\delta^i)).$  The first factor has highest weight $\eta = (r-i+1)\omega_{l-k+1} + (i-1)\omega_{l-k+2}.$
The second factor has a direct summand of highest weight 
$\xi =\omega_{k-2} +(2i-2)\omega_{k-1} + (2r-2i+1)\omega_k$ (the term $\omega_{k-2}$ is not present if $k = 2).$ As $k \le \frac{1}{2}(l+2)$ we can then apply Lemma \ref{abtimescd} to
the tensor product of these factors; setting  $(a,b) = (2i-2,2r-2i+1)$ and $(c,d) = (r-i+1, i-1)$ in Lemma \ref{abtimescd} we obtain  a direct summand of multiplicity at least $2$ and having highest weight  $\eta + \xi -(\alpha_{k-1} + 2\alpha_k + \dots +2\alpha_{l-k+1} + \alpha_{l-k+2}).$  Now, tensoring with the remaining factors $V_{C^0}(\mu^0), V_{C^r}(\mu^r)$  we conclude that there is a direct summand $V_{L_X'}(\rho)$ of $V^2$ having multiplicity at least 2 and satisfying $S(\rho) \ge  S(\mu^0 \downarrow L_X') + 3r-3 +S(\mu^r \downarrow L_X') $.  On the other hand the $S$-value of highest weight
of an $L_X'$ irreducible summand of $V^1(Q_Y)$ is $S(\mu^0 \downarrow L_X') + r +S(\mu^r \downarrow L_X')$ so this  contradicts Lemmas \ref{induct} and \ref{weightsum}.

If $\mu^i = \lambda_{r_i}^i$ we proceed as above, working with $\gamma_{i+1}.$ This leaves  the special
case  where $i = k = 2, r = 4$ and  $\mu^i = \lambda_2^i$ or its dual.  We give details for the first of these 
noting that the other case is essentially the same.  Then in $V^2(Q_Y)$ there is an $L_Y'$- irreducible summand of highest weight $\lambda -\gamma_i - \beta_1^i - \beta_2^i$ which affords the module 
$$ V_{C^0}(\mu^0) \otimes V_{C^{i-1}}(\lambda_{r_{i-1}}^{i-1})\otimes V_{C^i}(\lambda_3^i) \otimes V_{C^r}(\mu^r).$$
An easy weight count (or use Magma) shows that $V_{C^i}(\lambda_3^i) \downarrow L_X'$ contains the irreducible of highest weight
$(660 \dots 00) - \alpha_1-2\alpha_2 = (6320\dots 0)$ with multiplicity $2$.  Therefore, by Lemma \ref{wedgefactors}  there is a direct summand $V_{L_X'}(\rho)$ of 
$V^2$ having multiplicity at least 2 and satisfying $S(\rho) \ge  S(\mu^0 \downarrow L_X') + r + 11 +S(\mu^r \downarrow L_X'). $  An upper bound for the $S$-value of $V^1(Q_Y) \downarrow L_X'$ is $S(\mu^0 \downarrow L_X') + r + 8 +S(\mu^r \downarrow L_X')$, so this again contradicts Lemmas \ref{induct} and \ref{weightsum}.

(ii)  Suppose $1 < i \le r-1$ but $\la \lambda, \gamma_i \ra  \ne 0.$   In view of (i) we have
$$V^1(Q_Y) \downarrow L_Y'  = V_{C^0}(\mu^0) \otimes V_{C^1}(\mu^1)\otimes V_{C^{r-1}}(\mu^{r-1}) \otimes V_{C^r}(\mu^r),$$
noting that at most one of the middle two terms is nonzero.  First assume that $i \ne 2,r-1.$   Then $V_{\gamma_i}^2(Q_Y)$ contains 
$$V^1(Q_Y) \otimes V_{C^{i-1}}(\lambda_{r_{i-1}}^{i-1})\otimes V_{C^i}(\lambda_1^i).$$
Using Lemma \ref{abtimescd} we see that the tensor product of the last two terms contains an irreducible summand
of multiplicity two  with $S$-value at least $2r-2.$  Then Lemmas \ref{induct} and \ref{weightsum} give a contradiction. 

Now assume $i = 2$ or $i = r-1.$   Then $r \ge 3$ and if $r = 3$ then $i = 2 = r-1.$  For these cases an adjustment to the above argument is required if $i = 2$ and $\mu^1 \ne 0 $ or if $i = r-1$ and  $\mu^{r-1} \ne 0. $
We will work out the first case as the second case is essentially the same.  So suppose $i= 2$ and $\mu^1 \ne 0.$
The term $V_{C^{1}}(\lambda_{r_{1}}^1)$  must be replaced by $V_{C^{1}}(\mu^1 + \lambda_{r_{1}}^1).$ The possible choices for $\mu^1$ are $\l_1^1$, $\l_2^1$, $2\l_1^1$,  or the dual of one of these. In each case we can write down an explicit highest weight of a composition factor of $V_{C^{1}}(\mu^1 + \lambda_{r_{1}}^1) \downarrow L_X'.$  Indeed there is a maximal vector of weight $(\mu^1 +\l_{r_1}^1) \downarrow S_X.$ Lemma \ref{abtimescd} still applies to show that there is an $L_X'$-composition factor in 
$V_{C^{1}}(\mu^1 + \lambda_{r_{1}}^1) \otimes V_{C^i}(\lambda_1^i)$ having multiplicity 2 and
$S$-value at least $S((\mu^1 + \lambda_{r_{1}}^1) \downarrow L_X') + r - 2$.
 Moreover, Lemma \ref{p.29} shows that $S((\mu^1 + \lambda_{r_{1}}^1) \downarrow L_X') \ge S(\lambda_{r_{1}}^1 \downarrow L_X') + r.$   Therefore $V_{C^{1}}(\mu^1 + \lambda_{r_{1}}^1) \otimes V_{C^i}(\lambda_1^i)$ has a composition factor of multiplicity at least $2$ and $S$-value at least $2r-2$ which gives the same contradiction as before.  \hal

At this point we are reduced to  
$$V^1(Q_Y) \downarrow L_Y'  = V_{C^0}(\mu^0) \otimes V_{C^1}(\mu^1)\otimes V_{C^{r-1}}(\mu^{r-1}) \otimes V_{C^r}(\mu^r).$$

Moreover  Lemma \ref{atmost1} implies that at most one of the two middle terms is nontrivial.

\begin{lem}\label{mu_r}  Assume  $k < \frac{l+2}{2}$.  Then $\mu^r \ne \lambda_1^r, \lambda_2^r, \lambda_1^r+\lambda_{r_r}^r, \lambda_3^r$ or $2\lambda_1^r.$
\end{lem}

\pf Assume $k < \frac{l+2}{2}$.  Then the result follows from  Lemma \ref{dimineq} and consideration of the dual of $V$  as otherwise we contradict the Inductive Hypothesis.  
 For example, if $l =3$ and $r = k = 2$, then $C^0 = A_{19}$ and $C^r = A_{9}$.  If  $\mu^r = \lambda_1^r, \lambda_2^r, \lambda_1^r+\lambda_{r_r}^r, \lambda_3^r$ or $2\lambda_1^r, $ we would have $(\mu^*)^0 = \lambda_9^0, \lambda_8^0, \lambda_1^0+\lambda_9^0, \lambda_7^0$ or $2\lambda_9^0$, respectively.  \hal

\begin{lem}\label{mu1notomega1etc}   $\mu^1 \ne \lambda_1^1, \lambda_2^1$ or $2\lambda_1^1$. 
\end{lem} 

\pf  By way of contradiction assume the result is false. First assume $\mu^1 = \l_1^1$. We first claim $\mu^0$ and $\mu^r$ are each 0, a natural module, or the dual of a natural module.  Indeed this follows from Lemmas \ref{twolabs} and \ref{tensorprodMF} unless $k = r = 2$ and $\mu^r = 2\l_1^r$ or its dual.  In this exceptional case it follows from Theorem \ref{LR} that ($\mu^1 \otimes \mu^r)\downarrow L_X' \supseteq (2110 \ldots )^2$ or $(10\ldots 02)^2$, respectively, and this is a contradiction.  This establishes the claim.  Also at most one of $\mu^0$ and $\mu^r$ is nonzero, as otherwise $V^1$ is not MF.

Lemma \ref{no1<i<r-1} implies that $V^2(Q_Y) = V_{\gamma_1}^2(Q_Y) + V_{\gamma_2}^2(Q_Y) + V_{\gamma_r}^2(Q_Y).$  Now
  $\lambda - \gamma_1 -\beta_1^1$ affords the $L_Y'$-module
  $$V_{C^0}(\mu^0 + \lambda_{r_0}^0) \otimes V_{C^1}(\lambda_2^1) \otimes V_{C^r}(\mu^r)$$
  in $V_{\gamma_1}^2(Q_Y)$  and $V_{C^1}(\lambda_2^1) \downarrow L_X'$ has  irreducible summands of highest weights $\nu_1 = 2\delta^1 - \alpha_k,$  $\nu_2 = 2\delta^1 - \alpha_{k-1}$, and $\nu_3 = 2\delta^1 - (\alpha_k+ \alpha_{k-1}) = \om_{k-2} + \om_{k-1}+(2r-3)\om_k+\om_{k+1}$  (omit $\om_{k-2}$ if $k = 2$). These
  direct summands satisfy $S(\nu_1), S(\nu_2) \le  2r$ and $S(\nu_3) \le 2r$ ($2r-1$ if $k= 2$).   In fact all irreducible summands of $V_{C^1}(\l_2^1) \downarrow L_X'$ have $S$-value at most $2r$.  This is because all composition factors
have the form $2\delta^1-\sum c_i\a_i.$ 

  Suppose $\mu^0 = \lambda_1^0.$  Then $\mu^r = 0$ and $V^1(Q_Y) = V_{C^0}(\mu^0) \otimes V_{C^1}(\mu^1).$   We have
  $$V_{C^0}(\mu^0 + \lambda_{r_0}^0) \downarrow L_X' = (V_{L_X'}(r\omega_k) \otimes V_{L_X'}(r\omega_{l-k+1})) - 0.$$ 
  The tensor product contains a direct summand of highest weight 
  $$(r\omega_k + r\omega_{l-k+1}) - (\alpha_k + \dots + \alpha_{l-k+1}) = \omega_{k-1} + (r-1)\omega_k + (r-1)\omega_{l-k+1} + \omega_{l-k+2}.$$
  The rank $2$ case of Lemma \ref{abtimescd}  applies to the tensor product of irreducibles of highest weights
    $\omega_{k-1} + (r-1)\omega_k$ and $\omega_{k-1}+(2r-3)\omega_k$ (which occur as restrictions to an $A_2$ Levi factor of $C^0$ of the module   $V_{C^0}(\mu^0 + \lambda_{r_0}^0) \downarrow L_X'$ and $\nu_3$, respectively. This  then implies that $V_{C^0}(\mu^0 + \lambda_{r_0}^0) \otimes V_{C^1}(\lambda_2^1)$
    has a composition factor of multiplicity at least $2$ whose highest weight has $S$-value at least 
    $4r-1$.   Lemma \ref{weightsum} implies that this  bounded
    above by $r+r  + 1$, a contradiction.
    
   Therefore, $\mu^0 = 0$  or $\lambda_{r_0}^0.$  Also, it follows from the above discussion that an upper bound for the 
$S$-value of highest weights of irreducible $L_X'$ composition factors  on
  $V_{\gamma_1}^2(Q_Y)$ is $4r$ or $3r$ according to whether or not one of $\mu^0$ or $ \mu^r$ is nonzero.  The largest possible value occurs
  if $\la\lambda, \gamma_1\ra \ne 0.$  We obtain a similar conclusion for  $V_{\gamma_r}^2(Q_Y).$ And if $r > 2$,
  the largest $S$-value appearing in  the restriction to $L_X'$ of $V_{\gamma_2}^2(Q_Y)$ is $2r$ or  $r$
  according to whether or not one of $\mu^0$ or $ \mu^r$ is nonzero.
  
   We see that
 $V^3(Q_Y)$ has a direct summand with highest weight $\lambda - \beta_{r_0}^0 -2\gamma_1 -2\beta_1^1-\beta_2^1$
 which affords the irreducible module 
  $$V_{C^0}(\mu^0 + \lambda_{r_0-1}^0) \otimes V_{C^1}(\lambda_3^1) \otimes V_{C^r}(\mu^r).$$
  Using the fact that $V_{C^0}(\lambda_{r_0}^0 + \lambda_{r_0-1}^0) = V_{C^0}(\lambda_{r_0}^0) \otimes V_{C^0} (\lambda_{r_0-1}^0) - V_{C^0}(\lambda_{r_0-2}^0)$, we see that the restriction of the first factor  to  $L_X'$  contains an irreducible module of highest weight   
   $$\omega_{l-k}+(2r-2)\omega_{l-k+1}+\omega_{l-k+2}$$
   or
   $$\omega_{l-k}+(3r-2)\omega_{l-k+1}+\omega_{l-k+2},$$
   according to whether $\mu^0 = 0$ or $\lambda_{r_0}^0.$
   And the second tensor factor has an $L_X'$ irreducible module of highest weight 
   $$3\omega_{k-1} + 3(r-1)\omega_k - (\alpha_{k-1}+\alpha_k).$$
   Finally, as above the last factor restricts to an irreducible for $L_X'$ of highest weight $0, r\omega_{k-1},$ or $r\omega_{l-k+2}.$
   Applying Lemma \ref{abtimescd} we find that  $V^3(Q_Y)$    has an irreducible $L_X'$ composition factor of multiplicity at least 2 and having highest weight $\eta$ such that
   either $S(\eta) \ge (2r) + (3r-1)  + S(\mu^r \downarrow L_X') $ or $S(\eta) \ge (3r) + (3r-1) $, respectively, noting
   that $\mu^r = 0$ in case $\mu^0 \ne 0.$ It follows from the above paragraph and Lemma \ref{weightsum} that $5r-1 \le 3r+1$ or  $6r-1 \le 4r+1$, respectively.  In either case this is a
    contradiction.   Therefore, $\mu^1 \ne \lambda_1^1$.
   
   Next consider the case $\mu^1 = \lambda_2^1$ and assume for now that $r > 2$.   This case is considerably easier since
   $\lambda - \gamma_1 -\beta_1^1 -\beta_2^1$ affords
   $$V_{C^0}(\mu^0 + \lambda_{r_0}^0) \otimes V_{C^1}(\lambda_3^1) \otimes V_{C^r}(\mu^r)$$
   and there is  a composition factor in $V_{C^1}(\lambda_3^1) \downarrow L_X'$ having multiplicity at least 2 and of highest weight $(3\omega_{k-1} + 3(r -1)\omega_k) - (\alpha_{k-1}+2\alpha_k)$   which we see from the equality $V_{C^1}(\lambda_3^1) = \wedge^3(V_{C^1}(\lambda_1^1)).$    Lemma \ref{p.29} shows that $S(V_{C^0}(\mu^0 + \lambda_{r_0}^0) \downarrow L_X')
   = S(V_{C^0}(\mu^0) \downarrow L_X') + r$.   Now an $S$-value argument gives a contradiction.
   
   Now assume  $\mu^1 = \lambda_2^1$ with $r = 2$. This is a slight variation from the above.  Viewing $V_{C^1}(\lambda_3^1)$ as wedge cube we obtain a composition factor of multiplicity $2$ and of highest weight
   $(3\omega_{k-1}+3\omega_k) - (2\alpha_{k-1}+2\alpha_k+\alpha_{k+1}) = 2\omega_{k-2} +\omega_{k-1} +  2\omega_k+\omega_{k+2}$   (the term $\omega_{k-2}$ does not appear if $k = 2$).   As the maximal $S$-value of $V_{C^1}(\lambda_2^1) \downarrow L_X'$ is $4$ we again obtain a contradiction by comparing  $S$-values.     
   
    Finally consider the  case $\mu^1 = 2\lambda_1^1$.  There is an irreducible summand of $V^2(Q_Y)$ afforded
    by $\lambda - \gamma_1 - \beta_1^1.$  The corresponding irreducible for $L_Y'$ is
   $$V_{C^0}(\mu^0 + \lambda_{r_0}^0) \otimes V_{C^1}(\lambda_1^1 + \lambda_2^1) \otimes V_{C^r}(\mu^r).$$
   
   We claim $V_{C^1}(\lambda_1^1 + \lambda_2^1) \downarrow L_X'$  contains  the irreducible module  of highest weight $3\omega_{k-1} + 3(r-1)\omega_k - \alpha_{k-1}-\alpha_k$  with multiplicity at least $2$.  To see this we first note that  $V_{C^1}(\lambda_1^1 + \lambda_2^1) = V_{C^1}(\lambda_1^1) \otimes V_{C^1}(\lambda_2^1)  - V_{C^1}(\lambda_3^1).$  Viewing $V_{C^1}(\lambda_3^1)$ as $\wedge^3(\lambda_1^1)$ it is clear that  $3\omega_{k-1} + 3(r-1)\omega_k - \alpha_{k-1}-\alpha_k$ occurs within  $V_{C^1}(\lambda_3^1) \downarrow L_X'$ with multiplicity $1$.  On  other hand $V_{C^1}(\lambda_2^1) \downarrow L_X'$ contains the irreducibles with highest weights
   $2\omega_{k-1} + 2(r-1)\omega_k - \alpha_k,$  $2\omega_{k-1} + 2(r-1)\omega_k - \alpha_{k-1}$, and $2\omega_{k-1} + 2(r-1)\omega_k - \alpha_{k-1}-\alpha_k.$  Tensoring each of these with $V_{C^1}(\lambda_1^1) \downarrow L_X'$ yields a summand of highest weight $3\omega_{k-1} + 3(r-1)\omega_k - \alpha_{k-1}-\alpha_k$, so this establishes the claim.
   
 The arguments of the above paragraph show $V_{C^1}(\lambda_1^1 + \lambda_2^1) \downarrow L_X'$ contains
   an irreducible summand with multiplicity $2$ and whose highest weight has $S$-value at least $3r-1.$ From Lemma \ref{p.29} we see that the $S$-value of $V_{C^0}(\mu^0 + \lambda_{r_0}^0) \downarrow L_X'$ is at least $S(V_{C^0}(\mu^0) \downarrow L_X') +r$.  We therefore obtain a contradiction as the  maximal $S$-value of $V^1$ is $S(V_{C^0}(\mu^0) \downarrow L_X') + 2r + S(V_{C^r}(\mu^r)) \downarrow L_X'.$      \hal
   
\begin{lem}\label{mu1=0}  $\mu^1 = \mu^{r-1} = 0$. 
\end{lem} 

\pf  Assume $\mu^1 \ne 0.$  Then Lemmas \ref{mu1notomega1etc} and \ref{induc} imply that $\mu^1 = \l_{r_1}^1, 2\l_{r_1}^1$, or $ \l_{r_1 - 1}^1.$ If $r = 2$, then $C^1 = C^{r-1}$ and we apply Lemma \ref{dimineq}.  If $\dim V_{L_X'}(r\om_k) = \dim V_{L_X'}(r\om_{k-1})$, then we obtain a contradiction by applying  Lemma \ref{mu1notomega1etc} to $V^*$.	Otherwise $V^*$ contradicts Lemma \ref{induc}.  Therefore assume $r > 2.$ 
 Then Lemma \ref{atmost1} implies that $\mu^i = 0$ for all $i \ne 0, 1, r.$ 
 
 Assume  $\mu^1 = \lambda_{r_1}^1.$  Then $V_{\gamma_2}^2(Q_Y)$ contains an irreducible summand $V_{C^0}(\mu^0) \otimes V_{C^1} (\lambda_{r_1-1}^1) \otimes V_{C^2}( \lambda_1^2) +V_{C^r}( \mu^r).$  The assumption $r > 2,$ Lemma \ref{lemma1.5(i)(ii)}, and  Lemma \ref{abtimescd} show that this tensor product has an  $L_X'$ composition factor of multiplicity at least 2 and for which the 
$S$-value of the highest weight is at least  $S(\mu^0 \downarrow L_X') + 2r  + r  -2 + S(\mu^r \downarrow L_X').$  On the other hand  the $S$-values of highest weights of summands of $V^1$ are at most  $S(\mu^0 \downarrow L_X') + r +S(\mu^r \downarrow L_X').$   Therefore
Lemma \ref{weightsum} implies that $3r-2 \le r+1,$ a contradiction.  

If  $\mu^1 = 2\lambda_{r_1}^1$, then $V_{\gamma_2}^2(Q_Y) $ contains $V_{C^0}(\mu^0) \otimes V_{C^1} (\lambda_{r_1-1}^1 + \lambda_{r_1}^1) \otimes V_{C^2}( \lambda_1^2) +V_{C^r}( \mu^r)$
and as in the previous result the restriction to $L_X'$ of the second tensor factor contains an irreducible summand of multiplicity 2 and with
S-value at least $3r-1.$  At this point we obtain a contradiction as
above. Now assume $\mu^1 =\lambda_{r_1-1}^1.$   Here $V_{\gamma_2}^2(Q_Y) $ contains $V_{C^0}(\mu^0) \otimes V_{C^1} (\lambda_{r_1-2}^1) \otimes V_{C^2}( \lambda_1^2) +V_{C^r}( \mu^r).$  As in the proof of Lemma \ref{mu1notomega1etc} the restriction of the  second tensor factor contains $3(r-1)\omega_{l-k+1} + 3\omega_{l-k+1} - 2\a_{l-k+1} -\a_{l-k+2}$ with multiplicity 2 and we again obtain a contradiction from an $S$-value consideration.
Therefore $\mu^1 = 0.$

Now consider $\mu^{r-1}.$  If $\mu^{r-1} \in \{\lambda_1^{r-1}, \lambda_2^{r-1}, 2\lambda_1^{r-1}\},$ then  arguments analgous to those above yield a  contradiction.  Therefore, assume 
$\mu^{r-1} \in \{\lambda_{r_{r-1}}^{r-1},  \lambda_{r_{r-1}-1}^{r-1}, 2\lambda_{r_{r-1}}^{r-1} \}$.  Here we consider the dual representation. Lemma \ref{dimineq} shows that either $\dim V_{L_X'}(r\omega_k) = \dim V_{L_X'}(r\omega_{k-1})$ or else
$\dim V_{L_X'}(r\omega_k) > \dim V_{L_X'}(r\omega_{k-1}) + 4.$  In the first case
 the labelling of $V^*$ contradicts the fact that we have just shown that $\mu^1 = 0.$  And in the second case we contradict the inductive hypothesis.   \hal

At this point we have
$$V^1(Q_Y)   = V_{C^0}(\mu^0)  \otimes V_{C^r}(\mu^r).$$

\begin{lem}\label{lambdadelta=0}  $\la\lambda, \gamma_i \ra  = 0$ for  $1 \le i \le r$.
\end{lem}

\pf  Assume  $\la \lambda, \gamma_i \ra \ne 0.$   Then by Lemma \ref{no1<i<r-1}  we have $i = 1$ or $i = r$, and by 
Lemmas \ref{no1<i<r-1} and \ref{mu1=0} we have  $V^2(Q_Y) = V_{\gamma_1}^2(Q_Y) + V_{\gamma_r}^2(Q_Y).$ 

We will require an upper bound for $S$-values of irreducible $L_X'$ modules appearing in $V_{\gamma_r}^2(Q_Y).$  
Checking the possibilities for $\mu^r$ and using Lemma \ref{p.29} which shows that $S(V_{C^r}(\l_1^r+\mu^r) \downarrow L_X') =  S(V_{C^r}(\mu^r)\downarrow L_X')+ r $, we see  that a maximal $S$-value occurs when $\la\lambda, \gamma_r \ra \ne 0$ where the irreducible 
afforded by $\lambda - \gamma_r$ has the largest $S$-value.  This irreducible is
$$V_{C^0}(\mu^0)  \otimes V_{C^{r-1}}(\lambda_{r_{r-1}}^{r-1}) \otimes  
V_{C^r}(\lambda_1^r +\mu^r),$$
where (in view of Lemma \ref{p.29}(ii)) the $S$-value is  $S(V_{C^0}(\mu^0) \downarrow L_X') + r + r + S(V_{C^r}(\mu^r) \downarrow L_X').$

Suppose $i = 1$.  Here $\lambda - \gamma_1$
affords 
$$ V_{C^0}(\mu^0 + \lambda_{r_0}^0) \otimes V_{C^1}(\lambda_1^1)  \otimes V_{C^r}(\mu^r)$$
and the $S$-value of the highest weight upon restriction to $S_X$ is $S(\mu^0 \downarrow L_X') + S(\mu^r \downarrow L_X') + 2r$.  

First suppose $\mu^0 = 0.$  Here we pass to $V_{\g_1}^3$, where we obtain
 $$V_{C^0}(\lambda_{r_0-1}^0) \otimes V_{C^1}(\lambda_2^1)  \otimes V_{C^r}(\mu^r).$$
Restricting to $L_X'$ we have  a summand of the form $V_{L_X'}(\omega_{k-1} + (2r-2)\omega_k+\omega_{k+1})^* \otimes \wedge^2V_{L_X'}(\omega_{k-1} + (r-1)\omega_k) \otimes (V_{C^r}(\mu^r) \downarrow L_X').$ The second tensor factor contains a summand of highest weight $(2\omega_{k-1} + 2(r-1)\omega_k) -(\alpha_{k-1}+\alpha_k).$
Then  Lemma \ref{abtimescd}  yields a composition
factor of multiplicity at least 2 and having  $S$-value at least $(2r) + (2r) - 2 + S(\mu^r \downarrow L_X') = 4r-2 +S(\mu^r \downarrow L_X').$  But the largest $S$-value of an $L_X'$ composition factor of the summand of $V^2(Q_Y)$ is $2r + S(\mu^r \downarrow S_X),$ so this is a contradiction.  Therefore, $\mu^0 \ne 0.$

In the remaining cases $\mu^0 = \lambda_1^0, \lambda_2^0, 2\lambda_1^0, \lambda_1^0+\lambda_{r_0}^0, \lambda_3^0$ (here $k = r = 2$) or the dual of one of these.  We claim that in each case
there is a composition factor in $V_{\gamma_1}^2(Q_Y)$ having multiplicity at least $2$  and having $S$-value at least $S(\mu^0 \downarrow L_X') + 2r - 2 + S(\mu^r \downarrow L_X')$.   This will yield a contradiction.
Towards this end  we note that Lemma \ref{v2gamma1} implies
\begin{equation}\label{tensor}
V_{\gamma_1}^2(Q_Y) \supseteq V_{C^0}(\mu^0) \otimes V_{C^0}(\l_{r_0}^0) \otimes V_{C^1}(\l_1^1)  \otimes V_{C^r}(\mu^r).
\end{equation} 
We will establish the claim using the above tensor product together with Lemmas \ref{abtimesce}, \ref{abtimescd} and \ref{compfactor}.

If $\mu^0 = c\lambda_1^0$ for c = 1,2, then Lemma \ref{abtimesce} shows that the restriction of $V_{C^0}(\mu^0) \otimes V_{C^0}(\l_{r_0}^0)$
contains $\om_{k-1} + (cr-1)\om_k + (r-1)\om_{l-k+1}+ \om_{l-k+2}$.  Tensoring with the other two factors and applying  Lemma \ref{abtimescd}, we have  a composition factor of multiplicity 2 and S-value at least $S(\mu^0\downarrow L_X') + 2r -1+ S(\mu^r\downarrow L_X')$.  Now assume
$\mu^0 = c\lambda_{r_0}^0.$  The restriction of the tensor product of this with $V_{C^0}(\l_{r_0}^0)$ contains $\om_{l-k}+((c+1)r-2)\om_{l-k+1}+ \om_{l-k+2}$.  Applying Lemma \ref{abtimescd} we obtain a composition factor of multiplicity 2
and S-value  $(c+2)r-1 + S(\mu^r\downarrow L_X') = S(\mu^0\downarrow L_X') +2r-1 + S(\mu^r\downarrow L_X').$

The cases $\mu^0 = \lambda_2^0$ and its dual are very similar to those above.  Lemma \ref{wedgefactors} shows that the restriction of 
$ \lambda_2^0$ to $L_X'$ contains $\om_{k-1} + (2r-2)\om_k + \om_{k+1}$.  Lemma \ref{compfactor} shows that tensoring this with the restriction of $V_{C^1}(\l_1^1)$ produces a summand with multiplicity 2 and  highest weight $(2\om_{k-1} + (3r-3)\om_k + \om_{k+1}) - \alpha_{k-1}-\alpha_k$.  Tensoring with the remaining factors and considering S-values gives the claim.  And if $\mu^0 = \lambda_{r_0 -1}^0$, then applying Lemma \ref{abtimescd} to the restriction of $V_{C^0}(\mu^0)  \otimes V_{C^1}(\l_1^1)$ yields a composition factor with multiplicity at least 2 having highest weight at least $3r-1$.   Tensoring with the remaining factors yields the claim.

Next assume $\mu^0  = \lambda_1^0 + \lambda_{r_0}^0.$  The restriction of the first tensor factor of (\ref{tensor}) to $L_X'$
contains $r\om_k + r\om_{l-k+1}$ so that tensoring with the restriction of $V_{C^0}(\l_{r_0}^0)$ and applying  Lemma \ref{abtimesce} we obtain a summand with highest weight  $r\omega_k + 2r\omega_{l-k+1} -(\alpha_k + \dots + \alpha_{l-k+1}) = \omega_{k-1} + (r-1)\omega_k + (2r-1)\omega_{l-k+1} + \omega_{l-k+2}.$ 
Again Lemma \ref{abtimescd} yields the claim.

Now assume $\mu^0  = \lambda_3^0$ with $k = r = 2.$ If $l > 3$, then  $\mu^0 \downarrow L_X' \supseteq (2\om_1 + 3\om_2 + \om_4)$  and if $l = 3$ then $\mu^0 \downarrow L_X' \supseteq (3\om_1 + 3\om_3).$  In each case these summands have maximal $S$-value. Tensoring with $V_{C^1}(\l_1^1) \downarrow L_X' = \om_1 + \om_2$ and using Lemma \ref{compfactor}  gives a composition factor with multiplicity 2 and highest weight $2\om_1 + 3\om_2 + \om_3+ \om_4$ or $3\om_1 + \om_2 + 2\om_3$, respectively.   Then tensoring with $V_{C^0}(\l_{r_0}^0)  \otimes V_{C^r}(\mu^r)$ yields a composition factor with multiplicity 2 and S-value at least $7+2 + S(\mu^r \downarrow L_X')$ or  $6+2 + S(\mu^r \downarrow L_X'),$ respectively.   In either case the claim holds.

The final case is $\mu^0  = \lambda_{r_0-2}^0$ again with $ r = k  = 2.$  We first settle the case $l = 3.$ Then the restriction of 
$V_{C^0}(\mu^0) \otimes V_{C^0}(\l_{r_0}^0) \otimes V_{C^1}(\l_1^1)$ contains $(303) \otimes (020) \otimes (110)$ and this contains $(324)^2$.  Tensoring with $V_{C^r}(\mu^r)$ and comparing S-values yields the claim
here.  So now assume $l \ge 4.$  Here the restriction of $V_{C^0}(\mu^0) \otimes V_{C^0}(\l_{r_0}^0)$ contains $\om_{l-3} + \om_{l-2} + 3\om_{l-1} + 3\om_l$.  Now tensor with $V_{C^1}(\l_1^1) \downarrow L_X' = \om_1 + \om_2$ and apply Lemma \ref{abtimescd} to establish the claim.  Therefore we have a contradiction when $i = 1.$

Finally assume $i = r.$  Here we again use Lemma \ref{dimineq}. If Lemma \ref{dimineq}(i)  holds, then consideration of  $V^*$ reduces us to the case just considered.  And if Lemma \ref{dimineq}(ii) holds, then $\dim V_{L_X'}(r\omega_k) > \dim V_{L_X'}(r\omega_{k-1}) + 4$ and we contradict Lemma \ref{induc} in $V^*$. \hal

\begin{lem}\label{mu0mur} The possibilities for $\mu^0$ and $\mu^r$ are as follows:
\begin{itemize}
\item[{\rm (i)}] $\mu^0 \in \{ \lambda_1^0, \,\lambda_2^0, \,\lambda_3^0, \,2\lambda_1^0\}$;
\item[{\rm (ii)}]   $\mu^r \in \{ \lambda_{r_r}^0,\, \lambda_{r_r-1}^0,\, \lambda_{r_r-2}^0,\, 2\lambda_{r_r}^0\}.$
\end{itemize}
\end{lem}

\pf By taking duals it suffices to prove either (i) or (ii).  In view of previous lemmas we cannot have $\mu^0 = \mu^r = 0,$ for otherwise
$\lambda = 0.$  And by duality we may assume $\mu^0 \ne 0.$ In view of Lemma \ref{induc}, seeking a contradiction, we may assume $\mu^0 \in \{ \lambda_{r_0}^0, \lambda_{r_0-1}^0, \lambda_1^0 + \lambda_{r_0}^0, \lambda_{r_0-2}^0,  2\lambda_{r_0}^0 \}.$  In each case we obtain a contradiction from consideration of 
$V_{\gamma_1}^2(Q_Y).$   

First assume $\mu^0  =\lambda_{r_0}^0$.  Here $V_{\gamma_1}^2(Q_Y) = V_{C^0}(\lambda_{r_0-1}^0) \otimes  V_{C^1}(\lambda_1^1) \otimes V_{C^r}(\mu^r).$  Using Lemmas \ref{wedgefactors} and  
Lemma \ref{abtimescd}  we see that the tensor product of the first two terms contains  an irreducible $L_X'$ summand appearing with multiplicity at least $2$
and having $S$-value  at least $(2r) + r -1= 3r-1$.  On the other hand
 $V^1(Q_Y) = V_{C^0}( \lambda_{r_0}^0) \otimes V_{C^r}(\mu^r)$ which has $S$-value $r + S(\mu^r \downarrow L_X')$, giving a  contradiction.

Next consider $\mu^0  = \lambda_1^0 + \lambda_{r_0}^0.$  Here $V_{\gamma_1}^2(Q_Y)$ contains 
 $(V_{C^0}(\lambda_1^0 +\lambda_{r_0-1}^0) + V_{C^0}(\lambda_{r_0}^0)) \otimes  V_{C^1}(\lambda_1^1) \otimes V_{C^r}(\mu^r)$  which equals   $V_{C^0}(\lambda_1^0) \otimes V_{C^0}(\lambda_{r_0-1}^0)   \otimes  V_{C^1}(\lambda_1^1) \otimes V_{C^r}(\mu^r).$  We obtain a contradiction by applying Lemma \ref{abtimescd} to the restrictions of the tensor product of the second and third factors and then tensoring with the remaining factors. 
 
If $\mu^0 = 2\l_{r_0}^0$, then $V_{\gamma_1}^2(Q_Y)$ contains  $V_{C^0}(\lambda_{r_0-1}^0+ 
\lambda_{r_0}^0) \otimes  V_{C^1}(\lambda_1^1) \otimes V_{C^r}(\mu^r)$  and the restriction to $L_X'$ of the first tensor factor has a composition factor of highest weight $\om_{l-k} + (3r-2)\om_{l-k+1} + \om_{l-k+2}.$
Tensoring this with the restriction of the second factor and applying Lemma \ref{abtimescd}  we obtain a composition factor of multiplicity 2 and S-value $4r-2$. Now tensoring with the third factor yields a contradiction.
 
 If  $\mu^0 = \lambda_{r_0-1}^0$ we obtain a contradiction from Lemma \ref{abtimescd}  using the fact that the restriction of $\lambda_{r_0-2}^0$ contains $3\om_{l-k} + (3r-6)\om_{l-k+1}+3 \om_{l-k+2}.$  Finally assume $\mu^0 = \lambda_{r_0-2}^0$, so that $r=k=2.$  Then  $V_{C^0}(\lambda_{r_0-3}^0) \downarrow L_X'$  contains $\om_{l-k-1} + \om_{l-k} + 4\om_{l-k+1} + \om_{l-k+2}$ or $(214)$, according as $l \ge 4$ or $l = 3.$  At this point we tensor with $V_{C^1}(\lambda_1^1)$, apply  Lemma \ref{abtimesce}, and obtain the usual contradiction.
\hal

\begin{lem}\label{mu0murnonzero} If both $\mu^0$ and $\mu^r$ are nonzero, then
   $\lambda=  \lambda_1+\lambda_n.$
\end{lem}

\pf  Suppose $\mu^0 \ne 0 \ne \mu^r$.  As $V^1$ is MF, it is immediate
from  Proposition \ref{stem1.1.A}, Lemma \ref{wedgefactors} and Lemma \ref{mu0mur}  
that either $\mu^0 = \l_1^0$ or $\mu^r = \l_{r_r}^r$ and taking duals, if necessary, we may assume that $\mu^0 = \l_1^0$.  Assume $\l \ne \l_1 + \l_n$. By  Lemma \ref{lambdadelta=0} we have  $\la\lambda, \gamma_i \ra  = 0$ for  $1 \le i \le r$. 

Then Lemma \ref{mu0mur} implies that $\mu^r \in \{  \lambda_{r_r-1}^r, \lambda_{r_r-2}^r, 2\lambda_{r_r}^r\}.$  
 Proposition \ref{stem1.1.A} and Lemma \ref{wedgefactors} yield a contradiction unless $r = 2$ and either $\mu^r   =\lambda_{r_r-1}^0$ with $k = 2$ or $\mu^r   = 2\lambda_{r_r}^0.$   

First assume we have the latter case and consider the dual module, $V^*$.  If $k > 2$ and $l > 4$, then using Magma we see that $S^2(2\omega_k)$ has a composition of highest weight
$4\omega_k - (\alpha_{k-1} + 2\alpha_k + \alpha_{k+1}) = \omega_{k-2} + 2\omega_k + \omega_{k+2}$ and the tensor product
of this with $(\mu^*)^r = \lambda_{r_r}^r$ is not MF by  Proposition \ref{stem1.1.A}.  Next suppose $k>2$ and $l = 4$.
Here $k = 3$ and $(V^*)^1 = S^2(0020) \otimes (0020)$.  A Magma check shows that  there is a composition factor of highest weight $0222$ appearing with multiplicity at least 2.

Finally, consider $k = 2.$  If $l = 3$,  then $(V^*)^1 = S^2(020) \otimes (002)$ and a Magma check shows that this contains $(220)^2$,   a contradiction.
Now assume $l > 3$.   Here $(V^*)^1 = S^2(2\omega_2) \otimes (2\omega_l)$.
The first tensor factor  has  summands of highest weights   $4\omega_2$ and $2\omega_1+2\omega_3$.  Using Lemma \ref{aibj12} we see that each of these yields an irreducible module
of highest weight $2\omega_1+2\omega_2$ upon tensoring with $2\om_l$, so this is again a contradiction.

Now assume that $\mu^r = \lambda_{r_r-1}^r$ with $r = 2 = k.$  Here we again consider $V^*$ so that
$(V^*)^1 = \wedge^2(2\omega_2) \otimes (2\omega_l).$  The first factor has summands
of highest weights $\omega_1 + 2\omega_2 +\omega_3$ and $\omega_1 + \omega_3 +\omega_4$  (the  term $\omega_4$ does not occur if $l = 3$). Tensoring with the second factor we will argue that each of these yields a summand  with highest weight $\omega_1 + \omega_2 +\omega_3.$ This is based on  a hom argument together with a Magma computation.
For example,  for the first containment we have 
\[
\begin{array}{ll}
{\rm Hom}\left((1110 \ldots  0), \,(1210\ldots  0) \otimes (0 \ldots  02)\right) & \cong 
{\rm Hom}\left(0, \,(0\ldots  0111) \otimes (1210\ldots  0) \otimes (0 \ldots  02)\right) \\
  & \cong {\rm Hom}\left((0 \ldots  121),\, (0\ldots  0111) \otimes (0 \ldots  02)\right).
\end{array}
\]
We see that the last term is nonzero since $(0\ldots  0111) \otimes (0 \ldots  02)$ contains $(0\ldots  0113)-\a_l.$ Similarly for the second case.  So this contradics the fact that $(V^*)^1$ is MF. \hal

We can now complete the proof of  Theorem \ref{main} for $l > 2$. We must consider the cases where just one of $\mu^0$ and $\mu^r$ is nonzero.  Replacing $V$ by $V^*$ if necessary
we may suppose $\mu^0 \ne 0.$  Then Lemma \ref{mu0mur}  and the above results imply
$\lambda \in \{\lambda_1, \lambda_2, \lambda_3, 2\lambda_1\}$ and if $\lambda = \lambda_3$, then $r = 2.$  
So at this point we have completed the proof of Theorem \ref{main} in the case where $l>2$.


\section{Case $l=2$. }  We now prove Theorem \ref{main} in the case where $l = 2$. Here there are complications due to the fact that there are
 additional examples to consider in the induction hypothesis, especially when $r = 2.$
Recall that $W = V_X(\d)$ where $\d = 0r0 = r\omega_2$  and recall that  $Y = SL(W) = A_n$.

In this case the induction hypothesis in the statement of Theorem \ref{main} gives the following.

\begin{lem}\label{inducl2} Assume that $l=2$.
\begin{itemize}
\item[{\rm (i)}] If $1 < i < r-1,$ then $\mu^i = 0, \,\lambda_1^i$ or $\lambda_{r_i}^i.$
\item[{\rm (ii)}] If  $r>2$, then  $\mu^1 $ or $(\mu^*)^1$  is in  $\{0,\, \lambda_1^1, \, \lambda_2^1, \,2\lambda_1^1\}$.
\item[{\rm (iii)}] If $r=2$, then $\mu^1 $ or $(\mu^*)^1$ is in  $\{ 0,\, \lambda_1^1, \, \lambda_2^1,\,  \lambda_3^1,\, 2\lambda_1^1,\, 3\lambda_1^1\}$. 
\end{itemize}
\no Also for $i = 0$ or $r$, either
\begin{itemize}
\item[{\rm (iv)}] $\mu^i$ or $(\mu^i)^*$ is in  $\{ 0, \, j\l_1^i,  \, \l_j^i, \, \l_1^i + \l_{r_i}^i, \,
 \l_1^i + \l_2^i (r = 3)\}$, or
\item[{\rm (v)}] $r = 2$ and $\mu^i$ or $(\mu^i)^*$ is either in (iv) or is in
\[
\{ 0, \, a\l_1^i + \l_2^i (a \le 3), \,  \l_1^i + \l_3^i,  \, \l_2^i + \l_3^i,\, \l_1^i + \l_4^i, \,
\l_2^i + \l_4^i, \, \l_1^i + 2\l_5^i, \,  \l_1^i + 3\l_5^i, \, 2\l_2^i, 3\l_2^i \}.
\]
\end{itemize}
\end{lem} 
  
  Note that in (iv) above the values of $j$ depend on $r.$  Indeed, $j\lambda_1^0$ is allowed
  for any $j$ if $r = 2$;  $j = 4$ requires $r \le 3$; $j = 3$ requires $r \le 5;$ and $j = 2$ holds for all $r$. Similarly,
   $\l_j^0$ is allowed for any $j$ if $r= 2$; $j = 5$ requires $r \le 3$; $j = 4$ requires $r \le 4$; $j =3$ requires $r \le 6$; and $j = 2$ holds for all $r$.

We begin with some initial observations.  Firstly,  $W \downarrow L_X'$ is self-dual, so there is a duality
among the factors $(C^0, C^r), (C^1, C^{r-1}), \dots$ which will be exploited.  Lemmas \ref{atmost1} and \ref{no1<i<r-1} hold in this situation, noting for $r = 2$ they hold trivially, and the proofs are valid if $r > 2.$ 
We next rule out a special case of (iv) above when $r = 3.$

\begin{lem}\label{11000} Assume that $r = 3$ and $ i = 0$ or $r$.
If $\mu^i$ or $(\mu^i)^*$ is $\l_1^i + \l_2^i,$ then $V\downarrow X$ is not MF.
\end{lem}

\pf Assume false.  Then taking duals we can assume that $\mu^0 = (110 \ldots  0)$ or $(0\ldots  011).$  Using Magma we find that  $\mu^0\downarrow L_X'= 25+14+33+11+22+41+17$, respectively the  dual. Then by  Lemmas \ref{twolabs}, \ref{tensorprodMF}, and \ref{pieri} we find that $\mu^i = 0 $ for $i \ne 0.$

Assume that $\mu^0 = (110 \ldots  0).$ Then $V_{\g_1}^2 \supseteq ((20\ldots  0)+(010\ldots 0))\downarrow L_X') \otimes (12)$.  The first tensor factor is $((10\ldots  0)\otimes (10\ldots  0)) \downarrow L_X'$, so $V_{\g_1}^2 \supseteq (03)\otimes (03) \otimes (12)$ and this contains $(23)^6.$ Only three such factors can arise from $V^1,$ so this is a contradiction.

Now assume $\mu^0 = (0\ldots  011).$ Here  $V_{\g_1}^2 \supseteq ((0\ldots  020)+(0\ldots  0101))\downarrow L_X') \otimes (12).$   This time the first tensor factor is $((0\ldots  010)\otimes (0\ldots  010)) - (0 \ldots  01000))\downarrow L_X' $ and we find that $V_{\g_1}^2 \supseteq (33)^3 \otimes (12) \supseteq  (34)^6$ and this is again a contradiction.   \hal

\begin{lem}  Theorem $\ref{main}$ holds if $r > 2.$
\end{lem}

\pf Assume $r > 2.$  The proof  amounts to noting that most of the proof of  Theorem \ref{main} for $l>2$ goes through
in this case, often with simplified proofs, although we must take into account extra possibilities which appear in Lemma \ref{inducl2}(iv).  We will discuss the changes required to establish the preceding Lemmas \ref{mu1notomega1etc} -   \ref{mu0murnonzero}.

 \ref{mu1notomega1etc}:   Here we must show $\mu^1 \ne \lambda_1^1, \lambda_2^1,$ or $2\lambda_1^1.$ Assume $\mu^1 = \lambda_1^1.$ As before $\mu^0 = 0, \lambda_1^0$, or $\lambda_{r_0}^0.$ Similarly for $\mu^r,$ but we cannot have both $\mu^0$ and $\mu^r$ nonzero.  With $\nu_1$, $\nu_2$ and $\nu_3$ as before we have $S(\nu_1), S(\nu_2), S(\nu_3)$ are $2r-1, 2r-1, 2r-2,$ respectively, and  these are the irreducible summands of the restriction  $V_{C^1}(\l_2^1) \downarrow L_X'$ with the largest S-value. Assume $\mu^0 = \lambda_1^0.$  Then $\mu^r = 0$ and the irreducible summand of $V_{\gamma_1}^2(Q_Y)$ afforded by
$\lambda - \gamma_1 -\beta_1^1$ contains 
$V_{C^0}(\mu^0 + \lambda_{r_0}^0) \otimes V_{C^1}(\l_2^1)$ and the restriction to  $L_X'$ contains the irreducible of highest weight $(rr) \otimes (3,2r-4).$ Applying Lemma \ref{abtimescd} we obtain a contradiction as before, by considering $S$-values.  So
 $\mu^0 = 0$ or  $\lambda_{r_0}^0$.  For the latter case the restriction to $L_X'$ of  the irreducible module afforded by $\lambda - \beta_{r_0}^0 -\gamma_1 $  contains $(2r-2,1) \otimes (2,2r-2) \supseteq (2r-1,2r-2)^2$ by Lemma \ref{abtimescd}. From here $S$-values again give  
a contradiction at level 2.  Suppose $\mu^0 = 0.$  Here we obtain a contradiction from the summand of  $V^3(Q_Y)$ afforded by $\l - \b_{r_0}^0-2\gamma_1-2\b_1^1-\b_2^1$ which  affords a composition factor $\rho$ which is the tensor product of $\wedge^2(r0) \otimes \wedge^3(1(r-1))$ and the restriction of $V_{C^r}(\mu^r)$.   The largest $S$-value appearing among $L_X'$ irreducibles  of $V^2(Q_Y)$ is $3r$ or $4r$, according to whether or not $\mu^r = 0,$   while the $L_Y'$-irreducible afforded by $\rho$   contains
in its restriction to $L_X'$ the module $(2r-2,1) \otimes (2,3r-4) \otimes V_{C^r}(\mu^r) \downarrow L_X'$.  The usual arguments yield a contradiction if $r > 3.$  If $r = 3$, then a Magma computation shows that $V^3$ contains $(4,7)^2,$ which is again a contradiction.  The remaining cases of the lemma are $\mu^1 = \lambda_2^1$ and $2\lambda_1^1$ and these proceed just as in the lemma.  

\ref{mu1=0}:  We must show $\mu^1 = \mu^{r-1} = 0,$ and by duality it will suffice to show $\mu^1  = 0.$ Proceed as in the lemma. If $\mu^1 = \lambda_{r_1}^1,$ then  $V_{\gamma_1}^2(Q_Y)$ has an irreducible summand of the form   $V_{C^0}(\mu^0) \otimes V_{C^1}(\lambda_{r_1-1}^1) \otimes V_{C^2}(\lambda_1^2) \otimes V_{C^r}(\mu^r).$   Lemma \ref{abtimescd} shows that there is an irreducible $L_X'$ composition factor of multiplicity at least 2 and for which the $S$-value of the highest weight is at least  $S(\mu^0 \downarrow S_X) + (2r -1) + r - 2 + S(\mu^r \downarrow S_X)$ and this yields a contradiction. The same argument gives a contradiction if $\mu^1 = 2\lambda_{r_1}^1$ or $\lambda_{r_1-1}^1.$   Therefore $\mu^1 = 0$ and we get $\mu^{r-1} = 0$ by duality.

\ref{lambdadelta=0}: We must show that $\la\lambda, \gamma_i \ra = 0$ for $1 \le i \le r$.
By Lemma \ref{no1<i<r-1} we need only consider $\gamma_1$ and $\gamma_r $ and by duality we may will work with $\gamma_1.$  By way of contradiction assume  $\la\lambda, \gamma_1 \ra \ne 0.$

  Suppose $\mu^0 = 0$. Recall that $\mu^1 = 0.$ Then $\l -\b_{r_0}^0-2\gamma_1-\b_1^1$ affords the summand $V_{C^0}(\lambda_{r_0-1}^0) \otimes V_{C^1}(\lambda_2^1)  \otimes V_{C^r}(\mu^r)$ of $V^3(Q_Y).$ The restriction of the tensor product of the first two terms contains $((2r-2)1)\otimes (3(2r-4))$ which contains $(2r(2r-4))^2.$  An S-value argument gives a contradiction.
  So now suppose $\mu^0 \ne 0.$

 Here we obtain a contradiction in $V^2(Q_Y)$, but  we must consider all possibilities in (iv) of Lemma \ref{inducl2}. Now $\lambda-\gamma_1$ affords 
$$V_{C^0}(\mu^0 + \lambda_{r_0}^0) \otimes V_{C^1}(\lambda_1^1) \otimes V_{C^r}(\mu^r).$$
We claim that in each case there is a summand of $V_{C^0}(\mu^0 + \lambda_{r_0}^0) \downarrow L_X'$
with highest weight of the form $a\omega_1 +b\omega_2$ such that $ab\ne 0$ and $a+b \ge S(\mu^0 \downarrow L_X')  - 2 + r.$
This is easily checked if $\mu^0$ or $(\mu^*)^0$ is in  $\{ \l_1^0, \, 2\l_1^i,  \, \l_2^0, \, \l_1^0 + \l_{r_0}^0 \}.$ In the remaining cases $r$ is bounded and a Magma computation gives the assertion.

Then Lemma \ref{abtimescd} yields a composition factor of multiplicity 2 with $S$-value at least $(S(\mu^0 \downarrow L_X') -2 + r ) + r-2 + S(\mu^r \downarrow L_X')$ and we obtain the usual numerical contradiction.

\ref{mu0mur}:   We begin just as in the lemma. We can assume $\mu^0 \ne 0$ and the goal is
to show that $\mu^0 \in \{ \lambda_1^0, \,\lambda_2^0, \,\lambda_3^0, \,2\lambda_1^0\}. $ Assume this does not hold. 
We first claim that if $\mu^0$ or  $(\mu^*)^0$ is one of the exceptional cases  $\l_j$ for $3 \le j \le 5$ or $j\l_1^0$  for 
$j = 3,4$, then $\mu^r = 0$.  Indeed Lemmas \ref{twolabs} and \ref{tensorprodMF} imply that $\mu^r$ or its dual is in  
$\{0, \, \l_1^r \}$ and then a Magma computation shows
that $\mu^r = 0$, establishing the claim.   So if $\mu^0$ is one of the exceptional cases, then $\l = j\l_1$ or $\l_j$ and
a Magma computation shows that $V\downarrow X$ is not MF.

At this point we proceed as in the proof of Lemma \ref{mu0mur} where it is first necessary to rule out cases $\mu^0 = \lambda_1^0 + \lambda_{r_0}^0,$ $\lambda_{r_0}^0$,  $\lambda_{r_0-1}^0$, $\lambda_{r_0-2}^0$, and $2\lambda_{r_0}^0$.  In addition  we must rule out the extra cases $\mu^0 = j\lambda_{r_0}^0 (j \ge 3)$ and $\lambda_{r_0-j+1}^0 (3\le j \le 5).$ 

The case $\mu^0 = \lambda_1^0 + \lambda_{r_0}^0$ is handled just as in Lemma \ref{mu0mur}.
If $\mu^0 = \lambda_{r_0}^0$, then $V_{\gamma_1}^2(Q_Y)$ contains $V_{C^0}(\lambda_{r_0-1}^0) \otimes V_{C^1}(\l_1^1) \otimes V_{C^r}(\mu^r).$   The restriction to $L_X'$ of the tensor product of the first two terms
contains $((2r-2)1) \otimes (1(r-1)$ and this contains $((2r-2)(r-1))^2$.  An S-value argument gives a contradiction.  If $\mu^0 = 2\lambda_{r_0}^0,$ then $V_{\gamma_1}^2(Q_Y)$ contains  $(V_{C^0}(\lambda_{r_0-1}^0+ \lambda_{r_0}) \otimes  V_{C^1}(\lambda_1^1) \otimes V_{C^r}(\mu^r).$  The tensor product
of the restriction of the first two terms contains $((3r-2)(r-1))^2$ and we again have a contradiction.  Suppose 
$\mu^0 = \lambda_{r_0-1}^0.$ We have $\wedge^3(r0) \supseteq ((3r-6)3),$
so that the tensor product with  $(1(r-1))$ contains $((3r-6)(r+1))^2.$  We then get an S-value contradiction noting that the S-value of $\wedge^2(r0)$ is $2r-1.$  Similarly if $\mu^0 = \lambda_{r_0-2}^0,$ then
$\wedge^4(r0)$ contains $((4r-7)2)$ and the tensor product with $(1(r-1)$ contains $((4r-7),r)^2.$  We then obtain a contradiction noting that the S-value of $\wedge^3(r0)$ is $3r-3.$  

It remains to consider the exceptional cases $ j\l_{r_0}^0$ and $\l_{r_0-j+1}^0$, where we have shown
that $\mu^r = 0.$

Now assume $\mu^0 = j\lambda_{r_0}^0$ for $j \ge 3,$ which only holds for a few values of $r$.  Then $V_{\gamma_1}^2(Q_Y)$ contains an irreducible summand afforded by $\lambda - \beta_{r_0}^0 - \gamma_1$ and this affords $V_{C^0}(\lambda_{r_0 -1}^0 + (j-1)\lambda_{r_0}^0) \otimes V_{C^1}(\lambda_1^1).$  The restriction to $L_X'$ of the first tensor factor contains an irreducible summand of high
weight $((j-1)r),0)+ ((2r-2)1) = (((j+1)r-2),1).$  The tensor product of this with the restriction of $V_{C^1}(\lambda_1^1)$ contains $(((j+1)r-2)(r-1))^2$ which has S-value $(j+2)r-3.$  The S-value of $V^1$ is $jr$, so this is a contradiction.

 Now consider the cases where $\mu^0 = \lambda_{r_0-j+1}^0$ for $3\le j \le 5.$ If $r = 3$ a direct check with
Magma shows that the relevant wedge powers of $030$ are not MF.  Otherwise $r = 4,5,6$ and we argue
as usual to get a contradiction in  $V_{\gamma_1}^2(Q_Y).$  Indeed this summand contains the irreducible module $V_{C^0}(\lambda_{r_0 -j}^0)  \otimes V_{C^1}(\lambda_1^1) \otimes V_{C^r}(\mu_r).$  Using Magma we find an irreducible of the first term of form $(ab)$ where $ab \ne 0$
and $a+b > S(V_{C^0}(\mu^0) \downarrow L_X').$  Then Lemma \ref{abtimescd}  guarantees that there
is an irreducible summand of  $V^2$ with $S$-value at least $(a+b) + r -2,$  which gives
a contradiction.

\ref{mu0murnonzero}:  Suppose $\mu^0 \ne 0 \ne \mu^r$, but $\lambda \ne \lambda_1 + \lambda_n.$  If  $\mu^0 \ne \lambda_1^0$ and $\mu^r \ne \lambda_{r_r}^r$,   
then  Lemmas \ref{inducl2}(iv), \ref{twolabs} and Proposition \ref{tensorprodMF} contradict the fact that
$V^1$ is MF.  So, using the dual if necessary, we may assume $\mu^0 = \l_1^0$ but $\mu^r \ne \l_r^r$. By 
(\ref{mu0mur}) above we can assume $\mu_r \in \{ 2\l_{r_r}^r, \, \l_{r_r-1}^r, \, \l_{r_r-2}^r \}$.  In each case it is easy to argue
that $V^1$ is not MF.

We can now complete the proof of the lemma.  We consider the case where just one of $\mu^0$ and $\mu^r$ is nonzero which we may take to be $\mu^0$.  Then we have reduced to $\lambda = \lambda_1$, $\lambda_2$, or $2\lambda_1$, (the case $\lambda_3$ is out since $r > 2$.)   \hal


We are left with the case where $r=2$, so that $X = A_3$ embedded in $A_{19}$ via the representaton 
$\d = 2\omega_2.$  We proceed with a series of lemmas.  Here $C^0 = C^2 = A_5,$ while $C^1 = A_7.$  As usual $\gamma_1, \gamma_2$ are the fundamental nodes adjacent to $C^0, C^2$, respectively.

In the proofs to follow we often use Magma for computations.
In particular,  Lemma \ref{A2info}  records additional information regarding the cases in Lemma \ref{inducl2}(v) and we often use this without reference.

\begin{lem}\label{mu1is0} We have $\mu^1 = 0.$
\end{lem}

\pf  Suppose $\mu^1 \ne 0.$  Then from Lemma \ref{inducl2}(iii)  and taking duals if necessary we can assume  $\mu^1 \in  \{\lambda_1^1,  \lambda_2^1, \lambda_3^1,2\lambda_1^1, 3\lambda_1^1\}$  and we check using Magma that $V_{C^1}(\mu^1) \downarrow L_X'$ has
a summand from $(11)$, $(11)$, $(22)$, $(22)$, $(33),$ respectively.  Since $V^1$ is MF, it follows that
at most one of $ \mu^0$ or $\mu^2$ is nonzero and applying Lemma \ref{twolabs} and using Magma to check the special cases we find that $\mu^i = 0$, $\lambda_1^i$, or 
$\lambda_5^i$  for $i = 0$ or $2.$  For future reference we note that $\wedge^2(11) = (03)+(30)+(11)$ and
$\wedge^3(11) = (03)+(30)+(11)+(22)+(00).$

Consider $V_{\gamma_1}^2(Q_Y).$  
First suppose that $\mu^1 = \lambda_3^1.$    Then $V_{\gamma_1}^2(Q_Y)$
contains $V_{C^0}(\mu^0 +\lambda_5^0) \otimes V_{C^1}(\lambda_4^1)$.  Now  $V_{C^1}(\lambda_4^1) \downarrow L_X' = (22)^2 + (11)^2,$ so tensoring with the $(22)$ terms we immediately obtain summands  of multiplicity
two, and a  check of $S$-values (using Lemma \ref{p.29}) yields a contradiction.  The same argument works if $\mu^1 = 2\lambda_1^1$ since $V^2(Q_Y)$
contains $V_{C^0}(\mu^0 +\lambda_5^0) \otimes V_{C^1}(\lambda_1^1+ \lambda_2^1)$ and the restriction
of the second factor contains $(22)^2.$  Similarly for $\mu^1 = 3\lambda_1^1$, using Magma to see that $(33)^2$
appears in the restriction of $V_{C^1}(2\lambda_1^1+ \lambda_2^1).$

So we are left with the cases $\mu^1 = \lambda_1^1$ and $\lambda_2^1.$
First suppose $\mu^0 = 0$ so that $\mu^2 = 0, \lambda_1^2$ or $\lambda_5^2$ which restricts
to $L_X'$ as $(ab) = (00), (20),$ or $(02)$, respectively.   If $\mu^1 = \lambda_1^1$, then the restriction of $V_{\gamma_1}^2(Q_Y)$ contains $(20) \otimes \wedge^2(11) \otimes (ab)$ which contains $(12)^3 \otimes (ab).$  However, it follows from Corollary \ref{V^2(Q_X)}   that at most one such irreducible arises from   $V^1(Q_Y)$  a contradiction.
Suppose $\mu^1 = \lambda_2^1.$  Then $\mu^2 = 0$, as otherwise $V^1$ is not MF. Here the restriction of  $V_{\gamma_1}^2(Q_Y)$ contains $V_{L_X'}(20) \otimes \wedge^3(11)$
and $(12)^4$ appears. But Corollary \ref{V^2(Q_X)} shows that only $(12)^2$ can arise from  $V^1(Q_Y)$, again
a contradiction.  Therefore $\mu^0 = \lambda_1^0$ or $\lambda_5^0$ and hence $\mu^2 = 0.$

Suppose  $\mu^0 = \lambda_1^0.$  This forces $\mu^1 = \lambda_1^1$, as otherwise $V^1$ would not be MF. Then $V_{\gamma_1}^2(Q_Y)$ contains $V_{C^0}(\lambda_1^0+\lambda_5^0)\otimes V_{C^1}(\lambda_2^1)$  and the restriction to $L_X'$ contains $(33)^3.$ 
We get a contradiction by considering $S$-values.  

Finally, assume  $\mu^0 = \lambda_5^0.$  Again we have $\mu^1 = \lambda_1^1$ and so  $V_{\gamma_1}^2(Q_Y)$  contains $V_{C^0}(2\lambda_5^0)\otimes V_{C^1}(\lambda_2^1)$.  The restriction to $L_X'$ contains $(51)^2$ and $S$-values yield a contradiction. \hal

\begin{lem}\label{mu20}Suppose that $\mu^2 = 0$.  Then $\langle \lambda, \gamma_1  \rangle = 0 = \langle \lambda,\gamma_2 \rangle.$
\end{lem}

\pf  We have $V^1(Q_Y) = V_{C^0}(\mu^0).$  Let $(ab)$ be an irreducible constitutent of $V^1$ for which $a+b$ is maximal. First consider those cases where $ab \ne 0.$ This includes most of the cases
listed in Lemma \ref{A2info}. 

Suppose $\langle \lambda, \gamma_2 \rangle \ne 0.$ Then the summand of  $V_{\gamma_2}^2(Q_Y)$ afforded by $\lambda - \gamma_2$   contains  $V_{C^0}(\mu^0) \otimes V_{C^1}(\lambda_7^1) \otimes V_{C^2}(\lambda_1^2)$.  Now $(V_{C^1}(\lambda_7^1) \otimes V_{C^2}(\lambda_1^2)) \downarrow L_X' = (11) \otimes (20)\supseteq (31).$  Therefore, Lemma \ref{abtimescd} shows that $V^2$
contains $((a+2)b)^2$, and an $S$-value argument gives a contradiction.  Therefore, $\langle \lambda, \gamma_2  \rangle = 0.$

Now suppose $\langle \lambda, \gamma_1  \rangle \ne 0.$ Then the summand of  
$V_{\gamma_1}^2(Q_Y)$ afforded by $\lambda - \gamma_1$  contains $V_{C^0}(\mu^0 + \lambda_5^0) \otimes V_{C^1}(\lambda_1^0).$  The argument in the proof of Lemma \ref{p.29} shows that the dominant weight $((a+2)b)$ is subdominant to the highest weight of an irreducible summand of $V_{C^0}(\mu^0 + \lambda_5^0) \downarrow L_X' \subseteq(V_{C^0}(\mu^0) \otimes V_{C^0}(\lambda_5^0)) \downarrow L_X'.$  For the group $A_2$, proper subdominant weights  of irreducibles have $S$-values strictly less than the highest weight. So maximality of $a+b$ implies
that $V_{C^0}(\mu^0 + \lambda_5^0) \downarrow L_X' \supseteq ((a+2)b) $ and $(V_{C^0}(\mu^0 + \lambda_5^0) \otimes V_{C^1}(\lambda_1^0)) \downarrow L_X' \supseteq ((a+2)b) \otimes (11).$  Then Lemma \ref{abtimescd} again shows that this contains $((a+2)b)^2$ and we again have a contradiction by considering $S$-values.

Now consider the situations where the maximal value of $a+b$ occurs only for $a = 0$ or $b = 0$.
If this weight has the form $(0b)$ then assuming $\langle \lambda, \gamma_2 \rangle \ne 0, $  we see that $V_{\gamma_2}^2(Q_Y)$  contains $(0b)\otimes  (11)\otimes (20)$ which contains $(2b)^2.$
This gives a contradiction using $S$-values. And if $\langle \lambda, \gamma_1 \rangle \ne 0$, we argue as above to  get $(2b) \otimes (11)$ which yields the same contradiction.

Finally, assume the highest weight has the form $(a0)$, which only occurs for $\mu^0 = c\lambda_5$ with
$a = 2c.$ If $\langle \lambda, \gamma_2 \rangle \ne 0, $ then the restriction of $V_{\gamma_2}^2(Q_Y)$  contains $S^c(20) \otimes (11) \otimes (20) \supseteq  ((2c-2)2)^2 \otimes (11) \supseteq ((2c-1)3)^2,$ provided $c > 1$,  where an $S$-value argument gives a contradiction.  And if $c = 1$, we have $(20) \otimes (11) \otimes (20) \supseteq  (32)^2 $ and we again have a contradiction.  Therefore $\langle \lambda, \gamma_2 \rangle = 0.$ Now assume $\langle \lambda, \gamma_1  \rangle \ne 0.$ Here the restriction of  $V_{\gamma_1}^2(Q_Y)$
contains $S^{c+1}(20) \otimes (11) \supseteq ((2c-2)2) \otimes (11) \supseteq (2c-2,2)^2.$  Here we do not get a contradiction using $S$-values but an application of  Corollary \ref{V^2(Q_X)} does yield a contradiction.  \hal

\begin{lem}\label{mu0not} Suppose that $\mu^2 = 0. $ Then  $\mu^0$ is not one of the following weights:
\[
\begin{array}{l}
0,\, a\lambda_1^0+\lambda_2^0,\, \lambda_1^0+\lambda_3^0, \,\lambda_1^0+\lambda_4^0,\, \lambda_1^0+\lambda_5^0,\,\lambda_2^0+\lambda_3^0,\, \lambda_2^0+\lambda_4^0, \,\lambda_1^0+2\lambda_5^0,\, \lambda_1^0+3\lambda_5^0,\\
 2\lambda_2^0, \,3\lambda_2^0,\, \lambda_j^0 \,(j = 4, 5),\, j\lambda_1^0 (j>2). 
\end{array}
\] 
\end{lem}

\pf  Suppose $\mu^2 = 0$ and the lemma is false. The previous lemma shows that $\langle \lambda, \gamma_1  \rangle = 0 = \langle \lambda,\gamma_2 \rangle.$ It follows from the hypotheses and Lemma \ref{020} that  $\mu^0 \ne a\l_1 + \l_2$.   By way of contradiction assume $\mu^0$ is one of the other listed weights.
 In view of previous lemmas we have 
$\lambda \in  \{0,  \lambda_1+\lambda_3, \lambda_1+\lambda_4, \lambda_2+\lambda_3, \lambda_2+\lambda_4, \lambda_1+2\lambda_5, \lambda_1+3\lambda_5,  2\lambda_2, 3\lambda_2, \lambda_j (j = 4,5), j\lambda_1 (j>2)  \}$ and we aim for a contradiction.
As $V$ is nontrivial, $\lambda \ne 0.$

Several of the remaining  weights are settled  using Magma.  For example, if $\lambda = \lambda_j (j = 4,5)$, then $V$ is the corresponding wedge power of $(020)$ and a Magma computation shows that this is not MF.  
If $\l = \l_1 + \l_3$, then $V = V_Y(\lambda_1) \otimes V_Y(\lambda_3) - V_Y(\lambda_4)$ and a Magma computation shows that $V \downarrow X$ is not MF.  Similarly for  $\lambda_1+\lambda_4$ and $\lambda_1+\lambda_5.$ For $\lambda_2+\lambda_3$ we have $V = (V_Y(\lambda_2) \otimes V_Y(\lambda_3)) - (V_Y(\lambda_1) \otimes  V_Y( \lambda_4))$ and a Magma computation shows this is not MF.  Likewise, $V_Y(2\lambda_2) =  (V_Y(\lambda_2) \otimes V_Y(\lambda_2)) - (V_Y(\lambda_1) \otimes V_Y(\lambda_3))$ and we see that $(222)^2$ appears in the restriction.

For $\lambda = \lambda_2+\lambda_4$, we consider $V_{\gamma_1}^2(Q_Y).$    This summand contains $V_{C^0}( \lambda_2^0+\lambda_3^0) \otimes V_{C^1}(\lambda_1^1)$ so restricting to $L_X'$ and using Lemma \ref{A2info} we  see that the restriction contains  $(42)^4$.  However, Corollary \ref{V^2(Q_X)} shows that at most  two summands $(42)$ can arise from  $V^1(Q_Y),$  so this is a contradiction.  
 
 Suppose $\lambda = \lambda_1 + 2\lambda_5.$  Then the $V_{\gamma_1}^2(Q_Y)$    contains $V_{C^0}( \lambda_1^0+\lambda_4^0 + \lambda_5^0) \otimes V_{C^1}(\lambda_1^1).$  Using
 Magma we check that $ V_{C^0}( \lambda_1^0+\lambda_4^0 + \lambda_5^0) = (V_{C^0}( \lambda_1^0+\lambda_4^0) \otimes  V_{C^0}(\lambda_5^0)) - V_{C^0}( \lambda_1^0+\lambda_3^0) - V_{C^0}(\lambda_4^0).$
 The information in Lemma \ref{A2info} now shows that $(32)^2$ appears in $ V_{C^0}( \lambda_1^0+\lambda_4^0 + \lambda_5^0) \downarrow L_X'$ and hence $(32)^4$ appears in the restriction of $V^2(Q_Y)$, contradicting
Corollary  \ref{V^2(Q_X)}.
 
 Similarly, suppose $\lambda = \lambda_1 + 3\lambda_5.$  Then   $V_{\gamma_1}^2(Q_Y)$    contains $V_{C^0}( \lambda_1^0+\lambda_4^0 + 2\lambda_5^0) \otimes V_{C^1}(\lambda_1^1).$  Using
 Magma we check that $ V_{C^0}( \lambda_1^0+\lambda_4^0 + 2\lambda_5^0) = (V_{C^0}( \lambda_1^0+2\lambda_5^0) \otimes  V_{C^0}(\lambda_4^0)) - V_{C^0}( \lambda_1^0+\lambda_3^0 + \lambda_5^0) - V_{C^0}(3\lambda_5^0) - V_{C^0}(\lambda_4^0 + \lambda_5^0) .$
 The information in Lemma \ref{A2info}  shows that $(52)^3$ appears in the tensor product but $(52)$  has multiplicity at most 1  in the subtracted terms. Therefore  $(52)^4$ appears in the restriction of $V^2(Q_Y)$, contradicting  Corollary \ref{V^2(Q_X)}.
 
 Next consider $\lambda = 3\lambda_2.$ Here  $V_{\gamma_1}^2(Q_Y)$ contains 
 $V_{C^0}( \lambda_1^0+2\lambda_2^0) \otimes V_{C^1}(\lambda_1^1)$.  Using the
fact that $ V_{C^0}( \l_1^0+2\l_2^0) = (V_{C^0}( \l_1^0) \otimes V_{C^0}(2\l_2^0))  - V_{C^0}(\l_2^0 +\l_3^0)$ and the information in Lemma \ref{A2info} we find that the restriction to $L_X'$ contains
 $(31)^4,$ and Corollary \ref{V^2(Q_X)} shows that  only $(31)^2$ can arise from $V^1(Q_Y).$
 
Finally, consider the case $\lambda = j\lambda_1$ for $j > 2$ where the module restricts to $S^j(020).$  For $j = 3$ we use  Magma to see that $(020)$ appears with multiplicity $2$.  For  $4 \le j \le 7$  a Magma computation shows that there is an irreducible of highest weight $(2(2j-6)2)$ that appears with multiplicity  $2.$  We claim
that this holds for $j > 7$ as well.  

Let $j > 7$ and set $\psi = (0(2j)0)$, so that
$(2(2j-6)2) = \psi -(\alpha_1+4\alpha_2+\alpha_3).$  In $S^j(020)$ any symmetric tensor which results in a
weight which has the form $\psi - (c_1\a_1+c_2\a_2+c_3\a_3)$ with $c_1+c_2+c_3 \le 6$ must be the symmetric product of at least $j-6$ copies of $\d = 020 $ followed by terms where certain roots are subtracted from $\d$.
It follows that the multiplicity of $(2(2j-6)2) = \psi -(\alpha_1+4\alpha_2+\alpha_3)$ in $S^j(020)$ equals
the multiplicity of $(262)$ in $S^6(020)$ and a Magma computation shows that this is 2.  \hal
 
 \begin{lem}\label{mu0notomega5} $\mu^0 \ne \lambda_5^0$ and $\mu^2 \ne \lambda_1^2.$
 \end{lem}
 
 \pf  Using duals it suffices to  prove the first assertion. Suppose $\mu^0 = \lambda_5^0$.   By Lemma \ref{mu0not}, $\mu^2 \ne 0.$ 
 Then  from Lemmas \ref{inducl2},   \ref{A2info},  \ref{nonemf}, and a Magma computation we see that
 $\mu^2 = \lambda_1^2,$ $\lambda_5^2,$ $\l_2^2,$ or $\l_4^2.$  Indeed, otherwise, $V^1(Q_Y)$ is not MF.  Therefore, $V^1 = (20) \otimes (20),$ $(20) \otimes (02),$ $(20) \otimes (21)$, or 
$ (20) \otimes (12).$
 
Now $V_{\gamma_1}^2(Q_Y)$ contains the irreducible summand afforded by $\lambda -\beta_5^0 - \gamma_1$ and $L_Y'$
  acts as $V_{C^0}(\lambda_4^0) \otimes V_{C^1}(\lambda_1^1) \otimes V_{C^2}(\mu^2). $  
  Then a Magma computation shows that
 $V_{\gamma_1}^2(Q_Y) \downarrow L_X'$ contains $(41)^4,$  $(23)^4,$ $(42)^6$, or $(22)^6$, according to whether
  $\mu^2 = \lambda_1^2,$ $\lambda_5^2,$ $\l_2^2,$ or $\l_4^2.$  In each case this contradicts Corollary \ref{V^2(Q_X)}.   \hal

 \begin{lem}\label{lambda1lambdan} If $\mu^0 \ne 0 \ne \mu^2$, then $\lambda = \lambda_1 + \lambda_n.$
 \end{lem}
 
 \pf  Suppose $\mu^0 \ne 0 \ne \mu^2.$   Taking duals, if necessary, and applying the last lemma, the inductive hypothesis,  Lemma \ref{twolabs}, Lemma \ref{A2info}, and Magma we may assume  that $\mu^0 = \lambda_1^0, $ $\lambda_2^0$,  or $\lambda_4^0$ and $\mu^2 = \lambda_5^2.$ Indeed otherwise, $V^1$ is not MF.  If $\mu^0 = \lambda_2^0$  or $\lambda_4^0,$  then we obtain
 a contradiction within $V_{\gamma_1}^2(Q_Y)$ in the usual way.  So we now assume $\mu^0 = \lambda_1^0.$
 It remains to show that $\langle \lambda, \gamma_1  \rangle = 0 = \langle \lambda, \gamma_2  \rangle$ and
 by dualty it suffices to show that $\langle \lambda, \gamma_1  \rangle = 0.$  Otherwise $V_{\gamma_1}^2(Q_Y) \supseteq V_{C^0}(\l_1^0 + \l_5^0) \otimes V_{C^1}(\l_1^1) \otimes V_{C^2}(\l_5^2)$ and restricting to
 $L_X'$ we have a summand $(22) \otimes (11) \otimes (02)$ which contains $(24)^2.$    This is a contradiction since  $V^1 = (02) \otimes (02).$   \hal
 
\vspace{4mm}
 We now aim to complete the proof of  Theorem \ref{main}. Taking duals, if necessary,  and using  Lemmas \ref{mu1is0},  \ref{mu20}, and \ref{lambda1lambdan}  we may now assume
 that $V^1(Q_Y) = V_{C^0}(\mu^0),$ $\mu^0 \ne 0$, and   $\la \lambda, \gamma_i  \ra = 0$ for $i = 1,2.$
 Using previous lemmas, we are done unless $\lambda \in \{ j\lambda_5 (j>1)$, $\l_2$, $\l_3$, $\lambda_4+a\lambda_5$ ($a \le 3$), $\lambda_3+\lambda_5$, $\lambda_2+\lambda_5$, $\lambda_3+\lambda_4$, $2\lambda_1+\lambda_5$, $3\lambda_1+\lambda_5$, $2\lambda_4$, $3\lambda_4 \}. $  In view of Lemma \ref{020} we can rule
 out the case $\l = \lambda_4+a\lambda_5$, and both $\l_2$ and $\l_3$ are included in the conclusion of 
Theorem \ref{main}.

For the other cases we proceed as follows. For the moment exclude the case $\lambda = j\lambda_5,$ for $j \ge 3.$
Then  $V_{C^0}(\mu^0) \downarrow L_X'$ is given by taking duals in Lemma \ref{A2info}.  We
consider  $V_{\gamma_1}^2(Q_Y)$ in the usual way.    In most cases the dual of the restriction of this irreducible is also given in Lemma \ref{A2info} and using Lemma \ref{abtimescd} and
 Corollary \ref{V^2(Q_X)} we obtain a contradiction.  For example if $\lambda = 2\lambda_5$
the irreducible summand afforded by $\lambda - \beta_5 - \gamma_1$ is $V_{C^0}(\lambda_4^0 + \lambda_5^0) \otimes V_{C^1}(\lambda_1^1).$  The second tensor factor restricts to $(11)$, while the first
restricts to $(41)+(11)+(22).$  Taking tensor products with $V_{C^1}(\lambda_1^1) \downarrow L_X' = (11)$ we find that  $(33)^2$ appears in 
$V^2$ contradicting Corollary \ref{V^2(Q_X)}.  Using the same arguments with $\lambda = \lambda_3+\lambda_5, \lambda_2+\lambda_5, 2\lambda_4$ we find
that  $(32)^4$, $(22)^5$, $(51)^3$ occur, and we obtain a contradiction.

If $\lambda = \lambda_3 + \lambda_4$, then  $V_{\gamma_1}^2(Q_Y)$
contains $V_{C^0}(2\lambda_3^0) \otimes V_{C^1}(\lambda_1^1).$ Using Magma we check that  $V_{C^0}(2\lambda_3^0) = S^2(V_{C^0}(\lambda_3^0)) - V_{C^0}(\lambda_1^0 + \lambda_5^0).$  Another application of Magma
shows that $V_{C^0}(2\lambda_3^0) \downarrow L_X' \supseteq ((33) + (22)^2).$  Therefore, the restriction of $V^2(Q_Y)$
to $L_X'$ contains $((33) +(22)^2) \otimes(11) \supseteq (22)^5.$
Now Lemma \ref{A2info} gives the restriction of $V^1(Q_Y)$ and we contradict  Corollary \ref{V^2(Q_X)}. 

Suppose $\lambda = 2\lambda_1+\lambda_5$.  Here $V^2(Q_Y)$ contains an irreducible summand with highest weight afforded by $\lambda - \beta_5^0 -\gamma_1$
and this affords $V_{C^0}(2\lambda_1^0+\lambda_4^0) \otimes V_{C^1}(\lambda_1^1).$  Now
$V_{C^0}(2\lambda_1^0+\lambda_4^0) = (V_{C^0}(2\lambda_1^0) \otimes V_{C^0}(\lambda_4^0)) - 
V_{C^0}(\lambda_1^0 + \lambda_5^0).$ Restricting to $L_X'$ we have $(((04) + (20)) \otimes (21)) - ((22)+(11)+(00))$.  Using Magma we see that this contains $(22) + (03) +(33)$, so that $V^2$ contains $((22) + (03) +(33)) \otimes (11)$ which contains $(22)^4$.  This contradicts Corollary \ref{V^2(Q_X)}.  

Assume $\lambda = 3\lambda_1+\lambda_5$.  Then $V^2(Q_Y)$ contains an irreducible summand with highest weight afforded by $\lambda - \beta_5^0 -\gamma_1$
and this affords $V_{C^0}(3\lambda_1^0+\lambda_4^0) \otimes V_{C^1}(\lambda_1^1).$  Now
$V_{C^0}(3\lambda_1^0+\lambda_4^0) = (V_{C^0}(3\lambda_1^0) \otimes V_{C^0}(\lambda_4^0)) - 
V_{C^0}(2\lambda_1^0 + \lambda_5^0).$ Using Magma and  Lemma \ref{A2info} we see that the restriction
to  $V_{C^0}(3\lambda_1^0+\lambda_4^0) \downarrow L_X'$ contains
$(32)^5+(21)$ so that  $V^2\supseteq (32)^6,$ a contradiction.

Suppose $\lambda = 3\lambda_4.$ In this case $V^2(Q_Y) \supseteq  V_{C^0}(\lambda_3^0+2\lambda_4^0) \otimes V_{C^1}(\lambda_1^1).$  From Magma we check that $V_{C^0}(\lambda_3^0+2\lambda_4^0) = (V_{C^0}(\lambda_3^0) \otimes V_{C^0}(2\lambda_4^0)) - V_{C^0}(\lambda_2^0+\lambda_4^0 +\lambda_5^0) -
V_{C^0}(\lambda_1^0+\lambda_4^0 ).$ Using Lemma \ref{A2info}, we see
that $(V_{C^0}(\lambda_3^0) \otimes V_{C^0}(2\lambda_4^0)) \downarrow L_X' \supseteq (61)^2$ and $(61)$ does not appear in the summands deleted from the tensor product.  Therefore, $V^2 \supseteq (61)^2 \otimes (11) \supseteq (61)^4$ and this contradicts Corollary \ref{V^2(Q_X)}.

Finally we must consider $\lambda = j\lambda_5$ for $j > 2.$  We will consider the summand of $V^3(Q_Y)$ afforded by $\lambda  - 2\beta_5^0 -2\gamma_1$  which affords $V_{C^0}(2\lambda_4^0 + (j-2)\lambda_5^0) \otimes V_{C^1}(2\lambda_1^1).$   Now $\lambda_5^0$ restricts to $(20)$ and $\lambda_4^0$ restricts to $(21)$. It follows that there is a maximal vector in $V_{C^0}(2\lambda_4^0 + (j-2)\lambda_5^0)$ whose weight upon the restriction to the maximal torus $S_X$ of $L_X'$ is $((2j-4)0) + (4,2)) = ((2j)2).$

It now follows that 
$V^3 \supseteq ((2j)2) \otimes S^2(11) \supseteq ((2j)2) \otimes ((22) +(11))\supseteq ((2j+1)3)^3.$
 On the other hand $V^2(Q_Y)  = V_{C^0}(0001(j-1)) \otimes V_{C^1}(\lambda_1^1)$ and restricting to  $L_X'$ this is contained in $ S^{j-1}(20) \otimes \wedge^2(20) \otimes (11)$ which has highest weight $((2j+1)2)$ and $S$-value $2j+3.$ All other dominant weights have smaller $S$-values.  Consequently only one composition factor   $((2j+1)3)$ can possibly arise from $V^1(Q_Y)$. Therefore, we obtain a contradiction. This contradiction completes the proof of  Theorem \ref{main}.



\chapter{The case $\d = r\omega_1$, $r\ge 2$}\label{ro1}

In this chapter we prove Theorem \ref{MAINTHM} in the case where $\d = r\o_1$ with $r\ge 2$. 
Recall our basic notation: $X = A_{l+1}$, $W = V_X(\d)$, $Y = SL(W) = A_n$, and $V = V_Y(\l)$ such that $V \downarrow X$ is MF. A fundamental root system for $X$ is denoted $\Pi(X) = \{\alpha_1,\dots,\alpha_l,\alpha_{l+1}\}$, and $\Pi(L_X') = \{\alpha_1,\dots,\alpha_l\}$, with
corresponding fundamental dominant weights $\{\omega_1,\dots,\omega_{l+1}\}$ (viewed as 
weights for $L_X$ and  $L_X'$, as well).

We divide the analysis into two subsections -- the cases $r=2$ and $r\ge 3$.

\section{The case $\delta = 2\omega_1$}\label{sec:2omega1}

Set $\delta = 2\omega_1$. Here there are two levels $W^1(Q_X)$ and $W^2(Q_X)$ on which $L_X'$ acts irreducibly with highest weights $2\o_1$ and $\o_1$, respectively, so $L_Y' = C^0\times C^1 \cong A_{r_0}\times A_l$, where $r_0 = \frac{(l+1)(l+2)}{2}-1$.  We write $\Pi(Y) = \{\beta_1,\dots,\beta_n\}$ and 
$\Pi(C^0) = \{\beta_1^0,\dots,\beta_{r_0}^0\}$, so $\beta_i^0 = \beta_i$ for $1\leq i\leq r_0$. Set 
$\gamma_1 = \beta_{r_0+1}$ and $\gamma_2 =\beta_n$ and finally 
set $\beta_{l-i+1}^1=\beta_{n-i}$, for $1\leq i\leq l$, so $\Pi(C^1) = 
\{\beta_1^1,\dots,\beta_l^1\}$.
The corresponding fundamental dominant weights are denoted $\lambda_j^i$, $i=0,1$, 
$1\leq j\leq {\rm rank}(C^i)=r_i$.

We establish the following theorem. 
\begin{thm}\label{thm:2omega1} Let $X=A_{l+1}$ and $\delta=2\omega_1$. Suppose $V_Y(\lambda)\downarrow X$ is
 multiplicity-free, where $\l \ne 0,\l_1,\l_n$. Then $\lambda$ or $\lambda^*$ is as in Tables $\ref{TAB2}-\ref{TAB4}$ of Theorem $\ref{MAINTHM}$.
\end{thm}

 \subsection{Proof of Theorem~\ref{thm:2omega1}}

We now begin the proof of Theorem~\ref{thm:2omega1}. By Theorem \ref{a2r0}, the result holds for $l=1$,
 so we now assume 
$l\geq 2$ and that the result holds for any embedding  
$A_{m+1}\subseteq SL(V_{A_{m+1}}(2\omega_1))$, 
with $m<l$. We refer to the list of possibilities for $\l$ in smaller rank cases as the inductive list.

As in Theorem \ref{thm:2omega1}, suppose $V = V_Y(\l)$ and $V\downarrow X$ is MF.
 The proof is accomplished in a sequence of propositions, treating different 
configurations for the weight $\lambda$.

\begin{lem}\label{V1triv} Theorem $\ref{thm:2omega1}$ holds If $V^1(Q_Y)$ is the trivial $L_Y$-module.
\end{lem}

\pf Suppose $V^1(Q_Y)$ is trivial. Set $x=\langle \lambda,\gamma_1\rangle$ and $y=\langle \lambda,\gamma_2\rangle$. 
If $x=0$, then $\lambda^*$ is as in Table \ref{TAB3}, so we assume from now on that $x>0$;
in particular, $V^2$ has a summand 
$(2\o_1)^*\otimes \omega_1 = (\omega_1+2\omega_l)\oplus \omega_l$ and if $y\ne 0$,   then $V^2$ has an
additional summand $\omega_l$.

Consider first the case where $y=0$. If $x=1$ then $\l$ is in Table \ref{TAB3}, and so 
we may assume $x>1$. 
Applying the induction hypothesis to $W=V^*$,
we reduce to the case $l=2$. Then $V^2= (\omega_1+2\omega_2)\oplus \omega_2$, and $V^3(Q_Y)$ has $L_Y$-summands
afforded by $\lambda-2\gamma_1$, $\lambda-\gamma_1-\beta_1^1-\beta_2^1-\gamma_2$, and 
$\lambda-\beta_{r_0}^0-2\gamma_1-\beta^1_1$,
giving rise to $S^2(2\o_2)\otimes 2\omega_1$, respectively $2\o_2$, 
$\wedge^2(2\o_2)\otimes \omega_2$.
Decomposing the first tensor product produces two summands $2\o_2$, and the last 
tensor product produces a third such summand,  and hence $V^3 \supseteq (2\o_2)^4$. Only two of 
these summands can arise from summands of $V^2$ and so we obtain a contradiction.

We may now assume that $xy\ne 0$; if $x=y=1$, the induction hypothesis applied to $V^*$ 
shows that $l\leq 5$ and so $\lambda$ is as in Table \ref{TAB2}. So we
now assume at least one of $x,y$ is different from 1.
Recall  $V^2 = (\omega_1+2\omega_l)\oplus \omega_l \oplus \omega_l$.

We claim that $x=1$. For otherwise, $V^3(Q_Y)$ has summands afforded by $\lambda-2\gamma_1$, 
$\lambda-\gamma_1-\gamma_2$, and two further summands, each afforded by the weight
$\lambda-\gamma_1-\beta_1^1-\cdots-\beta_l^1-\gamma_2$. Restricting these summands to $L_X'$
produces five summands $2\omega_l$. However, only three such summands can arise from 
$V^2$, contradicting Proposition \ref{induct}. Hence $x=1$ as claimed and so we have $y>1$. But now applying the induction hypothesis  to $V^*$ produces a contradiction, thus completing the proof.\hal

In view of the previous lemma, we assume from now on that $V^1(Q_Y)$ is nontrivial.

\begin{lem}\label{V1nt_1} We have $\langle\lambda,\gamma_1\rangle\cdot\langle\lambda,\gamma_2\rangle=0.$
\end{lem}

\pf Set $\langle\lambda,\gamma_1\rangle = x$ and
$\langle\lambda,\gamma_2\rangle=y$ and suppose $xy\ne 0$. 
Applying the induction hypothesis to the dual module $V^*$, we find that the following conditions hold:
\begin{itemize}
\item[] $\mu^1=0$, and hence $\mu^0\ne 0$,
\item[] $x=y=1$, and
\item[] $\mu^0$ is supported on the first $l+2$ nodes, i.e. those nodes corresponding to the 
roots 
$\beta_1^0,\dots,\beta_{l+2}^0$.
\end{itemize}

Now $V^2(Q_Y)$ has a summand afforded by $\lambda-\gamma_2$, and Corollary~\ref{cover}
shows that this affords all $L_X'$-summands of $V^2(Q_Y)$ arising from $\sum_{n_i=0} V_i$.
Hence, if any other $L_Y$-summands of $V^2(Q_Y)$ have non-multiplicity-free restriction to 
$L_X'$, we obtain a 
contradiction. We will use Proposition~\ref{v2gamma1} throughout; that is, $V^2(Q_Y)$ has a submodule of the form $V_{C^0}(\mu^0)\otimes V_{C^0}(\lambda_{r_0}^0)\otimes V_{C^1}(\lambda^1_1)$. Restriction of this summand to $L_X'$ affords the following $L_X'$-summand of $V^2$:
\[
V_{C^0}(\mu^0)\downarrow {L_X'}\otimes 2\omega_l\otimes \omega_1 = 
V_{C^0}(\mu^0)\downarrow {L_X'}\otimes ((\omega_1+2\omega_l)\oplus \omega_l).
\]
By Proposition \ref{tensorprodMF} and the above remarks, if $V_{C^0}(\mu^0)\downarrow {L_X'}$ has any irreducible 
summand with two nonzero labels, 
we obtain a contradiction. 

We consider one-by-one, the possibilities for $\mu^0$, given by the inductive hypothesis. We apply  
Lemma~\ref{twolabs} and the above remarks, and reduce to the case $\mu^0=a\lambda_1^0$, for 
$a\leq 2$. But now
we have the summand $2\omega_1\otimes 2\omega_l\otimes \omega_1$, if $a=1$,
and otherwise a summand $S^2(2\omega_1)\otimes 2\omega_l\otimes \omega_1$.
Each of these can easily be seen to be non-MF, completing the proof of the lemma.\hal

\begin{lem}\label{V1nt_2} If $\langle\lambda,\gamma_1\rangle +
\langle\lambda,\gamma_2\rangle\ne 0$, then $\lambda$ or $\lambda^*$ is in Tables $\ref{TAB2}-\ref{TAB4}$ of Theorem $\ref{MAINTHM}$.
\end{lem}

\pf First set $x = \langle\lambda,\gamma_1\rangle$ and
$y=\langle\lambda,\gamma_2\rangle$. By Lemma~\ref{V1nt_1}, $xy=0$.
Applying the inductive hypothesis to the dual module $V^*$, we see that one of the following holds:
\begin{enumerate}[(A)]
\item $l\geq3$, $(x,y)=(1,0)$ and $\mu^1 = 0$,
\item $l\geq 3$ and $(x,y) = (0,a)$, or
\item $l=2$.
\end{enumerate}

We will treat these three cases separately.

\noindent{\bf{Case (A):}} Here we have $l\geq3$, $(x,y)=(1,0)$ and $\mu^1 = 0$.
The aim here is to show that $\mu^0=\lambda_1^0$ and $l\leq4$. 
Applying the inductive hypothesis
to the modules $V$ and $V^*$, we deduce that $\mu^0$ lies in the set of weights
\[
\begin{array}{rl}
b\lambda_1^0,\,& b\geq 1,\\
\lambda_j^0,\, &2\leq j\leq l+3,\\
\lambda_1^0+\lambda_t^0,\,&t\leq {\rm min}\{7,l+3\},\\
a\lambda_2^0, a\lambda_1^0+\lambda_2^0,\, &a=2,3, \\
\lambda_2^0+\lambda_3^0.&
\end{array}
\]
Moreover, as in the proof of the preceding lemma,
we apply Proposition~\ref{v2gamma1} to produce $L_X'$-summands of $V^2$ of the form
\[
V_{C^0}(\mu^0)\downarrow {L_X'}\otimes (\omega_1+2\omega_l),\mbox{ and }
V_{C^0}(\mu^0)\downarrow {L_X'}\otimes \omega_l.
\]
Now Corollary~\ref{cover} implies that the first of these summands must be multiplicity-free,
and Lemma~\ref{twolabs} and Proposition~\ref{tensorprodMF} then yield that 
$\mu^0=\lambda_1^0$ or $2\lambda_1^0$. In the latter case, we can easily see that the first summand 
is not multiplicity-free and so we have reduced to $\mu^0=\lambda_1$, as desired.
To see that $l\leq4$ we apply Lemma~\ref{nonMF12}.  This then completes the consideration of Case (A).

\noindent{\bf{Case (B):}} Here we have $l\geq 3$, $(x,y) = (0,a)$.

We first suppose $\mu^1\ne 0$.
Then considering $V^*$, we deduce that either (i) $a=1$ and $\mu^1=\lambda_j^1$ for some $j$,
or (ii) $a=2$ or $3$ and $\mu^1=\lambda^1_l$.
Moreover, by case (A), Lemma~\ref{V1nt_1}, and the induction hypothesis 
(applied to $V$ and $V^*$), we may assume  that the support of $\mu^0$ lies in the set 
$\{\beta_1^0,\ldots ,\beta_{l+1}^0\}$. 

Consider now the case (ii). Note that if $\mu^0=0$, then $\lambda^*$ is as in Tables \ref{TAB2}-\ref{TAB4}; so 
we assume $\mu^0\ne 0$. Here $\lambda-\gamma_2$ and $\lambda-\beta_l^1-\gamma_2$ afford 
$L_Y'$-summands of $V^2(Q_Y)$, the sum of which has restriction to $L_X'$ being 
$$V_{C^0}(\mu^0)\downarrow L_X'\otimes V_{C^1}(\mu^1)\downarrow L_X'\otimes \omega_l.$$ 
Then applying 
Corollary~\ref{cover}, we see that any other summand of $V^2(Q_Y)$ must be a multiplicity free 
$L_X'$-module. Now $\nu=\lambda-\gamma_1-\beta_1^1-\cdots-\beta^1_l$ affords a summand 
$V_{C^0}(\mu^0+\lambda_{r_0}^0)$, so applying the inductive hypothesis we deduce that 
$\mu^0=b\lambda_1^0$ for $1\leq b\leq 3$ or $\mu^0=\lambda^0_j$ for $2\leq j \le {\rm min}\{6,l+1\}$. In the first case,
 the weights $\lambda-\gamma_1-\beta_1^1-\cdots-\beta^1_l$ and 
$\lambda-\beta_1^0-\cdots-\beta^0_{r_0}-\gamma_1-\beta_1^1-\cdots-\beta^1_l$ afford summands the
 sum of which is isomorphic to $S^b(2\omega_1)\otimes 2\omega_l$ and then Lemma~\ref{nonemf} implies that $b=1$. Now adding to these two summands the summand afforded by $\lambda-\beta_1^0-\cdots-\beta^0_{r_0}-\gamma_1$, we again obtain a repeated $L_X'$-summand. So finally we have reduced to 
the case where $\lambda=\lambda_j+\lambda_{n-1}+a\lambda_n$, for $2\leq j\leq 6$ and $a=2,3$. 
The final contradiction in this case will come from considering the summands of $V^2(Q_Y)$ afforded 
by the weight $\nu$ given above and two further summands afforded by $\nu'=\lambda-\beta_j^0-\cdots-\beta_{r_0}^0-\gamma_1$ and 
$\nu''=\lambda-\beta_j^0-\cdots-\beta_{r_0}^0-\gamma_1-\beta_1^1-\cdots-\beta_l^1$. If $j>2$, 
Lemmas~\ref{twolabs} and \ref{tensorprodMF} imply that the summand afforded by $\nu'$ is non-MF.
And if $j=2$, the sum of the summands afforded by $\nu'$ and $\nu''$ is a non-MF $L_X'$-module. This completes the consideration of (ii).

Now we continue with the assumption that $\mu^1\ne0$ and consider case (i), where $a=1$ and  $\mu^1=\lambda_j^1$ for some $j$. Hence we have either $\mu^0=0$ or $\mu^0$ lies in the set
\[
\begin{array}{rl}
b\lambda_1^0,\,& b\geq 1,\\
\lambda_j^0 , \, & 2\leq j\leq l+1, \\
\lambda_1^0+\lambda_j^0, \, & 2\leq j\leq{\rm min}\{7,l+1\},\\ 
c\lambda_2^0,\,\lambda_2^0+\lambda_3^0,c\lambda_1^0+\lambda_2^0,\, & c=2,3.
\end{array}
\]
The cases where $\mu^0$ is $\l_1^0+\l_j^0$, $\l_j^0$ or $c\l_2^0$  are ruled out by Lemma~\ref{nonMF} (applied to $V^*$ in the second case); also $\l_2^0+\l_3^0$ is excluded by applying Lemmas~\ref{twolabs} and \ref{tensorprodMF} to $V^*$. 

Assuming $\mu^0\ne 0$, this leaves the possibilities $\mu^0=b\lambda_1^0$ or $c\lambda_1^0+\lambda_2^0$.

Suppose $\mu^0=b\lambda_1^0$, and set $M: = V^*$, so $M$ has highest weight 
$\l^*=\lambda_1+\lambda_i+b\lambda_n$, for some
$2\leq i\leq l+1$. Now $ \l^*-\gamma_2$ affords an $L_X'$-summand of $M^2$ of the form 
$M^1\otimes \omega_l$. So by Corollary~\ref{cover}, any remaining summands
of $M^2$ must be multiplicity-free. The weights
 $\l^*-\beta_i^0-\beta_{i+1}^0-\cdots-\gamma_1$ and
 $\l^*-\beta_1^0-\beta_2^0-\cdots-\gamma_1$ afford  
$L_X'$-summands of $M^2$ the sum of which is precisely the $L_X'$-module 
$2\omega_1\otimes \wedge^{i-1}(2\omega_1)\otimes \omega_1$. Since 
$2\omega_1\otimes \omega_1$ has a summand 
$(\omega_1+\omega_2)$, Proposition~\ref{tensorprodMF} and Lemma~\ref{twolabs} show that 
$2\omega_1\otimes \wedge^{i-1}(2\omega_1)\otimes \omega_1$ is not MF for $i>2$, 
 and for $i=2$ it is easy to check that 
$2\omega_1\otimes \wedge^{i-1}(2\omega_1)\otimes \omega_1$ is not MF. Hence $\mu^0\ne b\lambda_1^0$.

Now assume $\mu^0=c\lambda_1^0+\lambda_2^0$, $c=2,3$.
Here $\lambda-\gamma_2$ and $\lambda-\beta_j^1-\cdots-\gamma_2$ afford summands of $V^2(Q_Y)$, 
the sum of which, restricted to $L_X'$,  is $V^1\otimes \omega_l$. So then 
Corollary~\ref{cover} implies that any other summand of $V^2(Q_Y)$ must be a multiplicity 
free $L_X'$-module.
But $\lambda-\gamma_1-\beta_1^1-\cdots-\beta^1_{j}$ affords a non-multiplicity free summand by 
the induction hypothesis.

We have now shown that if $\mu^1\ne 0$, then $\mu^0=0$. But then we apply induction to $V^*$ and 
see that $\lambda^*$ is as in Tables \ref{TAB2}-\ref{TAB4} of Theorem \ref{MAINTHM}.

Henceforth, we will assume that $\mu^1=0$ (still in Case (B)).
Let us make a few general remarks.  By Case A, applied to $M=V^*$,
we may assume $\langle\lambda^*,\gamma_1\rangle=0$, that is 
$\langle\lambda,\beta_{l+2}\rangle=0$. Note that if 
$\langle\lambda,\beta_1\rangle\ne 0$, we may assume, by the first
part of the consideration of Case B applied to $M$, that 
$\langle\lambda,\beta_j\rangle=0$ for $2\leq j\leq l+1$. 

We now apply the induction hypothesis both to $V$ and $M=V^*$, and eliminate all
possibilities where $\lambda$ or $\lambda^*$ has been covered in Case A or 
by the above discussion. These considerations
 allow us to reduce to the following list: 
\begin{enumerate}[(1)]
\item $\mu^0 =b\lambda_1^0$, $b\geq 1$.
\item $\mu^0=\lambda_j^0$, for $2\leq j\leq l+1$.
\item $\mu^0=\lambda_{l+3}^0$ and $a\leq 3$.
\item $\mu^0=\lambda_j^0$, for $j>l+3$ and $a=1$.
\item $\mu^0=\lambda_1^0+\lambda_j^0$, for $j>l+3$ and $a=1$.
\item $\mu^0=\lambda_1^0+\lambda_{l+3}^0$ and $a\leq3$.
\item $\mu^0=b\lambda_2^0$, $b=2,3$.
\item $\mu^0=\lambda_2^0+\lambda_3^0$, and applying Lemma~\ref{twolabs}(1),
we deduce that $a\leq 2$.
\item $\mu^0 = \lambda_2^0+\lambda_{r_0-1}^0$ and $a=1$.
\item $\mu^0=b\lambda_1^0+\lambda_{r_0}^0$, $a=1$ and $b=2,3$.
\item $\mu^0 = \lambda_j^0+\lambda_{r_0}^0$, $a=1$ and $2\leq j\leq {\rm min}\{6,l+1\}$.
\item $\mu^0=b\lambda_1^0+\lambda_2^0$, $b=2,3$.
\end{enumerate}

Now we will repeatedly apply Corollary~\ref{cover} as in Case (A). In particular,
since $ \lambda-\gamma_2$ affords an $L_X'$-summand of $V^2$ of the form 
$V^1\otimes\omega_l$, any remaining summands
of $V^2$ must be multiplicity-free. This quickly rules out cases (9), (10) and (11), 
as well as (7) when $b=3$.

The configuration of case (12) is also quite easy; we have summands afforded by 
$\lambda-\beta_2^0-\cdots-\beta_{r_0}^0-\gamma_1$ and  
$\lambda-\beta_1-\beta_2^0-\cdots-\beta_{r_0}^0-\gamma_1$, the sum of which affords 
$S^b(2\omega_1)\otimes 2\omega_1\otimes \omega_1$, and Lemma~\ref{nonemf} 
provides the contradiction.

Cases (5) and (6) can be treated simultaneously; set $j=l+3$ in case (6). The weights 
$\lambda-\beta^0_j-\cdots-\beta^0_{r_0}-\gamma_1$ and $\lambda-\beta^0_1-\cdots-\beta^0_{r_0}-\gamma_1$ 
each afford 
irreducible $L_Y$-summands of $V^2(Q_Y)$, the sum of which
affords the module $V_{C^0}(\lambda_1^0)\otimes V_{C^0}( \lambda_{j-1}^0)\otimes 
V_{C^1}(\lambda_1^1).$ Restricting this to $L_X'$, we obtain the $L_X'$-module
$2\omega_1\otimes \wedge^{j-1}(2\omega_1)\otimes \omega_1$.
Since $(\omega_1+\omega_2)$ is an irreducible summand of the tensor product $2\omega_1\otimes \omega_1$, using Proposition~\ref{tensorprodMF} and Lemma~\ref{twolabs} we see that the three-fold tensor product is not multiplicity-free, hence ruling out these configurations.

For case (7) when $b=2$, we have a summand of $V^2(Q_Y)$ afforded by
the weight $\lambda-\beta_2-\cdots-\beta_{r_0}-\gamma_1$, which upon 
restriction to $L_X'$
is seen to be non multiplicity-free by Lemma~\ref{nonMF}(4).

For case (8), we note that $V^2(Q_Y)$ has a $L_Y'$-summand 
$V_{C^0}(\lambda_1^0+\lambda_3^0)\otimes V_{C^1}(\lambda_1^1)$.
This is isomorphic to $\wedge^2(V_{C^0}(\lambda_2^0))\otimes V_{C^1}(\lambda_1^1)$,
and upon restriction to $L_X'$ we obtain an $L_X'$-summand 
$\wedge^2(2\omega_1+\omega_2)\otimes \omega_1$, which is easily checked to be non-MF.

 We now consider case
(2) when $a>1$. First suppose that $2<j<l+1$ and  consider
the module $M=V^*$. Now 
$M^1 = S^a(2\omega_1)\otimes \omega_m$, where $2\leq m\leq l-1$,
while $M^2$ has summands 
\begin{itemize}
\item[(i)] $V_{C^0}(a\lambda_1^0+\lambda_{r_0}^0)\downarrow L_X'\otimes \omega_{m+1}$,
\item[(ii)] $S^{a-1}(2\omega_1)\otimes (\omega_1+\omega_m)$,
\item[(iii)] $S^{a-1}(2\omega_1)\otimes \omega_{m+1}$,
\item[(iv)] $S^a(2\omega_1)\otimes \omega_{m-1}$.
\end{itemize}
The sum of the first and third summands is isomorphic to 
$S^a(2\omega_1)\otimes 2\omega_l\otimes \omega_{m+1}$,
which in turn is isomorphic to 
$$\bigl(S^a(2\omega_1)\otimes (\omega_{m+1}+2\omega_l)\bigr)
\oplus \bigl(S^a(2\omega_1)\otimes (\omega_m+\omega_l)\bigr).$$  
On the other hand,
\begin{equation}\label{eq1}
M^1\otimes \omega_l = S^a(2\omega_1)\otimes \omega_m\otimes \omega_l
=\bigl(S^a(2\omega_1)\otimes (\omega_m+\omega_l)\bigr)
\oplus \bigl(S^a(2\omega_1)\otimes (\omega_{m-1})\bigr).
\end{equation}
It now suffices to see that 
the $L_X'$-module 
\[
\bigl(S^a(2\omega_1)\otimes (\omega_{m+1}+2\omega_l)\bigr)
\oplus \bigl(S^{a-1}(2\omega_1) \otimes (\omega_1+\omega_m)\bigr)
\]
 is not multiplicity-free, as  Corollary~\ref{cover} then produces the desired contradiction. 
This follows from Lemma~\ref{twolabs},
if $a\geq 4$ or if $a=3$ and $m<l-1$. Now if $a=3$ and $m=l-1$, the first summand
$S^3(2\omega_1)\otimes 3\omega_{l}$ contains $4\omega_1+\omega_l$ with
 multiplicity 2, 
and if $a=2$ and $m\leq l-1$,
the summand $S^2(2\omega_1)\otimes (\omega_{m+1}+2\omega_l)$ contains 
$(2\omega_1+\omega_{m+1})$
with multiplicity 2. This completes the case $2<j<l+1$.

We now consider the limit cases in (2), where $j=2$ or $j=l+1$ (and still with $a>1$). Here as well we consider the 
dual module $M=V^*$,
with highest weight $\lambda^*=a\lambda_1+\lambda_{n-1}$, respectively, 
$\lambda^* = a\lambda_1+\lambda_{n-l}$,
and $M^1 = S^a(2\omega_1)\otimes \omega_l$, respectively
$S^a(2\omega_1)\otimes \omega_1$.
We once again have the four summands listed in (i)- (iv) above, but where we interpret the 
subscripts $m\pm1$ accordingly. 
In case $j=l+1$ and so $m=1$, we rewrite the various summands and see that 
$M^2$ has  summands
 $S^a(2\omega_1)\otimes (\omega_2+2\omega_l)$, $S^a(2\omega_1)
\otimes (\omega_1+\omega_l)$, 
$S^{a-1}(2\omega_1)\otimes 2\omega_1$, and $S^a(2\omega_1)$.  As in (\ref{eq1}),
we have
$$M^1\otimes \omega_l = 
\bigl(S^a(2\omega_1)\otimes (\omega_1+\omega_l)\bigr)\oplus S^a(2\omega_1).$$ 
One now checks that 
$S^a(2\omega_1)\otimes (\omega_2+2\omega_l)$ is not multiplicity-free and 
 Corollary~\ref{cover} gives the desired contradiction.

The case where $j=2$ and $M$ has highest weight $a\lambda_1+\lambda_{n-1}$ is not quite as 
straightforward, though the
 arguments are similar. Again rewriting the four summands in 1- 4 given above 
(with $m=l$), we see that 
$M^2$ contains the submodule 
\[
 (S^a(2\omega_1)\otimes 2\omega_l) \oplus (S^a(2\omega_1)
\otimes \omega_{l-1}) \oplus (S^{a-1}(2\omega_1)\otimes (\omega_1+\omega_l)).
\]
 This time the contribution to 
$M^2$ from $M^1(Q_Y)$ is covered by 
$S^a(2\omega_1)\otimes 2\omega_l\oplus S^a(2\omega_1)\otimes \omega_{l-1}$.
So  Corollary~\ref{cover} gives the desired contradiction whenever 
$S^{a-1}(2\omega_1)\otimes (\omega_1+\omega_l)$ is not multiplicity-free. 
Lemmas~\ref{twolabs} and \ref{tensorprodMF}
show that this tensor product is not multiplicity-free if $a\geq 4$. So we must consider
 $a=2,3$. 
It is straightforward to see that $(2\omega_1+\omega_2)$ occurs with multiplicity two in 
$S^2(2\omega_1)\otimes (\omega_1+\omega_l)$, handling the case $a=3$. 
The case $a=2$ is then treated by applying Lemma~\ref{nonMF}(6) to the module $M$.

We now consider case (3) with $a=2,3$. Let $M = V^*$, with highest weight 
$a\lambda_1+\lambda_{r_0}$. 
Then $M^1 = V_{C^0}(a\lambda_1^0+\lambda_{r_0}^0)\downarrow L_X'$. The weights $\lambda-
\beta_{r_0}^0-\gamma_1$
and $\lambda-\beta_1^0-\cdots-\beta_{r_0}^0-\gamma_1$ afford summands of $W^2(Q_Y)$, the sum of 
which restricts to $L_X'$ as $S^a(2\omega_1)\otimes (\omega_{l-1}+2\omega_l)\otimes \omega_1$.
The latter contains $M^1\otimes \o_l \oplus (S^a(2\omega_1)\otimes (\omega_1+\omega_{l-1}+\omega_l))$.
Hence, by  Corollary~\ref{cover}, it suffices to show that $S^a(2\omega_1)\otimes
(\omega_1+\omega_{l-1}+\omega_l)$
is not multiplicity-free. This follows directly from Lemmas~\ref{twolabs} and \ref{tensorprodMF} if
 $a=3$. For $a=2$, it is an easy  check.

We now turn to case (1), where $\lambda=b\lambda_1+a\lambda_n$. 
We first consider the case where $a,b\geq 2$. Note that 
$$V_Y(b\lambda_1)\otimes V_Y(a\lambda_n)=V_Y(\lambda)\oplus 
(V_Y((b-1)\lambda_1)\otimes V_Y((a-1)\lambda_n)).$$
Using the fact that $S^d(V_X(2\omega_1))$ has summands $V_X(2d\omega_1)$ and 
$V_X((2d-4)\omega_1+2\omega_2)$, for $d\geq 2$,
it is easy to check that $S^b(2\omega_1)\otimes S^a(2\omega_{l+1})$ has three occurrences of 
$V_X((2b-2)\omega_1+(2a-2)\omega_{l+1})$. There is exactly one such summand in 
$(V_Y((b-1)\lambda_1)\otimes V_Y((a-1)\lambda_n))\downarrow X$ and so $V_Y(\lambda)$ is not
 multiplicity-free.

We now handle case (1) when either $a$ or $b$ is 1; by
 duality we may assume $a=1$. If $b\leq3$ then $\lambda$ is as in Tables \ref{TAB2}-\ref{TAB4}, so assume $b\geq4$. Here we have the isomorphism of $X$-modules
$S^b(2\omega_1)\otimes V_X(2\omega_{l+1})\cong V\downarrow X\oplus S^{b-1}(2\omega_1)$.
Now one checks that $S^b(2\omega_1)$ has summands 
$V_X((2b-4)\omega_1+2\omega_2)$,
$V_X((2b-8)\omega_1+4\omega_2)$ and $V_X((2b-6)\omega_1+2\omega_3)$. 
(Recall we are assuming that $l\geq 3$ here.) Now the Littlewood-Richardson rules (Theorem \ref{LR}) show that 
$V_X((2b-6)\omega_1+2\omega_2)$ occurs with multiplicity three in 
$S^b(2\omega_1)\otimes V_X(2\omega_{l+1})$. Since
$S^{b-1}(2\omega_1)$ is multiplicity-free by Theorem \ref{symc}, this gives the desired contradiction.

It remains to consider cases (2) and (3) when $a=1$, and case (4). That is, we have 
$\lambda = \lambda_j+\lambda_n$, where $2\leq j\leq r_0$ and $j\ne l+2$. If $j\leq 6$ or if $j\geq n-6$, then 
$\lambda^*$ is as in Tables \ref{TAB2}-\ref{TAB4},  so assume that $j\geq 7$.
For $3\le l \le 6$ we check the conclusion using Magma.
Hence we now assume that $l\geq 7$. In particular, $j\leq r_0\leq n-7$. 
So applying the inductive hypothesis to $V^*$, we deduce that one of the following holds:
\begin{enumerate}[i.]
\item $l+3\leq j\leq l+8$, or
\item $7\leq j\leq l+1$.
\end{enumerate}
As $l \ge 7$ we have $7\leq j\leq l+8$, in which case Lemma~\ref{nonMF12} applies to show that $V_Y(\lambda_j+\lambda_n)\downarrow X$ is not MF.
 This completes the proof of the proposition in Case (B).
 
\noindent{\bf{Case (C):}} Assume now that $l=2$, so $n=9$ and precisely
one of $\langle\lambda,\beta_6\rangle$, $\langle\lambda,\beta_9\rangle$ is 
nonzero. By Proposition~\ref{V1triv}, we may assume $V^1(Q_Y)$ non trivial and
 similarly for $V^*$. Set $\lambda = a\lambda_1+b\lambda_2+c\lambda_3+d\lambda_4+e\lambda_5+x\lambda_6+
w\lambda_7+z\lambda_8+y\lambda_9$.
In what follows, it is often helpful to consult Lemma \ref{A2info} to see that $V^1$ is not multiplicity free.
Applying Lemma~\ref{V1nt_1} to $V$ and $V^*$, we see that $ad=0=xy$. 

To simplify the exposition, let us also write $\Pi(Y) = \{\beta_1,\beta_2,\dots,\beta_9\}$ with 
$\Pi(L_Y) = \{\beta_1$,$\beta_2$,$\beta_3$,$\beta_4$,$\beta_5$, $\beta_7$,$\beta_8\}$.
We now define two finite sets  ${\mathcal E}$ and $\mathcal F$ of dominant weights for $Y$. 
The set ${\mathcal E}$ consists of the following weights:
\[
\begin{array}{l}
\lambda_1+\lambda_j+x\lambda_6, \,b\lambda_2+x\lambda_6, \,
 \lambda_j+\lambda_4+x\lambda_6,\, \lambda_3+x\lambda_6,\, d\lambda_4+x\lambda_6, \\
\lambda_1+\lambda_j+\lambda_6+\lambda_7,\, b\lambda_2+\lambda_6+\lambda_7,\,
\lambda_j+\lambda_4+\lambda_6+\lambda_7,\, \lambda_3+\lambda_6+\lambda_7, \\
d\lambda_4+\lambda_6+\lambda_7,  \,\lambda_5+\lambda_6,\,\lambda_r+\lambda_5+\lambda_6,\, a\lambda_1+\lambda_5+\lambda_6, \,\lambda_1+\lambda_j+\lambda_6+\lambda_8, \\
b\lambda_2+\lambda_6+\lambda_8,\,\lambda_j+\lambda_4+\lambda_6+\lambda_8,\,
\lambda_3+\lambda_6+\lambda_8, \\
d\lambda_4+\lambda_6+\lambda_8,\, \lambda_1+e\lambda_5+\lambda_6,\, e\lambda_5+\lambda_6, \lambda_4+e\lambda_5+\lambda_6, \\
a'\lambda_1+\lambda_2+x\lambda_6,\,a'\lambda_1+\lambda_2+\lambda_6+\lambda_7,\, a'\lambda_1+\lambda_2+\lambda_6+\lambda_8, \\
\end{array}
\]
where $a\leq3$, $a'=2,3$, $j=2,3$, $b\leq 3$, $d\leq 3$, $1\le r\le 4$, $1\leq x\leq 3$, and $e=2,3$.
And we define ${\mathcal F}$ to consist of the following weights:
\[
\begin{array}{l}
\lambda_1+\lambda_j+\lambda_r+\lambda_9,\, b\lambda_2+\lambda_r+\lambda_9,\,
 \lambda_3+\lambda_r+\lambda_9,a\lambda_1+\lambda_5+y\lambda_9, \\
\lambda_k+\lambda_5+y\lambda_9,\,e\lambda_5+\lambda_9,\lambda_1+e\lambda_5+\lambda_9,\,\lambda_8+y\lambda_9,\, a\lambda_1+\lambda_2+\lambda_8+y\lambda_9, \\
 \lambda_1+\lambda_j+\lambda_8+y\lambda_9,\,  \lambda_j+\lambda_8+y\lambda_9,\, b\lambda_2+\lambda_8+y\lambda_9, \,a\lambda_1+\lambda_2+\lambda_7+\lambda_9,
\end{array}
\]
where $e=2,3$, $a\leq3$, $y\leq3$, $j=2,3$, $r=7,8$, $b\leq3$, and $k=2,3,4$.

We begin by handling the case where $\l$ is in  ${\mathcal E} \cup {\mathcal F}$.
Let $1\leq x\leq 3$. The cases $\lambda=\lambda_r+\lambda_j+x\lambda_6$,
for $r\in\{1,4\}$
and $j\in\{2,3\}$, $\lambda=b\lambda_2+x\lambda_6$, $d\lambda_4+x\lambda_6$, with $b,d\leq 3$, $\lambda_3+x\lambda_6$, and  $a'\lambda_1+\lambda_2+x\lambda_6$,
$a'=2,3$, are all
treated in a similar manner. Here we may assume $b\geq 1$ and $d\geq 1$ by Proposition~\ref{V1triv}.
Using Proposition~\ref{v2gamma1}, we produce summands of $V^2(Q_Y)$ the sum of which is isomorphic to $V^1(Q_Y)\otimes V_{C^0}(\lambda_5^0) \otimes V_{C^1}(\lambda_1^1)$.
Restricting to $L_X'$ then gives
$V^1\otimes ((\omega_1+2\omega_2)\oplus \omega_2)$. 
Now we
conclude using Corollary~\ref{cover} and Lemmas \ref{twolabs} and \ref{tensorprodMF}, except when $V^1 = \wedge^3(\delta)$. But in this latter case, a Magma check shows that
  $V^1\otimes (\omega_1+2\omega_2)$ is not MF and Corollary~\ref{cover} again gives the result.

For the weights $\lambda_1+\lambda_j+\lambda_6+\lambda_7$,
$\lambda_j+\lambda_4+\lambda_6+\lambda_7$,
$b\lambda_2+\lambda_6+\lambda_7$, $\lambda_3+\lambda_6+\lambda_7$, $a'\lambda_1+\lambda_2+\lambda_6+\lambda_7$,  and $d\lambda_4+\lambda_6+\lambda_7$, with $j=2,3$, $b,d\leq 3$ and $a'=2,3$, $V^1$ is MF only if $\lambda\in\{\lambda_6+\lambda_7, \lambda_j+\lambda_6+\lambda_7,\lambda_4+\lambda_6+\lambda_7, j=2,3\}$.
(This relies on direct calculation using Lemma \ref{A2info}.) For the first case, it is a straightfoward Magma check to see that $V\downarrow X$ is not MF. For $\lambda=\lambda_2+\lambda_6+\lambda_7$, as in many of the previous cases, we compare the restriction of the $L_Y'$-summands of $V^2(Q_Y)$ afforded by $\lambda-\beta_6$, $\lambda-\beta_6-\beta_7$, $\lambda-\beta_2-\cdots-\beta_6$ and $\lambda-\beta_2-\cdots-\beta_7$ and the module 
$V^1\otimes \omega_2$, apply Corollary \ref{cover} and Lemma \ref{A2info} and obtain a contradiction. The argument is the same for $\lambda_4+\lambda_6+\lambda_7$, and indeed for the weight $\lambda_3+\lambda_6+\lambda_7$ as well.

For the weights $\lambda_1+\lambda_j+\lambda_6+\lambda_8$,
$\lambda_j+\lambda_4+\lambda_6+\lambda_8$, $j=2,3$, $b\lambda_2+\lambda_6+\lambda_8$, $d\lambda_4+\lambda_6+\lambda_8$, $b,d\leq 3$, $a'\lambda_1+\lambda_2+\lambda_6+\lambda_8$, $a'=2,3$,
   and $\lambda_3+\lambda_6+\lambda_8$, we see that $V^1$ is MF only if $\lambda\in\{\lambda_6+\lambda_8,\lambda_2+\lambda_6+\lambda_8,\lambda_4+\lambda_6+\lambda_8,\lambda_3+\lambda_6+\lambda_8\}$. 
The first case can be ruled out with a direct Magma check, and in the remaining three cases, 
Lemma \ref{v2gamma1} and Corollary \ref{cover} produce the desired contradiction.

Consider now all weights with
$\langle\lambda,\beta_5\rangle\langle\lambda,\beta_6\rangle\ne0$
and
  $\langle\lambda,\beta_j\rangle=0$ for $7\leq j\leq9$. For $\lambda=\lambda_5+\lambda_6$, a direct Magma calculation shows that $V\downarrow X$ is not MF. For the other cases, we consider the $L_Y'$-summands of $V^2(Q_Y)$ afforded by $\lambda-\beta_6$, $\lambda-\beta_5-\beta_6$ and $\lambda-\beta_i-\cdots-\beta_6$, where $i$ is minimal such that $\langle\lambda,\beta_i\rangle\ne 0$.
As usual, we compare the sum of these and the module $V^1\otimes \omega_2$ and apply  
Corollary~\ref{cover} to rule out every configuration.

At this point the weights in ${\mathcal E}$ have been dealt with.

Now consider the weights $\l_3+\eta$, $\lambda_1+\lambda_j+\eta$, $a\lambda_1+\lambda_2+\eta$, and $b\lambda_2+\eta$, where
  $j=2,3$, $a=2,3$,
$b\leq 3$ and $\eta=\lambda_8+y\lambda_9$ or $\lambda_7+y\lambda_9$, $1\leq y\leq 3$. 
As above, we see that $V^1$ is MF only if $\lambda=\lambda_i+\eta$, for $i=2,3$ or $\lambda=\eta$. The latter cases are handled with Magma. For the remaining cases, we can argue
 using Corollary~\ref{cover}, as follows. The weights $\lambda-\beta_9$ and $\lambda-\beta_8-\beta_9$, respectively $\lambda-\beta_7-\beta_8-\beta_9$, if $\eta=\lambda_8+y\lambda_9$, resp. 
$\lambda_7+y\lambda_9$, afford summands of $V^2(Q_Y)$ the sum of which has restriction giving an $L_X'$-summand of the form $V^1\otimes \omega_2$.  It is straightforward in each case to find a further non-MF summand of $V^2(Q_Y)$. 

For the weights $a\lambda_1+\lambda_5+y\lambda_9$ (with $a>0$),
$\lambda_j+\lambda_5+y\lambda_9$, and
$\lambda_1+e\lambda_5+y\lambda_9$, with $a\leq 3$, $j=2,3,4$, $e\leq 3$, we use the fact that $\lambda-\beta_9$ affords an $L_X'$-summand of $V^2$ which is isomorphic to $V^1\otimes \omega_2$ and then produce further summands in which there is a
  multiplicity and therefore rule these cases out using Corollary~\ref{cover}.

We must also
consider the weights $\lambda=e\lambda_5+y\lambda_9$, where $e\geq 1$ and one of $e$ and $y$ is equal to 1 and the other is at most 3. As usual $\lambda-\beta_9$ affords a summand of $V^2(Q_Y)$ and we analyze the remaining summands, apply Corollary~\ref{cover}, and deduce that $e=1$.
If $e=1$, induction gives that $y\leq 3$, and we may assume $y>1$ as otherwise $\lambda^*$ is in Table \ref{TAB2} of Theorem \ref{MAINTHM}. For $y=2$, we switch to the module $M=V^*$, with highest weight $\mu=2\lambda_1+\lambda_5$ and compare 
$M^1\otimes \omega_2$ and the sum of the summands of $M^2(Q_Y)$ afforded by $\mu-\beta_5-\beta_6$ and $\mu-\beta_1-\cdots-\beta_6$. We see that the latter upon restriction to $L_X'$ has the module $2\omega_1+\omega_2$ occurring with multiplicity 6, while this occurs only with multiplicity 3 in the former, contradicting  Corollary \ref{cover}. A similar argument rules out the case $y=3$.

This completes the analysis of the case where $\l$ is in the finite set ${\mathcal E} \cup {\mathcal F}$, So we assume from now on that 
$\l, \l^* \not \in{\mathcal E} \cup {\mathcal F}$.

Now suppose $x\ne 0$ and so $y=0$. Then considering $V^*$, we deduce that $(y,z,w,x,e)$ is one of
$(0,0,0,x,0), x\leq 3, 
(0,0, 1,1,0), (0,1,0,1,0)$, or $(0,0,0,1,e), 1\leq e\leq 3$. Moreover, each of the corresponding highest weights 
for $C^0$ restricts
to a weight 
for $L_X'$ having two nonzero coefficients, so using Proposition~\ref{tensorprodMF}, we further deduce that $bc=0$.
Now returning to $V$ and applying all of the above conditions, together with the assumption that 
$\l, \l^* \not \in{\mathcal E} \cup {\mathcal F}$, we deduce that $\lambda$ is one of the following:

\begin{enumerate}[(i)] 
\item $a\lambda_1+x\lambda_6$, $1\leq x\leq 3$, $a\geq 1$;
\item $a\lambda_1+\lambda_6+\lambda_7$, $a\geq 1$;
\item $a\lambda_1+\lambda_6+\lambda_8$, $a\geq 1$.
\end{enumerate}

Before considering the above infinite families, we turn to the case where $y\ne 0$ and $x=0$. 
We may assume $d=0$, else $V^*$ satisfies the conditions of the case $x\ne 0$.
Considering again $V^*$,
we deduce that $(y,z,w,x,e)$ is one of $(y,0,0,0,0), (1,0,1,0,0), (1,0,0,0,e),
1\leq e\leq3$, or
$(y,0,0,0,1), (y,1,0,0,0)$, $y\leq 3$. Arguing exactly as in the case $x\ne 0$, we deduce that in all cases 
except $(y,0,0,0,0)$,
we have $bc=0$.  Now returning to $V$, and using that $V^1$ is MF, we deduce
that either $\l^*$ is as in Tables \ref{TAB2}-\ref{TAB4}, 
or $\lambda$ or $\l^*$ is  one of the following:

\begin{enumerate}[]
\item{\rm (iv)} $a\lambda_1+y\lambda_9$;
\item{\rm (v)} $\lambda_1+\lambda_j+y\lambda_9, j=2,3$, $y\geq 1$;
\item{\rm (vi)} $\lambda_2+\lambda_3+y\lambda_9$, $y\geq 1$;
\item{\rm (vii)} $b\lambda_2+y\lambda_9$, $b\leq 3$, $y\geq 1$;
\item{\rm (viii)} $\lambda_3+y\lambda_9$, $y\geq 2$;
\item{\rm (ix)} $a\lambda_1+\lambda_8+y\lambda_9$, $1\leq y\leq 3$, $a\geq 1$;
\item{\rm (x)} $a\lambda_1+\lambda_2+y\lambda_9$, $a\leq 3$, $y\geq 1$.
\end{enumerate}

Take $\l$ to be as in (i)--(x), and let $M=V^*$. 
We shall treat cases (iv) and (viii) later. For all other cases, we produce in Table \ref{nuwts} weights $\nu_i, \eta_i$ of 
$V^2(Q_Y)$ or $M^2(Q_Y)$ such that $\sum \nu_i\downarrow L_X' = V^1 \otimes \o_2$ or $M^1\otimes \o_2$ and 
$\sum \eta_i\downarrow L_X'$ is non-MF; this gives a contradiction by Corollary \ref{cover}. (The non-MF assertions for 
$\sum \eta_i\downarrow L_X'$ are justified using Lemmas \ref{A2info}, \ref{nonMF} and \ref{twolabs}, as well as Theorems  \ref{rgreater} and \ref{a2r0}.)

\begin{table}[h]
\caption{}\label{nuwts}
\[
\begin{array}{|l|l|l|l|}
\hline
\hbox{Case} & \nu_i & \eta_i & \sum \eta_i\downarrow L_X' \supseteq \\
\hline
\hbox{(i)},\,x=2,3 & \l^*-\b_9 & \l^*-\b_4-\b_5-\b_6 & V_{C^0}(\l_3^0+(x-1)\l_4^0)\downarrow L_X' \otimes \o_1 \\
\hbox{(ii)} & \l^*-\b_9 & \l^*-\b_4-\b_5-\b_6 & V_{C^0}(2\l_3^0)\downarrow L_X' \otimes \o_1 \\
\hbox{(iii)} & \l^*-\b_9 & \l^*-\b_4-\b_5-\b_6 & V_{C^0}(\l_2^0+\l_3^0)\downarrow L_X' \otimes \o_1 \\
\hbox{(v)} & \l-\b_9 & \l-\b_j-\cdots -\b_6, & 2\o_1\otimes \wedge^{j-1}(2\o_1) \otimes \o_1 \\
     &              & \l-\b_1-\cdots -\b_6 &    \\
\hbox{(vi)} & \l-\b_9 & \l-\b_3-\cdots -\b_6 & V_{C^0}(2\l_2^0)\downarrow L_X' \otimes \o_1 \\
\hbox{(vii)},b=1 & \l^*-\b_8-\b_9, & \l^*-\b_1-\cdots -\b_6 & S^{y-1}(2\o_1) \otimes (\o_1+\o_2) \\
           \;\;\hbox{and }   y\ge 3           &  \l^*-\b_6-\b_7-\b_8, & & \\
                          & \l^*-\b_1-\cdots -\b_8 & & \\
\hbox{(ix)} & \l-\b_9, & \l-\b_6-\b_7 -\b_8, & (S^a(2\o_1) \otimes 2\o_2) \oplus \\
       &  \l-\b_8-\b_9 & \l-\b_1-\cdots -\b_6, &  (S^{a-1}(2\o_1)\otimes (\o_1+\o_2)) \\
       &                        & \l-\b_1-\cdots -\b_8 & \\
\hbox{(x)} & \l-\b_9 & \l-\b_2-\cdots -\b_6, & S^a(2\o_1) \otimes 2\o_1\otimes \o_1 \\
      &            & \l-\b_1-\cdots -\b_6 & \\
\hline
\end{array}
\]
\end{table}
A couple of special cases have been omitted in the table, namely case (i) with $x=1$, and (vii) for various small values of 
$b,y$.  In the former case, if $a=1$, $\lambda$ is in Table 2; we use Magma to exclude $a=2$,
and for $a\ge 3$ the weights $\l-\b_6$ and 
$\l-\b_1-\cdots -\b_6$ afford summands of $V^2(Q_Y)$ whose restriction to $L_X'$ contains 
$(V^1\otimes \o_2) \oplus (S^a(2\o_1)\otimes (\o_1+2\o_2))$; the latter summand is non-MF by Lemma \ref{twolabs} and Proposition \ref{tensorprodMF}. As for (vii), if $b,y\ge 2$ then $V^1$ is non-MF by Lemma \ref{nonemf}, and we 
use Magma to handle the remaining cases $y=1, b\le 3$ and $y=2,b=1$ (note that the case $y=b=1$ is an MF example in Table \ref{TAB2} of Theorem \ref{MAINTHM}).

Now consider case (viii). Here we use use following weights affording summands of $M^2(Q_Y)$: 
$\l^*-\b_7-\b_8-\b_9$, $\l^*-\b_6-\b_7$, $\l^*-\b_1-\cdots -\b_6$, $\l^*-\b_1-\cdots -\b_7$. The restriction of these to $L_X'$ contains 
\[
(M^1\otimes \o_2) \oplus (S^y(2\o_1)\otimes 3\o_2) \oplus (S^{y-1}(2\o_1)\otimes 2\o_1),
\]
and the sum of the last two summands is non-MF, giving a contradiction by Corollary \ref{cover}.

It remains to handle case (iv). Let us first assume that $y=1$ in which case we may assume $a\geq 4$ (as 
otherwise $\lambda$ is as in Table \ref{TAB2} of Theorem \ref{MAINTHM}). Then $V\downarrow X\oplus S^{a-1}(2\omega_1)$ 
is isomorphic to $S^a(2\omega_1)\otimes V_X(2\omega_3)$. Recall that $S^{a-1}(2\omega_1)$ is MF 
(by Theorem \ref{symc}), so 
it suffices to show that $S^a(2\omega_1)\otimes V_X(2\omega_3)$ has a summand of multiplicity 3. 
One first checks that $S^a(2\omega_1)$ has irreducible summands of highest weights 
$a\delta-2\alpha_1$,
$a\delta-4\alpha_1$ and $a\delta-4\alpha_1-2\alpha_2$. Then tensoring these with 
$V_X(2\omega_3)$ and 
using Proposition \ref{pieri},
we obtain three summands of highest weight $(2a-6)\omega_1+2\omega_2$. This completes the 
consideration of
the case where $y=1$. 

We now assume $y\geq 2$. If $a=1$ then $\lambda^*$ is in Tables \ref{TAB2}-\ref{TAB4}, so we 
assume as well $a\geq 2$. 
Here we have 
\[
V\downarrow X = \left(S^a(2\omega_1)\otimes S^y(2\omega_3)\right)-\left(S^{a-1}(2\omega_1)\otimes S^{y-1}(2\omega_3)\right).
\]
Now $S^a(2\omega_1)\supseteq 2a\omega_1\oplus (2a\omega_1-2\alpha_1)$ and 
$S^y(2\omega_3)\supseteq 2y\omega_3\oplus (2y\omega_3-2\a_3)$.
We count the occurrences of the summand
$\nu=2a\omega_1+2y\omega_3-2\alpha_1-2\alpha_2-2\alpha_3$:
we obtain one summand from $2a\o_1\otimes 2y\omega_3$, one from $2a\omega_1\otimes (2y\omega_3-2\alpha_3)$ and one from $(2a\omega_1-2\alpha_1)\otimes 2y\omega_3$.
Now consider $S^{a-1}(2\omega_1)\otimes S^{y-1}(2\omega_3)$. This has highest weight exactly 
$2a\omega_1+2y\omega_3-222$ with multiplicity 1, so has only one summand $\nu$. Hence we see that 
$V\downarrow X$ is not MF. \hal

\begin{lem}\label{V1nt_3} If  $\langle\lambda,\gamma_i\rangle=0$
and $\langle\lambda^*,\gamma_i\rangle=0$  for 
$i=1,2$ and $l\geq 3$, then $\lambda$ or $\lambda^*$ is as in Tables $\ref{TAB2}-\ref{TAB4}$ of Theorem $\ref{MAINTHM}$.
\end{lem}

\pf Consider first the case where $\mu^1\ne 0$. 
Applying the induction hypothesis to $V^*$, we see that 
 $$\mu^1\in\{\lambda_j^1, 2\lambda_l^1,3\lambda_l^1,\lambda_{l-1}^1+\lambda_l^1\}.$$
Also note that if $\mu^0=0$ (still with $\mu^1\ne0$), then the conclusion holds. 
Hence we assume both $\mu^1$ and $\mu^0$ are nonzero.
Applying the induction hypothesis to $V$ and $V^*$, as well as the initial assumptions on $\lambda$,
 we deduce that the support of $\mu^0$ lies in $\{\beta_2^0,\dots,\beta_{l+1}^0, \beta^0_{l+3}, \beta^0_{l+4}\}$. 
 Hence $\mu^0$ is one of the following:
\[
\begin{array}{l}
\lambda^0_k, \,2\leq k\leq l+4,\, k\ne l+2,\\
2\lambda_2^0,3\lambda_2^0, \lambda_2^0+\lambda_3^0.
\end{array}
\] 
Now apply the induction hypothesis and Lemmas~\ref{nonemf},  \ref{nonMF}, 
 \ref{twolabs} and \ref{tensorprodMF} to both $V$ and $V^*$
to reduce to $(\mu^0,\mu^1) = (\lambda_j^0,\lambda_k^1)$, where one of
\begin{enumerate}[(a)]
\item $2\leq j\leq l+1$,
\item $j=l+3$ and $l-4\leq k\leq l$,
\item $j=l+4$ and $k=l$.
\end{enumerate}
Now $V^1 = \wedge^j(2\omega_1)\otimes \omega_k$, while we have the following 
summands of $V^2$:\begin{enumerate}
\item $\wedge^{j-1}(2\omega_1)\otimes (\omega_1+\omega_k)$, afforded by 
$\lambda-\beta_j^0-\beta_{j+1}^0-\cdots-\gamma_1$,
\item $V_{C^0}(\lambda_j^0+\lambda_{r_0}^0)\downarrow L_X'\otimes ((1-\delta_{k,l})\omega_{k+1})$, 
afforded by $\lambda - \gamma_1-\beta_1^1-\cdots-\beta_k^1$,
\item $\wedge^j(2\omega_1)\otimes ((1-\delta_{k,1})\omega_{k-1})$, afforded by 
$\lambda-\beta_k^1-\beta_{k+1}^1-\cdots-\gamma_2$,
\item $\wedge^{j-1}(2\omega_1)\otimes ((1-\delta_{k,l})\omega_{k+1})$, afforded by 
$\lambda-\beta_j^0-\beta_{j+1}^0-\cdots-\gamma_1-\beta_1^1-\cdots-\beta_k^1$.
\end{enumerate}

Assume $k\ne l$. The the sum of the summands (1)-(4) contains $(V^1\otimes \o_l) \oplus Z$, where
\begin{equation}\label{V2Z}
Z = (\wedge^{j-1}(2\omega_1)\otimes (\omega_1+\omega_k))\oplus
(\wedge^j(2\omega_1)\otimes (\omega_{k+1}+2\omega_l)).
\end{equation} 
Hence it suffices to show that $Z$ is not MF. 
By Lemma~\ref{twolabs}, $\wedge^j(2\omega_1)$ has a summand with two nonzero labels and then
Proposition~\ref{tensorprodMF} implies that $\wedge^j(2\omega_1)\otimes 
(\omega_{k+1}+2\omega_l)$ is not 
multiplicity-free if $k\ne l-1$. If $k=l-1$, the same reasoning
shows that the 
summand $\wedge^{j-1}(2\omega_1)\otimes (\omega_1+\omega_k)$ is not multiplicity-free as 
long as $j\ne 2$. 
So suppose $j=2$ and $k=l-1$. Then one checks directly that the module $Z$  in (\ref{V2Z}) is not 
multiplicity-free.

Finally, suppose $k=l$; then we may assume $j\geq 3$ (else the conclusion holds). 
Here the sum of (1)-(4) above contains $(V^1\otimes \omega_l)\oplus N$, where 
$N=\wedge^{j-1}(2\omega_1)\otimes (\omega_1+\omega_l)$. But now since $j\geq 3$, 
Lemma~\ref{twolabs} shows that $N$ is non-MF, contradicting Corollary \ref{cover}.

We now turn to the case where $\mu^1 =0$. As usual, the inductive hypothesis
gives a list of possible $\mu^0$ (and hence $\lambda$). We consider each of 
these in turn and eliminate all cases where $\lambda^*$ has been handled by the above 
considerations
or where $V^*$ is inductively impossible.
We deduce that either $\lambda$ or $\lambda^*$ is as in the conclusion, or $l=3$ and 
$\lambda=\lambda_6+\lambda_9$ or $\lambda=\lambda_7+\lambda_8$. These final possibilities are handled using Magma.
\hal

\begin{lem}\label{V1nt_4} If $\langle\lambda,\gamma_i\rangle=0$ for 
$i=1,2$ and $l=2$, then $\lambda$ or $\lambda^*$  is as in Tables $\ref{TAB2}-\ref{TAB4}$ of Theorem $\ref{MAINTHM}$.
\end{lem}

\pf Throughout the proof, we will use the decompositions given in Lemma \ref{A2info} without  explicit 
reference.
Applying Lemmas~\ref{V1triv} and \ref{V1nt_2} to $V$ and $V^*$, we
 may assume $\lambda = a\lambda_2+b\lambda_3+c\lambda_5+
x\lambda_7+y\lambda_8$. 
If $(x,y)\ne (0,0)$, then $(x,y)\in\{(1,1),  (1,0), (0,y): y\leq 3\}$ 
and 
similarly for $(a,b)$. If $xy\ne 0$, then the induction hypothesis and
Proposition~\ref{tensorprodMF} implies that $c=0$ and
$V_{C^0}(a\lambda_2^0+b\lambda_3^0)\downarrow L_X'$ must 
have all summands with highest weight a multiple
of a fundamental dominant weight. Similarly, for $V_{C^0}(y\lambda_2^0+x\lambda_3^0)\downarrow L_X'$, 
if $ab\ne 0$.

We now  note that Lemma~\ref{twolabs} shows that if $(x,y) = (1,1)$ then either $a=0=b$ and so 
$\lambda^*$ is as in the conclusion, or $a=0$ and $b=1$. In the second case, setting $M=V^*$, 
it is a straightforward check to see that $M^1$ is not multiplicity-free.

Now suppose $(x,y) = (0,y)$ for some $0<y\leq 3$. Considering $V^*$, we deduce that 
$(y,c)\in\{(1,1), (1,0), (2,0), (3,0)\}$. We first show that $c=0$. If not, then $y=1$ and we see that 
$(a,b,c)\in\{(0,0,1), (1,0,1),(0,1,1)\}$. If $(a,b)\ne (0,0)$ it is 
straightforward to see that $V^1$ is not multiplicity-free. 
Hence, if $c\ne 0$, we have $\lambda = \lambda_5+\lambda_8$. But now a direct check shows 
that $V^*\downarrow X$ is not multiplicity-free. Hence if $(x,y)=(0,y)$ with $y\ne 0$, then $c=0$ as 
claimed.
Now if $(a,b)=(0,0)$ then $\lambda^*$ is as in the conclusion. So suppose 
$(a,b)\ne (0,0)$. The induction hypothesis, the preliminary remarks of the proof, and 
the first cases handled above, applied to $V$ and $V^*$, give that 
 $(a,b)\in\{(0,1),(a,0) : a\leq 3\}$. For the case $(a,b)=(0,1)$, one checks 
that $V^1$ is multiplicity-free only if 
$y=1$, that is $\lambda=\lambda_3+\lambda_8$.
Again one checks directly that $V\downarrow X$ is not multiplicity-free. So finally, we 
are left with $\lambda=a\lambda_2+y\lambda_8$,with $a,y\leq 3$. If $a=0$ or if $a=1=y$, then 
$\lambda^*$ is as in the conclusion, so without loss of 
generality we may assume $a\geq 2$. Now it is a direct check to see that 
$V^1$ is not multiplicity-free.

Suppose now that $(x,y) = (1,0)$. Consider first the case where $c\ne0$, where
by induction $(a,b)\in\{(1,0), (0,1), (0,0)\}$. In the first 
two cases, $V^1$ is not multiplicity-free. In the third case, 
$V^1 = S^c(2\omega_2)\otimes \omega_1$, while 
$V^2$ has summands 
\begin{enumerate}[(i)]
\item $V_{C^0}(\lambda_4^0+(c-1)\lambda_5^0)\downarrow L_X'\otimes 2\omega_1$, 
\item $S^{c+1}(2\omega_2)\otimes \omega_2$,
\item $S^c(2\omega_2)$, and 
\item $V_{C^0}(\lambda_4^0+(c-1)\lambda_5^0)\otimes \omega_2$.
\end{enumerate}
The sum of the summands in (ii) and (iv) gives rise to $L_X'$-summands 
$S^c(2\omega_2)\otimes 3\omega_2$ and $S^c(2\omega_2)
\otimes (\omega_1+\omega_2)$. The latter of these and the summand (iii) sum 
to $V^1\otimes \omega_2$. 
Now one checks that (i) together with $S^c(2\omega_2)\otimes 3\omega_2$ is
not MF, contradicting Corollary \ref{cover}.
Hence if $(x,y)=(1,0)$, then $c=0$. 
Moreover, using the previous cases, we reduce to $(a,b)\in\{(0,0), (0,1)\}$.
 If $(a,b)=(0,0)$,
 then $\lambda$ is as in the conclusion. In the remaining case, $\lambda=\lambda_3+\lambda_7$ and 
we can check  directly that $V_Y(\lambda)\downarrow X$ is not 
multiplicity-free.

The above considerations reduce us to the case $(x,y) = (0,0)$. By duality we
 may also assume $(a,b) = (0,0)$, so $\lambda=c\lambda_5$. If $c=1$, then 
$\lambda$ is as in the conclusion, so assume $c\geq 2$. We now apply Lemma~\ref{nonMF11}
to conclude.\hal

This completes the proof of Theorem \ref{thm:2omega1}.




\section {The case $\delta = r\omega_1$, $r\geq 3$}\label{sec:r>2}


Now assume $\delta = r\omega_1$ with $r\geq 3$. Then $L_Y' = C^0\times \cdots \times C^{r-1}$, and 
the embedding of $L_X'$ in $C^i$ is via the representation on $W^{i+1}(Q_X)$, with highest weight $(r-i)\omega_1$.
As usual we write $\Pi(Y) = \{\beta_1,\dots,\beta_n\}$ and 
$\Pi(C^i) = \{\beta_1^0,\dots,\beta_{r_i}^0\}$, for $0\leq i\leq r-1$. 

In this subsection we prove

\begin{thm}\label{thm:romega1} Let $X=A_{l+1}$ and $\delta=r\omega_1$, $r\geq 3$. Suppose $V_Y(\lambda)\downarrow X$ is 
 multiplicity-free, where $\l \ne 0,\l_1,\l_n$. Then  $\lambda$ or $\lambda^*$ is as in Tables $\ref{TAB2}-\ref{TAB4}$ of Theorem $\ref{MAINTHM}$.
\end{thm}






 \subsection{Proof of Theorem~\ref{thm:romega1}}\label{mainsec}

Let $X, Y, \d$ be as in the hypothesis, and assume $V\downarrow X$ is MF, where $V = V_Y(\l)$ and $\l\ne 0,\l_1,\l_n$. 
We assume $l\geq 2$, as the case $l=1$ has been treated in Chapter~\ref{casea2}.

\begin{lem}\label{V1_triv} Suppose $V^1(Q_Y)$ is the trivial module. Then the conclusion of Theorem $\ref{thm:romega1}$ holds.
\end{lem}

\pf We first suppose  $\langle\lambda,\gamma_i\rangle\ne 0$ for some $1\leq i\leq r-1$. Then 
$\lambda-\gamma_i$ affords a summand of 
$V^2$ isomorphic to $((r-i+1)\omega_l)\otimes ((r-i)\omega_1)$.
By Proposition~\ref{pieri}, the tensor product has summands $(\omega_1+2\omega_l)$ and
$\omega_l$. Applying  Corollary~\ref{cover}, we see that there exists at most one 
value of $i$ such that $1\leq i\leq r-1$ and $\langle\lambda,\gamma_i\rangle\ne 0$.

Consider now the case where $\langle\lambda,\gamma_{r}\rangle\ne 0$ and $\langle\lambda,\gamma_i\rangle\ne 0$ for a 
(unique) $1\leq i\leq r-1$. Consideration of the dual $\l^*$, together with the inductive hypothesis, shows that $i \ne r-1$.
Then 
$$V^2 = \omega_l \oplus \Bigl((r-i+1)\omega_l\otimes (r-i)\omega_1\Bigr).$$
Now $V^3$ has a  summand  
$\wedge^2((r-i+1)\omega_l)\otimes \wedge^2((r-i)\omega_1)$, afforded by 
$\lambda-\beta_{r_{i-1}}^{i-1}-2\gamma_i-\beta^{i}_1$.
The first tensor factor contains a summand $(\o_{l-1}+(2(r-i+1)-2)\o_l)$, and the second contains $((2(r-i)-2)\o_1+\o_2)$. By Lemma \ref{july4/17}, 
the tensor product of these has a repeated summand with $S$-value at least $4(r-i)-2$. Hence by Proposition \ref{sbd}, we have 
$4(r-i)-2 \le 2(r-i)+2$, which implies that $r-i=2$. Then 
\[
V^2 = \o_l^2 \oplus (2\o_1+3\o_l)\oplus (\o_1+2\o_l).
\]
The summand $\wedge^2(3\omega_l)\otimes \wedge^2(2\omega_1)$ of $V^3$ contains $(2\o_1+4\o_l)^2$. Also $V^3$ has a summand 
$3\o_l\otimes 2\o_1\otimes \o_l$, afforded by $\l-\g_i-\g_r$, and this also contains $(2\o_1+4\o_l)$. Hence $V^3 \supseteq (2\o_1+4\o_l)^3$. However, only one of these 
composition factors can arise from $V^2$, which is a contradiction by Corollary \ref{conseq}. We conclude that there exists a unique $i$ with $1\le i \le r$ such that 
$\langle \l,\g_i\rangle \ne 0$.

Suppose $\langle\lambda,\gamma_i\rangle\ne 0$ for some $1<i<r-1$, in particular $r\geq 4$. Then $V^2 = 
((r-i+1)\omega_l)\otimes ((r-i)\omega_1)$, and so the maximum $S$-value of an 
irreducible summand is $2r-2i+1$.
On the other hand, the weight $\lambda-\beta_{r_{i-1}}^{i-1}-2\gamma_i-\beta_1^{i}$ affords a summand $M_2$ 
(as above) of $V^3(Q_Y)$. But by Lemma~\ref{aibjnotzero}(ii), there is a repeated irreducible summand of the restriction to $L_X'$ 
which has $S$-value
$4(r-i+1)-6$. Hence $4(r-i+1)-6\leq 2(r-i)+2$ and we deduce that $i=r-2$. 

We now consider the module $M=V^*$, and we observe that 
$\langle\l^*,\beta_j\rangle\ne 0$ 
for $j= l+2+{l+2\choose 2}$. Note as well that $r_0={l+r\choose r}-1$ and $r_0-j>5$ as 
long as 
$l\geq 3$. But then $M^1$ is not MF and we have a contradiction.
Hence we reduce to the case $l=2$, still with $r\geq 4$. But now $r_0\geq14$ and $j=10$ and 
again we see 
that $M^1$ is not MF. Hence $\langle\lambda,\gamma_i\rangle= 0$ for  $1<i<r-1$.

Consider now the case where $\langle\lambda,\gamma_1\rangle\ne 0$. We  argue as above using
 $S$-values; the maximum $S$-value for $V^2$ is $2r-1$ and we have a 
repeated summand in $V^3$ with $S$-value $4r-6$. Hence $4r-6\leq 2r$ and so 
$r=3$. Repeating our analysis of $V^*$ as above, we deduce that $l\leq 3$. If $l=3$, a Magma 
check shows that $\wedge^2(003)\otimes \wedge^2(200)$ has a 
summand $103$ occurring with multiplicity 3. There is another
such summand afforded by $\lambda-\gamma_1-\beta_1^1-\cdots-\beta_{r_1}^1-\gamma_2$.
On the other hand,  $V^2 = (003) \otimes (200) = (203)+(102)+(001)$ so at most two summands $103$ arise from 
$V^2$, providing the desired contradiction.

If $\langle\lambda,\gamma_1\rangle=1$ and $l=2$, $V\downarrow X = \wedge^{10}(3\omega_1)$ and a Magma check 
shows this is not MF.  So we now
assume $\langle\lambda,\gamma_1\rangle>1$ and $l=2$. Then 
$V^2 = (2\omega_1+3\omega_2)\oplus (\omega_1+2\omega_2)
\oplus \omega_2$, while $\lambda-2\gamma_1$ affords a summand 
$S^2(3\omega_2)\otimes S^2(2\omega_1)$ of $V^3$ which has a summand 
$(2\omega_1+4\omega_2)$ occurring with 
multiplicity 4, only one of which can arise from $V^2$. This is the final contradiction in 
 case $\langle\lambda,\gamma_1\rangle\ne 0$.

Now consider the case where $\langle\lambda,\gamma_{r-1}\rangle\ne0$. Here we turn to $V^*$,
where our inductive hypothesis yields that $\langle\lambda,\gamma_{r-1}\rangle=1$ and 
$(l,r)\in\{(2,3), (2,4), (3,3)\}$.
In each of these configurations, $\lambda^*$ is as in the conclusion of Theorem \ref{thm:romega1}.

Finally, we are left to consider the case where $a:=\langle\lambda,\gamma_{r}\rangle\ne 0$. If 
$a\leq 2$,
then the conclusion holds. If $a>2$, then we apply the inductive hypothesis to $V^*$ and 
deduce that $(a,r)\in\{(3,3), (3,4), (3,5), (4,3)\}$.
In each case  the conclusion of Theorem \ref{thm:romega1} holds.

This completes the proof of the lemma.\hal

\begin{lem}\label{mu_r-2_1} If $l\geq3$, then $\mu^{r-2}=0$.
\end{lem}

\pf This will follow from a comparison of the ranks $r_0$, $r_{r-2}$ and $r_{r-1}$. Recall 
$r_0={l+r\choose r}-1$, while $r_{r-1}=l$ and $r_{r-2} = {l+2\choose 2}-1$.
One checks that for $r\geq 3$ and $l\geq5$,
\begin{equation}\label{ineq1}  r_{r-1}+r_{r-2}+2\leq r_0/2.
\end{equation}
Suppose $\mu_{r-2}\ne 0$. Then the highest weight of the dual $M=V^*$ has a nonzero coefficient
of $\lambda_j$ for some $l+3\leq j\leq r_{r-1}+2+r_{r-2}$. But by (\ref{ineq1}), this latter is at most $\frac{1}{2}r_0$ if $l\geq 5$, contradicting the inductive hypothesis. For $l=3,4$, (\ref{ineq1}) holds as long as 
$r\geq 4$. Finally, for $l=3,4$ and $r=3$, it is a direct check to see that $M^1(Q_Y)$ is not inductively 
allowed. \hal

\begin{lem}\label{gamma_r-1non0} If $\langle\lambda,\gamma_{r}\rangle\ne0$, then one of the following 
holds:
\begin{itemize}
\item[{\rm (i)}] $\mu^{r-1}=0=\mu^{r-2}$, $\langle\lambda,\gamma_{r-1}\rangle =0$ and if $l\geq 3$ 
then $\langle\lambda,\gamma_{r-2}\rangle=0$.
\item[{\rm (ii)}] $r=3$, $\langle\lambda,\gamma_{r}\rangle=1$, $\mu^{r-1}=\lambda_{r_{r-1}}^{r-1}$, 
$\mu^{r-2}=0$ and $\langle\lambda,\gamma_{r-1}\rangle=0$. Moreover, if $l\geq 3$ then 
$\langle\lambda,\gamma_{r-2}\rangle=0$.
\item[{\rm (iii)}] $r=3$, $l=2$, $\langle\lambda,\gamma_{r}\rangle=1$, $\mu^{r-1}=0$, 
$\langle\lambda,\gamma_{r-1}\rangle=0$, and $\mu^{r-2}=\lambda^{r-2}_1$.
\end{itemize}
\end{lem}

\pf Suppose first that $\mu^{r-1}=0$. By Lemma~\ref{mu_r-2_1}, either 
$\mu^{r-2}=0$ or $l=2$. If $\mu^{r-2}\ne 0$, applying the inductive hypothesis
to $V^*$ gives the conclusion of (iii).
Now suppose $\mu^{r-1}=0=\mu^{r-2}$. Again, applying the inductive hypothesis to $V^*$ gives the conclusion
 of (i).

Now consider the case where $\mu^{r-1}\ne 0$. The only possibility inductively allowed when considering 
the module $V^*$ is that $\langle\lambda,\gamma_{r}\rangle=1$ and $\mu^{r-1}=\lambda_{r_{r-1}}^{r-1}$, and 
this only for $r=3$. The remaining conditions of (ii) follow as in previous cases.\hal

\begin{lem}\label{mu_r-2_2} If $\mu^{r-2}\ne 0$, then $l=2$, $r=3$, $\mu^{m}=0$ for $m\ne r-2$, and 
$\lambda^*$ is as in Tables $\ref{TAB2}-\ref{TAB4}$ of Theorem $\ref{MAINTHM}$.
\end{lem}

\pf Assume $\mu^{r-2}\ne0$, so by Lemma~\ref{mu_r-2_1}, $l=2$. Then $\langle\lambda, 
\beta_j\rangle\ne0$ for some 
$n-8\leq j\leq n-4$ and applying the induction hypothesis to $V^*$ shows that $r=3$
and that $\mu^{r-1}=\mu^2=0$. Let us fix some
notation to be used in the rest of this proof. We now have $n=19$, $r_0=9$, $r_1=5$ and 
$r_2=2$.
We will use the notation $\Pi(Y)=\{\beta_i\ |\ 1\leq i\leq 19\}$ and 
$\Pi(C^0) = \{\beta_i\ |\ 1\leq i\leq 9\}$, $\Pi(C^1)=\{\beta_i\ |\ 11\leq i\leq 15\}$ and 
$\Pi(C^2)=\{\beta_{17},\beta_{18}\}$. Set 
$\mu^1=a\lambda_1^1+b\lambda_2^1+c\lambda_3^1+d\lambda_4^1+e\lambda_5^1$, 
$x=\langle\lambda,\beta_{16}\rangle$ and 
$y=\langle\lambda,\beta_{19}\rangle$. 

Now applying the inductive hypothesis to $V$ and $M=V^*$, we obtain the following list of possibilities for $(a,b,c,d,e,x,y)$:
\begin{enumerate}[(1)]
\item $(1,0,0,0,0,0,1)$
\item $(1,1,0,0,0,0,0)$
\item $(a,0,0,0,0,0,0)$, $a\leq4$
\item $(0,1,0,0,0,0,0)$
\item $(0,0,1,0,0,0,0)$
\item $(0,0,0,1,0,0,0)$
\item $(0,0,0,0,1,0,0)$
\end{enumerate}

\noindent{\bf {Case I:}} Suppose $\langle\l^*,\beta_{19}\rangle\ne 0$ (that is, $\langle\lambda,\beta_1\rangle\ne 0$).

Here we apply Lemma~\ref{gamma_r-1non0} to $V^*$ and find that 
\begin{enumerate}[(i)]
\item $\l^*=\nu+z\lambda_{19}$, or
\item $\l^*=\nu+\lambda_{18}+\lambda_{19}$ or 
\item $\l^*=\nu+\lambda_{11}+\lambda_{19}$,
\end{enumerate} where $\nu=y\lambda_1+x\lambda_4+e\lambda_5+d\lambda_6+c\lambda_7+b\lambda_8+
a\lambda_9+f\lambda_{10}$, for $(a,b,c,d,e,x,y)$ as above and $z\ne 0$.

The weight $\l^*-\b_{19}$ (together with $\l^*-\b_{18}-\b_{19}$ in case (ii)) affords an $L_X'$-summand of $M^2$ isomorphic to $M^1\otimes \o_l$. So by Corollary \ref{cover}, it suffices to produce an additional non-MF summand of 
$M^2$  to obtain a contradiction.  
Set $N=(2\omega_1)$, $(2\omega_1)\otimes \omega_2$, respectively 
$S^2(2\omega_1)$, according to whether $\lambda^*$ is as in (i), (ii), respectively 
(iii) above. Then for each of the configurations (1) to (7), we indicate below 
an additional $L_X'$-summand of $M^2$:
\begin{enumerate}[(1)]
\item $3\omega_1\otimes \wedge^2(3\omega_2)\otimes N$;
\item $\bigl((\wedge^2(3\omega_2)\otimes \wedge^2(3\omega_2))/\wedge^4(3\omega_2)\bigr)\otimes N$;
\item $\bigl((S^a(3\omega_2)\otimes 3\omega_2)/S^{a+1}(3\omega_2)\bigr)\otimes N$;
\item $\wedge^3(3\omega_2)\otimes N$;
\item $\wedge^4(3\omega_2)\otimes N$;
\item $\wedge^5(3\omega_2)\otimes N$;
\item $\wedge^6(3\omega_2)\otimes N$;
\end{enumerate}
Using Magma, we find that the above summand is non-MF in all cases
 except case (3), when $a=1$ and $\l^*=\nu+z\lambda_{19}$ or $\nu+\lambda_{18}+\lambda_{19}$.
 So here we have $\l^*=\lambda_9+f\lambda_{10}+z\lambda_{19}$ or 
$\l^*=\lambda_9+f\lambda_{10}+\lambda_{18}+\lambda_{19}$. If $f\ne 0$, $M^2$
has an additional summand $3\omega_2\otimes 3\omega_2\otimes N$, which is not MF.
So we have $f=0$. Now return to the module $V$, which has highest weight 
$\lambda=z\lambda_1+\lambda_{11}$ or $\lambda_1+\lambda_2+\lambda_{11}$. We see that $V^1$ is
MF only in the first case and for $z=1$. So $\lambda=\lambda_1+\lambda_{11}$. Now it is a 
Magma check to see that $V\downarrow X$ is not MF.
This completes the consideration of Case I.

\noindent{\bf {Case II:}} Suppose $\langle\l^*,\beta_{19}\rangle= 0$ (that is, $\langle\lambda,\beta_1\rangle=0$).

We first show $\mu^0=0$. Suppose not; then the inductive hypothesis gives 
$\mu^0\in\{\lambda_i^0,\,\lambda_8^0+\lambda_9^0,\,p\lambda_9^0 \,:\,2\leq i\leq 8$, 
$1\leq p\leq 4\}$. Recall from  (1) - (7) above that $\mu^1\in\{\lambda_1^1+\lambda_2^1,\,a\lambda_1^1,\,\lambda_j^1 \,:\, 1\leq a\leq 4, 2\leq j\leq 5\}$. Since $V^1$ is MF, using 
Lemma~\ref{twolabs}, Lemma~\ref{MFtensors} and Proposition~\ref{tensorprodMF}, we reduce to the following list of
 possible pairs $(\mu^0,\mu^1)$:
\begin {enumerate}[(a)]
\item $(\lambda^0_i,\lambda_1^1)$, $2\leq i\leq 8$;
\item $(\lambda_i^0,\lambda_j^1)$, $2\leq i\leq 8$, $j=3,5$;
\item $(\lambda_8^0+\lambda_9^0,a\lambda_1^1)$, $a=1,2$;
\item $(\lambda_8^0+\lambda_9^0,\lambda_j^1)$, $j=3,5$;
\item $(\lambda_9^0,\lambda_1^1+\lambda_2^1)$;
\item $(\lambda_9^0,a\lambda_1^1)$, $1\leq a\leq 4$;
\item $(\lambda_9^0,\lambda_j^1)$, $2\leq j\leq 5$.
\end{enumerate}

For the pairs in (a), (b), (c), (d) and (e), one checks (using Lemma \ref{A2info}) and Magma that $V^1$ is not MF, ruling out these configurations.
For (f), $V^1$ is MF only if $\mu^1=\lambda^1_1$
and for (g) only if $j\ne 3$. So we have $(\mu^0,\mu^1)=(\lambda_9^0,\lambda_j^1)$ for $j=1,2,4,5$. If 
$\langle\lambda,\beta_{10}\rangle\ne 0$, then 
one checks that $V^2(Q_Y)$ contains a submodule 
$V_{C^0}(\lambda_9^0)\otimes V_{C^0}(\lambda_9^0)\otimes 
V_{C^1}(\lambda_1^1)\otimes V_{C^1}(\lambda^1_j)$. In particular, $V^2$ contains a submodule
 $M = V^1\otimes \omega_2$ and it is then straightforward to 
check that $V^2/M$ is not MF, and we conclude by applying  Corollary~\ref{cover}. 
Hence in all cases, $\langle\lambda,\beta_{10}\rangle=0$.

The maximum $S$-value in $V^1$ is $5$, if $j=1$ or $5$ and $6$ if $j=2$ or $4$. In the case 
$\mu^1=\lambda^1_1$, 
combining various summands of 
$V^2$  we obtain a submodule $3\omega_2\otimes 3\omega_2\otimes 
(2\omega_1+\omega_2)$ which has a multiplicity 2 summand with $S$-value 7, contradicting 
Proposition \ref{sbd}. In each of the remaining cases, we argue similarly and 
find a repeated summand of $V^2$ whose $S$-value is larger than that 
allowed by Proposition \ref{sbd}. 
So finally, in Case II, we have shown that $\mu^0=0$. 

We now apply the induction hypothesis to $M=V^*$
and find that 
\begin{equation}\label{disp10}
\l^*\in\{\lambda_j+\lambda_9+f\lambda_{10},\; a\lambda_9+f\lambda_{10},\;
\lambda_k+f\lambda_{10}\,:\, j=1,8,\, a\leq 4,\, 5\leq k\leq 8\}.
\end{equation}
 Now if $f\ne 0$, we apply Proposition~\ref{v2gamma1} to see 
that $M^2(Q_Y)$ contains a summand $M^1(Q_Y)\otimes V_{C^0}(\lambda_{r_0}^0)\otimes V_{C_1}(\lambda_1^1)$.
The restriction of this summand to $L_X'$ then yields 
$M^1\otimes 3\omega_2\otimes 2\omega_1$. Since 
$3\omega_2\otimes 2\omega_1 = (2\omega_1+3\omega_2)\oplus (\omega_1+2\omega_2)\oplus \omega_2$,  Corollary~\ref{cover} shows that $M^1\otimes((2\omega_1+3\omega_2)\oplus (\omega_1+2\omega_2))$ must be
 MF. We now apply Lemmas~\ref{twolabs} and \ref{tensorprodMF} to reduce 
to the case $\l^* = \lambda_9+f\lambda_{10}$. Here it is a direct check to see that 
$3\omega_2\otimes ((2\omega_1+3\omega_2)\oplus (\omega_1+2\omega_2))$ is not MF. Hence we conclude that $f=0$ in all cases. 

At this point the list of possibilities in (\ref{disp10}) is bounded, and Magma computations complete the proof. \hal

\begin{lem}\label{mu_r-1_nonzero} If $\mu^{r-1}\ne 0$, then $\mu^{r-2}=0$ and one of the following holds:
\begin{itemize}
\item[{\rm (i)}] $\mu^{r-1} = \lambda_l^{r-1}$ and $\langle\lambda,\gamma_{r-1}\rangle=0$.
 Moreover, if $\langle\lambda,\gamma_{r}\rangle\ne 0$, then $r=3$ and $\langle\lambda,\gamma_{r}\rangle=1$.
\item[{\rm (ii)}] $3\leq r\leq 6$, $\mu^{r-1}=\lambda_{l-1}^{r-1}$, and $\langle\lambda,\gamma_k\rangle=0$ for 
$k=r-1,r$.
\item[{\rm (iii)}] $r=3,4$, $l\geq 3$, $\mu^{r-1}=\lambda^{r-1}_{l-2}$, and $\langle\lambda,\gamma_k\rangle=0$
for $k=r-1,r$.
\item[{\rm (iv)}] $r=3$, $l\geq4$, $\mu^{r-1}=\lambda^{r-1}_{l-3}$, and $\langle\lambda,\gamma_k\rangle=0$, 
for $k=r-1,r$.
\end{itemize}
\end{lem}

\pf That $\mu^{r-2}=0$ follows from Proposition~\ref{mu_r-2_2}. We consider the module $M=V^*$.
Then $\langle\l^*,\beta_j\rangle\ne 0$ for some $2\leq j\leq l+1<r_0/2$.
We deduce that if $r\geq 7$, then $\mu^{r-1}=\lambda^{r-1}_{l}$ and 
$\langle\lambda,\gamma_{r-1}\rangle=0$. The remaining statement in (i) also follows from 
the inductive hypothesis applied to $M^1(Q_Y)$.

For $r<7$, there are a few other possibilities allowed inductively, namely those described by parts 
(ii), (iii) and (iv). This is straightforward.\hal

\begin{lem}\label{threenonzero_mu} At most two of $\mu^i\, (0\leq i\leq r-1)$ are nonzero.
\end{lem}

\pf Suppose $\mu^i,\mu^j,\mu^k$ are all nonzero for $i,j,k$ distinct indices. Then 
Lemma~\ref{mu_r-2_2} implies that $i,j,k\ne r-2$. Then
at most one of $V_{C^m}(\mu^m)\downarrow L_X'$, $m\in\{i,j,k\}$ can have a summand with two 
nonzero labels 
(by Proposition~\ref{tensorprodMF}). Consider first the case where one of these restrictions, say 
$V_{C^i}(\mu^i)\downarrow L_X'$, 
has a summand with two nonzero labels and
 say $j<k$. Combining the results Lemmas~\ref{twolabs}, \ref{mu_r-2_2}, and 
\ref{mu_r-1_nonzero}, we see that one of the following holds:
\begin{enumerate}[(i)]
\item $\mu^m\in\{\lambda^m_1,\lambda^m_{r_m}\}$ for $m=j,k$,
\item $\mu^j\in\{\lambda^j_1,\lambda^j_{r_j}\}$ and $k=r-1$ and $\mu^{r-1} = \lambda^{r-1}_m$ for 
$m=l-3, l-2$, or $ l-1$, and $r\leq 6$.
\end{enumerate}
But now one checks that $V_{C^j}(\mu^j)\downarrow L_X'\otimes V_{C^k}(\mu^k)\downarrow L_X'$ has a summand with two nonzero labels and therefore tensoring with $V_{C^i}(\mu^i)$ yields a multiplicity in $V^1$. 

We have reduced to the case where all three restrictions $V_{C^m}(\mu^m)\downarrow L_X'$, $m\in\{i,j,k\}$ 
must have only summands with one nonzero 
label and applying Lemma \ref{twolabs}, we deduce that for $m\in\{i,j,k\}$, $V_{C^m}(\mu^m)$ is either the natural $C^m$-module or 
its dual, or $m=r-1$ and $\mu^{r-1}$ is as in (ii) above; indeed more precisely, $\mu^{r-1}$ satisfies 
condition (ii), (iii) or (iv) of Lemma~\ref{mu_r-1_nonzero}. We claim that in each case the three-fold 
tensor product $V_{C^i}(\mu^i)\downarrow L_X'\otimes V_{C^j}(\mu^j)\downarrow L_X'\otimes V_{C^k}(\mu^k)\downarrow L_X'$ is not MF.
The following assertions can be verified using the LR rules (Theorem \ref{LR}), or Lemmas~\ref{abtimesce} and \ref{aibj12}; 
this justifies the claim,
and completes the proof of the Lemma. Recall that $i,j,k\ne r-2$, and so in each of the following, 
we may assume $a,b,c\geq 1$, $a,b,c\ne 2$, and in cases (c), (d) and (e), $a,b\ne 1$:
\begin{enumerate}[(a)]
\item $a\omega_1\otimes b\omega_1\otimes c\omega_1 \supseteq ((a+b+c-2)\omega_1+\omega_2)^2$,
\item $a\omega_1\otimes b\omega_1\otimes c\omega_l \supseteq ((a+b-1)\omega_1+(c-1)\omega_l)^2$,
\item $a\omega_1\otimes b\omega_1\otimes \omega_m \supseteq ((a+b-1)\omega_1+\omega_{m+1})^2$, for  
$m=l-1,l-2,l-3$, $m>1$,
\item $a\omega_1\otimes b\omega_l\otimes \omega_m \supseteq ((a-1)\omega_1+\omega_{m}+(b-1)\omega_l)^2$, 
for $m=l-1,l-2,l-3$, $m>1$, 
\item  $a\omega_l\otimes b\omega_l\otimes \omega_m \supseteq (\omega_{m-1}+(a+b-1)\omega_l)^2$, for  
$m=l-1,l-2,l-3$, $m>1$.
\end{enumerate}\hal

\begin{lem}\label{muimuj_nonzero} Suppose $\mu^i\ne 0\ne \mu^j$ for $0\leq i<j\leq r-1$. Then 
$j=r-1$.
\end{lem}

\pf  Suppose $j<r-1$. Then Lemmas~\ref{twolabs}, \ref{MFtensors} and \ref{mu_r-2_2} 
imply that  $j<r-2$ and 
$\mu^m\in\{\lambda^m_1,\lambda^m_{r_m}\}$ for $m=i,j$. Note that $r\geq r-i>r-j\geq3$. 
Moreover,  Lemma~\ref{threenonzero_mu} 
shows that $\mu^m=0$ for $m\ne i,j$.   There are now four different cases to consider:
\begin{enumerate}[(i)]
\item $\mu^m=\lambda^m_1$, for $m=i,j$;
\item $\mu^m=\lambda^m_{r_m}$, for $m=i,j$;
\item $\mu^i=\lambda^i_1$ and $\mu^j=\lambda^j_{r_j}$; or
\item $\mu^i=\lambda^i_{r_i}$ and $\mu^j=\lambda^j_1$.
\end{enumerate}
Note that in all cases the maximum $S$-value in $V^1$ is $r-i+r-j=2r-i-j$.

In each case we will use Lemma \ref{nonMF_tensors} to produce a repeated summand of $V^2$ of $S$-value at least $2r-i+j+2$, contradicting Proposition \ref{sbd}. Table \ref{hwties} below gives the highest weights $\nu_t$ of $L_Y$-summands of $V^2(Q_Y)$, and the sum of their restrictions to $L_X'$, to which Lemma \ref{nonMF_tensors} applies. \hal 
\begin{table}[h]
\caption{}\label{hwties}
\[
\begin{array}{|l|l|l|}
\hline 
\hbox{Case} & \nu_t & \sum \nu_t\downarrow L_X' \\
\hline
\hbox{(i)},\,i\ne 0 & \l-\g_i-\b_1^i & (r-i+1)\o_l\otimes \wedge^2(r-i)\o_1 \otimes (r-j)\o_1 \\
\hbox{(i)},\,i=0,j>1 & \l-\g_j-\b_1^j & r\o_1\otimes (r-j+1)\o_l \otimes \wedge^2(r-j)\o_1 \\
\hbox{(i)},\,i=0,j=1 & \l-\g_1-\b_1^1, & r\o_1\otimes r\o_l \otimes \wedge^2(r-1)\o_1 \\
                                 & \l-\b_1^0-\cdots -\b_{r_0}^0-\g_1-\b_1^1 & \\
\hbox{(ii)} & \l-\b_{r_j}^j-\g_{j+1} & (r-i)\o_l\otimes \wedge^2(r-j)\o_l \otimes (r-j-1)\o_1 \\
\hbox{(iii)} & \l-\b_{r_j}^j-\g_{j+1} & (r-i)\o_1\otimes \wedge^2(r-j)\o_l \otimes (r-j-1)\o_1 \\
\hbox{(iv)},\,j>i+1 & \l-\b_{r_i}^i-\g_{i+1} & \wedge^2(r-i)\o_l \otimes (r-i-1)\o_1\otimes (r-j)\o_1 \\
\hbox{(iv)},\,j=i+1 & \l-\b_{r_i}^i-\g_{i+1}, & \wedge^2(r-i)\o_l \otimes (r-i-1)\o_1\otimes (r-i-1)\o_1 \\
                                & \l-\b_{r_i}^i-\g_{i+1}-\b_1^{i+1} & \\
\hline
\end{array}
\]
\end{table}

\begin{lem}\label{one_mui_nonzero} At most one $\mu^i$ ($0\le i\le r-1$)  is nonzero.
\end{lem}

\pf Assume the contrary. Then Lemma \ref{muimuj_nonzero} implies that $\mu^{r-1}\ne 0$, and so 
Lemma \ref{mu_r-1_nonzero} applies.
In particular, $\mu^{r-2}=0$ and by Lemma \ref{threenonzero_mu}, there exists a unique $i<r-1$ with $\mu^i\ne 0$.
 We must now consider 
the cases which are listed in 
Lemma~\ref{mu_r-1_nonzero}.

\noindent{\bf{Case~\ref{mu_r-1_nonzero}(i).}} Here $\mu^{r-1}=\lambda^{r-1}_l$ and note that 
$\lambda-\beta^{r-1}_l-\gamma_{r}$ and $\lambda-\gamma_{r-1}-\beta^{r-1}_1-\cdots-\beta^{r-1}_l$ 
afford 
summands of $V^2(Q_Y)$, the sum of which has restriction to $L_X'$ equal to  
$V^1\otimes \omega_l$. Hence by  Corollary~\ref{cover}, we will obtain
a contradiction if any other summand of $V^2$ is not MF. If $i>0$, then
we have summands of $V^2(Q_Y)$ afforded by $\lambda-\gamma_{i}-\eta$ and $\lambda-\chi-\gamma_{i+1}$,
where $\eta,\chi\in\Pi(C^i)$ are of minimal height such that $\lambda-\gamma_{i}-\eta$ and 
$\lambda-\chi-\gamma_{i+1}$ 
are weights of $V$. Upon restriction to $L_X'$ we obtain
summands of the form $((r-i+1)\omega_l)\otimes V_{C^i}(\mu^i-\eta+\lambda_1^i)\downarrow L_X'\otimes 
\omega_l$ and 
$((r-i-1)\omega_1)\otimes V_{C^i}(\mu^i-\chi+\lambda^i_{r_i})\downarrow L_X'\otimes \omega_l$. 
So we obtain a contradiction
if $ V_{C^i}(\mu^i-\eta+\lambda_1^i)\downarrow L_X'$ or  $V_{C^i}(\mu^i-\chi+\lambda_{r_i}^i)_{L_X'}$ has a 
summand with two 
nonzero labels. But this is easy to check essentially using the same arguments as in 
Lemma~\ref{twolabs}. 
Hence we deduce that $i=0$ and so $\mu^0\ne 0\ne \mu^{r-1}$.

We now show that  $\mu^0$ has support  among the first $l+2$ nodes. This can be seen by 
considering the dual module $M=V^*$. If $\mu^0$ has a nonzero label on one of
the nodes corresponding to some
root in the set $\{\beta^0_{l+3},\ldots ,\beta^0_{r_0}\}$,
then we apply the induction hypothesis, together with Lemmas~\ref{mu_r-2_2} and \ref{muimuj_nonzero} 
to the module $M$. In particular, we deduce that $(r,l)\ne (3,2)$,
$M^1(Q_Y) = V_{C^0}(\lambda^0_2)$ or $V_{C^0}(\lambda_1^0+\lambda_2^0)$, and
there exists $\gamma_k$, $1\leq k< r-1$ such that $\langle\l^*,\gamma_k\rangle\ne 0$. If $k>1$,
$\l^*-\gamma_k$ affords a summand of $M^2(Q_Y)$ 
whose  restriction to $L_X'$ is $M^1\otimes ((r-k+1)\omega_l)\otimes ((r-k)\omega_1)$.
This then produces a summand $M^1\otimes \omega_l$ 
as well as a summand $M^1\otimes ((r-k)\omega_1)\otimes((r-k+1)\omega_l)$.
The second tensor 
product is not MF, by Lemma~\ref{twolabs} and then  Corollary~\ref{cover} provides the
 desired contradiction. If $k=1$, then by Lemma~\ref{v2gamma1}, $M^2(Q_Y)$ contains a summand
of the form $M^1(Q_Y)\otimes V_{C^0}(\lambda^0_{r_0})\otimes V_{C^1}(\lambda^1_1)$. Now restricting this to $L_X'$, we have $M^1\otimes r\omega_l\otimes ((r-1)\omega_1)$. Decomposing
$r\omega_l\otimes ((r-1)\omega_1)$ affords summands $\omega_l$ and 
$((r-1)\omega_1+r\omega_l)$. In addition we note that $M^1$ has a summand 
with two nonzero labels and so 
applying Proposition~\ref{tensorprodMF} and  Corollary~\ref{cover} we obtain a contradiction in this case as well.

We now return to our consideration of $V$, having shown that $\mu^0$ has support among the first 
$l+2$ nodes. Moreover, $V_{C^0}(\mu^0)\downarrow L_X'\otimes \omega_l$ must be MF. So we 
apply Lemma~\ref{MFtensors} and deduce that
$\mu^0=\lambda_1^0,2\lambda_1^0$ or $\lambda_2^0$. Now $\lambda-\beta_l^{r-1}-\gamma_{r}$ and 
$\lambda-\gamma_{r-1}-\beta^1_{r-1}-\cdots-\beta_{r-1}^l$ provide summands of $V^2(Q_Y)$ the sum of which has 
restriction to $L_X'$
giving $V^1\otimes \omega_l$. Hence by Corollary \ref{cover}, any other summand of $V^2$ 
must be a multiplicity-free $L_X'$-module. In particular, $\langle\lambda,\gamma_j\rangle=0$ for all
 $j<r$.
If $\langle\lambda,\gamma_{r}\rangle= 0$, then we can appeal to Lemma~\ref{nonMF2} to conclude. If 
$\langle\lambda,\gamma_{r}\rangle\ne 0$, then by 
Lemma~\ref{mu_r-1_nonzero}, we have $r=3$ and $\langle\lambda,\gamma_{r}\rangle=1$. 
Now turn again to the module $M=V^*$, which has highest weight $\l^*=\lambda_1+\lambda_2+a\lambda_n$, $a=1,2$ 
or $\lambda_1+\lambda_2+\lambda_{n-1}$. In this case, we note that the weight $\l^* -\gamma_3$, in the first two 
cases,
or $\l^*-\beta_l^2-\gamma_3$ together with $\l^*-\gamma_2-\beta_1^2-\cdots-\beta_l^2$
in the third case,  afford summands of $M^2(Q_Y)$ whose
restriction to $L_X'$ affords $M^1\otimes \omega_l$. 
But then we also have a summand $S^2(3\omega_1)\otimes 2\omega_1\otimes \eta$, afforded by
$\l^*-\beta_2^0-\cdots-\beta^0_{r_0}-\gamma_1$, where
$\eta\in\{0,\omega_l\}$. This three-fold tensor product is not multiplicity-free,
contradicting  Corollary~\ref{cover}. This 
completes the consideration of Case~\ref{mu_r-1_nonzero}(i).

\noindent{\bf{Case~\ref{mu_r-1_nonzero}(ii).}} Here $3\leq r\leq 6$, $\mu^{r-1} = 
\lambda^{r-1}_{l-1}$ and 
$\langle\lambda,\gamma_{r}\rangle=0=\langle\lambda,\gamma_{r-1}\rangle$. 
Note that $V^2$ has summands $V_{C^i}(\mu^i)\downarrow L_X'\otimes 2\omega_l
\otimes \omega_l$ and
 $V_{C^i}(\mu^i)\downarrow L_X'\otimes \omega_{l-2}$ afforded by 
$\lambda-\gamma_{r-1}-\beta_1^{r-1}-\cdots-\beta^{r-1}_{l-1}$, respectively 
$\lambda-\beta_{l-1}^{r-1}-\beta^{r-1}_l-\gamma_{r}$. The sum of these is equal to 
$(V^1\otimes \omega_l)\oplus (V_{C^i}(\mu^i)\downarrow L_X'\otimes 3\omega_{l}).$ 
Recalling that $i<r-2$, apply
 Corollary~\ref{cover} and Lemma~\ref{MFtensors} to deduce that $\mu^i=\lambda^i_1$ or 
$\lambda^i_{r_i}$. If $\mu^i=\lambda_1^i$, we have an additional summand of $V^2(Q_Y)$ afforded
by $\lambda-\beta_1^i-\cdots-\beta^i_{r_i}-\gamma_{i+1}$; the restriction to $L_X'$ is 
$((r-i-1)\omega_1)\otimes \omega_{l-1}$. We now find that  
 $V^2/(V^1\otimes \omega_l) \supseteq ((r-i-2)\omega_1+\omega_l)^2$,   contradicting 
 Corollary~\ref{cover}. If $\mu^i=\lambda_{r_i}^i$, we have a summand of $V^2/(V^1\otimes \omega_l)$ afforded by 
$\lambda-\beta^i_{r_i}-\gamma_{i+1}$, 
whose restriction to $L_X'$ gives $\wedge^2((r-i)\omega_l)\otimes
 ((r-i-1)\omega_1)\otimes \omega_{l-1}$ which is not MF by Proposition~\ref{tensorprodMF} and Lemma~\ref{twolabs}.
This completes the case~\ref{mu_r-1_nonzero}(ii).

\noindent{\bf{Case~\ref{mu_r-1_nonzero}(iii).}} Here $r=3,4$, $l\geq 3$, 
$\mu^{r-1}=\lambda^{r-1}_{l-2}$, and $\langle\lambda,\gamma_{r}\rangle=0=\langle\lambda,\gamma_{r-1}\rangle$.
The argument is very similar to that of the previous case. Note that 
$V^2$ has summands $V_{C^i}(\mu^i)\downarrow L_X'\otimes 2\omega_l
\otimes \omega_{l-1}$ and
 $V_{C^i}(\mu^i)\downarrow L_X'\otimes \omega_{l-3}$ afforded by 
$\lambda-\gamma_{r-1}-\beta_1^{r-1}-\cdots-\beta^{r-1}_{l-2}$, respectively 
$\lambda-\beta^{r-1}_{l-2}-\beta_{l-1}^{r-1}-\beta^{r-1}_l-\gamma_{r}$. The sum of these is
 equal to $\bigl(V^1\otimes \omega_l\bigr)\oplus 
\bigl(V_{C^i}(\mu^i)\downarrow L_X'\otimes (\omega_{l-1}+2\omega_{l})\bigr).$ 
Recalling that $i<r-2$, apply
 Corollary~\ref{cover} and Lemma~\ref{twolabs} to deduce that $\mu^i=\lambda^i_1$ or 
$\lambda^i_{r_i}$. In the first case, we have an additional summand of $V^2(Q_Y)$ afforded
by $\lambda-\beta_1^i-\cdots-\beta^i_{r_i}-\gamma_{i+1}$; the restriction to $L_X'$ is 
$((r-i-1)\omega_1)\otimes \omega_{l-2}$. We now find that  
 $V^2/(V^1\otimes \omega_l) \supseteq ((r-i-2)\omega_1+\omega_{l-1})^2$, contradicting 
 Corollary~\ref{cover}. In the second case, we have a summand afforded by 
$\lambda-\beta^i_{r_i}-\gamma_{i+1}$, whose restriction to $L_X'$ gives 
$\wedge^2((r-i)\omega_l)\otimes ((r-i-1)\omega_1)\otimes \omega_{l-2}$ 
which is not MF by Proposition~\ref{tensorprodMF} and Lemma~\ref{twolabs} (and a direct calculation if 
$l=3$). This completes the consideration of Case~\ref{mu_r-1_nonzero}(iii).

\noindent{\bf{Case~\ref{mu_r-1_nonzero}(iv).}} Here $r=3$, $l\geq4$, 
$\mu^{r-1}=\lambda^{r-1}_{l-3}$ and 
$\langle\lambda,\gamma_{r-1}\rangle=0=\langle\lambda,\gamma_{r}\rangle$; in particular we see 
that $i=0$. Precisely as in the previous two cases, we deduce that $\mu^0=\lambda^0_1$ or 
$\lambda^0_{r_0}$. If $\mu^0 = \lambda^0_{r_0}$, we argue as above and easily produce a contradiction. 
If $\mu^0=\lambda^0_1$, we note that $\langle\lambda,\gamma_1\rangle=0$, else 
we contradict Corollary \ref{cover}.
So we now have $\lambda=\lambda_1+\lambda_{n-4}$. Now consider the dual module $M=V^*$.
Here $\l^*-\gamma_{r}$ affords a summand of $M^2$ of the form 
$M^1\otimes \omega_l$ and so it suffices to note that 
$\l^*-\beta^0_5-\beta^0_6-\cdots-\beta^0_{r_0}-\gamma_1$ affords an additional summand 
$\wedge^4(3\omega_1)\otimes 2\omega_1$, which by Lemma~\ref{MFtensors} is not 
MF and this gives the final contradiction.\hal

\begin{lem}\label{mu_r-1_nonzero_2} If $\mu^{r-1}\ne 0$, then the conclusion of Theorem $\ref{thm:romega1}$ holds.
\end{lem}

\pf We apply Lemma~\ref{one_mui_nonzero} and see that 
$V^1(Q_Y) = V_{C^{r-1}}(\mu^{r-1})$. We consider the cases
listed in Lemma~\ref{mu_r-1_nonzero}. In each case, if $\langle\lambda,\gamma_k\rangle=0$ for 
all $1\leq k\leq r-1$, then the result holds.
So we assume there exists $k\leq r-1$ such that $\langle\lambda,\gamma_k\rangle\ne 0$
and produce a contradiction in each case. In fact, Lemma~\ref{mu_r-1_nonzero}
shows that $k\leq r-2$. Hence $\lambda-\gamma_k$ affords a summand of $V^2$ of the 
form $((r-k+1)\omega_l)\otimes ((r-k)\omega_1)\otimes V^1$.
The decomposition of the first two tensor factors yields summands 
$(2\omega_1+3\omega_l)$, $(\omega_1+2\omega_l)$, and 
$\omega_l$, and now we conclude using Corollary \ref{cover}. \hal

\begin{lem}\label{mu0_nonzero} If $\mu^0\ne 0$, then 
the conclusion of Theorem $\ref{thm:romega1}$ holds.
\end{lem}

\pf Lemma~\ref{one_mui_nonzero} implies that $V^1(Q_Y) = V_{C^0}(\mu^0)$.
We first treat the case where the support of $\mu^0$ is entirely contained in the first $l+1$ nodes. If 
$\langle\lambda,\beta^0_j\rangle\ne 0$ for some $2\leq j\leq l+1$, the result follows from 
Lemma~\ref{mu_r-1_nonzero_2} applied to $V^*$. So we may assume $\mu^0=a\lambda^0_1$, and so 
$V^1=S^a(r\omega_1)$.
Moreover, if $\langle\lambda,\gamma_k\rangle=0$ for all $k$, the result follows from the induction 
hypothesis. So we may assume that $\langle\lambda,\gamma_k\rangle\ne 0$ for some $k$.

 We claim that $\langle\lambda,\gamma_1\rangle=0$. For otherwise Proposition~\ref{v2gamma1}
 shows that 
$V^2$ has a submodule of the form $S^a(r\omega_1)\otimes 
r\omega_l\otimes (r-1)\omega_1$. But $r\omega_l\otimes (r-1)\omega_1$ 
contains $(2\o_1+3\o_l)\oplus (\omega_1+2\omega_l)\oplus \omega_l$ and an application of 
Lemma \ref{twolabs} and Corollary \ref{cover} yields a contradiction.

Hence we now have $\langle\lambda,\gamma_k\rangle\ne 0$ for some $k>1$. Then $\lambda-\gamma_k$ 
affords a 
summand of $V^2$ of the form $M=S^a(r\omega_1)\otimes 
((r-k+1)\omega_l)\otimes ((r-k)\omega_1)$, which contains 
$S^a(r\omega_1)\otimes \omega_l$. In addition, 
$\lambda-\beta_1^0-\cdots-\beta_{r_0}^0-\gamma_1$ affords a summand 
$N=S^{a-1}(r\omega_1)\otimes (r-1)\omega_1$. Now 
$(M\oplus N)/(S^a(r\omega_1)\otimes \omega_l)$ must be MF by  Corollary~\ref{cover}. 
We then deduce that $k=r$ and $N$ must be MF.
Since $r\geq 3$, 
applying Lemma \ref{MFtensors} gives that $a\leq 2$.
 By considering $V^*$, we reduce to the following 
configurations: up to duals, $\lambda \in \{\lambda_1+\lambda_n,\,2\lambda_1+\lambda_n,\, 2\lambda_1+2\lambda_n\}$. 
The first possibility is in Table \ref{TAB2} of Theorem \ref{MAINTHM}.
 So it remains to rule out the second and third possibilities. 

If $\l = 2\l_1+\l_n$, we use the fact that
$S^2V_X(r\omega_1)\otimes V_X(r\omega_{l+1}) = V_Y(2\lambda_1+\lambda_n)\downarrow X\oplus V_X(r\omega_1)$, 
while the tensor product on the left-hand side of the equality has a summand 
$V_X((2r-2)\omega_1+(r-2)\omega_{l+1})$ with multiplicity 2. 

For $\l = 2\l_1+2\l_n$, we note that $S^2V_X(r\omega_1)\otimes S^2V_X(r\omega_{l+1})=
V_Y(2\lambda_1+2\lambda_n)\downarrow X\oplus V_X(r\omega_1)\otimes V_X(r\omega_{l+1})$. The tensor product on 
the left has a summand $V_X((2r-2)\omega_1+(2r-2)\omega_{l+1})$ with multiplicity at least 2
which does not occur in $V_X(r\omega_1)\otimes V_X(r\omega_{l+1})$. 
This completes the consideration of the case where $\mu^0$ is supported on some of the first $l+1$ nodes. 


Henceforth we will assume that $\mu^0$ has a nonzero label on a node corresponding to some root in the set
 $\{\beta_{l+2}^0,\dots,\beta_{r_0}^0\}$. 
The inductive hypothesis then gives that $\mu^0$ is one of the following: 
\[
\begin{array}{l}
a\lambda_{r_0}^0 \\
\lambda_1^0+\lambda_{r_0}^0 \\
\lambda_{r_0-i}^0, \,i=1,2,3,4, \hbox{ where } r\leq 6,4,3, \hbox{ when } i=2,3,4 \hbox{ respectively} \\
\lambda^0_{r_0-1}+\lambda^0_{r_0} \hbox{ and } r=3 \\
\end{array}
\]
If $\mu^0=\lambda^0_{r_0-1}+\lambda^0_{r_0}$ and $r=3$, then 
we consider $M=V^*$.
By Lemma~\ref{mu_r-2_2}, $l\geq 3$, so one checks that $r_1+r_2+4\leq r_0-1$ which contradicts the 
induction hypothesis as
$M^1$ is not MF.

We now treat each of the remaining cases. For each $\mu^0$, 
we will indicate below the maximum $S$-value of $V^1$, a
set of weights $\nu_i$ affording summands of $V^2(Q_Y)$ and a summand of $(\sum\nu_i)\downarrow {L_X'}$ having a
repeated summand of $S$-value exceeding the given $S$-value by at least 2, contradicting Corollary \ref{Svalues}.
We will use without reference the  LR rules, and Lemmas \ref{tensor_1}, \ref{nonMF_tensors}, \ref{mu_r-2_2}. Note as well that Lemma~\ref{mu_r-2_2} applied to $V^*$ covers the case $\mu^0=\lambda_{r_0-4}^0$ when $l=2$ and also shows that we may assume $r\geq 4$ when $\mu^0=\lambda_{r_0-i}^0$ and $l=2$.


\bigbreak
\noindent Case $\mu^0=a\lambda^0_{r_0}$, with $a\geq 2$: 
$S$-value $ar$, 
$\nu_i=\lambda-\beta_{r_0}^0-\gamma_1$, $(\sum\nu_i)\downarrow {L_X'}\supseteq
V_{C^0}(\lambda_{r_0-1}^0+(a-1)\lambda_{r_0})\downarrow L_X'\otimes (r-1)\omega_1$,
repeated summand
$(r-2)\omega_1+\omega_{l-1}+((a+1)r-3)\omega_l$\bigbreak

\noindent Case $\mu^0=\lambda_{r_0}^0$, $r\geq 4$: $S$-value $r$, $\nu_i=\lambda-\beta_{r_0}^0-\gamma_1$,
$(\sum\nu_i)\downarrow {L_X'}\supseteq\wedge^2(r\omega_l)\otimes (r-1)\omega_1$, repeated summand
$(r-3)\omega_1+\omega_{l-1}+(2r-4)\omega_l$.
\bigbreak

\noindent Case $\mu^0=\lambda^0_1+\lambda^0_{r_0}$: $S$-value $2r$, $\nu_1 = \lambda-\beta_{r_0}^0-\gamma_1$
and $\nu_2=\lambda-\beta_1^0-\cdots-\beta_{r_0}^0-\gamma_1$, $(\sum\nu_i)\downarrow {L_X'}\supseteq 
(r\omega_1)
\otimes (r-1)\omega_1\otimes \wedge^2(r\omega_l)$, repeated summand
$(2r-2)\omega_1+(2r-1)\omega_l$ 

\bigbreak
\noindent Case $\mu^0=\lambda^0_{r_0-1}$, $l\geq 3$, $r\geq 4$:  $S$-value $2r-1$, $\nu_1=\lambda-\beta_{r_0-1}^0-\beta_{r_0}^0-\gamma_1$, $(\sum\nu_i)\downarrow {L_X'}\supseteq \wedge^{3}(r\omega_l)\otimes (r-1)\omega_1$, repeated summand
 $(r-3)\omega_1+\omega_{l-1}+(3r-4)\omega_l$.\bigbreak

\noindent Case $\mu^0=\lambda_{r_0-2}^0$, $l\geq 3$, $r\geq 4$: $S$-value $3r-2$, $\nu_i=\lambda-\beta_{r_0-2}^0-\beta_{r_0-1}^0-\beta_{r_0}^0-\gamma_1$, $(\sum\nu_i)\downarrow {L_X'}\supseteq \wedge^{4}(r\omega_l)\otimes (r-1)\omega_1$,
repeated summand $(r-3)\omega_1+\omega_{l-2}+(4r-5)\omega_l$.\bigbreak

\noindent Case $\mu_0=\lambda_{r_0-i}^0$, $i=1,2,3$, $l=2$: $S$-value $2r-1$, $3r-3$, $11$, respectively (if $i=3$ then $r=4$). $\nu_i=\lambda-\beta_{r_0-i}^0-\beta_{r_0-i+1}^0-\cdots-\beta_{r_0}^0-\gamma_1$, $(\sum\nu_i)\downarrow {L_X'}\supseteq \wedge^{i+2}(r\omega_2)\otimes (r-1)\omega_1$, repeated summand $(r-2)\omega_1+(3r-4)\omega_2$, $r\omega_1+(4r-8)\omega_2$, $4\omega_1+9\omega_2$, for $i=1,2,3$ respectively.

\bigbreak
Next we handle a couple of special cases, omitted above.
 First consider $\mu^0=\lambda_{r_0-3}^0$, $r=4$. If  $l\geq 3$, then $n-r_0+1\leq r_0$ and so the $Y$-module $V^*$ 
with highest weight $\lambda^*$ 
also has a nonzero 
restriction to $C^0$. Indeed, $\langle\lambda^*,\beta_{l+{{l+2}\choose2}+3}\rangle\ne 0$.
But then the induction hypothesis yields a
 contradiction. So $l=2$. Now we apply Lemma~\ref{mu_r-2_2} to $V^*$ to conclude.

Now consider $\mu^0=\lambda_{r_0}^0$, $r=3$. If $l\geq 4$, then $n-r_0+1\leq r_0-5$ and the induction hypothesis
applied to $V^*$ gives a contradiction. If $l=2$, Lemma \ref{mu_r-2_2} applied to $V^*$ gives a contradiction.

The arguments of the proof so far have reduced us to the case where $r=3$ and $l\geq 3$, and where $l=3$ if $\mu^0=\lambda_{r_0}^0$. We now replace $V$ by $V^*$ and see that the only inductively allowed weights $\lambda^*$ occur for $l=3$. So we now have $r_0=19$, $r_1=9$ and $r_2=3$, while $n=34$, and 
$\lambda=\lambda_{19-i}+x\lambda_{20}+y\lambda_{30}+z\lambda_{34}$, where $0\leq i\leq 4$.
 Now if $y$ or $z$ is nonzero, then 
$\lambda-\beta_{30}$, respectively $\lambda-\beta_{34}$ provides a summand of $V^2$ of the form 
$V^1\otimes \omega_3$. But we also have a summand of the form 
$\wedge^{i+2}(3\omega_3)\otimes 2\omega_1$
which has a repeated summand (by Lemma~\ref{MFtensors}), contradicting  Corollary~\ref{cover}. 
Hence $y=0=z$. If $x\ne 0$, 
considering the dual module, we see that $x=1$ and $\l^*=\lambda=\lambda_{15}+\lambda_{20}$. Finally, applying as well
Lemma~\ref{V1_triv} and the previously considered cases to the dual module we reduce to the following
 possible configurations (up to duals):

\begin{enumerate}[(a)]
\item $\lambda=\lambda_{15}+\lambda_{20}$
\item $\lambda=\lambda_{17}$
\item $\lambda=\lambda_{16}$.
\end{enumerate}

\noindent In all of these cases, a Magma check gives the conclusion.
 This completes the 
proof of the lemma.\hal

\begin{lem}We have $\mu^i=0$ for $0<i<r-2$.
\end{lem}

\pf Assume $\mu^i\ne 0$ for some $0<i<r-2$.
In particular, $r\geq 4$. Moreover, by Lemma~\ref{one_mui_nonzero}, $\mu^j=0$ for all $j\ne i$.

\noindent {\bf{Case A: $0<i<r-3$.}} Here we have $r\geq 5$ and $r-i\geq 4$.
We consider successively the possibilities for $\mu^i$ arising from the inductive hypothesis.

\noindent Case: $\mu^i=\lambda_1^i+\lambda_{r_i}^i$.  

Then the maximum $S$-value of 
$V^1$ is $2(r-i)$. Now $V^2$ has a summand
$V_{C^{i-1}}(\lambda^{i-1}_{r_{i-1}})\otimes V_{C^i}(\lambda_2^i+\lambda_{r^i}^i)$, afforded by 
$\lambda-\gamma_{i}-\beta_1^i$ and a second summand 
$V_{C^{i-1}}(\lambda^{i-1}_{r_{i-1}})\otimes V_{C^i}(\lambda_1^i)$
afforded by $\lambda-\gamma_{i}-\beta_1^i-\cdots-\beta^i_{r_i}$. The sum of 
these then restricts to  give the $L_X'$-module 
$((r-i+1)\omega_l)\otimes \wedge^2((r-i)\omega_1))\otimes 
((r-i)\omega_l)$. 
Applying Lemma~\ref{nonMF_tensors} and Proposition \ref{sbd}, an $S$-value comparison implies 
that $4(r-i)-1\leq 2(r-i)+1$, a contradiction.

\noindent Case: $\mu^i=2\lambda^i_1$ or $2\lambda^i_{r_i}$.

Once again the maximum $S$-value in $V^1$ is $2(r-i)$. Here,
$V^2$ has a summand $(a\omega_l)\otimes V_{C^i}(\lambda_1^i+\lambda_2^i)\downarrow L_X'$, respectively,
 its dual, for $a=r-i+1$, 
respectively $r-i-1$. By Lemma~\ref{tensor_1}, $V_{C^i}(\lambda_1^i+\lambda_2^i)\downarrow L_X'$ has 
$L_X'$-summands $((3(r-i)-2)\omega_1+\omega_2)$ 
and $((3(r-i)-4)\omega_1+2\omega_2)$.
This then produces in the tensor product two summands 
$((3(r-i)-3)\omega_1+\omega_2+(a-1)\omega_l)$ (or the dual). Hence, we compare $S$-values and
 see that $4(r-i)-4\leq 2(r-i)+1$, giving a contradiction.

\noindent Case: $\mu^i=\lambda^i_2$ or $\lambda^i_{r_i-1}$.

The maximum $S$-value in $V^1$ is $2(r-i)-1$. 
Here we have a summand of $V^2$ of the form 
$\wedge^3((r-i)\omega_1)\otimes a\omega_l$, respectively the dual, for $a=r-i+1$, 
respectively $r-i-1$. Here one checks that there is a repeated summand
$((3(r-i)-4)\omega_1+\omega_2+(a-2)\omega_l)$, (resp. its dual), and comparing $S$-values gives that 
$4(r-i)-6\leq 2(r-i)$, again contradicting the assumptions of case A.

\noindent Case: $\mu^i=3\lambda^i_1$ or $3\lambda^i_{r_i}$.

The maximum
$S$-value in $V^1$ is $3(r-i)$. Let us treat only the first case, the second being entirely 
similar. We have a summand of $V^2$ of the form 
$V_{C^i}(2\lambda_1^i+\lambda_2^i)\downarrow L_X'\otimes ((r-i+1)\omega_l)$. Now
note that $$V_{C^i}(2\lambda^i_1+\lambda_2^i) = \wedge^2V_{C^i}(2\lambda_1^i) = 
\wedge^2(S^2(V_{C^i}(\lambda^i_1))).$$
Hence restricting to $L_X'$ we have a summand of $V^2$ of the form 
$\wedge^2(S^2((r-i)\omega_1)\otimes ((r-i+1)\omega_l)$. Now one checks
that $((4(r-i)-3)\omega_1+\omega_2+(r-i)\omega_l)$ is a repeated summand of the tensor product 
 and an $S$-value comparison gives again a 
contradiction.

The cases $\mu^i=\lambda^i_1$ or $\lambda^i_{r_i}$ are  
entirely similar, and we omit the details.

\noindent Case: $\mu^i=\lambda^i_3$ or $\lambda^i_{r_i-2}$.

This is handled in a similar manner. 
Here we use the summand of $V^2$ of the form $\wedge^4((r-i)\omega_1)
\otimes ((r-i+1)\omega_l)$ or $\wedge^4((r-i)\omega_l)\otimes ((r-i-1)\omega_1)$. 
We find that $V^2$ has a repeated summand 
$((4(r-i)-9)\omega_1+3\omega_2+(a-3)\omega_l)$, with $a=r-i+1$ or the dual with $a=r-i-1$. 
So comparing $S$-values gives that $5(r-i)-10\leq 3(r-i)-2+1$ and
so we deduce that $r-i=4$. We will treat this case together with the case $\mu^i=\lambda^i_4$ or 
$\lambda^i_{r_i-3}$, as in these cases we also have by the induction hypothesis that $r-i=4$. 

In Case A, we have reduced to $i=r-4$ and $\mu^i=\lambda^i_s$, for 
$s\in\{3,4,r_i-2,r_i-3\}$.
In particular, $r\ge 5$. Recall that $r_a = {{l+r-a}\choose {r-a}}-1$. 

First note  that if $s=3$ or 4, then
\[
\langle\lambda^*,\beta^0_j\rangle\ne 0\mbox{ for }j=r_{r-1}+r_{r-2}+r_{r-3}+r_{r-4}+5-s,
\]
while if $s=r_i-t$, $t=2,3$, then
\[
\langle\lambda^*,\beta^0_j\rangle\ne 0\mbox{ for }j=r_{r-1}+r_{r-2}+r_{r-3}+5+t.
\]
As easy inductive argument shows that for $l\ge 4$, $r\ge 5$, 
\begin{equation}\label{binommy}
r_{r-1}+r_{r-2}+r_{r-3}+r_{r-4}-2 \le r_0.
\end{equation}
It follows that $\l^*\downarrow C^0 \ne 0$, so we can apply Lemma \ref{mu0_nonzero} to the dual $V^*$ to obtain the conclusion.

To complete the consideration of Case A, we must treat the remaining possibilities 
$\mu^i\in\{\lambda_3^i,\lambda_4^i$, $\lambda_{r_i-2}^i$,$\lambda_{r_i-3}^i\}$ when $i=r-4$ and $l=2,3$. 
For $l=2$, one checks that the inequality (\ref{binommy}) holds as long as $r\geq 7$, and for 
$l=3$ as long as $r\geq 6$.  So it remains to consider the pairs $(l,r)=(2,5), (2,6), (3,5)$. 
The third case is the easiest; 
here $\langle\lambda^*,\beta_m\rangle\ne 0$ for $m$ one of $38,39,66,67$, while $r_0=55$ and $r_1=34$,
 and this contradicts the 
induction hypothesis. For the case $(l,r)=(2,6)$, by applying the inductive hypothesis to $V^*$, 
we deduce that $\mu^i=\lambda^i_4$. But then in $M=V^*$, we have that $C^1$ acts nontrivially on 
$M^1(Q_Y)$ with a highest weight which has been handled above.

So it now remains to consider the four cases arising when $(l,r)=(2,5)$ and $\mu^i=\mu^1=\lambda_3^1$,
$\lambda_4^1$, $\lambda_{r_1-2}^1$, $\lambda_{r_1-3}^1$.
The second case is again ruled out by considering the dual module. The 
dual of the third case has been treated previously. 
In the first case, the maximum $S$-value in $V^1$ is $9$, while $V^2$ has a summand
of the form $5\omega_2\otimes \wedge^4(4\omega_1)$ which produces a repeated summand with
$S$-value $13$, ruling out this case. Finally, in the fourth case we find that $V^2$ contains 
$(5\o_1+7\o_2)^2$, whereas this composition factor cannot arise from $V^1$.

This completes the consideration of case A.

\noindent{\bf{Case B:}} $i=r-3$, so in particular $r\geq 4$. The argument here is very similar to that
 used for the case $i=r-4$ 
above. We first claim that for $l\geq 2$ and $r\geq 5$, or  for $l\geq 3$ and $r\geq 4$,
we have
$\langle\lambda^*,\beta_j\rangle\ne 0$ for some $j\leq r_0$. 
Consequently, the module $V^*$ has been treated in Lemma~\ref{mu0_nonzero}, giving the conclusion.
To prove the claim, it suffices to show that 
\[
r_{r-1}+r_{r-2}+r_{r-3}+3 \le r_0.
\]
But this is easily proved by induction, under the given conditions.

Finally, to complete the consideration of Case B, we must consider the configurations where 
$(l,r)=(2,4)$. Here, we have $n=34$ and $r_0=14$ and $r_1=9$. If 
$\mu^1 = \mu^{i}\in\{\lambda^1_1+\lambda^1_{r_1}, \lambda^1_{r_1-1}+\lambda^1_{r_1}, a\lambda^1_{r_1}, \lambda^1_{r_1-k}, k=1,2,3\}$, then $\lambda^*$ has  non trivial
restriction to $C^0$ and so this case has been covered by Lemma~\ref{mu0_nonzero}. The possibilities $\mu^1=\lambda^1_1+\lambda^1_2$ and $\mu^1=a\lambda_1^1$, for $a>1$, are ruled
out by considering the dual module and seeing that these are not inductively allowed. Hence we are left with 
$\mu^1=\lambda^1_j$ for $j=1,2,3,4,5$. The maximal $S$-value in $V^1$ is, respectively, 
$3, 5, 6, 7, 7$. In each case, there is a summand of $V^2$ of the form 
$4\omega_2\otimes \wedge^{j+1}(3\omega_1)$. Now one checks using Magma that there is a repeated 
irreducible summand
with $S$-value, respectively, $5, 8, 9, 9, 10$, providing the desired contradiction. 
This completes the consideration of Case B and the proof of the lemma.\hal

\vspace{4mm}
The previous lemma, together with  Lemmas~\ref{V1_triv}, \ref{mu_r-2_2}, \ref{one_mui_nonzero}, 
\ref{mu_r-1_nonzero_2}, and \ref{mu0_nonzero},  complete the proof of Theorem~\ref{thm:romega1}.

\chapter{The case $\d = \omega_i$ with  $i\ge 3$}\label{oige3}

Let $X = A_{l+1}$ with $l\ge 2$, let $W = V_X(\d)$ and $Y = SL(W)$. Suppose $V = V_Y(\l)$ is an irreducible $Y$-module such that $V\downarrow X$ is multiplicity-free and $\l$ is not $\l_1$ or its dual.
Let $L_X' < L_Y' = C^0\times \cdots \times C^k$ as in Chapter \ref{levelset} and let $\mu^i$ be the restriction of $\l$ to $T_Y\cap C^i$, so that $V^1(Q_Y)\downarrow L_Y' = V_{C^0}(\mu^0) \otimes \cdots \otimes V_{C^k}(\mu^k)$. 

In this chapter we handle the case where $\d = \om_i$, $3\le i \le \frac{l+2}{2}$. Note that this forces $l \ge 4.$ As in Chapter \ref{rkge2}, we adopt the inductive assumption that Theorem \ref{MAINTHM} holds in general for groups of rank smaller than $l+1$.

\begin{theor}\label{omilist} Let $X=A_{l+1}$, $\d = \om_i$ with $3\le i \le \frac{l+2}{2}$, and assume that the conclusion of Theorem $\ref{MAINTHM}$ holds for groups $A_m$ of rank $m<l+1$. Then $\l,\d$ are as in Tables $\ref{TAB2}-\ref{TAB4}$ of Theorem $\ref{MAINTHM}$.
\end{theor}

It turns out that the most complicated case of the proof is when $i= \frac{l+2}{2}$. We treat this in a separate subsection, after the general case.

\section{The case where $i< \frac{l+2}{2}$}

Assume $3\le i \le \frac{l+1}{2}$, so that $l\ge 5$. We summarize some preliminary information. Let $\g_1 = \b_{r_0+1}$ be the node between $C^0$ and $C^1$.

\begin{lem}\label{prelimomi}
{\rm (i)} $k=1$.

{\rm (ii)} As $L_X'$-modules, $W^1(Q_X) \cong \om_i$ and $W^2(Q_X) \cong \om_{i-1}$.

{\rm (iii)} $C^0 \cong A_{r_0}$ and $C^1\cong A_{r_1}$, where $r_0 = {{l+1}\choose i}-1$, $r_1 = {{l+1}\choose {i-1}}-1$.

{\rm (iv)} $\la \l,\g_1\ra = 0$.
\end{lem}

\pf Parts (i), (ii) and (iii) are immediate from Theorem \ref{LEVELS}. 

Now we prove (iv). Assume first that $l\ge 6$. 
Suppose $\la \lambda, \gamma_1 \ra  \ne 0.$  We have  $C^0 = A_{r_0}$ and $C^1 = A_{r_1}$ with $r_0,r_1$ as in (iii).  Notice that $r_0 > r_1.$  Now consider $V^*.$  Then $\la (\mu^*)^0, \beta_{r_1+1}^0\ra \ne 0.$ Since $l\ge 6$, both $r_1+1 > 5$ and $r_0-(r_1+1) > 5.$  Therefore, the induction assumption in Theorem \ref{omilist} implies that $V_{C^0}((\mu^*)^0) \downarrow L_X'$ is not MF, which is a contradiction. 

Hence $l=5$, $i=3$, $r_0 = 19$ and $r_1 = 14$. 
Assume that $x: =  \la \l,\g_1\ra \ne 0$.  For the dual $V^*$, the highest weight $\l^*$ has coefficient $x$ on $\l_{15} = \l_{r_0-4}$. 
Hence from the induction assumption we must have $(\mu^*)^0 = \l_{15}^0$, so that $x=1$ and $\mu^1=0$. 

If $\mu^0 = 0$ then $\l^*  = \l_{15}$, so $V^* = \wedge^{15}(V_{A_6}(\om_3))$. A Magma computation shows that this is not MF.  Hence $\mu^0\ne 0$. 

We work with the dual $V^*$. We know that $(\mu^*)^0 = \l_{15}^0$. Suppose $(\mu^*)^1 =0$. Then $\l^* = \l = \l_{15}+\l_{20}$, 
so   Lemma \ref{v2gamma1} implies that $V^2(Q_Y)\supseteq V_{C^0}(\lambda_{15}^0)\otimes V_{C^0}(\lambda_{19}^0)\otimes
V_{C^1}(\lambda_1^1)$, whose restriction to $L_X'$ is $\wedge^5(\om_3) \otimes \om_3 \otimes \om_2$ and this has a multiplicity 2 irreducible summand not in $\wedge^5(\om_3)\otimes \om_5$, contradicting Corollary \ref{cover}. Therefore $(\mu^*)^1 \ne 0$.

We claim that $(V^*)^1$ is not MF.  A Magma computation shows that $\wedge^5(\om_3) \supseteq (20102)$.  So the claim follows from Lemmas \ref{twolabs}, \ref{tensorprodMF} and \ref{stem1.1.A}, unless $(\mu^*)^1 = \l_1^1,2\l_1^1$ or the dual of one of these.  But here use of Magma gives the claim. This is a contradiction.  \hal

Note that since $V_{C^i}(\mu^i) \downarrow L_X'$ is MF for each $i$, the induction assumption in Theorem \ref{omilist} gives us a list of possibilities for each of the weights $\mu^i$. We refer to such a list as the inductive list for $\mu^i$.

It is convenient next to deal with the smallest case, namely $i=3,l=5$.

\begin{lem}\label{i3l5}
If $i=3,l=5$ then $\l,\d$ are as in Tables $\ref{TAB2}-\ref{TAB4}$ of Theorem $\ref{MAINTHM}$.
\end{lem}

\pf Suppose $i=3,l=5$. We have $r_0=19,r_1=14$, and $\la \l,\g_1\ra = 0$ by Lemma \ref{prelimomi}. 
If $\mu^0$ and $(\mu^*)^0$ are both in $\{0, \l_1^0\}$, then $\l = \l_1, $ $\l_n$, or $\l_1+\l_n.$ The first two are out  as in the introduction to this section, and the third case is in Table \ref{TAB2} of Theorem \ref{MAINTHM}.
So we may assume that $\mu^0 \ne 0,\l_1^0$.  We have $V^1 = (V_{C^0}(\mu^0) \otimes V_{C^1}(\mu^1))\downarrow L_X'$. Hence from the inductive assumption, $\mu^0$ is one of the following: 
\[
\begin{array}{l}
\l_j^0\,(j>1),\,a\l_1^0\,(2\le a\le 5),\,a\l_{19}^0\,(2\le a\le 5), \\ \l_1^0+\l_2^0,\,\l_{18}^0+\l_{19}^0,\,\l_1^0+\l_{19}^0,\,\l_1^0+\l_{18}^0,\,\l_2^0+\l_{19}^0.
\end{array}
\]

If  $\mu^0 = a\l_{19}^0\,(a\ge 2)$ or $\l_{18}^0+\l_{19}^0$, then $(\mu^*)^0 = a\lambda_{16}^0 +\cdots$ or $\l_{16}^0+\l_{17}^0+\cdots$, neither of which is on the inductive list of possibilities. So these cases do not occur.

Now assume $\mu^0= \l_1^0+\l_{19}^0$. If $\mu^1\ne 0$ then $(\mu^*)^0$ does not belong to the inductive list. Hence $\mu^1=0$ and so $\lambda = \lambda_1 + \lambda_{19}.$  Then  $V^* = \lambda_{16} + \lambda_{34}$ and a Magma check shows that $(V^*)^1$ is not MF, a contradiction.  A similar argument shows that $\mu^0$ is not $\l_1^0+\l_{18}^0$ or  $\l_2^0+\l_{19}^0.$ 

Thus $\mu^0 = \l_j^0\,(1<j\le 19)$, $a\l_1^0\,(2\le a\le 5)$ or $\l_1^0+\l_2^0$. If $\mu^1 = 0$ then $\l = \l_j$, $a\l_1$ or $\l_1+\l_2$, and a Magma check shows that in the first case $V\downarrow X$ is MF only when $j \le 6$. Hence $\l,\d$ are as in Tables \ref{TAB2}-\ref{TAB4}.
Now assume $\mu^1\ne 0$. From the inductive list, and ruling out cases where $(\mu^*)^0$ is not on the list, we see that $\mu^1$ must be $b\l_{14}^1\,(b\ge 2)$, $\l_{13}^1+\l_{14}^1$ or $\l_i^1$ for some $i$. Also $j<19$, if $\mu^0 = \l_j^0$, since the case $j=19$ is ruled out by considering $V^*$. Lemma \ref{twolabs} shows that  $V_{C^0}(\mu^0) \downarrow L_X'$ has a composition factor with 2 nonzero labels, so $V_{C^1}(\mu^1)\downarrow L_X'$ can have no such factor. Another application of  Lemma \ref{twolabs} shows  that $\mu^1 = \l_1^1$, $\l_{14}^1$ or $2\l_{14}^1$.

At this point we know that $V_{C^0}(\mu^0) \downarrow L_X'$ is $\wedge^j(\om_3)$, $S^a(\om_3)\,(2\le a\le 5)$ or $(\om_3\otimes \wedge^2(\om_3))/\wedge^3(\om_3)$, and
$V_{C^1}(\mu^1) \downarrow L_X' = \om_2$, $\om_2^*$ or $S^2(\om_2^*)$. Now a Magma check shows that in each of these cases $ (V_{C^0}(\mu^0) \otimes V_{C^1}(\mu^1))\downarrow L_X'$ is not MF, a final contradiction. \hal

\vspace{4mm}
From now on assume that $l\ge 6$. Arguing as at the beginning of the proof of the previous lemma, we can suppose that
 $\mu^0 \ne 0,\l_1^0$. 

At this point we make a definition that will be referred to here and in later sections.

\begin{defn}\label{innerdef}
{\rm Suppose the weight $\mu^0$ is either in the inductive list given by Tables \ref{TAB1}--\ref{TAB4} of Theorem \ref{MAINTHM}, or is the dual of one of these. 
We call $\mu^0$ {\it outer} if it is one of the weights in the tables, and {\it inner} if it is the dual of one of these.
The case $\l_1^0 + \l_{r_0}^0,$ will be considered both inner and outer.} \end{defn}

\begin{lem}\label{mu0omi}
$\mu^0$ is one of the following:
\[
\begin{array}{l}
2\l_1^0, \,3\l_1^0\,(i\le 5), \,4\l_1^0\,(i=3), \,5\l_1^0\,(i=3), \\ \l_2^0, \,\l_3^0\,(i\le 6),  \,\l_3^0\,(i=7, l = 13), \,\l_4^0\,(i\le 4),\,\l_5^0\,(i=3,l\le7),\\ \l_6^0\,(i=3,l=6), \,\l_1^0+\l_2^0\,(i=3).
\end{array}
\]
\end{lem}

\pf Since $V_{C^0}(\mu^0) \downarrow L_X' $ is MF, $\mu^0$ is in the inductive list given by Tables \ref{TAB2}--\ref{TAB4}.
If $\mu^0$ is inner, then comparing ranks as in the proof of Lemma \ref{prelimomi}, we see that $(\mu^*)^0$ is not on the inductive list. Hence $\mu^0$ is outer, and these are the possibilities listed in the conclusion. \hal

\begin{lem}\label{mu1ominon0}
If $\mu^1=0$ then $\l,\d$ are as in Tables $\ref{TAB2}-\ref{TAB4}$.
\end{lem}

\pf Suppose $\mu^1=0$. Since $\la \l,\g_1\ra =0$ by Lemma \ref{prelimomi}, it then follows immediately that $\l,\d$ are in the tables, except for $\d=\om_3$, $\l = \l_5$ or $\l_6$ with $l=7$ or 6, respectively, or $\d = \om_7$ with $l = 13.$ However a Magma check shows that $\wedge^5(\om_3)$ (resp. $\wedge^6(\om_3)$, $\wedge^3(\om_7)$) is not MF for $A_8$ (resp. $A_7$, $A_{14}$), as required. \hal

\vspace{4mm}
Hence we assume now that $\mu^1\ne 0$.

\begin{lem}\label{mu1omi}
We have $\mu^1 = \l_{r_1}^1$ or $2\l_{r_1}^1$. In the latter case, $i=3$.
\end{lem}

\pf Since $V_{C^1}(\mu^1)\downarrow L_X'$ is MF, $\mu^1$ is given by the inductive list. Since $(\mu^*)^0$ must also be on the list, it follows from Lemma \ref{mu0omi} that $\mu^1$ is $a\l_{r_1}^1$ with $a\le 5$, or $\l_{r_1-a+1}^1$ with $a\le 5$, $\l_{r_1-5}^1$ ($i = 3, l = 6$) or $\l_{r_1}^1+\l_{r_1-1}^1$ with $i=3$.

The possibilities for $\mu^0$ are given by Lemma \ref{mu0omi}, and using  Lemma \ref{twolabs} we see that in each case
$V_{C^0}(\mu^0) \downarrow L_X'$ has a composition factor with 2 nonzero labels, so $V_{C^1}(\mu^1)\downarrow L_X'$ can have no such factor. Hence it follows using  Lemma \ref{twolabs} that $\mu^1$ is as in the conclusion. \hal

\vspace{4mm}
The proof of Theorem \ref{omilist} in the cases where $i \le \frac{l+1}{2}$ now follows quickly: the possibilities for $\mu^0$ and $\mu^1$ are given by Lemmas \ref{mu0omi} and \ref{mu1omi}. In every case however, Proposition \ref{newprop} shows that
$(V_{C^0}(\mu^0) \otimes V_{C^1}(\mu^1))\downarrow L_X'$ is not MF. This final contradiction completes the proof. 

\section{The case where $i=\frac{l+2}{2}$}

Suppose now that $\d = \om_i$ where $i = \frac{l+2}{2}$ and $l \ge 4$ is even. As in Lemma \ref{prelimomi}, $k=1$,
$C^0 \cong C^1 \cong A_{r_0}$ where $r_0 = {l+1 \choose i} - 1$, and $W^1(Q_X) \cong \om_i$, $W^2(Q_X) \cong \om_{i-1} \cong \om_i^*$. Again let $\g_1$ be the node between $C^0$ and $C^1$.

\begin{lem}\label{no1}
If $\mu^0=\mu^1=0$, then $i=3$, $l=4$ and $\l = \l_{10}$, as in Table $\ref{TAB4}$.
\end{lem}

\pf  Suppose $\mu^0=\mu^1=0$, and let $x = \la \l,\g_1\ra$ (so that $\l = x\l_{r_0+1}$). Now $V^1(Q_Y) = 0$, while in $V^2(Q_Y)$, the weight $\l - \g_1$ affords $V_{C^0}(\l_{r_0}^0) \otimes V_{C^1}(\l_1^1)$, whence \[ \begin{array}{rl}
V^2  = & \om_{i-1}\otimes \om_{i-1} \\
 = & 2\om_{i-1}\oplus (\om_{i-2}+\om_i) \oplus (\om_{i-3}+\om_{i+1}) \oplus \cdots \oplus (\om_1+\om_{l-1}) \oplus \om_l.
\end{array}
\]
Indeed, this can be checked using Theorem \ref{LR}.  In $V^3(Q_Y)$, the weight $\l-\b_{r_0}^0-2\g_1-\b_{1}^1$ affords $V_{C^0}(\l_{r_0-1}^0) \otimes V_{C^1}(\l_2^1)$, which restricts to $L_X'$ as $\wedge^2(\om_{i-1}) \otimes \wedge^2(\om_{i-1})$. Lemma \ref{compfactor}(ii) implies that this has a multiplicity 2 summand of highest weight $ \om_{i-3} + \om_{i-2} + \om_i + \om_{i+1}$, where the first
	term is absent if $i = 3.$ If $i+1<l$ then by Theorem \ref{LEVELS} this summand does not appear in $\sum _{i, n_i = 0} V_i^3(Q_X) + \sum _{j, n_j = 1} V_j^2(Q_X).$ This gives a contradiction by Proposition \ref{induct}.

Hence $i+1 = l$, which implies that $i=3,l=4$. If $x>1$ then $V^3(Q_Y)$ has a further summand afforded by $\l-2\g_1$, which restricts to $L_X'$ as $S^2(\om_2) \otimes S^2(\om_2)$. This has a summand $(2\om_2 + \om_4)^3$ but just one of these appears in the next level of a summand of $V^2$. So again Proposition \ref{induct} gives a contradiction. Hence $x=1$.

At this point we have $i=3$, $l=4$ and $\l = \l_{10}$, as in the conclusion. \hal

\vspace{2mm}
The case where $l=4, i=3$ requires special treatment, and we postpone this until the end of the proof in Section \ref{i3l4}. Until then, and in view of Lemma \ref{no1}, we assume that \begin{equation}\label{munotzero}
\mu^0 \ne 0, \, i \ge 4,\hbox{ and } l \ge 6.
\end{equation}

 \subsection{The case where $\mu^1\ne 0$}

Assume in this subsection that $\mu^1 \ne 0$.

\begin{lem}\label{assum}
Replacing $V$ by $V^*$ if necessary, we may assume that $\mu^1$ is $\l_1^1$ or $\l_{r_1}^1$.
\end{lem}

\pf From the inductive list, $\mu^1$ is $\l_j^1\,(j\le 5)$, $\l_6 (l = 6,i=4)$, $c\l_1^1\,(c\le 5)$, $\l_1^1+\l_{r_1}^1$, $\l_1^1+\l_2^1\,(i=4)$, or the dual of one of these.
Using Lemma \ref{twolabs}, we see that either $V_{C^1}(\mu^1) \downarrow L_X'$ has a composition factor with at least two nonzero labels, or $\mu^1$ is as in the conclusion; by duality, the same observation applies to $\mu^0$. Since $V_{C^1}(\mu^1) \downarrow L_X'$ and
$V_{C^0}(\mu^0) \downarrow L_X'$ cannot both have such a composition factor, the result follows. \hal

\begin{lem}\label{no2}
$\mu^0$ is $\l_1^0$ or $\l_{r_0}^0$.
\end{lem}

\pf Suppose false. Then $\mu^0$ is $\l_j^0\,(2\le j\le 5)$, $\l_6^0\,(i=4)$, $c\l_1^0\,(2\le c\le 5)$, $\l_1^0+\l_2^0\,(i=4)$, $\l_1^0+\l_{r_0}^0$ or the dual of one of these, noting that certain weights are possible only for certain values of $i$.
Apart from the last case, it follows from Lemma \ref{newprop} that $V^1$ is not MF, a contradiction. Finally, when $\mu^0 = \l_1^0+\l_{r_0}^0$, we have  $V_{C^0}(\mu^0)\downarrow L_X'
= (\om_i \otimes \om_{i-1})/0$, and Theorem \ref{LR} implies that  $V^1$ has a composition factor $\om_{i-2}+\om_{i-1}+\om_{i+1}$ (respectively $\om_{i-2}+\om_{i}+\om_{i+1}$)  with multiplicity 2 if $\mu^1 = \l_1^1$ (respectively  $\mu^1 = \l_{r_0}^1$),  again a contradiction. \hal

\begin{lem}\label{no3}
We have $\mu^1 \ne \l_1^1$.
\end{lem}

\pf Suppose $\mu^1 = \l_1^1$. We know that $\mu^0 = \l_1^0$ or $\l_{r_0}^0$ by Lemma \ref{no2}.

Assume $\mu^0 = \l_{r_0}^0$. Then $V^1 = \om_{i-1} \otimes \om_{i-1}$, and this has $S$-value 2. In $V^2$, the weight $\l-\b_{r_0}^0-\g_1$ gives a summand $\wedge^2(\om_{i-1}) \otimes S^2(\om_{i-1})$, while the weight $\l-\g_1-\b_{1}^1$ gives a summand $S^2(\om_{i-1}) \otimes \wedge^2(\om_{i-1})$. Hence $V^2$ has a multiplicity 2 summand of highest weight $\om_{i-2}+2\om_{i-1}+\om_i$. But this has $S$-value 4, giving a contradiction by Proposition \ref{sbd}. 

Now assume $\mu^0 = \l_1^0$. Then $V^1 = \om_i \otimes \om_{i-1}$, again of $S$-value 2.
In $V^2(Q_Y)$, the weights $\l-\b_1-\cdots -\b_{r_0}^0-\g_1$, $\l-\g_1-\b_{1}^1$ and $\l-\b_1^0-\cdots -\b_{r_0}^0-\g_1-\b_{1}^1$ afford summands for $C^0\times C^1$ of highest weights $0 \otimes 2\l_1^1$, $(\l_1^0+\l_{r_0}^0) \otimes \l_2^1$ and
$0 \otimes \l_2^1$. The sum of these, restricted to $L_X'$, is $(\om_i \otimes \om_{i-1} \otimes \wedge^2(\om_{i-1})) \oplus S^2(\om_{i-1})$, and this contains a multiplicity 2 summand of highest weight $2\om_{i-2}+\om_i+\om_{i+1}$. This has $S$-value 4, so we again have a contradiction by Proposition \ref{sbd}. \hal

\begin{lem}\label{no5}
Suppose $\mu^1 = \l_{r_1}^1$. Then $\l = \l_1+\l_n$, as in Table $\ref{TAB2}$.
\end{lem}

\pf If $\mu^0 = \l_{r_0}^0$, then in $V^*$ we have $(\mu^*)^1 = \l_1^1$, contrary to the previous lemma. Hence $\mu^0 = \l_1^0$. Then
$V^1 = \om_i\otimes \om_i$, which has $S$-value 2.

Suppose $\la \l,\g_1\ra \ne 0$. Then $\l-\g_1$ affords the summand $V_{C^0}(\l_1^0+\l_{r_0}^0) \otimes V_{C^1}(\l_1^1+\l_{r_0}^1)$, and the restriction of this to $L_X'$ contains $(\om_{i-1}+\om_i) \otimes (\om_{i-1}+\om_i)$, which has a multiplicity 2 summand of highest weight $\om_{i-2}+\om_{i-1}+\om_i+\om_{i+1}$ (see Lemma \ref{compfactor}(ii)). This has $S$-value 4, giving a contradiction by Proposition \ref{sbd}.

Hence  $\la \l,\g_1\ra = 0$, which means that $\l = \l_1+\l_n$, as in the conclusion. \hal

\vspace{2mm}
This completes the proof of Theorem \ref{omilist} in the case where $\mu^1\ne 0$ (assuming that $i = \frac{l+2}{2}$ and, as in (\ref{munotzero}), that $\mu^0\ne 0$, $l\ge 6$).

 \subsection{The case where $\mu^1=0$}

Suppose now that $\mu^1=0$, and continue to assume that $i=\frac{l+2}{2}$ and, as in (\ref{munotzero}) that $\mu^0\ne 0$ and $l\ge 6$.

\begin{lem}\label{lposs} $\mu^0$ is one of the following, or its dual:
\[
\begin{array}{l}
\l_1^0,\,2\l_1^0, \,3\l_1^0\,(i\le 6), \,4\l_1^0\,(i=4), \,5\l_1^0\,(i=4), \\ \l_2^0, \,\l_3^0\,(i\le 7), \,\l_4^0\,(i\le 5),\,\l_5^0\,(i=4),\,\l_6^0\,(i=4), \\ \l_1^0+\l_2^0\,(i=4),\,\l_1^0+\l_{r_0}^0. \\ \end{array} \] \end{lem}

\pf This is immediate from the inductive list of possibilities, noting that we need to replace $\om_i$ by $\om_{i-1} = \om_i^*$ when referring to the list, and also that $i-1 \ge 3$ as $l\ge 6$. \hal

In the next lemma, recall Definition \ref{innerdef}.

\begin{lem}\label{louter}
 $\mu^0$ is outer and also $\mu^0 \ne \l_1^0+\l_{r_0}^0$.
\end{lem}

\pf Suppose $\mu^0$ is inner, which means that it is the dual of one of the members of the list in Lemma \ref{lposs}.

First consider the case where $\mu^0 = \l_{r_0}^0$. Here $V^1 = \om_{i-1}$. The weight $\l-\b_{r_0}^0-\g_1$ gives a summand $V_{C^0}(\l_{r_0-1}^0) \otimes V_{C^1}(\l_1^1)$ of $V^2(Q_Y)$, which restricts to $L_X'$ as $\wedge^2(\om_{i-1}) \otimes \om_{i-1}$.  For $i \ge 5$ this contains $(\om_{i-4} + \om_{i-1} + \om_{i+2})$ with multiplicity 2. Indeed $\wedge^2(\om_{i-1})$ has summands $(\om_{i-2} + \om_i)$ and  $(\om_{i-4} + \om_{i+2})$, so using Theorem \ref{LR} we see that tensoring each of these with $\om_{i-1}$ produces a summand  $(\om_{i-4} + \om_{i-1} + \om_{i+2}).$ This summand has $S$-value 3, so it cannot arise from  $V^1$, a contradiction. 

Hence $i = 4$, so $l = 6$. Recall that $V^1 = \om_3$. 
 If $\la \l,\g_1\ra \ne 0$ then the weight $\l-\g_1$ gives a summand $S^2(\om_{i-1}) \otimes \om_{i-1}$ in $V^2$, and this contains $(101010)^2$. Therefore $\la \l,\g_1\ra = 0$. 
At this point we have $\l= \l_{r_0} = \l_{34}$. Now $V^2 = \wedge^2(\om_3) \otimes \om_3$, which has $S$-value 3. On the other hand the weight $\l-\b_{r_0-1}^0-2\b_{r_0}^0-2\g_1-\b_{1}^1$ gives a summand $\wedge^3(\om_3)\otimes \wedge^2(\om_3)$ in $V^3$, and this has a multiplicity 3 summand of highest weight $120101$, of $S$-value 5. This is a contradiction by Proposition \ref{sbd}.

Next assume that $\mu^0 = c\l_{r_0}^0$ with $2\le c \le 5$ and recall the restriction  on $i$ as in Lemma \ref{lposs} and its proof.
If $c\ge 3$ then $V^1 = S^c(\o_{i-1})$, of $S$-value $c$, and a Magma check shows that $V^2$ has a summand of multiplicity at least 2 with $S$-value $c+2$, contradicting Proposition \ref{sbd}.
Now suppose $c = 2.$ Here $V^1 = S^2(\om_{i-1})$. In $V^2(Q_Y)$, the weight $\l-\b_{r_0}^0-\g_1$  affords
$V_{C^0}(\l_{r_0-1}^0+\l_{r_0}^0) \otimes V_{C^1}(\l_1^1)$ This restricts to $L_X'$ as
$((\wedge^2(\om_{i-1}) \otimes \om_{i-1}) - \wedge^3(\om_{i-1})) \otimes \om_{i-1}$.  One checks that this contains $(\om_{i-3} + \om_{i-2} + \om_i + \om_{i+1})^2.$  This irreducible summand cannot arise from $S^2(\om_{i-1})$, so this is a contradiction.

Consider now $\mu^0 = \l_{r_0-j+1}^0$ with $2\le j\le 5$ or $j=6\,(i=4)$. Note that by Lemma \ref{lposs}, we have $i\le 7$ when $j\ge 3$.
Then $V^1  = \wedge^j(\om_{i-1})$, while $V^2$ contains
$\wedge^{j+1}(\om_{i-1}) \otimes \om_{i-1}$. When $j=2$, the $S$-value of $V^1$ is 2, while $V^2$ has a multiplicity 2 summand of highest weight $\om_{i-3}+\om_{i-2}+\om_i+\om_{i+1}$, of $S$-value 4.  And when $3 \le j \le 5$ (so that $i \le 7$) a Magma checks show that
$\wedge^{j+1}(\om_{i-1}) \otimes \om_{i-1}$ has a summand with multiplicity at least 2 that cannot arise from $\wedge^j(\om_{i-1})$.

This leaves the cases where $\mu^0 = \l_1^0+\l_{r_0}^0$ or $\mu^0 = \l_{r_0}^0+\l_{r_0-1}^0\,(i=4)$. In the former case,  $V^1$ is the quotient of $\om_i\otimes \om_{i-1}$ by a 1-dimensional trivial module. In $V^2(Q_Y)$, the weights $\l-\b_{r_0}^0-\g_1$ and $\l-\b_1^0-\cdots -\b_{r_0}^0-\g_1$ afford modules, the sum of which restricts to $L_X'$ as $ \om_i \otimes \wedge^2(\om_{i-1}) \otimes \om_{i-1}$ and this contains 
$\om_{i-3} + \om_{i-1} + \om_i + \om_{i+1}$ with multiplicity at least 2.   Indeed this can
be checked from Magma for the case $i = 4$, $l = 6$ and it then follows for larger values of $i$ and $l.$  This summand has $S$-value 4, whereas the $S$-value of $\om_i \otimes \om_{i-1}$ is 2.  Hence we have a contradiction. A similar argument rules out the case where
$\mu^0 = \l_{r_0}^0+\l_{r_0-1}^0$ with $i=4$.  \hal

\begin{lem}\label{linner}
$\la \l,\g_1\ra = 0$ and $\l,\d$ are as in Tables $\ref{TAB1}-\ref{TAB4}$.
\end{lem}

\pf  By the previous two lemmas, $\mu^0$ is $\l_1^0+\l_2^0\,(i=4)$,  $\l_j^0\,(2\le j\le 7)$, or $c\l_1^0\,(1\le c\le 5)$, with various restrictions on $i$ when $j\ge 3$ or $c\ge 3$. If we show that $\la \l,\g_1\ra = 0$, then $\l$ is $\l_1+\l_2\,(i=4)$,  $\l_j$, or $c\l_1$. 
The possibilities $\l_1^0+\l_2^0$ and $\l_6^0$ (both with $i=4, l=6$) are ruled out by a Magma computation, so this leaves the cases  $\l = \l_j$, or $c\l_1$ where we check that either  $\l,\d$ are in Tables $\ref{TAB1}-\ref{TAB4}$ or we have one of the following configurations:  
$l = 6$ with $\l = 4\l_1$, $5\l_1$, $\l_5$, or $\l_6;$  $l = 8$ with $\l = \l_4$; $l = 10$ with $\l = 3\l_1$. A Magma check shows that none of the exceptional configurations is MF.

Consequently  we  suppose  that $\la \l,\g_1\ra \ne 0$ and aim for a contradiction. We work through the possibilities with the aid of Lemma \ref{v2gamma1}.  Indeed, in each case the lemma shows that $V^2(Q_Y) \supseteq V_{C^0}(\mu^0) \otimes V_{C^0}(\l_{r_0}^0) \otimes V_{C^1}(\l_1^1)$.  Restricting to $L_X'$ this becomes $V^1 \otimes \om_{i-1} \otimes \om_{i-1}$.  Since $\om_{i-1} \otimes \om_{i-1} = S^2(\om_{i-1}) + \wedge^2(\om_{i-1})$ and one of these summands contains $\om_l,$ we see that  $V^2$ contains $V^1 \otimes \om_l$ together with $V^1$ tensored with the remaining summand, say $J$. 

Consider first  $\mu^0 = \l_j^0$. Here $V^1  = \wedge^j(\om_i)$ and so $V^1 \otimes J$ is not
MF, as it contains the tensor product of two irreducibles each with at least two nonzero labels. Thus we have a contradiction by Corollary \ref{cover}. The same argument applies to  $\mu^0=\l_1^0+\l_2^0\,(i=4)$ and $\mu^0 = c\l_1^0$ with $c \ge 2.$    And if $c = 1$ the second summand is $\om_i \otimes \wedge^2(\om_{i-1})$ or $\om_i \otimes S^2(\om_{i-1})$ and as $i \ge 4$, neither of these is MF. \hal

\vspace{4mm}
This completes the proof of Theorem \ref{omilist} in the case where $i = \frac{l+2}{2}$ and $l \ge 6$.

 \subsection{The case $i=3,l=4$}\label{i3l4}

Now we complete the proof of Theorem \ref{omilist} by considering the case where $i=3,l=4$. Here $X = A_5$, $Y=A_{19}$, $C^0\cong C^1 \cong A_9$, and $W^1(Q_X) \cong \om_3$, $W^2(Q_X) \cong \om_2 \cong \om_3^*$. As usual let $\g_1$ be the node between $C^0$ and $C^1$.

Using Lemma \ref{no1} and replacing $V$ by $V^*$ if necessary, we may assume that \begin{equation}\label{muassump} \mu^0\ne 0.
\end{equation}

\begin{lem}\label{i3l41}
Suppose $\mu^1\ne 0$. Then $\l=\l_1+\l_{19}$, $\l_1+\l_{18}$ or $\l_2+\l_{19}$, as in Tables $\ref{TAB2}, \ref{TAB4}$ of Theorem $\ref{MAINTHM}$. 
\end{lem}

\pf From the inductive list, each $\mu^j$ ($j=0,1$)  is one of the following, or its dual:
\begin{equation}\label{mu1lis}
\begin{array}{l}
c\l_1^j \,(c\ge 2) \\
a\l_2^j \,(2 \le a \le 5) \\
a\l_1^j + \l_2^j\, (a \ge 2)\\
a\l_1^j + \l_9^j\, (a \ge 2)\\
\l_i^j,\\
2\l_3^j,2\l_4^j,\\
\l_1^j+\l_i^j\,(i\ge 2), \\
\l_2^j+\l_3^j, \\
\l_2^j+\l_8^j,\\
\l_1^j + 2\l_2^j.
\end{array}
\end{equation}
It follows from Lemma \ref{twolabs} that $V_{C^1}(\mu^1)\downarrow L_X'$ has a composition factor with at least  two nonzero labels in all cases except for $\mu^1 = \l_1^1,\,2\l_1^1$ or the dual of one of these.
Hence, replacing $V$ by $V^*$ if necessary, we may assume that $\mu^1 = \l_1^1$, $2\l_1^1$ or the dual of one of these.

Now consider $\mu^0$. For all cases in (\ref{mu1lis}) where $\mu^0$ is not $c\l_1^0$, $a\l_2^0\,(a\ge 1)$, $ a\l_1^0 + \l_2^0\, (a \ge 2)$, $a\l_1^0+ \l_9^0\, (a \ge 2)$ or the dual of one of these, we check using Magma that $(V_{C^0}(\mu^0)\otimes V_{C^1}(\mu^1))\downarrow L_X'$ is not MF, a contradiction. 

We next rule out some configurations under the assumption that  $\mu^1=2\l_1^1$ or the dual.  If  $\mu^0 = c\l_1^0 \,(c \ge 2)$ or the dual Proposition \ref{newprop1} shows that $(V_{C^0}(\mu^0)\otimes V_{C^1}(\mu^1))\downarrow L_X'$ is not MF.  And if $\mu^0 = \l_1^0$ or $\l_2^0$ or the dual  then a Magma computation shows that $(V_{C^0}(\mu^0)\otimes V_{C^1}(\mu^1))\downarrow L_X'$ is not MF.

Suppose $\mu^0 = a\l_1^0+ \l_9^0\,(a\ge 2).$ Then $\mu^0 \downarrow L_X' = (S^a(\om_3) \otimes \om_2) - S^{a-1}(\om_3)$ and as in the proof of Lemma \ref{a4(a0...01)} we see that $S^a(\om_3) \downarrow L_X'$ contains $(00a0)$ and $(10(a-2)0)$. Tensoring with $\om_2$ we get composition factors $(10(a-1)1)$, $(11(a-2)0)$ which
are not contained in $S^{a-1}(\om_3).$   If $\mu^1 = \l_1^1$ or $ \l_9^1$, then tensoring $\mu^1 \downarrow L_X'$ with each of these composition factors produces $(01(a-1)0)^2$ or $(20(a-2)1)^2,$ respectively.
On the other hand if $\mu^1 = 2\l_1^1$ or $ 2\l_9^1$, then $\mu^1 \downarrow L_X'$ is $S^2(\om_2)$ or $S^2(\om_3)$, respectively and tensoring with $(10(a-1)1)$ produces $(10a0)^2$ or $(20(a-1)1)^2$, respectively.
Therefore $(V_{C^0}(\mu^0)\otimes V_{C^1}(\mu^1))\downarrow L_X'$ is not MF.  We have allowed for $\mu^1$ or its dual in the above, so the same holds if $\mu^0 = \l_1^0+ a\l_9^0.$ 

Next consider $\mu^0 = a\l_1^0+ \l_2^0\,(a\ge 2)$. 
Consider the embedding of $L_X'=A_4$ in $A_9$ acting on $\delta=\omega_3$.
Then we can work out the following restrictions:
If $\{\a_1,\a_2,\a_3,\a_4\}$ are the simple roots of the $A_4$ and $\{\b_i^0\ |\ 1\leq i\leq 9\}$ are the simple roots of the $A_9$, then we have $\b_1^0\downarrow T_X = \a_3$, $\b_2^0\downarrow T_X = \a_4$,
  $\b_3^0\downarrow T_X = \a_2-\a_4$,  $\b_4^0\downarrow T_X = \a_4$,
  $\b_5^0\downarrow T_X = \a_3$,  $\b_6^0\downarrow T_X = \a_1-\a_3-\a_4$,
  $\b_7^0\downarrow T_X = \a_4$,  $\b_8^0\downarrow T_X = \a_3$,  $\b_9^0\downarrow T_X = \a_2$.

Now the $T_Y$-weight $ a\l_1^0+ \l_2^0$ restricts to $A_4$ as $\nu=\omega_2+a\omega_3+\omega_4$ and affords the highest weight of a summand. We check that the weights $\nu -\a_2-\a_3$, $\nu-\a_3-\a_4,$ and $\nu - \a_3$ cannot be the highest weights of $L_X'$-summands by using the above restrictions and noticing that if $(a\l_1^0 + \l_2^0) - \sum a_j\b_j^0$ is a weight of $V_{C^0}(a\l_1^0+\l_2^0)$ and if $a_i = 0$ for some $i > 2$, then $a_j = 0$ for all $j \ge i.$

But the weight $\nu-\a_2-\a_3-\a_4$ occurs with multiplicity 4 in the summand already found while  the weights $\mu^0-\b_1^0-\b_2^0-\b_3^0-\b_4^0$,
$\mu^0-\b_2^0-\b_3^0-\b_4^0-\b_5^0$, and $\mu^0-\b_1^0-2\b_2^0-\b_3^0$ all restrict to this weight and have respective multiplicities 2, 1, and 2. Therefore the restriction also has an irreducible  summand with  highest weight $\nu-\a_2-\a_3-\a_4$ which is $\om_1 + a\omega_3$.  That is
$\mu^0 \downarrow L_X' \supseteq (01a1) + (10a0).$

As $\mu^1 =  \l_1^1$, $\l_9^1$, $2\l_1^1$, or $2\l_9^1$ one checks using Theorem \ref{LR} that $V^1$
contains $(11a0)^2$, $(11(a-1)1)^2$, $(01(a+1)0)^2$, or $(00(a+1)1)^2$, respectively, and so it is not MF.   The possibilities for $\mu^1$ considered are closed under taking duals, so we see that $V^1$ is  not MF if $\mu^0 = \l_8^0 + a\l_9^0 \,(a\ge 2).$

Now consider $\mu^0=a\l_2^0\,(a\ge 2)$. Arguing as above, we see that $\mu^0\downarrow L_X' \supseteq (0a0a)+(1(a-2)2(a-2))$. 
If $\mu^1=2\l_1^1$ or $2\l_9^1$, then $V^1 = (\mu^0\otimes \mu^1)\downarrow L_X'$ contains $(2(a-1)1(a-1))^2$ or $(1(a-1)2(a-1))^2$, a contradiction. And if $\mu^1 = \l_1^1$ or $\l_9^1$, we again check that $V^1$ is not MF: for $a=2$ this is a Magma check, and for $a\ge 3$, in the respective cases $V^1$ contains $(0(a-1)1(a-1))^2$ or $(1(a-1)1(a-1))^2$. 

Hence we conclude that
\[
\mu^1=\l_1^1 \hbox{ (or dual)}, \mu^0=\l_1^0, \l_2^0 \hbox{ (or dual)}.
\]
 
Suppose $\mu^1=\l_1^1$. If $\mu^0 = \l_9^0$ or $\l_1^0$, we obtain a contradiction exactly as in the proof of Lemma \ref{no3}. If $\mu^0 = \l_8^0$ then $V^1 = \wedge^2(\om_2)\otimes \om_2$, which has $S$-value 3. On the other hand, the weights $\l-\g_1-\b_{1}^1$ and $\l-\b_8^0-\b_9^0-\g_1-\b_{1}^1$ give summands of $V^2(Q_Y)$, the sum of whose restrictions to $L_X'$ equals $\wedge^2(\om_2)\otimes \om_2 \otimes \wedge^2(\om_2)$. This contains a multiplicity 2 summand $3011$ of $S$-value 5, contradicting Proposition \ref{sbd}. And if $\mu^0 = \l_2^0$ we get a similar contradiction: again $V^1$ has $S$-value 3, while the weights $\l-\g_1-\b_{1}^1$, $\l-\b_2^0-\cdots -\b_9^0-\g_1-\b_{1}^1$ lead to a multiplicity 2 summand $2102$ of $S$-value 5 in $V^1$.

It follows that $\mu^1 = \l_9^1$. If $\mu^0$ is $\l_9^0$ or $\l_8^0$ then $V^*$ has $(\mu^*)^1 = \l_1^1$ or $\l_2^1$.
The first possibility was ruled out in the previous paragraph, and the second is handled in the same way. Hence $\mu^0 = \l_1^0$ or $\l_2^0$. Then the $S$-value of $V^1$ is 2 or 3, respectively.
If $\la \l,\g_1\ra \ne 0$, then Proposition \ref{v2gamma1} shows that $V^2(Q_Y) \supseteq V_{C^0}(\mu^0) \otimes V_{C^0}(\l_9^0) \otimes V_{C^1}(\l_1^1+\l_9^1)$ and the restriction of this to $L_X'$ has a summand $1111$ or $1202$ of multiplicity at least 2 and $S$-value 4 or 5, respectively, contrary to Proposition \ref{sbd}.

Consequently $\la \l,\g_1 \ra = 0$. At this point we have established that $\l=\l_1+\l_{19}$ or $\l_2+\l_{19}$, as in the conclusion. \hal 

\vspace{2mm}
In view of the previous lemma we assume from now on that $\mu^1=0$. We also know that $\mu^0$ or its dual is as in the list (\ref{mu1lis}). 

\begin{lem}\label{nextstep} One of the following holds:
\begin{itemize}
\item[(i)] $\mu^0 = \l_i^0$, $c\l_1^0\,(c\ge 2)$,\,  or $c\l_9^0\,(c\ge 2)$.
\item[(ii)] $\l=\l_1+\l_2$, as in Table $\ref{TAB2}$.
\end{itemize}
\end{lem}

\pf Suppose (i) does not hold. Then $\mu^0$ or its dual is as in (\ref{mu1lis}), excluding the first  and fifth rows.
In view of Lemma \ref{00100} we can exclude the cases $\mu^0 = a\l_1^0 + \l_2^0$ and $\l_8^0 + a\l_9^0$ for $a \ge 2.$ 

Consider first the case where $\mu^0$ is inner. Then $\mu^0$ is $2\l_7^0$, $2\l_6^0$, $\l_i^0+\l_9^0\,(i\le 8)$, $\l_7^0+\l_8^0$, $2\l_8^0 + \l_9^0$,  $\l_2^0+\l_8^0$, or $a \l_8^0\,(2 \le a \le 5).$  For all but the last possibility we argue as follows:  we use Magma to compute the restriction $V^1 = V_{C^0}(\mu^0)\downarrow L_X'$.
We then use one or two weights of the form $\l-\b_i^0-\cdots -\b_9^0-\g_1$ to find summands of $V^2(Q_Y)$ in the usual way, and compute that the restriction of these to $L_X'$ contains a submodule $M \cong V^1 \otimes \om_4$, and that $V^1/M$ is not MF, which is a contradiction by Corollary \ref{cover}.  We omit the details of  these computations.

We next consider  $\mu^0 = a \l_8^0\,(2 \le a \le 5)$ which requires  more work.  We will make use of the restrictions of the $\b_1^0, \ldots , \b_9^0$ to $T_X$ that were given in the proof of Lemma \ref{i3l41}. As $\l_9^0 \downarrow L_X' = (0100)$ we see from these restrictions   that $\l_8^0 = 2\l_9^0 - \b_9^0 $ and $\l_7^0 = 3\l_9^0 - 2\b_9^0-\b_8^0 $ restrict to $(1010)$ and $(2001)$, respectively.  

Set  $b = a-1$ so that $V^2(Q_Y) =  V_{C^0}(0 \ldots  0 1b0) \otimes V_{C^1}(\l_1^1).$ Let $M$ denote the first tensor factor and $\g$ its highest weight. If $v$ is a maximal vector of $M$ then $M$ is spanned by vectors of the form $f_{\g_1}\cdots f_{\g_r}v$ where each $f_{\g_i}$ is a root vector  for a root with negative coefficient of $\b_7^0$ or $\b_8^0.$ It follows that weights of $M$ have the form $ \g - \sum c_i\b_i^0$ and $c_7 \ge c_6 \ge \cdots \ge c_1.$

We claim that $M \downarrow L_X'$ contains summands of highest weights  $\xi_1 = ((b+2)0b1)$ and $\xi_2 = (b1b1).$ This is immediate for $\xi_1$ since the above shows that $v$ affords a summand with this highest weight.  Next note that $\xi_2 = \xi_1 - \a_1$ and  this weight occurs as the restriction of $\g-\beta_5-\beta_6-\beta_7$ and  also $\g-\beta_6-\beta_7-\beta_8$ and so has multiplicity at least 3 in $M.$

If the claim is false there must be a summand of highest weight $\xi_3 $ such that  $\xi_2$ is a subdominant to $\xi_3$ and $\xi_3 \ne \xi_1, \xi_2.$ From the restrictions in Lemma \ref{i3l41} we see that writing $\xi_3$ as the restriction of $ \g - \sum c_i\b_i^0$, then $\xi_3$ has $S$-value $(2b+3)-c_7-c_4+c_3-c_2 \ge 2b+2$. Therefore  $c_7 + c_4 + c_2 \le c_3 + 1$. From the above $c_7 \ge c_4 \ge c_3$ and this implies  that $c_7 \le 1.$ At this point it is an easy check to see that the only possibility is where  $\xi_3$ is the restriction of $\nu - \b_6^0 -\b_7^0$ so that  $\xi_3 = (b0(b+2)0) = \xi_2 + \a_3$. This weight has multiplicity 1, so this establishes the claim.

Therefore   $V^2 \supseteq ((b+2)0b1) + (b1b1)) \otimes (0100).$ An application of Magma shows this tensor product contains $((b+1)1(b-1)2)^2 = (a1(a-2)2)^2$ which has  
$S$-value $2a+1.$ 
We next observe that $V^1 \subseteq S^a(1010)$ and another application of Magma shows that  $S^a(1010)$ equals $(a0a0)$ together with  a sum of irreducibles having strictly smaller  of $S$-values. As $(a0a0)$ cannot contribute a term $(a1(a-2)2)$ to $V^2$ it follows from Corollary \ref{cover}   that $V\downarrow X$ is not MF.

Now suppose $\mu^0$ is not inner, so it is $a\l_2^0\,(2 \le a \le 5)$, $\l_1^0+\l_i^0,$  $\l_2^0+\l_3^0, $ $2\l_3^0$, $2\l_4^0$, or $\l_1^0 + 2\l_2^0$.  Lemma \ref{00100}  rules out the cases  $a\l_2^0$.  
In each of the last three cases we obtain a contradiction to Corollary \ref{cover}   by considering composition factors of $V^2(Q_Y)$ of the form $V_{C^0}(\l_2^0+\l_3^0) \otimes V_{C^1}(\l_1^1),  V_{C^0}(\l_3^0+\l_4^0) \otimes V_{C^1}(\l_1^1), (V_{C^0}(2\l_1^0+\l_2^0) + V_{C^0}(2\l_2^0))  \otimes V_{C^1}(\l_1^1),$ respectively. 

So now assume that $\mu^0 = \l_1^0 + \l_i^0$ or $\l_2^0 + \l_3^0$.  First suppose that $\la \l,\g_1 \ra \ne 0$. Then Lemma \ref{v2gamma1} shows that  $V^2 \supseteq V^1 \otimes (0100) \otimes (0100)$ which equals $V^1 \otimes (0001) + (1010) + (0200),$ and we again contradict Corollary \ref{cover}. 
Hence $\la \l,\g_1 \ra = 0$, and so $\l = \l_1+\l_i\,(i\le 9)$ or $\l_2+\l_3$.  At this point we verify using Magma that the only one of these weights for which $V _{A_{19}}(\l)$ is MF on restriction to $X = A_5$ is $\l = \l_1+\l_2$, as in conclusion (ii).
\hal

\begin{lem}\label{not9}
$\mu^0$ is not $c\l_9^0\,(c\ge 2)$.
\end{lem}

\pf Suppose $\mu^0 = c\l_9^0$ with $c\ge 2$.  We have
$V^1 = S^c(\om_2)$, while $V^2(Q_Y)$ contains $V_{C^0}(\l_8^0 + (c-1)\l_9^0) \otimes V_{C^1}(\l_1^1)$. Taking duals in the discussion of the case $a\l_1^0 + \l_2^0$  in the proof of Lemma \ref{i3l41} shows that $V_{C^0}(\l_8^0 +(c-1)\l_9^0) \downarrow L_X' \supseteq (1(c-1)10) + (0(c-1)01)$, so that $V^2$ contains $(1(c-2)11)^2$.  On the other hand  $V^1 = (0c00) + (0(c-2)01)+ \cdots$, so $(1(c-1)11)^2$ cannot arise from $V^1(Q_Y),$ a contradiction.  \hal

\begin{lem}\label{final}
We have $\l=\l_i$,  $c\l_1\,(c\le 5)$ or $\l_1+\l_2$, as in Tables $\ref{TAB2}-\ref{TAB4}$ of Theorem $\ref{MAINTHM}$.
\end{lem}

\pf Assume $\l \ne \l_1+\l_2$. Then by the previous two lemmas we have $\mu^0=\l_i^0$ or  $c\l_1^0$. 

Assume $\la \l,\g_1 \ra \ne 0$. For $\mu^0=\l_i^0$, the proof of Lemma \ref{linner} shows that $V^1 
= \wedge^i(\om_3)$, while $V^2$ contains $\wedge^i(\om_3)\otimes \om_2\otimes \om_2$. Now $\om_2\otimes \om_2 \supseteq \om_4 + (\om_1 + \om_3)$ so we get the usual contradiction from Corollary \ref{cover}.   Similarly,
for $\mu^0 = c\l_1^0$, $V^2$ contains $S^c(\om_3)\otimes \om_2\otimes \om_2$, and this contains $S^c(\om_3)\otimes \om_4$ with non-MF quotient. Thus we have a contradiction by Corollary \ref{cover}.
Hence $\l=\l_i$ or $c\l_1$. In the latter case $c\le 5$, since  $V = S^c(V_X(\om_3))$ is not MF for $c\ge 6$ by Proposition \ref{newprop1}(ii).  \hal

\vspace{4mm}
This completes the proof of Theorem \ref{omilist}.


\chapter{The case $\d = \omega_2$}\label{o2sec}

Let $X = A_{l+1}$.  In this chapter we consider the case where  $\d = \om_2$.  Set $W = V_X(\d)$ and $Y = SL(W) = A_n$. Suppose $V = V_Y(\l)$ is an irreducible $Y$-module such that 
$V\downarrow X$ is multiplicity-free and $\l$ is not $\l_1$ or its dual.  

Recall the usual notation. There are two levels $W^1(Q_X) \cong V_{L_X'}(\o_2)$ and $W^2(Q_X) \cong V_{L_X'}(\o_1)$,
so $L_X' < L_Y' = C^0 \times C^1,$  where $C^0 = A_{r_0}$ for $r_0 = \frac{(l+1)l}{2}-1$ and $C^1 = A_l.$ Let $\mu^i$ be the restriction of $\l$ to $T_Y\cap C^i$, so that 
$V^1(Q_Y)\downarrow L_Y' = V_{C^0}(\mu^0) \otimes V_{C^1}(\mu^1)$. 

We break up the analysis into the three cases: $l=2$ (Subsection \ref{leq2}), $l=3$ (Subsection \ref{leq3}), and  $l\ge 4$ (Subsection \ref{lge4}).

\section{$X = A_3$, $\delta = \omega_2$}\label{leq2}

In this section we establish the following result.

\begin{thm}\label{om2a3}
Let $X = A_3$, $W = V_X(\om_2)$ and $Y = SL(W) \cong A_5$. Suppose $V = V_Y(\l)$ is an irreducible $Y$-module such that $V\downarrow X$ is MF. Then up to duals, $\l$ is one of the following weights:
\begin{itemize}
\item[{\rm (i)}] $a\l_i\,(a \ge 1)$, 
\item[{\rm (ii)}] $\l_i+a\l_j\,(a \ge 1)$,
\item[{\rm (iii)}] $11100$ or $11001$.
\end{itemize}
Each of these is in Table $\ref{TAB4}$ of Theorem $\ref{MAINTHM}$.
\end{thm}

We prove the theorem  in a series of lemmas. Let $X,W,\l$ be as in the hypothesis of the theorem, and write 
\[
\l=abcde.
\]
We have $A_2 \cong L_X' < L_Y' = C^0\times C^1 \cong A_2A_2$, and as $L_X'$-modules, $W^1(Q_X) \cong 01$ and $W^2(Q_X) \cong 10$. Hence 
\begin{equation}\label{v1res}
V^1 = ba \otimes de,
\end{equation}
where as usual, we abbreviate $V^i(Q_Y)\downarrow L_X'$ by $V^i$.

\begin{lem}\label{abnon}
Theorem $\ref{om2a3}$ holds when $ab\ne 0$ and $(d,e) \ne (0,0)$.
\end{lem}

\pf Suppose $ab\ne 0$ and $(d,e) \ne (0,0)$. Since $V^1$ is MF, it follows by 
Proposition \ref{tensorprodMF} that either $d$ or $e$ is 0.

Assume first that $d\ne 0$, so that $e=0$. In $V^2$  the following summands appear, afforded by the given weights:
\[
\begin{array}{l|l}
\hline
\hbox{summand of }V^2  & \hbox{afforded by } \\
\hline
(b-1,a+1) \otimes (d+1,0) & \l-\b_2-\b_3 \\
(b,a-1) \otimes (d+1,0) & \l-\b_1-\b_2-\b_3 \\
(b-1,a+1) \otimes (d-1,1) & \l-\b_2-\b_3-\b_4 \\
(b,a-1) \otimes (d-1,1) & \l-\b_1-\b_2-\b_3-\b_4 \\
(b+1,a) \otimes (d-1,1) & \l-\b_3-\b_4 \\
\hline
\end{array}
\]
The first four summands sum to 
\[
((b-1,a+1) + (b,a-1))\otimes d0 \otimes 10.
\]
This contains $(ba \otimes 01 \otimes d0)+((b-1,a)\otimes d0)$, the first term of which is $V^1 \otimes 01$. Hence, including the fifth summand in the table above, we have
\[
(V^1 \otimes 01)+ ((b-1,a)\otimes d0) +  ((b+1,a) \otimes (d-1,1)) \subseteq V^2.
\]
Now $ (b+1,a) \otimes (d-1,1)$ is not MF unless $d=1$, in which case both $(b-1,a)\otimes d0$ and $(b+1,a) \otimes (d-1,1)$ have a summand $(b,a)$. So in any case we have a contradiction by Corollary \ref{cover}(ii). 

Hence $d=0$, and so $e\ne 0$ by hypothesis. Then $V^2$ has the following summands:
\[
\begin{array}{l|l}
\hline
\hbox{summand of }V^2  & \hbox{afforded by } \\
\hline
(b-1,a+1) \otimes (1,e) & \l-\b_2-\b_3 \\
(b,a-1) \otimes (1,e) & \l-\b_1-\b_2-\b_3 \\
(b-1,a+1) \otimes (0,e-1) & \l-\b_2-\b_3-\b_4-\b_5 \\
(b,a-1) \otimes (0,e-1) & \l-\b_1-\b_2-\b_3-\b_4-\b_5 \\
(b+1,a) \otimes (0,e-1) & \l-\b_3-\b_4-\b_5 \\
\hline
\end{array}
\]

Suppose $c\ne 0$. Then $V^2$ also has a summand $(b+1,a) \otimes (1,e)$, afforded by $\l-\b_3$, and the sum of this together with the summands in the above table is the tensor product $ba \otimes 10 \otimes 10 \otimes 0e$, which is equal to 
\[
ba \otimes (01+20)\otimes 0e = (V^1\otimes 01) + (ba\otimes 20\otimes 0e).
\]
Since $ba\otimes 20\otimes 0e$ is not MF, this contradicts Corollary \ref{cover}(ii).

Hence $c=0$. The first four summands in the table above sum to 
\[
((b-1,a+1) + (b,a-1))\otimes 0e \otimes 10.
\]
This contains $V^1\otimes 01 + ((b-1,a)\otimes 0e) + Z$, where 
\[
Z= \left\{ \begin{array}{l}
(b-2,a+2)\otimes 0e, \hbox{ if }b\ge 2 \\
(b,a-2)\otimes 0e, \hbox{ if }a\ge 2 \\
0, \hbox{ otherwise.}
\end{array}
\right.
\]
If $b\ge 2$ then $((b-1,a)\otimes 0e) + Z$ contains $(b-1,a+e)^2$; and if $a\ge 2$ then 
$((b-1,a)\otimes 0e) + Z$ contains $(b,a+e-2)^2$. Hence $a,b\le 1$ by Corollary \ref{cover}(ii).

We now have $\l = 1100e$. When $e=1$, this is conclusion (iv) of Theorem \ref{om2a3}. So assume $e\ge 2$.
Then including the fifth summand in the above table, $V^2$ contains \linebreak $(V^1\otimes 01) + (01\otimes 0e) + (21\otimes (0,e-1))$.
Since $ 01\otimes 0e + 21\otimes (0,e-1)$ contains $(1,e-1)^2$, this again contradicts  Corollary \ref{cover}(ii).
\hal 

\begin{lem}\label{abnonlem}
Theorem $\ref{om2a3}$ holds when $ab\ne 0$ and $d=e=0$.
\end{lem}

\pf Assume the hypothesis, so that $\l = abc00$ with $a,b\ne 0$.

Suppose first that $c=0$ and $a,b\ge 2$. Then $V^1 = ba$ and the weights $\l-\b_2-\b_3, \l-\b_1-\b_2-\b_3$ give  
\begin{equation}\label{v2facs}
\begin{array}{ll}
V^2 & = ((b-1,a+1)\otimes 10) + ((b,a-1)\otimes 10) \\
       &  = (b,a+1)+(b-2,a+2)+(b-1,a)^2+(b+1,a-1)+(b,a-2).
\end{array}
\end{equation}
Now $V^3$ contains the following summands (see Proposition \ref{pieri}):
\[
\begin{array}{l|l}
\hline
\hbox{summand of }V^3  & \hbox{afforded by } \\
\hline
(b-2,a+2) \otimes 20 & \l-2\b_2-2\b_3 \\
(b,a-2) \otimes 20 & \l-2\b_1-2\b_2-2\b_3 \\
(b-1,a) \otimes 20 & \l-\b_1-2\b_2-2\b_3 \\
\hline
\end{array}
\]
In addition, the weight $\l-\b_1-2\b_2-2\b_3-\b_4$ has multiplicity 3 in $V^3$ (it is conjugate to $\l-\b_1-2\b_2-\b_3$), whereas it appears only twice in the above summands (the first and third in the table). Hence $V^3$ has a summand 
$(b-1,a)\otimes 01$ in addition to those in the table above. Each summand in the table, and also $(b-1,a)\otimes 01$, has a composition factor $(b-2,a)$, so $(b-2,a)^4 \subseteq V^3$. However, of the composition factors of $V^2$ in (\ref{v2facs}), only $(b-1,a)$ has $(b-2,a)$ in its first level by Corollary \ref{V^2(Q_X)}. This contradicts Corollary \ref{conseq}(ii).

It follows that if $c=0$, then $a$ or $b$ is 1, and so $\l$ is as in Theorem \ref{om2a3}(ii).

Now suppose that $c=1$ and $a,b\ge 2$. Then in addition to the summands in (\ref{v2facs}), $V^2$ has a summand $(b+1,a)\otimes 10$ afforded by $\l-\b_3$, and so 
\begin{equation}\label{v2facs1}
V^2 = (b,a+1)^2+(b+1,a-1)^2+(b-1,a)^2+(b+2,a)+(b,a-2)+(b-2,a+2).
\end{equation}
As in the previous case, the only composition factors among these that have $(b-2,a)$ in their first level are $(b-1,a)^2$. However we see similarly that $V^3$ contains $(b-2,a)^4$, contradicting Corollary \ref{conseq}(ii).

Hence if $c=1$ then $a$ or $b$ is 1. Hence by Proposition \ref{d1100}, $\l$ is $11100$ as in Theorem \ref{om2a3}(iii).

Finally, the case where $c \ge 2$ is ruled out by Proposition \ref{abc00,c>1}. \hal

\vspace{2mm}
In view of the previous lemmas, we can now assume that $ab=0$, and also by duality that $de=0$.

\begin{lem}\label{adnon}
Theorem $\ref{om2a3}$ holds when $a\ne 0$ and $d\ne 0$.
\end{lem}

\pf By the above remarks $\l = a0cd0$, $V^1 = 0a \otimes d0$ and $V^2$ has the following summands:
\[
\begin{array}{l|l}
\hline
\hbox{summand of }V^2  & \hbox{afforded by } \\
\hline
(0,a-1) \otimes (d+1,0) & \l-\b_1-\b_2-\b_3 \\
(1,a) \otimes (d-1,1) & \l-\b_3-\b_4 \\
(0,a-1) \otimes (d-1,1) & \l-\b_1-\b_2-\b_3-\b_4 \\
\hline
\end{array}
\]
If $c\ne 0$ then $V^2$ also has a summand $(1,a)\otimes (d+1,0)$ afforded by $\l-\b_3$, and the sum of this together with the summands in the above table is $10\otimes 0a\otimes d0\otimes 10$, which is equal to 
$(V^1 \otimes 01) + ( 0a\otimes d0\otimes 20)$. The latter summand is not MF as it contains $(d+1,a-1)^2$, so this contradicts 
Corollary \ref{cover}(ii).

Hence $c=0$. Suppose $a\ge 2$ and $d\ge 2$. The second and third entries in the above table sum to 
$10\otimes 0a\otimes (d-1,1)$, which is equal to 
\[
0a \otimes (d1+(d-1,0)+(d-2,2)) = (0a \otimes d0\otimes 01) + (0a\otimes (d-2,2)).
\]
Adding in the first summand in the above table, we see that $V^2$ contains 
\[
(V^1 \otimes 01) + (0a\otimes (d-2,2)) + ((0,a-1) \otimes (d+1,0)).
\]
However the latter two summands both contain $(d,a-2)$, so this contradicts Corollary \ref{cover}(ii).

It follows that $a$ or $d$ is equal to 1, so $\l = a0010$ or $100d0$, as in (i) of Theorem \ref{om2a3}. \hal

\begin{lem}\label{aenon}
Theorem $\ref{om2a3}$ holds when $a\ne 0$ and $e\ne 0$.
\end{lem}

\pf Here $\l = a0c0e$, and replacing $V$ by $V^*$ if necessary we can assume that $a\ge e$.
We have $V^1 = 0a \otimes 0e$, and $V^2$ has the following summands:
\[
\begin{array}{l|l}
\hline
\hbox{summand of }V^2  & \hbox{afforded by } \\
\hline
(0,a-1) \otimes (1,e) & \l-\b_1-\b_2-\b_3 \\
(1,a) \otimes (0,e-1) & \l-\b_3-\b_4-\b_5 \\
(0,a-1) \otimes (0,e-1) & \l-\b_1-\b_2-\b_3-\b_4-\b_5 \\
\hline
\end{array}
\]
If $c\ne 0$ then $V^2$ also has a summand $(1,a)\otimes (1,e)$ afforded by $\l-\b_3$, and the sum of this together with the summands in the above table is $10\otimes 0a\otimes 0e\otimes 10$, which is equal to 
$(V^1 \otimes 01) + ( 0a\otimes 0e\otimes 20)$. The latter summand is not MF as it contains $(1,a+e-1)^2$, so this contradicts 
Corollary \ref{cover}(ii).

Hence $c=0$ and $\l = a000e$. Assume $a\ge e \ge 2$. Now
\[
V^1 =  0a \otimes 0e = (0,a+e)+(1,a+e-2)+(2,a+e-4)+\cdots +(e,a-e).
\]
However, all three of the summands in the above table have $(0,a+e-2)$ as a composition factor, whereas only the composition factor $(1,a+e-2)$ of $V^1$ has $(0,a+e-2)$ in its first level (by Corollary \ref{V^2(Q_X)}). This contradicts 
Corollary \ref{conseq}(i).

It follows that under the assumption that $a\ge e$, we must have $e=1$, so that $\l = a0001$, as in (i) of the theorem. \hal

\begin{lem}\label{anondeo}
Theorem $\ref{om2a3}$ holds when $a\ne 0$ and $d=e= 0$.
\end{lem}

\pf Here $\l = a0c00$. Assume $a,c \ge 2$. Then the weights $ \l-\b_1-\b_2-\b_3, \l-\b_3$ give
\[
\begin{array}{ll}
V^2 & = 1a \otimes 10 + (0,a-1)\otimes 10 \\
        & = (2,a)+(0,a+1)+(1,a-1)^2+(0,a-2).
\end{array}
\]
Of these composition factors only $(1,a-1)$ has $(2,a-2)$ in its first level. However, we argue in the usual way that $V^3$ contains $(2,a-2)^4$, arising from the following summands:
\[
\begin{array}{l|l}
\hline
\hbox{summand of }V^3  & \hbox{afforded by } \\
\hline
2a\otimes 20 & \l-2\b_3 \\
(1,a-1)\otimes 20 & \l-\b_1-\b_2-2\b_3 \\
(0,a-2)\otimes 20 &  \l-2\b_1-2\b_2-2\b_3 \\
(1,a-1)\otimes 01 &  \l-\b_1-\b_2-2\b_3-\b_4 \\
\hline
\end{array}
\]
This contradicts Corollary \ref{conseq}(ii).  Therefore $c = 0$ or $c=1$ or $a=1$.  In each case Theorem \ref{om2a3} holds. \hal

\vspace{4mm}
\no {\bf Completion of the proof of Theorem \ref{om2a3}}

In view of the previous few lemmas we can assume that $a=0$, and also $e=0$ (by duality). Hence $\l = 0bcd0$.
If $b=d=0$ then $\l = c\l_3$ is as in (i) of the theorem, so (by duality) we may assume that $b\ne 0$ and also $b\ge d$.

Suppose $d\ne 0$. Then 
\[
V^1 = b0\otimes d0 = (b+d,0)+(b+d-2,1)+\cdots +(b-d,d),
\]
 and $V^2$ has the following summands: 
\[
\begin{array}{l|l}
\hline
\hbox{summand of }V^2  & \hbox{afforded by } \\
\hline
(b-1,1) \otimes (d+1,0) & \l-\b_2-\b_3 \\
(b+1,0) \otimes (d-1,1) & \l-\b_3-\b_4 \\
(b-1,1) \otimes (d-1,1) & \l-\b_2-\b_3-\b_4 \\
\hline
\end{array}
\]
If $c\ne 0$ then in addition $\l-\b_3$ affords a summand $(b+1,0)\otimes (d+1,0)$ in $V^2$, and the sum of this and the  above three summands is $b0\otimes 10\otimes d0\otimes 10$.  This is equal to $(V^1\otimes 01)+(b0\otimes d0\otimes 20)$, and the latter summand is not MF, contradicting Corollary \ref{cover}(ii).

Hence $c=0$ and $\l = 0b0d0$. If $b \ge d \ge 2$, then each of the three summands in the above table has $(b+d-2,2)$ as a composition factor, whereas only the composition factor $(b+d-2,1)$ of $V^1$ has $(b+d-2,2)$ in its first level; this is a contradiction  by Corollary \ref{V^2(Q_X)}. It follows that $d = 1$, and so $\l$ is as in Theorem \ref{om2a3}(ii).

It remains to consider the case where $d=0$, so that $\l = 0bc00$. Assume $b,c\ge 2$. Then 
\[
\begin{array}{ll}
V^2 & = ((b-1,1)\otimes 10)+((b+1,0)\otimes 10) \\
        & = (b+2,0)+(b,1)^2+(b-2,2)+(b-1,0),
\end{array}
\]
while $V^3$ contains the following summands:
\[
\begin{array}{l|l}
\hline
\hbox{summand of }V^3  & \hbox{afforded by } \\
\hline
(b+2,0)\otimes 20 & \l-2\b_3 \\
(b,1)\otimes 20 & \l-\b_2-2\b_3 \\
(b-2,2)\otimes 20 &  \l-2\b_2-2\b_3 \\
(b,1)\otimes 01 &  \l-\b_2-2\b_3-\b_4 \\
\hline
\end{array}
\]
Each of these summands has a composition factor $b2$, so $V^3$ contains $b2^4$. On the other hand only the composition factor $b1$ of $V^2$ has $b2$ in its first level, so this contradicts Corollary \ref{conseq}(ii). It follows that either $b\le 1$ or $c\le 1$, so $\l$ is as in Theorem \ref{om2a3}(i, ii).

This completes the proof of Theorem \ref{om2a3}.

\section{$X = A_4$, $\delta=\omega_2$}\label{leq3}

In this subsection we establish the following result.

\begin{thm}\label{om2a4}
Let $X = A_4$, $W = V_X(\om_2)$ and $Y = SL(W) \cong A_9$. Suppose $V = V_Y(\l)$ is an irreducible $Y$-module such that $V\downarrow X$ is MF. Then up to duals, $\l$ is one of the following weights:
\[
\begin{array}{l}
a\l_1, \\
a\l_1+\l_2, \\
a\l_1+\l_9, \\
\l_i, \\
\l_1+\l_i, \\
a\l_2\,(a\le 5), \\
2\l_3,\,2\l_4, \\
\l_1+2\l_2,\,\l_2+\l_3,\,\l_2+\l_8.
\end{array}
\]
Each of these is in Tables $\ref{TAB2} - \ref{TAB4}$ of Theorem $\ref{MAINTHM}$.
\end{thm}

We prove the theorem  in several further subsections. Let $X,W,\l$ be as in the hypothesis of the theorem. We have 
 $A_3 \cong L_X' < L_Y' = C^0\times C^1 \cong A_5A_3$, and as $L_X'$-modules, $W^1(Q_X) \cong 010$ and $W^2(Q_X) \cong 100$. As usual, let $\mu^i = \l\downarrow T_Y \cap C^i$ for $i=0,1$. 

As usual we shall abbreviate $V^i(Q_Y)\downarrow L_X'$ by $V^i$. In particular, 
\[
V^1 \cong (V_{C^0}(\mu^0) \otimes V_{C^1}(\mu^1))\downarrow L_X'
\]
is MF. Hence the possibilities for $\mu^0$ are given by Theorem \ref{om2a3}.

 \subsection{The case where $\mu^1=0$} \label{mu10a4}

Assume in this subsection that the hypotheses of Theorem \ref{om2a4} hold, and that $\mu^1=0$. Note that $\b_6$ is the fundamental root lying between $C^0$ and $C^1$ in the Dynkin diagram of $Y=A_9$.

\begin{lem}\label{lnu2}
Suppose that $\mu^0\ne 0$ and that every composition factor of $V_{C^0}(\mu^0) \downarrow L_X'$ is a multiple of a fundamental dominant weight. Then $\mu^0 = a\l_1^0$, $a\l_5^0$ or $\l_3^0$.
\end{lem}

\pf This is immediate fro Lemma \ref{twolabs}.  \hal

\begin{lem}\label{a4g0}
If $\mu^0 \ne 0$, then either $\la \l,\b_6 \ra = 0$ or $\l=\l_1+\l_6$ (as in the conclusion of Theorem $\ref{om2a4}$).
\end{lem}

\pf Suppose $\mu^0 \ne 0$ and $x = \la \l,\b_6 \ra \ne 0$. By Proposition \ref{v2gamma1},
\begin{equation}\label{v2cont}
V^2 \supseteq V^1 \otimes 010 \otimes 100 = V^1 \otimes (001 + 110).
\end{equation}
Hence $V^1 \otimes 110$ is MF by Corollary \ref{cover}(ii). Now Proposition \ref{tensorprodMF} shows that each composition factor of $V^1$ must be a multiple of a fundamental dominant weight, so Lemma \ref{lnu2} shows that 
\begin{equation}\label{mu0possib}
\mu^0 = a\l_1^0,\,a\l_5^0 \hbox{  or } \l_3^0.
\end{equation}

Assume first that $\mu^0=a\l_1^0$. If  $a\ge 2$ then Theorem \ref{SOpowers} shows that 
$V^1 = S^a(010)$ contains $0a0 + 0(a-2)0$, and hence Theorem \ref{LR} implies
\[
V^1 \otimes 110 \supseteq (0a0 + 0(a-2)0)\otimes 110 \supseteq (1(a-1)0)^2,
\]
contradicting the fact that $V^1 \otimes 110$ is MF. Hence $a=1$. Now $V^2(Q_Y)$ has two summands, afforded by the weights $\l-\b_6$ and $\l-\b_1-\cdots -\b_6$; these summands are $V_{C^0}(\l_1^0+\l_5^0)\otimes V_{C^1}(100)$ and $V_{C^0}(0) \otimes V_{C^1}(100)$ respectively. Hence
\begin{equation}\label{v21}
V^2 = 010 \otimes 010 \otimes 100 = 011^2+100^2+120+201.
\end{equation}
If $x\ge 2$ then $V^3(Q_Y)$ has the following summands:
\[
\begin{array}{l|l}
\hline
\hbox{summand of }V^3(Q_Y)  & \hbox{afforded by } \\
\hline
V_{C^0}(\l_1^0+2\l_5^0)\otimes V_{C^1}(200) & \l-2\b_6 \\
V_{C^0}(\l_1^0+\l_4^0)\otimes V_{C^1}(010) & \l-\b_5-2\b_6-\b_7 \\
\hline
\end{array}
\]
Restricting these summands to $L_X'$, we see that $V^3 \supseteq 121^3$. However, of the composition factors of $V^2$ in (\ref{v21}), only one has 121 in its first level (namely, the composition factor 120). This contradicts Corollary \ref{conseq}(ii).
Hence $x=1$, and now we have $\l = \l_1+\l_6$, as in the conclusion. 

Next assume that $\mu^0 = a\l_5^0$, the second case of (\ref{mu0possib}). As above we deduce that $a=1$ and that the composition factors of $V^2$ are as in (\ref{v21}) (this time coming from summands afforded by $\l-\b_6$ and $\l-\b_5-\b_6$). If $x\ge 2$, then $V^3(Q_Y)$ has the following summands:
\[
\begin{array}{l|l}
\hline
\hbox{summand of }V^3(Q_Y)  & \hbox{afforded by } \\
\hline
V_{C^0}(3\l_5^0)\otimes V_{C^1}(200) & \l-2\b_6 \\
V_{C^0}(\l_4^0+\l_5^0)\otimes V_{C^1}(200) & \l-\b_5-2\b_6 \\
V_{C^0}(\l_4^0+\l_5^0)\otimes V_{C^1}(010) & \l-\b_5-2\b_6-\b_7 \\
\hline
\end{array}
\]
Restricting to $L_X'$, we see that $V^3 \supseteq 121^3$ and as before this contradicts  Corollary \ref{conseq}(ii).
Hence $x=1$ and $\l = \l_5+\l_6$. However a Magma computation shows that $V_{A_9}(\l_5+\l_6) \downarrow X$ is not MF. 

Now assume  that $\mu^0 = \l_3^0$, the final case of (\ref{mu0possib}). Here $V^1 = \wedge^3(010) = 200+002$.
It follows that  $V^1 \otimes 110 \supseteq 011^2$, contradicting the fact that $V^1\otimes 110$ is MF. \hal

\vspace{2mm}
We now suppose that $\mu^0\ne 0$, so that we may assume $\la \l,\b_6\ra =  0$ by Lemma \ref{a4g0}. In particular, $\l = \mu^0$. The possibilities for $\mu^0$ are given  by Theorem \ref{om2a3}, and we divide these up as in Table \ref{mu0posss}, according to various cases involving the dual $\l^*$.  In the table, $(\mu^*)^i$ is the restriction of $\l^*$ to $C^i$ for $i=0,1$, we let $x^* = \la \l^*,\b_6\ra$, and $a>0$. Note that in the table, the cases in (6) and (7) have been separated off from the others for convenience of reference, even though they satisfy the conditions in the second column for previous cases.

In the following, we shall freely use information given in the results of Section \ref{a3a5} concerning the composition factors of the (MF) restrictions of the modules $V_{A_5}(\mu^0)$ to $L_X' = A_3$.

\begin{table}[h]
\caption{}\label{mu0posss}
\[
\begin{array}{|l|l|l|}
\hline
\hbox{case} & \hbox{condition} & \hbox{possible }\mu^0 \\
\hline
(1) & (\mu^*)^0\ne 0,\,(\mu^*)^1=0, \,x^*\ne 0 & 0001a,\,000a1 \\
\hline
(2) & (\mu^*)^0= 0,\,x^*= 0 & 1a000,\,a1000,\, \\
  & & 10a00,\,a0100, \\
  & & 01a00,\,0a100 \\
\hline
(3)  & (\mu^*)^0= 0,\,x^*= 1 & a0010,\,0a010,\,00a10 \\
\hline
(4)  & (\mu^*)^0= 0,\,x^*= a & 100a0,\,010a0,\,001a0 \\
\hline
(5) & (\mu^*)^0\ne 0,\,(\mu^*)^1\ne 0 & 1000a,\,a0001, \\
  & & 0100a,\,0a001, \\
  & & 0010a,\,00a01 \\
\hline
(6) & & 11100,\,11001, \\
   & & 00111,\,10011 \\
\hline
(7) & & a\l_i^0 \\
\hline
\end{array}
\]
\end{table}

\begin{lem}\label{t91}
Theorem $\ref{om2a4}$ holds when  $\mu^0$ is as in case $(1)$ of Table $\ref{mu0posss}$.
\end{lem}

\pf In this case $(\mu^*)^0\ne 0$, $(\mu^*)^1=0$ and $x^* = \la \l^*,\b_6\ra \ne 0$. So Lemma \ref{a4g0} applied to $\l^*$ gives the conclusion. \hal

\begin{lem}\label{t92}
Theorem $\ref{om2a4}$ holds when  $\mu^0$ is as in case $(2)$ of Table $\ref{mu0posss}$.
\end{lem}

\pf Suppose $\mu^0$ is as in (2) of Table \ref{mu0posss}. If $\mu^0 = a1000$ then $\l = a\l_1+\l_2$, which is in the list in the conclusion of Theorem \ref{om2a4}; and we shall postpone the case where $\mu^0 = a0100$ until the end of this proof.  In the other cases we have
\[
(\mu^*)^0= 0,\,x^*= 0,\, (\mu^*)^1 = a01,\,0a1,\,a10 \hbox{ or }1a0.
\]
If $a=1$, or if $(\mu^*)^1  = 0a1$ with $a=2$, then $\l$ is the dual of a weight in the conclusion of Theorem \ref{om2a4},   so assume from now on that $a\ge 2$, and also that $a\ge 3$ in the case where $(\mu^*)^1  = 0a1$.

We now work out the full list of composition factors of $(V^*)^2$. In each case this has precisely two summands afforded by two weights in the set $\{\l^*-\b_6-\b_7,\,\l^*-\b_6-\b_7-\b_8,\,\l^*-\b_6-\b_7-\b_8-\b_9\}$:
\[
\begin{array}{|l|l|l|}
\hline
(\mu^*)^1   & \hbox{summands of }(V^*)^2 &  \hbox{comp. factors of }(V^*)^2 \\
\hline
a01 & 010 \otimes (a-1)11,\,010\otimes a00 & a10^2,\,(a-1)01^2,\,(a-1)21, \\
    & & a02,\,(a-2)12,\,(a-2)20 \\ 
\hline
0a1 & 010\otimes 0(a-1)2,\,010\otimes 0a0 & 1(a-1)1^2,\,0a2,\,1(a-2)3, \\
  & & 0(a-2)2,\,0(a+1)0,\,0(a-1)0 \\
\hline
a10 & 010\otimes (a-1)20,\,010\otimes a01 & a11^2,\,(a-1)10^2,\,(a-1)30, \\
   & & (a-2)21,\,(a-1)02,\,(a+1)00 \\
\hline
1a0 & 010\otimes 0(a+1)0,\,010\otimes 1(a-1)1 & 1a1^2,\,0a0^2,\,0(a+2)0, \\
   & & 2(a-2)2,\,0(a-1)2,\,2(a-1)0, \\
  & & 1(a-2)1 \\
\hline
\end{array}
\]
Next we work out some summands of $(V^*)^3(Q_Y)$. In the table below, we denote by $\l^*-abc\ldots$ the weight $\l^*- a\b_5-b\b_6-c\b_7-\cdots$:
\[
\begin{array}{|l|l|l|}
\hline
(\mu^*)^1  & \hbox{summands of }(V^*)^3(Q_Y) &  \hbox{afforded by } \\
\hline
a01 &2\l_5^0\otimes (a-2)21 & \l^*-022 \\ 
 &2\l_5^0\otimes (a-1)10 & \l^*-02211 \\ 
 &\l_4^0\otimes (a-1)02 & \l^*-1221 \\ 
 &\l_4^0\otimes (a-1)10 & \l^*-12211 \\ 
\hline
0a1 &2\l_5^0\otimes 0(a-2)3 & \l^*-0222 \\ 
 &2\l_5^0\otimes 0(a-1)1 & \l^*-02221 \\ 
 &\l_4^0\otimes 0(a-1)1 & \l^*-12221 \\ 
\hline
a10 &2\l_5^0\otimes (a-2)30 & \l^*-022 \\ 
 &2\l_5^0\otimes (a-1)11 & \l^*-0221 \\ 
 &\l_4^0\otimes (a-1)11 & \l^*-1221 \\ 
\hline
1a0 &2\l_5^0\otimes 0a1 & \l^*-0221 \\ 
 &2\l_5^0\otimes 1(a-2)2 & \l^*-0222 \\ 
 &\l_4^0\otimes 0a1 & \l^*-1221 \\ 
 &\l_4^0\otimes 1(a-1)0 & \l^*-12221 \\ 
\hline
\end{array}
\]
Now $V_{C^0}(2\l_5^0) \downarrow L_X' = 020+000$ and $V_{C^0}(\l_4^0) \downarrow L_X' = 101$. So restricting the above summands to $L_X'$, we deduce that $(V^*)^3$ contains the following composition factors:
\[
\begin{array}{ll}
(\mu^*)^1  = a01: & (V^*)^3 \supseteq ((a-2)21)^5 \\
(\mu^*)^1  = 0a1: & (V^*)^3 \supseteq (1(a-3)2)^3 \\
(\mu^*)^1  = a10: & (V^*)^3 \supseteq ((a-2)30)^4 \\
(\mu^*)^1  = 1a0: & (V^*)^3 \supseteq (1(a-2)2)^5.
\end{array}
\]
Let $k=5,3,4,5$ denote the multiplicity of the composition factor in the above list in the respective cases. We now check that at most $k-2$ of these composition factors are in the first level of the composition factors of $(V^*)^2$. This is a contradiction by Corollary \ref{conseq}(ii). 

It remains to handle the postponed case where $\mu^0 = a0100$. We deal with this directly with $V$, not the dual.
The case $a = 1$ is in the conclusion of Theorem \ref{om2a4}.  So now assume $a \ge 2$. Here  $V^2(Q_Y)$ is the sum of  $ (a-1\,0100) \otimes (100)$ and $ (a1000) \otimes (100)$, 
and these add to  $(a0000) \otimes (01000) \otimes (100).$  Restricting to $L_X'$, this is 
$S^a(010) \otimes (101) \otimes (100)$, which is
\[
\left(0a0 + 0(a-2)0+ \cdots\right) \otimes \left(201 + 011 + 100\right).
\]
Using Theorem \ref{LR} we find that $0a0$ tensored with each of $100$, $011$, and $201$ produces a term $0(a-1)1$.  Another such summand is contained in  $(0(a-2)0) \otimes (011).$  Therefore, $V^2$  contains $(0(a-1)1)^4$, and this highest weight has $S$-value $a$.

Now $V^1(Q_Y) = a0100 = (a0000 \otimes 00100) - ((a-1)0010)$. Restricting to $L_X'$ and using Lemma \ref{a3om2}, we see that $V^1$ is the difference of
\[
\left(0a0 + 0(a-2)0 + \cdots\right) \otimes (200 + 002)
\]
and
\[
\left((1(a-1)1) +(1(a-3)1) + \cdots \right) + \left((2(a-2)0)^+ + (2(a-4)0)^+ + \cdots\right).
\]
 Any summands of $V^1$ that contribute a composition factor $0(a-1)1$ in $V^2$ must have $S$-value 
$a-1$, $a$ or $a+1$.  We easily see that all such summands lie in $\left(0a0 + 0(a-2)0\right) \otimes (200 + 002)$.
Now $S$-value considerations show that the possible summands are $0(a-2)2$, $1(a-1)1$, $2(a-2)0$ and $1(a-3)1$, each of which occurs
with multiplicity at most 1 in $V^1$.  Of these, only $0(a-2)2$ and $1(a-1)1$ can actually yield a factor $0(a-1)1$. So this yields at most
two terms $0(a-1)1$ and we conclude that $V_Y(\l) \downarrow X$ is not MF.  \hal

\begin{lem}\label{t93}
Theorem $\ref{om2a4}$ holds when  $\mu^0$ is as in case $(3)$ of Table $\ref{mu0posss}$.
\end{lem}

\pf This is similar to the previous lemma. Suppose $\mu^0$ is as in (3) of Table \ref{mu0posss}. Here
\[
(\mu^*)^0= 0,\,x^*= 1,\, (\mu^*)^1  = 00a,\,0a0 \hbox{ or }a00.
\]
If $a=1$ then a Magma computation rules out $\l = \l_2+\l_4$ or $\l_3+\l_4$, and $\l = \l_1+\l_4$ is in the conclusion of 
Theorem \ref{om2a4}.  So assume from now on that $a\ge 2$.

We work out the composition factors of $(V^*)^2$: this has two summands, afforded by $\l^*-\b_6$ and  $\l^*-\b_6-\b_7-\cdots -\b_m$ with $m \in \{7,8,9\}$:
\[
\begin{array}{lll}
(\mu^*)^1  = 00a:\;(V^*)^2 & = & 010\otimes (10a+00(a-1)) \\
                                               & = & 01(a-1)^2,\,11a,\,20(a-1), \\
                                               &  & 00(a+1),\, 10(a-2) \\
(\mu^*)^1  = 0a0:\;(V^*)^2 & = & 010\otimes (1a0+0(a-1)1) \\
                                               & = & 1(a-1)0^2,\,0a1^2,\,1(a+1)0, \\
                                               &  & 2(a-1)1,\, 1(a-2)2,\,0(a-2)1 \\
(\mu^*)^1  = a00:\;(V^*)^2 & = & 010\otimes ((a+1)00+(a-1)10) \\
                                               & = & a01^2,\,(a+1)10,\,(a-1)20, \\
                                               & & (a-1)00,\,(a-2)11.
\end{array}
\]
Next, in $(V^*)^3(Q_Y)$ we have the following summands, where $\l^*-abc\ldots$ denotes $\l^*-a\b_5-b\b_6-\cdots$: 
\[
\begin{array}{|l|l|l|}
\hline
(\mu^*)^1   & \hbox{summands of }(V^*)^3(Q_Y) &  \hbox{afforded by } \\
\hline
00a &2\l_5^0\otimes 10(a-1) & \l^*-02111 \\ 
 &2\l_5^0\otimes 00(a-2) & \l^*-02222 \\ 
 &\l_4^0\otimes 01a & \l^*-121 \\ 
 &\l_4^0\otimes 10(a-1) & \l^*-12111 \\ 
\hline
0a0 &2\l_5^0\otimes 1(a-1)1 & \l^*-0211 \\ 
 &2\l_5^0\otimes 0(a-2)2 & \l^*-0222 \\ 
 &\l_4^0\otimes 0(a+1)0 & \l^*-121 \\ 
 &\l_4^0\otimes 1(a-1)1 & \l^*-1211 \\ 
\hline
a00 &2\l_5^0\otimes a10 & \l^*-021 \\ 
 &2\l_5^0\otimes (a-2)20 & \l^*-022 \\ 
 &\l_4^0\otimes a10 & \l^*-121 \\ 
\hline
\end{array}
\]
Restricting to $L_X'$, we see that $(V^*)^3$ contains $(02(a-2))^k$, $(0a2)^k$ or $((a-1)21)^k$ with $k=4$, 4 or 3 in the respective cases for $(\mu^*)^1$. However, only $k-2$ of these composition factors are in the first level of composition factors of $(V^*)^2$, contradicting Corollary \ref{conseq}(ii). \hal

\begin{lem}\label{t94}
Theorem $\ref{om2a4}$ holds when  $\mu^0$ is as in case $(4)$ of Table $\ref{mu0posss}$.
\end{lem}

\pf  This is again similar to the previous lemmas. Suppose $\mu^0$ is as in (4) of Table \ref{mu0posss}. Then
\[
(\mu^*)^0= 0,\,x^*= a,\, (\mu^*)^1 = 001,\,010 \hbox{ or }100.
\]
The case $a=1$ follows as in the previous lemma, so assume from now on that $a\ge 2$.

We work out the composition factors of $(V^*)^2$: this has two summands, afforded by $\l^*-\b_6$ and one of  $\l^*-\b_6-\b_7$, $\l^*-\b_6-\b_7-\b_8$ and $\l^*-\b_6-\b_7-\b_8-\b_9$:
\[
\begin{array}{lll}
(\mu^*)^1  = 001:\;(V^*)^2 & = & 010\otimes (101+000) \\
                                               & = & 010^2,\,111,\,200,\,002 \\
(\mu^*)^1  = 010:\;(V^*)^2 & = & 010\otimes (110+001) \\
                                               & = & 011^2,\,120, 201,100^2 \\
(\mu^*)^1  = 100:\;(V^*)^2 & = & 010\otimes (200+010) \\
                                               & = & 210,\,101^2,\,020,\,000. \\
\end{array}
\]
In $(V^*)^3(Q_Y)$ we have the following summands, where $\l^*-abc\ldots$ denotes $\l^*-a\b_5-b\b_6-\cdots$:
\[
\begin{array}{|l|l|l|}
\hline
(\mu^*)^1   & \hbox{summands of }(V^*)^3(Q_Y) &  \hbox{afforded by } \\
\hline
001 &2\l_5^0\otimes 201 & \l^*-02 \\ 
 &2\l_5^0\otimes 100 & \l^*-02111 \\ 
 &\l_4^0\otimes 011 & \l^*-121 \\ 
\hline
010 &2\l_5^0\otimes 210 & \l^*-02 \\ 
 &2\l_5^0\otimes 101 & \l^*-0211 \\ 
 &\l_4^0\otimes 020 & \l^*-121 \\ 
\hline
100 &2\l_5^0\otimes 300 & \l^*-02 \\ 
 &2\l_5^0\otimes 110 & \l^*-021 \\ 
 &\l_4^0\otimes 110 & \l^*-121 \\ 
\hline
\end{array}
\]
Restricting to $L_X'$, we see that $(V^*)^3$ contains $120^k$, $210^k$ or $211^k$ with $k=3,\,4$ or 3 in the respective cases for $(\mu^*)^1$. 
However, only $k-2$ of these composition factors are in the first level of composition factors of $(V^*)^2$, contradicting Corollary \ref{conseq}(ii). \hal

\begin{lem}\label{t95}
Theorem $\ref{om2a4}$ holds when  $\mu^0$ is as in case $(5)$ of Table $\ref{mu0posss}$.
\end{lem}

\pf  Suppose $\mu^0$ is as in (5) of Table \ref{mu0posss}. We postpone the cases where $\mu^0 = 1000a$, $0100a$ or $0010a$ until later in the proof. For the other cases we have
\[
(\mu^*)^0= 00001,\,x^*= 0,\, (\mu^*)^1 = 00a,\,0a0 \hbox{ or }a00.
\]
The case $a=1$ follows as in previous lemmas, so assume from now on that $a\ge 2$.

Now $(V^*)^1 = ((\mu^*)^0 \otimes (\mu^*)^1)\downarrow L_X'$, so the composition factors of $(V^*)^1$ are as follows:
\[
\begin{array}{l}
(\mu^*)^1 = 00a:\;(V^*)^1 = 01a,\,10(a-1) \\
(\mu^*)^1 = 0a0:\;(V^*)^1 = 0(a+1)0,\,1(a-1)1, 0(a-1)0 \\
(\mu^*)^1 = a00:\;(V^*)^1 = a10,\,(a-1)01. \\
\end{array}
\]
In $(V^*)^2$ we have the following summands, where $\l^*-abc\ldots$ denotes $\l^*-a\b_5-b\b_6-\cdots$:
\[
\begin{array}{|l|l|l|}
\hline
(\mu^*)^1   & \hbox{summands of }(V^*)^2 &  \hbox{afforded by } \\
\hline
00a &(020+000)\otimes 00(a-1) & \l^*-01111 \\ 
 & 101\otimes 10a & \l^*-11 \\ 
 & 101 \otimes 00(a-1) & \l^*-11111 \\ 
\hline
0a0 &(020+000)\otimes 0(a-1)1 & \l^*-0111 \\ 
 & 101\otimes 1a0 & \l^*-11 \\ 
 & 101 \otimes 0(a-1)1 & \l^*-1111 \\ 
\hline
a00 &(020+000)\otimes (a-1)10 & \l^*-011 \\ 
 & 101\otimes (a+1)00 & \l^*-11 \\ 
 & 101 \otimes (a-1)10 & \l^*-111 \\ 
\hline
\end{array}
\]
Restricting to $L_X'$, we see that $(V^*)^2$ contains $(00(a-1))^k$, $(0(a-1)1)^k$ or $((a-1)10)^k$ with $k=3,\,5$ or 4 in the respective cases for $(\mu^*)^1$. However, at most $k-2$ of these composition factors are in the first level of composition factors of $(V^*)^1$, contradicting Corollary \ref{conseq}(i). 

It remains to deal with the cases where $\mu^0 = 1000a$, $0100a$ or $0010a$. 

First suppose $\mu^0 = 1000a$. By Lemma \ref{a3om2}(i), 
\[
V^1 = (0(a+1)0) + (1(a-1)1) + (0(a-1)0) +(1(a-3)1) +(0(a-3)0) +   \cdots.
\]
On the other hand $V^2(Q_Y)$ contains summands $0000a \otimes 100$ and $(1001(a-1)) \otimes (100)$.  Since 
$10000 \otimes (0001(a-1)) = (1001(a-1)) + (0000a) + (0001(a-2))$, it follows that 
\[
V^2(Q_Y) \supseteq \left((10000) \otimes (0001(a-1)) - (0001(a-2))\right) \otimes 100.
\]
Using Lemma \ref{a3om2}(iv) we restrict the  first tensor product to $L_X'$ and obtain
\[
010 \otimes\left(1(a-1)1 + 1(a-3)1 + \cdots+  0(a-1)0 + 0(a-3)0 + \cdots\right) \otimes 100,
\]
although the term $000$ does not occur in the middle tensor factor if $a$ is odd.
Now $010 \otimes 0(a-1)0$ and $010 \otimes 1(a-3)1$  each contain $1(a-2)1$,  while $010 \otimes 1(a-1)1$ 
contains $2(a-2)2$ and  $010 \otimes 1(a-1)1$ contains $2(a-1)0$. Therefore, tensoring with $100$ we obtain 
 $(2(a-2)1)^4$.  Another application of Lemma \ref{a3om2}(iv) shows that this irreducible appears in the subtracted
tensor product with multiplicity 1.  Therefore $V^2$ contains $(2(a-2)1)^3.$  On the other
hand,  at most one such term can arise from $V^1$.  Therefore $V_Y(\l) \downarrow X$ is not MF in this case.  

Next suppose $\mu^0 = 0100a$. The case $a=1$ is handled with a Magma calculation, so assume $a\ge 2$.  
 By Lemma \ref{a3om2}(ii), 
\[
V^1 = (1a1) + (1(a-2)1) + \cdots + (2(a-1)0)^+ + (2(a-3)0)^+ \cdots.
\]
 We will be counting  summands of highest weight  $2(a-1)1$ in $V^2$, 
and we note that at most two such summands can arise from $V^1$.

Now $V^2(Q_Y)$ has  summands $1000a \otimes 100$ and  $(0101(a-1))\otimes 100$.  Moreover
 $01000 \otimes (0001(a-1)) = (0101(a-1)) + (1000a) + (1001(a-2)) + (0000(a-1))$. Therefore $V^2(Q_Y)$ contains
\[
 \left((01000 \otimes 0001(a-1)) \otimes 100\right) - \left((1001(a-2) + 0000(a-1)) \otimes 100\right).
\]
 Restricting the first summand to $L_X'$ and using Lemma \ref{a3om2}(iv), we obtain
\[
101 \otimes \left(1(a-1)1 + 1(a-3)1 + \cdots + 0(a-1)0+0(a-3)0 + \cdots\right) \otimes 100,
\]
although the term $000$ does not occur in the middle tensor product if $a$ is odd.  Using Theorem \ref{LR} we see that 
$101 \otimes 1(a-1)1$ contains $(1(a-1)1)^3$, $2(a-1)2$ and $3(a-2)1$.  Also $101 \otimes 0(a-1)0$ contains $1(a-1)1$.  So tensoring with $100$ we obtain $(2(a-1)1)^6$.

So it remains to consider the multiplicity of $2(a-1)1$ in the restrictions of $(1001(a-2)) \otimes 100$ and
 $ 0000(a-1) \otimes 100$.  The first tensor product is contained in $\left(10000 \otimes 0001(a-2)\right) \otimes 100$ and
 using Lemma \ref{a3om2}(iv) again, we see that $2(a-1)1$ occurs with multiplicity 1 in this.  It does not occur in the
 second tensor product.  Therefore, $V^2$ contains $(2(a-1)1)^5$ and $V_Y(\l) \downarrow X$ is not MF. 

Finally, suppose $\mu^0 = 0010a$. 
Consider $V^*$, where $(\mu^*)^0 = 0000a$ and $(\mu^*)^1 = 100$.  Here $(V^*)^1 = S^a(010) \otimes (100)$. 
Also $(V^*)^2(Q_Y)$ contains $((0000(a+1)) + (0001(a-1)))\otimes (010)$, and the restriction of this to $L_X'$ is 
$S^a(010) \otimes (010) \otimes (010)$; this contains $((V^*)^1\otimes 001) + (S^a(010)\otimes 020)$. By Theorem \ref{SOpowers} we have 
$S^a(010) \supseteq (0a0)+(0(a-2)0)$, so $S^a(010)\otimes 020$ contains $(0a0)^2$ and is not MF, contradicting Corollary \ref{cover}.
  \hal

\begin{lem}\label{t96}
Theorem $\ref{om2a4}$ holds when  $\mu^0$ is as in case $(6)$ of Table $\ref{mu0posss}$.
\end{lem}

\pf 
In this case $\mu^0 = 11100$, $11001$, $00111$, or $10011$, and $\l = \mu^0$. To establish the assertion we must show that in each case $V \downarrow X$ fails to be MF.  We will do this by considering $V^2$.

  First assume  $\mu^0 = 11100.$  Passing to $V^*$ we have $(\mu^*)^0 = 0$ and $(\mu^*)^1 = 111.$  Therefore $(V^*)^2 \supseteq 010 \otimes (021 +102+ 110)$.  This contains  $(011)^3$ but only one such  summand can arise from $(V^*)^1$, a contradiction.

Next suppose that $\mu^0 = 00111$. Here $(\mu^*)^0 = 00001$,   $(\mu^*)^1 = 100$, and  $\la \l^*,\b_6\ra = 1.$
Then $(V^*)^1  = 010 \otimes 100 = 110 + 001$. In $(V^*)^2$ there are composition factors afforded by $\l^*  -\b_6$, $\l^*  -\b_5-\b_6,$ and $\l^* -\b_6-\b_7.$  Restricting to $L_X'$ we see that these contain $(111)^3$.  Only one such composition factor can arise from  $(V^*)^1$ so this is a contradiction.

If $\mu^0 = 11001,$ then $(\mu^*)^0 = 00001$ and $(\mu^*)^1 = 011.$ Then
$(V^*)^1 = 010 \otimes 011 = 021 + 102 + 110 + 001.$  In $(V^*)^2$ there are composition factors afforded by 
$\l^* - \b_5-\b_6$, $\l^*  -\b_6-\b_7-\b_8$ and $\l^*  -\b_6-\b_7-\b_8-\b_9.$  Restricting to $L_X'$ we see that these contain $(111)^5$ and only three of these summands can arise from $(V^*)^1$.  So here too we have a contradiction.

Finally assume that $\mu^0 = 10011$ so that $(\mu^*)^0 = 00001,$  $(\mu^*)^1 = 001,$ 
and $\la \l^*,\b_6\ra = 1.$   In $(V^*)^2$ there are composition factors afforded by $\l^*  -\b_6$ and $\l^*  -\b_5-\b_6.$  Restricting to $L_X'$ we see that these contain $(101)^4$.  However only two such composition factors can arise from   $(V^*)^1  = 010 \otimes 001 = 011 + 100$.  So again we have a contradiction.  \hal

\vspace{2mm}
To complete this subsection (the case where $\mu^1=0$), it remains to handle the cases where $\mu^0 = a\l_i^0$ (as in (7) of Table \ref{mu0posss}), and also where $\mu^0 = 0$.

\begin{lem}\label{t97l3}
Theorem $\ref{om2a4}$ holds when  $\mu^0=a\l_3^0$.
\end{lem}

\pf Suppose that $\mu^0 = a\l_3^0$ (so $\l = a\l_3$). For $a\le 2$ this is in the list
in the conclusion of Theorem \ref{om2a4}, and for $a=3$ a Magma computation shows that $V_Y(\l)\downarrow X$ is not MF.   So we assume that $a \ge 4$.

 It is convenient to work with $\l^* = a\l_7$. For this we have $(\mu^*)^0=0$, $x^*=0$, $(\mu^*)^1 = a00$. 

Consider $(V^*)^3(Q_Y)$. This has summands $2\l_5^0 \otimes (a-2)20$ and $\l_4^0\otimes a-1\,01$, afforded by $\l^*-022$ and $\l^*-1221$ respectively (where $\l^*-abc\ldots = \l^*-a\b_5-\cdots$). Restricting to $L_X'$, these give the following composition factors of $(V^*)^3$:
\begin{equation}\label{v3cfs1}
\begin{array}{l}
(a-2)20^3,\,(a-1)01^3,\,a02^2,\,(a-2)12^2,\,(a-3)11^2, \\
(a-1)21,\,(a-3)31,\,(a-2)40,\,(a-4)22,\,a10,( a-2)00^2. 
\end{array}
\end{equation}
Note that our final contradiction will come from bounding the number of summands in $(V^*)^3$ which 
could give rise to the summand $(a-3)20$ in $(V^*)^4$ and applying Corollary \ref{icover}. This summand can arise from 
summands of the form $(a-3)11$, $(a-4)30$ or $(a-2)20$, which are summands of $S$-value at least $a-1$. 

We claim that the list (\ref{v3cfs1}) contains all composition factors of $(V^*)^3$ with $S$-value at least $a-1$.
To prove this, 
let $\eta$ afford the highest weight of an $L_Y'$-summand of $(V^*)^3(Q_Y)$. Then $\eta$ is subdominant to 
the weight $\lambda^*-2\beta_6-2\beta_7$. So $\eta$ is of the form 
$\lambda^*-2\beta_6-x\beta_7-y\beta_8-z\beta_9$ or
 $\lambda^*-\beta_5-2\beta_6-x\beta_7-y\beta_8-z\beta_9$, for $x,y,z$ non negative integers satisfying:
$x\geq 2$, $y\geq 2z$, $x+z\geq 2y$ and $y+a+2\geq 2x$. Moreover, since this weight should afford an $L_X'$ 
summand with $S$-value at least $a-1$, we have the additional constraint that $x+z\leq 5$. So we see that 
there are a finite number of possible triples $(x,y,z)$, and we will consider them in turn.

It is easy to see that the two listed weights, corresponding to the triples $(2,0,0)$  and $(2,1,0)$, 
afford summands. 
Let us call these two weights $\mu = \lambda^*-2\beta_6-2\beta_7$ and 
$\nu = \lambda^*-\beta_5-2\beta_6-2\beta_7-\beta_8$. If $x=2$, these are the only possible weights affording 
summands. Note that for all remaining cases, we may assume $y\ne 0$, as if $y=0$ then $z=0$ and the multiplicity
of the given weight is 1 and already occurs in the summand afforded by $\mu$ or $\nu$.

So now suppose $x=3$; then we have $z\in\{0,1\}$ and $y\leq \frac{3+z}{2}$. In the case 
$(x,z) = (3,0)$, we must consider the weights
$\lambda^*-2\beta_6-3\beta_7-\beta_8$ and  $\lambda^*-\beta_5-2\beta_6-3\beta_7-\beta_8$.
 Since $a\geq 4$, we may apply Proposition \ref{cav1} to see that the multiplicity of the above 
weights are respectively 2, 3. It is then straightforward to see that the multiplicity in the 
summand afforded by $\mu$ is 2, and the second weight occurs also in the summand 
afforded by $\nu$ with multiplicity 1. Hence none of these weights affords an $L_Y'$-summand of 
$(V^*)^3(Q_Y)$. For the case $(x,z)=(3,1)$, consider the weight $\lambda^*-\beta_5-2\beta_6-3\beta_7-2\beta_8-\beta_9$.
Using Proposition \ref{cav1}, we see that this weight has multiplicity 6 in $V^*$, while in the summand afforded by $\mu$ it occurs with 
multiplicity 3, and with the same multiplicity in the summand afforded by $\nu$, so does not afford
an additional $L_Y'$-summand. The other case is similar, though easier.

Suppose now that $x=4$, so $z\in\{0,1\}$. If $(x,z)=(4,0)$, then $1\leq y\leq 2$.  The corresponding 
weights are \begin{enumerate}
\item $\lambda^*-2\beta_6-4\beta_7-\beta_8$, 
\item $\lambda^*-2\beta_6-4\beta_7-2\beta_8$,  
\item $\lambda^*-\beta_5-2\beta_6-4\beta_7-\beta_8$, and 
\item $\lambda^*-\beta_5-2\beta_6-4\beta_7-2\beta_8$.
\end{enumerate}

For all of the above weights, since $a\geq 4$, we may use 4.9 to calculate their multiplicities,
 replacing $a$ by $4$ in the weight $\lambda^*$. In each case, we see that the weight occurs with the same 
multiplicity in the sum of the two summands afforded by $\mu$ and $\nu$ as in the module $V^*$ and so does not afford a summand.  We argue similarly for $(x,z)=(4,1)$, where $y=2$.

Consider the pair $(x,z)=(5,0)$, where we must have $y\in\{1,2\}$. Then the corresponding 
weights are \begin{enumerate}
\item $\lambda^*-2\beta_6-5\beta_7-\beta_8$, 
\item $\lambda^*-2\beta_6-5\beta_7-2\beta_8$,  
\item $\lambda^*-\beta_5-2\beta_6-5\beta_7-\beta_8$, and 
\item $\lambda^*-\beta_5-2\beta_6-5\beta_7-2\beta_8$.
\end{enumerate}

 For the first
 weight 
we have $y=1$ and since $y+a+2\geq 2x=10$, we have $a\geq 7$ and so we may use Proposition \ref{cav1} to 
determine the multiplicity of this weight by considering its multiplicity in the module where we 
replace $a$ by 5,
which is then seen to be 2. On the other hand, its multiplicity in the summand afforded by $\mu$ 
is also 2 (again using 4.9 to replace $a$ by  5).
For the second weight, we have $a\geq 6$ and can again calculate the multiplicity of the weight
in the module where we replace $a$ by 5, giving multiplicity 3, while again this is the 
multiplicity in the summand afforded by $\mu$. The other two cases are entirely similar; for 
example, the last
weight has multiplicity 4 while its multiplicity in the summand afforded by $\mu$ is 
3 and in the summand afforded by $\nu$ is 1.

This proves our claim that the list (\ref{v3cfs1}) contains all composition factors of $(V^*)^3$ with $S$-value at least $a-1$.

Now consider $(V^*)^4(Q_Y)$. This has summands as follows, writing $\l^*-abc\ldots $ for the weight $\l^*-a\b_4-b\b_5-\cdots$:
\[
\begin{array}{|l|l|}
\hline
\hbox{summand of }(V^*)^4(Q_Y) &  \hbox{afforded by } \\
\hline
3\l_5^0\otimes (a-3)30 & \l^*-0033 \\
(\l_4^0+\l_5^0)\otimes (a-2)11 & \l^*-01331 \\
\l_3^0\otimes (a-1)00 & \l^*-123321 \\
\hline
\end{array}
\]
Restricting to $L_X'$, we calculate that $(V^*)^4 \supseteq  ((a-3)20)^7$. However only five of these composition factors can appear in the first level of the list (\ref{v3cfs1}), and since this list contain all 
composition factors of $(V^*)^3$ of $S$-value at least $a-1$, this contradicts Corollary \ref{conseq}(ii). \hal


\begin{lem}\label{t97l4}
Theorem $\ref{om2a4}$ holds when  $\mu^0=a\l_4^0$.
\end{lem}

\pf Suppose that $\mu^0 = a\l_4^0$ (so $\l = a\l_4$). For $a\le 2$ this is in the list
in the conclusion of Theorem \ref{om2a4}, and for $a=3$ a Magma computation shows that $V_Y(\l)\downarrow X$ is not MF. So we assume that $a \ge 4$.

Now $(V^*)^4$ has precisely three summands $3\l_5^0 \otimes 300$,  $(\l_4^0+\l_5^0)\otimes 110$ 
and $\l_3^0\otimes 001$, afforded by $\l^*-003$ and $\l^*-0131$ and $\l^*-12321$ respectively (where $\l^*-abc\ldots = \l^*-a\b_4-\cdots$). Hence the full list of composition factors of $(V^*)^4$ is:
\begin{equation}\label{v4list}
201^5,\,011^4,\,003^3,\,100^3,\,120^3,\,112^3,\,221^2,\,310^2,\,302,\,031,\,330.
\end{equation}
In $(V^*)^5(Q_Y)$ we have the following summands (where this time $\l^*-abc\ldots = \l^*-a\b_3-b\b_4-\cdots$):
\[
\begin{array}{|l|l|}
\hline
\hbox{summand of }(V^*)^5(Q_Y) &  \hbox{afforded by } \\
\hline
4\l_5^0\otimes 400 & \l^*-0004 \\
(\l_4^0+2\l_5^0)\otimes 210 & \l^*-00141 \\
2\l_4^0\otimes 020 & \l^*-00242 \\
(\l_3^0+\l_5^0)\otimes 101 & \l^*-012421 \\
\hline
\end{array}
\]
Restricting to $L_X'$, we calculate that $(V^*)^5 \supseteq (311)^7$. However only five of these composition factors can appear in the first level of the list (\ref{v4list}). This contradicts Proposition \ref{induct}. \hal

\begin{lem}\label{t97l5}
Theorem $\ref{om2a4}$ holds when  $\mu^0=a\l_5^0$.
\end{lem}

\pf   Here $\l = a\l_5$. For $a = 1$ this is an MF example in Table \ref{TAB3} by Theorem \ref{wedgec},  and for $a=2$ a Magma check shows that $V_Y(\l)\downarrow X$ is not MF; so  assume $a \ge 3$.
  Then  $V^2(Q_Y) = (0001(a-1)) \otimes 100$, so restricting to $L_X'$ and using Lemma \ref{a3om2}(iv) we have 
\[
V^2 = \left(1(a-1)1 + 1(a-3)1 + \cdots + 0(a-1)0+0(a-3)0 + \cdots\right) \otimes 100,
\]
although $000$ does not appear in the first tensor factor when $a$ is odd.  We will be concerned
with $L_X'$-composition factors in $V^3$ of highest weight $3(a-2)1$.  This can arise from only one composition factor of $V^2$, namely $2(a-1)1$.

Now $V^3(Q_Y)$ contains $(0002(a-2)) \otimes 200$, and
\begin{equation}\label{subt}
0002(a-2) = \left(00010 \otimes (0001(a-2))\right) - (0010(a-1)) - (0100(a-2)) - (0011(a-3)),
\end{equation}
Restricting the tensor product to $L_X'$ and using Lemma \ref{a3om2}(iv), we obtain
\[
101 \otimes \left(1(a-2)1 + 1(a-4)1 + \cdots + 0(a-2)0 + \cdots\right).
\]
Now $101 \otimes (1(a-2)1) \supseteq (1(a-2)1)^3 + (2(a-2)2) + (2(a-1)0)+(3(a-3)1)$, and 
$101 \otimes (0(a-2)0) \supseteq 1(a-2)1$.  Tensoring with $200$ produces $(3(a-2)1)^7$.

In addition  $V^3(Q_Y)$ contains $(0010(a-1)) \otimes 010$ (afforded by
$\l - \b_4 - 2\b_5 - 2\b_6 - \b_7$), and using Lemma \ref{a3om2}(iii) we see that this produces
an additional term $3(a-2)1$, giving  $(3(a-2)1)^8$.

To count the multiplicity of $3(a-2)1$ in $V^3$, we must consider the subtracted 
 terms  $(0010(a-1))$, $(0100(a-2))$, $(0011(a-3))$ in (\ref{subt}), and look for the multiplicity of $3(a-2)1$ after tensoring with $200$.   The second and third terms are contained in $00100 \otimes (0001(a-3))$, 
and restricting to $L_X'$ this is
\[
(200 + 002) \otimes  \left((1(a-3)1) + (1(a-5)1) + \cdots + (0(a-3)0) + \cdots\right).
\]
The only possibilities here arise from $(200 + 002) \otimes  (1(a-3)1)$, and the composition factors of this of $S$-value at least $a$ are 
$(1(a-2)1)^2+(3(a-3)1)^++(2(a-4)2)^2$. Tensoring this with $200$ produces just $(3(a-2)1)^3$.  The remaining subtracted term is $(0010(a-1))$,
and using Lemma \ref{a3om2}(iii) we see that the restriction to $L_X'$ of $(0010(a-1)) \otimes 200$
contains just $(3(a-2)1)^2$.  

We conclude that $V^3$ contains $(3(a-2)1)^3$, and hence $V_Y(\l) \downarrow X$ is not MF.  \hal

\begin{lem}\label{t97both0}
Theorem $\ref{om2a4}$ holds when  $\mu^0=0$.
\end{lem}

\pf In this case $\mu^0=\mu^1=0$, so $\l = a\l_6$. Then $\l^* = a\l_4$, so this is covered by Lemma \ref{t97l4}. \hal

\vspace{2mm}
By the previous lemmas, the only remaining cases are $\mu^0 = a\l_1^0$ or $a\l_2^0$, in which case $\l = a\l_1$ or 
$a\l_2$. In the latter case, $a\le 5$ by Lemma \ref{alambda2}. Hence $\l$ is as in the conclusion of Theorem \ref{om2a4}, as required.

This completes the analysis of this subsection, the case where $\mu^1=0$.

 \subsection{The case where $\mu^1\ne 0$}\label{mu1ne0a4}

In this subsection we complete the proof of Theorem \ref{om2a4} by handling the case where $\mu^1\ne 0$.

By the previous Subsection \ref{mu10a4}, we may assume that $(\mu^*)^1\ne 0$, and hence also $\mu^0 \ne0$.

\begin{lem}\label{mu0mu1posss}
The possibilities for $\mu^0$ and $\mu^1$ are as follows:
\[
\begin{array}{|l|l|}
\hline
\mu^0 & \mu^1 \\
\hline
a\l_1^0 & \hbox{any} \\
\l_3^0 & \hbox{any} \\
a\l_3^0\,(a\ge 2) & 100,\,001 \\
\l_2^0 & d00,\,0d0,\,00d \\
a\l_2^0\,(a\ge 2) & 100,\,001 \\
\hline
\end{array}
\]
\end{lem}

\pf First note that $\la \mu^0,\b_i^0\ra \ne 0$ for some $i\le 3$, since otherwise $(\mu^*)^1\ne 0$.

Suppose now that every composition factor of $V_{C^0}(\mu^0)\downarrow L_X'$ is a multiple of a fundamental dominant weight. Then $\mu^0 = a\l_1^0$ or $\l_3^0$ by Lemma \ref{lnu2}.  Hence $\mu^0$ is as in the first two lines of the table in the conclusion.

So assume now that  $V_{C^0}(\mu^0)\downarrow L_X'$ has a composition factor with at least two nonzero labels (so $\mu^0 \ne a\l_1^0,\l_3^0$). As $V^1 = (\mu^0\otimes \mu^1)\downarrow L_X'$ is MF,  Proposition \ref{tensorprodMF} implies that $\mu^1 = d00$, $0d0$ or $00d$.

The possiblities for $\mu^0$ are given by Theorem \ref{om2a3}. In Table \ref{mu0mu1cfs} we list 
(using results from Section \ref{taby4}) some composition factors of 
$\mu^0\downarrow L_X'$ and also of $(\mu^0\otimes \mu^1)\downarrow L_X'$. Since the latter restriction is MF, it follows from the table that $\mu^0,\mu^1$ are as in the conclusion of the lemma. \hal

\begin{table}[h]
\caption{}\label{mu0mu1cfs}
\[
\begin{array}{|l|l|l|l|}
\hline
\mu^0 & \hbox{some cfs of }\mu^0\downarrow L_X' & \mu^1 & \hbox{some cfs of } \\
&&& (\mu^0 \otimes \mu^1)\downarrow L_X' \\
\hline
a\l_3^0& (2a)00,\,(2a-2)02, & d00\,(d\ge 2) & ((2a+d-4)20)^2 \\
                             \,(a\ge 2)                     & (2a-4)20,\,00(2a),  & 0d0\,(d\ge 2) & ((2a-2)(d-2)2)^2 \\
                                                                              & 20(2a-2),\,02(2a-4)  & 010 & ((2a-1)01)^2 \\
\hline
a\l_2^0 & 202,\,020,\,000 & d00\,(d\ge 2) & d00^2 \\
                                          \,(a\hbox{ even})                               && 0d0 & 1d1^2 \\
a\l_2^0 & 303,\,121,\,101 & d00\,(d\ge 2) & ((d+1)01)^2 \\
                                                 \,(a\ge 3 \hbox{ odd})                            &   & 0d0 & 1d1^2 \\
\hline
a\l_1^0+\l_5^0 & 0(a+1)0,\,1(a-1)1, & d00 & ((d-1)a1)^2 \\
                                                                               &  0(a-1)0    & 0d0 & (1(a+d-1)1)^2 \\
\hline
a\l_1^0+\l_4^0 & 1a1,\,2(a-1)0, & d00 & ((d+1)(a-2)1)^2\,(a\ge 2) \\
                                                                  & 1(a-2)1\,(a\ge 2),\, && ((d-2)20)^2\,(a=1,d\ge2) \\
                                                                  &0(a-1)2 && 110^2\,(a=d=1) \\
                                                       && 0d0 & (1(a+d-2)1)^2\,(a\ge 2) \\
                                                         &         && 2d0^2\,(a=1) \\
\hline
\l_1^0+a\l_4^0 & a1a,\,(a-1)0(a+1), & d00 & ((a+d-1)0(a+1))^2 \\
                           & (a+1)0(a-1)  & 0d0 & (a(d-1)a)^2 \\
\hline
a\l_1^0+\l_3^0 & 1(a-1)1,\,2a0, & d00 & ((d+1)(a-1)1)^2 \\
                 &0a2 & 0d0 & (1(a+d-1)1)^2 \\
\hline
\l_1^0+a\l_3^0 & (2a)10,\,(2a-1)01, & d00 & ((2a+d-1)01)^2 \\
                           & 01(2a),\,10(2a-1) & 0d0 & ((2a-1)d1)^2 \\
\hline
a\l_1^0+\l_2^0  & 1a1,\,0a0  & d00 & da0^2 \\
                      && 0d0 & (0(a+d)0)^2 \\
\hline
\l_1^0+a\l_2^0 & 111,\,010 & d00 & d10^2 \\
                         && 0d0 & (0(d+1)0)^2 \\
\hline
a\l_2^0+\l_3^0 & (a+2)0a,\,a1a,\,  & d00 & ((a+d)1a)^2 \\
                           & a0(a+2),\,(a-1)1(a-1) & 0d0 & ((a-1)(d+1)(a-1))^2 \\
\hline
\l_2^0+a\l_3^0  & (2a-1)11,\,11(2a-1), & d00 &  ((2a+d-3)21)^2 \\
                           & (2a-3)21,\,12(2a-3),  & 0d0 & ((2a-2)d2)^2 \\
\hline
a\l_2^0+\l_4^0 & (a+1)0(a+1),\,(a-1)1(a+1), & d00 & ((a+d-1)1(a+1))^2 \\
                           & (a+1)1(a-1),\,a0a & 0d0 & (ada)^2 \\
\hline
11001 & 121,\,202,\,012,\,210 & d00 & (d12)^2 \\
                             && 0d0 & (2d2)^2 \\
\hline
11100 & 113,\,311,\,121,\,202 & d00 & ((d+1)21)^2 \\
                 &&  0d0 & (2d2)^2 \\
\hline
\end{array}
\]
\end{table}

\begin{lem}\label{x0obv}
We have $\la \l,\b_6\ra = 0$.
\end{lem}

\pf This follows by applying Lemma \ref{mu0mu1posss} to the dual $\l^*$. \hal

\begin{lem}\label{al1}
Theorem $\ref{om2a4}$ holds when  $\mu^0=a\l_1^0$.
\end{lem}

\pf Suppose $\mu^0=a\l_1^0$, and let $\mu^1 = yzw$. Applying Lemma \ref{mu0mu1posss} to the dual $\l^*$, we see that just one of $y,z,w$ is nonzero.

Assume first that $y\ne 0$, so $\l = a\l_1+y\l_7$. Suppose that $y\ge 2$. Then $a=1$ by Lemma \ref{mu0mu1posss}
applied to $V^*$, and so 
\[
V^1 = 010 \otimes y00 = y10 + (y-1)01;
\]
but $V^2(Q_Y)$ contains $\l_1^0 \otimes \l_5^0 \otimes ((y-1)10)$, whence
\[
V^2 \supseteq 010 \otimes 010 \otimes  ((y-1)10) \supseteq ((y-2)21)^2.
\]
Since $(y-2)21$ is not in the first level of any composition factor of $V^1$, this contradicts  Corollary \ref{conseq}(i). 

Hence $y=1$. Consider $\l^* = \l_3+a\l_9$. If $a=1$ this is in the list in Theorem \ref{om2a4}, and if $a=2$ a Magma computation shows that $V_Y(\l)\downarrow X$ is not MF. So assume $a\ge 3$. Now $V_{C^0}(\l_3^0) \downarrow L_X' = 200+002$, so 
\[
(V^*)^1 = (200+002)\otimes 00a =  20a + 00(a+2) + 10(a-1) + 01a + 00(a-2) + 02(a-2).
\]
Next, $(V^*)^2(Q_Y)$ has the following summands:
\[
\begin{array}{|l|l|}
\hline
\hbox{summand of }(V^*)^2(Q_Y) &  \hbox{afforded by } \\
\hline
\l_2^0\otimes 10a & \l^*-\b_3-\cdots-\b_6 \\
(\l_3^0+\l_5^0)\otimes 00(a-1) & \l^*-\b_6-\cdots - \b_9 \\
\l_2^0\otimes 00(a-1) & \l^*-\b_3-\cdots-\b_9 \\
\hline
\end{array}
\]                                                           
Restricting to $L_X'$, the first summand is $101 \otimes 10a$, and the second and third sum to 
\[
(\l_3^0\otimes \l_5^0)\downarrow L_X' \otimes 00(a-1) = (200+002)\otimes 010\otimes (00(a-1)).
\]
Hence we see that $(V^*)^2 \supseteq (11(a-2))^4$. However only two of these composition factors are in the first level of a composition factor of $(V^*)^1$, so this contradicts Corollary \ref{conseq}(i).  

Next assume $z\ne 0$, so $\l=a\l_1+z\l_8$. Lemma \ref{mu0mu1posss} applied to $\l^*$ shows that either $a=1$ or $z=1$.

Suppose $a=1$. Then we can assume that $z\ge 2$, as $\l_1+\l_8$ is on the list in Theorem \ref{om2a4}. So 
\[
V^1 = 010 \otimes 0z0 = (0(z+1)0) + (1(z-1)1) + (0(z-1)0).
\]
Also 
\[
\begin{array}{ll}
V^2 & \supseteq (\l_1^0 \otimes \l_5^0)\downarrow L_X' \otimes (0(z-1)1)  \\
        & = (010\otimes 010\otimes 0(z-1)1) \\
        & \supseteq (0(z-1)1)^4.
\end{array}
\]
However only two of these composition factors are in the first level of a composition factor of $V^1$, contradicting 
Corollary \ref{conseq}(i).  

Now suppose $z=1$. As above, we may assume that $a\ge 2$. Here $\l^* = \l_2+a\l_9$, and
\[
(V^*)^1 = 101 \otimes 00a = (10(a+1)) + (11(a-1)) + (00a) + (01(a-2)).
\]
Also 
\[
\begin{array}{ll}
(V^*)^2 & \supseteq (\l_2^0 \otimes \l_5^0)\downarrow L_X' \otimes (00(a-1))\;+\; (\l_1^0\downarrow L_X' \otimes 10a)  \\
        & = (101\otimes 010\otimes 00(a-1))\;+\; 010\otimes 10a \\
        & \supseteq (20(a-1))^3.
\end{array}
\]
However only one of these composition factors is in the first level of a composition factor of $V^1$, contradicting 
Corollary \ref{conseq}(i).  

Finally, assume $w\ne 0$, so $\l = a\l_1+w\l_9$. If $a$ or $w$ is 1 then $\l$ is in the list in Theorem \ref{om2a4}, so assume $a,w\ge 2$. Then 
\[
V^1 = S^a(010)\otimes 00w = (0a0\,+\,0(a-2)0 +\cdots )\otimes 00w.
\]
The composition factors of this having $S$-value at least $a+w-3$ are
\begin{equation}\label{aw3}
0aw,\,1(a-1)(w-1),\,2(a-2)(w-2),\,3(a-3)(w-3),\,0(a-2)w,\,1(a-3)(w-1),
\end{equation}
all occurring in $V^1$ with multiplicity at most 1. Also
\[
\begin{array}{ll}
V^2 & \supseteq (a\l_1^0 \otimes \l_5^0)\downarrow L_X' \otimes (00(w-1))\;+\; ((a-1)\l_1^0\downarrow L_X' \otimes 10w)  \\
        & \supseteq (0(a-1)(w-1))^4.
\end{array}
\]
However only two of these composition factors are in the first level of a composition factor in the list (\ref{aw3}), contradicting Corollary \ref{conseq}(i). \hal

\begin{lem}\label{al3}
Theorem $\ref{om2a4}$ holds when  $\mu^0=a\l_3^0$.
\end{lem}

\pf Suppose first that $\mu^0 = \l_3^0$.   Let $\mu^1 = dbc$. Applying Lemma \ref{mu0mu1posss} to the dual $\l^*$, we see that only one of $d,b,c$ is nonzero. If $c \ne 0$ then $(\mu^*)^0 = c\l_1^0$, a case covered by Lemma \ref{al1}. If $b\ne 0$ then 
\[
V^1 = \l_3^0\downarrow L_X' \otimes 0b0 = (200+002)\otimes 0b0 \supseteq (1(b-1)1)^2,
\]
contradicting the fact that $V^1$ is MF. Finally, assume $d\ne 0$. If $d=1$ then $\l=\l_3+\l_7$ and a Magma check shows that $V_Y(\l)\downarrow X$ is not MF. So take $d\ge 2$. We have 
\[
\begin{array}{ll}
V^1 &= (200+002)\otimes d00 \\
       & =  ((d+2)00) + (d10) + ((d-2)20)+ (d02) + ((d-1)01) + ((d-2)00).
\end{array}
\]
In $V^2$, the weights $\l-\b_6-\b_7$, $\l-\b_3-\cdots -\b_7$ afford summands adding to $(\l_3^0\otimes \l_5^0)\downarrow L_X' \otimes ((d-1)10)$, so
\[
V^2 \supseteq (200+002)\otimes 010 \otimes ((d-1)10) \supseteq (d11)^4.
\]
However only two composition factors $d11$ are in the first level of composition factors of $V^1$, so this contradicts  
Corollary \ref{conseq}(i). 

Now suppose that $\mu^0 = a\l_3^0$ with $a\ge 2$. By Lemma \ref{mu0mu1posss}, $\mu^1 = 100$ or $001$. In the latter case the dual $\l^*$ has $(\mu^*)^0 = \l_1^0$, a case covered by Lemma \ref{al1}. In the former case, $(\mu^*)^0=\l_3^0$, which was handled above. \hal

\begin{lem}\label{al2}
Theorem $\ref{om2a4}$ holds when  $\mu^0=a\l_2^0$.
\end{lem}

\pf Suppose $\mu^0 = a\l_2^0$. If  $\mu^1 = d00$ or $00d$ then $(\mu^*)^0 = d\l_3^0$ or $d\l_1^0$, cases covered by Lemmas \ref{al1}, \ref{al3}. Hence Lemma \ref{mu0mu1posss} implies that $\mu^1 = 0d0$ and also $a=1$, so that $\l = \l_2+d\l_8$. A further application  of Lemma \ref{mu0mu1posss} to the dual $\l^*$ shows that $d=1$. It follows that $\l=\l_2+\l_8$, which is in the list in the conclusion of Theorem \ref{om2a4}. \hal

\vspace{4mm}
This completes the proof of Theorem \ref{om2a4}.


\section{$X=A_{l+1}$ with $l\ge 4$, $\delta=\omega_2$}\label{lge4}

 In this subsection we consider the case $X=A_{l+1}$ with $l\ge 4$ and $\d=\o_2.$  We establish
 the following result
 
\begin{thm}\label{deltai}  Assume $X = A_{l+1}$ with $l\ge 4$,  let $W = V_X(\om_2)$ and let $Y = SL(W)$. Suppose $V = V_Y(\l)$ is an irreducible $Y$-module such that $V\downarrow X$ is MF.
Then up to duals, $\l$ is as in Tables $\ref{TAB2}$, $\ref{TAB3}$ of Theorem $\ref{MAINTHM}$.
  \end{thm}

We will prove the theorem by induction on $l$. Assume that the hypotheses of  the theorem hold. Note that by replacing $V$ by the dual $V^*$ if necessary, we may assume that $\mu^0\ne 0$. The induction hypothesis provides us with a list of possibilities for the weight $\mu^0$, since $V_{C^0}(\mu^0) \downarrow L_X'$ is MF and $L_X'$ embeds in $C^0$ via the weight $\o_2$. We record the possiblities in the next lemma.


\begin{lem}\label{mu0posso2} The possibilities for $\mu^0$ are as follows, listed up to duals:
\[
\begin{array}{rl}
\hbox{ any }l\ge 4: & \l_i^0,\,c\l_1^0,\,2\l_2^0,\,3\l_2^0, \\
                                & \l_1^0+\l_i^0\,(i\le 7),\,\l_1+\l_{r_0+2-i}\,(2\le i\le 7), \\
                                & c\l_1^0+\l_{r_0}^0\,(c\le 3),\,c\l_1^0+\l_2^0\,(c\le 3),\, \\
                                & \l_2^0+\l_3^0,\,\l_2^0+\l_{r_0-1}^0,\\
\hbox{extras for }l=4: & 2\l_3^0,\,2\l_4^0,\, c\l_2 \,(c=4,5)\\
                                      & c\l_1^0+\l_9^0,\,c\l_1^0+\l_2^0 \,(\hbox{any }c \ge 1) \\
                                      & \l_1^0+2\l_2^0 \end{array} \] \end{lem}

We now work through each possibility for $\mu^0$ in the list.

\begin{lem}\label{mu0=cl2} Assume $\mu^0 =  c\l_2^0$ or $c\l_{r_0-1}^0$ with $c \ge 2.$  Then $\l = c\l_2$ and $c \le 3$. \end{lem}

\pf  First suppose $\mu^0 = c\l_{r_0-1}^0.$ Then the induction hypothesis implies that $(V^*)^1$ is not MF
unless $l = 4$, $c = 2$, $\mu^1=0$ and $\la \l,\g_1\ra = 0$.   Assume this occurs. Then 
\begin{equation}\label{dispy}
V^1 = S^2(0101) - \wedge^4(0010) = (0011)+(1020)+(0202)+(2000).
\end{equation}
 Now $V^2(Q_Y) = V_{C^0}(\l_7^0+\l_8^0) \otimes V_{C^1}(\l_1^1)$.
The first tensor factor restricts to $L_X'$ as $(\wedge^3(0010) \otimes \wedge^2(0010)) - (\wedge^4(0010) \otimes (0010)).$ Tensoring this with $1000$ we obtain $(1110)^3.$  On the other hand at most one factor of $(1110)$ can arise from $V^1$.  This is a contradiction.
Therefore, from now on we assume that
$\mu^0 = c\l_2^0$.  

We first claim that $\mu^1 = 0.$ Assume otherwise.  First suppose that $l \ge 5$ so that the induction hypothesis implies $c \le 3.$ 
By Proposition \ref{2l2,3l2}, $V_{C^0}(c\l_2^0)  \downarrow L_X' $ contains $(\om_1 + \om_2 + \om_5) \oplus (2\om_2 + \om_4)$ or $(2\om_2 + \om_3 + \om_5) \oplus 
(\om_1 + 2\om_2 + \om_3 + \om_4)$, according as $ c = 2$ or $c = 3.$ 

Now consider the possibilities for $\mu^1.$  Since each of the above summands contains a constituent with 3 nonzero labels, it follows from Proposition \ref{stem1.1.A} that $\mu^1 = \l_i^1, d\l_1^1$, or $d\l_l^1$ with  $d > 1$ in the latter two cases. The second case is impossible as  $(V^*)^1$ is not MF by the induction hypothesis. Therefore $\mu^1 = \l_i^1$ or $d\l_l^1.$

Assume $\mu^1 = \l_i^1$. If $i >5$, then Lemma \ref{abtimesce}(i) shows that $V^1$ contains $(\om_1 + \om_2 + \om_4 + \om_{i+1})^2$ or $(2\om_2 + \om_3 + \om_4 + \om_{i+1})^2,$ in the respective cases $c=2$ or $c=3$,  a contradiction.  And if $i\le 5$ we can use Theorem \ref{LR} to see that $V^1$ is not MF.
Now assume  $\mu^1 = d\l_l^1$. If $l > 5,$  then Lemma \ref{abtimesce}(i) shows that  $V^1$ contains $(\om_1 + \om_2 + \om_4 + (d-1)\om_l)^2$ or $(2\om_2 + \om_3 + \om_4 + (d-1)\om_l)^2,$ in the respective cases $c=2$ or $c=3$; and if $l = 5$ then Proposition \ref{pieri} shows that  $V^1$ is not MF. 
  
 To complete the proof of the claim we must still settle the case $l = 4$ where now we must allow  values of $2\le c \le 5.$ First assume $c = 2.$ Here Proposition \ref{stem1.1.A} 
together with the dual of (\ref{dispy}) shows that $\mu^1 = d\l_j^1$ for some $j$  and it is easy to see that $V^1$ is not MF using either Lemma \ref{abtimesce}(i) or Theorem \ref{LR}.  Indeed, $V^1$ contains  $((d-1)200)^2, (1(d+1)00)^2, (11d0)^2$ or $ (110d)^2$ according as $ j = 1,2,3$ or $4$.  Hence $\mu^1=0$ in this case ($c=2$).

Now assume $c \ge 3$, so $\mu^0 = c\l_2^0$ with $c = 3,4$ or 5. The composition factors of $V_{C^0}(c\l_2^0)\downarrow L_X'$ are given in Table \ref{comple} of Lemma \ref{a5,a6,a7}. Also, using  Proposition \ref{stem1.1.A}
and Lemma \ref{twolabs}, we have $\mu^1 = d\l_1^1$, $d\l_l^1$, $\l_2^1$ or $\l_3^1$. Hence we find a repeated composition factor in $V^1$ as in the following table:
\[
\begin{array}{c|ccc}
\hline
\mu^1 & \hbox{repeated cf, }c=3 & c=4 & c=5 \\
\hline
d\l_1^1 & (d210)^2 & (d211)^2 & ((d+1)111)^2 \\
d\l_4^1 & (111\,d-1)^2 & (220d)^2 & (212d)^2 \\
\l_2^1 & (1201)^2 & (1121)^2 & (1012)^2 \\
\l_3^1 & (1111)^2 & (2011)^2 & (2101)^2 \\
\hline
\end{array}
\]
In particular, $V^1$ is not MF, a contradiction. This establishes the claim that $\mu^1=0$.

 Now assume that $\langle \l, \gamma_1 \rangle >0. $  Then by Proposition \ref{v2gamma1}, $V^2(Q_Y)$  contains
 $V_{C^0}(c\l_2^0) \otimes V_{C^0}(\l_{r_0}^0) \otimes V_{C^1}(\l_1^1)$.
Restricting to $L_X'$ this becomes $V^1 \otimes \om_{l-1} \otimes \om_1.$ Now $\om_{l-1} \otimes \om_1 =  \om_l \oplus  (\om_1 + \om_{l-1}).$  Therefore
$$V^2  \supseteq  (V^1 \otimes  \om_l)  + (V^1 \otimes (\om_1 + \om_{l-1})).$$ As $V^1$ has a constituent with two nonzero labels, Corollary \ref{cover} now gives a contradiction.

Hence $\langle \l, \gamma_1 \rangle = 0$, and so $\l = c\l_2$. If $l\ge 5$ then the inductive hypothesis implies that $c\le 3$; and the same holds for $l=4$ by Lemma \ref{donnanewa5lem}.  \hal



The next result covers certain possibilities that only occur for $l = 4.$

\begin{lem}\label{2l32l4}  Assume $l = 4.$ Then $\mu^0 \ne 2\l_3^0$, $2\l_4^0$, $2\l_6^0$, $2\l_7^0$, $\l_1^0+2\l_2^0$ or 
$\l_9^0+2\l_8^0$. 
\end{lem}

\pf  Assume first that $\mu^0 = 2\l_3^0,2\l_4^0, 2\l_6^0,$ or $2\l_7^0.$   Identifying irreducible representations with their highest weights we first 
check that $2\l_3^0 = S^2(\l_3^0) - (\l_1^0 + \l_5^0) = S^2(\l_3^0) - (\l_1^0 \otimes \l_5^0) + \l_6^0.$
Restricting to $L_X'$ this is
$$(4002)+(0003)+(0040)+(1020)+(2021)+(2000)+(0202)+(2110)+(1101).$$
Similarly, $2\l_4^0 = S^2(\l_4^0) - (\l_2^0 + \l_6^0) -\l_8^0= S^2(\l_4^0) - (\l_2^0 \otimes \l_6^0) + (\l_1^0\otimes\l_7^0)-\l_8^0$ and restricting  $L_X'$ this is
$$(2022)+(0020)+(4011)+(3002)+(1102)+(1021)+(2111)+$$
$$(0004)+(0203)+(2001)+(1110)+(0130)+(2200)+(6000). $$ 

The remaining cases are duals of the above.  
It therefore follows that in each case   $\mu^0 \downarrow L_X'$ has a summand with at least 3 nonzero labels.
We first claim that $\mu^1 = 0.$ Suppose false.  We will see that in each case $V^1$
fails to be MF.

It follows from Proposition \ref{stem1.1.A}  that $\mu^1= d\l_1^1$, $d\l_4^1$, $\l_2^1$, or $\l_3^1.$
If $\mu^0 = 2\l_3^0$, then 
\[
\begin{array}{l}
(d000) \otimes ((2000) + (1101)) \supseteq (d100)^2, (000d) \otimes ((2110) + (1101)) \supseteq (111(d-1))^2,\\
 (0100) \otimes ((2110) + (1101)) \supseteq (2100)^2, \hbox{ and} \\
(0010) \otimes ((2110) + (1101)) \supseteq (1111)^2.
\end{array}
\]
  Therefore $V^1$ is not MF, establishing the claim for $\mu^0 = 2\l_3^0.$ The case $\mu^0 = 2\l_7^0$ follows by duality.
Similarly, if $\mu^0 = 2\l_4^0$, we have 
\[
\begin{array}{l}
(d000) \otimes ((4011) + (2111)) \supseteq ((d+2)111)^2,\\
 (000d) \otimes ((1110) + (1102)) \supseteq (111d)^2,\\
 (0100) \otimes ((1102) + (1021)) \supseteq (0112)^2,\\
(0010) \otimes ((1021) + (1102)) \supseteq (2011)^2.
\end{array}
\] 
 So the claim holds for $\mu^0 = 2\l_4^0$ and by duality $\mu^0 = 2\l_6^0.$
So we now have $\mu^1 = 0.$

We now study $V^2(Q_Y)$ and obtain a contradiction in each case.  First assume that $\mu^0 = 2\l_3^0.$
Then $V^2(Q_Y) $ contains $V_{C^0}(\l_2^0 + \l_3^0) \otimes V_{C^1}(\l_1^1).$   The first tensor factor
restricted to $L_X'$ is $(\wedge^2(0100) \otimes \wedge^3(0100)) - ((0100) \otimes \wedge^4(0100)).$
This decomposes as $(1111) +(0110)+(1030)+(1200)+(1001) + (3011)+(2010)+(0102)$
and tensoring with $(1000)$ we see that $V^2 \supseteq (0101)^3.$  On the other hand
from the first paragraph of the proof we see that at most one such factor arises from $V^1$.  This is  a contradiction.

Next assume $\mu^0 = 2\l_4^0.$  Here we see that  $V^2(Q_Y)$ contains $V_{C^0}(\l_3^0 + \l_4^0) \otimes V_{C^1}(\l_1^1).$ The first tensor factor restricts to $L_X'$ as $(\wedge^3(0100) \otimes \wedge^4(0100)) - (\wedge^2(0100) \otimes \wedge^5(0100)),$ which contains $(1112)+(1201)+(0111)+(1120)+(2011)^2.$  Tensoring with $(1000)$
we obtain $(1111)^6.$  From the first paragraph we see that at most $(1111)^4$ can arise from $V^1$,  so this is  a contradiction.

Now assume $\mu^0 = 2\l_7^0.$  Then $V^1$ is the dual of the restriction of $2\l_3$ given in the first paragraph of the proof.  On the other hand $V^2(Q_Y)$ contains $V_{C^0}(\l_6^0 + \l_7^0) \otimes V_{C^1}(\l_1^1).$ The first tensor factor restricts to $L_X'$ as the dual of $(\wedge^3(0100) \otimes \wedge^4(0100)) - (\wedge^2(0100) \otimes \wedge^5(0100)),$ which contains $(2111)+(1301)+(0211),$ 
so the tensor product of this with $(1000)$ contains $(1211)^3.$  At most one summand $(1211)$ can arise from 
$V^1$,  so we again have a contradiction.

Finally, assume $\mu^0 = 2\l_6^0.$  Then $V^1$ is the dual of the restriction of $2\l_4$ given in the first paragraph of the proof.  Then 
$V^2(Q_Y)$ contains $V_{C^0}(\l_5^0 + \l_6^0) \otimes V_{C^1}(\l_1^1).$  Restricting the first tensor factor to $L_X'$ yields $(\wedge^5(0100) \otimes \wedge^6(0100)) - (\wedge^4(0100) \otimes \wedge^7(0100))$ and one checks that this contains $(0011)+(1012)^2+(1020)+(1101)^2.$  The tensor product of this with $(1000)$ contains $(1011)^6.$ However $V^1$ yields at most $(1011)^3$, again a contradiction. 

It remains to consider the cases $\mu^0 =  \l_1^0 + 2\l_2^0$ or $2\l_8^0 + \l_9^0$. The latter case is ruled out by considering $V^*$, so assume $\mu^0 = \l_1^0 + 2\l_2^0$.  A Magma check shows that 
\[
V_{C^0}(\mu^0)\downarrow L_X’ = 2010 + 0102 + 2120 + 1001+ 1111 + 0110 +1200 +0301.
\]
We claim $\mu^1 = 0$.  Assume not. Then from Proposition \ref{stem1.1.A} we see that  
$\mu^1 = d\l_1^1$, $d\l_4^1$, $\l_2^1$, or $\l_3^1$. But we find that 
\[
\begin{array}{l}
(d000) \otimes ((2010) + (0110)) \supseteq (d110)^2, \\
(000d) \otimes ((0102) + (0110)) \supseteq (011d)^2, \\
(0100) \otimes ((2010) + (1200)) \supseteq (2110)^2, \\
(0010) \otimes ((1200) + (1111)) \supseteq (1210)^2,
\end{array}
\]
showing that  $V^1$ is not MF in each case, a contradiction.

Hence $\mu^1=0$. 
Observe that  $\la \l,  \g_1 \ra = 0$, since otherwise Lemma \ref{twolabs} together Proposition \ref{tensorprodMF} shows that 
$(V^*)^1$ is not MF. Hence $\l = \l_1+2\l_2$. However, a Magma check shows that $V_Y(\l_1+2\l_2)\downarrow X$ is not MF.
This contradiction completes the proof.
 \hal
 
\begin{lem}\label{muequalslambdai}  Assume $\mu^0 = \lambda_i^0$ for some $i$ with $1<i\le r_0$.  
If $\lambda$ or its dual is not in 
 Table $\ref{TAB2}$ or $\ref{TAB3}$ (of Theorem $\ref{MAINTHM}$), 
then  $\lambda  = \lambda_i + c\lambda_n.$
\end{lem}

\pf    If $i = r_0 =\frac{(l+1)l}{2} - 1$, then taking duals we see that all nonzero labels
on $\l^*$ appear in $(\mu^*)^0.$  Then the induction hypothesis forces $\l^* = \l_{l+2}$ or $\l_1 + \l_{l+2}$
and the assertion follows.   So from now on we assume $i < r_0$.

 First assume that $\mu^1 = 0.$  
If $\la \lambda, \gamma_1 \ra = 0$, then 
$\lambda = \lambda_i$ which is an example in Table \ref{TAB3}.  Now suppose $\la \lambda, \gamma_1\ra \ne 0.$
Set $\nu^0 = \lambda - \beta_i^0 - \cdots - \gamma_1.$ 
Using Lemma \ref{v2gamma1}, we have $V^2(Q_Y) \supseteq (\lambda_i^0 + \lambda_{r_0}^0) \otimes \lambda_1^1) + (\nu^0 \otimes \lambda_1^1).$ These two terms sum to $\lambda_i^0 \otimes \lambda_{r_0}^0 \otimes \lambda_1^1$
which restricts to $L_X'$ as $\wedge^i(\om_2) \otimes  \om_{l-1} \otimes \om_1. $  As $\om_{l-1} \otimes \om_1 = (\om_1 + \om_{l-1}) + \om_l$, the result is $(\wedge^i(\om_2) \otimes \om_l) + (\wedge^i(\om_2) \otimes (\om_1 +\om_{l-1})).$ 
Now combining Lemma \ref{twolabs} with Corollary \ref{cover} gives a contradiction.

Therefore, from now on we  assume $\mu^1 \ne 0.$  Another application of Lemma \ref{twolabs} together with 
Proposition \ref{tensorprodMF} implies that $\mu^1 = c\lambda_j^1.$  Taking duals, applying the induction hypothesis and Lemmas \ref{mu0=cl2} and \ref{2l32l4}, we see that $\mu^1 = \lambda_j^1$ or $c\lambda_l^1.$

We next claim that $\la \lambda, \gamma_1 \ra = 0.$  By way of contradiction assume $\la \lambda, \gamma_1 \ra > 0.$ First assume $\mu^1 = \lambda_j^1$.  Then 
$$V^2(Q_Y) \supseteq ((\l_i^0 + \l_{r_0}^0) \otimes (\l_1^1 + \l_j^1)) +   (\l_{i-1}^0 \otimes (\l_1^1 + \l_j^1)) + ((\l_i^0 + \l_{r_0}^0) \otimes \l_{j+1}^1)  + (\l_{i-1}^0 \otimes \l_{j+1}^1),$$ 
where we set $\l_{j+1}^1 = 0$ in case $j = l.$ Combining the first two summands and the last two summands we see that
 $$V^2(Q_Y) \supseteq ((\l_i^0 \otimes \l_{r_0}^0) \otimes (\l_1^1 + \l_j^1)) +   ((\l_i^0 \otimes \l_{r_0}^0) \otimes \l_{j+1}^1).$$
 And these can be combined to yield
$$V^2(Q_Y) \supseteq (\l_i^0 \otimes \l_{r_0}^0) \otimes (\l_1^1 \otimes \l_j^1).$$
Restricting to $L_X'$ we have
$$V^2 \supseteq \wedge^i(\om_2) \otimes \om_{l-1} \otimes \om_1 \otimes \om_j.$$
Expanding the tensor product of the middle two terms we find that
$$V^2 \supseteq (\wedge^i(\om_2) \otimes \om_j \otimes \om_l)  + (\wedge^i(\om_2) \otimes (\om_1 + \om_{l-1}) \otimes \om_j),$$ 
and Corollary \ref{cover} gives a contradiction.
The argument for the case $\mu^1 = c\lambda_l^1$ is essentially the same.

So from now we assume that $\la \lambda, \gamma_1 \ra = 0.$  Notice that if $\mu^1 = c\l_l$, then
$\l = \l_i + c\l_l$ as in the statement of the Proposition.
Therefore we now may assume  $\mu^1 = \l_j^1$ for $j < l.$  Then
$$V^2(Q_Y) \supseteq   (\l_{i-1}^0 \otimes (\l_1^1 + \l_j^1)) + ((\l_i^0 + \l_{r_0}^0) \otimes \l_{j+1}^1) + (\l_{i-1}^0 \otimes \l_{j+1}^1).$$
Combining the second two terms we have
$$V^2(Q_Y) \supseteq   (\l_{i-1}^0 \otimes (\l_1^1 + \l_j^1)) + ((\l_i^0 \otimes \l_{r_0}^0) \otimes \l_{j+1}^1).$$
Restricting to $L_X'$ we have
$$V^2 \supseteq \wedge^{i-1}(\om_2) \otimes (\om_1 + \om_j) + \wedge^i(\om_2) \otimes \om_{j+1} \otimes \om_{l-1}.$$
As $\om_{j+1} \otimes \om_{l-1} \supseteq (\om_j + \om_l) \oplus \om_{j-1} = \om_j \otimes \om_l$ we see that
$$V^2 \supseteq \wedge^{i-1}(\om_2) \otimes (\om_1 + \om_j) + \wedge^i(\om_2) \otimes \om_j \otimes \om_l.$$
Again, Lemma \ref{twolabs} and Corollary \ref{cover} give a contradiction, provided $i > 2.$   So assume $i = 2.$  Replacing $V$ by $V^*$ and repeating the above argument we see that $j = l-1$ or $l$. Therefore, $\lambda = \lambda_2 + \lambda_{n-1}$ or $\lambda_2 + \lambda_n.$  The former is on the list of examples in Theorem \ref{MAINTHM}, as is the dual of the latter.  \hal

\begin{lem}\label{mu0=l1+li}  Assume $\mu^0 = \l_1^0 + \l_i^0$ for $2 \le i \le 7.$  Then $\lambda$ is as in 
Tables $\ref{TAB2}$, $\ref{TAB4}$ of Theorem $\ref{MAINTHM}$.
\end{lem}

\pf  First assume that $\mu^1 = 0.$  If $\la \l, \gamma_1 \ra = 0$, then $\l$ is on the list of examples in  Tables \ref{TAB2}, \ref{TAB3}.  So suppose $\la \l, \gamma_1 \ra > 0.$  
Then by Lemma \ref{v2gamma1}, $V^2$ contains $((\l_1^0+\l_i^0)\downarrow L_X') \otimes \o_{l-1}\otimes \o_1$.
As $\om_{l-1} \otimes \om_1 = (\om_1 + \om_{l-1}) \oplus \om_l$, it follows that $V^2$ contains 
$(V^1\otimes \o_l) + ((\l_1^0+\l_i^0)\downarrow L_X' \otimes (\o_1+ \o_{l-1}))$. Now we obtain a contradiction 
in the usual way using Proposition \ref{tensorprodMF}, Corollary \ref{cover} and Lemma \ref{twolabs}. 

 So from now on we assume $\mu^1 \ne 0.$ 
As $V^1$ is MF, it follows from Lemma \ref{twolabs}  that
$\mu^1 = c\l_j^1$ for some $j$.  First suppose that $j < l-1.$  Consider $V^*$. 
By induction together with Lemma \ref{2l32l4}, we see that $c = 1.$    If $i\le l$, then $(V^*)^1(Q_Y)$ is the tensor product
of two modules and on restriction to $L_X'$, each has a summand with with at least two nonzero labels. Hence this is a contradiction by Proposition \ref{tensorprodMF}. 
Now assume $i > l.$ Then $7 \ge i > l$, so that $4 \le l \le 6.$  At this point Magma calculations show that $V^1$ is not MF, a contradiction.  

So we can now assume $j \ge l-1.$ If $j = l-1$, then Lemma \ref{mu0=cl2} applied to $V^*$ implies that $c \le 3.$  Arguing as in the above paragraph we get a contradiction.  Namely,  $(V^*)^1$ is not MF if $i \le l$
and if $i > l,$ then $l = 4,5,$ or $6$ and Magma checks show that $V^1$ is not MF.
Therefore we may now assume that $\mu^1 = c\lambda_l^1.$  We argue as in the first paragraph of the proof that $\la \l,\g_1\ra = 0$. 
If $c = 1$, then $\l = \l_1+\l_i+\l_n$, which is not possible by Lemma \ref{nolab1}. So  assume $c > 1.$

Consider $V^*.$  Suppose first that  $i =  l+1$. Then $\l^* = c\l_1 + \l_{r_0+1} + \l_n$ and we argue as above that 
$$(V^*)^2  \supseteq (S^c(\om_2) \otimes \om_l \otimes \om_l ) + (S^c(\om_2) \otimes (\om_1 + \om_{l-1}) \otimes \om_l).$$
Using Proposition \ref{stem1.1.A} we see  that  the second summand is not MF, which is a contradiction by Corollary \ref{cover}.  Next,  assume $i > l+1$, so $l \le 5$ (as $i\le 7$).  Then $(\mu^*)^0 = c\l_1^0 + \l_k^0$ for  
some $k$ with $2 < k \le r_0$, and the inductive list of examples implies that $k = r_0$.  If $l = 5$ then the induction hypothesis implies that $c = 2$ or $3$ and Magma can be used
to show that $(V^*)^1$ is not MF.  Now suppose $l = 4.$  Then $\mu^0 = \l_1^0  + \l_6^0$ or $ \l_1^0  + \l_7^0$, and 
$\mu^1 = c\l_4^1.$ Then $\mu^0 \downarrow L_X'$ contains $(1010) + (1002)$ or $(1110)+(1102)$,  
and tensoring with $(000c)$ we get $(101c)^2$ or $(111c)^2$.  Therefore  $V^1$ is not MF, a contradiction.  Hence $i \le l.$ 

We now have $(\mu^*)^0 = c\l_1^0$ and $(\mu^*)^1 = \l_k^1 + \l_l^1,$  where $k = l-i+1.$  Therefore, 
$(V^*)^1 = S^c(\om_2) \otimes (\om_k + \om_l).$   If $c > 3$, then the first tensor
factor has an irreducible summand with at least two nonzero labels and hence $(V^*)^1$
is not MF, a contradiction.  Therefore, $c = 2.$

Here too, we claim that $(V^*)^1$
is not MF.  First note that $(\mu^*)^0 \downarrow L_X' = 2\om_2 \oplus \om_4$.  If $k = 1$  then Theorem \ref{LR} implies that tensoring each of $2\om_2$ and $\om_4$ with
$\om_1 + \om_l$ produces a summand $\om_1+ \om_3$.   Finally, for $k>1$, Theorem \ref{LR} implies that
$(V^*)^1$ contains $ (\om_2 + \om_{k+1})^2.$ \hal

\begin{lem}\label{mu0=l1+ln-+2-i}  Assume $\mu^0 = \l_1^0 + \l_{r_0+2-i}^0$ with $2 \le i \le 7$ and $r_0+2-i>7$.  
Then $\l =  \l_1 + \l_{r_0+2-i}.$  
\end{lem}

\pf    The argument in the first few paragraphs of the proof of Lemma \ref{mu0=l1+li} shows that $\mu^1 = c\l_j^1$ for some $j$, and we can assume that $c\ne 0$. Note also that application of the induction hypothesis to $V^*$ gives $\la \l,\g_1\ra=0$.
 
We will work with $V^*$  which has highest weight $\l^* = c\l_k + \l_t + \l_n$, where $1 \le k \le l$ and
$l+2 \le t \le l+7$.  Therefore, $(V^*)^1(Q_Y) = ( c\l_k^0 + \l_t^0) \otimes (\l_l^1).$
 If $l \ge 6,$ then  we see from the inductive list that the restriction to $L_X'$ of the first factor of the tensor product is not MF.

This leaves the cases where $l = 4$ or 5.  Suppose that $l = 5$ so that $n = 20$ and $r_0=14$.  Here $7 \le t \le 12$ and from the inductive list of examples we see that $c =k=1$ and $\l^* = \l_1 + \l_t + \l_{20}$.  A Magma computation shows that in each case 
$(V^*)^1$ fails to be MF. 

Finally, assume that $l = 4.$  Here we have $\l^* = c\l_k + \l_t + \l_{14}$ and from the first paragraph $t = 6$ or $7$.   
Also $k \le 4.$  Then 
$(V^*)^1(Q_Y) = (c\l_k^0 + \l_t^0) \otimes (\l_4^1)$ and from the inductive list of examples we 
see that  $c = k = 1.$   A Magma computation shows that none of these are MF on restriction to $L_X'$.  \hal

\begin{lem}\label{mu0=l2+l3} Assume $\mu^0 =  \l_2^0 + \l_3^0$ or $\l_{r_0-2}^0+\l_{r_0-1}^0.$ Then $\l = \l_2 + \l_3$, as in Table $\ref{TAB2}$ of Theorem $\ref{MAINTHM}$.
\end{lem}

\pf  First assume $\mu^0 =  \l_{r_0-2}^0+\l_{r_0-1}^0.$  If $l \ge 5$, then taking duals we see that $(\mu^*)^0$ does not satisfy the induction hypothesis.  This argument also yields a contradiction if  $l = 4$ unless   $\l = \l_7 + \l_8$.
And if $l = 4$ and $\l = \l_7 + \l_8,$ a Magma calculation shows that $V \downarrow X$ is not MF.

So now assume $\mu^0 =  \l_2^0 + \l_3^0.$  Now  $V^1(Q_Y) = V_{C^0}(\mu^0) \otimes V_{C^1}(\mu^1)$ and $V_{C^0}(\mu^0) \downarrow L_X' =  (\wedge^2(\om_2) \otimes \wedge^3(\om_2)) -
 (\om_2 \otimes \wedge^4(\om_2)).$  Since all the nonzero labels in highest weights of the relevant modules occur in the first several nodes, we can apply Theorem \ref{LR} together with some Magma checks to show that $V_{C^0}(\mu^0) \downarrow L_X' \supseteq (30110 \ldots  0) + (11110 \ldots  0).$ 
 
If $\mu^1\ne 0$ then arguing as before, we have $\mu^1 = c\l_j^1$ for some $j$, and 
Lemma \ref{abtimesce} shows that $V^1$ contains
  $(2\om_1 + \om_3 +\om_4 + (c-1)\om_j+\om_{j+1})^2$,  $(2\om_1 + \om_3 +c\om_4 + \om_5)^2,$ $(2\om_1 + c\om_3 +2\om_4)^2$, $(2\om_1 + (c-1)\om_2 +2\om_3 +\om_4)^2,$ or $((c+1)\om_1 + \om_2 +\om_3 +\om_4)^2,$ according as $j > 4, j = 4, j = 3, j = 2,$ or $j = 1, $ where the terms $\om_{j+1}, \om_5$ 
  are omitted in the first two cases if $j = l$ or $l = 4$, respectively.   Therefore $\mu^1 = 0.$
 
Now if $ \langle \l, \g_1 \rangle = 0,$  then the conclusion holds,   so suppose this is not the case.
Then  $(V^*)^1(Q_Y)= V_{C^0}(d\l_{l+1}^0) \otimes V_{C^1}(\l_{l-2}^1 + \l_{l-1}^1).$ Now  
Lemma \ref{twolabs} shows that $(V^*)^1$ is not MF, a contradiction.  \hal
 
 \begin{lem}\label{mu0=cl1+ll} Assume $\mu^0 =  c\l_1^0 + \l_{r_0}^0$ or $\l_1^0 + c\l_{r_0}^0.$   Then   
$\l = c\l_1 + \l_{r_0}$, and $c \le 3$ if $l\ge 5$. 
\end{lem}
 
 \pf  Applying the induction hypothesis to $(V^*)^1$, we see that $\mu^0 =  c\l_1^0 + \l_{r_0}^0$, $ \langle \l, \g_1 \rangle = 0$ and that
 $\mu^1 = 0$ or $\mu^1 = \l_l.$ In the former case the conclusion holds.  So suppose $\mu^1 = \l_l.$ Then $(\mu^*)^0 = \l_1^0 + \l_{l+2}^0$ and the induction hypothesis implies
 that either $l+2 \le 7$ or $l+2 \ge \frac{(l+1)l}{2} -6.$  Therefore $l \le 5.$  Now apply Lemma \ref{mu0=l1+li} to $V^*$ to obtain 
a contradiction.  \hal

 \begin{lem}\label{mu0=cl1} Assume $\mu^0 =  c\l_1^0$ with $c \ge 1$, or $\mu^0 = c\l_{r_0}^0 $ with $c > 1.$ Then, up to duals, either $\l$ is as in the conclusion of Theorem $\ref{MAINTHM}$, or $\l$ is one of  $c\l_1 + b\l_n$, $ \l_1 + \l_{r_0+1}$, 
$\l_1 + \l_i + \l_n$ ($i\le 7)$ or $c\l_1 + \l_i$ $(i > r_0+1)$.
\end{lem}

\pf If  $\mu^0 = c\l_{r_0}^0$ with $c > 1,$  then the induction hypothesis implies that $(V^*)^1$ is not MF.  Therefore, from now on we assume that  $\mu^0 = c\l_1^0$. 

We first consider the case where $\langle \l, \g_1 \rangle \ne 0.$  Applying the induction hypothesis to
 $(V^*)^1$ we see that $\langle \l, \g_1 \rangle = 1$ and $\mu^1 = 0$ or $ \l_l^1$ (with $l\le 6$ in the latter case).
Then by Proposition \ref{v2gamma1},
  $$V^2(Q_Y) \supseteq V_{C^0}(c\l_1^0) \otimes  V_{C^0}(\l_{r_0}^0) \otimes V_{C^1}(\l_1^1 + \mu^1).$$
 Suppose $\mu^1 = 0.$ Then restricting the above  to $L_X'$ this becomes $S^c(\om_2) \otimes  \om_{l-1} \otimes \om_1$, 
 which equals $S^c(\om_2) \otimes \om_l$ + $S^c(\om_2) \otimes  (\om_1 + \om_{l-1}).$
 If $c \ge 3$ then the first tensor factor of the second summand has an irreducible constituent with highest weight having at least  two nonzero labels by Lemma \ref{twolabs}, so this summand fails to be MF by Proposition \ref{stem1.1.A}, giving a contradiction by Corollary \ref{cover}.  Now suppose  $c = 2.$ Then Theorem \ref{LR} implies that $S^c(\om_2) \otimes (\om_1 + \om_{l-1})$ contains $ (\o_1 + \o_3 + \o_l)^2,$ again a contradiction.  And if $c =1$ then $\l = \l_1 + \l_{r_0+1}$ which is one of the exceptions listed in the conclusion. 
  
Now suppose $\mu^1 = \l_l^1$ (still with $\langle \l, \g_1 \rangle \ne 0$).  Here $V^2(Q_Y)$ has  additional summands,  $V_{C^0}(c\l_1^0 + \l_{r_0}^0) \otimes 0$ and $V_{C^0}((c-1)\l_1^0) \otimes 0.$  Combining these and adding the result to the above we see   that 
$$V^2(Q_Y) \supseteq V_{C^0}(c\l_1^0) \otimes  V_{C^0}(\l_{r_0}^0) \otimes V_{C^1}(\l_1^1) \otimes V_{C^1}(\l_l^1).$$
 Restricting to $L_X'$ this becomes
 $$S^c(\om_2) \otimes \om_{l-1} \otimes \om_1 \otimes \om_l.$$
Combining the middle two terms   gives 
 $$(S^c(\om_2) \otimes \om_l \otimes \om_l) + (S^c(\om_2) \otimes (\om_1 + \om_{l-1}) \otimes \om_l).$$
Proposition \ref{stem1.1.A} shows that
 the second term is not MF, giving a contradiction by Corollary \ref{cover}.  

Hence we can now assume that $\langle \l, \g_1 \rangle = 0.$
 
  Suppose $c \ge 3.$  Then $S^c(\om_2) $ contains $c\om_2$  and $ (c-2)\om_2 + \om_4.$  As $V^1$ is MF, it follows from Proposition \ref{stem1.1.A} that $\mu^1 = d\l_j^1$ for some $j.$  Now taking duals we see from the inductive
 hypothesis and Lemma \ref{2l32l4} that either $d = 1$ or $j \in \{l-1,l \}.$
If $4 \le j < l$, then applying Lemma \ref{aibjnotzero}
 we see that the tensor product of $\mu^1\downarrow L_X'$ with each of the summands contains $(c-1)\om_2 +(d-1)\om_j + \om_{j+2}$ (omit the last term if $j = l-1$) and hence
 $V^1$ is not MF.  Similarly, for $j = 2$ or $3$  using either Lemma \ref{aibjnotzero} or easy weight
 considerations.  Therefore, for $c \ge 3$ we can assume $\mu^1 = \l_1^1$ or $d\l_l^1.$  In either case the
 conclusion of the lemma holds. 
 
 Now assume $c \le 2.$  First suppose that $\mu^1$ is not of the form $d\l_j^1.$  Then consideration of $(V^*)^1$ using the induction hypothesis implies that
 $(\mu^*)^0$ is one of the following:
\[
\begin{array}{l}
\l_2^0 + \l_3^0, \\
\l_1^0 + \l_i^0 \hbox{ with } i \le {\rm max}(7,l), \\
\l_1^0+2\l_2^0, \\
a\l_1^0+\l_2^0\,(a\ge 2).
\end{array}
\]
The first and third possibilities are excluded by Lemmas \ref{mu0=l2+l3} and \ref{2l32l4} respectively.
Next suppose $(\mu^*)^0 = \l_1^0 + \l_i^0$ with $i \le {\rm max}(7,l)$.  If $c = 1$, then this is listed in the conclusion.
And if $c = 2$ then we obtain a contradiction by applying Lemma \ref{mu0=l1+li} to $V^*$.
Now assume that $(\mu^*)^0 = a\l_1^0 +\l_2^0$. 
Observe that $a\l_1^0 +\l_2^0 = ((a+1)\l_1^0\otimes \l_1^0) - (a+2)\l_1^0$, which restricts to $L_X'$ as
$(S^{a+1}(\o_2)\otimes \o_2)-S^{a+2}(\o_2)$. 
By Lemmas \ref{cfssc} and \ref{LR_om2},
\[
\begin{array}{ll}
S^{a+1}(\o_2)\otimes \o_2& \supseteq \left(((a-1)\o_2+\o_4)\oplus ((a-3)\o_2+2\o_4)\right) \otimes \o_2 \\
             & \supseteq (\o_1+(a-2)\o_2+\o_3+\o_4) \oplus ((a-2)\o_2+2\o_4)^2.
\end{array}
\]
The composition factors of $S^{a+2}(\o_2)$ of $S$-value at least $a$ are $(a+2)\o_2$, $a\o_2+\o_4$, $(a-2)\o_2+2\o_4$, all with multiplicity 1 (again by Lemma \ref{cfssc}). Hence $(a\l_1^0 +\l_2^0)\downarrow L_X'$ contains
$(\o_1+(a-2)\o_2+\o_3+\o_4) \oplus ((a-2)\o_2+2\o_4)$. Tensoring this with $c\o_l$ gives 
$((a-2)\o_2+\o_3+\o_4+(c-1)\o_l)^2 \subseteq (V^*)^1$, a contradiction.

 It remains to consider the case $\mu^1 = d\l_j^1$ (still with $c\le 2$).  As above this implies that either $d = 1$ or $j \in \{l-1,l \}.$
 
 If $d = 1$ or if $j = l$, then $\l$ is one of the cases listed in the conclusion. So we are left with the case $\mu^1 = d\l_{l-1}^1$ and $d > 1.$  The induction hypothesis applied
 to $(V^*)^1$ implies that $d \le 5$.  
 
 If $c = 2$, then we argue that in each case $V^1 = S^2(\om_2) \otimes d\om_{l-1}$ fails to be MF.  Indeed, $S^2(\om_2) = 2\om_2 \oplus \om_4$  and Lemma \ref{aibjnotzero} implies that $V^1$ contains $(\om_2 + (d-1)\om_{l-1})^2$, a contradiction. 
 
 Finally assume that $c = 1$ and consider $V^*.$  Here we check that $V_{C^0}((\mu^*)^0) \downarrow L_X'$
contains $ (2\om_2+\om_4) \oplus (\om_1+\om_2+\om_5)$ if $d=2$,  and contains 
$(\om_1+\om_2+\om_4+ \om_5) \oplus (2\om_1+\om_2+\om_3+ \om_5)$ if $d=3$ (delete the $\om_5$ terms  if $l = 4$.)
Then we find that $(V^*)^1$ contains $(\om_1 + \om_2 +\om_4)^2$ or $(\om_1+\om_2+\om_3+ \om_5)^2$ in the respective cases. And if $d=4$ or 5, then $l=4$ and we  see that $V^1$ is not MF, using the decomposition of $(d\l_{l-1}^1)\downarrow L_X'$ given by Table \ref{comple}. This is a final contradiction.  \hal
 
 \begin{lem}\label{age2lem} Assume  $a \ge 2$ and $\mu^0 = a\l_1^0+\l_2^0$ or $\l_{r_0-1}^0+a\l_{r_0}^0$.  Then $\l = a\l_1 +\l_2$ and $a \le 3.$  
\end{lem} 

\pf   Inductively we have $a \le 3$ if $l >4.$  If $\mu^0 = \l_{r_0-1}^0+a\l_{r_0}^0,$ then taking duals we see that $(\mu^*)^0 \downarrow L_X'$ is not MF, a contradiction.  So from now on assume $\mu^0 = a\l_1^0+\l_2^0.$
In order for $V^1$ to be MF we must have $\mu^1  = d\l_i^1$ for some $i$. 
Now pass to $V^*$.  Here $(\mu^*)^1 = \l_{l-1}^1+a\l_l^1$.  In order for $(V^*)^1$ to be MF it is necessary that all composition factors of $(\mu^*)^0 \downarrow L_X'$ to have at most one nonzero label.   
By Lemma \ref{twolabs} the only possibilities are $(\mu^*)^0 = 0, \l_1^0$ or $2\l_1^0$.  In the last case Lemma \ref{aibjnotzero}  implies that $(V^*)^1 \supseteq (\om_2 + a\om_l)^2,$  a contradiction.  Therefore $(\mu^*)^0 = 0$ or $ \l_1^0$.

Suppose that $(\mu^*)^0 = \l_1^0$ so that $V^1\downarrow L_X' = ((a10\ldots  0) \downarrow L_X') \otimes (0\ldots  01).$  We claim that this is not MF.  From the proof of Lemma \ref{mu0=cl1} we see that the first tensor factor contains $(1(a-2)110 \ldots  0)$ and $(0(a-2)020 \ldots  0),$ so tensoring with $(0 \ldots  01)$ we get $(0(a-2)110 \ldots  0)^2$.  Thus the claim holds and so $(\mu^*)^0 \ne \l_1^0.$

Finally suppose $(\mu^*)^0 = 0.$   If $ \la \l, \g_1 \ra  \ne 0$, then using Lemma \ref{twolabs} and Proposition \ref{tensorprodMF} we see that $(V^*)^1$ is not MF, a contradiction.  Therefore $\l = a\l_1 + \l_2$, so the 
 result holds unless $l = 4.$  So now assume $l = 4$ and $a \ge 4.$  The usual arguments show that $V^2 = S^a(0100) \otimes (0100) \otimes (1000).$  As $a \ge 4,$ $S^a(0100) \supseteq (0(a-2)01) + (0(a-4)02).$  We then check that each summand tensored  with $(0100) \otimes (1000)$ contains $(0(a-3)01)^2.$  Therefore $V^2$ contains $(0(a-3)01)^4.$  We claim that only two such summands can arise from $V^1$.
The factor $(0(a-3)01)$ in $V^2$  can only arise from factors $(1(a-3)01),$ $(0(a-3)00),$ and $(0(a-4)11)$
in $V^1.$  As $V^1$ is MF it will suffice to show that $(0(a-3)00)$ does not appear in $V^1$.

We have $a\l_1^0+\l_2^0 = (S^{a+1}(\l_1^0) \otimes \l_1^0) - S^{a+2}(\l_1^0)$. 
 So it will suffice to show that $S^{a+1}(0100) \otimes (0100)$ does not contain $(0(a-3)00)$. Towards this end we recall the proof of Lemma \ref{a4(a0...01)} where it was noted that the composition factors of  $S^{a+1}(0100)$  have highest weights of the form $(0x0y)$ subject to $x + 2y = a+1$. On the other hand Lemma \ref{LR_om2} shows that $(0x0y) \otimes (0100) \supseteq (0(a-3)00)$ only if $(0x0y) = (0(a-4)00)$, which does not satisfy the equality  $x + 2y = a+1$.  Therefore $(0(a-3)00)$ does not appear in $V^1$, completing the proof.  \hal

 \begin{lem}\label{newpossib} We have $\mu^0 \ne \l_2^0+\l_{r_0-1}^0$. 
\end{lem} 

\pf  Write $\mu^1 = a_1\l_1^1 + \cdots +a_l\l_l^1$. For $V^*$, $(\mu^*)^0$ has labelling $(a_l, a_{l-1}, \ldots , a_1,x,0,1,0, \ldots , 0)$, where the last nonzero entry appears at node $l+3 \le r_0+2.$ The induction hypothesis  implies that $(a_l, a_{l-1}, \ldots , a_1,x) = (0\ldots  0)$ or $(10\ldots  0).$ Returning to $V$ we see that $\la \l,\g_1\ra= 0$, and $\mu^1 = 0$ or $\l_l^1.$ In the former case we can immediately apply Lemma \ref{newestlem}, so assume $\la \l,\g_1\ra= 0$ and $\mu^1 = \l_l^1.$

Then $(\mu^*)^0 = \l_1^0 + \l_{l+3}^0$ and $(\mu^*)^1= \l_{l-1}^1.$  If  $l \ge 5$, then $l+3 > 7$ and $l +3 < r_0 -5$ contradicting the induction hypothesis.  Therefore $l = 4.$ At this point a Magma check shows that $(V^*)^1$ is not MF, a contradiction.  \hal

 \begin{lem}\label{cl1+li}  Suppose that $\l = c\l_1 + \l_i$  with $c > 1$.
 \begin{itemize}
\item[{\rm (i)}]  If $i >r_0$, then either  $i =r_0+2$ or $i = n.$
 \item[{\rm (ii)}]  If $i = n$ then $c \le 3.$
\end{itemize}
 \end{lem}
 
 \pf  (i) Recall that Lemma \ref{mu0=cl1} shows that $i \ne r_0+1.$ Assume that $r_0+2 < i < n$ and set $j = i - (r_0+1).$ Then $1 < j < l.$ We will show that $V^1 = S^c(010\ldots 0)  \otimes \om_j$ is not MF, which is a contradiction.   Note that $S^c(010\ldots 0) \supseteq (0c0\ldots 0) + (0(c-2)01\ldots 0).$  
 
 Suppose that $4 \le  j < l-1.$
 In this case Lemma \ref{aibjnotzero} shows  that $S^c(010\ldots 0)  \otimes \om_j \supseteq ((c-1)\om_2 + \om_{j+2})^2$, where we delete the last term if $j = l-1$.   
  If $j = 3$, then using  Lemmas \ref{aibjnotzero}, \ref{abtimesce} and \ref{LR}, we see that
 both $(0c0\ldots 0) \otimes \om_3$ and  $(0(c-2)01\ldots 0) \otimes \om_3$ contain $(0(c-1)0010\ldots 0)$,
 where we omit the last term if $l = 4.$
 And if  $j = 2$, Lemma \ref{LR_om2} shows that $(0c0\ldots 0) \otimes \om_2$ and  $(0(c-2)01\ldots 0) \otimes \om_2$ both contain $(0(c-1)010\ldots 0).$    So in each case we have a contradiction. 
 
  (ii) Assume  $\l = c\l_1 + \l_n.$  By way of contradiction assume that $c \ge 4.$ By Lemma \ref{cfssc}, 
  $S^c(\om_2)$ contains irreducible summands of highest weights $\nu_1 = (c-2)\om_2 + \om_4, \nu_2 = (c-4)\om_2 + 2\om_4$, and $\nu_3 = (c-3)\om_2 + \om_6,$ where in the last weight we omit  $\om_6$ if $l = 4.$
  And applications of Lemma \ref{aibjnotzero} show that the tensor product  of each of $\nu_1, \nu_2,$ and $\nu_3$ with $\om_{l-1}$ contains the irreducible of highest weight $\nu = (c-3)\om_2 + \om_4$.   
It follows that $ S^c(\om_2) \otimes \om_{l-1} \supseteq ((c-3)\om_2 + \om_4)^3.$   On the other hand
working in $A_n$ we have $c\l_1 \otimes \l_n = (c\l_1 + \l_n) + (c-1)\l_1$ and another application of Lemma \ref{cfssc}  shows that $(c-3)\om_2 + \om_4$ occurs in $S^{c-1}(\om_2)$
with multiplicity $1$. It therefore follows that $(c\l_1 + \l_n) \downarrow X$ is not MF.   \hal

At this point, Lemmas \ref{mu0posso2} -- \ref{cl1+li} show that either $\l$ is as in Tables \ref{TAB2}, \ref{TAB3} of Theorem \ref{MAINTHM}, or $\l$ is one of the following possibilities, up to duals:
\begin{enumerate}[]
\item{\rm (1)} $c\l_1+b\l_n$ ($b,c\ge 1$)
\item{\rm (2)} $c\l_1+\l_i$ ($c\ge 1$, $i>r_0$)
\item{\rm (3)} $c\l_1+\l_{r_0}$ ($c\le 3$ if $l\ge 5$)
\item{\rm (4)} $\l_1+\l_{r_0+2-i}$ ($2\le i\le 7$)
\item{\rm (5)} $\l_i+c\l_n$ ($i\le r_0$, $c\ge 1$).
\end{enumerate}

Consider case (1). If $b=1$ then $c\le 3$ by Lemma \ref{cl1+li}, and so $\l$ is in Table \ref{TAB2} of Theorem \ref{MAINTHM}. Hence we can assume that $b,c\ge 2$. But then $V\downarrow X$ is not MF by Lemma \ref{nolab2}, a contradiction.

Now suppose $\l$ is as in (2). If $c\ge 2$, then $i=r_0+2$ or $n$ by Lemma \ref{cl1+li}. In the first case $\l^* = \l_l+c\l_n$, for which $V\downarrow X$ is shown to be non-MF by Lemma \ref{(s+2)cge2}; and in the second, $c\le 3$ by 
 Lemma \ref{cl1+li}, and so $\l$ is as in Table \ref{TAB2}. Hence $c=1$ and $\l = \l_1+\l_i$ with $i>r_0$. 
If $i \ge n-5$ then $\l$ is as in Table \ref{TAB2}, so we can assume $i\le n-6$. Then  
$V\downarrow X$ is not MF by Lemma \ref{nolab3}.

If $\l$ is as in (3), then $\l^* = \l_{l+2}+c\l_n$, and $V\downarrow X$ is not MF by Lemma \ref{cl1+ls}.
Likewise, $V\downarrow X$ is not MF for $\l$ as in (4), by Lemma \ref{nolab4}.

Finally, suppose $\l = \l_i+c\l_n$ with $i\le r_0$, $c\ge 1$, as in (5). Then $\l^* = c\l_1+\l_j$ where $j=n-i+1$. If $j > r_0$ this is as in case (2), already dealt with. And if $j\le r_0$ then Lemmas \ref{mu0posso2} and \ref{age2lem} applied to $(\mu^*)^0$ imply that either $\l^*$ is as in Table \ref{TAB2}, or it is as in one of cases (3) and (4) above, hence already dealt with.

At this point the proof of Theorem \ref{deltai} is complete.

\chapter{The case $\d = \omega_1+\omega_{l+1}$}\label{o1ol1}

In this chapter we prove the following result.

\vspace{4mm}
\begin{theor}\label{om1oml+1}
Let $X = A_{l+1}$ ($l\ge 1$), let $W = V_X(\om_1+\om_{l+1})$ and let $Y = SL(W)$. Suppose $V = V_Y(\l)$ is an irreducible $Y$-module such that $V\downarrow X$ is MF. Then up to duals, $\l$ is $\l_1$, $2\l_1$, $\l_2$, $\l_3$ or $3\l_1\,(l=1)$, as in Table 
$\ref{TAB1}$ of Theorem $\ref{MAINTHM}$.
\end{theor}

Before the proof of the theorem, here is a preliminary result that we will need.

\begin{lemma}\label{mults}
The following hold for $A_l$-modules, $l\ge 2$:
\begin{itemize}
\item[{\rm (i)}] $\wedge^2(\om_1+\om_l) = (\om_1+\om_l) \oplus (\om_2+2\om_l)\oplus (2\om_1+\om_{l-1})$
\item[{\rm (ii)}] $S^2(\om_1+\om_l) = (2\om_1+2\om_l)\oplus (\om_1+\om_l)\oplus (\om_2+\om_{l-1}) \oplus 0$
\item[{\rm (iii)}] If $a\ge 2,\,b\ge 1$ then $((a-2)\om_1+\om_2))\otimes \om_l \otimes b\om_1$ contains  $((a+b-3)\om_1+\om_2)^2$.
\item[{\rm (iv)}] $2\om_l \otimes (\om_1+\om_l)$, 
and $2\om_l\otimes (\om_2+\om_{l-1})$ both have a composition factor $\om_1+\om_{l-1}+\om_l$; if $l\ge 3$, 
so do $\om_{l-1}\otimes (\om_2+2\om_l)$ and $\om_{l-1}\otimes (2\om_1+\om_{l-1})$; and if $l=2$ then $10\otimes 30$ contains $21$.
\end{itemize}
\end{lemma}

\pf Parts (i) and (ii) follow from Lemma \ref{wedge2adj}. Parts (iii), (iv) are proved using Corollary \ref{tensorwithoml}, \ref{LR_om2} and Proposition \ref{pieri}. \hal

\vspace{4mm} 
We now embark on the proof of Theorem \ref{om1oml+1}. The proof goes by induction on $l$, the case $l=1$ being covered by Chapter \ref{casea2}. 
 Let $X$, $\d$, $W$ and $V = V_Y(\l)$ be as in the hypothesis of the theorem. Assume for a contradiction that $\l$ is not as in the conclusion.

The Levi subgroup $L_Y' = C^0C^1C^2$, and as $L_X'$-modules we have 
\begin{equation}\label{wemb}
W^1(Q_X)\cong \om_1, \;W^2(Q_X)\cong (\om_1+\om_l)\oplus 0,\; W^3(Q_X) \cong \om_l,
\end{equation}
so that $C^0\cong C^2 \cong A_l$ and $C^1 \cong A_{(l+1)^2-1}$. As usual we let $\mu^i$ be the restriction of $\l$ to $T_Y\cap C^i$, so that $V^1(Q_Y) \downarrow L_Y' = V_{C^0}(\mu^0) \otimes V_{C^1}(\mu^1) \otimes V_{C^2}(\mu^2)$,  and we write
$V^i = V^i(Q_Y)\downarrow L_X'$.

\begin{lemma}\label{mu10}
We have $\mu^1 = 0$.
\end{lemma}

\pf Suppose $\mu^1\ne 0$. Now $V^1$ is MF. In particular $V_{C^1}(\mu^1) \downarrow L_X'$ is MF, and by (\ref{wemb}), $L_X' = A_l$ is embedded in $C^1$ via the representation $(\om_1+\om_l)\oplus 0$. 
So $L_X'$ embeds into a maximal Levi subgroup of $C^1$, and all composition factors of this Levi subgroup on $V_{C^1}(\mu^1)$  must have  MF restriction  to $L_X'$. Hence, considering levels within $C^1$ for the Levi subgroup, and using the inductive hypothesis, we see that 
\[
\mu^1 = \l_1^1, 2\l_1^1,\l_2^1,\l_3^1 \hbox{ or }3\l_1\,(l=2).
\]
By Lemma \ref{mults}(i,ii), both $\wedge^2(\om_1+\om_l)$ and $S^2(\om_1+\om_l)$ have a composition factor $\om_1+\om_l$; so does $\wedge^3(\om_1+\om_l)$ (see Proposition \ref{wedgeadj}). Hence $\wedge^2((\om_1+\om_l)\oplus 0)$,  
$S^2((\om_1+\om_l)\oplus 0)$ and $\wedge^3((\om_1+\om_l)\oplus 0)$ are all non-MF; and so is $S^3(11+00)$ for $l=2$. It follows that $\mu^1 = \l_1^1$.

Suppose $\mu^0 \ne 0$. Then by Proposition \ref{tensorprodMF}, $\mu^0 = a\l_i^0$ for some $i$ and some $a\ge 1$. Hence
\[
\begin{array}{ll}
V^1 & \supseteq a\om_i \otimes ((\om_1+\om_l)\oplus 0) \\
        & = a\om_i \otimes \om_1 \otimes \om_l \\
        & \supseteq ((\om_1+a\om_i) \oplus  ((a-1)\om_i+\om_{i+1}))\otimes \om_l. \\
\end{array}
\]
By Corollary \ref{pieri} this contains $(a\om_i)^2$, contradicting the fact that $V^1$ is MF.

Hence $\mu^0 = 0$, and similarly $\mu^2=0$. So $V^1 = (\om_1+\om_l)\oplus 0$. In $V^2$, the weight 
$\l -\b_{l+1}- \b_{l+2}$ affords a summand 
\[
\om_l \otimes \wedge^2((\om_1+\om_l)\oplus 0).
\]
By Lemma \ref{mults}(i),   $\wedge^2((\om_1+\om_l)\oplus 0)$ has composition factors $(\om_1+\om_l)^2$ and $\om_2+2\om_l$; and by Corollary \ref{pieri} the tensor product of $\om_l$ with each of these has a composition factor $\om_1+2\om_l$.
Hence $V^2$ contains $(\om_1+2\om_l)^3$, whereas by Corollary \ref{V^2(Q_X)} only one such composition factor is in the first level of a composition factor of $V^1$. This contradicts Proposition \ref{induct}. \hal

\begin{lemma}\label{mu0non}
Either $\mu^0 \ne 0$ or $\mu^2\ne 0$.
\end{lemma}

\pf Suppose false; then using Lemma \ref{mu10} we have $\mu^0=\mu^1=\mu^2=0$, and so $V^1 = 0$. Let $\g_1 = \b_{l+1}$  and $\g_2 = \b_{n-l}$ (the nodes not in $C^0,C^1$ or $C^2$), and let $x=\la \l,\g_1\ra$, $y=\la \l,\g_2\ra$.

By duality we may assume that $x \ne 0$. If also $y\ne 0$, then in $V^2$ the weights $\l-\g_1$ and $\l-\g_2$ each affords $((\om_1+\om_l)\oplus 0) \otimes \om_l$, which by Corollary \ref{pieri} contains $\om_l^2$; hence $V^2 \supseteq \o_l^4$, which contradicts Corollary \ref{cover}.
Hence $y=0$. It follows that $V^2$ is afforded by the weight $\l-\g_1$, so that by Corollary \ref{pieri},
\[
V^2 = \om_l \otimes ((\om_1+\om_l)\oplus 0) = (\om_1+2\om_l)\oplus (\om_1+\om_{l-1})\oplus \om_l^2.
\]

If $x\ge 2$ then in $V^3$ we have the following summands:
\[
\begin{array}{|l|l|}
\hline
\hbox{summand of }V^3 & \hbox{afforded by } \\
\hline
2\om_l\otimes S^2((\om_1+\om_l)\oplus 0) & \l-2\g_1 \\
\om_{l-1}\otimes \wedge^2((\om_1+\om_l)\oplus 0) & \l-\b_l-2\g_1-\b_{l+2} \\
\hline
\end{array}
\]
By Lemma \ref{mults}(i,ii,iv), between them these summands contain $(\om_1+\om_{l-1}+\om_l)^4$. However  $(\om_1+\om_{l-1}+\om_l)$ only appears in the first level of two of the composition factors of $V^2$, by
 Corollary \ref{V^2(Q_X)}, so this contradicts Proposition \ref{induct}.

Hence $x=1$ and so $\l = \l_{l+1}$. If $l=2$ then $\l=\l_3$ as in Table \ref{TAB1}, contrary to our initial assumption, so $l\ge 3$. In $V^3$, again the weight $\l-\b_l-2\g_1-\b_{l+2}$ affords the second summand in the table above, and this contains 
$(\om_1+\om_{l-1}+\om_l)^4$ by Lemma \ref{mults}, giving a contradiction to Corollary \ref{V^2(Q_X)} as before. This completes the proof of the lemma. \hal

\vspace{4mm}
Recall that if $\l = \sum c_i\l_i$ then $L(\l)$ is defined to be the number of values of $i$ such that $c_i \ne 0$, with similar definitions for $L(\mu^i)$.

\begin{lemma}\label{L02} We have $L(\mu^0) \le 1$ and $L(\mu^2)\le 1$.
\end{lemma}

\pf Suppose $L(\mu^0) \ge 2$, and write $\mu^0 = \mu'+b\l_j^0+a\l_k^0$, where $a,b\ne 0$, $j<k$ and $j,k$ are maximal subject to having nonzero coefficients. Now 
\[
V^1 = (\mu^0 \otimes \mu^2)\downarrow L_X',
\]
while in $V^2$, the weight $\l-\b_k-\cdots -\b_{l+1}$ affords a summand $\nu \otimes ((\om_1+\om_l)\oplus 0) \otimes \mu^2\downarrow L_X'$, where 
\[
\nu = \mu'+ b\om_j+\om_{k-1}+(a-1)\om_k
\]
(where we have written also $\mu'$ for $\mu'\downarrow L_X'$). 
Thus
\[
V^2 \supseteq \nu \otimes \om_1 \otimes \om_l \otimes (\mu^2\downarrow L_X').
\]
Now  using Corollary \ref{pieri} we see that $\nu \otimes \om_1$ has composition factors $\nu_1,\nu_2,\nu_3$, where
\[
\begin{array}{l}
\nu_1 = \mu'+b\om_j+a\om_k, \\
\nu_2 = \mu'+\om_1+b\om_j+\om_{k-1}+(a-1)\om_k,\\
\nu_3 = \mu'+(b-1)\om_j+\om_{j+1}+\om_{k-1}+(a-1)\om_k.
\end{array}
\]
Then $\nu_1 \otimes \om_l \otimes (\mu^2\downarrow L_X') = V^1\otimes \om_l$, while $\nu_2\otimes \om_l$ and $\nu_3\otimes \om_l$ both contain a composition factor $ \mu'+b\om_j+\om_{k-1}+(a-1)\om_k$. This is a contradiction by 
Corollary \ref{cover}(ii). Thus $L(\mu^0) \le 1$, and by duality, $L(\mu^2) \le 1$ also. \hal

\begin{lemma}\label{muposs1} We have $\mu^0 = a\l_1^0$ and $\mu^2 = b\l_l^2$ for some $a,b\ge 0$.
\end{lemma}

\pf Suppose false. Then using duality and Lemma \ref{L02}, we may assume that $\mu^2 = b\l_i^2$ with $b\ge 1$ and $i<l$. Note that  
by Lemma \ref{mu10},
\[
V^1 = (\mu^0 \otimes \mu^2)\downarrow L_X' = (\mu^0 \downarrow L_X') \otimes b\om_{l-i+1}.
\]
In $V^2$, the weight $\l-\b_{n-l}-\cdots -\b_{n-l+i}$ affords a summand
\[
\begin{array}{ll}
& (\mu^0\downarrow L_X') \otimes ((\om_1+\om_l)\oplus 0) \otimes (\om_{l-i}+(b-1)\om_{l-i+1})\\
= & (\mu^0\downarrow L_X') \otimes \om_1\otimes \om_l \otimes (\om_{l-i}+(b-1)\om_{l-i+1}).
\end{array}
\]
Assume $b\ge 2$. Then using Corollary \ref{pieri} we see that $(\om_{l-i}+(b-1)\om_{l-i+1}) \otimes \om_1$ has summands $\nu_1,\nu_2,\nu_3$, where
\[
\begin{array}{l}
\nu_1 = b\om_{l-i+1}, \\
\nu_2 = \om_1+\om_{l-i}+(b-1)\om_{l-i+1}, \\
\nu_3 = \om_{l-i}+(b-2)\om_{l-i+1}+\om_{l-i+2}
\end{array}
\]
(where there is no $\om_{l-i+2}$ term in $\nu_3$ if $i=1$). However, $(\mu^0\downarrow L_X') \otimes \nu_1 \otimes \om_l = V^1 \otimes \om_l$, while (again using Corollary \ref{pieri}) $\nu_2 \otimes \om_l$ and $\nu_3 \otimes \om_l$ both have a composition factor $\om_{l-i}+(b-1)\om_{l-i+1}$, so this contradicts Corollary \ref{cover}(ii).

Hence $b=1$ and $\mu^2 = \l_i^2$. Now $V^1 = (\mu^0\downarrow L_X') \otimes \om_{l-i+1}$, while 
$V^2 \supseteq (\mu^0\downarrow L_X') \otimes \nu$, where
\[
\nu = \om_1\otimes \om_l \otimes \om_{l-i} = \om_l \otimes ((\om_1+\om_{l-i})\oplus \om_{l-i+1}).
\]
It follows by  Corollary \ref{cover}(ii) that $(\mu^0\downarrow L_X') \otimes \om_l \otimes (\om_1+\om_{l-i})$ is MF. By Proposition \ref{stem1.1.A}, this forces $\mu^0 = 0$ or $\l_l^0$.

If $\mu^0 = \l_l^0$, then $(\mu^0\downarrow L_X') \otimes \om_l \otimes (\om_1+\om_{l-i}) = \om_l \otimes \om_l \otimes (\om_1+\om_{l-i})$ contains $(\om_{l-i}+\om_l)^2$, a contradiction. Hence $\mu^0=0$.

We now have $\mu^0=\mu^1=0$, $\mu^2=\l_i^2$ with $i<l$. Also 
\[
V^2 \supseteq \om_1\otimes \om_l\otimes \om_{l-i} = (V^1\otimes \om_l)+((\om_1+\om_{l-i})\otimes \om_l).
\]
Let $\g_1 = \b_{l+1}$  and $\g_2 = \b_{n-l}$ (the nodes not in $C^0,C^1$ or $C^2$), and let $x=\la \l,\g_1\ra$, $y=\la \l,\g_2\ra$. 

If $x\ne 0$, then the weight $\l-\g_1$ affords a further summand $\om_l \otimes ((\om_1+\om_l)\oplus 0) \otimes \om_{l-i+1}$ in $V^2$. This must be MF, so Proposition \ref{stem1.1.A} implies that $i=1$; however $\om_l \otimes (\om_1+\om_l) \otimes \om_l$ is not MF as it contains $(\om_1+\om_{l-1}+\om_l)^2$. Therefore $x=0$.

If $y\ne 0$, then the weight $\l-\g_2$ affords a further summand $((\om_1+\om_l)\oplus 0) \otimes (\om_{l-i+1}+\om_l)$ in $V^2$; however this is again not MF. Hence $y=0$.

At this point we have $\l = \l_{n-l+i}$ with $i<l$. Replace $V$ by its dual to take $\l = \l_j$ with $1<j\le l$. By our initial assumption that $\l$ is not in Table \ref{TAB1}, we have $j\ge 4$, and so also $l\ge 4$. Now
\[
\begin{array}{ll}
V^1  & = \om_j,  \\
V^2 & = \om_{j-1}\otimes \om_1\otimes \om_l \\
        & = (\om_1+\om_{j-1}+\om_l) \oplus  (\om_1+\om_{j-2})\oplus  (\om_j+\om_l)\oplus  \om_{j-1}^2, 
\end{array}
\]
and in $V^3$ we have the summand $\om_{j-2}\otimes \wedge^2((\om_1+\om_l)\oplus 0)$, afforded by 
$\l-\b_{j-1}-2\b_j-2\b_{j+1}-\cdots -2\g_1-\b_{l+2}$. 
 Using Lemma \ref{mults} together with Section \ref{LRres}, we see that between them, the two summands in $V^3$ contain $(\om_1+\om_{j-3})^3$. However by Corollary \ref{V^2(Q_X)} only one of these is in the first level of a composition factor of $V^2$, so this contradicts Proposition \ref{induct}. This completes the proof of the lemma. \hal

\vspace{2mm}
By the previous lemmas, we now have 
\[
\mu^1=0, \;\mu^0=a\l_1^0, \; \mu^2 = b\l_l^2
\]
 for some $a,b$, not both zero. Replacing $V$ by its dual if necessary, we can assume that $a\ge b$. In particular, $a\ne 0$, and 
\[
V^1 = a\om_1 \otimes b\om_1.
\]
As before, define $\g_1 = \b_{l+1}$  and $\g_2 = \b_{n-l}$, and let $x=\la \l,\g_1\ra$, $y=\la \l,\g_2\ra$. 

\begin{lemma}\label{muposs2} We have $x=y=0$. 
\end{lemma}

\pf In $V^2$, the weight $\l-\b_1-\b_2-\cdots -\g_1$ affords a summand 
\begin{equation}\label{summ100}
(a-1)\om_1 \otimes \om_1\otimes \om_l \otimes b\om_1,
\end{equation}
and this contains $V^1 \otimes \om_l$. If $x \ne 0$, then in addition the weight $\l-\g_1$ affords a summand
\[
(a\om_1+\om_l) \otimes \om_1\otimes \om_l \otimes b\om_1
\]
of $V^2$, and this is not MF, contradicting Corollary \ref{cover}(ii). Similarly, if $y\ne 0$ the weight $\l-\g_2$ afford a summand
\[
a\om_1 \otimes \om_1\otimes \om_l\otimes (b\om_1+\om_l)
\]
of $V^2$, and this is also not MF. Hence $x=y=0$. \hal

\begin{lemma}\label{muposs2} We have $b=0$ and $a\ge 3$.
\end{lemma}

\pf If $a=1$ then also $b=1$ (as $\l \ne \l_1$ by assumption), and so $\l = \l_1+\l_n$. Then 
\[
V\downarrow X = ((\om_1+\om_l) \otimes (\om_1+\om_l))-0.
\]
However  this is not MF, a contradiction. Hence $a\ge 2$.

Now suppose that $b\ne 0$. As in (\ref{summ100}) above, $V^2$ has a summand 
\[
\begin{array}{ll}
& (a-1)\om_1 \otimes \om_1\otimes \om_l \otimes b\om_1 \\
= & (a\om_1 \oplus  ((a-2)\om_1+\om_2))\otimes \om_l \otimes b\om_1.
\end{array}
\]
However, $((a-2)\om_1+\om_2))\otimes \om_l \otimes b\om_1$ contains $((a+b-3)\om_1+\om_2)^2$ by Lemma \ref{mults}(iii), which contradicts Corollary \ref{cover}(ii). 

Hence $b=0$. As $\l\ne 2\l_1$ by assumption, we also have $a\ge 3$. \hal

\vspace{4mm}
\no {\bf Completion of the proof}

At this point we have $\l = a\l_1$ with $a\ge 3$. In this case $V^2$ is equal to the summand in (\ref{summ100}), so 
using Corollary \ref{pieri}, 
\[
\begin{array}{ll}
V^2 & = (a-1)\om_1 \otimes \om_1\otimes \om_l  \\
       & = (a\om_1 \oplus  ((a-2)\om_1+\om_2))\otimes \om_l \\
       & = (a\om_1+\om_l) \oplus  ((a-2)\om_1+\om_2+\om_l)\oplus  ((a-3)\om_1+\om_2)\oplus  ((a-1)\om_1)^2.
\end{array}
\]
Now in $V^3$ we have the following summands:
\[
\begin{array}{l|l}
\hline
\hbox{summand of }V^3  & \hbox{afforded by } \\
\hline
(a-2)\om_1 \otimes S^2((\om_1+\om_l)\oplus 0) & \l-2\b_1-2\b_2-\cdots -2\g_1 \\
(a-1)\om_1 \otimes \om_l & \l-\b_1-\b_2-\cdots -\g_2 \\
\hline
\end{array}
\]
Now $S^2((\om_1+\om_l)\oplus 0) \supseteq (\om_1+\om_l)^2 \oplus  0^2$ by Lemma \ref{mults}(ii), so using Corollary \ref{pieri} we see that the first summand in the above table contains $((a-2)\om_1)^4$; the second also contains $(a-2)\om_1$. So $V^3$ contains 
$(a-2)\om_1$ with multiplicity at least 5. However, by Corollary \ref{V^2(Q_X)} there are only three composition factors of $V^2$ that have $(a-2)\om_1$ in their first level, namely $((a-3)\om_1+\om_2)$ and $((a-1)\om_1)^2$.  
This is a contradiction by Proposition \ref{induct}.

This completes the proof of Theorem \ref{om1oml+1}.

\chapter{Proof of Theorem \ref{MAINTHM}, Part I: $V_{C^i}(\mu^i)$  is usually trivial}\label{pfstart}

We adopt the hypotheses of Theorem \ref{MAINTHM}. Note that for $X=A_1$, the theorem was proved in \cite{LST}; and the case where $X=A_2$ was covered in Chapter \ref{casea2}. So let $X = A_{l+1}$ with $l\ge 2$, let $W = V_X(\d)$ and $Y = SL(W) = A_n$. Suppose $V = V_Y(\l)$ is a nontrivial  irreducible $Y$-module such that $V\downarrow X$ is multiplicity-free and $\l$ is not $\l_1$ or its dual.

Let $L_X' < L_Y' = C^0\times \cdots\times C^k$ as in Chapter \ref{levelset} and let $\mu^i$ be the restriction of $\l$ to $T_Y\cap C^i$, so that 
$V^1(Q_Y)\downarrow L_Y' = V_{C^0}(\mu^0) \otimes \cdots \otimes V_{C^k}(\mu^k)$. 

Let $\d = \sum_{i=1}^{l+1}d_i\om_i$. By Theorem \ref{LEVELS}, $L_X'$ is irreducible on the levels $W^1(Q_X)$ and $W^{k+1}(Q_X)$, with highest weights $\d',\d''$ respectively, where 
\[
\d' = \sum_{i=1}^l d_i\om_i, \;
\d'' = \sum_{i=1}^l d_{i+1}\om_i.
\]
For the other levels $W^{i+1}(Q_X)$ ($0<i<k$), let  $\delta_1^i, \ldots , \delta_{k_i}^i$ denote the highest weights of
the irreducible $L_X'$-summands; that is, 
\[
W^{i+1}(Q_X)\downarrow L_X' = \sum_{j=1}^{k_i}V_{L_X'}(\d_j^i).
\]
For each $i$ the projection of $L_X'$ to $C^i$ corresponds to the
action of $L_X'$ on the $i$th level $W^{i+1}(Q_X)$. We know that $V^1$ is MF; in particular
$V_{C^i}(\mu^i) \downarrow L_X'$ is MF for all $i$. 

Throughout the proof we adopt the inductive hypothesis that Theorem \ref{MAINTHM} holds for groups of type $A_m$ with $m<l+1$. The induction starts with $X = A_2$, as that case of Theorem \ref{MAINTHM} has been established in Chapter \ref{casea2}. We can also assume that $\d \ne r\o_j$ for any $r,j$ by the main results of Chapters \ref{rkge2}, \ref{ro1}, \ref{oige3}, \ref{o2sec}. As a consequence, Theorem \ref{LEVELS} implies that $k_i\ge 2$ for all $i \ne 0,k$.

For each $i$  let $C^i$ have  fundamental
roots  $\Pi(C^i) = \{\beta_1^i, \ldots , \beta_{r_i}^i\}$  and corresponding fundamental dominant weights
 $\{ \lambda_1^i, \ldots , \lambda_{r_i}^i \}.$  Write 
\[
\mu^i = \sum_{j=1}^{r_i} c_j^i\l_j^i.
\]
Let $\g_i$ denote the fundamental root between $C^{i-1}$ and $C^i$ for $i = 1, \ldots , k.$ 

There are two main results in this chapter.
Recall our notation that for a dominant weight $\nu$, the number of nonzero coefficients in the expression for $\nu$ as a sum of fundamental dominant weights is denoted by $L(\nu)$.

\begin{theor}\label{NORDU}  Assume that $\delta$ is not of the form $r\omega_s$ (that is, $L(\d)\ge 2$) and that the  induction hypothesis holds.  
\begin{itemize}
\item[{\rm (i)}]  For $1 \le  i \le k-1$,  $V_{C^i}(\mu^i)$  is a trivial, natural, or dual of natural module.
\item[{\rm (ii)}]  If $L(\delta_j^i)\le 1$  for each $1\le j \le k_i$, then  $i = 1$ or $k-1,$ the embedding of $L_X'$ in $C^i$ is given by $V_{L_X'}(2\omega_1) + V_{L_X'}(\omega_2)$ or its dual, and $\delta = d_1\omega_1 + \omega_2$
or $\omega_l + d_{l+1}\omega_{l+1}.$
\end{itemize}
\end{theor}

The theorem has strong consequences.  Namely in part (i), except perhaps for the first and last levels, the $\mu^i$ each afford a trivial, natural, or dual of natural module for $C^i$.  Therefore,
 $V_{C^i}(\mu^i)\downarrow L_X' $ is either trivial, $ \sum_j V_{L_X'}(\delta_j^i)$ or the dual of this module.
 Now part (ii) shows that except in very special situations, this restriction always has an irreducible
for which $L(\d_j^i)\ge 2$.  So in view of Proposition \ref{tensorprodMF}, ignoring the
 special cases in (ii) there can be at most one nontrivial $\mu^i$ for $i \ne 0,k.$

Theorem \ref{NORDU} and some additional arguments
yield the following key result which shows that except for some very special configurations
the only possible nonzero $\mu^i$ occur for $i = 0, k$.

 \begin{theor}\label{MUIZERO}  Assume that  the  induction hypothesis holds and $\delta \ne r\omega_j$ or    $\omega_1 + \omega_{l+1}.$ If $0 < i < k$,  then  $\mu^i = 0.$ 
 \end{theor}
 
 Before beginning the proof of these theorems we state and prove a useful
 corollary of Theorem \ref{MUIZERO}.
 
 As usual we adopt the following notation, for each $i$:
\[
V^i = V^i(Q_Y) \downarrow L_X'.
\]
  By Corollary \ref{cover} we know that $A := V^1 \otimes V_{L_X'}(\omega_l)$ covers
  those irreducible summands in $V^2$ that arise from $V^1$.
  The following result shows that in certain  situations both $A$ and a specified additional summand occur in  $V^2$.  In order for $V\downarrow X$ to be MF this additional summand must also be  MF.  The summand appears in the restriction of $V_{\gamma_1}^2(Q_Y)$ (notation as in Chapter \ref{notation}).  Further summands will appear if $\la\lambda, \gamma_1\ra \ne 0$ and also possibly in $V_{\gamma_i}^2(Q_Y)$ for $i > 1.$

  \begin{coro}\label{coversums}  Assume $V^1(Q_Y) = V_{L_Y'}(\mu^0 \otimes \mu^k)$, and that $\mu^0 \ne 0$ and $d_{l+1} \ne 0$.  Let $i$ be minimal with $c_i^0 \ne 0$ and let $\nu^0$ denote the restriction of $ \mu^0 - \beta_i^0 - \cdots - \b_{r_0}^0-\gamma_1$   to $C^0.$  Then
  $V^2 \supseteq A \oplus B$, where $A := V^1 \otimes V_{L_X'}(\omega_l)$ and $B$ is as follows:
  \begin{itemize}
\item[{\rm (i)}] if $i > 1$, then $B = V_{L_Y'} ((\lambda_1^0 + \nu^0) \otimes \mu^k) \downarrow L_X' \otimes V_{L_X'}(\omega_l)$;
  \item[{\rm (ii)}] if $i = 1,$  $c_1^0 \ge 2,$ then $B =  V_{L_Y'}((\mu^0 -\beta_1^0)) \otimes \mu^k) \downarrow L_X'  \otimes V_{L_X'}(\omega_l)$;
  \item[{\rm (iii)}] if $i = 1= c_1^0$, then either  $\mu^0 =  \lambda_1^0$ and $B = 0$; or 
  $\mu^0 \ne  \lambda_1^0$ and  $B = V_{L_Y'}((\mu^0 -\beta_1^0- \cdots - \beta_j^0) \otimes \mu^k) \downarrow L_X'  \otimes V_{L_X'}(\omega_l)$, where $j >1$  is minimal with  $c_j^0 > 0.$
 \end{itemize}
  \end{coro}

 \pf  By hypothesis $V^1(Q_Y) = V_{L_Y'}(\mu^0 \otimes \mu^k)$. The assumption $d_{l+1} \ne 0$ together with 
Corollary \ref{V^2(Q_X)} and Theorem \ref{LR} imply that 
\[
V_{C^1}(\lambda_1^1) \downarrow L_X' = V_{L_X'}(\delta') \otimes V_{L_X'}(\omega_l) = 
V_{C^0}(\lambda_1^0) \downarrow L_X' \otimes V_{L_X'}(\omega_l).
  \]
 Therefore Theorem \ref{MUIZERO} implies that  $V^2(Q_Y) \supseteq V_{L_Y'}(\nu^0 \otimes \lambda_1^1 \otimes \mu^k)$, so we have  
\begin{equation}\label{2ndp}
V^2  \supseteq V_{C^0}(\nu^0 \otimes \lambda_1^0) \downarrow L_X' \otimes V_{L_X'}(\omega_l) \otimes 
V_{C^k}(\mu^k) \downarrow L_X'.
\end{equation}
    Assume $i > 1$.  Then   $V_{C^0}(\nu^0 \otimes \lambda_1^0)  \supseteq 
V_{C^0}(\lambda_1^0 + \nu^0)  \oplus V_{C^0}(\mu^0),$  since $\mu^0 = 
(\lambda_1^0 + \nu^0) - \beta_1^0 - \cdots - \beta_{i-1}^0.$
  Therefore 
\[
V^2  \supseteq \left(V_{L_Y'}(\mu^0 \otimes \mu^k) \downarrow L_X'
 \otimes   V_{L_X'}(\omega_l)\right)  \oplus   \left((V_{L_Y'}((\lambda_1^0 + \nu^0) \otimes \mu^k)\downarrow L_X')
 \otimes V_{L_X'}(\omega_l)\right).
\]
   The first summand is $ A$ and the second is $B$ so this yields (i).

     Now assume $i = 1.$  Note that $\lambda_1^0 + \nu^0 = \mu^0$.  If $c_1^0 > 1$, then 
$V_{C^0}(\lambda_1^0 \otimes \nu^0) \supseteq V_{C^0}(\mu^0) \oplus V_{C^0}(\mu^0 - \beta_1^0)$.  And if 
$c_1^0 = 1$, then one of the following holds:
\begin{itemize}
\item[(1)]  $\mu^0 = \lambda_1^0$, $\nu^0 = 0$ and  $V_{C^0}(\lambda_1^0 \otimes \nu^0) = V_{C^0}(\mu^0)$, or   \item[(2)]  $\mu^0 \ne \lambda_1^0$, $j$ exists as in (iii), and 
$V_{C^0}(\lambda_1^0 \otimes \nu^0) \supseteq V_{C^0}(\mu^0) \oplus V_{C^0}(\mu^0 - \beta_1^0 - \cdots - \beta_j^0)$.
\end{itemize}    
At this point we argue as before using (\ref{2ndp}) to get conclusion (ii) or (iii) of the statement.  \hal

\section{Proof of Theorem \ref{NORDU}}

The proof of Theorem \ref{NORDU} will follow from a series of lemmas which we now begin.
Adopt the hypotheses of the theorem, and the above notation.

\begin{lem}\label{murestriction}  Fix a level $i$ and a weight $\delta_j^i$.  Let $D = V_{L_X'}(\delta_j^i).$
\begin{itemize}
\item[{\rm (i)}]  If the only nontrivial irreducible representations of $SL(D)$  that are MF upon restriction to the image of $L_X'$ are $D$ and $D^*$, then $\mu^i \in \{ 0, \lambda_1^i, \lambda_{r_i}^i \}.$ 
\item[{\rm (ii)}] If the only nontrivial irreducible representations of $SL(D)$  that are MF upon restriction to the image of $L_X'$ are $D$, $\wedge^2(D)$, $S^2(D)$ and duals of these,   then $\mu^i \in \{ 0, \lambda_1^i,  \lambda_2^i, 2\lambda_1^i, \lambda_{r_i}^i, \lambda_{r_i-1}^i,2\lambda_{r_i}^i \}.$ 
\end{itemize}
\end{lem}

\pf  To simplify notation set $C = C^i$ and $\mu = \mu^i$. If $k_i = 1$, then $C = SL(D)$ and (i) and (ii) are  immediate from the hypotheses.  So assume $k _i > 1.$    Then there is a proper Levi factor of $C$ of type $A_x \times A_y$, where $L_X'$ embeds into $A_x$ via $V_{L_X'}(\delta_j^i)$ and  into $A_y$ by the sum of the irreducible modules $V_{L_X'}(\delta_m^i)$ for $m \ne j.$  The hypothesis in each of (i) and (ii) imply that
$D$ is not a natural or dual module for $L_X'$, as otherwise $L_X'$ and $SL(D)$ both have type $A_l$ and any irreducible representation of $SL(D)$ is irreducible (hence MF) for $L_X'$.  Consequently, $x \ge 5$ (as $l\ge 2$)

Dropping  superscripts and subscripts $i,$ let  $\Pi(C) = \{\beta_1, \ldots, \beta_r\}$ with
corresponding fundamental weights $\{\l_1, \ldots, \l_r\}.$ Let $\beta_k$ be  the node outside $\Pi(A_x) \cup \Pi(A_y)$ and reorder, if necessary, so that $\Pi(A_x) = \{\beta_1, \ldots, \beta_x\}$. Write $\mu = \sum a_j\l_j.$

Consider the levels of $V_C(\mu)$ with respect to the above Levi subgroup. We know that 
 $V_C(\mu)\downarrow L_X'$ is MF.  Each composition factor of $V_C(\mu)\downarrow (A_x \times A_y)$ is a tensor product of irreducibles for the factors, and such an irreducible for $A_x$, when restricted to the projection of $L_X'$, must be MF.  Hence the irreducible is either trivial or it must be one of the modules indicated in (i) or (ii) (i.e. $D$ or $D^*$ in (i), or $D,\wedge^2(D), S^2(D)$ or a dual in (ii)).
Therefore, all composition factors of $A_x$ on $V_C(\mu)$ are among these.  Consequently, the same must hold
for any Levi factor of $C$ conjugate to $A_x$.

(i)  If $\mu$ is not $\l_1, \l_r,$ or $\l_1 + \l_r,$ then we claim that we can
find a Levi factor conjugate to $A_x,$ for which there is a composition factor which is not a trivial, natural, or
dual module.  To see this when $a_j \ne 0$ for some $j \ne 1, r$,  just take a rank $x$ Levi with fundamental system
a subset of $\Pi(C)$ and for which $\beta_j$ is not an end node.  The claim is also clear for $c\l_1 + d\l_r$ with $c>1$ or $d>1$; hence the claim is proved. It remains to rule out the case where $\mu = \l_1 + \l_r.$  But here $A_x$ has a composition factor with highest weight $\l_1 + \l_x$, contradicting the hypothesis of (i).
 
Now consider (ii) and proceed as  in the last paragraph. We easily get a contradiction unless $\mu = a_1\l_1 + a_2\l_2 + a_{r-1}\l_{r-1} + a_r\l_r.$  Moreover, restricting to $A_x$ we see that
$(a_1, a_2) \in \{ (0,0), (1,0), (0,1), (2,0) \}$ and restricting to the conjugate of $A_x$ with fundamental system
$\beta_{r-x+1}, \ldots , \beta_r,$ we find that  $(a_{r-1}, a_r) \in \{ (0,0), (0,1), (1,0), (0,2) \}.$
If either $(a_1, a_2) = (0,0)$ or $(a_{r-1}, a_r) = (0,0)$, then we obtain (ii).  Therefore assume
neither pair is $(0,0).$  Let $j \in \{r-1, r \}$ with $a_j \ne 0$.  Then $\mu - \beta_{x+1} - \cdots -\beta_j$
affords a nontrivial module at level $1$ in $V_C(\mu)$ for $A_x \times A_y$ and the highest weight restricts to $A_x$
as $a_1\l_1 + a_2\l_2  + \l_x$.  This contradicts the hypothesis of (ii)
and completes the proof of the lemma. \hal

The inductive hypothesis  shows that the hypothesis of Lemma \ref{murestriction}(i) holds if
there exists $j$ such that $\delta_j^i$ has at least three nonzero labels.  The next three lemmas focus
on situations where the irreducible module corresponding to the largest monomial at a given level has just one or two nonzero labels.  These results do not require the inductive hypothesis.
 
Recall that $W = V_X(\delta)$ with $\delta = d_1\omega_1 + \cdots + d_{l+1}\omega_{l+1}.$  Fix a level of $W$, say level $x$.  Theorem \ref{LEVELS} describes a filtration of the level under the action of $L_X'$.  
Recall the ordering of the monomials $f_1^{a_1} \cdots  f_{l+1}^{a_{l+1}}$ introduced in Section \ref{pflevel}. 
Assume the highest weight of the factor with the largest monomial at level $x$ is  $\nu$, where $\nu$ is afforded by $f_j^{c_j} \cdots f_k^{c_k}v,$  where $c_j \ne 0 \ne c_k$ and $1 \le j \le k \le l+1.$  We will often abuse terminology and simply say that $\nu$ is afforded by $f_j^{c_j} \cdots f_k^{c_k}.$
Then $x = c_1 + \cdots + c_k$ and $c_i \le d_i$ for each $i.$  To simplify notation we will identify  $\nu$  and
other weights with their restriction to  $L_X'$.   Also $f_j^{c_j} \cdots f_k^{c_k}$ is a weight vector for $T_X$ in the algebra $N$ defined in Section \ref{mainreslev}, and we will often identify it with the corresponding weight.

We first make a few observations from the fact that $\nu$ has the largest monomial.  First, $d_i = 0$ for all $i > k.$  Also,  if $j < k$ then $c_i=d_i$ for all $i>j$, since otherwise the monomial 
  $ f_j^{c_j-1} \cdots f_{i-1}^{c_{i-1}}f_i^{c_i+1}f_{i+1}^{c_{i+1}} \cdots$ is strictly larger than $f_j^{c_j} \cdots f_k^{c_k}$, a contradiction. Hence the label of $\nu$ at node $i-1$ is $d_{i-1} - c_{i-1} + d_i.$ Since $\d \ne r\o_j$ by hypothesis, it follows that $\nu \ne 0$.

\begin{lem}\label{1label}  Assume that $L(\d)\ge 2$, and that $\nu$ is the highest weight afforded by the largest monomial  $f_j^{c_j}\cdots f_k^{c_k}$   at some level of $W$ that is not the top or bottom level.  If $L(\nu)=1$, then 
the $S$-value $S(\nu)\ge 2$ and one of the following holds:
\begin{itemize}
\item[{\rm (i)}] $f_j^{c_j}\cdots f_k^{c_k} = f_k^{c_k}$, $\d = d_{k-1}\omega_{k-1} + d_k\omega_k$   and $\nu =  (d_{k-1}+c_k)\omega_{k-1}.$ Also, $c_k = d_k$ unless $k = l+1$;
\item[{\rm (ii)}] $f_j^{c_j}\cdots f_k^{c_k} = f_1^{c_1}f_2^{d_2}$,  $c_1 < d_1$, $\d = d_1\omega_1 + d_2\omega_2$   and $\nu =  (d_1-c_1+d_2)\omega_1.$
\end{itemize}
\end{lem}

\pf  If $k = 1$, then by the above $\d = d_1\omega_1$ and we are assuming this is not the case.  Therefore $k > 1$   and
the hypothesis implies $\nu = (d_{k-1}-c_{k-1} + c_k)\omega_{k-1}.$  
 
 There are two cases.  First assume $j = k$, so that $\nu $ is afforded by $f_k^{c_k}.$  By hypothesis there exists $i < k$ with $d_i \ne 0.$  If  $i < k-1$, then $\nu$ also has a nonzero coefficient of $\omega_i$, a contradiction.
 Therefore, $i = k-1,$  $\nu = (d_{k-1} + c_k)\omega_{k-1}.$    
 If $k < l+1$, then we must have $c_k = d_k$ to avoid a nonzero multiple of $\omega_k.$

Now suppose $j < k.$  This forces $j = 1,$ as otherwise $\nu$ would have a nonzero multiple of $\omega_{j-1}$.  So here
the earlier remarks imply $\nu$ is afforded by  $f_1^{c_1}f_2^{d_2} \cdots f_k^{d_k}.$  Also $c_1 < d_1$ as otherwise the monomial would be the largest
among all monomials and this only occurs at the last level. Then the coefficient of $\omega_1$ in $\nu$ is $d_1-c_1 + d_2>0$  and this forces $1 = k-1,$ whence $k = 2.$  Therefore $\d = d_1\omega_1 + 
d_2\omega_2$ and $\nu = (d_1-c_1+d_2)\omega_1.$    \hal
  
\begin{lem}\label{2labels}  Assume $L(\d)\ge 2$ and that $\nu$ is the highest weight afforded by the largest monomial  $f_j^{c_j}\cdots f_k^{c_k}$ at some level of $W$ that is not the top or bottom level.  Assume  that $L(\nu)=2$.  
Then $k > 1$ and one of the following holds:
\begin{itemize}
\item[{\rm (i)}]  $f_j^{c_j}\cdots f_k^{c_k} = f_k^{c_k}, $ $k < l+1, c_k < d_k,$ $\d = d_{k-1}\omega_{k-1} + d_k\omega_k$, and 
$\nu = (d_{k-1}+c_k) \omega_{k-1} + (d_k-c_k)\omega_k$;
\item[{\rm (ii)}]  $f_j^{c_j}\cdots f_k^{c_k} = f_k^{c_k}$,   $\d = d_i\omega_i + d_{k-1}\omega_{k-1} + d_k\omega_k,$ for $i < k-1$ and   $\nu = d_i\omega_i + (d_{k-1} + c_k)\omega_{k-1}.$  Also, $d_i \ne 0$ and either $c_k = d_k$ or $k = l+1$;
\item[{\rm (iii)}] $f_j^{c_j}\cdots f_k^{c_k} = f_1^{c_1}f_2^{d_2}f_k^{d_k},$  $c_1 < d_1$, $d_2 \ge 0$, $\d = d_1\omega_1 +d_2\omega_2 + d_k\omega_k$ and $\nu = (d_1-c_1+d_2)\omega_1 + d_k\omega_{k-1}$; 
\item[{\rm (iv)}] $f_j^{c_j}\cdots f_k^{c_k} = f_j^{d_j}f_k^{d_k}$,  $1<j<k$, $d_{j-1} \ne 0 \ne d_j,$ $\d = d_{j-1}\omega_{j-1} + d_j\omega_j + d_k\omega_k$ and $\nu = (d_{j-1} + d_j)\omega_{j-1} + d_k\omega_{k-1}$;
\item[{\rm (v)}] $f_j^{c_j}\cdots f_k^{c_k} = f_{k-1}^{c_{k-1}}f_k^{d_k}$, $2 < k$, $c_{k-1} \le d_{k-1}, $ $\d = d_{k-2}\omega_{k-2} + d_{k-1}\omega_{k-1} + d_k\omega_k$ and $\nu = (d_{k-2} + c_{k-1})\omega_{k-2} + (d_{k-1}-c_{k-1} + d_k)\omega_{k-1}.$
\end{itemize}
\end{lem} 

\pf By assumption, $L(\nu)=2$. As above $k > 1,$ $c_i = d_i$ for $i>j$,  and $d_i = 0$ for $i > k.$

The coefficient of $\omega_{k-1}$ in $\nu$ is $d_{k-1} - c_{k-1} + c_k >0,$  so there is
precisely one other nonzero coefficient of $\nu.$  This could occur at $\omega_k$ if $k< l+1,$  $c_k < d_k,$
and $\nu = (d_{k-1} -c_{k-1} + c_k)\omega_{k-1} + (d_k - c_k)\omega_k.$  The comments preceding Lemma \ref{1label}
imply $j =k$ in this situation.  Also, $d_i = 0$ for $i < k-1$, as otherwise $\nu$ would also have a nonzero multiple of $\omega_i.$  So here (i) holds.

Suppose $j = k.$ Excluding the above  case,  we have 
$\nu = d_i\omega_i + (d_{k-1} + c_k)\omega_{k-1}$ for some $i < k-1,$  $\d =  d_i\omega_i + d_{k-1}\omega_{k-1} + d_k\omega_k,$   and (ii) holds. Also either $c_k = d_k$ or $k = l+1.$

Now assume that $j < k.$  Suppose in addition that $c_j < d_j.$ Then $\nu$ has  nonzero coefficients of $\omega_j,$ $\omega_{k-1}$, and also $\omega_{j-1}$ provided $j > 1.$ All other coefficients must be $0$.  If $j > 1$, then
we must have $j = k-1$.  In this case $\d = d_{k-2}\omega_{k-2} + d_{k-1}\omega_{k-1} + d_k\omega_k$ and 
$\nu = (d_{k-2} + c_{k-1})\omega_{k-2} + (d_{k-1}-c_{k-1} + d_k)\omega_{k-1},$ giving (v).
Now suppose $j = 1.$   If $d_s \ne 0$  for some $j < s < k$, then the coefficient of $\omega_{s-1}$ is nonzero, a contradiction unless $s = 2$.   So in this case $\d = d_1\omega_1 +d_2\omega_2 + d_k\omega_k$ and $\nu = 
(d_1-c_1+d_2)\omega_1 + d_k\omega_{k-1}$ is afforded by $f_1^{c_1}f_2^{d_2}f_k^{d_k}$, giving (iii),   which allows for the possibility that $d_2 = 0.$

Now suppose $c_j = d_j$.  We are assuming that the level is not the first or last level and we have $c_i = d_i$ for $i \ge j.$
Hence there must exist $i < j$ with $d_i \ne 0$.  Therefore, $j >1$ and for any such $i$ the label of $\omega_i$ in $\nu$
is nonzero.  But the coefficient of $\omega_{j-1}$  is also nonzero. Therefore, $i = j-1,$ and hence $i$ is unique, and $\d = d_{j-1}\omega_{j-1} + d_j\omega_j + d_k\omega_k.$ This gives (iv).  \hal

\begin{lem}\label{notw2ors2}  Assume that $L(\d)\ge 2$, and that we are in one of the cases of Lemma $\ref{1label}(i),(ii)$ or Lemma $\ref{2labels}(i)-(v)$ at
level $i$, where this is not the top or bottom level.  Write $D$ for the natural module for $C^i$. Then $V_{C^i}(\mu^i) $ is not $\wedge^2(D)$, $S^2(D)$, or the dual of one of these modules. 
\end{lem}

\pf  We will work with the weight $\nu$ given in the lemmas together with the next highest weight of a composition factor of the level.

Consider Lemma \ref{1label}(i). Here $\nu  =  (d_{k-1}+c_k)\omega_{k-1}$ is afforded by $f_k^{c_k},$  $\d = d_{k-1}\omega_{k-1} + d_k\omega_k,$ and $d_{k-1} \ne 0.$   At level $i$ there is also a composition factor of highest weight $\nu' $ afforded by  $f_{k-1}f_k^{c_k-1}.$
Notice that $\nu' = \nu - \beta_{k-1}^i.$
The restriction to $L_X'$ of $\wedge^2(D)$ contains 
$$\wedge^2(V_{L_X'}(\nu)) \oplus (V_{L_X'}(\nu) \otimes V_{L_X'}( \nu')).$$  The first summand contains the
irreducible of highest weight $2\nu - \beta_{k-1}^i$, and this  is also 
precisely the highest weight of the tensor product. 
Now consider $S^2(D)\downarrow L_X'$, which contains
$$S^2(V_{L_X'}(\nu)) \oplus S^2(V_{L_X'}( \nu')).$$
In this sum the first summand has a submodule of highest weight $2\nu - 2\beta_{k-1}^i$ (as $d_{k-1}+c_k \ge 2$) and this is
precisely the highest weight of the second summand.  So in each case the full module fails to be MF.  So neither 
$\wedge^2(D)\downarrow L_X'$ nor $S^2(D)\downarrow L_X'$  is MF, and this also holds for the duals of these modules.

The other cases are similar.  We will consider the case in Lemma  \ref{2labels}(iii) and leave the remaining cases to the reader. Here 
$\nu  = (d_1-c_1+d_2)\omega_1 + d_k\omega_{k-1}$ is afforded by $ f_1^{c_1}f_2^{d_2}f_k^{d_k}$ and $d_1 > c_1 \ne 0 \ne d_k.$ The argument differs somewhat according to whether or not $d_2 = 0.$ If $d_2 \ne 0,$ then there is a
 second irreducible at level $i$ with highest weight $\nu'  = \nu - \beta_1^i$ afforded by 
$f_1^{c_1+1}f_2^{d_2-1}f_k^{d_k}.$ So the argument proceeds just as above.  
 
 Now suppose $d_2 = 0.$   Here $\nu$ is afforded by  $f_1^{c_1}f_k^{d_k}$ and the next highest weight is  $\nu'  = \nu - \psi^i$, afforded by $f_1^{c_1+1}f_k^{d_k-1}$, where $\psi^i = \beta_1^i + \cdots + \beta_{k-1}^i.$  
 The restriction to $L_X'$ of $\wedge^2(D)$ contains 
$$\wedge^2(V_{L_X'}(\nu)) \oplus (V_{L_X'}(\nu) \otimes V_{L_X'}( \nu - \psi^i)).$$
As $\nu = (d_1-c_1)\omega_1 + d_k\omega_{k-1},$ it follows from  Lemma \ref{lemma1.5(i)(ii)}(ii) that $\wedge^2(V_{L_X'}(\nu))$ contains an irreducible summand of highest weight $2\nu - \psi^i$, and as before the above sum is not MF.

Now consider  $S^2(D)\downarrow L_X'$, which contains
$$S^2(V_{L_X'}(\nu)) \oplus (V_{L_X'}( \nu) \otimes V_{L_X'}(\nu - \psi^i)).$$
It follows from  Lemma \ref{lemma1.5(i)(ii)}(ii) that the first summand has an irreducible summand of highest weight $2\nu - \psi^i$,   and this is the highest weight of the second summand. \hal

\begin{lem}\label{trivnatordual}  Assume that $L(\d)\ge 2$, fix a level $i$ for $W=V_X(\delta)$ other than the top or bottom level, and assume  the inductive hypothesis. 
 Then either $V_{C^i}(\mu^i)$ is a trivial, natural, or dual of natural module, 
or the level restricts to $L_X'$  as one of the following or its dual:
\[
\begin{array}{l}
2\omega_1 \oplus \omega_2,\\
30\oplus  11\; (L_X' = A_2),\\
020 \oplus 101\; (L_X' = A_3).
\end{array}
\]
\end{lem}

\pf  Consider a highest weight $\psi$ for $L_X'$ at  level $i$, and let $N = V_{L_X'}(\psi)$.   
If $L(\psi)\ge 3$, then the induction hypothesis implies that the only nontrivial irreducibles for $SL(N)$ that are MF upon restriction
to $L_X'$ are $N$ and $N^*$.  Then Lemma  \ref{murestriction}(i) gives the assertion.
So now assume that all such $\psi$ have at most two nonzero coefficients.  Then 
Lemmas \ref{1label} and \ref{2labels} apply.

If $L(\psi)=2$,  our induction hypothesis  says that, with one family of exceptions,  the only nontrivial  irreducibles for $SL(N)$ that are MF upon restriction to $L_X'$ are $N$, $\wedge^2N$, $S^2N$, and duals of these.   The exceptions  occur for $\psi = 10 \ldots  01$, where $\wedge^3N$ is MF, and also $S^3N$ 
when  $L_X' = A_2$. Therefore, excluding these exceptions, we can apply  Lemma \ref{murestriction} to restrict the possibilities for $\mu^i.$  Then Lemma \ref{notw2ors2} shows  that $V_{C^i}(\mu^i)$ is a trivial, natural, or dual of natural module.

We must now consider situations where all irreducibles appearing at level $i$ have highest weight $\psi$ with at most one nonzero label or else  $\psi = 10 \ldots  01$.  

We begin with the situation of Lemma \ref{2labels}.  Then $\nu = 10\ldots  01$.
This cannot happen for cases (i) or (iv) of Lemma \ref{2labels}. For case (v)  this would force $d_{k-2} = 0$ and $c_{k-1} = d_{k-1}$.  But then, $i$ is the bottom level, against our hypothesis. Now consider cases (ii) and (iii) with $\delta = 10\ldots  0d_{l+1}$ or $ d_10\ldots  01,$ respectively.  In both cases  the other irreducible at this level is the trivial module $0\ldots  0$.
Here $L_X'$ embeds into a Levi subgroup of rank $1$ less than the rank of $C^i $ and it is easy to argue
as in Lemma \ref{murestriction} that 
$\mu^i \in \{ 0, \lambda_1^i,  \lambda_2^i, 2\lambda_1^i, \lambda_3^i, 3\lambda_1^i \}$ or the dual of one of these, 
where $3\lambda_1^i$ (or its dual) can only occur if $L_X' = A_2.$  

Now $\wedge^3(10\ldots  01 + 0\ldots  0) = \wedge^3(10\ldots  01) + \wedge^2(10\ldots  01).$   It is  shown in Lemma \ref{wedge3adj}  that $\wedge^3(10\ldots  01)\supseteq  (010\ldots  02) $  (or $(03)$ if $l=2$)  and this irreducible is also  clearly present in
 $\wedge^2(10 \ldots  01)$.  Therefore, $\wedge^3(10\ldots  01 + 0\ldots  0)$ is not MF.  And $S^3(11 + 00) \supseteq 22^2.$ Therefore neither  of these representations
nor their duals are MF.  It then follows from Lemma \ref{notw2ors2} that $V_{C^i}(\mu^i)$ is a trivial, natural, or dual module.

Next consider Lemma \ref{1label}.  Let $\nu'$ be the next highest weight at level $i$.  First assume 
that $L(\nu')=1$. Suppose Lemma \ref{1label}(i) holds with $\nu = (d_{k-1}+c_k)\omega_{k-1}.$  If $k \le l,$  then $c_k = d_k$ and $\nu' =  \omega_{k-2} +(d_{k-1}+d_k-2)\omega_{k-1} + \omega_k$  is afforded by $ f_{k-1}f_k^{d_k-1}$  so this forces $k-2 = 0 = d_{k-1}+d_k-2.$ Therefore $k = 2$ and $d_{k-1} = d_k = 1.$  So $i = 1$ and  the level has the form 
$2\omega_1 \oplus \omega_2.$  The other alternative in Lemma \ref{1label}(i) is that $k = l+1.$  Here  $\nu' =  \omega_{l-1} +(d_l+c_{l+1}-2)\omega_l $ so that $d_l+c_{l+1}-2 = 0$, whence $d_l = c_{l+1} = 1$ and again $i = 1.  $ The level is the dual of the previous one, namely $2\omega_l \oplus \omega_{l -1}.$

Now assume that Lemma \ref{1label}(ii) holds, still under the assumption that $L(\nu')=1$. Then $k = 2$ and $\nu = (d_1 - c_1 + d_2)\omega_1.$  Here $\nu' = 
 (d_1 - c_1 + d_2 - 2)\omega_1 + \omega_2$  is afforded by $f_1^{c_1+1}f_2^{d_2-1}.$ Therefore, $d_1 - c_1 + d_2 - 2 = 0$ forcing $c_1 = d_1-1$ and $d_2 = 1.$  So in this case we are in the next to last level and the full level
decomposes as $2\omega_1 \oplus \omega_2$.  

At this point we are left with the case where $L(\nu')=2$, so by the above   $\nu' = 10\ldots  01.$  If Lemma \ref{1label}(i) holds with  $1 < k \le l $, then $k = l$ and  $\nu' $ is afforded by
  $ f_{l-1}f_l^{c_l-1}.$  Then   $\nu'  = \om_{l-2} + (d_{l-1}+c_l-2)\omega_{l-1} + \omega_l.$  This forces $l = 2$ or $3$.  In the latter case $d_2+ d_3 = 2$, so $d_2 = d_3 = 1$ and the level decomposes as $(020 + 101)$.  And in the former case $d_1+c_2 = 3$  and $(d_1,c_2) = (1,2)$ or $(2,1).$ In either case the level decomposes as $30 + 11.$
Now suppose $k = l+1.$  Here $\nu' $ is afforded by $f_l^1f_{l+1}^{c_{l+1}-1}$, so $\nu' =  \omega_{l-1} +
(d_l+c_{l+1}-2)\omega_l$ so  again $L_X' = A_2$  and again
$d_2+c_3 = 3.$ This time the level decomposes as $03 + 11$.

Finally, suppose Lemma \ref{1label}(ii) holds.  Then $\delta = d_1\omega_1 + d_2\omega_2$, so this forces $L_X' = A_2$ and $\nu' = 11$.  As $\nu' =  (d_1 - c_1 + d_2 -2)\omega_1 + \omega_2$
we have $d_1-c_1 + d_2 = 3$ and $(d_1-c_1,d_2) = (1,2)$ or $(2,1)$.  Here the
full decomposition of the level is $30 + 11.$  \hal

\begin{cor}\label{maintwo} Part {\rm (ii)} of Theorem $\ref{NORDU}$ holds.
\end{cor}

\pf Assume the hypothesis of the theorem. As $i\ne 0,k$, there are at least two $L_X'$-summands $\d_j^i$ in the $i$th level of $W$. As in the previous proof, let the two highest such weights be $\nu$ and $\nu'$. By assumption, $L(\nu)=L(\nu')=1$. At this point the arguments in paragraphs 6 and 7 of the previous proof give the conclusion. \hal

We now establish  lemmas that settle the  cases described at the end of Lemma \ref{trivnatordual}.
The next two lemmas concern the situation where $L_X'$ embeds in $C^i$ via the representation $2\om_1 \oplus \om_2$ or its dual.  We assume the former and introduce some extra terminology.  We have $L_X'$ embedded in the Levi factor $A_x \times A_y$ via representations with high  weights $2\om_1$ and $\om_2$, respectively. Set $\Pi(A_x) =
 \{\beta_1, \ldots, \beta_x\}$  and $\Pi(A_y) =
 \{\beta_{x+2}, \ldots, \beta_{x+y+1}\}$ with corresponding fundamental dominant weights $\lambda_1^i, \ldots, \lambda_x^i$ and  $\lambda_{x+2}^i, \ldots, \lambda_{x+y+1}^i.$
 
  In proving these results we will use the term  ``exceptional weight"   for $A_x$ or $A_y$ to indicate one of a certain number of specific weights or the dual of such a weight.  The weights are listed below.
\[
\begin{array}{ll}
\hbox{Excep. weights for }A_x: &     
\lambda_1^i + \lambda_j^i \,(j\le 7 \hbox{ or }j\ge x-5),\, 2\lambda_1^i + \lambda_x^i,\,3\lambda_1^i + \lambda_x^i,\,
 \\
 & \lambda_2^i + \lambda_{x-1}^i,\,2\lambda_2^i ,\, 3\lambda_2^i,\, \lambda_2^i + \lambda_3^i,\, 2\lambda_1^i + \lambda_2^i,\,3\lambda_1^i + \lambda_2^i \\
\hbox{Excep. weights for }A_y \,(l \ge 4): & \hbox{same as }A_x, \hbox{ with obvious change of notation,  but} \\
& \hbox{for }l = 4 \hbox{ add } \l_{x+2}^i + 2\l_{x+3}^i,\,  a\l_{x+2}^i + \l_{x+3}^i\, (\hbox{all } a),\\
&  a\l_{x+2}^i + \l_{x+y+1}^i\, (\hbox{all }a),\,a\l_{x+3}\, (a\le 5),\,2\l^i_{x+4},\, 2\l^i_{x+5}.\\     
\hbox{Excep. wts for }A_y = A_5\, (l = 3):&  a\l_j^i,\,a\l_j^i+\l_k^i \,(\hbox{all }a,j,k),\,11100,\,11001.
\end{array}
\]

\begin{lem}\label{longlem} Suppose  the embedding   $L_X' < C^i$ corresponds to $2\omega_1 \oplus \omega_2$ or its dual and $L_X' \ne A_2$.  Assume the  inductive hypothesis. Then $V_{C_i}(\mu^i)$ is a trivial module, the natural module,  or the dual of the natural module.
\end{lem}

\pf   Note that since $l\ge 3$ by hypothesis, $x\ge 9$; also $x>y$. 
Write  $\mu^i = c_1\lambda_1^i + \cdots + c_{r_i}\lambda_{r_i}^i.$   In view of Lemma \ref{notw2ors2} it will suffice to show that $\mu^i$ is  $0,\,\l_1^i,\,\l_2^i,\,2\l_1^i$ or the dual of one of these.

Since $V_{C^i}(\mu^i) \downarrow L_X'$ is MF, Proposition \ref{tensorprodMF} implies that  in each composition factor of $A_x \times A_y$ on $V^i$, one of the tensor factors  has the property that upon restriction to $L_X'$,  all the irreducibles have  highest weights with at most one nonzero label.
 
Let $\e_1, \e_2$ denote the restrictions of $\mu^i$ to the maximal tori of $A_x, A_y.$ The inductive hypothesis implies that 
the possibilities for $\e_1,\e_2$ are as follows, listed up to duals:
\begin{equation}\label{posse1e2}
\begin{array}{l}
\hbox{possibilities for }\e_1:\;c\l_1^i,\,\l_j^i, \hbox{ or exceptional weight} \\
\hbox{possibilities for }\e_2:\;c\l_{x+2}^i,\,\l_j^i, \hbox{ or exceptional weight.}
\end{array}
\end{equation}
This also holds for conjugates of $A_x$ and $ A_y$.
  
We next note that 
\begin{equation}\label{nextnote} 
(c_x, c_{x+1}, c_{x+2}) \in \{(0,0,0),\, (0,0,1),\, (0,1,0),\,(1,0,0)\}.
\end{equation}
 Indeed, otherwise,  the conjugate of $A_x$ with base $\{ \beta_5^i, \ldots, \beta_{x+4}^i \}$ will have a composition factor for which the restriction to the corresponding conjugate of $L_X'$ is not MF as it is not one of the possibilities  described in (\ref{posse1e2}). 
Note that this rules out  the possibilities $\e_1 = c\l_x^i$, $\e_2 = c\l_{x+2}^i$ for $c \ge 2.$ 

Next observe that by Lemma \ref{twolabs}, if $\e_1$ is as in (\ref{posse1e2}) and  is not $c\l_1^1$ ($c\le 2$) or the dual of this,  then $V_{A_x}(\e_1)\downarrow L_X'$ has a summand of highest weight $\nu$ with $L(\nu)\ge 2$. The same conclusion holds for $V_{A_y}(\e_2)\downarrow L_X'$ provided $l\ge 4$ and $\e_2$ is not $c\l_{x+2}^i$ ($c\le 2$) or the dual, by Lemma \ref{twolabs}.

\vspace{2mm}
\no {\bf Case $\e_1,\e_2 \ne 0$}

Suppose that $\e_1,\e_2 \ne 0$. Then the previous paragraph together with Proposition \ref{tensorprodMF} implies that one of the following holds:
\begin{itemize}
\item[(1)] $\e_1 = c\l_1^i$ or $c\l_x^i$ with $c\le 2$,
\item[(2)] $\e_2 = d\l_{x+2}^i$ or $d\l_{r_i}^i$ with $d \le 2$,
\item[(3)] $l=3$ and $\e_2= \l_3^i,\,c\l_{x+2}^i$ or $c\l_{x+6}^i$.
\end{itemize}
By assumption, there exist $i,j$ with  $i < x+1 < j$ such that $c_i, c_j \ne 0$.  If $c_{x+1} \ne 0,$ then by (\ref{nextnote}) we have $(c_x, c_{x+1}, c_{x+2}) =  (0,1,0).$ Then at level 1 for $V_{C^i}(\mu^i)\downarrow A_xA_y$ there is an irreducible module  $V_{A_x}(\e_1 + \lambda_x^i) \otimes V_{A_y}( \lambda_{x+2}^i + \e_2).$ This must be MF upon restriction to $L_X'$, so from the induction hypothesis and the previous paragraph we obtain a contradiction using  Proposition \ref{tensorprodMF}.  Therefore,  $c_{x+1} = 0.$ 
 
 We first rule out the following specific possibilities for $\mu^i$:
\begin{equation}\label{listmui}
\l_1^i + \l_{r_i}^i,\,  2\l_1^i + \l_{r_i}^i,\, 3\l_1^i + \l_{r_i}^i,\, \l_1^i + \l_{r_i-1}^i, \,2\l_1^i + \l_{r_i-1}^i,\,
2\l_1^i +2\l_{r_i}^i,\,\l_2^i + \l_{r_i-1}^i.
\end{equation}
In these cases writing $\l_1^i \otimes \l_{r_i}^i$ to denote $V_{C^i}(\l_1^i) \otimes V_{C^i}(\l_{r_i}^i)$ and so on, $V_{C^i}(\mu^i)$ has the form 
\[
\begin{array}{l}
(\l_1^i \otimes \l_{r_i}^i) - 0,\; (2\l_1^i \otimes \l_{r_i}^i) - \l_1^i,\; (3\l_1^i \otimes \l_{r_i}^i) - 2\l_1^i,\;
(\l_1^i \otimes \l_{r_i-1}^i) - \l_{r_i}^i,\\
 (2\l_1^i \otimes \l_{r_i-1}^i)-(\l_1^i\otimes \l_{r_i}^i)+0,\; (2\l_1^i \otimes 2\l_{r_i}^i) - (\l_1^i \otimes \l_{r_i}^i), \;
(\l_2^i \otimes \l_{r_i-1}^i) - (\l_1^i \otimes \l_{r_i}^i),
\end{array}
\]
 respectively.  It is straightforward to see that in each case the restriction to $L_X'$ fails to be MF. For example, 
$(\l_1^i \otimes \l_{r_i}^i)\downarrow L_X' = (2\om_1\oplus \om_2) \otimes(2\om_l\oplus \om_{l-1}) \supseteq (\om_1+\om_l)^2.$ Also $\l_2^i \downarrow L_X' \supseteq (2\om_1+\om_2)^2$ and $3\l_1^i \downarrow L_X' \supseteq (\om_1+\om_2 + \om_3)^2. $ The duals of these  configurations also fail to be MF upon restriction to $L_X'$.
 
We now work our way through the list of possibilities for $\e_2$ given by (\ref{posse1e2}), bearing in mind the restrictions implied by (\ref{nextnote}) and (\ref{listmui}).

First suppose $\e_2 = \l_j^i$ with $j > x+2.$ Then at level 1 of $V_{C^i}(\mu^i)$, 
there is an irreducible summand  afforded by
 $\mu^i - \psi - \beta_{x+1}^i$, where $\psi$ is a sum of fundamental roots in $\Pi(A_x),$ which must restrict 
 to $A_y$ with an exceptional highest weight.  Excluding $\e_1 = c\l_1^i \,(c \le 3)$ or $\l_2^i$, Proposition \ref{tensorprodMF} implies that $\psi$ can be chosen such that the resulting tensor product fails to be MF upon restriction to $L_X'$. Consider  the excluded cases.  Here level 1 contains a composition factor afforded by $\mu^i - \psi - \beta_{x+1}^i$, where this time $\psi$ is a sum of fundamental roots in $\Pi(A_y)$, and Proposition \ref{tensorprodMF} gives a contradiction unless $\e_2 = \l_{r_i-1}^i,\,\l_{r_i}^i$ or $\l_{x+3}^i$, with $l=3$ in the last case. In the last case, a Magma computation shows that level 1 is not MF for $L_X'$. In the other cases, Proposition \ref{tensorprodMF} shows that level 0 is not MF for $L_X'$ if $\e_1  = 3\l_1^i$ and  $\e_2 = \l_{r_i-1}^i$. All the remaining cases for $\mu^i$ are in the list (\ref{listmui}). Therefore $\e_2 \ne \l_j^i$ for $j > x+2.$
 
 Next suppose $\e_2 = d\l_{r_i}^i$ with $d > 1.$ We allow $l = 3$  where all irreducibles in the restriction to $L_X’$  have at most 1 nonzero label.  At level 1 there is an irreducible afforded by $\mu^i - \psi - \beta_{x+1}^i$, where $\psi$ is a sum of fundamental roots in $\Pi(A_x)$ and this restricts to $A_y$ as $\l_{x+2} + d\l_{r_i}^i.$ This will not be MF unless
$\e_1 = c\l_1^i$ for $c \le 3$. So we have reduced to the case $\mu^i = c\l_1^i + d\l_{r_i}^i.$ Taking duals and using the same argument we find that $d \le 3$. Therefore in view of (\ref{listmui}) and taking duals, if necessary, we can assume that $(c,d) =(3,2)$ or $(3,3)$. In both
 cases we see that  the restriction of $V_{A_x}(\e_1) \otimes V_{A_y}(\e_2)$ to $L_X'$ fails to be MF.  Indeed
   the restriction contains $((60\ldots 0) + (220\ldots 0)) \otimes (0 \ldots  0d0)$ and two  applications of Lemma \ref{aibj12} shows that this contains $(40\ldots 0(d-2)2)^2.$

   Now suppose $\e_2 = c\lambda_{x+2}^i$. Then at level 1 there is an irreducible afforded by $\mu^i - \beta_{x+1}^i - \beta_{x+2}^i$ which yields a contradiction to Proposition \ref{tensorprodMF} unless $\e_1 = \l_x^i$.  In the latter case,  level 1 contains $\wedge^2(2\o_l) \otimes S^c(\o_2)\otimes \o_2$, which can be seen to be non-MF using Proposition 
\ref{tensorprodMF}.

 We have now reduced to the case where $\e_2$ is an exceptional weight for $A_y$.  Then Proposition \ref{tensorprodMF} together with Lemma \ref{twolabs} gives a contradiction at level 0 unless $\e_1 = c\l_1^i\,(c\le 2)$ or $c\l_x^i$.   In the first case
  we can take duals to reduce to the case $\e_2 = c\l_{r_i}^i$ previously considered.   In the last case Proposition \ref{tensorprodMF} yields a contradiction  at level 1  using  the weight $\mu^i- (\beta_x^i+  \beta_{x+1}^i)$.  

\vspace{2mm}
\no {\bf Case $\e_1\ne 0$, $\e_2=0$}

Assume $\e_1\ne 0$, $\e_2=0$, and consider the possibilities for $\e_1$ given by (\ref{posse1e2}). 
As mentioned in the first paragraph, the conclusion of the lemma holds if $\e_1 = \l_1^i$, $\l_2^i$ or $2\l_1^i$ provided $c_{x+1} = 0.$  And if one of these cases occurs with    $c_{x+1} \ne 0$, then taking duals we are back in the case $\e_1,\e_2\ne 0$ already considered. So assume $\e_1 \ne \l_1^i$, $\l_2^i$ or $2\l_1^i$.
 
 Let   $s$ be maximal such that  $s \le x$ and $c_s \ne 0.$  Suppose $s > 1.$  

If $c_s > 1$  then we can find a conjugate of $A_x$ with fundamental system a subset of $\Pi(C^i)$ for which $\mu^i$ restricts to a highest weight which is not in the list of possibilities for $\e_1$ given in (\ref{posse1e2}), 
unless $s=2$ and  $c_s = 2$ or $3$. In these cases $\mu^i = c\l_2^i$ ($c = 2,3$) (note that $c_{x+1} = 0$, as above). First suppose $\mu^i = 3\l_2^i. $ Then at
 level 3 the module $V_{A_x}(3\l_1^i) \otimes V_{A_y}(3\l_{x+2}^i)$ appears, and we can use Proposition \ref{tensorprodMF} (together with Magma for $l=3$) to get  a contradiction. 
Likewise, if $\mu^i = 2\l_2^i$ then Lemma \ref{nonMF}(iv) shows that level 1 is not MF for $L_X'$.  Therefore  $c_s = 1$.  
  
From (\ref{nextnote}) we know that $c_{x+1} \le 1$  and that if equality holds, then $c_x = 0,$ forcing $s < x.$
Suppose $c_{x+1} = 1.$  There  is an irreducible at level $2$ afforded by $\mu^i - \beta_x^i - 2\beta_{x+1}^i - \beta_{x+2}^i.$   This affords $V_{A_x}(\e_1 + \l_{x-1}^i) \otimes V_{A_y}(\l_{x+3}^i).$  Restricting to $L_X'$ the second tensor factor affords $(1010\ldots  0).$  In addition, either $V_{A_x}(\e_1 + \l_x^i)$ is    not one of the modules in the list of possibilities for $\e_1$ in (\ref{posse1e2}), or if it  is, then the restriction to $L_X'$ involves an irreducible with at least two nontrivial labels.  Either way we have a contradiction.  Therefore $c_{x+1} = 0$ (still assuming 
$s > 1.)$

If $c_{s-1}=0$ then there is an irreducible at level 2 afforded
 by $\mu^i - (\beta_{s-1}^i + 2\beta_s^i + \cdots + 2\beta_{x+1}^i + \beta_{x+2}^i)$ and we obtain the same
 contradiction unless $\e_1 = \lambda_2^i$, $\lambda_3^i$ or $\l_1^i+\l_3^i$.
If $V_{C^i}(\mu^i)$ is the wedge square of the natural module, then the restriction to $L_X'$ is not MF as $\wedge^2(2\o_1\oplus \o_2) \supseteq (2\o_1+\o_2)^2$.   Also the wedge cube of the natural module  contains a submodule which restricts to $L_X'$ as the sum of  $\wedge^2(2\o_1) \otimes \o_2 = (2\o_1+\o_2)\otimes \o_2$  and $2\o_1 \otimes \wedge^2(\o_2) =2\o_1 \otimes (\o_1+\o_3)$.  Each of these tensor products contains a summand of highest weight $\o_1+\o_2+\o_3$,   a contradiction.  Finally, if $\e_1 = \l_1^i+\l_3^i$, then 
$V_{C^i}(\mu^i) \supseteq V_{A_x}(\l_1^i)\otimes V_{A_x}(\l_2^i)\otimes V_{A_y}(\l_{x+2}^i)$ in level 1 for $A_x\times A_y$, which is not MF.

Hence $c_{s-1}\ne 0$, and so $\mu^i$ is either 
$a\l_1^i + \l_2^i$ or $\l_2^i+\l_3^i$ (note that it is not the dual of either of these, as can be seen by considering a suitable conjugate of $A_x$). In the first case, level 1 restricted to $L_X'$ contains a summand $S^a(2\o_1)\otimes 2\o_1 \otimes \o_2$, which is not MF by Lemma \ref{nonemf} for $a\ge 2$, and is also non-MF for $a=1$ as it contains $(3\o_1+\o_3)^2$. And if $\mu^i = \l_2^i+\l_3^i$, level 1 contains a summand $V_{A_x}(\l_1^i+\l_3^i) \otimes V_{A_y}(\l_{x+2}^i)$, whose restriction to $L_X'$ is not MF by Lemma \ref{nonMF}(4).

Therefore $s = 1$. It follows that $\e_1 = c\l_1^i$.  If $c_{x+1} \ne 0$, then taking duals we are back in
the case where $\e_1,\e_2 \ne 0$, so $c_{x+1} = 0$ and $\mu^i = c\l_1^i$. 
If $c\le 2$ then the conclusion of the lemma holds, so assume $c\ge 3$.
Then if $M$ denotes the natural module for $C^i$, $V_{C^i}(\mu^i)$ is the symmetric power $S^cM$.   
If $c = 3$, then as noted earlier $S^3M\downarrow L_X'$ contains $(\o_1+\o_2+\o_3)^2$ as a submodule, a contradiction. 
 
Hence  $c > 3$. Then $S^cM\downarrow L_X'$  contains submodules of the form 
$S^{c-1}(2\o_1) \otimes \o_2$ and also $S^{c-3}(2\o_1) \otimes S^3(\o_2)$.  The former tensor product contains 
  $((2c-6)\o_1+2\o_2) \otimes \o_2$, which contains an irreducible submodule of highest weight 
$(2c-6)\o_1+3\o_2$.  The latter tensor product contains $(2c-6)\o_1 \otimes  3\o_2$ which also contains 
$(2c-6)\o_1+3\o_2$. This is a  contradiction.  

This completes the case where $\e_1\ne 0, \e_2=0$. 

Finally, if either $\e_1= 0, \e_2\ne 0$ or $\e_1= \e_2=0$, then dualizing gives one of the cases that we have already considered -- namely, $\e_1,\e_2\ne 0$ or $\e_1\ne 0, \e_2=0$.  This completes the proof of the lemma. \hal

The next lemma handles the case where $L_X' = A_2$, excluded in the previous result.

\begin{lem} Suppose $L_X' = A_2$ and   $L_X' < C^i$ corresponds to $20 + 01$ or its dual.
Then $V_{C_i}(\mu^i)$ is a trivial module, the natural module, or the dual of the natural module. 
\end{lem}

\pf  We may assume the embedding to be $20 + 01$ so that $L_X'$ embeds into the Levi subgroup $A_x \times A_y$ of $C^i$,  where $A_x = A_5$ and $A_y = A_2.$  Write $\mu^i = c_1\lambda_1^i + \cdots + c_8\lambda_8^i$ and let $\e_1$ and $\e_2$  be the restrictions of $\mu^i$ to $A_x$ and $A_y,$ respectively. 

Assume first that $\e_1$ is one of the exceptional weights listed before the statement of Lemma \ref{longlem}.
The restrictions of these to $L_X'$ are given by Lemma \ref{A2info}. Let $\e_2 = ab$. 

Note that $\e_1 \ne 10001$ or $10002$, since otherwise at level 1 we would have 
$((10000\otimes 00010)\downarrow L_X') \otimes (b,a+1)$ or $(10011\downarrow L_X')\otimes (b,a+1)$, and neither of these is MF. And if $\e_1 = 00011$ then at level 1 we would have $((00101) + (00020)) \otimes (b,a+1)$ which contains $(13)^2 \otimes (b,a+1)$, which is again a contradiction.  
Also, if $c_6,c_7,c_8$ are not all zero, then there
is an irreducible module at level $1$ for which the highest weight restricts to $A_x$ as $\e_1 + \lambda_5^i.$ 
But then the restriction of this to $L_X'$ is not MF, a contradiction.  Therefore, $c_6 = c_7 = c_8 = 0.$

For the remaining exceptional weights and their duals, we can again 
produce irreducible $A_xA_y$-summands at level 1, the sum of whose restrictions to $L_X'$ is not MF; in most cases a single summand suffices, but in a few, two summands are needed. Here is an example of a case where two are required, $\e_1 = 21000$. At level 1 there are $A_xA_y$-summands $11000\otimes 10$ and $30000\otimes 10$, and these sum to 
$(20000\otimes 10000)\otimes 10$. The restriction of this to $L_X'$ contains $22^2$. All other cases for $\e_1$ an exceptional weight are similar.

 At this point the exceptional weights and their duals have been ruled out as possibilities for $\e_1$.  We can
 also assume that the exceptional modules do not occur in $V_{C^i}(\mu^i)^*$.
Now consider the remaining possibilities for $\e_1$. 
There is no restriction on $\e_2$, while $\e_1 \in \{0, c\lambda_1^i, c\lambda_5^i,
\lambda_2^i, \lambda_3^i, \lambda_4^i \}$.  Taking duals we obtain certain restrictions.  In particular, the fact that the exceptional cases for $\e_1$ have been excluded, when applied to  $V_{C^i}(\mu^i)^*$, shows that  
$(c_6,c_7,c_8) \in \{ (0,0,0), (0,1,0),  (0,0,d), (1,0,0) \}.$  

First suppose that $\e_1 = \e_2 = 0.$   If $c_6 \ne 0$ then by the above $c_6 = 1.$  But this implies
that  $V_{C^i}(\mu^i)^*$ is  the third wedge of the natural module, which is not MF when restricted to $L_X'$, 
since $(30)^2$ appears.  Therefore we may assume that either $\e_1$ or $\e_2$ is nonzero.  

Next assume both $\e_1 , \e_2 \ne 0$.  
First suppose that $(c_6,c_7,c_8) = (0,1,0).$  Suppose $\e_1 = c\l_5^i$ and consider   $V_{C^i}(\mu^i)^*$. 
 If $c > 1$ the first tensor factor is not MF when restricted to $L_X'$.  And if $c = 1$, then we have one of the
exceptional weights and we again have a contradiction. Therefore, $\e_1 \ne c\lambda_5^i.$
The same argument shows that $\e_1 \ne \lambda_4^i.$  If $\e_1 = \lambda_3^i$, then  $V_{C^i}(\mu^i)^*$
 has highest weight $\lambda_2^i + \lambda_6^i$ and at level $1$ we have the irreducible $(01001) \otimes 10.$  The
restriction to $L_X'$ contains $(13)^2$, so this is impossible.  If $\e_1 = \lambda_2^i,$ then
$V_{C^i}(\mu^i) = \lambda_2^i + \lambda_7^i =  (\lambda_2^i \otimes \lambda_7^i) - (\lambda_1^i \otimes \lambda_8^i).$  But this is not
MF upon restriction to $L_X'$: indeed, $(33)^4$ occurs.  Now suppose $\e_1 = c\lambda_1^i.$
 Then at level $1$ the module $(c0001) \otimes 01$ appears.   If $c > 3$, the first factor is not MF upon restriction to $L_X'$.  If $c = 3$, then the restriction to $L_X'$ contains $(23)^2.$ If $c = 2$ the restriction to $L_X'$ contains $(03)^2$.  
 So suppose $c = 1$ where $V_{C^i}(\mu^i) = \lambda_1^i + \lambda_7^i = (\lambda_1^i \otimes \lambda_7^i) - \lambda_8^i.$  Here the restriction contains $(32)^2$ and so is not MF.   Therefore $(c_6,c_7,c_8) = (0,0,d).$
 
Suppose $\e_1 = c\lambda_1^i.$  Taking duals we may assume $c \ge d.$ The argument in the preceding paragraph implies that $c \le 3$.  At level $d$ the module $(c000d) \otimes 00$ appears and the first factor is not MF upon restriction to $L_X'$ unless 
$(c,d) = (1,1), (2,1),$  or $ (3,1).$  If $(c,d) = (1,1)$, then $V_{C^i}(\mu^i)\downarrow L_X'$ contains $(11)^2$.  
If $c = 2$, then 
$\mu^i = 2\lambda_1^i + \lambda_8^i = (2\lambda_1^i \otimes  \lambda_8^i) - \lambda_1^i$ and the restriction to $L_X'$
contains $(04)^2$, a contradiction.  And if $c = 3$, then $\mu^i = 3\lambda_1^i + \lambda_8^i = (3\lambda_1^i \otimes  \lambda_8^i) - 2\lambda_1^i$ and the restriction to $L_X'$ contains $(13)^2,$ again a contradiction.

Now suppose $\e_1 = c\lambda_5^i.$  Then taking duals we see that
$c = d = 1$ and we have one of the exceptional cases treated above. If $\e_1 = \lambda_4^i$, then 
$V_{C^i}(\mu^i)^* = d0001000$, a case already treated.   If
$\e_1 = \lambda_3^i$, then $(V_{C^i}(\mu^i))^* = 
d0000100$ and at level $1$ we have the module $(d0001) \otimes (10).$  The main result for $A_2$ (proved in Chapter \ref{casea2}) forces
$d \le 3$. And for $d = 1,2,3$,  restricting to $L_X'$ we get $(12)^2$, $(32)^2,$  $(52)^2$, respectively.  The final case here is $\e_1 = \lambda_2^i.$ Then $(V_{C^i}(\mu^i))^* = 
d0000010$, and again this has already been handled.

Thus we may now assume just one of $\e_1$ and $\e_2$ is nonzero, and taking duals we may take it that $\e_1 \ne 0.$  As above we have $\e_1 \in \{0, c\lambda_1^i, c\lambda_5^i,
\lambda_2^i, \lambda_3^i, \lambda_4^i \}$ and $c_6 = 0$ or $1$.  If $\e_1 = c\l_5^i$, then 
$V_{C^i}(\mu^i)^* = 00c_6c0000$
so $c = 1,$  as otherwise either the restriction to $L_X'$ is not MF or we are in an exceptional case. The case $c_6 = 1$ is the dual of an exceptional case, so this is not possible.  And if $c_6 = 0$,
then $V_{C^i}(\mu^i)^*$ is the fourth wedge of $\lambda_1^i$ and the restriction to $L_X'$ contains $(12)^2.$  Therefore,
$\e_1 \ne c\lambda_5^i.$  Suppose $\e_1 = c\lambda_1^i.$  If $c_6 = 1,$ then taking duals we are back in the case where $\epsilon_1 \ne 0 \ne \epsilon_2.$  So we can assume $c_6 = 0. $ The lemma follows if $c = 1,$ so assume $c > 1.$
If $c = 2$, then the restriction contains $(02)^2$ and if $c > 2$ the restriction contains
 $S^c(20) + (S^{c-2}(20) \otimes S^2(01))$ which contains $((2c-4)2)^2$, a contradiction. 
 
  The remaining cases are $\e_1 = \lambda_2^i, \lambda_3^i,$ or $\lambda_4^i.$   If $c_6 = 0$, then  
$V_{C^i}(\mu^i)$ is the corresponding wedge of the natural module and an easy check shows that these wedges contain $(21)^2, (30)^2, (31)^2$, respectively, a contradiction.  
Suppose $c_6 = 1.$  If $\e_1 = \lambda_2^i$, then taking duals we
 again have a case with $\epsilon_1 \ne 0 \ne \epsilon_2.$ And if $\e_1 = \lambda_4^i,$ the dual involves an exceptional module.  Finally, assume $\e_1 = \lambda_3^i.$  Then at level 1 the module $(00101) \otimes (10)$
 appears and the restriction to $L_X'$ contains $(22)^2$, a  final contradiction. \hal

The next lemma deals with the last two cases given at the end of Lemma \ref{trivnatordual}.

\begin{lem}\label{finaltriv} Suppose that either  $L_X' = A_2$ and the embedding $L_X' < C^i$ corresponds to $30 + 11$ or  $03 + 11$, or $L_X' = A_3$ and $L_X' < C^i$ corresponds to $101 + 020$.
Then $V_{C^i}(\mu^i)$ is a trivial module, a natural module, or the dual of a natural module.
\end{lem}

\pf First consider $L_X' = A_2$.  By assumption, the embedding is $30 + 11.$ Here $C^i = A_{17}$ and $L_X'$ is inside the Levi subgroup $A_9 \times A_7.$  By the main result for $A_2$ (proved in Chapter \ref{casea2}), 
 the possible restrictions of $\mu^i$ to $A_9$ and $A_7$ are as follows:
\begin{equation}\label{posswts}
\begin{array}{ll}
\hbox{to }A_9:& \;(c,0,\ldots,0)\,(c\le 4),\,(0,1,0,\ldots,0),\,(0,0,1,\ldots,0),\,(0,0,0,1,\ldots,0),\\
                         & \; (0,0,0,0,1,\ldots,0),\,(1,1,0,\ldots,0),\,(1,0,\ldots,0,1)\hbox{ and duals} \\
\hbox{to }A_7:& \;(c,0,\ldots,0)\,(c\le 3),\,(0,1,0,\ldots,0),\,(0,0,1,0,\ldots,0) \hbox{ and duals}.
\end{array}
\end{equation}

 In view of Lemma \ref{notw2ors2} it will suffice to show that 
$V_{C^i}(\mu^i)$ is   trivial, the natural module, the wedge square of the natural module, the symmetric square
of the natural module, or the dual of one of these.

 As above, we know that
all irreducible summands of conjugates of these Levi factors $A_9$ and $A_7$ are MF when restricted to the appropriate
$A_2$ subgroup.  Write $\mu^i = c_1\lambda_1^i + \cdots + c_{17}\lambda_{17}^i$.  If $c_j \ne 0$ for
some $4 \le j \le 14$, then there exists a conjugate of $A_7$ built from the same fundamental system and having a composition factor which is not in the list (\ref{posswts}).  Therefore $c_j = 0$ for $4 \le j \le 14.$  

Next suppose that there exist $c_j,c_k \ne 0$ with $j \le 3$ and $k\ge 15$. Taking $j$ maximal for this we see that there is a maximal vector at level $1$ of $V_{C^i}(\mu ^i)$ for which the restriction to $A_7$ has highest weight $(1000c_{15}c_{16}c_{17})$ and this contradicts the list (\ref{posswts}).  Therefore the restriction of $\mu^i$ to one or the other of the Levi factors $A_9,A_7$ is  trivial.  Taking dual modules we may assume the restriction to $A_9$ is trivial.  

The restriction to the $A_7$ factor is $(0000c_{15}c_{16}c_{17}).$  If the $3$-tuple $(c_{15},c_{16},c_{17})$ is one of $(0,0,0)$, $(0,0,1)$, $(0,1,0)$, or $(0,0,2)$, then, as noted in the second paragraph, the result follows. So assume  $(c_{15},c_{16},c_{17}) = (1,0,0)$ or $(0,0,3)$. Then 
$\wedge^3(30 +11) \supseteq 22^2$ and $S^3(30 +11) \supseteq 33^2,$ so these are not MF, and similarly for their duals.  This is a contradiction. 

Finally suppose $L_X' = A_3$ and the embedding is $020 +101$.  Here $C^i = A_{34}$ and $L_X'$ is contained in the Levi subgroup
$A_{19} + A_{14}.$  The argument is as above with fewer special cases to consider.  Ultimately
we  check that $\wedge^3(020 + 101) \supseteq   (210)^2$ and hence is not MF.   \hal

\vspace{2mm}
The proof of Theorem \ref{NORDU}(i) follows from Lemmas \ref{trivnatordual} -- \ref{finaltriv}.  
And part (ii) follows from Corollary \ref{maintwo}.

\section{Proof of Theorem \ref{MUIZERO}}

We now begin the proof of Theorem \ref{MUIZERO}.
Assume the hypotheses  of the theorem and suppose
$\mu^i \ne 0$, where $i\ne 0,k$.  Then Theorem \ref{NORDU} shows that  $\mu^i $ affords the natural or dual module for $C^i$.

We first  assume $\mu^i = \lambda_1^i.$   Later we indicate the changes required for the dual case.

We first consider the case  $i \ge 2.$  We claim $\mu^{i-1} = 0.$  For otherwise,
$V_{C^{i-1}}(\mu^{i-1})$ is a natural or dual module, and part (ii) of Theorem  \ref{NORDU} implies that either
the restrictions to $L_X'$ of $V_{C^i}(\mu^i)$ and $V_{C^{i-1}}(\mu^{i-1})$ each contain an irreducible module
whose highest weight has at least two nontrivial coefficients, or we have one of the following exceptional cases: 
\begin{itemize}
\item[(a)] $\delta = \omega_l + d_{l+1}\omega_{l+1}$, $i = 2$ and $W^2(Q_X) = V_{L_X'}(2\om_l) + V_{L_X'}(\om_{l-1})$, or 
\item[(b)] $\delta =d_1 \omega_1 + \omega_2,$ $i = d_1$, and  $W^{d_1}(Q_X) = V_{L_X'}(2\om_1) + V_{L_X'}(\om_2)$. 
\end{itemize}
 Excluding the exceptional cases, Proposition \ref{tensorprodMF} implies that $V^1$ is not MF, a contradiction.  

Now assume (a) holds.  The assumption $i \ge 2$ implies $d_{l+1} > 1.$  Therefore $W^2(Q_X) = (0 \ldots  02) + (0 \ldots  010)$ and $W^3(Q_X) = (0 \ldots  03) + (0 \ldots  011).$  If $\mu^1 = \lambda_1^1$, then
$V_{C^1}(\mu^1) \otimes V_{C^2}(\mu^2)  \supseteq ((0 \ldots  02) \otimes (0 \ldots  011)) + ((0 \ldots  010) \otimes (0 \ldots  03)) \supseteq (0 \ldots  0 13)^2$, a contradiction.  On the other hand if $\mu^1 = \lambda_{r_1}^1$, then
$ V_{C^1}(\mu^1) \otimes V_{C^2}(\mu^2)  \supseteq ((20\ldots  0) \otimes(0 \ldots  03)) + ((010 \ldots  0) \otimes (0 \ldots  011)) \supseteq (10\ldots  02)^2$, again a contradiction.
So here, $\mu^{i-1} = 0.$  If (b) holds, then we must  have $d_1 > 1$ in view of our assumption $i \ge 2.$  This time 
$W^{d_1 -1}(Q_X) = (110\ldots 0) + (30\ldots 0)$ and the situation is dual to the one just considered and hence
again gives a contradiction.  This establishes the claim that $\mu^{i-1} = 0.$ 

Let ${\kappa}$ be the fundamental root between $C^{i-1}$ and $C^i$ and consider $V_{\kappa}^2(Q_Y)$.
Let $\xi = \bigotimes_{j \ne i}V_{C^j}(\mu^j)$, so that $V^1(Q_Y) = \xi \otimes V_{C^i}(\lambda_1^i).$  Then
$V_{{\kappa}}^2(Q_Y) \supseteq \xi \otimes V_{C^{i-1}}(\lambda_{r_{i-1}}^{i-1}) \otimes V_{C^i}(\lambda_2^i).$

At this point we prove a lemma that applies in this situation and at several points in Chapters 15 -17.

Suppose that $\nu$ is the highest weight of the irreducible  at level $i$ (for the action of $L_X'$ on $W$) with largest monomial in the
ordering.  As in the discussion preceding  Lemma \ref{1label}, $\nu = \delta + f_j^{c_j} \cdots f_k^{c_k}$ and
if $j<k$ then $c_s = d_s$ for all $s > j.$  In particular, $d_s = 0$ for $s > k.$

If $c_j = d_j$ or if $j = k$, then there must exist $e < j$ with $d_e \ne 0.$   For otherwise in the first case the $i$th level is the last level and in the second case $\delta = d_k\omega_k$, both of which contradict the hypothesis.
Taking $e$ maximal there is  another irreducible at level $i$ with highest weight $\nu' = \delta + f_ef_j^{d_j-1}f_{j+1}^{d_{j+1}} \cdots f_k^{d_k}$  or $\delta + f_ef_k^{c_k-1}$, respectively. Then $\nu' = \nu - \psi$, where $\psi = f_{\beta_e^i + \cdots +\beta_{j-1}^i}.$
Now suppose $j < k$ and $c_j < d_j.$ Let $j < e \le k$ be minimal with $d_e \ne 0.$  Here
there is  an irreducible with highest weight $\nu' = \delta + f_j^{c_j+1}f_e^{d_e-1}f_{e+1}^{d_{e+1}} \cdots f_k^{d_k}$  and  $\nu' = \nu - \psi$ with $\psi = f_{\beta_j^i + \cdots  + \beta_{e-1}^i}.$

\begin{lem}\label{mult2}
Assume that $\d \ne r\o_j$, and let $1 \le i \le k-1$. Let $\nu$ be the highest weight of the 
irreducible  at level $i$  with largest monomial. As above, there is a second largest monomial at level 
$i$, affording highest weight of the form $ \nu-\psi,$ for $\psi \in \Sigma^+(L_X')$. Moreover,
\begin{itemize}
\item[{\rm (i)}] $\wedge^2(W^{i+1}(Q_X))$ has a composition factor of multiplicity at least $2$ with highest weight $2\nu-\psi$, 
and $S(2\nu-\psi) \ge 2S(\nu)-1$ or $2S(\nu)-2,$ where the latter occurs only if $\delta = a\omega_1 + b\omega_{l+1}$ and $\nu = \delta + f_1^{c_1}f_{l+1}^{c_{l+1}}.$
\item[{\rm (ii)}] If $\psi$ is not a fundamental root, then $S^2(W^{i+1}(Q_X))$ has a composition factor of multiplicity at least $2$ with highest weight $2\nu-\psi$, and $S(2\nu-\psi) \ge 2S(\nu)-1$ or $2S(\nu)-2,$ where the latter occurs only if $\delta = a\omega_1 + b\omega_{l+1}$ and $\nu = \delta + f_1^{c_1}f_{l+1}^{c_{l+1}}.$
\item[{\rm (iii)}] If $\psi$ is a fundamental root, then $S^2(W^{i+1}(Q_X))$ has a composition factor of multiplicity at least $2$ with highest weight $2\nu-2\psi$, and $S(2\nu-2\psi) = 2S(\nu)$ or $2S(\nu)-2$, the latter only if $\psi$ is an end node.
\end{itemize}
\end{lem}

\pf  We have $\nu=\d+f_j^{c_j}\cdots f_k^{c_k}$ where each $c_j \le d_j$ and $\sum c_j = r.$  The hypothesis on $i$ implies that  $\nu \ne \d$ and $\nu \ne \d+f_1^{d_1}\cdots f_{l+1}^{d_{l+1}}.$  

(i)  Now $\wedge^2(W^{i+1}(Q_X)) \supseteq \wedge^2( V_{L_X'}(\nu) + V_{L_X'}(\nu - \psi))  \supseteq 
\wedge^2( V_{L_X'}(\nu) + ( V_{L_X'}(\nu) \otimes V_{L_X'}(\nu - \psi))).$  We claim that the first
summand contains an irreducible of highest weight $2\nu - \psi$.  The argument varies only slightly
in the two cases above. First assume  $c_j = d_j$ or $j = k.$  Let $e< j$ be maximal
with $d_e \ne 0.$  If $e < j-1$, then   $\nu = \ldots  d_e 0 \ldots  0 c_j \ldots $, where the $d_e, c_j$
appear at nodes $e, j-1$ respectively.  And if $e = j-1,$ then $\nu = \ldots  (d_{j-1}+c_j) \ldots .$
The claim is obvious in the latter situation and follows from a simple weight count in the former, noting that the only dominant weights of the wedge square strictly above $2\nu - \psi$ are $2\nu, 2\nu- \alpha_e$ and 
$2\nu - \alpha_{j-1}.$ 
  This gives the claim, which implies that $\wedge^2(W^{i+1}(Q_X)) \supseteq V_{L_X'}(2\nu - \psi)^2.$
Now consider $S(2\nu - \psi).$  The discussion preceding the proof shows that $\psi$ is a root, but not
the highest root of $\Sigma(L_X')$ unless $\nu =\d+f_1^{c_1}f_{l+1}^{c_{l+1}}$.  Therefore, $S(2\nu - \psi) \ge S(2\nu)-1 = 2S(\nu)-1$ or $2S(\nu)-2$ in the exceptional case.
The second case is similar.  If $j < e-1$ $\nu = \ldots   c_j (d_j-c_j)0\ldots  0 d_e \ldots $, where the $d_j-c_j, d_e$ appear at nodes $ j, e-1,$ respectively.  From here the argument is the same.

(ii), (iii)  Here we argue exactly as above, using the symmetric square rather than the wedge square.
The argument is  the same and we get (ii), unless $e = j-1$ (if $c_j = d_j$ or $j = k$)  or $e = j+1$ (otherwise).  In these cases  $\psi = \alpha_{j-1}$ or $\alpha_j$, respectively, and $S^2(V_{L_X'}(\nu))$ does not contain $V_{L_X'}(2\nu - \psi).$  However, here  $S^2(V_{L_X'}(\nu)) \supseteq V_{L_X'}(2\nu - 2\psi)$ which is also the high
weight of $S^2(V_{L_X'}(\nu - \psi))$, so that here $V_{L_X'}(2\nu - 2\psi)$ appears with multiplicity $2$.
As $\psi$ is a simple root we check $S$-values and obtain (iii). \hal

Lemma \ref{mult2} implies that  $V_{C^i}(\lambda_2^i) \downarrow L_X' \supseteq (2\nu - \psi)^2.$ Now consider $V_{C^{i-1}}(\lambda_{r_{i-1}}^{i-1}) = V_{C^{i-1}}(\lambda_1^{i-1})^*.$  Let $\gamma$ be
the highest weight of the irreducible summand of  $W^i(Q_X)$ corresponding to the largest monomial
in the ordering. Note that $\g\ne 0$ by the discussion following the proof of Lemma \ref{murestriction}. 
By Corollary \ref{V^2(Q_X)}, $S(\g) = S(V_{C^{i-1}}(\lambda_1^{i-1}) \downarrow L_X')$, and also 
$S(\gamma) \ge S(\nu) - 1$. 
 The $S$-values for a highest weight and the highest weight of the dual module
are equal.   Therefore, if $\epsilon$ is the highest weight of an irreducible appearing in $\xi \downarrow L_X'$ with 
largest $S$-value, then $S(V^1) = S(\e)+S(\nu)$ and also $V_{\kappa}^2(Q_Y) \downarrow L_X'$ contains an irreducible module appearing
with multiplicity at least $2$ and
whose highest weight has $S$-value  $S(\epsilon)  + S(\gamma) + S(2\nu - \psi).$
On the other hand, by Corollary \ref{V^2(Q_X)} the largest $S$-value among  highest weights of irreducibles in $V^2$ arising from $V^1(Q_Y)$  is at most $S(\epsilon) + S(\nu) + 1.$  So Proposition \ref{sbd} implies that 
$S(\epsilon) + S(\nu) + 1 \ge  S(\epsilon)  + S(\gamma) + S(2\nu - \psi).$   As 
$S(2\nu - \psi) = S(\nu) +S(\nu - \psi)$ this yields  $S(\nu - \psi) + S(\gamma) \le 1.$
As $\nu - \psi$ is a dominant weight, it follows that $S(\gamma) \le 1$.  This 
implies that $\gamma$ is a fundamental dominant weight, which contradicts Lemma \ref{1label}.
 
 Now assume $i = 1.$  Here we have no information on
 $\mu^0.$   Let $\gamma$ be the highest weight of an irreducible summand of $V_{C^0}(\mu^0) \downarrow L_X'$ having maximal 
$S$-value.  Also, let $\delta' $ be the restriction of $\delta$ to $S_X, $  so that $\delta' = \l_1^0 \downarrow S_X.$
 This time set $\xi = \bigotimes_{j \ne 0,1}V_{C^j}(\mu^j)$, so that $V^1(Q_Y) = 
\xi \otimes V_{C^0}(\mu^0) \otimes V_{C^1}(\lambda_1^1)$  and
$V_{\kappa}^2(Q_Y) \supseteq \xi \otimes V_{C^0}(\mu^0 + \lambda_{r_0}^0) \otimes V_{C^1}(\lambda_2^1).$
Therefore, with $\epsilon$, $\nu$ and $\psi$ as before, Proposition \ref{sbd} implies that 
$S(\epsilon) + S(\gamma) + S(\nu) + 1 \ge S(\epsilon) + 
S(V_{C^0}(\mu^0 + \lambda_{r_0}^0) \downarrow L_X') + S(2\nu - \psi)$; that is,
\begin{equation}\label{ineqalhorr}
S(\g)+1 \ge S(V_{C^0}(\mu^0 + \lambda_{r_0}^0)) + S(\nu - \psi).
\end{equation}
 Lemma \ref{p.29} shows that $S((V_{C^0}(\mu^0 + \lambda_{r_0}^0)\downarrow L_X') \ge S(\gamma) + S(\delta')$ (as $V_{C^0}(\lambda_{r_0}^0) \downarrow L_X'$ is irreducible).  
 By hypothesis $S(\delta') \ge 1$, and so the inequality (\ref{ineqalhorr}) implies that $S(\nu - \psi) = 0$.
As $\psi$ is a root this forces $\nu = \omega_1 + \omega_l.$ 
By hypothesis $\d \ne r\o_j$, and hence the fact that 
$\nu$ corresponds to the highest weight with the largest monomial at level 1 implies that 
$\delta = \omega_1 + d_{l+1}\omega_{l+1}.$  Also $\d \ne \o_1+\o_{l+1}$ by hypothesis, so $d_{l+1} \ge 2$. 
Note that $2\nu - \psi = \omega_1 + \omega_l$.  Further, $\delta' = \omega_1$, so that $C^0$ and $L_X'$ both have type $A_l$ and $V_{C^0}(\mu^0) \downarrow L_X'$ is irreducible, say with highest weight $\mu = \mu^0 \downarrow S_X.$ 

 Now    $L_X'$ acts on level 1 of $W$  as the sum of two irreducibles with highest weights  $\om_1+\om_l$ and $0$.  Similarly for level 2 of $W$, although the highest weights are   $\omega_1 + 2\omega_l$ and   $\omega_l$.  Theorem \ref{NORDU}(i) implies that $\mu^2$ affords the trivial, natural or dual module for $C^2$, so  Proposition \ref{tensorprodMF} forces $\mu^2 = 0.$

Let $\xi \downarrow L_X' = V_{L_X'}(\epsilon_1) + \cdots + V_{L_X'}(\epsilon_r) +V_{L_X'}(\epsilon_{r+1})+ \cdots$, where $\epsilon =\epsilon_1, \ldots , \epsilon_r$ is the complete set of  highest weights of irreducibles in the restriction with maximal $S$-value.   Note that each occurs with multiplicity $1$ since $V^1$ is MF.

 Let ${\kappa'}$ be the node between $C^1$ and $C^2.$  Then $V^2(Q_Y) \supseteq V_{\kappa}^2(Q_Y) + V_{{\kappa}'}^2(Q_Y).$  As mentioned above,  $V_{\kappa}^2(Q_Y)$ contains $V_{C^0}(\mu^0 + \lambda_l^0) \otimes V_{C^1}(\lambda_2^1) \otimes \xi$, and  as $V_{C^1}(\lambda_2^1) \downarrow L_X' \supseteq 
(V_{L_X'}(2\nu - \psi))^2  = (V_{L_X'}(\omega_1+\omega_l))^2$,   the restriction of $V_{\kappa}^2(Q_Y)$ to $L_X'$ contains $\sum_jV_{L_X'}((\mu + \omega_l ) + (\omega_1 + \omega_l) + \epsilon_j)^2.$ On the other hand 
$V_{{\kappa}'}^2(Q_Y)$ contains 
$V_{C^0}(\mu^0) \otimes V_{C^2}(\lambda_1^2) \otimes \xi$ and the restriction to $L_X'$
contains $\sum_jV_{L_X'}(\mu + (\omega_1 + 2\omega_l) + \epsilon_j).$  Therefore, the irreducible $L_X'$-modules with highest weights $\mu + \omega_1 + 2\omega_l + \epsilon_j $ each occur 
with multiplicity at least $3$ in $V^2$.  

On the other hand,  $V^1 = \sum_j\mu \otimes ((\om_1 + \om_l) +0) \otimes \epsilon_j $,  so by Corollary \ref{cover}, $\sum_{i,n_i=0}V_i^2(Q_X)$ is a submodule of  
\[
\begin{array}{rl}
V^1\otimes V_{L_X'}(\omega_l) =&  
\sum_j\mu \otimes \left((\om_1 + 2\om_l)\oplus  \om_l \oplus \om_l\oplus (\o_1+\o_{l-1})\right) \otimes \epsilon_j  \\
 = &  \left(\sum_j\mu \otimes (\om_1 + 2\om_l)\otimes \epsilon_j\right) + 
\left( \sum_j\mu \otimes \om_l \otimes \epsilon_j\right)^2 + \\
& \left(\sum_j\mu \otimes (\om_1 + \om_{l-1})\otimes \epsilon_j\right) .
\end{array}
\]
Fix $j\le r$. It is clear from $S$-value considerations that the second and third summands have no irreducible of highest weight $\mu + \om_1 + 2\om_l+ \epsilon_j$, while the irreducible with this highest weight appears precisely once in the first summand.  So this contradicts Proposition \ref{induct}(ii).  This completes the analysis when $\mu^i  = \lambda_1^i.$

 If  $\mu^i  = \lambda_{r_i}^i$ we use essentially the same argument, although we work to the right rather than the left and consider  $C^{i+1}$ rather than $C^{i-1}.$ First assume $i \le k-2,$ where $k = \sum_j d_j$ (see Theorem \ref{LEVELS}). We then show $\mu^{i+1} = 0$ using Lemmas \ref{NORDU} and  \ref{tensorprodMF}.  Letting $\kappa$ be the node between $C^i$ and $C^{i+1},$ we get a contradiction by studying $V_{\kappa}^2(Q_Y).$  The final case is $i = k-1.$  Here
 we can replace $W$ and $V$ by their duals and proceed as above.  \hal

  \chapter{Proof of Theorem \ref{MAINTHM}, Part II: $\mu^0$ is not inner}\label{ch15}

We continue with the notation introduced at the beginning of the previous chapter. In particular, $l\ge 2$.  
Theorem \ref{MUIZERO} shows that if we exclude a small number of possibilities for $\d$, then $\mu^i = 0$ for $0<i<k$. The main result of this section restricts the possibilities for the weight $\mu^0$, showing that it is not inner  (in the sense of Definition \ref{innerdef}) under certain additional  hypotheses. 
 
Recall that $\g_1$ denotes the node between $C^0$ and $C^1$. 
Recall also from Chapter \ref{notation} that $ V_{\gamma_1}^j(Q_Y)$ denotes the sum of weight spaces in 
$V^j(Q_Y)$ afforded by weights of the form $\l -\psi -\gamma_1$, where $\psi$ is a sum of positive roots in 
$\Sigma (L_Y')$. Define 
\[
S_1^2 = S(V_{\g_1}^2(Q_Y)\downarrow L_X').
\]

  \begin{theor}\label{munotinner}  Assume  the induction hypothesis holds, that $\la \lambda, \gamma_1 \ra = 0$, and that  $S_1^2 = S(V^2)$.  In addition, assume that either $L(\delta') \ge 2$  or that $\delta = a\omega_1 + b\omega_{l+1}$ with $a\ge 3$ and $a \ge b > 0$.   Then $\mu^0$ is not equal  to one of the following weights:
 \begin{itemize}
\item[{\rm (i)}]  $c\lambda_{r_0}^0$, with $c \le 4$.
 \item[{\rm (ii)}] $\lambda_{r_0 -c}^0$, with $c \le {\rm min}\left(4,\,\frac{1}{2}(r_0-1)\right)$.
 \item[{\rm (iii)}] $\lambda_1^0 + \lambda_{r_0}^0$.
 \item[{\rm (iv)}] $\lambda_{r_0 -1}^0 + \lambda_{r_0}^0$.
\end{itemize}
  \end{theor}
 
The hypothesis $S_1^2 = S(V^2)$ in the theorem just means that the summand of largest $S$-value in $V^2$ is afforded by a weight of the form $\l -\psi -\gamma_1$, as above. 
 The hypothesis on $\d$ and $\d'$  implies that $\delta \ne r\omega_j$ or    $\omega_1 + \omega_{l+1}$.  
Consequently, Theorem \ref{MUIZERO} implies that  $\mu^i = 0$ for $0 < i < k$ and so $V^1(Q_Y) = V_{C^0}(\mu^0) \otimes V_{C^k}(\mu^k).$ 
 
In the proof we adopt further notation as follows. Set  
\[
J^k = V_{C^k}(\mu^k) \downarrow L_X', \hbox{ and } S^k = S(J^k).
\]
 Let $\gamma = \gamma_1$ and  for $j \ge 2$ write 
\[
V_1^j = V_{\gamma_1}^j(Q_Y).
\]
Finally let $\nu$ be the highest weight of the $L_X'$-composition factor of $W$ corresponding to the largest monomial at level 1.
 
 \vspace{2mm} 
 We now work through the cases (i)-(iv) of Theorem \ref{munotinner}.
 
 \begin{lemma}\label{l1inner}
We have $\mu^0 \ne \lambda_{r_0}^0$.
\end{lemma}

\pf Suppose $\mu^0 =  \lambda_{r_0}^0$. 
Then $V^1 = (\d')^* \otimes J^k$ (abbreviating $V_{L_X'}(\d')$ by just $\d'$ as usual).  As $\la\lambda, \gamma_1\ra = 0$, 
$V_1^2$ is the irreducible afforded by $\l-\beta_{r_0}^0 - \g$, which is $V_{C^0}(\lambda_2^0)^* \otimes V_{C^1}(\lambda_1^1) \otimes V_{C^k}(\mu^k)$. Hence 
\begin{equation}\label{v2stuff}
V_1^2\downarrow L_X' = \wedge^2(\d')^* \otimes W^2(Q_X) \otimes J^k.
\end{equation}
We have $S( \wedge^2(\d')^* ) \le 2S(\d')$. Also Corollary \ref{V^2(Q_X)} implies that $\nu$ 
is a weight of highest $S$-value in $W^2(Q_X)$. Therefore 
\begin{equation}\label{seqn1}
S(V^2) = S_1^2 \le 2S(\d')+S(\nu) + S^k.
\end{equation}
Now consider $V_1^3$, where $\l-\beta_{r_0-1}^0 -2\beta_{r_0}^0-2\g-\beta_1^1$ affords $V_{C^0}(\lambda_{r_0-2}^0)\otimes V_{C^1}(\lambda_2^1) \otimes V_{C^k}(\mu^k)$, and hence 
\[
V_1^3\downarrow L_X' \supseteq \wedge^3(\d')^* \otimes \wedge^2(W^2(Q_X)) \otimes J^k.
\]
By Lemma \ref{mult2}, $\wedge^2(W^2(Q_X))$ has a composition factor of multiplicity at least 2 of highest weight $2\nu -\psi$, where $S(2\nu-\psi) \ge 2S(\nu)-1$ or $2S(\nu)-2$, the latter only if $\delta' = a\omega_1$. If $\d'$ has distinct nonzero coefficients $d_i,d_j$, then $\wedge^3(\d')$ contains the wedge of three vectors of weights $\d',\d'-\a_i,\d'-\a_j$, and hence has a composition factor of $S$-value at least $3S(\d')-2$. And if $\delta' = a\omega_1$, we see that there is a composition
factor of $S$-value at least $3S(\d')-3$.
 It follows that 
$V^3$ has a composition factor of multiplicity at least 2 and $S$-value at least $3S(\d')-2+2S(\nu)-1 + S^k,$ respectively $3S(\d')-3+2S(\nu)-2+S^k$. Hence by Proposition \ref{sbd} and (\ref{seqn1}),
\[
3S(\d')+2S(\nu)-c +S^k \le 2S(\d')+S(\nu)+1+S^k,
\]
where $c = 3$ or $ 5$, respectively.
It follows that $S(\d')+S(\nu)\le c+1$. If $d_{l+1}\ne 0$ then $S(\nu) = S(\d')+1$ (see Corollary \ref{V^2(Q_X)}(iii)), so $2S(\d')\le c$, which is impossible. Hence $d_{l+1}=0$ and so by hypothesis    $L(\d')\ge 2.$ Then Corollary \ref{V^2(Q_X)} shows that $S(\nu) = S(\d')$, whence 
\[
S(\nu) = S(\d') = 2.
\]
Thus $\d = \om_i+\om_j$, where $1\le i < j\le l$. If $i>1$ or $j<l$ then $\wedge^3(\d')$ has a composition factor of 
highest weight $3\d'-\a_i-\a_j$ which has $S$-value at least $3S(\d')-1$, which as above yields 
$S(\d')+S(\nu)\le 3$, a contradiction. Hence $i=1, j=l$ and we have
\[
\d = \om_1+\om_l.
\]
If $l = 2$, then a Magma computation shows that  $V^2 \supseteq (31)^3\otimes J^k$
and none of these summands can arise from $V^1$.
 So assume $l > 2.$ Now $W^2(Q_X)\downarrow L_X' = (\om_1+\om_{l-1})\oplus \om_l$ by Corollary \ref{V^2(Q_X)}, and 
$\wedge^2(\d')^* \supseteq (2\omega_1 + \omega_{l-1}) \oplus (\omega_2 + 2\omega_l).$  One then checks
that $\wedge^2(\d')^* \otimes (W^2(Q_X)\downarrow L_X') \supseteq (2\omega_1 + \omega_{l-1}+ \omega_l)^2.$ 
Therefore, $V^2$ has a repeated composition factor of  $S$-value $4 + S^k,$  which  contradicts 
Lemma \ref{weightsum} as the $S$-value of $V^1$ is $2 + S^k.$  \hal

\begin{lemma}\label{l2inner}
We have $\mu^0 \ne 2\lambda_{r_0}^0$.
\end{lemma}

\pf Suppose $\mu^0 = 2\lambda_{r_0}^0$. 
Note that $V^1 = S^2(\d')^* \otimes J^k$, which has $S$-value $2S(\d') + S^k$.

Now, $V_1^2$  is  the irreducible afforded by $\l-\beta_{r_0}^0-\g$, so that  
$V_1^2 = V_{C^0}(\lambda_{r_0}^0+\lambda_{r_0-1}^0) \otimes V_{C^1}(\lambda_1^1) \otimes V_{C^k}(\mu^k).$ 
The restriction of this to $L_X'$ is contained in  $(\d'^* \otimes \wedge^2(\d'^*)) \otimes W^2(Q_X) \otimes J^k$, 
so has $S$-value at most $3S(\d')+S(\nu) + S^k$.

In $V_1^3$ the weight $\l-2\beta_{r_0}-2\gamma$ affords $V_{C^0}(2\lambda_{r_0-1}^0) \otimes V_{C^1}(2\lambda_1^1) \otimes V_{C^k}(\mu^k)$.  
If $U=V_{C^0}(\lambda_1^0),$ then $S^2(\wedge^2(U)) = 2\lambda_2^0 \oplus \lambda_4^0$.  Now $\wedge^2(\delta')$ has a composition factor of highest weight $2\delta' - \alpha_j,$  so that $S^2(\wedge^2(U)) \downarrow L_X'$ has a composition factor of highest weight
$4\delta' - 2\alpha_j$ and this is not a weight of $\wedge^4(\delta').$  Therefore  $V_{C^0}(2\lambda_2^0)\downarrow L_X' \supseteq 4\delta' - 2\alpha_j$ which has $S$-value $2S(\wedge^2(\d'))$.
Taking duals it follows that $V_{C^0}(2\lambda_{r_0-1}^0) \downarrow L_X'$ has a composition factor with this $S$-value.
 
Also  $ V_{C^1}(2\lambda_1^1)$ restricts to $L_X'$ as $S^2(W^2(Q_X))$, which by Lemma \ref{mult2}  has a repeated composition factor of $S$-value at least $2S(\nu)-2$. Consequently $V^3$ has a composition factor of multiplicity at least 2 of $S$-value at least $2S(\wedge^2(\d'))+2S(\nu)-2 + S^k$.
Now Lemma \ref{weightsum} gives
\[
2S(\wedge^2(\d'))+2S(\nu)-2 + S^k  \le 3S(\d')+S(\nu) + S^k+1.
\]
Since $S(\wedge^2(\d')) \ge 2S(\d')-1$, it follows that $S(\d')+S(\nu) \le 5$. 

If $\delta = a\o_1 + b\o_{l+1}$, then $S(\d') = a \ge 3$ and $S(\nu) = a+1$, so this is impossible.
Therefore $L(\delta') \ge2.$ Here $S(\nu) = S(\delta')$ or $S(\d')+1$, so we must have $S(\delta') = 2$ and $\d = \om_i + \om_j.$
If $\d' \ne \om_1+\om_l$, then $S(\wedge^2(\d')) = 2S(\d')$ since the weight of the wedge of two vectors of weights $\d'$ and $\d'-\a_i$ or $\d'-\a_j$ has  $S$-value equal to that of $\d'$; now the above inequality becomes $S(\d')+S(\nu) \le 3$, which is impossible. Hence $\d' = \om_1+\om_l$.

Here we work with $V_1^2 = V_{C^0}(\lambda_{r_0}^0+\lambda_{r_0-1}^0)  \otimes V_{C^1}(\lambda_1^1) \otimes V_{C^k}(\mu^k).$
With $U$ as above, we have $U \otimes \wedge^2U \cong V_{C^0}(\lambda_1^0+\lambda_2^0) \oplus
\wedge^3U$.  Restricting the tensor product  to $L_X',$ we see that there is a composition factor of highest weight
$3\delta' - \alpha_1$  and this is not a weight of $\wedge^3\delta'.$   Therefore $V_{C^0}(\lambda_1^0+\lambda_2^0) \downarrow L_X' \supseteq 3\d'-\a_1.$

First suppose $l \ge 3.$ Taking duals we see that $V_{C^0}(\lambda_{r_0}^0+\lambda_{r_0-1}^0) \downarrow L_X'$ has a composition factor of highest weight $(30 \dots 011).$
Also $W^2(Q_X)$ contains $(\om_1+\om_{l-1})\oplus \om_l$ (see Corollary \ref{V^2(Q_X)}).  Lemma \ref{abtimesce} implies that $(30 \ldots  011) \otimes (10\ldots  010)$ has a composition factor of highest weight $(40 \ldots  021) -(\alpha_1 +\cdots + \alpha_{l-1}) = (30\ldots 012).$ Therefore 
$(30 \ldots  011) \otimes ((10\ldots  010) + (0 \ldots  01)) \supseteq (3 0 \ldots  012)^2$
and $V_1^2 \downarrow L_X'$ has a composition factor of multiplicity at least 2 and $S$-value at least $6 + S^k$. Since $S(V^1) = 2S(\d') + S^k = 4 + S^k$, this contradicts Lemma \ref{weightsum}.  

 Now suppose $l = 2.$  Here we find that $(41)$ is a composition factor  of
$V_{C^0}(\lambda_{r_0}^0+\lambda_{r_0-1}^0) \downarrow L_X'.$  Also $W^2(Q_X)$ contains $(20) + (01)$ so that
 $V_1^2 \downarrow L_X' \supseteq (41) \otimes ((20) + (01)) \otimes J^k  \supseteq (42)^2 \otimes J^k$ and we have a repeated composition factor of $S$-value $6 + S^k.$  As above $S(V^1) =  4 + S^k,$  and so we again contradict Lemma \ref{weightsum}. \hal

\begin{lemma}\label{l3inner}
We have $\mu^0 \ne \lambda_{r_0-1}^0$.
\end{lemma}

\pf Suppose $\mu^0 = \lambda_{r_0-1}^0$.  Here $V^1 = \wedge^2(\d')^* \otimes J^k. $  We also have $V_1^2 \downarrow L_X' = \wedge^3(\d')^* \otimes W^2(Q_X) \otimes J^k$ (afforded by $\l-\beta_{r_0-1}-\beta_{r_0}-\g$).

First assume that $\d = a\om_1 + b\om_{l+1}.$  Then $\d' = a\om_1$
and $S(V^1) = 2a-1 + S^k$.  We have $W^2(Q_X) \supseteq a\om_1 + \om_l$
and $\wedge^3(\d')^* \supseteq 3\om_{l-1} + (3a-6)\om_l.$  Therefore, 
$V_1^2 \downarrow L_X' \supseteq (3\om_{l-1} + (3a-6)\om_l) \otimes (a\om_1 + \om_l) \otimes J^k$
which by Lemma \ref{compfactor}(ii) contains $((a-1)\om_1 + 3\om_{l-1}+ (3a-6)\om_l) \otimes J^k$ with
multiplicity 2.  This has $S$-value $4a-4 + S^k$ so we have a contradiction provided $4a-4 + S^k > 2a + S^k$,
which holds as $a \ge 3.$

So from now on we assume $L(\d') \ge 2.$  In $V_1^3\downarrow L_X'$, the weight $\l- \beta_{r_0-2}^0-2\beta_{r_0-1}^0-2\beta_{r_0}^0-2\g-\beta_1^1$ affords 
$\wedge^4(\d')^* \otimes \wedge^2(W^2(Q_X))\otimes J^k$. Lemma \ref{mult2} shows that in the second tensor factor there is a multiplicity 2 summand with $S$-value at least $2S(\nu)-1$  and checking weights in the fourth wedge we see that the first factor has a summand with $S$-value at least  $4S(\d')-d$, where $d = 3$ unless $l=2$, in which 
case $d = 4$.  

Now $S(\wedge^3(\d')^* \otimes W^2(Q_X) \otimes J^k) \le  3S(\d')-e + S(\nu) + S^k$, where
$e =1$ unless $l=2$, in which case $e=2$.  Therefore we can use Lemma \ref{weightsum} to see that
\begin{equation}\label{ineqal}
4S(\d')-d+2S(\nu)-1+S^k \le 3S(\d') - e +S(\nu) + S^k +1.
\end{equation}
   Thus $S(\d')+S(\nu)\le d - e+2.$  As, $L(\delta') \ge 2$ we have $S(\nu) \ge S(\delta')$ and hence $2S(\d') \le 
   d - e+2 \le 4$.  Therefore $S(\delta') = L(\delta') = 2$  
and $\delta' = \omega_i + \omega_j$ for some $ i < j.$

If  $i>1$ and $j < l$, then we claim $S(\wedge^4 \delta') = 4S(\delta').$   Indeed, the   wedge of $\delta', \delta' - \alpha_i, \delta' - \alpha_j$, and $\delta' - \alpha_i - \alpha_j$ has this value and  is subdominant  to the highest weight of a composition factor whose $S$-value must be at least as large. This improved $S$-value gives a contradiction. Therefore $\d' = \om_1+\om_j$ or $\om_i+\om_l$.

Assume  that $l\ge 4$ and $\d'\ne \om_1+\om_l$. We claim that $\wedge^4(\d')$ has a composition factor of $S$-value at least 7 $(=4S(\d')-1)$; given this, the above argument improves to  $S(\d')+S(\nu)\le 2$, a contradiction. It suffices to find a dominant weight of $S$-value $7$. If  $\d' = \om_1+\om_j$ with $j<l$, the wedge of four vectors of weights $\d'$, $\d'-\a_j$, $\d'-\a_j-\a_{j+1}$, $\d'-\a_{j-1}-\a_j$ is a dominant weight with $S$-value at least 7.  A similar argument 
applies  if  $\d' = \om_i+\om_l.$

It remains to consider the cases where either $\d' = \om_1+\om_l$ or $l= 3$ and $\d' \in \{\om_1+\om_2, \om_2+\om_3\}$.
Consider the first case, $\d' = \om_1+\om_l$. Here $S(V^1) = S(\wedge^2(\d') \otimes J^k)  = 3 + S^k$. Now
\begin{equation}\label{awf}
V_1^2\downarrow L_X' \supseteq  \wedge^3(\d')^* \otimes W^2(Q_X) \otimes J^k.
\end{equation}
If $l \ge 4$, then Proposition \ref{compfactor} implies that 
\[
 \wedge^3(\d')^* \otimes W^2(Q_X) \supseteq (110\ldots  011) \otimes \left((10\ldots  010) + (0 \ldots  01)\right) \supseteq (110 \ldots  012)^2.
\]
If $l = 3$ the tensor product contains $(121) \otimes ((110) + (001))  \supseteq (122)^2.$  And if $l = 2$ the tensor product contains $(22) \otimes ((20) + (01)) \supseteq (23)^2.$ So for each of these 
the  tensor product of the first two factors in (\ref{awf})  has a composition factor of multiplicity at least 2 and $S$-value at least 5, contradicting Lemma \ref{weightsum}. If $l=3$ and $\delta' = \om_1+\om_2$, then $S(\wedge^2(\d') \otimes J^k)=4 + S^k$ and $V_1^2\downarrow L_X' \supseteq (122) \otimes ((200) + (010)) \otimes J^k \supseteq 
(132)^2 \otimes J^k$ and thus has a composition factor of multiplicity at least 2 and $S$-value at least $6 +S^k$, again a contradiction. Similar arguments apply if $l=3$ and $\delta' = \om_2+\om_3.$ \hal

\begin{lemma}\label{l4inner}
We have $\mu^0 \ne \lambda_{r_0-2}^0$ with $r_0 \ge 5.$
\end{lemma}

\pf Suppose $\mu^0 = \lambda_{r_0-2}^0$ with $r_0 \ge 5.$  First assume that $L(\delta') \ge 2.$  Then the induction hypothesis  implies $\d' = \om_1+\om_l$. Observe that $V^1 = \wedge^3(\d')^* \otimes J^k$, which has $S$-value $4 + S^k$. On the other hand, in $V_1^2\downarrow L_X'$, the weight 
$\l-\beta_{r_0-2}^0-\beta_{r_0-1}^0-\beta_{r_0}^0-\g$  affords $\wedge^4(\d')^* \otimes W^2(Q_X)\otimes J^k$, and $W^2(Q_X)$ contains 
$(\om_1+\om_{l-1})\oplus \om_l$. There is a composition factor of $\wedge^4(\delta')$ of highest weight given by the wedge of $\d'$, $\d'-\a_1$, $\d'-\a_l$, $\d'-\a_{l-1}-\a_l$.   Taking the duals we have a composition factor of highest weight $(1010\ldots  012)$, $(1022),$ $(113)$, or $(22)$ according as $l \ge 5, l=4, l=3, $ or $l=2$. Therefore,  if $l \ge 5$,  Lemma \ref{compfactor} implies that $V_1^2\downarrow L_X' \supseteq (1010\ldots  012) \otimes ((10\ldots 010) + (0\ldots  01)) \otimes J^k \supseteq (1010\ldots  013)^2 \otimes J^k.$ The repeated composition factor has $S$-value $6 + S^k$, which contradicts Lemma \ref{weightsum}.   For $l = 4,3,2$ we get repeated composition factors
$(1023)^2, (114)^2$, or $(23)^2,$ respectively.  In the first two cases we  again have a contradiction.  

Suppose $l = 2$, so that the $S$-value of the above repeated factor is only $5 + S^k.$  Here $\delta = 11x.$  If $x > 0$, then $W^2(Q_X) \supseteq (12 + 20)$ and   $(42)^2 \otimes J^k$  appears, hence there exists a repeated composition factor of $S$-value $6 + S^k$ and we obtain the same contradiction.  So finally assume $x = 0$ so that $\d = 110$, 
$C^k = A_2$ and $J^k = (st)$ is an irreducible module. Now $\wedge^4(\delta') \supseteq (22)^2$, so that the repeated factor  $(23) \otimes  J^k$ in $V_1^2\downarrow L_X'$ occurs with multiplicity $4.$  On the other hand, 
$V^1 = ((22) + (30)+(03)+(11)+(00)) \otimes (st)$ which has $S$-value $4+s+t.$  Moreover, 
$V^1$ can contribute at most 1 composition factor of $S$-value $5 + s + t$ to $V^2$ and so we have a contradiction.

Now suppose $\delta = a\omega_1 + b\omega_{l+1}$ with $a \ge 3.$  Here the $S$-value of  
$V^1$ is $3a-2 + S^k$, or $3a-3 +S^k$ if $l = 2$.  On the other hand, we argue as above that $\wedge^4(\d')^*$ contains an irreducible summand of highest weight $(0\ldots 0100(4a-4)) $, $(00(4a-4))$, or $(2(4a-7))$, according as $l \ge 4, l =3,$ or $l=2.$ 
If $l \ge 4$ then Lemma \ref{compfactor}(ii)  implies that 
\[
\begin{array}{ll}
V_1^2\downarrow L_X' & \supseteq (0\ldots 0100(4a-4)) \otimes ((a0\ldots 01) + ((a-1)0 \ldots  0)) \otimes J^k \\
                                       & \supseteq ((a-1)0\ldots  0100(4a-4))^2 \otimes J^k.
\end{array}
\]
  Hence there is a repeated composition factor with 
$S$-value $5a-4+S^k.$ But then Lemma \ref{weightsum} implies that $5a-4+S^k \le (3a-2 + S^k)+1, $ which is not the case.  If $l = 3$ or $l = 2$ we get repeated composition factors 
$((a-1)0(4a-4))^2$ or $((a+1)(4a-7))^2$, respectively, and once again this yields a contradiction.
\hal

\begin{lemma}\label{l5inner}
We have $\mu^0 \ne 3\lambda_{r_0}^0 $ or $4\lambda_{r_0}^0. $
\end{lemma}

\pf Suppose $\mu^0 = 3\lambda_{r_0}^0$. The induction hypothesis implies that either $l=2$ and $\d' = \om_1+\om_2$, or $\delta = a\omega_1 + b\omega_{l+1}.$ Then $V^1 = S^3(\d')^* \otimes J^k$, which has $S$-value $6+S^k$ or $3a+S^k$, respectively. In $V_1^2(Q_Y)$
the weight $\lambda-\beta_{r_0}-\g$ affords $V_{C^0}(\lambda_{r_0-1}^0+2\lambda_{r_0}^0) \otimes V_{C^1}(\lambda_1^1) \otimes V_{C^k}(\mu^k)$. If $\d' = \om_1+\om_2$  we compute that the restriction of the first tensor factor to $L_X' = A_2$ contains $(33)^2$, while the restriction to the second factor is $W^2(Q_X)$ which contains $20\oplus 01$. Hence $V_1^2(Q_Y)\downarrow L_X'$ contains $53^2 \otimes J^k$, of $S$-value $8 + S^k$. This contradicts Lemma \ref{weightsum}. 

Now assume that $\delta = a\omega_1 + b\omega_{l+1}.$  
Now $\l_{r_0}^0\otimes 3\l_{r_0}^0 = 4\l_{r_0}^0\oplus (\l_{r_0-1}+2\l_{r_0}^0)$. Moreover, $\d'^*\otimes S^3(\d'^*)$ contains $(0\ldots 01\,4a-2)$, while 
$S^4(\d'^*)$ does not. Therefore,
$V_1^2 \downarrow L_X' \supseteq (0\ldots  0 1(4a-2)) \otimes ((a0 \ldots  01) + ((a-1)0 \ldots  0)) \otimes J^k$ and Lemma \ref{compfactor}(ii) shows that the tensor product of the first two terms contains $((a-1)0\ldots  01(4a-2))^2$ (or $(a(4a-2))^2$ if $l = 2.$) Therefore we have a repeated composition factor in $V_1^2 \downarrow L_X'$ of $S$-value $5a-2 + S^k$ and this contradicts Lemma \ref{weightsum}. 

Now assume $\mu^0 = 4\lambda_{r_0}^0. $ The induction hypothesis implies that $\delta = a\omega_1 + b\omega_{l+1}$ with $a =3.$  Here  $V^1 = S^4(\d')^* \otimes J^k$, which has $S$-value $4a + S^k.$ In $V_1^2$ the weight $\l-\beta_{r_0}$ affords $V_{C^0}(\lambda_{r_0-1}^0+3\lambda_{r_0}^0) \otimes V_{C^1}(\lambda_1^1) \otimes V_{C^k}(\mu^k)$. Arguing as in the last paragraph we see that the restriction of the first term to $L_X'$ contains 
$(0 \ldots  01(5a-2))$
and Lemma \ref{compfactor}(ii) implies that there is a repeated composition factor of highest weight $((a-1)0\ldots  01(5a-2)) \otimes J^k$ (or $(a(5a-2)) \otimes J^k$ if $l = 2$) and  $S$-value 
$6a-2 + S^k.$  This contradicts Lemma \ref{weightsum}.  \hal

\begin{lemma}\label{l7inner}
We have $\mu^0 \ne \lambda_{r_0-3}^0 \,(r_0 \ge 7) $ or $\lambda_{r_0-4}^0\, (r_0 \ge 9)$.
\end{lemma}

\pf  Assume false.  Then the induction hypothesis implies that $\delta = a\omega_1 + b\omega_{l+1}$ with
$a\le 4$ or $a = 3$, respectively. 
We have  $V^1 = \wedge^4(\delta')^* \otimes J^k$ or $\wedge^5(\delta')^* \otimes J^k,$
respectively. 

In order to avoid special cases we  first assume that $l >  4$ or $l > 5$, respectively.  Then easy checks show that  $S(V^1) = 4a-3 + S^k$ or $5a-4 + S^k$, respectively.
Now $V_1^2$ has a composition factor afforded by $\lambda - \beta_{r_0-3}^0 - \cdots -\beta_{r_0}^0 -\gamma $ (respectively   $\lambda - \beta_{r_0-4}^0 - \cdots -\beta_{r_0}^0 -\gamma ).$   This affords  $V_{C^0}(\lambda_{r_0-4}^0) \otimes V_{C^1}(\lambda_1^1) \otimes V_{C^k}(\mu^k)$ (respectively $V_{C^0}(\lambda_{r_0-5}^0) \otimes V_{C^1}(\lambda_1^1) \otimes V_{C^k}(\mu^k)$) which restricts to $L_X'$ as 
$\wedge^5(\delta')^* \otimes ((a0 \ldots  01) + ((a-1)0 \ldots  0)) \otimes J^k$ (respectively $\wedge^6(\delta')^* \otimes ((a0 \ldots  01) + ((a-1)0 \ldots  0)) \otimes J^k$) and these
have composition factors $(0 \ldots  01000(5a-5)) \otimes  ((a0 \ldots  01) + ((a-1)0 \ldots  0))\otimes J^k$ (respectively $(0 \ldots  010000(6a-6))  \otimes  ((a0 \ldots  01) + ((a-1)0 \ldots  0))  \otimes J^k$.) 
Therefore Lemma \ref{compfactor}(ii) implies that there is   a repeated composition factor $((a-1) 0 \dots 01000(5a-5)) \otimes J^k$ (respectively $((a-1) 0 \dots 010000(6a-6)) \otimes J^k$) and once again we contradict Lemma \ref{weightsum}. 

This leaves the cases where $l \le 4$ (respectively $l \le 5$) which were excluded earlier. But since $a \le 4$ these can  be handled using Magma.  For example assume $a = l = 3.$   If $V^1 = \wedge^4(\delta')^* \otimes J^k$ then a Magma computation shows   that $S(\wedge^4(\delta')^* \otimes J^k) = 8 + S^k$ and that $\wedge^5(\delta')^* \supseteq  (027).$  Therefore $\wedge^5(\delta')^* \otimes  ((301) + (200))  \supseteq (027) \otimes  ((301) + (200)) \supseteq (227)^2.$ So there is a repeated composition factor of $S$-value $11 + S^k$ which is  a contradiction.  And if $V^1 = \wedge^5(\delta')^* \otimes J^k$  then  $S(V^1) = 9 + S^k$  while
$\wedge^6(\delta')^* \otimes  ((301) + (200))  \supseteq (127) \otimes  ((301) + (200)) \supseteq (327)^2$
again a contradiction.

 The remaining cases are left to the reader.  \hal

\begin{lemma}\label{l8inner}
We have $\mu^0  \ne \lambda_1^0 + \lambda_{r_0}^0$.
\end{lemma}

\pf Assume $\mu^0 = \lambda_1^0 + \lambda_{r_0}^0,$ so that $\delta = a\omega_1 + b\omega_{l+1}$
and $V^1 = ((a0 \ldots  0) \otimes (0 \ldots  0a) -0) \otimes J^k.$  Here $\lambda - \beta_{r_0}^0 -\gamma $ affords the  irreducible  $V_{C^0}(\lambda_1^0 + \lambda_{r_0 -1}^0) \otimes V_{C^1}(\lambda_1^1) \otimes V_{C^k}(\mu^k).$   The first factor is  $V_{C^0}(\lambda_1^0) \otimes V_{C^0}(\lambda_{r_0 -1}^0) - V_{C^0}(\lambda_{r_0 }^0)$ so
the restriction to $L_X'$ contains an irreducible of highest weight $(a0 \ldots  0 1(2a-2))$ if $l \ge 3$ and $((a+1)(2a-2))$ if $l = 2.$  So if $l \ge 3$, then Lemma \ref{compfactor}(i) implies that
$V_1^2 \downarrow L_X' \supseteq (a0 \ldots  01(2a-2)) \otimes ((a0 \ldots  01) + ((a-1)0 \ldots  0)) \otimes J^k \supseteq 
((2a-1)0 \ldots  01(2a-2))^2 \otimes J^k.$ Hence there  is  a repeated composition factor with $S$-value
$4a-2 + S^k$ and this contradicts Lemma \ref{weightsum}. And if $l = 2$ we have a repeated composition factor 
$(2a(2a-2)) + J^k$ and $S$-value $4a-2 +S^k,$ again a contradiction. \hal

\begin{lemma}\label{l9inner}
We have $\mu^0  \ne \lambda_{r_0-1}^0 + \lambda_{r_0}^0$.
\end{lemma}

\pf  Assume $\mu^0 = \lambda_{r_0-1}^0 + \lambda_{r_0}^0.$  Then by hypothesis and using the induction hypothesis we have $\delta = 3\omega_1 + b\omega_{l+1}.$  Therefore $V^1$ is contained in  $(\delta')^* \otimes \wedge^2(\delta')^* \otimes J^k$ which has $S$-value $3 +(2\cdot3 -1) + S^k = 8 + S^k.$  Next note that $\lambda - \beta_{r_0}^0-\g$ affords the irreducible $V_{C^0}(2\lambda_{r_0 -1}^0)  \otimes V_{C^1}(\lambda_1^1) \otimes V_{C^k}(\mu^k).$  The first tensor factor is contained in 
$V_{C^0}(\lambda_{r_0 -1}^0) \otimes V_{C^0}(\lambda_{r_0 -1}^0).$  
A maximal vector of $V_{C^0}(2\l_{r_0-1}^0)$ restricts to $L_X'$ as $(0\ldots  028)$. 
Therefore Lemma \ref{compfactor}(ii) implies that $V_1^2 \downarrow L_X' \supseteq
  (0 \ldots  028) \otimes ((30 \ldots  01) + (20 \ldots  0)) \otimes J^k \supseteq (20\ldots 028)^2 \otimes J^k$ (or $(48)^2$ if $l = 2$) with $S$-value $12 + S^k.$  This contradicts Lemma \ref{weightsum}.  \hal
  
\vspace{4mm}
  Theorem \ref{munotinner} follows from the Lemmas \ref{l1inner} -- \ref{l9inner}.

\chapter{Proof of Theorem \ref{MAINTHM}, Part III: $\langle \lambda, \gamma \rangle = 0$}\label{ch16}

Continue with the notation of the previous two chapters. In particular, recall the following notation:
\[
\begin{array}{l}
\d = \sum_1^{l+1}d_i\o_i, \;\d' = \sum_1^{l}d_i\o_i, \;\d'' = \sum_1^{l}d_{i+1}\o_i, \\
\g_i = \hbox{ node between } C^{i-1} \hbox{ and }C^i,\\
V_i^j = V_{\g_i}^j(Q_Y),\,S_i^j = S(V_i^j \downarrow L_X'), \\
J^i = V_{C^i}(\mu^i)\downarrow L_X',\, S^i = S(J^i).
\end{array}
\]
In this chapter we show that under certain hypotheses $\la \lambda, \gamma_i \ra = 0$ for all $i$.
The main result is as follows.

\begin{theor}\label{L(delta')ge2}  Assume the induction hypothesis.  Then $\la \lambda, \gamma_i \ra = 0$ for $1\le i\le k$, provided one of the following holds:
\begin{itemize}
\item[{\rm (i)}] $L(\delta') \ge 2$;  
\item[{\rm (ii)}] $\delta = a\omega_1 + b\omega_{l+1}$ with $a\ge 3$ and $a \ge b > 0.$
\end{itemize}
\end{theor}

Assume throughout this chapter that (i) or (ii) of  Theorem \ref{L(delta')ge2} holds. We aim to show that 
$\la \lambda, \gamma_i \ra = 0$ for all $i$.

Choose $i$ such that $S_i^2$ is maximal --  that is, $S_i^2 = S(V^2)$.

We will proceed in a series of lemmas.  The first just records  information that follows from the inductive hypothesis. 
 
\begin{lemma}\label{indhyp}  {\rm (i)}  If $L(\delta') \ge 2$, then $\mu^0$ or $(\mu^*)^0$ is in
$\{ 0,\, \lambda_1^0,\, \lambda_2^0,\, \lambda_3^0,\, 2\lambda_1^0,\, 3\lambda_1^0 \}.$

{\rm (ii)} If $\delta = a\omega_1 + b\omega_{l+1}$ with $a\ge 3$ and $a \ge b > 0$, then $\mu^0 $ or $(\mu^*)^0$ is one of the following:
\[
\begin{array}{l}
 0, \\
 c\lambda_1^0\, (c \le 4),\\
 \lambda_c^0\, (c \le 5),\\
 \lambda_1^0 + \lambda_{r_0}^0,\\
 \lambda_1^0 + \lambda_2^0 \,(a = 3).
\end{array}
\]
\end{lemma}

We remark that Theorem \ref{munotinner} shows that with  additional hypotheses,  the
above lemma  can be improved so as to delete the terms $(\mu^*)^0$ in both (i) and (ii).  We will have that situation later.

Recall that $\nu^j$ denotes the highest weight of an $L_X'$-composition factor in the $j$th level of $W$ arising from the largest monomial.

\begin{lemma}\label{S(nui)} Assume that  $1 \le j \le k-1.$  Then
\begin{itemize}
\item[{\rm (i)}] $\mu^j = 0$, and 
\item[{\rm (ii)}] $S(\nu^j) \ge 2$. 
\end{itemize}
\end{lemma}

\pf (i) This is Theorem \ref{MUIZERO}.  

(ii) Note first that $\nu^j \ne 0$ (as observed just before Lemma \ref{1label}). If $L(\nu^j) \ge 2$, then the assertion of (ii) is obvious.  Otherwise,  $L(\nu^j) = 1$ and the conclusion follows from Lemma \ref{1label}.  \hal

\begin{lemma}\label{gammai20} We have $V_{\gamma_j}^2(Q_Y) = 0$ unless one of the following holds:
\begin{itemize} 
\item[{\rm (i)}]$\la \lambda, \gamma_j \ra \ne 0$;
\item[{\rm (ii)}] $j = 1$ and $\mu^0 \ne 0$; 
\item[{\rm (iii)}]$j = k$ and $\mu^k \ne 0$.
\end{itemize}
\end{lemma}

\pf This follows  from Lemma \ref{S(nui)}(i).   \hal

\begin{lemma}\label{i=1,k}  {\rm (i)} We have $i = 1$ or  $i =  k.$  

{\rm (ii)}  Assume $\la\lambda, \gamma_j\ra \ne 0$ for some $j$ with $1<j<k$. Then the following hold.
\begin{itemize}
\item[{\rm (a)}]  If $L(\delta') \ge 2,$  then $S_j^2 \ge S^0 + S^k + 4.$
\item[{\rm (b)}]  If  $\delta = a\omega_1 + b\omega_{l+1}$, then $S_j^2 \ge S^0 + S^k + 5.$
\end{itemize}
\end{lemma}

\pf  (i) Suppose $i\ne 1,k$.  As $\l \ne 0$, $S(V^2) \ne 0$, so by Lemma \ref{gammai20} we have $\la \l, \g_i \ra \ne 0.$ Then $V_i^2 = V_{C^0}(\mu^0) \otimes V_{C^{i-1}}(\lambda_{r_{i-1}}^{i-1}) \otimes V_{C^i}(\lambda_1^i) \otimes V_{C^k}(\mu^k),$ which is afforded by $\lambda - \gamma_i$. 

The irreducible afforded by the largest monomial has maximal $S$-value, and the $S$-value of an irreducible module and its dual are equal.  Therefore $S(V_i^2 \downarrow L_X')  = S_i^2 = S^0 + S(\nu^{i-1}) + S(\nu^i) + S^k.$  

Now consider $V_i^3$  (the summand of $V^3(Q_Y)$ involving just $-2\gamma_i$) which contains an irreducible summand
$ V_{C^0}(\mu^0) \otimes V_{C^{i-1}}(\lambda_{r_{i-1} -1}^{i-1}) \otimes V_{C^i}(\lambda_2^i) \otimes V_{C^k}(\mu^k),$ afforded by $\lambda - \beta_{r_{i-1}}^{i-1}  -2\gamma_i - \beta_1^i$.    
 Lemma \ref{mult2} implies that the restriction to $L_X'$ contains 
a summand $(V_{C^0}(\mu^0)\downarrow L_X') \otimes (V_{L_X'}(2\nu^{i-1}-\psi^{i-1})^*)^2 \otimes (V_{L_X'}(2\nu^i - \psi^i))^2 \otimes (V_{C^k}(\mu^k)\downarrow L_X')$, where $\psi^i,\psi^{i-1}$ are as in Lemma \ref{mult2}.
  So we get a composition factor of
multiplicity $4$.  Taking $S$-values and using Lemma \ref{mult2}, we find that $S(V_i^3 \downarrow L_X') \ge S^0 + S(2\nu^{i-1}-\psi^{i-1}) + S(2\nu^i - \psi^i) + S^k \ge S^0 + 2S(\nu^{i-1}) - c + 2S(\nu^i) - c + S^k$, 
where $c =1$ or $2$, the latter only if $\delta = a\omega_1 + b\omega_{l+1}.$
Therefore Lemma \ref{nextlevelSvalue} implies that   $S(\nu^{i-1}) + S(\nu^i) \le 2c+1$. 
If $L(\delta') \ge 2,$ then $c = 1$ and this  contradicts Lemma \ref{S(nui)}.  

For the case $\delta = a\omega_1 + b\omega_{l+1}$ we note that the weights $\nu^s$ for $s = 0, \dots, k$ are as follows: \[
a\omega_1,\, a\omega_1 + \omega_l, \dots, a\omega_1 + b\omega_l,\, (a-1)\omega_1 + b\omega_l, \dots, \omega_1 + b\omega_l,\, b\omega_l.
\]
  So we have $S(\nu^{i-1}) + S(\nu^i) \ge 3 + 2b$.  Notice that this argument gives (ii)(b).  If $b \ge 2$  this is larger than $2c+1$ and we have a contradiction.  This will also hold for $b = 1$ unless  $i = k-1.$ 

In this last case we will obtain a contradiction in  $V^2(Q_Y).$  Indeed, as noted earlier, 
$V_i^2 = V_{C^0}(\mu^0) \otimes V_{C^{i-1}}(\lambda_{r_{i-1}}^{i-1}) \otimes V_{C^i}(\lambda_1^i) \otimes V_{C^k}(\mu^k).$  Restricting the middle two tensor factors to $L_X'$ we have $  ((20 \ldots  01) + (10 \ldots  0)) \otimes ((10 \ldots  01) + (0 \ldots  0))$   and this contains $(20 \ldots  01)^2$.  Therefore, $V_i^2 \downarrow L_X'$ contains a repeated composition factor with $S$-value $S^0 + S^k + 3$, a contradiction.

(ii) This follows from the above proof.  As noted above, (ii)(b) holds. For (ii)(a), we see as in the second paragraph that $S(V_j^2 \downarrow L_X') = S^0 + S(\nu^{j-1}) + S(\nu^j) + S^k.$  Lemma \ref{S(nui)} shows that the two middle terms are each at least 2, which gives the result.  \hal

\begin{lemma}\label{V2} The following hold.
\begin{itemize}
\item[{\rm (i)}] $V_1^2$ is isomorphic to a summand of $ V_{C^0}(\mu^0) \otimes V_{C^0}(\lambda_{r_0}^0) \otimes V_{C^1}(\lambda_1^1) \otimes V_{C^k}(\mu^k).$ Therefore, $S_1^2 \le S^0 + S(\delta') + S(\nu^1) + S^k.$
\item[{\rm (ii)}] $V_k^2 $ is isomorphic to a summand of $ V_{C^0}(\mu^0) \otimes  V_{C^{k-1}}(\lambda_{r_{k-1}}^{k-1}) \otimes V_{C^k}(\mu^k) \otimes V_{C^k}(\lambda_1^k)$.  Therefore, $S_k^2\le S^0 +  S(\nu^{k-1}) + S^k + S(\delta'').$
\item[{\rm (iii)}]  The $S$-value inequalities in (i) (resp. (ii)) are equalities if  $\la \lambda, \gamma_1 \ra \ne 0$ 
(resp. $ \la \lambda, \gamma_k \ra \ne 0$). 
 \item[{\rm (iv)}] If $\la \lambda, \gamma_1 \ra \ne 0,$ then $S_1^2 \ge S^0 + S^k + 4.$
\end{itemize}
\end{lemma}

\pf  We first set up some temporary notation. For each $j$ let   $(Q_Y)_{\gamma_j} = Q_jQ_Y'/Q_Y',$ where $Q_j$ is the product of all root subgroups for negative roots which involve $\gamma_j$. We can regard $(Q_Y)_{\gamma_j}$ as the direct product of root groups for roots of the form $-\gamma_j - \eta,$ where $\eta \in \Z\Sigma^+ (L_Y').$  Then $Q_Y/Q_Y' $ is the direct product of the quotients $(Q_Y)_{\gamma_j}$ and each of these is invariant under $L_Y'.$  

It follows from  \cite[(2.3)(i)]{mem1} that $[V,Q_Y^2]$ is a sum of weight spaces of level at least $2$, so a
 consideration of weights shows that $[V, Q_Y'] \le [V,Q_Y^2].$  Therefore we can abuse notation and
think of $V^2(Q_Y)$ as $ [V^1(Q_Y),Q_Y/Q_Y'] $  and $V_j^2 =  [V^1(Q_Y), (Q_Y)_{\gamma_j}] .$
By doing this we  avoid continually writing  quotients in the arguments to follow.

(i), (ii)  It will suffice to prove (i), so we now take $j = 1$ in the above. 
 Then $V_1^2 =  [V^1(Q_Y), (Q_Y)_{\gamma_1}] $  and we claim that there is a surjective map
 from $V^1(Q_Y) \otimes (Q_Y)_{\gamma_1}$ to $V_1^2$, commuting with the action of $L_Y$.  
 Let $x \in V$ and $a \in Q_1$ and consider the map  $\bar x \otimes \bar a \rightarrow [\bar x,\bar a].$
 This is well defined since $V$ and $Q_Y$  both act trivially on 
 $[V^1(Q_Y), (Q_Y)_{\gamma_j}]$. Moreover  the trivial action together with  the commutator identities $[xy,a] = [x,a]^y[y,a]$ and $[x,ab] = [x,b][x,a]^b$ imply that the map is linear in both coordinates.  This establishes the claim.
 As  $V^1(Q_Y) = V_{C^0}(\mu^0) \otimes V_{C^k} (\mu^k)$ and  $(Q_Y)_{\gamma_1}$  affords $V_{C^0}(\lambda_{r_0}^0) \otimes V_{C^1}(\lambda_1^1), $ this establishes the first assertion in (i).  The second assertion in (i) follows by taking $S$-values, noting that 
$V_{C^0}(\lambda_{r_0}^0)\downarrow L_X' = V_{L_X'}(\delta')^*.$

(iii)  Suppose $\la \lambda, \gamma_1 \ra \ne 0.$   Then Lemma \ref{v2gamma1}
shows $V_{\g_1}^2(Q_Y) \supseteq V_{C^0}(\mu^0) \otimes V_{C^0}(\l_{r_0}^0)  \otimes V_{C^1}(\l_1^1) \otimes V_{C^k}(\mu^k)$.  The result follows.

(iv) This follows from (i) and  (iii) since $S(\delta')  \ge 2$ by hypothesis, and $S(\nu^1) \ge 2$  by 
Lemma \ref{S(nui)}.  \hal

\begin{lemma}\label{notinner} Assume $i = 1$ and $\la \lambda, \gamma_1\ra = 0.$  
\begin{itemize}
\item[{\rm (i)}] $\mu^0$ is not inner.  
\item[{\rm (ii)}]  If $\mu^0 \ne \lambda_4^0, \lambda_5^0$ with $l \le 3,4$, respectively, 
then  $S_1^2 \le S^0 + S^k + 2.$
\item[{\rm (iii)}] If $\mu^0 = \lambda_4^0, \lambda_5^0$ with $l \le 3,4$, respectively, then 
$\delta = a\omega_1 + b\omega_{l+1}$ and $S_1^2 \le S^0 + S^k + 3,$
 unless $l = 2$ and $\mu^0 = \lambda_5^0$, in which case  $S_1^2 \le S^0 + S^k + 4.$
\end{itemize}
\end{lemma}

\pf  Assume the hypotheses of the lemma.  Then the hypotheses of Theorem \ref{munotinner} are satisfied and
part (i) follows. 

 (ii) Consider $S_1^2$. We will work through the possibilities in Lemma \ref{indhyp}.  By (i) we have $\mu^0 = 0, c\lambda_1^0, \lambda_c^0,$ or $\lambda_1^0 + \lambda_2^0.$ If $\mu^0 = 0$, then $V_1^2 = 0$ and the assertion holds.  Now suppose
$\mu^0 \ne 0.$ Then $V_1^2 = M^0 \otimes V_{C^1}(\lambda_1^1) \otimes V_{C^k}(\mu^k)$ where 
$M^0 = (c-1)\lambda_1^0$, $\lambda_{c-1}^0$, or $2\lambda_1^0 \oplus \l_2^0$, respectively. 

 We have $S(V_{C^1}(\lambda_1^1) \downarrow L_X') = S(\nu^1) \le S(\delta') + 1.$ If $\mu^0 = c\lambda_1^0$,  then $S^0 = cS(\delta')$, whereas $S(M^0 \downarrow L_X') = (c-1)S(\delta')$.  So 
$S_1^2 \le (c-1)S(\delta') +  (S(\delta')  +1) + S^k = S^0 +  S^k +1$, as required for (ii).
 
If $\mu^0 = \lambda_c^0$ ($c=2,3$), then $V_{C^0}(\mu^0) \downarrow L_X' =  \wedge^c(\d')$, and 
$M^0 \downarrow L_X' =  \wedge^{c-1}(\d')$. A consideration of weight vectors in $\wedge^c(\d'), \wedge^{c-1}(\d')$
 shows that $S(M^0 \downarrow L_X') \le S^0 - S(\delta') + 1$.  Therefore $S_1^2 \le ( S^0 - S(\delta') + 1) + (S(\delta')+1) + S^k = S^0 + S^k + 2$, as required.

 Next assume $\mu^0 = \lambda_1^0 + \lambda_2^0$, which only occurs for $\delta' = 3\omega_1.$  Here
$M^0 = 2\lambda_1^0 \oplus \l_2^0$ so that $S(M^0 \downarrow L_X') = 6.$  Now $V_{C^0}(\lambda_1^0 + \lambda_2^0) = V_{C^0}(\lambda_1^0) \otimes V_{C^0}(\lambda_2^0) - V_{C^0}(\lambda_3^0).$  It follows
  that $V_{C^0}(\mu^0) \downarrow L_X'$ has an irreducible of highest weight $3\delta'- \alpha_1,$ so that
  $S^0 = 8.$  Therefore,  $S_1^2 = 6 + (3+1) + S^k  = S^0 + S^k + 2.$
   
Next suppose $\mu^0 =  \lambda_c^0,$  for $c = 4, 5$, which only occurs for   $\delta = a\omega_1 + b\omega_{l+1}$ with $3\le a \le 4$ or $a=3$, respectively. Then $\delta' = a\omega_1$ and $S(\nu^1) = a+1.$  Here $V_{C^0}(\mu^0)$ affords $\wedge^c(a\omega_1)$ for $L_X'$.  Assume that  $l \ge c$.  
 Then there is a composition factor of highest weight $(ca-c)\o_1+\o_c$, 
 which has $S$-value $ca-c+1.$   Therefore, $S(M^0 \downarrow L_X') = (c-1)a - (c-1) +1$ and $S_1^2 \le ((c-1)a - (c-1) +1) + (a+1) + S^k = ca -c +3 + S^k = S^0 + S^k + 2.$
 This completes the argument for (ii).
 
 (iii)  Here we again suppose $\mu^0 =  \lambda_c^0,$  for $c = 4, 5,$ and consider the special
 cases $l = 3,4$, respectively.  In these cases the weights of the last paragraph exist but
 yield  $S$-values reduced by $1$.  Therefore the resulting inequality becomes $S_1^2 \le S^0 + S^k + 3.$ 
 
It remains to deal with the cases where $c = 4$ with $l = 2$, or $c = 5$ with $l =2,3.$  For $l = 2$ we find that $\wedge^3(a0)$, $\wedge^4(a0)$, $\wedge^5(a0)$ have composition factors  $((3a-3)0), ((4a-7)2), ((5a-10)2),$ respectively, and these have maximal $S$-values. Arguing as above we see that (iii) holds. Finally, suppose $l = 3$ with
 $c = 5.$  Here  $\wedge^4(a00)$ and  $\wedge^5(a00)$ have composition factors  $((4a-4)00)$ and $ ((5a-8)20),$ respectively, and these have maximal $S$-value.   Again we get (iii).   \hal

\begin{lemma}\label{gamma1} {\rm (i)} Assume $i = 1$. Then  $\la \lambda, \gamma_j\ra = 0$  for all $j \ne k.$

{\rm (ii)}  Assume $i = k$ and $\la \lambda, \gamma_k\ra \ne 0$. Then  $\delta'' = \omega_s$,   $\delta = a\omega_1 + \omega_{s+1}$ for some $s$, and also  $S_k^2 \le S^0 + S^k + 3,$  and  $\la \lambda, \gamma_j \ra = 0$ for $j < k.$
\end{lemma}

\pf (i) In view of Lemmas \ref{notinner} and \ref{i=1,k}(ii), it will suffice to show $\la \lambda, \gamma_1\ra = 0.$  By way of contradiction suppose $\la \lambda, \gamma_1\ra = c > 0.$  We claim that $V_1^3$ contains an
irreducible summand for $L_Y'$ with highest weight $\rho = \lambda - \beta_{r_0}^0 - 2\gamma_1 - \beta_1^1.$
To simplify notation set  $\alpha = \beta_{r_0}^0, \gamma = \gamma_1, \beta = \beta_1^1$, respectively.
 Let $a = \la \lambda, \alpha\ra$.  We work through the possibilities, noting that  we can use Magma for small rank groups to check certain weight space dimensions.  If  $a = 0$ and $c = 1$ then  $\lambda - \alpha - 2\gamma -\beta$ affords a highest weight of $V_1^3$ and the assertion is immediate.
 If $a = 0$ and $c > 1,$ then $\lambda -2\gamma $ affords an irreducible summand for $L_Y'$, but  $\rho$
 only occurs once in the corresponding irreducible, whereas $\rho$ has multiplicity two in $V$.
 Now suppose $a> 0.$  If $c = 1$, then the highest weight in $V_1^3$ is  $\lambda - \alpha - 2\gamma$ and the next is $\lambda - \alpha - 2\gamma - \beta$ which is conjugate to $\lambda -\alpha - \gamma - \beta,$ which has multiplicity 2.  So again there must be a summand of $V_1^3$ with highest weight $\lambda - \alpha - 2\gamma - \beta.$  Finally, assume $a > 0$ and $c > 1.$  There is a summand of highest weight $\lambda - 2\gamma.$ As
 $\lambda - \alpha - 2\gamma$ has multiplicity $2$ there is also a summand of this highest weight.  The next
 one is $\rho$ which has multiplicity $3$.  Again we have the claim.
  
From the claim we see that there is composition factor of $V_1^3$ for which the action of $L_Y'$
is $V_{C^0}(\mu^0 + \lambda_{r_0-1}^0) \otimes V_{C^1}(\lambda_2^1) \otimes V_{C^k}(\mu^k).$ Consider
the first factor. The proof of  Lemma \ref{p.29} shows that $S(V_{C^0}(\mu^0 + \lambda_{r_0-1}^0)\downarrow L_X') \ge 
S(V_{C^0}(\mu^0)\downarrow L_X') + S(\sigma)$ where $\sigma = \l_{r_0-1}^0\downarrow S_X$. An appropriate choice of Borel subgroups $B_X<B_Y$ gives $\sigma = 2(\d')^*-\a_{l-j+1}$ for any $j\le l$ satisfying $d_j\ne 0$.
 Using this and an application of Lemma \ref{mult2} we find that the restriction of 
 $V_{C^0}(\mu^0 + \lambda_{r_0-1}^0) \otimes V_{C^1}(\lambda_2^1) \otimes V_{C^k}(\mu^k)$ to $L_X'$ contains an irreducible summand with multiplicity at least $2$ and $S$-value  at least $ S^0 + 2S(\delta') - 1 + 2S(\nu^1) -c + S^k,$ where $c = 1$ or $2$, the latter only if $\delta = a\omega_1 + b\omega_{l+1}.$    On the other hand, by 
Lemmas \ref{V2} and \ref{weightsum} this must be at most $(S^0 + S(\delta') + S(\nu^1) + S^k) + 1.$  This reduces to $S(\delta') +S(\nu^1) \le c+2.$
 By hypothesis, $S(\delta') \ge 2.$   Also $S(\nu^1) \ge 2$ by Lemma \ref{S(nui)}(ii), and  $S(\nu^1) = S(\delta') + 1$ if $\delta = a\omega_1 + b\omega_{l+1}.$
 So this is a contradiction and we have shown that   $\la \lambda, \gamma_1\ra = 0.$
 
 (ii)  Assume that $i=k$ and $\la \lambda, \gamma_k\ra \ne 0.$  Replace both $W$ and $V$ by their duals.  So $C^0$ and $C^k$ are replaced by $\tilde C^0$ and $\tilde C^k$ where $L_X'$ is embedded via $(\d'')^*$ and $(\d')^*$, respectively. And we reverse the labelling on $V$. 
  Following  the  argument of (i), we obtain a contradiction unless $S(\delta'') \le 1$.  So assume $S(\delta'') \le 1$. 
Now $\delta'' \ne 0$, as otherwise $\delta = a\omega_1$ against our hypothesis.  Therefore,
  $\delta'' = \omega_s$ and $\delta = a\omega_1 + \omega_{s+1}.$  Now return to $W$ and $V$. Lemma \ref{V2}(iii) implies that $S_k^2 = S^0 + S^k + S(\nu^{k-1}) + S(\delta'').$ Also $2 \le S(\nu^{k-1}) \le S(\delta'') + 1,$  (the first inequality from Lemma \ref{1label}).  This implies that $S_k^2 \le S^0 + S^k + 3$, giving the inequality in (ii).  Also 
Lemmas \ref{i=1,k} and \ref {V2}(iv) force $\la \lambda, \gamma_j \ra = 0$ for $j < k.$  \hal

  \begin{lemma}\label{aomega1omega(l+1)}    If $\delta = a\omega_1 + \omega_{l+1}$, then  
$\la \lambda, \gamma_k\ra = 0.$    
  \end{lemma}
  
  \pf   By way of contradiction assume $\delta = a\omega_1 + \omega_{l+1}$ and $\la \lambda, \gamma_k\ra  \ne 0.$
Recall that $a\ge 3$ (by the hypothesis of Theorem \ref{L(delta')ge2}), and note that $W^k(Q_X) = \d'' = \o_l$. 
  Then Lemma \ref{i=1,k} together with   Lemma \ref{gamma1}(i),(ii) imply that $\la \lambda, \gamma_j\ra = 0$  for $1 \le j < k.$ Consider $V^*$ where $(\mu^*)^0$ has a nonzero label at node $l+1.$
  Applying the induction hypothesis we find that $(\mu^*)^0 = \lambda_{l+1}^0$, $\la \lambda, \gamma_k\ra = 1,$  $\mu^k = 0,$  $l =2,3, 4$, and
  $a \le 6,4,3,$ respectively.  In view of Theorem \ref{MUIZERO}   we know that $(\mu^*)^j = 0$ for $j \ne 0,k.   $
  
  Next we consider the possibilities for $\mu^0$, given by Lemma \ref{indhyp}(ii).
We shall first rule out the cases where $\mu^0$ is inner (i.e. the dual of one of the nonzero weights listed in 
 Lemma \ref{indhyp}(ii)). So suppose $\mu^0$ is inner. Then Theorem \ref{munotinner} implies that $i\ne 1$, hence $i=k$. 
Since $\mu^k=\mu^{k-1}=0$ and $\la \l,\g_k\ra = 1$, we have $S_k^2 = S^0+3$, and therefore
\begin{equation}\label{s0s1}
S^0+3 > S_1^2.
\end{equation}
On the other hand, for each possibility for $\mu^0$ we can compute a lower bound for $S_1^2$ using a suitable summand of $V_1^2$, as in Table \ref{follo}.
\begin{table}[h]
\caption{}\label{follo}
\[
\begin{array}{llll}
\hline
\mu^0 & S^0 & V_1^2 \supseteq & S_1^2 \ge \\
\hline
c\l_{r_0}^0 & ac & (\l_{r_0-1}^0+(c-1)\l_{r_0}^0)\otimes \l_1^1 & (c+2)a \\
\l_{r_0-c}^0 & \le (c+1)a-c & \l_{r_0-c-1}^0\otimes \l_1^1 & 4a-2\,(c=1) \\
(1\le c\le 3) &  &  & 5a-4\,(c=2) \\
& & & 6a-7\,(c=3) \\
\l_{r_0-4}^0 & \le 11\,(l\ge 3) & \l_{r_0-5}^0\otimes \l_1^1 & 14\,(l\ge 3) \\
(a=3) & 7\,(l=2) & & 11\,(l=2) \\
\l_{r_0-1}^0+\l_{r_0}^0 & 3a-1 & 2\l_{r_0-1}^0\otimes \l_1^1 & 5a-1 \\
\l_1^0+\l_{r_0}^0 & 2a & \l_1^0 \otimes \l_{r_0-1}^0\otimes \l_1^1 & 4a \\
\hline
\end{array}
\]
\end{table}
In all cases the inequality (\ref{s0s1}) is violated, except for the case where $\mu^0 = \l_{r_0-3}^0$, $a=3$. In this case we compute that $S^0\le 7$ (resp. 9) for $l=2$ (resp. $l\ge 3$), while $S_1^2 \ge 11$ (resp. 13), again contradicting (\ref{s0s1}). 

Hence $\mu^0$ is not inner. Now Lemma \ref{indhyp}(ii) implies that
\[
\mu^0 = 0, \,c\lambda_1^0 \,(c \le 4),\, \lambda_c^0 \,(2\le c \le 5), \hbox{ or } \lambda_1^0 + \lambda_2^0\, (a = 3).
\]
We will work through these possibilities with the aid of Magma.
  
  If $\mu^0 = 0,$ then $\lambda = \lambda_{n-l}.$ Then 
$V \downarrow X = (\wedge^{l+1}(a\omega_1 + \omega_{l+1}) )^*$ and a Magma computation shows that this fails to be MF.  Next assume that $\mu^0 = c\lambda_1^0 \, (c \le 4).$  Then
  $(V^*)^1 = \wedge^{l+1}(a\omega_1) \otimes (c\omega_1)$ and a Magma computation shows that
  this is not MF.  
  
  Now suppose $\mu^0 =  \lambda_c^0 \,(c \le 5).$  As $(\mu^*)^{k-1} = 0$, the bounds on $l$ restrict
  the possible values of $c$.  Indeed,  $c \le 3,4,5$ if $l = 2,3,4,$ respectively.  First assume $c = l+1.$  Now $V^2 = (V_1^2 + V_k^2)\downarrow L_X'$.  
The first summand is $\wedge^l(a\om_1) \otimes (a\om_1 + \om_l) + ((a-1)\om_1) = \wedge^l(a\om_1) \otimes (a\om_1) \otimes  (\om_l).$   Now $\wedge^l(a\om_1) \otimes (a\om_1)$ contains   $\wedge^{l+1}(a\om_1)$ as a direct
  summand.  Therefore $V_1^2\downarrow L_X'$ contains $\wedge^{l+1}(a\om_1) \otimes \om_l.$  On the other hand 
$V_k^2\downarrow L_X'$ contains  $\wedge^{l+1}(a\omega_1) \otimes (\omega_1 + \omega_l) \otimes \omega_l$ and  a Magma computation shows that this is not MF. This contradicts Corollary \ref{cover}.
  
    Now assume $c \le l.$  Therefore,
  $(V^*)^1 = \wedge^{l+1}(a\omega_1) \otimes (\omega_c)$ and again Magma shows that this is not MF.  
  
   Finally, assume that  $\mu^0 =  \lambda_1^0 + \lambda_2^0$ with $a = 3.$  Then
 $(V^*)^1 = \wedge^{l+1}(3\omega_1) \otimes (\omega_1+\omega_2).$  The
 first tensor factor has an irreducible summand with two nonzero labels, so the result is not MF by
Proposition \ref{tensorprodMF}.  \hal

 \begin{lemma}\label{inot1}  If $i = 1,$ then $\la\lambda, \gamma_j\ra = 0$ for all $j .$
 \end{lemma}
 
 \pf  Suppose $i = 1$.  By Lemma \ref{gamma1}(i), $\la\lambda, \gamma_j\ra = 0$ for $j < k$, so we need only show $\la\lambda, \gamma_k\ra = 0.$  Suppose
 false.  Then Lemma \ref{V2} shows that  $S_k^2 = S^0 +  S(\nu^{k-1}) + S^k + S(\delta'').$   On the other hand,  Lemma \ref{notinner} shows that $S_1^2 \le S^0 +  S^k + c,$  
where $c=2$ if $L(\d')\ge 2$, and $c=4$ if $\d = a\o_1+b\o_{l+1}$.
This forces $S(\nu^{k-1}) + S(\delta'') \le c$ respectively.  
If $L(\delta'') \ge 2,$  then Lemma \ref{1label} implies $S(\nu^{k-1}) \ge 2$.  But as $c = 2$ here, we have  $S(\delta'') = 0$, which is impossible.  Now assume $\delta = a\omega_1 + b\omega_{l+1}.$  Then   as in the proof of 
Lemma \ref{i=1,k}(ii), $S(\nu^{k-1}) = b+1.$  Therefore,  $S(\nu^{k-1}) + S(\delta'') = 2b+1$ which forces $b = 1.$
But this contradicts Lemma \ref{aomega1omega(l+1)}.    \hal
 
 \begin{lemma}\label{Vk2bound}  Assume that $\la\lambda, \gamma_k\ra = 0$, that $\delta'' = b\omega_s$ with 
$b > 1$ and $1 < s < l$,  and also that $(b,s) \ne (2,2).$ 
Then  $\mu^k \ne \lambda_1^k, \,2\lambda_1^k, \,\lambda_2^k,\, \lambda_3^k\, (b=2)$, or 
$\lambda_1^k + \lambda_{r_k}^k.$ 
 \end{lemma}
 
 \pf  The hypothesis on $\delta''$ implies that $\delta = a\omega_1 + b\omega_{s+1}$.  Therefore, the embedding
 of $L_X'$ in $C^{k-1}$ corresponds to the  representation
 $(\omega_1 +b\omega_s) \oplus ((b-1)\omega_s + \omega_{s+1}).$  We will work through the various cases with the aid of Lemmas \ref{aibjnotzero} and \ref{july4/17}.  
 
 First assume $\mu^k = \lambda_1^k$.  Then $V_k^2 = V_{C^{k-1}}(\lambda_{r_{k-1}}^{k-1}) \otimes 
V_{C^k}(\lambda_2^k) \otimes V_{C^0}(\mu^0) .$  Restricting
to $L_X'$ this contains $(\omega_{l-s} + (b-1)\omega_{l-s+1}) \otimes (\omega_{s-1} + (2b-2)\omega_s + \omega_{s+1}) \otimes J^0.$  Therefore we can apply Lemma \ref{july4/17} to see that this contains a composition factor
of multiplicity 2 and $S$-value at least $3b-2 + S^0$.  On the other hand, $S(V^1) = b + S^0$,
so this is a contradiction.

Next suppose $\mu^k = \lambda_2^k.$  Here $V_k^2 \downarrow L_X' \supseteq (\omega_{l-s} + (b-1)\omega_{l-s+1}) \otimes (3b\omega_s - (2\alpha_s + \alpha_{s \pm 1})) \otimes J^0$,  where we use $\alpha_{s + 1}$ if $s \le \frac{l}{2}$ and $\alpha_{s-1}$ otherwise.  An application of Lemma \ref{july4/17} implies that there is a repeated composition factor of $S$-value at least $4b-2 + S^0$ if $l > 3$ and $4b-3$ if $l = 3$.  On the other hand 
$S(V^1) = 2b + S^0,$ so this is a contradiction unless $l =3$ and $b = 2$.  But 
this  is ruled out by hypothesis.

If $\mu^k = \l_3^k$, then the induction hypothesis and the hypothesis of the lemma force 
$\delta'' = 2\omega_{l - 1}$  or $2\omega_2$, and so $\d = a\o_1+2\o_l$ or $a\o_1+2\o_3$. 
The argument is very similar to the previous cases.  Indeed in the first case $V_k^2 \downarrow L_X' \supseteq (\omega_1 + \omega_2) \otimes (\o_{l-3}+\o_{l-2}+4\o_{l-1}+\o_l) \otimes J^0.$  An application of Lemma \ref{aibjnotzero}(ii)   yields a repeated composition factor of $S$-value $8 + S^0$, which is a contradiction since $S(V^1) \le 6 + S^0.$  Similar reasoning applies in the second case.

Now suppose $\mu^k = 2\lambda_1^k.$  Then $V_k^2 \supseteq V_{C^{k-1}}(\lambda_{r_{k-1}}^{k-1}) \otimes V_{C^k}(\lambda_1^k + \lambda_2^k) \otimes V_{C^0}(\mu^0)$.  Now $V_{C^k}(\lambda_1^k + \lambda_2^k) =  (V_{C^k}(\lambda_1^k) \otimes V_{C^k}(\lambda_2^k)) - V_{C^k}(\lambda_3^k)$ so a weight consideration 
shows that $V_{C^k}(\lambda_1^k + \lambda_2^k) \downarrow L_X' \supseteq (\omega_{s-1} + (3b-2)\omega_s + \omega_{s + 1}).$  Therefore, $V_k^2 \downarrow L_X' \supseteq (\omega_{l-s} + (b-1)\omega_{l-s+1}) \otimes (\omega_{s-1} + (3b-2)\omega_s + \omega_{s + 1}) \otimes J^0$,  and Lemma \ref{aibjnotzero}(ii)  implies that there is a repeated composition factor with $S$-value at least $4b-2 + S^0.$  On the other hand, $S(V^1) = 2b + S^0$ 
and we get the usual contradiction.

Finally, assume  $\mu^k = \lambda_1^k + \lambda_{r_k}^k.$ Then $V_k^2 \supseteq V_{C^{k-1}}(\lambda_{r_{k-1}}^{k-1}) \otimes V_{C^k}(\l_2^k) \otimes V_{C^k}(\l_{r_k}^k) \otimes V_{C^0}(\mu^0).$   Restrict to $L_X'$. Then Lemma \ref{aibjnotzero}(ii) 
 applies to the product of the first two terms in the 4-fold tensor product to yield a repeated composition factor of $S$-value at least $3b-1$ and tensoring this with the other two terms we obtain a repeated factor with $S$-value at least $4b-1 + S^0.$   As $S(V^1) = 2b + S^0$,
this is a contradiction.  \hal
 
\vspace{4mm}
\noindent  {\bf Proof of Theorem \ref{L(delta')ge2}} 
 
In view of Lemma \ref{inot1} we may assume that $i \ne 1,$   so $i=k$ and $S_k^2 > S_1^2$ by Lemma \ref{i=1,k}.  
By way of contradiction, assume  that $\la\lambda, \gamma_j\ra \ne 0$  for some $j$.

If $\la\lambda, \gamma_k\ra \ne 0,$ then Lemma   \ref{gamma1}(ii) shows
 that $\delta'' = \omega_s$ and  $\delta = a\omega_1 + \omega_{s+1}.$  Lemma \ref{aomega1omega(l+1)} implies that 
$s < l$, and  Lemma \ref{diminequal} gives $r_0>r_k+2$. 
Consider $V^*$ and the corresponding restriction $(\mu^*)^0 = \lambda^* \downarrow C^0.$ 
Apply the induction hypothesis to $(\mu^*)^0$, which has a nonzero label at node  $1+r_k$ and possibly others. As  $L(\delta') \ge 2$ this 
forces $(\mu^*)^0 = \l_3^0$, $s=l-1$ and $a=1$, and this can only occur if  $l = 2$ and $\d = 110$. 
  Returning to $V$ this gives $\mu^k = 0$ and hence $S_k^2 = 3+S^0$.  
Inductively, we have $\mu^0 = 0$, $\l_1^0$, $\l_2^0$, $\l_3^0$, $2\l_1^0$, or the dual of one of these. 
In each case we calculate $S^0+3$, and we find that the inequality $S_k^2 > S_1^2$ forces $\mu^0 = 0,$ $\l_1^0$, $\l_2^0$, $\l_3^0$ or $\l_4^0 (=(\l_3^0)^*)$.  At this point a Magma computation shows that $V \downarrow X$ is not MF.

Therefore  $\la\lambda, \gamma_k\ra = 0$ and $j < k.$  This forces $\mu^k \ne 0$, since otherwise 
Lemma \ref{gammai20} shows that $V_k^2 = 0$ and $S_k^2 = 0$, which is a contradiction as $i=k$.
Further, it follows from Lemmas \ref{i=1,k}(ii) and \ref{V2}(iv) that  $S_k^2 \ge S^0 + S^k + 4.$

  If $L(\delta'') \ge 2$ then taking the duals of $V$ and $W,$  Theorem \ref{munotinner} shows that $\mu^k$ is outer, and then Lemma \ref{notinner}(ii) shows that $S_k^2 \le S^0 + S^k + 2,$  a contradiction. Therefore  $L(\delta'') =1$ so that $\delta'' = b\omega_s$ for some $s$, and hence $\delta = a\omega_1 + b\omega_{s+1}.$
It then follows from Lemma \ref{V2}(ii) that $S_k^2 \le S^0 + S^k + (b+1) + b = S^0 + S^k + 2b+1.$
 
Next we claim that $\la\lambda, \gamma_1\ra = 0.$  Otherwise  Lemma \ref{V2}(iii) implies that 
$S_1^2 = S^0 + S^k + (a+b) + (a+b)$ if $s \le l$   or $S_1^2 = S^0 + S^k + (a) + (a+1) = S^0 + S^k + 2a+1 $
if $s = l+1.$  In either case this violates the fact that $S_k^2 > S_1^2$ (in the latter case $a \ge b$ by hypothesis).
Therefore, $\la\lambda, \gamma_1\ra = 0.$

We summarize the information we have at this point: 
\begin{equation}\label{summ1}
\begin{array}{l}
\d = a\o_1+b\o_{s+1} \hbox{ and }i=k, \\
\la \l,\g_1\ra = \la \l,\g_k\ra = 0, \\
\mu^k \ne 0,\,\mu^t = 0 \hbox{ for }1\le t < k, \\
S_k^2 \ge S^0+S^k+4.
\end{array}
\end{equation}

Suppose first that $s\le l-1$, so that $L(\d')\ge 2$. Consider the dual $V^*$. The inductive hypothesis implies that $(\mu^*)^0$ or its dual is in $\{ 0,\, \lambda_1^0,\, \lambda_2^0,\, \lambda_3^0,\, 2\lambda_1^0,\, 3\lambda_1^0 \,(l=2)\}$. 

Assume that $r_0\ge r_k+3$. Then $(\mu^*)^0$ cannot be inner (otherwise $\mu^k$ would be 0). Hence 
\[
\mu^k \in \{\l_{r_k}^k, \, 2\l_{r_k}^k, \, 3\l_{r_k}^k\,(l=2), \, \l_{r_k-1}^k, \, \l_{r_k-2}^k \}.
\]
But now arguing just as in the proof of Lemma  \ref{notinner}(ii), we see that $S_k^2\le S^0+S^k+2$, contradicting the bound in (\ref{summ1}).

Now assume $r_0<r_k+3$. Then Lemma \ref{diminequal}
implies that $b > 1$  and also that $(b,s) \ne (2,2).$ It then follows from Lemma \ref{Vk2bound} that $\mu^k =  \lambda_{r_k}^k$, $2\lambda_{r_k}^k$, $\lambda_{r_k-1}^k$, or $\lambda_{r_k-2}^k\, (b=2).$  As above we argue as in the proof of Lemma \ref{notinner}(ii) to see that $S_k^2 \le S^0 + S^k + 2$, which is a contradiction. 

We have now established that $s=l$, that is, $\d = a\o_1+b\o_{l+1}$. By hypothesis, $a\ge 3$ and $a\ge b\ge 1$.

Consider $V^*$. The possibilities for $(\mu^*)^0$ are given by Lemma \ref{indhyp}(ii). 
We claim that $(\mu^*)^0$ is outer. If $a=b$, then $\d=\d^*$ and $S_1^2$ is maximal for $V^*$, so the claim follows from Theorem \ref{munotinner}. And if $a>b$, then $r_0 = \dim V_{A_l}(a\o_1)-1$, $r_k = \dim V_{A_l}(b\o_{l})-1$, and we compute that $r_0 \ge r_k+d$, where $d = 10$ if $l\ge 3$, $d = 5$ if $l=2$ and $a\ge 4$ or $(a,b)=(3,1)$, and $d=4$ if $l=2$ and $(a,b)=(3,2)$. Since $\mu^k\ne 0$ it follows that $(\mu^*)^0$ is outer.

 It follows that 
\[
\mu^k = \l_{r_k-c}^k\,(c\le 4),\,c\l_{r_k}^k\,(c\le 4), \hbox{ or }\l_{r_k-1}^k+\l_{r_k}^k\,(a=3).
\]

At this point we return to the proof of  Lemma \ref{notinner}(ii),(iii), but working with $S_k^2$ rather than $S_1^2$. The
aim is again to show that $S_k^2 \le S^0 + S^k + 3$, contradicting (\ref{summ1}). 
Note that this inequality fails in the case where $l=2$, $b=2,3$ and $\mu^k = \l_{r_k-4}^k$. These will be treated later, so assume for now that we are not in these cases.

The argument in the proof of Lemma \ref{notinner}(ii),(iii) shows that for $b\ge 3$, the inequality 
$S_k^2 \le S^0 + S^k + 3$ holds; and minor modifications yield the same conclusion for $b\le 2$.

It remains to handle the cases excluded above, where $l=2$, $b=2,3$ and $\mu^k = \l_{r_k-4}^k$. 
Suppose first $b=2$. Here $V^1 = (V_{C^0}(\mu^0) \downarrow L_X') \otimes (02)$, so $S(V^1) = S^0+2$, while 
$V^2 \supseteq (V_{C^0}(\mu^0) \downarrow L_X') \otimes ((12) + (01)) \otimes (12) \supseteq V_{C^0}(\mu^0) \downarrow L_X' \otimes (13)^2$.  Hence $V^2$ has a repeated composition factor of $S$-value $S(V^1)+2$, a contradiction. Similarly for $b=3$, $S(V^1) = S^0+7$ and $V^2$ has a repeated composition factor of $S$-value $S^0+9$, again a contradiction.

This completes the proof of Theorem \ref{L(delta')ge2}.

\chapter{Proof of Theorem \ref{MAINTHM}, Part IV: Completion}\label{pfend}

In this final chapter, we complete the proof of Theorem \ref{MAINTHM}. We continue with the notation of the previous three chapters. Here is the first main result of the chapter.

\begin{theor}\label{SUPP} Assume the inductive hypothesis, and suppose that $L(\d')\ge 2$. 
 Then $\l,\d$ are as in Table $\ref{TAB1}$ of Theorem $\ref{MAINTHM}$.
\end{theor}

In the next result we assume that $L(\d') = L(\d'')=1$, 
in which case we have $\d = a\om_1+b\om_{l+1}$, where $a,b>0$.

\begin{theor}\label{OM1OML}
Assume the inductive hypothesis, and suppose $\d = a\om_1+b\om_{l+1}$ with $a\ge b\ge 1$ and $a\ge 2$.  Then either $b=1$ or $(a,b)=(2,2)$. Moreover, up to duals we have $\lambda = 2\lambda_1$ or $ \lambda_2$,  and also $\l \ne 2\l_1$
if $(a,b)=(2,2)$, as in Table $\ref{TAB1}$ of Theorem $\ref{MAINTHM}$.
\end{theor}

These results complete the proof of Theorem \ref{MAINTHM}: the case $l = 1$ is handled in Chapter \ref{casea2}; the case where $L(\d)=1$ is in Chapters \ref{rkge2}--\ref{o2sec}; the case $L(\d') \ge 2$ (or $L(\d'')\ge 2$) is done in Theorem \ref{SUPP}; and finally the case where $\d = a\om_1+b\om_{l+1}$ is handled by Theorem \ref{OM1OML} together with Chapter \ref{o1ol1}.

\section{Proof of Theorem \ref{SUPP}}

In this section we prove Theorem \ref{SUPP}. Assume the hypothesis of the theorem.  Then
Theorems \ref{MUIZERO} and \ref{L(delta')ge2} imply that every nonzero $\l$-label is on $\mu^0$ or $\mu^k$.
That is, $\mu^i = 0$ for $i \ne 0,k$, and $\la \l,\gamma_j\ra = 0$ for all $j.$

 \subsection{The case where $\mu^0\ne 0, \mu^k =0$}\label{sec1711}

We assume in this subsection that every nonzero $\l$-label is on $\mu^0$.

By Proposition \ref{induct} we know that $V^1 = V_{C^0}(\mu^0)\downarrow L_X'$ is MF. The natural module for $C^0$ is $W^1(Q_X) \cong V_{L_X'}(\d')$, so using the inductive hypothesis we see that $\mu^0,\d'$ are one of the pairs in Table \ref{muposs}. 

\begin{table}[h]
\caption{}\label{muposs}
\[
\begin{array}{|l|l|}
\hline
\mu^0\;(\hbox{up to duals}) & \d'\; (\hbox{up to duals}) \\
\hline
\l_1^0 & \hbox{any} \\
2\l_1^0, \,\l_2^0 & \om_1+c\om_i,\,c\om_1+\om_i,\,\om_i+c\om_{i+1},\,c\om_i+\om_{i+1} \\
\l_2^0 & 2\om_1+2\om_2, 2\om_1+2\om_l , \om_2 + \om_{l-1}(l\ge4), \om_2 + \om_4 \\
\l_3^0 & \om_1+\om_l \\
3\l_1^0 & \om_1+\om_2\,(l=2) \\
\hline
\end{array}
\]
\end{table}

\no Note that $\mu^0$ is not $\l_1^0$, since otherwise $\l=\l_1$, contrary to assumption.

\begin{lem}\label{innerquote} We have $\mu^0 \ne \l_{r_0}^0,\, 2\l_{r_0}^0,\,\l_{r_0-1}^0,\, \l_{r_0-2}^0$ or $3\l_{r_0}^0$.
\end{lem}

\pf  By Lemma \ref{gammai20}, all $S_i^2$ are 0 except for $S_1^2$. Hence the conclusion follows from Theorem \ref{munotinner}. \hal

\begin{lem}\label{easyfirst}
If $\mu^0 = 2\l_1^0,\,\l_2^0,\,\l_3^0$ or $3\l_1^0$ as in Table $\ref{muposs}$, then $\l,\d$ are as in Table $\ref{TAB1}$ of Theorem $\ref{MAINTHM}$.
\end{lem}

\pf In this case we have $\l =  2\l_1,\l_2,\l_3$ or $3\l_1$ and $\d = \d'+x\om_{l+1}$ for some $x$.   

First we claim that $x=0$. Suppose false. Then $L(\d'')=2$ by the inductive hypothesis, so $\d' \ne \o_2+\o_{l-1}$ or $\o_2+\o_4$. For the remaining possibilities, $V_Y(\l)\downarrow X$ is not MF by Lemmas \ref{donna1}, \ref{wed3ex} and \ref{lemma1.5(iii)(iv)}(ii), a contradiction.

Hence $x=0$. Suppose now that 
\[
\d \in \{\om_1+c\om_i,\,c\om_1+\om_i,\,\om_i+c\om_{i+1},\,c\om_i+\om_{i+1},\,2\o_1+2\o_2,\,\o_2+\o_4\}.
\]
Then with the exception of $\d = \o_1+\o_l$, the inductive hypothesis implies that $\l=2\l_1$ or $\l_2$ and $(\l,\d)$ are as in Table \ref{TAB1}, as required. For the exceptional case, $\l=2\l_1$, $\l_2$, $\l_3$ or $3\l_1\,(l=2)$, and the last two do not occur by Lemmas \ref{wed3ex} and \ref{lemma1.5(iii)(iv)}(ii).

The remaining possibilities for $\d$ are 
\begin{equation}\label{list71}
\begin{array}{l}
\o_i+c\o_l\,(i\ne 1,l-1),\; c\o_i+\o_l\,(i\ne 1,l-1),\; \o_{l-3}+\o_{l-1}\,(l>5), \\
2\o_{l-1}+2\o_l\,(l\ge 3),\; \o_2+\o_{l-1}\,(l>5),\;  2\o_1+2\o_l\,(l\ge 3).
\end{array}
\end{equation}
The induction hypothesis implies that $\l=2\l_1$ or $\l_2$, although only $\l_2$ is possible in the last four cases. 

In the first two cases of (\ref{list71}) with $i \ge 4$ or $i=3$, $c>1$, consideration of $V^*$ gives  a contradiction using the induction hypothesis. The same contradiction applies in the third case with $l\ge 7$, and in the fourth with $l\ge 4$. This leaves the following possibilities for $\d$:
\begin{itemize}
\item[{\rm (1)}] $\d = \o_2+c\o_l,\,c\o_2+\o_l$,
\item[{\rm (2)}] $\d = \o_3+\o_l$,
\item[{\rm (3)}] $\d = \o_2+\o_{l-1}\,(l>5)$,
\item[{\rm (4)}] $\d = 2\o_{1}+2\o_l\,(l\ge 3)$,
\item[{\rm (5)}] $\d = \o_{l-3}+\o_{l-1}\,(l=6)$ or $2\o_{l-1}+2\o_l\,(l=2,3)$.
\end{itemize}
If $\l=\l_2$ and $\d$ is as in  case (1) with $c=1$,  or case (5) with $l=2$, then  $\l,\d$ are as in Table \ref{TAB1} of Theorem \ref{MAINTHM}. 
In all other cases we find that $V_Y(\l)\downarrow X$ is not MF: using Magma for (5), and using Lemmas \ref{0c0...010}, \ref{om3oml}, \ref{010...0100} and \ref{donna1}(ii) for cases (1),(2),(3),(4) respectively. \hal

 \subsection{The case where $\mu^0\ne 0, \mu^k \ne 0$}\label{bothne0}

Assume in this subsection that both $\mu^0$ and $\mu^k$ are nonzero.

\begin{lem}\label{1both}
We have $L(\d'')=1$, so that $\d = a\om_1+b\om_j$ and $\d''=b\om_{j-1}$ (where $2\le j\le l$).
\end{lem}

\pf If $L(\d'')\ge 2$ then by Proposition \ref{twolabs}, $V_{C^0}(\mu^0) \downarrow L_X'$ and 
$V_{C^k}(\mu^k) \downarrow L_X'$ have composition factors with highest weights $\nu_0$ and $\nu_k$ (respectively) such that $L(\nu_0)\ge 2$ and $L(\nu_k)\ge 2$. But $V^1$ contains $\nu_0\otimes \nu_k$,
which is not MF by Proposition \ref{tensorprodMF}, contradicting Proposition \ref{induct}(i). \hal

\begin{lem}\label{2both}
$\mu^k$ and $\d''$ are as in Table $\ref{mukposs}$.
\end{lem}

\begin{table}[h]
\caption{}\label{mukposs}
\[
\begin{array}{l|l}
\hline
\mu^k & \d'' \\
\hline
\l_1^k \hbox{ or } \l_{r_k}^k & \hbox{any }b\om_{j-1} \\
2\l_{r_k}^k    & \om_1,\,2\om_1,\,\om_2,\,\om_{l-1} \\
\l_{r_k-1}^k   & \om_1 \\
\l_{r_k-2}^k & \om_2\,(l=3) \\
3\l_{r_k}^k & \om_1\,(l=2) \\
\hline
\end{array}
\]
\end{table}

\pf  Assume that $\mu^k \ne \l_1^k$  or $\l_{r_k}^k$.
Since $V_{C^k}(\mu^k) \downarrow L_X'$ is MF, we can inductively assume that $\mu^k$ and $\d''$ are as in Tables \ref{TAB2}--\ref{TAB4} of Theorem \ref{MAINTHM}. 

As in the previous proof,  $V_{C^k}(\mu^k) \downarrow L_X'$ can have no composition factor with $L$-value greater than 1. Using Lemma \ref{twolabs}, this implies that 
\[
\d'' = \o_2,\,\o_{l-1},\,\o_1 \hbox{ or }2\o_1.
\]
Correspondingly, $\d' = a\om_1+\om_3$, $a\om_1+\om_l$, $a\o_1+\o_2$ or $a\om_1+2\om_2$. Hence Proposition \ref{diminequal} shows that $\dim V_{L_X}(\d') > \dim V_{L_X}(\d'')+2$ -- that is, $r_0>r_k+2$. It follows that the possibilities for $(\mu^*)^0$ are $2\l_1^0$, $\l_2^0$, $\l_3^0$ ($\d'=\o_1+\o_l$), $3\l_1^0$ ($l=2, \d'=\o_1+\o_2$) or the dual of one of these. This implies that $\mu^k$ is outer, and now a further application of Lemma \ref{twolabs} shows that $\mu^k$ is one of the following possibilities:
\[
\begin{array}{l}
2\l_{r_k}^k, \\
\l_{r_k-1}^k\,(\d''=\o_1), \\
\l_{r_k-2}^k\,(l=2,\d''=2\o_1 \hbox{ or }l=3,\d''=\o_2), \\
3\l_{r_k}^k\,(l=2, \d''=\o_1 \hbox{ or }\o_2).
\end{array}
\]
We can exclude the possibility $\mu^k = \l_{r_k-2}^k$ with $l=2,\d''=2\o_1$, as here we would have $\d = a20$, $(\mu^*)^0 = \l_3^0$, but $\wedge^3(a2)$ is not MF for $L_X'=A_2$, by Chapter \ref{casea2}; we can also exclude $3\l_{r_k}^k$ with $\d''=\o_2$, since $L(\d')>1$. It follows that $\mu^k, \d''$ are as in Table \ref{mukposs}. 
\hal

\begin{lem}\label{notlik}
$\mu^k$ is not $\l_1^k$.
\end{lem}

\pf Suppose $\mu^k = \l_1^k$. 

We first show that $r_0\ge r_k$. Assume false. 
 Since $\mu^k=\l_1^k$, for the dual $V^*$ we have 
$(\mu^*)^0 = 0$, and hence all the support of $\l^*$ is on $(\mu^*)^k$. As $\mu^0$ and $\mu^k$ are nonzero,
$(\mu^*)^k$ therefore has at least two nonzero labels, and it is MF on restriction to $L_X'$ via the module $b\om_{j-1}$.
Now  
\[
\dim V_{L_X}(\d') = r_0+1 < r_k+1 = \dim V_{L_X}(\d''),
\]
so Proposition \ref{diminequal} implies that $l>2$, $b\ge 2$ and $b\om_{j-1}$ is not $2\om_1$ or $3\om_1$. Hence the inductive list 
implies that $(\mu^*)^k$ is the adjoint representation $\l_1^k+\l_{r_k}^k$. It follows that $\l = \l_1+\l_{n-r_k+1}$. 

We now show that this cannot occur, under the slightly weaker assumption that $r_k\ge r_0-2$. First observe that $V^1 = (a\om_1+b\om_j)\otimes (b\om_{j-1})$ has $S$-value $a+2b$. Next, $V^2$ is the sum of $W^2(Q_X) \otimes (b\om_{j-1})$ and $(a\om_1+b\om_j)\otimes W^k(Q_X)^* \otimes \wedge^2(b\om_{j-1}).$  By Theorem \ref{LEVELS},  $W^2(Q_X)$ and $W^k(Q_X)$ have $S$-values at most $a+b+1$ and $b+1$ respectively. Hence 
$V^2$ has $S$-value at most $a+4b+1$. Now consider the next level 
$V^3$. If we write $\g = \b_{n-r_k}$ and let $\a$ be the adjacent node in $C^{k-1}$ and $\b,\d$ the two nearest nodes in $C^k$, then $\l-\a-2\g-2\b-\d$ affords the level 2 summand
\[
(a\om_1+b\om_j)\otimes \wedge^2(W^k(Q_X))^* \otimes \wedge^3(b\om_{j-1}).
\]
Now $W^k(Q_X) = (\om_1+b\om_{j-1})\oplus ((b-1)\om_{j-1}+\om_j)$ has $S$-value $b+1$, so 
by Lemma \ref{mult2}, $ \wedge^2(W^k(Q_X))$ has a summand of multiplicity at least 2 and $S$-value at least 
$2(b+1)-1$. Therefore $V^3$ has a repeated summand of $S$-value at least $a+b+(2b+1)+(3b-2)$. It now follows from Proposition \ref{sbd} that $a+6b-1 \le a+4b+2$, hence $b=1$. But then $r_0>r_k+2$ by Proposition \ref{diminequal}, a contradiction.

We have now shown that $r_0\ge r_k$. Consider the dual $V^*$. The weight $(\mu^*)^0$ involves $\l_{r_k}^0$ with coefficient 1, so as its restriction to $L_X'$ is MF, we see from the list that 
\[
(\mu^*)^0 = \l_{r_k}^0\; \hbox{ and } \; r_k \in \{2,3,r_0-2,r_0-1,r_0\}.
\]

Suppose $r_k=2$. Then $l=2$ and $\d''=\om_1$. The weight $\mu^0$ is in Table \ref{muposs}; consideration of $V^*$
shows that in fact   
$\mu^0$ must be one of $\l_1^0,\,2\l_1^0,\,\l_2^0,\,3\l_1^0$ (not the dual). Hence $\l = \l_1+\l_{n-1}$, $2\l_1+\l_{n-1}$, 
$\l_2+\l_{n-1}$ or $3\l_1+\l_{n-1}$. In the first case $\l^* = \l_2+\l_n$, so in $V^*$,  level 0 for $L_X'$ is $\wedge^2(a\om_1+\om_2) \otimes \om_2$; in the other cases level 0 in $V$ is $S^2(a\om_1+\om_2) \otimes \om_1$, $\wedge^2(a\om_1+\om_2) \otimes \om_1$ or $S^3(a\om_1+\om_2) \otimes \om_1$.
None of these are MF by Proposition \ref{a2one}, which is a contradiction. Hence $r_k\ne 2$.

If $r_k=3$ then $l=3$ and $\d'' = \om_1$. Then $(\mu^*)^0 = \l_3^0$. But then in $V^*$, level 0 for $L_X'$ is 
$\wedge^3(a\om_1+\om_2)$, which is not MF by  the inductive hypothesis.

Hence $r_0 \ge r_k \ge r_0-2$. If $r_k = r_0-2$ then in the dual $V^*$ we have $(\mu^*)^0 = \l_{r_0-2}^0$, so that $\wedge^3(\d')^*$ must be MF. Then by inspection of the inductive list we have $\d' = \om_1+\om_l$.  But then $r_0>r_k+2$ by Proposition \ref{diminequal}.

Therefore $r_k = r_0$ or $r_0-1$. 
Also $(\mu^*)^k$ and $\d''$ must be as in Table \ref{mukposs}. If $(\mu^*)^k \ne \l_1^k$ or $\l_{r_k}^k$, then the table shows that $b=1$ or $(b,j-1) = (2,1)$, so Proposition \ref{diminequal} implies that $r_0>r_k+2$, a contradiction. Therefore 
$(\mu^*)^k = \l_1^k$ or $\l_{r_k}^k$.

If $(\mu^*)^k = \l_{r_k}^k$ then for $V$ we have $\lambda = \l_1+\l_{n-r_k+1}$. This case was handled above. Hence 
$(\mu^*)^k = \l_1^k$ and so for $V$ we have 
\[
\l = \l_{r_0}+\l_{n-r_k+1}\;\hbox{ or }  \l_{r_0-1}+\l_{n-r_k+1}.
\]
Consider the first case, $\l = \l_{r_0}+\l_{n-r_k+1}$. Here $V^1 = (\d')^*\otimes \d''$ has $S$-value $a+2b$. At level 1, $V^2$ is the sum of $\wedge^2(\d')^* \otimes W^2(Q_X) \otimes \d''$ and 
$(\d')^* \otimes W^{k}(Q_X)^* \otimes \wedge^2(\d'')$. Hence 
$V^2$ has $S$-value at most the maximum of $2S(\d')+S(W^2(Q_X))+S(\d'')$ and 
$S(\d')+S(W^{k}(Q_X))+2S(\d'')$, which is $3a+4b$. Finally, at level 2,  $V^3$ has a summand 
$\wedge^3(\d')^* \otimes \wedge^2(W^2(Q_X)) \otimes \d''$. By Lemma \ref{mult2}, 
$\wedge^2(W^2(Q_X)) $ has a repeated summand of $S$-value at least $2(a+b)-1$, and so at level 2
there is a multiplicity 2 summand of $S$-value at least $(3S(\d')-2)+(2a+2b-1)+b$, which is at least $5a+6b-3$.
Therefore by Proposition \ref{sbd} we have $5a+6b-3 \le 3a+4b+1$, hence $a+b\le 2$. But then $\d' = \om_1+\om_j$ which implies that $r_0>r_k+2$ by Proposition \ref{diminequal}, a contradiction.

The last case, $\l =  \l_{r_0-1}+\l_{n-r_k+1}$, is very similar: here, $V^2$ has $S$-value equal to that of
$\wedge^3(\d')^* \otimes W^2(Q_X) \otimes \d''$, which is at most $4a+5b$, while 
$V^3$ has a summand $\wedge^4(\d')^* \otimes \wedge^2(W^2(Q_X)) \otimes \d''$, hence has 
a composition factor of multiplicity at least 2 and $S$-value at least $(4S(\d')-3)+(2a+2b-1)+b$. Now Proposition \ref{sbd} yields $6a+7b-4 \le 4a+5b+1$, which gives a contradiction by Proposition \ref{diminequal} as before.  \hal

Now we complete this section with

\begin{lem}\label{final both}
Theorem $\ref{SUPP}$ holds when $\mu^0$ and $\mu^k$ are both nonzero.
\end{lem}

\pf By Lemmas \ref{2both} and \ref{notlik}, we have $\mu^k = c\l_{r_k}^k$, $\l_{r_k-1}^k$ or  $\l_{r_k-2}^k$, 
with $\d''$ as in 
Table \ref{mukposs}. Therefore, considering $V^*$ we have $(\mu^*)^0 = c\l_1^0$, $\l_2^0$ or $\l_3^0$ as in Table \ref{muposs}. As $\lambda^*$ has at least $2$ nonzero labels, $(\mu^*)^k\ne 0$ and $(\mu^*)^k$ must also be as in Table \ref{mukposs}.

Assume first that $\mu^k = 3\l_{r_k}^k$, so that $\d'' = \om_1$, $l=2$. Then $V^*$ has $(\mu^*)^0 = 3\l_1^0$ and $\d' = \om_1+\om_2$. So 
$(V^*)^1$ is $S^3(\om_1+\om_2) \otimes S^c(\om_i)$ with $c\le 3$ and $i=1$ or 2, and this is not MF by Proposition \ref{mags}(i).

Next, let $\mu^k = \l_{r_k-2}^k$, so that $\d'' = \om_2$, $l=3$. Then $V^*$ has $(\mu^*)^0 = \l_3^0$ and $\d' = \om_1+\om_3$. But then 
$(V^*)^1$ is not MF by Proposition \ref{mags}(ii).

Thus $\mu^k$ is $\l_{r_k}^k$, $2\l_{r_k}^k$ or  $\l_{r_k-1}^k$. By considering $V^*$ it follows also that 
$(\mu^*)^0$ is $\l_1^0$, $2\l_1^0$ or  $\l_2^0$. So by  Table \ref{mukposs} and the above applied to $(\mu^*)^k$,  
$\l$ is one of the following or its dual:
\[
\l_1+\l_n,\,\l_2+\l_n,\, \l_2+\l_{n-1},\, 2\l_1+\l_n,\, 2\l_1+\l_{n-1},\, 2\l_1+2\l_n.
\]

In the first case $V\downarrow X$ is the adjoint module $V_X(\d) \otimes V_X(\d)^*/V_X(0)$, and this is not MF by 
Proposition \ref{tensorprodMF}  (noting that the trivial module $V_X(0)$ has multiplicity 1 in the tensor product). 

In the second case, the dual $\l^*=\l_1+\l_{n-1}$ and so Table \ref{mukposs} implies that $\d'' = \om_1$, hence $\d' = a\om_1+\om_2$. But then in $V$, we have $V^1 = \wedge^2(a\om_1+\om_2) \otimes \om_l$, which is not MF by Proposition \ref{morenonmf}(i).

Next consider $\l=\l_2+\l_{n-1}$ or $2\l_1+\l_{n-1}$. Here Table \ref{mukposs} gives $\d'' = \om_1$. Then 
$V^1 = \wedge^2(a\om_1+\om_2) \otimes \wedge^2(\om_l)$ or 
$S^2(a\om_1+\om_2) \otimes \wedge^2(\om_l)$, neither of which is MF by Proposition \ref{morenonmf}(i, iv).

Finally, suppose $\l=2\l_1+i\l_n$ with $i = 1$ or 2. Then $\l^* = i\l_1+2\l_n$ and Table \ref{mukposs} gives $\d'' = \om_1, 2\om_1, \om_2$ or $\om_{l-1}$. 
In each case Proposition \ref{morenonmf}(ii,iii) shows that $V^1$ is not MF. \hal

 \subsection{The case where $\mu^0 = 0, \mu^k \ne 0$}

In this subsection we complete the proof of Theorem \ref{SUPP} by handling the case where $\mu^0 = 0, \mu^k \ne 0$.
Assume that this holds. Considering $V^*$, we may also assume by the previous sections that 
$(\mu^*)^0 = 0, (\mu^*)^k \ne 0$. In particular this means that $r_0<r_k$. If $L(\d'')\ge 2$, then replacing $\d$ by $\d^*$ and $V$ by $V^*$ leads to the case dealt with in Section \ref{sec1711}. Thus we assume that $L(\d'')=1$, so that for some $j$ with $2\le j\le l$ we have $\d = a\o_1+b\o_j$, hence
\[
\d' = a\om_1+b\om_j, \;\d'' = b\om_{j-1}.
\]

\begin{lem}\label{init} The following hold.
\begin{itemize}
\item[{\rm (i)}] $(\mu^*)^0 = 0, (\mu^*)^k \ne 0$ and $\dim V_{L_X'}(\d') < \dim V_{L_X'}(\d'')$.
\item[{\rm (ii)}] $j\ge 3$, $b\ge 2$ and $(j,b)\ne (3,2)$.
\item[{\rm (iii)}] $\mu^k$ is $\l_1^k$, $2\l_1^k$, $\l_2^k$ or $\l_3^k$.
\item[{\rm (iv)}] If $\mu^k = \l_3^k$ then $b=2$ and $j=l$.
\end{itemize}
\end{lem}

\pf  (i) This just  records the assumptions in the preamble above.

 (ii) If $j=2$, or $b=1$, or $(j,b)=(3,2)$, then $\dim V(\d') > \dim V(\d'')$ by Proposition \ref{diminequal}, contradicting (i).

(iii) Since $\d'' = b\om_{j-1}$ with $b\ge 2$ and $2\le j-1 \le l-1$, the inductive list implies that 
 $\mu^k$ is $\l_1^k$, $2\l_1^k$, $\l_2^k$, $\l_3^k$ or the dual of one of these. However it is not a dual, since 
$(\mu^*)^0 =0$. Hence (iii) holds.

(iv) From the inductive list and (ii), if $\mu^k = \l_3^k$ then $\d''$ must be $2\om_{l-1}$. \hal

We deal with each of the possibilities in Lemma \ref{init}(iii) in turn. Recall that 
\[
W^{k}(Q_X) = (\om_1+b\om_{j-1}) \oplus ((b-1)\om_{j-1}+\om_j).
\]
Write $\a,\g,\b,\e$ for the roots $\b_{n-r_k-1},\,\b_{n-r_k},\,\b_{n-r_k+1},\,\b_{n-r_k+2}$ respectively.

\begin{lem}\label{l1kposs}
$\mu^k$ is not $\l_1^k$.
\end{lem}

\pf Suppose $\mu^k = \l_1^k$. Then $V^1 = \d''$ has $S$-value $b$. Also $V^2(Q_Y)$ is afforded by
$\l-\g-\b$, hence $V^2 = W^{k}(Q_X)^* \otimes \wedge^2(\d'')$, which has 
$S$-value at most $(b+1)+2b = 3b+1$.

Now consider $V^3(Q_Y)$. This has a summand afforded by $\l-\a-2\g-2\b-\e$, and hence $V^3$ has a summand $\wedge^2(W^{k}(Q_X))^* \otimes \wedge^3(\d'')$. By Lemma \ref{mult2} the first tensor factor has a repeated composition factor of $S$-value at least $2b+1$. Also since $\d'' = b\om_{j-1}$ with $2\le j-1\le l-1$, $\wedge^3(\d'')$ has a composition factor of $S$-value at least $3b-1$. Hence  $V^3$ has a summand of multiplicity at least 2 and $S$-value at least $(2b+1)+(3b-1)$. It now follows from Proposition \ref{sbd} that $5b \le (3b+1)+1$, whence $b=1$, contrary to Proposition \ref{init}(ii). \hal

\begin{lem}\label{l2kposs}
$\mu^k$ is not $\l_2^k$.
\end{lem}

\pf This is similar to the previous proof. Suppose $\mu^k = \l_2^k$. Then $V^1 = \wedge^2(\d'')$ has $S$-value $2b$. Also $V^2(Q_Y)$ is afforded by $\l-\g-\b-\e$, hence $V^2 = W^{k}(Q_X)^* \otimes \wedge^3(\d'')$, which has $S$-value at most $(b+1)+3b = 4b+1$. Finally,  $V^3$ has a summand $\wedge^2(W^{k}(Q_X))^* \otimes \wedge^4(\d'')$, which by Lemma \ref{mult2} has a repeated summand of $S$-value at least $(2b+1)+(4b-2)$. Hence by Proposition \ref{sbd} we have $6b-1\le 4b+2$, so $b=1$, a contradiction. \hal

\begin{lem}\label{2l1kposs}
$\mu^k$ is not $2\l_1^k$.
\end{lem}

\pf Suppose $\mu^k = 2\l_1^k$, so $V^1 = S^2(\d'')$ has $S$-value $2b$. Then $V^2(Q_Y)$ is afforded by $\l-\g-\b$, so that 
$V^2$ is contained in  $W^{k}(Q_X)^* \otimes (\d'' \otimes \wedge^2(\d''))$, which has $S$-value at most $4b+1$.

Now $V^3(Q_Y)$ has a summand afforded by $\l-2\g-2\b$, and this restricts to $C^{k-1}\times C^k$ as 
$2\l_{r_{k-1}}^{k-1} \otimes 2\l_2^k$. Also $2\l_{r_{k-1}}^{k-1} \downarrow L_X' = S^2(W^{k}(Q_X))^* $, which by Lemma \ref{mult2}  has a repeated composition factor of $S$-value at least $2b+1$. Consequently $V^3$ has a repeated composition factor of $S$-value at least $6b+1$. Now Proposition \ref{sbd} yields $6b+1 \le 4b+2$, a contradiction. \hal

\begin{lem}\label{l3kkposs}
$\mu^k$ is not $\l_3^k$.
\end{lem}

\pf Assume $\mu^k = \l_3^k$, so that $\d'' = 2\om_{l-1}$ by Lemma \ref{init}(iv). Then 
$V^2= W^{k}(Q_X)^* \otimes \wedge^4(\d'')$ has $S$-value at most 11.
Also $V^3$ contains $\wedge^2(W^{k}(Q_X))^* \otimes \wedge^5(\d'')$.
The first tensor factor contains a repeated summand of $S$-value 5 by Lemma \ref{mult2}, while the second has
a summand of $S$-value at least 8. Hence Proposition \ref{sbd} gives $13 \le 11+1$, a contradiction. \hal

\vspace{6mm} 
This completes the proof of Theorem \ref{SUPP}.

\section{Proof of Theorem \ref{OM1OML}: case $a\ge 3$}

In this section we assume that  
\begin{equation}\label{a3b2}
\d = a\om_1+b\om_{l+1} \hbox{ with } a\ge 3 \hbox{ and } a \ge b \ge 1.
\end{equation}
We aim to prove Theorem \ref{OM1OML} in this case. 

Theorems \ref{MUIZERO} and \ref{L(delta')ge2} imply that every nonzero $\l$-label is on $\mu^0$ or $\mu^k$.
That is, $\mu^i = 0$ for $i \ne 0,k$ , and $\la \l,\gamma_j\ra  = 0$ for all $j.$  This also holds for the dual $V^*.$
Moreover, as $a\ge b$ we have $r_0\ge r_k$, and hence, replacing $V$ by $V^*$ if necessary we may assume that 
\[
\mu^0 \ne 0.
\]
Moreover, if $a = b$, then we may assume that $S_1^2 \ge S_k^2$. Note that 
\[
\d' = a\om_1,\;\d'' = b\om_l.
\]
Therefore $\dim V_{L_X'}(\d') < \frac{1}{2}\dim V_X(\d)$, and it follows that $r_0<\frac{1}{2}n$  (recall that $n$ is the rank of $Y = SL(W)$).

\begin{lem}\label{mu0pos}
$\mu^0$ is one of the following, up to duals:
\[
\begin{array}{l}
c\l_1^0\;(c\le 4) \\
\l_c^0\; (c\le 5) \\
\l_1^0+\l_{r_0}^0 \\
\l_1^0+\l_2^0\; (a=3)
\end{array} 
\]
\end{lem}

\pf We know that $V_{C^0}(\mu^0)\downarrow L_X'$ is MF and that the natural module for $C^0$ restricts to $L_X'$ as $\d' = a\om_1$ with $a\ge 3$. Hence the conclusion is immediate from the inductive list of examples. \hal

Recall that we call $\mu^0$ an {\it outer} weight for $C^0$ if it is one of the weights listed in the conclusion of Lemma \ref{mu0pos}, and an {\it inner} weight if it is the dual of one of these. Note that $\l_1^0+\l_{r_0}^0$ is both outer and inner according to this definition.

\begin{lem}\label{mukpossib}
Assume $\mu^k\ne 0$.
\begin{itemize}
\item[{\rm (i)}] If $a>b$ then $\mu^k$ is one of the following:
\[
\begin{array}{l}
c\l_{r_k}^k\; (c\le 4) \\
\l_{r_k-c}^k\; (c\le 4) \\
\l_{r_k-1}^k+\l_{r_k}^k\; (a=3)
\end{array} 
\]
\item[{\rm (ii)}] If $a=b$ then $\mu^k$ is one of the weights in {\rm (i)} or its dual, or $\l_1^k+\l_{r_k}^k$.
\end{itemize}
\end{lem}

\pf  In the dual $V^*$, the weight $(\mu^*)^0$ must be as in Lemma \ref{mu0pos}. Suppose $a>b$ and $(\mu^*)^0$ is an inner weight. Then $(\mu^*)^0$ has a nonzero label on a node $\b_i^0$ for some $i\ge r_0-4$. Then $\mu^k$ has a nonzero label on node $\b_{n-i+1}$, and it follows that $r_k\ge i$. Hence $r_0-r_k\le 4$, so that 
\[
\dim V_{L_X'}(a\om_1) - \dim V_{L_X'}(b\om_l) \le 4.
\]
The only possibility is that $a=3,b=2$ and $l=2$. Here $r_0=9$, $r_k=5$ and also $(\mu^*)^0$ must be $\l_5^0$. However we check using Magma that $\wedge^5(3\om_1)\otimes (V_{C^k}((\mu^*)^k)\downarrow L_X')$ is not MF for all possibilities for $(\mu^*)^k$ (noting that $(\mu^*)^k \ne 0$).

Hence we have shown that if $a>b$ then $(\mu^*)^0$ is not inner.  The conclusions follow. \hal

\begin{lem}\label{deltaabnotinner}  $\mu^0$ is not inner.
\end{lem}

\pf Suppose $\mu^0$ is inner.
Assume first that $a>b$. Here we argue similiarly to the previous proof. For the $L_X'$-modules $\d',\d''$ we have $\dim (a\om_1) \ge 10$ and 
$\dim (a\om_1)-\dim (b\om_l) \ge 4$, with equality if and only if $l=2$, $a=3$ and $b=2$.
Now $\mu^0$ has a nonzero label on a node $\b_i^0$ for some $i \ge r_0-4$. It follows that for the dual $V^*$, 
either $\lambda^*$ has a nonzero label which is not on $C^0$ or $C^k$, or $l=2$, $a=3$, $b=2$ and $\mu^0 = \l_5^0$.
The former is impossible, so the latter holds. Then $\d = 3\om_1+2\om_3$ and 
\[
V^1 = \wedge^5(3\om_1) \otimes (V_{C^k}(\mu^k)\downarrow L_X').
\]
This must be MF. If $\mu^k = 0$ then $\l = \l_5$ and $V\downarrow X = \wedge^5(\d)$ which is not MF by a Magma computation.
And if $\mu^k\ne 0$ then this case was ruled out in the proof of Lemma \ref{mukpossib}.

Now assume $a=b$. Here we  are assuming that $S_1^2 \ge S_k^2,$ so that $S_1^2 = S(V^2)$ and  Theorem \ref{munotinner} gives the conclusion.    \hal

\begin{lem}\label{muknon}
If $\mu^k=0$, then $\l,\d$ are as in Table $\ref{TAB1}$ of Theorem $\ref{MAINTHM}$.
\end{lem}

\pf Suppose $\mu^k=0$. Then Lemmas \ref{mu0pos} and \ref{deltaabnotinner} imply that $\l$ is one of 
\[
c\l_1\;(2\le c\le 4),\;\;\l_c\;(2\le c\le 5),\;\;\l_1+\l_2\;(a=3).
\]
Now Lemmas \ref{wedgesquareab} and \ref{manymany} show that the only possibility is that $\delta = a\omega_1 + \omega_{l+1}$ with $\lambda = 2\lambda_1$ or $\lambda_2$. These possibilities are in Table \ref{TAB2}.  \hal

In view of the previous result, we assume now that $\mu^k\ne 0$. 

\begin{lem}\label{mukout}
 $\mu^k$ is as in Lemma $\ref{mukpossib}${\rm (i)}.
\end{lem}

\pf Assume false. Then by Lemma \ref{mukpossib}, $a=b \ge 3$ and $(\mu^*)^0$ is inner.  Recall that we are assuming that 
 $S_1^2 \ge S_k^2.$  If equality holds, then we can apply Theorem \ref{munotinner} to $V^*$ to get the result.  So suppose the
inequality is strict.  Then $\mu^0 \ne \l_1^0$, as otherwise $S_1^2 = a+1,$  while $S_k^2 \ge a+1.$ Therefore, Lemma \ref{twolabs} implies that $\mu^k = \l_1^k$ or $\l_{r_k}^k.$  In the latter case, the conclusion holds.  So assume
that $\mu^k = \l_1^k$.  Then $S_k^2 = S^0 + (a+1) + S(\wedge^2(a\omega_l)) = S^0 + (a+1) + (2a-1) = S^0 + S^k + 2a.$  On the other hand, considering the possibilities in Lemma \ref{mu0pos} we see that $S_1^2 \le S^0 + S^k +4$.  Indeed, this is worked out in the proof of Lemma \ref{notinner}.  So in all cases $S_1^2 < S_k^2$, a contradiction. \hal

\begin{lem}\label{finalposs}
If $b \ge 2$, then one of the following holds:
\begin{itemize}
\item[{\rm (i)}] $\mu^0 = \l_1^0$
\item[{\rm (ii)}] $\mu^k = \l_{r_k}^k$, $2\l_{r_k}^k$ (with $b=2$) or $\l_{r_k-2}^k$ (with $b=2,l=2$).
\end{itemize}
\end{lem}

\pf   The possibilities for $\mu^0,\mu^k$ are given by the previous lemmas. As $V^1 = (V_{C^0}(\mu^0) \otimes V_{C^k}(\mu^k))\downarrow L_X'$ is MF,  Proposition \ref{tensorprodMF} shows that one of the tensor factors must have all its $L_X'$-composition factors of $L$-value at most 1. Hence Lemma \ref{twolabs} gives the conclusion. \hal

\begin{lem}\label{not0}
If $b \ge 2,$ then $\mu^0$ is not $\l_1^0$.
\end{lem}

\pf Assume $\mu^0 = \l_1^0$. By Lemma \ref{mukout}, $\mu^k$ is as in Lemma \ref{mukpossib}(i).
Hence $V^1$ is one of the following possibilities:
\[
\begin{array}{l}
(a\om_1) \otimes S^c(b\om_1)\;\; (c \le 4) \\
(a\om_1) \otimes \wedge^c(b\om_1) \;\;(c\le 5) \\
(a\om_1) \otimes ((b\om_1)\otimes \wedge^2(b\om_1)/\wedge^3(b\om_1))\;\;(a=3,\; b=2 \hbox{ or } 3).
\end{array}
\]
However, with the exception of $c = 1$ in the first case, or $(b,c,l) = (2,2,l)$, $(2,4,2)$, $(2,5,2)$ in the second case, none of these is MF, by Lemma \ref{nonemf}. In the exceptional cases, the possibilities for $\l,\d$ are:
\[
\begin{array}{ll}
\hline
\l & \d \\
\hline
\l_1+\l_n & a\om_1+b\om_{l+1} \\
\l_1+\l_{n-1}  & a\om_1+2\om_{l+1} \\
\l_1+\l_{n-3} & a\om_1+2\om_3\;(l=2) \\
\l_1+\l_{n-4} & a\om_1+2\om_3\;(l=2) \\
\hline
\end{array}
\]
For the first line of the table, we have 
$V \downarrow X = ((a\omega_1 + b\omega_{l+1}) \otimes (b\omega_1 + a\omega_{l+1})) - 0$ which is not MF by Proposition \ref{tensorprodMF}.  The other cases are not MF by Lemma \ref{donnathreeseven}.  \hal

\begin{lem}\label{ifb2} If $b \ge 2$, then $\l,\d$ are as in Table $\ref{TAB2}$ of Theorem $\ref{MAINTHM}$.
\end{lem}
 
\pf  By Lemmas \ref{mu0pos},  \ref{deltaabnotinner}, and \ref{not0}, $\mu^0$ is $c\l_1^0\,(2\le c\le 4)$, 
$\l_c^0\,(2\le c\le 5)$ or $\l_1^0+\l_2^0\,(a=3)$.
Also, by the previous two lemmas, 
$\mu^k$ is  $\l_{r_k}^k$, $2\l_{r_k}^k$ (with $b=2$) or $\l_{r_k-2}^k$ (with $b=2,l=2$).

If $\mu^k =\l_{r_k}^k$, then the dual $V^*$ has $(\mu^*)^0= \l_1^0$ and Lemma \ref{not0} gives a contradiction.
Hence $\mu^k$ is $2\l_{r_k}^k$ (with $b=2$) or $\l_{r_k-2}^k$ (with $b=2,l=2$). 
Write $Z = V_{C^k}(\mu^k)\downarrow L_X'$, so that $Z$ is $S^2(2\om_1)$ or $\wedge^3(2\om_1)$ (with $l=2$ in the latter case). Then the possibiliities for  $V^1$ are
\[
\begin{array}{l}
S^c(a\om_1) \otimes Z\;\; (2\le c \le 4) \\
\wedge^c(a\om_1) \otimes Z \;\;(2\le c\le 5) \\
((3\om_1)\otimes \wedge^2(3\om_1)/\wedge^3(3\om_1)) \otimes Z. 
\end{array}
\]
However, none of these is MF, by Lemma \ref{nonemf}. \hal

  At this point we have $\delta = a\omega_1 + \omega_{l+1} = (a0 \ldots  01)$ with $a\ge 3$.  
Note that this implies that $C^k = A_l.$

Recall that both $\mu^0$ and $\mu^k$ are nonzero, the possibilities being given by  Lemmas \ref{mu0pos}, \ref{mukpossib}, and also $\mu^0$ is not inner (Lemma \ref{deltaabnotinner}).

\begin{lem}\label{naturals} $\mu^0 \ne \l_1^0$ and $\mu^k \ne \l_l^k.$
\end{lem}

\pf  Assume false.  Then taking duals, if necessary, we may assume $\mu^k = \l_l^k.$
If $\mu^0 = \l_1^0$, then  $V \downarrow X = (a0 \ldots  01) \otimes (10 \ldots  0a) - 0$.  But then Lemma \ref{compfactor}(iii) implies that there is a repeated composition factor of
weight $(a0 \ldots  0a),$ a contradiction.   Therefore,  $\mu^0 \ne \l_1^0.$

Suppose $\mu^0 = c\l_1^0$ for $2 \le c \le 4.$  Then $\lambda = c\lambda_1 + \lambda_n$ and it follows that 
 $V = (S^c(W) \otimes W^*) - S^{c-1}(W).$  Restricting to $X$ we have 
$S^c(a0 \ldots  01) \otimes (10\ldots 0a) - S^{c-1}(a0 \ldots  01).$  The tensor product contains
$((ca)0 \ldots  0c) \otimes (10\ldots 0a) \supseteq ((ca)0\ldots 0(c+a-1))^2$ by Lemma \ref{compfactor}(iii). An 
$S$-value comparison shows   that this does not appear in $S^{c-1}(a0 \ldots  01),$  so $V \downarrow X$
is not MF.

Now assume $\mu^0 = \l_c^0$ for $2 \le c \le 5.$  Then $\lambda = \lambda_c + \lambda_n$ and it follows that 
$V = (\wedge^c(W) \otimes W^*) - \wedge^{c-1}(W).$  Restricting to $X$ we have 
$\wedge^c(a0 \ldots  01) \otimes (10\ldots 0a) - \wedge^{c-1}(a0 \ldots  01).$
First assume $c = 2,3,4.$  Then $\wedge^c(a0 \ldots  01) \supseteq ((2a-2)10 \ldots  02), ((3a-6)30\ldots 03),$ or $((4a-12)60\ldots 04),$ respectively.  Therefore, Lemma \ref{compfactor}(i) implies that the tensor product contains
a repeated composition factor of highest weight $((2a-2)10 \ldots  0(a+1)), ((3a-6)30\ldots 0(a+2)),$ or $((4a-12)60\ldots 0(a+3)),$ respectively.  These have $S$-values $3a, 4a-1,5a-3,$ respectively. On the other hand
$S(\wedge^{c-1}(a0 \ldots  01))  = a+1, 2a+1, 3a+1$  respectively.  
Therefore, an $S$-value comparison shows that the repeated composition factor cannot appear in $\wedge^{c-1}(a0 \ldots  01)$ if $c = 2,3$ or $4$ and $V \downarrow X$ is not MF. 

Now assume $\mu^0 = \l_5^0$ where $a = 3$. Here the argument is slightly different. Consideration of $V^*$ implies  that $l \ge 5.$ It follows 
that $\wedge^5(30 \ldots  01) \ge (360\ldots  013) + (2510\ldots 05),$  so tensoring with $(10\ldots 03)$ we
have $(3510\ldots 016)^2.$  Therefore we have a repeated composition factor of $S$-value $16$ which is larger than the $S$-value of any composition factor of $\wedge^4(30\ldots 01)$. So again $V \downarrow X$ fails to be MF.

The final case is where $\mu^0 = \l_1^0 + \l_2^0$ with $a = 3.$  
Here we will argue that $V^1(Q_Y) = V_{C^0}(\mu^0) \otimes V_{C^k}(\mu^k)$ is not MF.  We have 
$V_{C^0}(\mu^0) = (V_{C^0}(\l_1^0) \otimes V_{C^0}(\l_2^0)) - V_{C^0}(\l_3^0).$
Restricting to $L_X'$ this is $((30\ldots  0) \otimes \wedge^2(30\ldots 0)) - \wedge^3(30\ldots 0).$
The tensor product contains $(30\ldots  0) \otimes (410\ldots 0). $ We use Magma to see
that this contains  $ (520\ldots 0) + (4110\ldots 0)$ ($ (52) + (41)$ if $l = 2$).   Therefore, $V^1(Q_Y) \downarrow L_X' \supseteq
((520\ldots 0) + (4110\ldots 0)) \otimes (10\ldots  0) \supseteq (5110\ldots 0)^2,$  ($(51)^2$ if $l = 2$), a contradiction.
\hal
 
 \begin{lem}  $\mu^0 \ne d\l_1^0$ and  $\mu^k  \ne d\l_l^k$ for $d \ge 2.$ 
 \end{lem}
 
 \pf  Assume false.  By taking duals we can assume $\mu^k  = d\l_l^k.$  We work through the possibilities in Lemma \ref{mu0pos} other than the case $\mu^0 \ne \l_1^0 + \l_{r_0}^0,$ which was ruled out in Lemma \ref{deltaabnotinner}. 
In all cases, $V^1$ is not MF by Lemma \ref{nonemf}. \hal

\begin{lem}\label{c3514} 
 It is not the case that $\mu^0 =  \l_c^0$ for $3 \le c \le 5$  and  $\mu^k   =\l_{l-d}^k$ for $1 \le d \le 4.$
\end{lem}

\pf Suppose  the assertion is false.   Then   $l \ge d+1$  and considering   $V^*$ we have $l \ge c$.   In all cases  
$V^1$ fails to be MF by Lemma \ref{nonemf}, a contradiction.  \hal

\vspace{4mm}
We can now complete the proof of Theorem \ref{OM1OML} in the case where $a\ge 3$. By Lemma \ref{ifb2} we can  assume that $b=1$. The previous lemmas imply that the possibilities for $\mu^0,\mu^k$ are as follows:
\[
\begin{array}{ll|l}
& \mu^0 & \mu^k \\
\hline
{\rm (1)} & \l_2^0 & \l_{l-d}^k\,(1\le d\le 4),\;\l_{l-1}^k+\l_l^k\,(a=3) \\
{\rm (2)} & \l_c^0\,(3\le c\le 5) & \l_{l-1}^k+\l_l^k\,(a=3) \\
{\rm (3)} &\l_1^0+ \l_2^0\,(a=3) & \l_{l-d}^k\,(1\le d\le 4),\;\l_{l-1}^k+\l_l^k \\
\hline
\end{array}
\]
In case (3), either $V^1$ or $(V^*)^1$ is not MF by Lemma \ref{nonemf}, together with Proposition \ref{tensorprodMF}. In case (2), the dual $V^*$ is as in (3). Finally, consider case (1). Here either $V^*$ is as in (2) or (3), or we have 
\[
\mu^0=\l_2^0,\;\mu^k= \l_{l-1}^k.
\]
Here $\l=\l_2+\l_{n-1}$ and $V \downarrow X =  (\wedge^2(a0\ldots 01) \otimes \wedge^2(10\ldots 0a)) - ((a0\ldots 01) \otimes (10\ldots 0a)).$ Then the
first tensor product contains $((2a-2)10 \ldots  02) \otimes (010\ldots 0(2a))$  and this contains $((2a-1)10\ldots 0(2a+1))^2$ (see Lemma \ref{compfactor}(iii)), so this is a final contradiction.  

\section{Proof of Theorem \ref{OM1OML}: case $a=2$}

In this section we complete the proof of Theorem \ref{OM1OML}  by handling the case where $a=2$. So assume that 
\begin{equation}\label{a3b2}
\d = 2\om_1+b\om_{l+1},\; b =1 \hbox{ or } 2.
\end{equation}
Note that the results of Chapters \ref{ch15} and \ref{ch16} do not apply here. 
We summarize some preliminary information for this case. As usual we abbreviate an $L_X'$-module by its highest weight (so that $2\om_1$ stands for $V_{L_X'}(2\om_1)$ and so on).

\begin{lem} \label{prela2} The following hold.
\begin{itemize}
\item[{\rm (i)}] $k=3$ if $b=1$, and $k=4$ if $b=2$.
\item[{\rm (ii)}] As $L_X'$-modules, $W^1(Q_X)\cong 2\om_1$, $W^{k+1}(Q_X) \cong b\om_l$ and $W^2(Q_X) \cong 
 2\om_1 \otimes \om_l$.
\item[{\rm (iii)}] If $b=2$ then $W^3(Q_X) \cong (2\om_1+2\om_l) \oplus (\om_1+\om_l)\oplus 0$ and 
$W^4(Q_X) \cong  \om_1 \otimes 2\om_l$.
\item[{\rm (iv)}] $\mu^i = 0$ for $1\le i \le k-1$.
\end{itemize}
\end{lem}

\pf Parts (i)-(iii) follow from Theorem \ref{LEVELS} and part (iv) from Theorem \ref{MUIZERO}. \hal

 As usual, for $1 \le i \le k$ let $\g_i$ be the node in the Dynkin diagram of $Y$ between $C^{i-1}$ and $C^i$, and let $x_i = \la \l,\g_i\ra$. 

\begin{lem}\label{mu0muk} At least one of $\mu^0$, $(\mu^*)^0$, $\mu^k$ and $(\mu^k)^*$ is nonzero. 
\end{lem}

\pf Suppose $\mu^0=\mu^k=0$. Then by Lemma \ref{prela2}(iv), the support of $\l$ is on the $\g_i$ only.
First consider the case where $b=1$. Here $k=3$ and the natural modules for $C^0,C^1,C^2, C^3$ have dimensions 
$\frac{1}{2}(l+1)(l+2)$, $\frac{1}{2}(l+1)^2(l+2)$, $(l+1)^2$, $l+1$, respectively. 
If any $x_i \ne 0$ then in $V^*$ we have $(\mu^*)^j \ne 0$ for some $j$; then $j=0$ or $k$ by Lemma \ref{prela2}(iv), giving the conclusion.

Now consider $b=2$. As $\mu^j=0$ for all $j$,  $V^1$ is a trivial module. If $x_1\ne 0$ then $\l-\g_1$ affords the highest weight of a summand of $V^2$ isomorphic to 
$W^1(Q_X)^* \otimes W^2(Q_X) \cong (2\om_l) \otimes  2\om_1 \otimes \om_l$. This has a summand $\om_1+2\om_l$ of multiplicity 2, which contradicts Proposition \ref{sbd}. Similarly, if $x_2\ne 0$ then $\l-\g_2$ affords a summand of 
$V^2$ isomorphic to $(\om_1+2\om_l)\otimes (2\om_1+2\om_l)$, which has a multiplicity 2 summand $2\om_1+3\om_l$, again contradicting Proposition \ref{sbd}. Hence $x_1=x_2=0$, and dualizing, $x_3=x_4=0$. But now we have forced $\l=0$, which is a contradiction. \hal

Replacing $V$ by $V^*$ if necessary, we may assume from now on that $\mu^0\ne 0$. 

\begin{lem}\label{mu0poss2}
{\rm (i)} If $b=1$ then $\mu^0$ is one of the following weights:
\[
\begin{array}{l}
c\l_1^0\,(c\ge 1),\, \l_i^0\,(i\ge 2),\,a\l_2^0\,(a\le 3), \\
\l_1^0+\l_i^0\,(2\le i \le 7),\,\l_2^0+\l_3^0,\,a\l_1^0+\l_2^0\,(a\le 3).
\end{array}
\]

{\rm (ii)} If $b=2$ then up to duals, $\mu^0$ is one of the weights in part (i), or one of the following:
\[
\begin{array}{l}
\l_1^0+\l_i^0\,(i\ge r_0-5),\,a\l_1^0+\l_{r_0}^0\,(a\le 3),\,\l_2^0+\l_{r_0-1}^0.
\end{array}
\]
\end{lem}

\pf The possibilities for $\mu^0$ can be listed using the inductive hypothesis for the weight $2\o_1$. 
The list consists of all the weights in (i) and (ii), together with their duals. 
If $b=1$ then a weight listed in part (ii), or the dual of one of those in either part, is not possible for $\mu^0$, since then the dual $V^*$ would satisfy $(\mu^*)^2 \ne 0$, contrary to Lemma \ref{prela2}(iv). \hal

\begin{table}[t]
\caption{Possibilities for $N$, $b=1$ or $2$}\label{summand}
\[
\begin{array}{|l|l|}
\hline
\mu^0 & \hbox{highest weights of } N \\
\hline
c\l_1^0\,(c\ge 2) & (c-2)\l_1^0+\l_2^0 \\
\l_i^0\,(i\ge 2) & \l_1^0+\l_{i-1}^0 \\
c\l_2^0\,(2\le c\le 3) & 2\l_1^0+(c-1)\l_2^0 \\
\l_1^0+\l_i^0\,(i\ge 3) &  2\l_1^0+\l_{i-1}^0 \\
\l_1^0+\l_2^0 &  3\l_1^0 \oplus \l_3^0 \\
\l_2^0+\l_3^0 & 2\l_1^0+\l_3^0 \\
a\l_1^0+\l_2^0\, (a=2,3) & (a+2)\l_1^0 \oplus ((a-2)\l_1^0+2\l_2^0) \\
\hline
\end{array}
\]
\end{table}

\begin{table}[t]
\caption{Extra possibilities for $N$, $b=2$}\label{summ2}
\[
\begin{array}{|l|l|}
\hline
\mu^0 & \hbox{highest weights of } N \\
\hline
c\l_{r_0}^0\,(c\ge 2) &  \l_1^0+\l_{r_0-1}^0+(c-1)\l_{r_0}^0 \\
c\l_{r_0-1}^0\,(2\le c \le 3) & \l_1^0+\l_{r_0-2}^0+(c-1)\l_{r_0-1}^0 \\
\l_i^0+\l_{r_0}^0\,(3\le i\le r_0-1) & \l_1^0+\l_{i-1}^0+\l_{r_0} \\
\l_2^0+\l_{r_0}^0 & \l_1^0+\l_2^0+\l_{r_0-1}^0 \\
a\l_1^0+\l_{r_0}\,(1\le a \le 3) &  (a+1)\l_1^0+\l_{r_0-1}^0 \\
\l_1^0+a\l_{r_0}^0\,(2\le a\le 3) & 2\l_1^0+\l_{r_0-1}^0+(a-1)\l_{r_0}^0 \\
\l_2^0+\l_{r_0-1}^0 & \l_1^0+\l_2^0+\l_{r_0-2}^0 \\
\l_{r_0-2}^0+\l_{r_0-1}^0 & \l_1^0+\l_{r_0-3}^0+\l_{r_0-1}^0 \\
\l_{r_0-1}^0+a\l_{r_0}^0\,(a=2,3) & \l_1^0+\l_{r_0-2}^0+a\l_{r_0}^0 \\
\hline
\end{array}
\]
\end{table}

\begin{lem}\label{v2sum}
{\rm (i)} $V^2$ has a submodule $M \cong V^1 \otimes \om_l$, and $V^2/M$ is MF.

{\rm (ii)} $V^2/M$ has a summand 
\[
(N\downarrow L_X') \otimes \om_l \otimes V_{C^k}(\mu^k)\downarrow L_X',
\] 
where $N$ is an irreducible $C^0$-module (or sum of two such) with highest weight as in 
Table $\ref{summand}$ if $b=1$, and is in Table $\ref{summand}$ or Table $\ref{summ2}$ if $b=2$.
\end{lem}

\pf This follows from Corollary \ref{coversums} and its proof.  The term in (i) is the
summand $A$ described prior to the statement of  Corollary \ref{coversums}, and $V^2/M$ is MF by Corollary \ref{cover}.  Now consider (ii). The statement of the Corollary covers all cases in Tables \ref{summand} and \ref{summ2} except for rows 4,5,7 and 4,5,6,7, respectively.  These cases are settled  using a slight adjustment to the proof.  For these  cases set $\nu' = \mu^0 - \beta_j - \cdots -\gamma_1, $  where  $j = i, 2, 2$  or 
$j =r_0, r_0,r_0,r_0-1$ respectively.  Then using $\nu'$ rather than $\nu$ in the proof of part (ii) of 
Corollary \ref{coversums}, we obtain the assertion, with the two exceptions of rows 5 and 7 of Table \ref{summand}.
For these cases we actually get two summands, one from $\nu$ and one
from $\nu' $ and this gives the assertion. \hal

\begin{lem}\label{mu0is}
$\mu^0$ is $ \l_1^0,\l_2^0$ or $2\l_1^0$.
\end{lem}

\pf Suppose false. Then Lemma \ref{v2sum}(i, ii) applies to give a summand $(N\downarrow L_X') \otimes V_{L_X'}(\om_l) \otimes V_{C^k}(\mu^k)\downarrow L_X'$ of $V^2/M$, which must be MF. 
The highest weight of $N$ must be on the inductive list of examples for $2\o_1$. Hence we see that $\mu^0$ is one of the following:
\[
2\lambda_1^0,\, 3\l_1^0,\,4\l_1^0,\,5\l_1^0,\, \l_i^0,\,2\l_2^0,\,\l_1^0+\l_2^0,\,2\l_1^0+\l_2^0,\,\l_1^0+\l_3^0.
\]

If $\mu^0 = 3\l_1^0$ then $N = 110\ldots 0 \cong (100\ldots 0)\otimes (010\ldots 0)/(001\ldots 0)$, so that $(N\downarrow L_X') \otimes V_{L_X'}(\om_l) \cong (2\om_1 \otimes \wedge^2(2\om_1)/\wedge^3(2\om_1))\otimes \om_l$, and it is readily seen that this has a summand $3\om_1+\om_2$ of multiplicity 2, so is not MF.

Next suppose $\mu^0 = \l_i^0$ with $i\ge 3$. Here $N = \l_1^0+\l_{i-1}^0$, and the inductive list shows that this can be MF only if either $i-1\le 7$ or $i-1\ge r_0-5$. In the first case we can use Magma for $l \le 4$ and otherwise Lemma \ref{nonMF}(3) to see that $(N\downarrow L_X') \otimes V_{L_X'}(\om_l) $ is not MF.  And in  the second case apply Lemma \ref{nonMF}(8) to get the same conclusion.  

If $\mu^0 = 2\l_2^0$, $4\l_1^0$ or $\l_1^0+\l_3^0$, then $N = 210\ldots 0$.  Then $N \downarrow L_X' \supseteq 
(6\o_1+\o_2) \oplus (4\o_1+2\o_2)$, and tensoring with $\om_l$ we obtain $(5\o_1+\o_2)^2$.

The remaining cases $\mu^0 = \l_1^0+\l_2^0$, $5\l_1^0$, $2\l_1^0+\l_2^0$ are similar: we find that 
$(N\downarrow L_X') \otimes \om_l$ has a multiplicity 2 summand of highest weight 
$3\om_1 + \om_2$, $7\o_1+\o_2$, $5\o_1+\o_2$, respectively. \hal

\vspace{4mm}
We can now complete the proof of Theorem \ref{OM1OML} in the case $a=2$. By Lemma \ref{mu0is}, $\mu^0 = \l_1^0$, $\l_2^0$ or $2\l_1^0$. 

Assume first that $\mu^0 = \l_1^0$. Then $\l-\b_1^0-\cdots - \b_{r_0}^0-\g_1$ affords the highest weight of a summand $V_{C^1}(\l_1^1) \otimes V_{C^k}(\mu^k)$ of $V^2(Q_Y)$, and the restriction of this to $L_X'$ is 
$M=2\om_1 \otimes \om_l \otimes (\mu^k)\downarrow L_X'$. Hence $V^2/M$ is MF, by Lemma \ref{v2sum}(i). 

We next argue that $x_i=0$ for all $i$, where $x_i  =\la\l,\g_i\ra$. 
If $x_1 \ne 0$, then $\l-\g_1$ affords a summand $V_{C^0}(\l_1^0+\l_{r_0}^0) \otimes 
V_{C^1}(\l_1^1) \otimes V_{C^k}(\mu^k)$ of $V^2(Q_Y)$. In the restriction of this to $L_X'$, the first two tensor factors have composition factors with highest weights having two nonzero labels, so by Proposition \ref{tensorprodMF} the restriction is not MF, a contradiction. Hence $x_1=0$. Similarly, if $x_i\ne 0$ for another value of $i$, then $\l-\g_i$ affords a non-MF summand of 
$V^2$. Hence $x_i=0$ for all $i$, as asserted. 

Since $\l\ne \l_1$ by the hypothesis of the theorem, we have $\mu^k\ne 0$. Let $j$ be minimal such that $\mu^k$ has a nonzero coefficient of $\l_j^k$. Then $\l-\b_j^k-\cdots -\b_1^k-\g_k$ affords the highest weight of a summand  $V_{C^0}(\l_1^0) \otimes 
V_{C^{k-1}}(\l_{r_{k-1}}^{k-1}) \otimes V_{C^k}(\mu^k-\l_j^k+\l_{j+1}^k)$ of $V^2(Q_Y)$, where the term $\l_{j+1}^k$ is not present if $j = r_k.$ The restriction of this to $L_X'$ is 
\[
2\om_1 \otimes (b\om_1\otimes \om_l) \otimes (\mu^k-\l_j^k+\l_{j+1}^k)\downarrow L_X'.
\]
 Using Proposition \ref{stem1.1.A}, we see that the only way this can be MF is if $\mu^k-\l_j^k+\l_{j+1}^k = 0$ and $\mu^k = \l_{r_k}^k$. This means that $\l = \l_1+\l_n$, the adjoint module for $Y$. But then $V\downarrow X = V(\d)\otimes V(\d^*)/V(0)$, and this is not MF by Proposition \ref{tensorprodMF}, a contradiction.

Next suppose $\mu^0 = 2\l_1^0$. By Lemma \ref{v2sum}(ii), the quotient $V^2/M$  has a summand $\wedge^2(2\om_1) \otimes \om_l \otimes (\mu^k)\downarrow L_X'$. This must be MF by Corollary \ref{cover}. Since 
$\wedge^2(2\om_1) $ has a summand with highest weight having 2 nonzero labels, this can only be the case if $\mu^k=0$ or $\l_l^k$. In the latter case the tensor product is not MF, so $\mu^k=0$.  
We claim that also $x_i= 0$ for all $i$. Indeed, if $x_k>0$ then $\l-\g_k$ affords a summand $V_{C^0}(2\l_1^0) \otimes 
V_{C^{k-1}}(\l_{r_{k-1}}^{k-1}) \otimes V_{C^k}(\l_1^k)$ of $V^2(Q_Y)$, and the restriction of this to $L_X'$ is 
$S^2(2\om_1)\otimes (b\om_1\otimes \om_l) \otimes (b\om_l)$, which is not MF. Similarly, if $x_i \ne 0$ for some other value of $i$, then $\l-\g_i$ affords a non-MF summand of $V^2$, proving the claim. We have now established that 
$\l = 2\l_1$, which for $\d = 2\om_1+\om_{l+1}$ is in the conclusion of Theorem \ref{OM1OML}. On the other hand
Lemma \ref{wedgesquareab} shows that $S^2(\d)$ is not MF if $b = 2.$

Finally, consider $\mu^0 = \l_2^0$. Suppose $\mu^k\ne 0$, and let $j$ be minimal such that $\mu^k$ has a nonzero coefficient of $\l_j^k$. Then $\l-\b_j^k-\cdots -\b_1^k-\g_k$ affords the highest weight of a summand  
$V_{C^0}(\l_2^0) \otimes V_{C^{k-1}}(\l_{r_{k-1}}^{k-1}) \otimes V_{C^k}(\mu^k-\l_j^k+\l_{j+1}^k)$ of $V^2(Q_Y)$, where once again we omit the term $\l_{j+1}^k$ if $j = r_k.$  The restriction of this to $L_X'$ is $\wedge^2(2\om_1)\otimes  (b\om_1\otimes \om_l) \otimes (\mu^k-\l_j^k+\l_{j+1}^k)\downarrow L_X'$, and the tensor product of the first 3 factors is not MF. Hence $\mu^k=0$. And $x_i = 0$ for all $i$ by the usual argument. Hence $\l=\l_2$, as in the conclusion of the theorem.

\vspace{4mm}
This completes the proof of Theorem \ref{OM1OML}.

\end{document}